\documentclass[a4paper]{article}
\usepackage{amsmath,amsfonts,amssymb}
\usepackage{graphicx, epsfig}

\def\captionskipAG{\vskip-10pt\penalty0}
\def\figskip{\vskip-10pt\penalty0}
\def\Figskip{\vskip-0pt\penalty0}

\def\F{Fig.\,}

\newtheorem{theorem}{Theorem}[section]
\newtheorem{Pseudo-Theorem}{Pseudo-Theorem}[section]

\newtheorem{lemma}[theorem]{Lemma}
\newtheorem{cor}[theorem]{Corollary}
\newtheorem{prop}[theorem]{Proposition}
\newtheorem{defn}[theorem]{Definition}
\newtheorem{ques}[theorem]{Question}
\newtheorem{conj}[theorem]{Conjecture}

\newtheorem{rem}[theorem]{Remark}

\newtheorem{prob}[theorem]{Problem}
\newtheorem{Numerology}[theorem]{Numerology}

\newtheorem{Scholium}[theorem]{Scholium} %%%Scholie in French
\newtheorem{Notation}[theorem]{Notation}

\newtheorem{quota}[theorem]{Quote}
\newtheorem{principle}[theorem]{Principle}

\newenvironment{proof}[1][Proof]{\textbf{#1.} }
{\hfill\rule{0.5em}{0.5em}\medskip}
\newenvironment{proof*}[1][Proof]{\textbf{#1.} }{}

\def\epsilon{\varepsilon} %%%%%%epsilon
\def\loccit{{\it loc.\,cit.}}
\def\RR{{\Bbb R}} %reals
\def\CC{{\Bbb C}} %complexes
\def\PP{{\Bbb P}} %projective space

\def\ZZ{{\Bbb Z}} %ganze Zahlen

\def\disc{{\frak D}} %discriminant

\def\mH{{\vert mH \vert}} %hyperspace of curves

\def\la{\langle} %%%%for the Viro symbols (but Gudkov better?)
\def\ra{\rangle}
\def\vc{\sqcup}

\def\al{\alpha}
\def\be{\beta}
\def\ga{\gamma}
\def\de{\delta}

\def\ep{\varepsilon}

\def\lam{\lambda} %%%%%%NOT TO BE CONFUSED WITH leftangle

\def\Qbar{\overline{\Bbb Q}}

\begin{document}

%%%%%%%%%%%%%%%%%%%%%OLD TITLE
%%%%%%%%%%%%%%%%%%%%%%%%%%
%\title{Ahlfors circle maps:
%historical ramblings}

%%%%%%%%NEW TITLE %%%13.03.13
%\title{Ahlfors circle maps and
%total reality:\\
%from Riemann to Rohlin, plus Viro}

%%%%%%%%NEW TITLE %%%13.03.13
\title{Notes on Hilbert's 16th:\\
%%%problem:\\
%%an experimental account of
experiencing Viro's theory}

%%%%borrowed by Abhyankar 1976
%excursion errante, divagation (vagabonde et décousue)

\author{Alexandre Gabard}
\maketitle

\vskip-15pt

\newbox\quotation
\setbox\quotation\vtop{\hsize
%%6.8cm
7.8cm \noindent

\footnotesize {\it Nihil est in infinito quod non prius fuerit
in finito.}

\noindent Andr\'e Bloch 1926
%in {\it of his papers
{\rm \cite{Bloch_1926_EM}, \cite{Bloch_1926_Mem}}.

\iffalse
%
CENSURE

%%%addded 25.07.13

\footnotesize {\it Ruhe giebt es nicht, bis zum Schlu{\ss}.}
%

\noindent Klaus Mann, Sohnemann von Thomas Mann (ca. 1939
%in {\it of his papers
{\rm \cite{Mann_1939_EM}}.

%%%%%%%added 08.08.13

\footnotesize {\it Ich bewundere dich, wie du flei{\ss}ig weiter
dampfst. (Sprichwort des Sohnemannes Alexander an der Mutter
Christa, during Juli-August 2013, eingeschrieben im Computer am
08.08.13, 02h48 AM=am Morgen!)}
%

\fi

%\smallskip
%{\it Ich musste die ver\"ucktesten Kurvenscharen integrieren.}

%\noindent Oswald Teichm\"uller writing to his Ph. D. student.

%\smallskip
%{\it Is there a life without a metric?}

%\noindent Misha Gromov in Spaces and questions, 1999.
}

\hfill{\hbox{\copy\quotation}}

\medskip
\newbox\abstract
\setbox\abstract\vtop{\hsize 12.2cm \noindent

%%%%%%%%%%%CLEAN ABSTRACT (ARXIV)

%%%\footnotesize
%\small
\noindent\textsc{Abstract.} This text is intended to become in the
long run Chapter~3 of our long saga dedicated to Riemann, Ahlfors
and Rohlin. Yet, as its contents evolved as mostly independent
(due to our inaptitude to interconnect  both trends as strongly as
we wished), it seemed preferable to publish it separately. More
factually, our account is an attempt to get familiarized with the
current consensus about Hilbert's 16th in degree 8. This is a
nearly finished piece of mathematics, thanks heroic breakthroughs
by Viro, Fiedler, Korchagin, Shustin, Chevallier, Orevkov, yet
still leaving undecided  six tantalizing bosons among a menagerie
of 104 logically possible distributions of ovals (respecting
B\'ezout, Gudkov periodicity, and the Fiedler-Viro imparity law
sieving away 4+36 schemes). This quest inevitably involves
glimpsing deep into the nebula referrable to as Viro's
patchworking, and the likewise spectacular obstructional laws of
Fiedler, Viro, Shustin, Orevkov. In the overall, the game is much
comparable to a pigeons hunting video-game, where 144 birds are
liberated in nature, with some of them strong enough to fly higher
and higher in the blue sky as to rejoin the stratospheric paradise
of eternal life (construction of a scheme in the algebro-geometric
category). Some other, less fortunate, birds were killed (a long
time ago) and crashed down miserably over  terrestrial crust
(prohibition). Alas, the hunt is unfinished with still six birds,
apparently too feeble to rally safely the paradise, yet too
acrobatic for any homosapiens being skillful enough to shoot them
down.

 }
%We hope at least
%to have demonstrated that Riemann surfaces is a game with little
%concessions for sane life.}

\centerline{\hbox{\copy\abstract}}

\iffalse
\bigskip

%2000 {\it Mathematics Subject Classification.} {\rm 57N99,
%57R30, 37E35.}

\bigskip\fi

%\normalsize

{\small \tableofcontents}

\section{Hilbert's 16th in degree 8}

[29.09.13] Hilbert's  problem in degree 8
%%%(the 1st unsettled case)
(the present frontier of knowledge)
%$H_{16}(8)$
is a puzzling piece of mathematics much revolutionized by Oleg
Yanovich Viro (1948--$\infty$)
 in the early 1980's,
%%'s revolution
yet by now stagnating along a fairly complex motion of
vortex-tube,
%%%mysterious state-of-affairs
especially when it comes to the six {\it bosons\/} not yet
elucidated. Those elusive objects---not yet detected nor
prohibited---are
%%%%represented
located by
%red-dashed rectangles
black holes of concentric circles on
Fig.\,\ref{SIMPLIFIED-TABLE_gurus-COPY:fig} below, which
%%%concentrates
synthesizes---in one {\it \"Uberblick\/}---all what mankind (and
machines?) knows about the problem (of predicting
%how $M$-octics
%distribute their ovals).
the different shapes of $M$-octics). The puzzle (albeit sembling
close to completion) seems resisting all assaults since circa one
decade, when the last progresses were scored by Orevkov ca. 2002.
It is presently unclear if the field requests completely novel
ideas, or just sharpening old weapons already available.

%%%metaphor of 05.10.13

Hilbert' 16th video-game is much akin to a pigeons or ducks
hunting: from 144 birds left free in nature, many could escape,
while flying higher and higher up to reaching the paradise of
eternal life (when constructed). Some
%few
other (less chanceful) birds  were killed, and felt dramatically
on terrestrial ground (when definitively prohibited). Alas, the
hunt is still open, as it remains in the blue sky, six birds
seemingly too weak to reach safely the paradise, yet too vigourous
 for any homosapiens being  skillful enough to shoot them
down.

Our own level of understanding (after ca. 4 months of work) stays
very low, in part because the combinatorics is
%quite
pretty overwhelming, but also because the main results (especially
those of prohibitive character) are quite hard stuff to digest.
Personally, we confess to have not yet been able
%%%to find enough time and energy
to digest even
%%those of
the first generation (4+36=forty many) prohibitions coming from
the Fiedler-Viro era in the early 80's. Besides, the proofs of
%%several
eight pivotal (sporadic) Viro obstructions (ca. 1984/86) were
apparently never clearly published\footnote{Added in proof
[06.10.13].---According to a recent e-mail by Viro, cf.
Sec.\,\ref{Correspondence-NEW}, it seems that those sporadic
prohibitions (by Viro) were integrated in a paper by Shustin.
Alas, for the moment, we lack a precise reference.}. Thus, a
certain aptitude, somewhat between prediction
%%power
ability \`a la Nostradamus and deep devotion to the
algebro-geometric crystal, seems requested to crack the problem.
This will surely not occur as an isolated prowess
%progress
specific to degree $8$
%%%
%prouesse checke in DICO
%%
but probably as a novel method (constructional or prohibitionary,
possibly a cocktail of both) spreading over all other degrees,
while offering new insights about the enigmatic morphogeneses of
{\it all\/} algebraic curves, conceived as an organic
(indestructible) entity. Besides, through their physical
incarnations, either as trajectories of planetary systems
(Kepler-Newton for conics), caustics in optics
(Huyghens-Newton-Thom, etc.), periodic atomic systems (as posited
in Gabard 2012 \cite{Gabard_2012/13}), or whatever else (spectral
curves in ferromagnetic percolation \`a la
Kenyon-Oukunkov-Cimasoni), algebraic curves seem to interact
widely with the natural world (assuming its existence of course).
%
%field of applications requesting highbrow skills in visualization
%that should now slowly
%%%penetrate the public domain
%become reality after computerization of the world.

% A certain
%aptitude, somewhat between mechanical drawing ability and artistic
%drawing, is required for successful flux plotting, and with
%practice people possessing such aptitude can learn to draw flux
%plots for a great variety of cases with relative ease.

\subsection{List of questions for Russian
or Toulouse
%%experts
colleagues}

[29.09.13] This section gathers some basic questions,  we could
not settle alone. Our
%own
jargon is
hopefully  self-explanatory (otherwise browse quickly through
Sec.\,\ref{jargon:sec}).

$\bullet$ (VSO)=(Viro's sporadic obstructions).---{\it Is a proof
of Viro's eight sporadic obstructions
 available in print?\/} Those were first announced in Viro 1986
 survey \cite{Viro_1986/86-Progress}, but  as far as we know never published in detail. The best
information we could
%%find
glean
%%%%glaner
is from Orevkov 2002
\cite{Orevkov_2001/02-classif-flexible-M-curves-degree-8}, where
it is remarked that the proofs
%of those obstructions
are similar to those implemented in Korchagin-Shustin 88/89
\cite{Korchagin-Shustin_1988/89}. If unavailable in print, who is
competent enough to write them down explicitly?
%%to democratize this knowledge?
Obvious candidates: Viro, Fiedler, Shustin, Korchagin, Orevkov, Le
Touz\'e, who else? Can someone
%take
include this didactic duty in his
%%%top-priority
agenda as to make the technique available to a broader
%class
spectrum of workers? If difficult, is it legitimate
%starting
to start doubting about all (or at least some) of those sporadic
prohibitions?
%, or are they so clear-cut to be hundred percent sure
%terrain?

\begin{figure}[h]\Figskip
%\vskip-1.2cm\penalty0
%\centering
\hskip-3.7cm\penalty0
\epsfig{figure=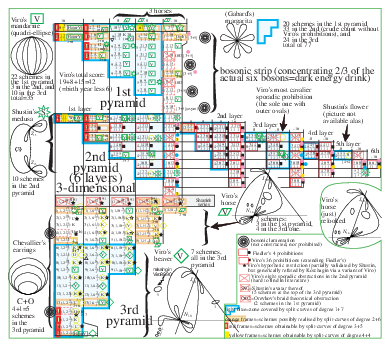,width=192mm}
\captionskipAG
  \caption{\label{SIMPLIFIED-TABLE_gurus-COPY:fig}%
The periodic table of elements (of all octic $M$-schemes).
%%%
\iffalse The game is a akin to pigeon or duck hunt: from 144 birds
left free in nature, many could escape and fly up to reach the
paradise to live for eternity (when geometrically constructed),
some few other (less chanceful) birds  were killed, and falling
dramatically on terrestrial ground (whence definitively
prohibited). Yet, the hunt is not finished, it remains in the sky,
six birds as yet to weak to fly up to the high stratospheric
paradise, and still to vigourous for any homosapiens being able to
kill it.
%
\fi
%
} \figskip
\end{figure}

If politically correct, the present (mostly Soviet) census reports
83 octic $M$-schemes constructed (by the following authors listed
chronologically, plus counting their prolixity: Harnack 1876=Ha=2,
Hilbert 1891=Hi=4, Wiman 1923=W=1, Gudkov 1971=G=2, Korchagin
78=K78=1, Viro 80=V=42(!), Shustin 87/89{\allowbreak}=S=6+1=7,
Korchagin 89=K=19, Chevallier 02=C=4, Orevkov 02=O=1) and 15
prohibited (Viro 84/86=V=8, Shustin 89=S=5 , Orevkov 99/02=O=2).
This leaves {\it six\/} cases undecided among the universe of 104,
as it stabilized after the Fiedler-Viro law of imparity (reigning
in the trinested case). This striking result  is quite akin to a
2nd law of thermodynamics, being the {\it first\/} nontrivial
prohibition beyond  Gudkov-Rohlin periodicity. Without
Fiedler-Viro, the universe would include 144 schemes (like the
Tupolev), all respecting Gudkov periodicity. All
%%of our
%%human
terrestrial knowledge in degree 8 is summed up in the following
catalogue (Fig.\,\ref{SIMPLIFIED-TABLE_gurus-COPY:fig}), whose
ground architecture (mostly governed by Gudkov periodicity) is
stable despite possible errors in the decoration of the building.

%%{\footnotesize
{

{\it Notational conventions}.---Each symbol of the table encodes a
distribution of 22 ovals (the maximum possible in degree $8$)
respecting evident B\'ezout obstructions for line and conics (but
not more!\footnote{Do not forget asking S\'everine (Le Touz\'e) if
her cubics technique for $M$-nonics adapts to octics, and if so
which sort of results does it reproduce.}), and Gudkov periodicity
$\chi \equiv k^2 \pmod 8$, where $\chi$ is the characteristic of
the Ragsdale membrane bounding the curve from inside, and $k$ the
semi-degree (here $k=4$). Our symbolism is basically that of
Gudkov (as opposed say to that of Viro), which in our opinion
overuses brackets, whence a typographical handicap when it comes
to represent all symbols on a single page. It was often
 criticized that Gudkov symbols have the sibylline drawback of
not standing on  a single line. Yet the usual typographical trick
(used e.g. with continuous fractions) permits one to write down
the subnested scheme  also on a single line. For instance, instead
of writing $1\frac{2\frac{17}{1}}{1}$ we may use the condensed
notation $1(1,2\frac{17}{1})$ (also due to Gudkov) as encoding the
distribution with one outer oval, one large oval
%%$(1,\dots)$
enclosing simultaneously 2 ovals plus one oval
%%%enclosing
surrounding 17 mini-ovals. As indicated by our table this scheme
was first constructed by Hilbert (and then re-accessed by Viro's
method).

}

Of course it is hoped, that all those (Soviet) results are
correct. Having personally not yet assimilated all of them,
%%(despite full Leningrad(ian) origins),
we keep open the option of
some mistakes needing
%revisionism.
revision. This is merely a subjective incertitude
%principle
allied
to our own incompetence. Evidently, we have no  serious objections
against the actual consensus, which looks logically robust
(noncontradictory) and plausible   (yet not
%%%omnipotent,
omniscient, hence unsatisfactory). Further, as we already said,
one of the deplorable issue is that some of the most
%%%key
formidable results (especially Viro's eight sporadic laws (=eight
commandments) are not yet available in print. This seems
especially deplorable, as it is the natural sequel to Oleg
Yanovich's heroic saga, implying in particular some limits to the
patchworking method as a purely combinatorial/random fabric
%generator
of curves.
%along any pedestrian
%fashion.
Of course, not all published results are true (nor are
all truths published), yet the public-domain seems a prerequisite
toward checking (resp. assimilating) truths.
%
%
%
\iffalse This looks extremely deplorable, since any (mediocre)
freethinker will quickly founders down
%%%(=sombrer)
along a hibernating digestive process of wasting his energy in
rediscovering the proof of such prohibitions.\fi

\begin{figure}[h]\Figskip
%\vskip-1.2cm\penalty0
%%%\centering
\hskip-2.7cm\penalty0
\epsfig{figure=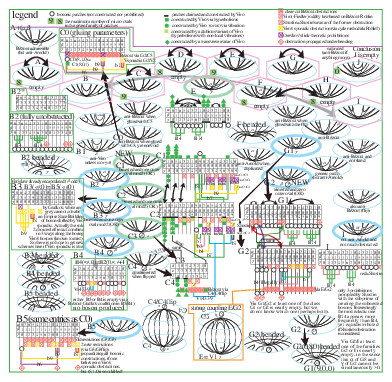,width=172mm}
\captionskipAG
  \caption{\label{ViroDEGREE8_exotic_patches0_BEND-COPY:fig}%
Catalogue of all patches under bending duality} \figskip
\end{figure}

$\bullet$ Patch mirabilis
%%%%CHECKED ON GOOGLE Gauss annus mirabilis
C2$(9,0,0)$.---On studying Viro's method
with extended parameters (counting micro-ovals), one notices the
%%%special
strategical role of the patch C2$(9,0,0)$ (i.e. lateral binested
lune with 9 ovals in the interstice, cf.
Fig.\,\ref{ViroDEGREE8_exotic_patches0_BEND-COPY:fig}). If this
patch exists, then two new bosons are
%%%(miraculously)
(spontaneously) created without any pain. So our question is {\it
whether anybody on this planet ever succeeded to prohibit this
patch.} As a small indicator of the difficulty,
%remark that
neither Viro's imparity law nor Orevkov's deep
obstructions accomplish this job. Cross-link to our
text=Question~\ref{patch-mirabilis-V2(9,0,0):ques}.

UPGRADE [06.10.13].---Actually this patch (C2(9,0,0)) was
prohibited in Shustin's PhD Thesis, compare his recent e-mail in
Sec.\,\ref{Correspondence-NEW} which supplies also another
argument. We acknowledge very much Shustin for communicating us
this precious information.

$\bullet$ (Census of all E-patches for $X_{21}=:F4$.)---In the
same vein there is the question whether the patch family E=V1
(trinested lune, cf.
%%%%%Fig.\,\ref{ViroDEGREE8_exotic_patches0_SYS:fig})
Fig.\,\ref{ViroDEGREE8_exotic_patches0_BEND-COPY:fig}) is complete
under Viro's theory. Here again we are not aware of prohibitions,
leaving open the opportunity to construct both subnested bosons
via an enriched family of patches extending  Viro's. This does not
necessarily mean that there is a method more puissant than Viro's,
but rather that it might be
%%%varied
strengthened along  twists not yet explored (e.g., by futurist
artists of the Moscow school).

We could then solve Hilbert's 16th purely in the vicinity of the
quadri-ellipse at least if Orevkov obstructions are false. (If
Orevkov is wrong, all four binested bosons could be small
perturbations of the quadri-ellipse, cf.
Fig.\,\ref{ViroDEGREE8_extended_BIS:fig}.) In contrast, it may be
true that Viro's actual census of patches has already reached its
%eternal
ultimate crystallization (i.e. no more patches than those
constructed by Viro are available), in which case we really need
Viro's zoo (beaver+horse), Shustin's art (medusa+?), plus
eventually some of your own do-it-yourself creatures to create
additional $M$-schemes. In this scenario, Hilbert's puzzle of
isotopic classification requests insufflating more artistic
freedom (which as we know since Sebastian Bach often reduces to
finitary combinatorics intermingled with sensations of infinity).
Paraphrasing, if Viro's theory of $X_{21}$-patches is already
frozen as it stands, then more flexibility may come from  global
patterns of artwork, rather than through (optical) dissipations of
the (quadruple) rainbow $X_{21}$ as a rigid (unimaginative) form
of patchwork exploiting only the quadri-ellipse.

$\bullet$ (VCP)=(Viro's census of patches is trivially
uncomplete?).---While studying Viro's proof via his own and our
pictures, we noted slight divergences summarizable by saying that
there is in our opinion besides Viro's double-lunes with ovals
injected in the lateral simple lune (=C2 in  our
%catalogue=Fig.\,\ref{ViroDEGREE8_exotic_patches0_SYS:fig}, or
catalogue=Fig.\,\ref{ViroDEGREE8_exotic_patches0_BEND-COPY:fig}
below), also perfect avatars where  micro-ovals are pullulating
instead in the inner simple-lune (cf. class-C1 in the same
catalogue). We would like to ask
%%Oleg Yanovich
Viro (or some other expert), if he believes our patches being also
legal despite not catalogued in Viro 89/90
\cite{Viro_1989/90-Construction}.

$\bullet$ (BS)=Bending symmetry.---A somewhat related question is
an experimental observation that the (extended) catalogue of all
patches looks stable under the symmetry of bending amounting to
invaginate the patch via a motion of horseshoe (compare
Fig.\,\ref{ViroDEGREE8_exotic_patches0_BEND:fig} or its copy right
below=Fig.\,\ref{ViroDEGREE8_exotic_patches0_BEND-COPY:fig}). Can
someone (probably fluent with hyperbolisms or other Cremona
transformations) explain theoretically the presence of such a
symmetry if
%%%it exists
real at all? If yes,  this could be an important tool toward
finishing the exact classification of all $X_{21}$-patches.

{\it Upgrade} [06.10.13]. An e-mail by Viro (cf.
Sec.\,\ref{Correspondence-NEW}) seems to answer this question
completely, via a very simple quadratic transformation.

{\it Prohibitions}.---Another possibility is that all six bosons
(or a portion thereof) will never materialize so that actually no
more constructions are possible, but pure art of prohibition is
requested at some highbrow level of excellence \`a la Viro,
Orevkov, etc. How will those look alike is difficult for us to
predict.

$\bullet$ One may guess elementary techniques of interpolation by
``adjoint'' curves salesman travelling through the deepest ovals
in a sufficiently complicated way as to
%overwhelm
exasperate B\'ezout. This we  call basically the method of DEePest
Penetration (DEPP, like the German word for ``idiot'' or the
US-amerindian actor). Alas, we were not able as yet to implement
it in any successful way, but expect a giant spectrum of
applicability and variability of this method. Very crudely put,
the intuition amounts saying that algebraic curves can nest,
%but
yet the
%complexity
intricacy of the nesting is inherently bounded by the degree (as
we daily experiment with lines and conics). In this respect, it
may be noted that all known prohibitions (Fiedler's four, Viro's
36, Viro's 8 sporadic, Shustin's five, Orevkov's two) all pertain
to curves which are somehow over-nested. For instance, all but one
of Viro's eight  sporadic obstructions pertain to trinested
schemes without outer ovals, and Shustin's five kill exactly the
subnested schemes lacking outer ovals.

$\bullet$ Another promising method (but also miserably ineffective
in our fingers) is that of {\it total reality\/} allied to the
conformal maps of Riemann-Ahlfors. Here, we claim that our short
article Gabard 2013B \cite{Gabard_2013B-Riemann's-flirt} (inspired
by Le Touz\'e) shows how to render purely synthetical (via
M\"obius-von Staudt) the Riemann-Ahlfors map in the schlichtartig
case corresponding to Harnack-maximality (due say to
Schottky-Bieberbach-Grunsky). So we can proudly speak about le
{\it th\'eor\`eme de Riemann rendu synth\'etique}, since the
r\^ole of analysis has been completely
%eliminated for plane $M$-curves.
banished. Alas, it must be confessed that this trivial result is
probably just the first stone toward understanding concretely how
it restricts the distribution of ovals. In the case at hand
(octics), this involves a pencil of sextics  perhaps already
difficult to visualize. Loosely put, the Riemann-Ahlfors map may
be interpreted as a dynamique de l'\'el\'ectron (or
dextrogyration), and so brings a certain dynamical flavor into the
static loci traced by the algebraic orthosymmetric equation. Hence
at the vegetative level at least, one gets another sort of
symbiosis between both parts of Hilbert's 16th. More recently (see
(\ref{fake-medusa:lem})) we noted the usefulness of a sort of
Poincar\'e-Bendixson trapping argument to rule out the {\it
fake-medusa\/} (as a hypothetical avatar of Shustin's medusa,
whose existence would have seriously
%conflicted with
jeopardized Viro's sporadic rules). So again,  there is a germ of
loose overlap between the qualitative theory of differential
equations and the arrangements of algebraic ovals.

$\bullet$ Of course one may also imagine that new prohibitions may
follow patterns of ``old'' methods by Fiedler, Viro, Shustin,
Orevkov, where tools like Fiedler's alternating rule, complex
orientations, clever 2-cycles in $\CC P^2 $, or Orevkov link
theory came to the forefront. And what about Hilbert-Rohn? Is  it
in full desuetude, or does it request a serious rehabilitation as
suggested by Orevkov-Shustin in a relatively recent paper (ca.
2002).  Alas, we confess as yet to have not studied nor
assimilated any of those deep works, and so have nothing more
precise to say.

$\bullet$ (Cubics as tool to prohibit $M$-schemes).---As noted
above, the periodic table of 144 octic $M$-schemes already takes
into account basic B\'ezout obstructions by lines (standard bound
on the depths of nests) and conics (impeding quadri-nested
schemes). One can wonder about the r\^ole of cubics, especially in
reference to S\'everine Le Touz\'e prowess in degree 9
(prohibition of 223 many $M$-nonics, see Le Touz\'e 2002
\cite{Fiedler-Le-Touzé_2002-Cubics}, as well as the subsequent
2009 paper). So the question is whether cubics afford novel proofs
of (old) octic prohibitions originally found by Fiedler, Viro,
Shustin, Orevkov. In our opinion, this natural passage to degree 3
does not exclude the option that already conics may supply
prohibitions (old or novel) if one is able to implement colorful
interpolations by conics salesman-travelling between the nests.
Besides, the method of total reality suggests that degree 6 is the
critical degree where to arrange total reality of a pencil. Yet,
as we said, it looks extremely hard to extract concrete
prohibitive statements.

Besides, it seems hard to claim that all constructions have been
explored. Crudely put, we see a confrontation between patchwork
(using the quadri-ellipse) and artwork (using more imaginative
curves as initiated by Viro, Shustin, etc.), the former being
rather local while the 2nd being more global, yet ultimately also
depending on semi-local patches, yet possibly for other
singularities than $X_{21}$.

\subsection{More shameful points that I could not yet clarify}

[02.10.13]

$\bullet$ As yet I never managed to trace (if possible at all?)
the ground curve (if any?) employed by Shustin to create his
last-discovered scheme $4\frac{5}{1}\frac{5}{1}\frac{5}{1}$. If my
question is meaningful (amounting roughly to say that Shustin does
not employ a more general method that the naive version of Viro's
method which I am able to understand, where the Newton polygon is
banished!?), I would be extremely grateful if someone can send me
the photo-portrait of this Shustin curve (in scanned pdf format).

\subsection{Oracle Orevkov: periodicity modulo 3}

$\bullet$ {\it Toward an Orevkov periodicity modulo
3}.---[03.10.13] As already noted a long time ago, but deciphered
more deeply this night (ca. 03h12), it seems that there is a
certain periodicity by 3 for binested bosons. This is primarily
motivated by Orevkov's obstructions of $b3$ and $b6$, where
$bn=1\frac{n}{1}\frac{19-n}{1}$. If we extend the passage from
$b3$ to $b6$ along a progression by 3, we get $b9$ (a boson not
yet constructed, and now posited as nonexistent), and then
$b_{12}=b_{7}$ (another boson that would also dematerialize), next
$b_{15}=b_{4}$ (also a boson not yet constructed, but that we
posit as nonexistent by propagating Orevkov), and finally
$b_{18}=b_{1}$, the last boson in this series, which would also
not exist. The moral is that what looks fairly chaotic in the
bosonic strip of Fig.\,\ref{SIMPLIFIED-TABLE_gurus-COPY:fig}
becomes perfectly regular under the palindromic symmetry allied to
the permutability of both Gudkov fractions due to the evident
shuffling isotopy of both nests.

The sole (very little) trouble is that extending the Orevkov pair
$b3\to b6$ backwardly we get
$b0=1\frac{0}{1}\frac{19-0}{1}=2\frac{19}{1}$, which is
constructed by Viro (via the horse curve). There is surely a way
to exclude this case, as it really lives outside the binested
realm.

On the basis of this experimental observation (extrapolating
widely Orevkov) there is two questions:

(1).---We seems to be in a situation a bit akin to the love-story
Gudkov/Arnold(+Rohlin), where an experimentalist
%observes
detect a periodicity that only abstract freethinkers are able to
%% suppute
reckon as deep 4D-topology. If an Orevkov-Gabard periodicity
modulo 3 governs the extremal bosonic strip of Gudkov's pyramid,
its seems natural to wonder which topology governs it (vague
guess: link theory, slice genuses, and the well-known yoga between
3D- and 4D-topology via membranes bounding the link). In the
situation of an algebraic (plane) curve there must be (starting
with the link of a singularity) varied ways to get such a setting.
More seriously, look also at Orevkov proofs.

(2).---It would be nice if this threefold periodicity (or another
of another period) also holds in the trinested context as to
supply more order in the apparent chaos reigning around Viro's
sporadic obstructions. Here the story may start with Fiedler's
scheme $f6:=1\frac{6}{1}\frac{12}{1}$, where in general we set
$fn:=\frac{1}{1}\frac{n}{1}\frac{18-n}{1}$. Propagating this
modulo 3, gives $f3$ prohibited by Viro sporadic, and forwardly
gives $f9$ (also Viro sporadic), then $f12=f6$ (Fiedler regular),
and $f15=f3$ (again the same Viro sporadic), and finally $f18=b1$.

One objection is that if we start instead with Fiedler's
$f2=1\frac{2}{1}\frac{16}{1}$, propagation modulo 3 gives $f5$
(Viro sporadic), then $f8$ (Fiedler regular), but then $f11$
(constructible by Shustin's medusa), hence periodicity looks
disrupted. Yet, continuing gives $f14$ (Fiedler), and finally
$f17=f1$, which is also constructible by Shustin. So one could
argue that for this initial condition there is no lovely
coincidence with the binested case, hence we cannot posit
threefold periodicity.

However starting with $f2$, we can join the binested realm via a
period of 2 transforming $f2$ backwardly to
$\frac{1}{1}\frac{0}{1}\frac{18}{1}=b1$ which we assume now as
prohibited (by our postulation of threefold periodicity). Next,
propagating $f2$ forwardly by this novel twofold periodicity gives
$f4$, $f6$, $f8$, $f10=f8$, $f12=f6$, $f14=f4$, etc., and those
guys are simply prohibited by Fiedler.

So there is a biperiodic structure with period 2,3 regulating at
least the 1st layer of the trinested pyramid, and this
bi-periodicity explains all of Fiedler and Viro (otherwise
sporadic) prohibitions. To extend the picture one should inspect
the other layers of the trinested pyramid.

We may start with Viro's $\frac{3}{1}\frac{3}{1}\frac{13}{1}$. To
rally (join) the binested realm, suffice to apply 3-periodicity
leading back to $\frac{3}{1}\frac{0}{1}\frac{16}{1}=b3$, which is
Orevkov's anti-scheme.  It may aver useful setting
$v_k(\ell)=\frac{k}{1}\frac{\ell}{1}\frac{19-(k+\ell)}{1}$, yet
let us try to avoid boring notation (where our $v$ obviously
stands for Viro). Our scheme above
$\frac{3}{1}\frac{3}{1}\frac{13}{1}=v_3(3)$ propagates modulo 3 to
$\frac{3}{1}\frac{6}{1}\frac{10}{1}=v_3(6)$ (Viro regular), and
then $\frac{3}{1}\frac{9}{1}\frac{7}{1}=v_3(9)=v_3(7)$ (Viro
sporadic), and next
$\frac{3}{1}\frac{12}{1}\frac{4}{1}=v_3(12)=v_3(4)$ (Viro
regular), and finally $\frac{3}{1}\frac{15}{1}\frac{1}{1}$, where
we rally the 1st layer, exactly at the place of Viro's ``1st''
sporadic obstruction ($f3$). Besides, the symbol may also
bifurcate when operating on the 1st fraction, to
$\frac{0}{1}\frac{15}{1}\frac{4}{1}$ rallying/fitting/joining
thereby the binested realm along prohibited schemes (under our
hypothetic threefold periodicity).

So the philosophy is a bit as follows: as we know Gudkov fourfold
periodicity is governed by fourfold extended (spin) manifolds \`a
la Rohlin, while it seems that on the boundary case of the Gudkov
pyramids there is reigning a periodicity modulo 3 (suggested by
Orevkov, but also Viro, and even Fiedler), which by analogy might
be governed by 3D-topology.

We may now experiment 3-fold periodicity (3-periodicity) higher in
the telescopic layers forming the trinested pyramid. Specifically,
we may start from $\frac{4}{1}\frac{4}{1}\frac{11}{1}$ (which is
Viro regularly prohibited). If we move ``backwardly'' along 3-fold
periodicity may give Shustin's (constructible) scheme
$\frac{1}{1}\frac{7}{1}\frac{11}{1}$, and one gets a bad feeling
of periodicity breaking. But as we know this Shustin's scheme is
not in phase with the binested obstructions, and therefore our
starting place ($\frac{4}{1}\frac{4}{1}\frac{11}{1}$) is not
adequate to 3-fold periodicity. Actually, it must be subsumed to a
2-fold periodicity explaining all nearby prohibitions, mostly of
Viro's regular sort.

So where is the next place to look for 3-fold propagation? If we
start with say $4\frac{2}{1}\frac{3}{1}\frac{10}{1}$, we get
$4\frac{5}{1}\frac{0}{1}\frac{10}{1}=5\frac{5}{1}\frac{10}{1}$
(constructed by Viro) or
$4\frac{2}{1}\frac{0}{1}\frac{13}{1}=5\frac{2}{1}\frac{13}{1}$
(likewise constructed by Viro). This is troubling but maybe again
explainable by arguing that our starting position is not adequate
for triple periodicity (triperiodicity).

Of course the biggest puzzle is to explain why Viro's most
sporadic obstruction $4\frac{3}{1}\frac{3}{1}\frac{9}{1}$ ought to
be predicted in terms of the (2,3)-biperiodicity.

[10h43:03.10.13] If we start with
$\frac{5}{1}\frac{5}{1}\frac{9}{1}$ we find under 3-periodicity
$\frac{5}{1}\frac{8}{1}\frac{6}{1}$ (Viro regularly prohibited),
next---while jumping over Shustin's construction
$\frac{5}{1}\frac{7}{1}\frac{7}{1}$---we get
$\frac{5}{1}\frac{11}{1}\frac{3}{1}$ (Viro sporadic). This becomes
next $\frac{5}{1}\frac{14}{1}\frac{0}{1}=
1\frac{5}{1}\frac{14}{1}$ (constructed by Viro). If we operation
modulo 3 on the last numerator of the initial scheme
$\frac{5}{1}\frac{5}{1}\frac{9}{1}$, we also arrive at
$\frac{5}{1}\frac{14}{1}\frac{0}{1}$, which is the same scheme of
Viro. Thus we look again disphased (out phased). So it seems that
the sole way to salvage our postulation of periodicity (while
staying politically correct, i.e. without conflicting with the
actual consensus) is to operate on the 1st numerator to get
$\frac{2}{1}\frac{8}{1}\frac{9}{1}$ (Viro regularly prohibited),
but then again acting mod 3 on the 3rd coefficient bring us to
$\frac{2}{1}\frac{17}{1}\frac{0}{1}$, which is Viro-constructed.
Hence we are again
%out-phased.
out of phase.

Finally, starting with $\frac{6}{1}\frac{6}{1}\frac{7}{1}$
triperiodicity bring us to $\frac{0}{1}\frac{12}{1}\frac{7}{1}$ or
$\frac{0}{1}\frac{6}{1}\frac{13}{1}$ (when trading in favor of the
3rd numerator), and now we are in good prohibitive phase.

So what is the conclusion? Obviously our understanding of
triperiodicity is nor perfectly translucid, yet this is maybe our
due to our own stubbornness.

Another related way to pose the question of a periodicity would be
to define another invariant $\varphi$ beside the Euler-Ragsdale
characteristic $\chi=p-n$, somehow transverse to the latter in the
sense that it would vary along the vertical strips of the pyramid
(Fig.\,\ref{SIMPLIFIED-TABLE_gurus-COPY:fig}). Then we
%%%would posit
expect $\varphi$ to be predestined  modulo 3, while explaining all
Fiedler, Viro, Orevkov prohibitions (plus all the
dematerialization of all 4 binested bosons).

More naively, we may return to our old viewpoint, yet more
systematically by starting from the bottom while elevating
progressively into higher layers, as opposed to starting from a
random high-position but often landing down along an out phased
regime (jet-lag).

So the story starts in the bosonic strip (binested with one outer
oval) where it seems perfectly coherent to posit triperiodicity.
This postulation---basically anchored in Orevkov---would kill all
four bosons. Next we move to the 2nd pyramid, starting with
Fiedler's (anti)-scheme $\frac{1}{1}\frac{2}{1}\frac{16}{1}$.
Modulo two, while acting on the 2nd coefficient, this may be
connected to the boson $b1$, and also be dragged down as to cover
all four Fiedler's prohibitions. Besides, acting on the 1st and
3rd coefficient gives $\frac{3}{1}\frac{2}{1}\frac{14}{1}$
(aV:=anti-Viro regular, i.e. imparity law),
$\frac{5}{1}\frac{2}{1}\frac{12}{1}$ (aV),
$\frac{7}{1}\frac{2}{1}\frac{10}{1}$ (aV),
$\frac{9}{1}\frac{2}{1}\frac{8}{1}$ (aV),
$\frac{11}{1}\frac{2}{1}\frac{6}{1}$ (aV),
$\frac{13}{1}\frac{2}{1}\frac{4}{1}$ (aV), and finally
$\frac{15}{1}\frac{2}{1}\frac{2}{1}$, sweeping thereby the full
extreme-right row of the 2nd layer (in accordance Viro regular).
Next, turning back to $\frac{3}{1}\frac{2}{1}\frac{14}{1}$, we may
move right  to $\frac{3}{1}\frac{4}{1}\frac{12}{1}$ (aV), and then
$\frac{3}{1}\frac{6}{1}\frac{10}{1}$ (aV),
$\frac{3}{1}\frac{8}{1}\frac{8}{1}$ (aV), followed by a
palindromic repetition. Further, from
$\frac{3}{1}\frac{4}{1}\frac{12}{1}$ we may expand the prohibited
territory right to $\frac{5}{1}\frac{4}{1}\frac{10}{1}$ (aV), and
its downwards companion $\frac{7}{1}\frac{4}{1}\frac{8}{1}$, or
alternatively climb further to
 $\frac{5}{1}\frac{6}{1}\frac{8}{1}$ (aV), from which place it remains
 only the possibility to reach the very summit of the pyramid with
$\frac{7}{1}\frac{6}{1}\frac{6}{1}$ (aV). In conclusion we got a
big armada of obstruction explained by periodicity mod 2 starting
from Fiedler's 1st scheme $\frac{1}{1}\frac{2}{1}\frac{16}{1}$
(aF), in turn allied to the first boson $b1$.

By analogy starting from $b1$, triperiodicity explains Orevkov
(plus killing all binested bosons), and then elevates to the 1st
layer as $\frac{1}{1}\frac{3}{1}\frac{15}{1}$ (aVs=Viro sporadic
obstruction), and then propagating properly in this layer, while
jumping correctly over the construction of Eugenii Shustin. The
next sort of transformation mod 3 makes the move
$\frac{1}{1}\frac{3}{1}\frac{15}{1}\mapsto
\frac{4}{1}\frac{3}{1}\frac{12}{1}$ from which position we may
explain many prohibitions of the right-row of the 3rd layer. At
this stage, we traced orange versus lilac bubbles on
Fig.\,\ref{SIMPLIFIED-TABLE_gurus-COPY:fig} as to show the
propagation of prohibition under resp. double and triple
periodicity (with initial condition $b1$). This permits to keep in
memory schemes already prohibited, and we note that only Viro's
$\frac{3}{1}\frac{5}{1}\frac{11}{1}$ is missed by our
biperiodicity. Finally, via triperiodicity we may climb to the 4th
and even 6th layer via the moves
$\frac{3}{1}\frac{4}{1}\frac{12}{1} \mapsto
\frac{4}{1}\frac{6}{1}\frac{9}{1}\mapsto
\frac{6}{1}\frac{6}{1}\frac{7}{1}$.

All this is pretty coherent and  noncontradictory with the actual
census, but alas some prohibitions are missed by our recursive
bi-propagation, namely $\frac{1}{1}\frac{5}{1}\frac{13}{1}$ in the
1st layer, and some other few schemes in the higher layers
(compare the pyramid-figure to get the exact enumeration). As just
said, we miss $\frac{1}{1}\frac{5}{1}\frac{13}{1}$ and it seems
natural getting it by the principle of reductionism to the 1st
boson $b1$. This forces introducing another fivefold periodicity
effecting the move $b1=\frac{1}{1}\frac{0}{1}\frac{18}{1}\mapsto
\frac{1}{1}\frac{5}{1}\frac{13}{1}$. On the diagram (always in
reference to Fig.\,\ref{SIMPLIFIED-TABLE_gurus-COPY:fig}), this
prompt introducing a third blue color for this fivefold
periodicity. (Of course one is pleased to appeal to 5 after 2 and
3 since we seems moving along the natural sequence of prime
numbers, aber Hallo!) At any rate, five-propagation along the 1st
layer looks in phase with the consensus. Actually, we may climb to
$\frac{3}{1}\frac{6}{1}\frac{10}{1}$ which maybe altered to
$\frac{3}{1}\frac{11}{1}\frac{5}{1}$ (which we missed as yet).
From either one of those positions we may elevate to
$\frac{8}{1}\frac{6}{1}\frac{5}{1}=\frac{5}{1}\frac{6}{1}\frac{8}{1}$
(alas not new but already covered by 2-periodicity).

As yet we still do not have a completely recursive law explaining
all Viro prohibitions (we miss
$\frac{4}{1}\frac{4}{1}\frac{11}{1}$ and
$\frac{5}{1}\frac{5}{1}\frac{9}{1}$, not to mention the sporadic
obstruction with outer ovals
$4\frac{3}{1}\frac{3}{1}\frac{9}{1}$).

Meanwhile, we may ask if we overlooked the period four, which
albeit not a prime number seems in accordance with the thesis of
reductionism to the 1st boson $b1$, in view of Fiedler's 2nd
prohibition, i.e. $\frac{1}{1}\frac{4}{1}\frac{14}{1}$. Ouh sorry,
of course such a periodicity will be coarser than that by 2, and
so will add no novel information.

So the next step is a periodicity by seven, but this looks
incompatible with Shustin (construction of
$\frac{1}{1}\frac{7}{1}\frac{11}{1}$).

Of course we could inspect higher primes period, but  first make a
puzzling remark. If we look at the schemes not yet covered by our
triple periodicity (with period 2,3,5), we have first
$\frac{4}{1}\frac{4}{1}\frac{11}{1}$. Applying a direct move to
the bosonic strip, we get either
$\frac{8}{1}\frac{0}{1}\frac{11}{1}$, or
$\frac{4}{1}\frac{0}{1}\frac{15}{1}$, which is puzzling as the 1st
is constructed while the other we expect prohibited.

Actually, if we take as ground principle that of starting
systematically from $b1$ we see that there is no periodicity by
seven as it would land on $b8$ (constructed by Viro).

But working so, we get a problem with 2-periodicity already.
Indeed the passage $b1\to b3$ (the latter prohibited by Orevkov)
suggests periodicity by 2, yet when propagated further to b5 this
conflicts with Viro's construction of this scheme. So one must
dictate somewhat artificially an absence of periodicity in the 1st
pyramid.

Next, for $\frac{4}{1}\frac{4}{1}\frac{11}{1}$, we could tabulate
on a periodicity mod 11, transforming it into
$\frac{4}{1}\frac{15}{1}\frac{0}{1}$, which is the boson $b4$.
Further the later is self-dual under 11-periodicity. The moral of
this (and the previous paragraph) it that periodicity might not
always be starting from $b1$, but maybe for certain other initial
conditions.

Last for $\frac{5}{1}\frac{5}{1}\frac{9}{1}$, we may attempt a
reduction to the binested realm by emptying a nest via 5-fold
periodicity. Alas, this leads either to
$\frac{10}{1}\frac{0}{1}\frac{9}{1}$ (boson forbidden by
triperiodicity), or to $\frac{5}{1}\frac{0}{1}\frac{14}{1}$
(constructed by Viro).

In conclusion, we were close to decipher a hidden periodicity
explaining all the apparent reigning around the  Viro, Fiedler,
Orevkov prohibitions while jumping acrobatically over the
constructivist mines posed by Shustin, but alas this turns out to
cover not exactly all cases, while being a somewhat ad-hoc
rafistolage. We expect to address this question more successfully
at the occasion, yet it must be confessed that the overall
approach is not extremely deep, while being logically founded on
the truth of several deep works, we had not yet the occasion to
check the geometric foundation.

[04.10.13] Actually, by our periodicities (mod 2,3,5) we could
explain all prohibitions safe two with $\chi=-16$, namely
$\frac{4}{1}\frac{4}{1}\frac{11}{1}$ and
$\frac{6}{1}\frac{6}{1}\frac{7}{1}$. The former may be excluded as
it is covered by Viro's imparity law which is very regular. So it
seems that in the all ($\chi$ arbitrary) there is only  two Viro
prohibitions not covered by periodicity, namely
$4\frac{3}{1}\frac{3}{1}\frac{9}{1}$ and
$\frac{5}{1}\frac{5}{1}\frac{9}{1}$. Of course one expediting
solution would be that those two sporadic obstructions to be
wrong, yet this request a construction that presently nobody is
able to make.

%%%%%%QUESTIONS FOR VIRO ETC

\subsection{Fixing our jargon}\label{jargon:sec}

[29.09.13] Here we define some terminology and abbreviations of
our own cooking, that we shall constantly use without extra
reference. This requests special boring attention as our
terminology is nonstandard and highly improvised.

$\bullet$ the main contributors to the field are the following
geometers whose name are often abridged by sole initials
especially on combinatorial tables, where little room is left to
write in extenso the original constructor (or prohibitor):
Harnack=Ha (1876), Hilbert=Hi (1891), Wiman=W (1923), G=Gudvov
(1969--72), K78=Korchagin 1978 (variant of Brusotti), F=Fiedler,
Viro=V, S=shustin, K=Korchagin (later), C=Chevallier, O=Orevkov.
On tables an ``a''-privative, like aB, aG, aF, aV, aVs usually
means anti-B\'ezout, anti-Gudkov, anti-Fiedler, anti-Viro, or
anti-Viro sporadic.

$\bullet$ a {\it nest\/} is an oval (of an algebraic curve or an
abstract scheme) which is nonempty in the sense that inside the
unique bounding disc for the oval there appear other ovals of the
curve. An {\it egg\/} is an oval which is empty, i.e. no other
oval inside of it.

$\bullet$ (uninested, binested, trinested, subnested).---As a
trivial consequence of B\'ezout, schemes in degree 8 can have one,
two, three or four nests, but not more than that. Actually, in
case of four nests the configuration collapse to the doubled
quadrifolium $\frac{1}{1}\frac{1}{1}\frac{1}{1}\frac{1}{1}$. The
corresponding Gudkov symbols are $x\frac{y}{1}$ ({\it
uninested\/}), $x\frac{y}{1}\frac{z}{1}$ ({\it binested\/}),
$x\frac{y}{1}\frac{z}{1}\frac{w}{1}$ ({\it trinested\/}). Besides
an octic scheme is the first degree where apart from the trivial
deep nest case appears the option of a subnest where a little bird
(moineau=Jack Sparrow) constructs a little nest squatting one of a
larger bird (crow). Such schemes we call subnested and the have a
Gudkov  symbol of the form $x(1,y\frac{z}{1})$ where $x$ counts
the outer ovals, the first $1$ stands for the big nest of the
crow, $y$ is the number of {\it big eggs\/} in the big ``crow''
nest, and $z$ is the number of (little) sparrow eggs at depth 2.

$\bullet$ Viro's {\it imparity law\/} (or oddness/oddity law)
refers to the result
%of Viro
ca. 1980/83  (extending an earlier one by Fiedler)  that each nest
of a trinested $M$-scheme contain an odd number of ovals. This is
a spectacular result yet not the {\it dernier mot\/} of the story,
since there is also:

$\bullet$ Viro's {\it sporadic obstruction\/} of eight schemes
first listed in Viro's seminal survey on progresses over the past
six years. Alas, those do not seem to have been explicitly proven
in literature and therefore we deliberately---to accuse this
lamentable state-of-affairs---adhere to a doubtful attitude
against them, whence our terminology {\it sporadic}. To caricature
at the extreme (at the risk of being unfair), one could say that
Viro stated those as prohibitions not so much because he was able
to prohibit them, but merely because he was not able to construct
them. As a sibylline avatar, remind that Viro apparently baffled
%%%% dejouer derouter
 himself when
conjecturing nonexistence of many octic $M$-schemes subsequently
constructed by Korchagin. Our ``apparently'' is just an
incertitude principle allied to Korchagin's construction, which we
could not follow.

$\bullet$ (Boson).---A {\it boson\/} is one of the six $M$-schemes
in degree 8 not yet known to exist. This ignorance is a dramatic
%feature
cliff between human brains versus arithmetical (capitalistic)
unpitying
%unpitoyable
law of higher computation. Nonetheless the game looks still worth
studying as algebra seems to embody a principle of economical
depiction of marvellous algebraic drawings. In fact 4 bosons are
binested and two are subnested. So sometimes we abuse terminology
by calling boson, any binested scheme. This lack of imagination
from our side allied to a parsimony of jargon should not cause any
confusion.

$\bullet$ (Macro, quantum and micro-ovals).---We distinguish along
any method of patchworking three kind of ovals: macro-ovals
visible on the curve right after the dissipation, quantum-ovals
placed on the singular artwork but whose exact location is not
exactly known, and micro-oval which are those pullulating in the
vicinity of the singularity damping during the dissipation
process.

$\bullet$ (Quadri-ellipse).---Viro's method of gluing
(patchworking) primarily involves a very basic singular octic
composed of four ellipses each pairwise bitangent at the same two
points.

\subsection{A loose idea of maltese singularities}

[30.09.13] So for instance new schemes not realized (realizable?)
nearby the quadri-ellipse (F4+F4) went constructed via Viro's
beaver (O5+F3), horse (O5+F3), and Shustin's medusa (C4+C4) (cf.
Fig.\,\ref{SIMPLIFIED-TABLE_gurus:fig}). Those employ different
kind of singularities catalogued by Arnold, but we use here or own
coding where e.g. O5 means an ordinary quintuple point with five
distinct tangents, F$k$ means a ``flat'' point (or rainbow) with
$k$ branches entertaining 2nd order tangency between themselves,
and C4 is a candelabrum consisting of F3 plus a fourth transverse
branch.

As a striking example of ``artwork'' Shustin's medusa gave six new
$M$-schemes (in degree $m=8$), and one may wonder about the
existence of other singular curve creating the bosons not yet
known to exist. Reading all of our text there should be candidate,
yet we should at the occasion make a census of all candidate. For
the moment we just mention one example of the margarita curve
(Fig.\ref{SIMPLIFIED-TABLE_gurus:fig}) with a septuple ordinary
point O7, which could create the boson $14(1,2\frac{4}{1})$ upon
using a suitable affine $M$-septic.

As it stands, the example of Shustin's medusa raises the question
if the more transverse version of the candelabrum where two flat
points (F2) of order 2 cross transversally (maltese cross) also
lead to new $M$-curves. First, by comparison to the case of F4 and
C4 where the number of micro-ovals pullulating is three unit less
that the numbers of crossings of a generic perturbation of the
singularity, we infer (loosely?) that the maltese M4 should
produce 5 micro-ovals. Next on using the singular Harnack bound we
see that a curve with two M4 (each eating eight units to the
genus) have 5 quantum ovals (whose exact exact location is not yet
known.) After few trials, we also arrived at the shape of the
medusa, yet using different singularities. Then maybe the
obstruction plaguing fake-medusas (see \ref{fake-medusa:lem})
vanish by chance, and so an opportunity to corrupt Viro's sporadic
rules. Specifically, arranging the quantum ovals as 1+2+2
(central+peripheral+peripheral) we get the singular octic termed
the langouste, whose dissipation (admittedly improvised along mere
arithmetic Gudkov periodicity) could create some of Viro's
sporadic obstructions. In defense of Viro, it may however noticed
that different pasting yields schemes violating the more
established Viro (regular) imparity law. So our construction is
only a pseudo-counterexample to Viro, but it seems to us still
interesting as the trapping obstruction disappears in contrast to
the fake-medusa based on 2 candelabrums C4 (see again
(\ref{fake-medusa:lem})).

\begin{figure}[h]\Figskip
%\vskip-1.2cm\penalty0
%\centering
\hskip-1.7cm\penalty0
\epsfig{figure=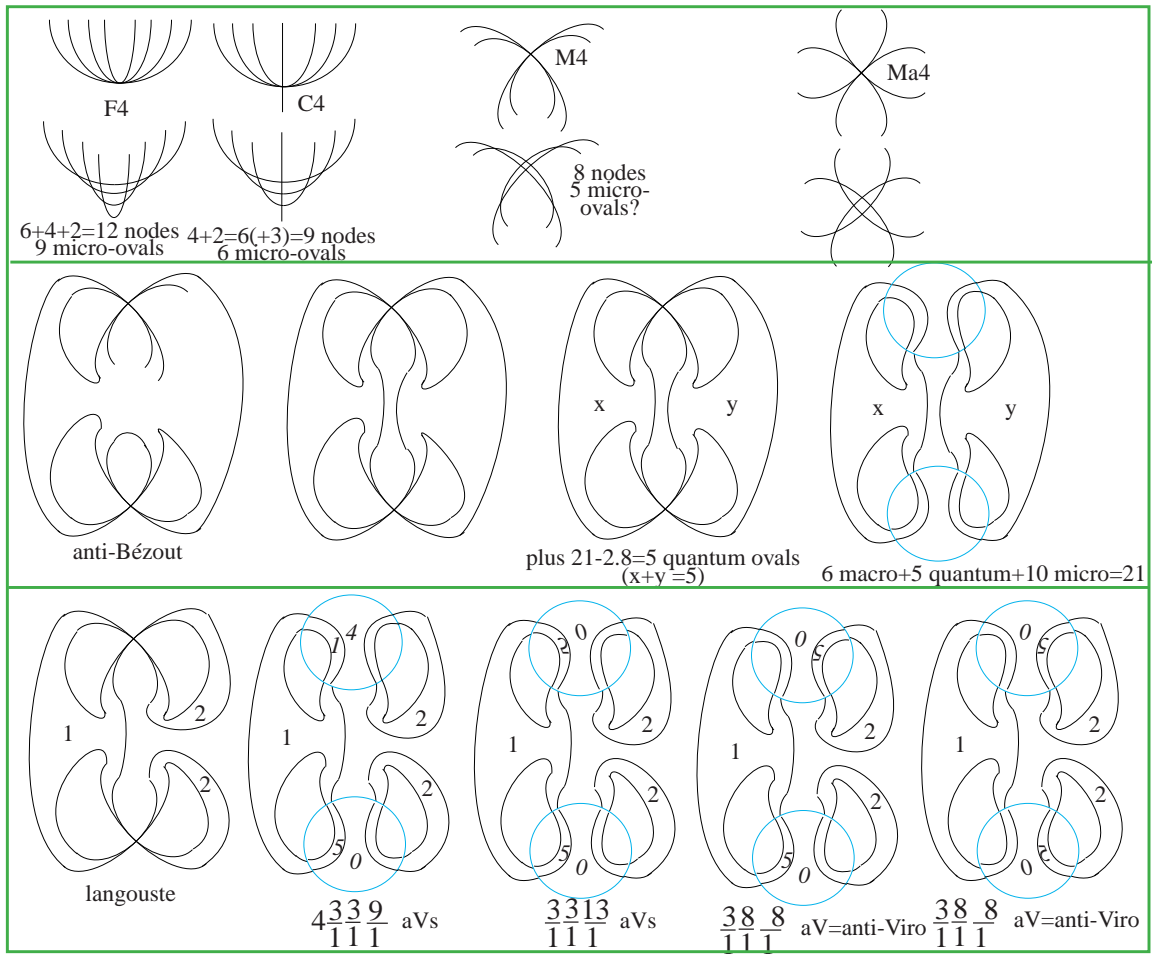,width=152mm} \captionskipAG
  \caption{\label{ViroDEGREE8_maltese:fig}%
  Maltese as a transversalization of Shustin's medusa}
\figskip
\end{figure}

$\bullet$ Another loose question concerns the ubiquity of the
patchwork method as we conceive it. In a perhaps limited sense
(yet broad enough to include all construction by Viro, and
Shustin) we may interpret patchwork as the data of a guru (GroUnd
singulaR cUrve), e.g. quadri-ellipse, etc., plus some patch
prescriptions yielding an $M$-scheme by gluing the patches in
place of the singularities. One naive question is whether all the
more recent constructions of Korchagin, Chevallier, Orevkov (using
as a rule more the story of the Newton polygon) are likewise
interpretable in our narrow sense. If so, make explicit in each
case which is the GroUnd singulaR cUrve (GURU).

\subsection{Correspondence}\label{Correspondence-NEW}

$\bullet$ [03.10.13] Some news from Alex and question on the
Hilbert-Viro problem in degree=8

Dear Colleagues,

I was much fascinated to explore during now ca. 4 months the
fantastic achievements of Oleg Yanovich  and all the other experts
around it. I focused especially on Hilbert's problem in degree 8,
but feel now somehow blocked by certain questions, notably on the
patches for $X_{21}$ (quadruple point with 2nd order tangency as
appearing on the quadri-ellipse). I send you only the new portion
of my text as it is already quite heavy. (418 pages and many
figures: this material is intended to be v3 of my long arXiv
article on Ahlfors, Rohlin and now Viro).

In Sec. 2.2 there is a list of questions (summarized below) which
I would be very happy to submit to your attention, in case you are
not overwhelmed by other duties. Alas, on my side I must start
anew a boring editorial duty for our Math. Journal in Geneva for a
period of about 2 weeks.

Maybe the most clear-cut question I have is whether there is a
complete classification of all dissipations of $X_{21}$, or at
least more detailed information than I was presently able to
compile. Referring for concreteness to my Figure 2 on page 7 (the
catalogue of ``all'' graphically possible  patches), I think that
I found with the class C1 nearby the center of that plate ``new''
patches that were not explicitly mentioned in the article Viro
1989/90 (Leningrad Math. J.). Of course those patches are so much
akin to Viro's (C2 in my notation) that they yield no new schemes
(and were therefore probably deliberately omitted by Oleg 89/90
for page-making convenience). Yet I was still wondering if they
really do exist. (I think that they can be constructed by slight
variants of Viro's purely geometric methods.)

Somewhat more conceptually, it seems to me that there is a global
symmetry of bending regulating the catalogue of all patches, and
amounting to invaginate a patch along a horseshoe motion inverting
all curvatures so-to-speak. (See for instance F and F bended on
that plate). Observationally, all of the available data (to me,
via the reading of Viro) seem to respect this symmetry. It would
be nice if there is a theoretical justification of this ``hidden''
symmetry (maybe via a sort of hyperbolism?). This bending symmetry
could perhaps help in classifying all patches.

As a more specific question, I wonder if Viro's E-patch involving
a trinested lune (near the middle-top of the figure 2) has been
fully classified (especially since it gives opportunities to
create new bosons=[M-schemes]). Likewise I wonder if anybody ever
succeeded to prohibit the patch C2(9,0,0) [which after Orevkov
seems to be the sole undecided patch]. This patch would permit the
creation of the bosons b1 and b7, where $bn$ is the binested
scheme $1 n/1 (19-n)/1$ having n ovals inside a nest.

Last but not least, it seems to me that Viro's eight sporadic
prohibitions (1984, but first announced in the 1986 survey) have
not yet been fully demonstrated in print. So I wonder if anyone is
still able to write down the details. As a question for
S\'everine, I wonder if those (or other octic prohibitions) follow
from your methods with cubics (as implemented to kill ca. 220
nonics in the 2002 and 2009 papers).

Very finally, being not so skilled with Newton polygons I missed
to understand Shustin's last construction of 4 5/1 5/1 5/1. If
someone know the picture of the singular curve leading to this
scheme, I would be very grateful to receive a photo-portrait of
this curve. Sorry if this question is ill-posed???

 Many thanks for all your attention, and I apologize for disturbing
you with my modest questions. Albeit 250 pages long my text
contains nothing original (apart boring tabulations showing all
possible patchworks). So I merely attempted to get familiarized
with the basic combinatorial aspects of patchworking, yet without
being able to discuss properly the prohibitions. So let me know if
there is any readable source available, especially on Viro
sporadic. (I heard by citation about a Texas survey by Korchagin
ca. 1998, which I could not find as yet. If there is an electronic
copy available, it would be excellent...)

All the best, Alex

PS: I added Korchagin and Chevallier to the mailing list, and send
them my best greetings, while apologizing that my text does not
properly reflect their deep contributions.

$\bullet$$\bullet$$\bullet$ [05.10.13] Viro's answer

Dear Alexandre,

C1 and C2 are diffeomorphic by a quadratic transformation
$(x,y)\mapsto (x,y+ax^2)$ with appropriate $a$. This is local
diffeomorphism. This is probably what you are looking for in your
attempt to understand ``duality''.

Shustin's PhD was devoted to smoothings of X21. So it's better to
ask him. He wrote about ``sporadic'' prohibitions. My prohibitions
have been published in his paper. I have no references right away.

All the best, Oleg

$\bullet$ [06.10.13] Gabard's reply

Dear Oleg and the other colleagues,

Many thanks, Oleg, for your prompt answer about the duality I was
expecting. This certainly explains everything I was looking for.
Of course, there is still some questions about special patches
(like C2(9,0,0)), and so I should really take a look to Shustin's
PhD as you suggest (published in VINITI?). This looks extremely
exciting, but alas I need now to work for two weeks for our
Journal L'Enseign. Math. in Geneva.

I will send to the arXiv my notes today, right after integrating
your answer.

All the very best, Alex

$\bullet$$\bullet$$\bullet$$\bullet$ [06.10.13] Shustin's
brilliant answer:

Dear Alex,

The patch C2(9,0,0) you asked for indeed has boon prohibited in my
PhD thesis by means of a version of the Hilbert-Rohn method.
Another way to prohibit it is (as you mention in the table in page
7) is to glue up it with itself and come to a non-existing curve.
Indeed, scanning the patch by a pencil of real lines through the
singular point one can join each of the 9 odd ovals to its
neighbors and one extreme oval with the interior loop by imaginary
discs. In the pathworking procedure these discs persist, and
hence, in the double cover of the plane branches along the
obtained non-singular curve of degree 8 and type 2+1(19), one
observes a sequence of 37 spheres realizing a sublattice $A_{37}$
in the invariant part of the second homology contrary to the
negative signature $s_{-}=35$. I believe a similar prohibition can
be obtained for the original patch when considering it on $F_2$,
and as I remember, Oleg Viro did this a long ago.

With best wishes, Eugenii

$\bullet$ [06.10.13] Gabard's (modest) reply

Dear Eugenii (and the other colleagues), many thanks for this
beautiful answer I look forward to digest in two weeks, but which
sounds extremely elegant. Meanwhile, my arXiv submission is still
much jeopardized due to the arXiv compilator being very sensitive
to coherent cross-links between all the figures. So I had to work
hard adjusting all this. This difficulty leaves me now the
opportunity to include your letter in the text. I hope my
subsequent attempt to compile the file within the arXiv interface
will succeed.

Many many thanks again to Eugenii for this fantastic answer! All
the best, Alex

\section{Getting started}

{\it Terminology} [28.08.13].---Quite typical to Hilbert's 16th
[in its Russian cultivation] is the issue that several
configurations (of the theory, or the Praxis if you prefer German
realism) are not yet constructed nor prohibited (despite the
intrinsic triviality of the algebraico-arithmetical realm making
the whole ``Zeuthen-Klein-Hilbertian'' question a Godd-given
tautology alike). This puzzling state-of-affairs reminisces the
high-energy quest of fundamental particles in natural sciences
(CERN[=Centre Urop\'een pour la Recherche Nucl\'eaire en
Gasp\'esie], etc., i.e. the US and Chinese competitors if any?).
By analogy, call any {\it undecided\/} distribution of ovals
(i.e., not yet known to exist nor to be prohibited) a {\it
boson\/}. As we shall, explain in the sequel the actual Russian
census posits (or demonstrates? if one is clever enough) that
there is actually (after Orevkov 2002) six $M$-bosons in degree 8,
where $M$ refers as usual to Harnack maximality (in the
stenography of Academician Ivan Georgievich Petrovskii). The
question of knowing of many bosons live at the other levels may be
considered as anecdotic but is probably essentially subsumed to
the $M$-case (granting Hilbertio-Russian-Chevalleresque
superstition of longstanding, where the latter refer to
Chevallier, a well-known expert from Toulouse).

Digressing a bit this invites to the following:

\begin{defn}
A mathematical problem is trivial if it merely requests
immortality of the investigator. For instance, the distribution of
primes is a trivial problem (Eratosthenes's crible), and---by way
of consequence---so is probably Riemann's hypothesis on the zeroes
of the zeta function $\zeta(s)=\sum_{n=1}^\infty n^{-s} $).
Probably Hilbert's 16 th problem is likewise trivial, but as yet
we lack (despite the efforts of Klein, Hilbert, Petrovskii,
Arnold, Rohlin, Fiedler, Viro, Korchagin, Shustin, Gabard???) any
clear-cut algorithm reducing this geometrical story to a boring
matter of arithmetics. Of course  it is also permitted to dream of
a world (\`a la Riemann-Klein-Thom-Gabard) where geometry is
stronger than arithmetics (anti-Gauss-Kervaire, etc.), since the
latter discipline is merely discrete geometry alias combinatorics.
This is why Christian W\"utrich (just to name one among many
anti-geometer) is not the master of the world, as he likes to
joke, about English pseudo-Scholars capitalizing his own
progeniture.
\end{defn}

[03.06.13] This section presents a self-contained essay to reach
the actual frontier of knowledge when it comes to octics which is
the first degree where the problem of the distribution of ovals of
algebraic curves is not yet completely settled.

As we shall see, Viro's method of construction is one of the big
breakthrough allowing one to get fairly close to a complete
understanding. The power of his  method supersedes violently all
what what was possible by older perturbation methods \`a la
Harnack, Hilbert or even Gudkov! Interestingly Viro's method is
not a closed machine but one developing further through the
fingers of other experts like Shustin, Korchagin, Chevallier,
Orevkov. So many new curves were obtained by reworking (or
twisting) Viro's method.

At the level of prohibitions also, Viro's role is again fairly
prominent though being based on ideas of Rohlin, Fiedler, etc.
Developing all this in full details should occupy at least
something like 50 pages, so we ask for the reader's patience. Our
intention is to present the theory in its full details so we shall
start slowly by the trivialities, and progressively try to carry
out the big rocks.

Maybe a fairly original result of us (but based on Shustin) is
Theorem~\ref{RMC:cter-example-via-Shustin} below which seems to
disprove Rohlin's maximality conjecture (even the sense thereof
which remained open after the Shustin/Polotovski disproof).

\subsection{Hilbert's 16th for $M$-octics: Harnack 76, Hilbert 91,
Wiman 23, Gudkov 71, Rohlin 72, Fiedler 79, Viro 1980(=Gold
medal=42 schemes), Shustin 85/87/88 (=Bronze medal=7 schemes),
Korchagin 78/88/89 (Silver medal=20 schemes), Chevallier 02 (4
schemes), Orevkov 02 (1 scheme and 2 prohibited)}
\label{Degree8-M-curve:sec}

[29.04.13] This is a big story, not yet completely elucidated.
Among a menagerie of 104 permissible in the Rohlin-Fiedler-Viro
era of complex orientations, it remains now 6 schemes left
undecided after the last advances due to Orevkov 2002
\cite{Orevkov_2001/02-classif-flexible-M-curves-degree-8}. This
still crystallizes the present frontier of knowledge as about of
2013. This affords a nearly complete solution, yet it is  not
clear how much time consuming it will be until completing the full
programme. The 6 bastions de r\'esistance seems quite hard to
crack (and now resisting 11 years of efforts). Looping back to
Viro 1989/90 \cite[p.\,1126]{Viro_1989/90-Construction}, the
%big
incontestable master of the theory wrote: ``The isotopy
classification of nonsingular real projective algebraic plane
curves of degree 8 has not yet been completed, although it is
reasonable to think that it will be completed within the next few
years. In any case, the last ten years have seen much progress
[starting with Viro 1980], and no diminishing of the intensity of
work on the subject.'' Now in 2013, it seems evident that Viro's
agenda was a bit overoptimistic but in substance he might be right
that we are close to the goal, primarily  thanks to his
revolutionary insights and many aficionados (Shustin, Korchagin,
Chevallier, Orevkov, etc.) that joined the br\`eche which he
created. However at the methodological level, Fiedler, Viro or
Shustin deep prohibitions are often poorly published and some work
is required to make their results more accessible to the grand
public.

Our naive idea could be that the method of total reality
(involving here pencil of sextics could help to crack the
problem). Alas, presently we are not even able to tackle the
Hilbert-Rohn obstruction in degree 6, cf. Gabard 2013B
\cite{Gabard_2013B-Riemann's-flirt} for a failing attempt. However
it may be suspected that this was due to lack of cleverness of us,
while Riemann's method of total reality could be the key to the
problem.

[28.04.13] One get a nearly clean view of what happens in degree 8
(say focusing first on $M$-curves) via Viro's seminal breakthrough
1980 \cite{Viro_1980-degree-7-8-and-Ragsdale}. There a table of 52
isotopy types (=schemes) are effectively constructed by Viro's
method of perturbation of complicated singularities
%(Viro's method[s])
when not already realized by a more ancient device, like Harnack
1876, Hilbert 1891, Wiman 1923, Gudkov 1971, Korchagin 1978. On
this table of 52 schemes the ubiquity of Viro's method is already
demonstrated, since the older methods only realize a marginal
portion of isotopy types, namely each method scores so many
schemes as tabulated below:

$\bullet$ Harnack=2; Hilbert=4\footnote{[29.04.13] In Orevkov 2002
\cite[p.\,726,
table]{Orevkov_2001/02-classif-flexible-M-curves-degree-8} Hilbert
only scores 3 schemes. We do not know who between Viro and Orevkov
is right.}; Wiman=1; Gudkov=2; Korchagin=1; Viro=the rest=42.
(This yields a total of $52$ schemes that were known to be
realized in 1980, among a total of 104 logically possible schemes
permissible by the Gudkov congruence and the advanced
B\'ezout-style prohibition of Viro-Fiedler 1980, cf. Viro 1983/84
\cite{Viro_1983/84-new-prohibitions} for the proof.) To understand
why 104 schemes are logically possible, compare our
Fig.\,\ref{Degree8-M-curve-TABLE:fig}.

This little count explains Viro's prose announcing proudly in 1980
(\loccit): ``In this article we formulate a definitive answer for
$m=7$ and some new results on curves of higher degree. Among these
results are the construction of $M$-curves refuting the well-known
Ragsdale conjecture [6](=Ragsdale 1906 \cite{Ragsdale_1906}), the
realization of 42 new isotopy types of $M$-curves of degree $8$
([only] 10 types were realized earlier), and a theorem on
$M$-curves of degree 8 with three nests that excludes 36 isotopy
types not previously excluded.'' To visualize the 36 schemes
prohibited by Viro (and the 4 prohibited by Fiedler somewhat
earlier) compare our
Fig.\,\ref{Degree8-M-curve-TABLE-FIEDLER_VIRO:fig}.)

All this is excellent but it does not tell really what remains to
be done. After this tour-de-force of Viro (and other workers),
exactly 6 types of $M$-octics remains now undecided.
%(cf. e.g.
%Kharlamov-Viro XXXX(=undated) \cite{Kharlamov-Viro_XXXX-UNDATED}).
It would be interesting to see if those schemes can be prohibited
by the method of total reality. Alas, presently we are not even
clever enough to recover the basic Hilbert-Rohn prohibition in
degree 6 (cf. Gabard 2013B \cite{Gabard_2013B-Riemann's-flirt}).
So the case of degree 8, is the ideal terrain in the long-run to
test our philosophy (sketched in the Introd.) that all
obstructions of Hilbert's 16th can be explained via the method of
total reality (and the felicity of pure orthosymmetry \`a la Felix
Klein and concomitantly Rohlin's maximality conjecture).

In fact on reading better Viro 1980 (p.\,569) one sees that the
obvious restrictions (i)-(ii)-(iii) listed on p.\,568 (namely (i)
the definition of an $M$-curve; (ii) Gudkov's congruence
$\chi\equiv_8 k^2=16\equiv 0$); (iii) the obvious consequences of
B\'ezout's theorem restricting the schemes to one of the following
list (given in Gudkov's notation):
$$
\al\frac{\be}{1}, \quad \al\frac{\be}{1}\frac{\ga}{1}, \quad
\al\frac{\be}{1}\frac{\ga}{1}\frac{\de}{1}, \quad \al (1, \be
\frac{\ga}{1}).
$$
This conjointly with  Viro's extension of Fiedler's prohibition
(cf. (\ref{Viro-Fiedler-prohibition:thm})) stating that in the
case of 3 nests (the 3rd kind listed above) all parameters $\be,
\ga, \de$ have to be odd (if nonzero), leaves 104 logically
possible schemes, of which Viro's method (with forerunners---like
Harnack, Hilbert, Wiman, Gudkov---probably always
%%%encompassed by
subsumed to Viro's method as in degree 6) realize 52 types, i.e.
exactly the half number, so that the {\it question of the
realizability of 52 types remains open.} This was the state of the
art in Viro 1980. Meanwhile  only six $M$-schemes are in suspense.
%
%(cf. again Kharlamov-Viro XXXX \cite{Kharlamov-Viro_XXXX-UNDATED},
%which alas does not supply accurate cross-reference). The same
%information (``only 6 isotopy types are questionable'') is to be
%found in Viro 2008
%\cite[p.\,198]{Viro_2008-From-the-16th-Hilb-to-tropical}, again
%without accurate cross-reference, nor a listing of the 6
%problematic cases.
The key reference, surveying all what was done previously, is
Orevkov 2001/02
\cite{Orevkov_2001/02-classif-flexible-M-curves-degree-8} where it
is also supplied a complete classification in the case of
pseudo-holomorphic $M$-curves.

\begin{theorem}
\label{Hilbert-16th-deg-8-M-curve(Viro-Orevkov):thm} {\rm (Hilbert's 16th nearly settled for $m=8$
safe 6 questionable schemes)}.---Since  2002 (Orevkov) and still
in 2008 (e.g., Viro's Japanese survey {\rm
\cite{Viro_2008-From-the-16th-Hilb-to-tropical}}) (and probably
still in 2013, April) there remains exactly six $M$-schemes of
degree~$8$ (among the $104$ logically possible)
%%%%%MY MISTAKE VIRO IS RIGHT
%[WARNING: in our
%opinion this should be $102$, and not $104$ as stated in Viro
%1980, 1983/84, compare our Fig.\,\ref{Degree8-M-curve-TABLE:fig}]
spoiling the completion of Hilbert's 16th. Those are the following
$6$ schemes called (by us) the Hilbert-Viro bosons which are not
yet known to be prohibited nor
%known
to be constructible:
$$4(1,2\frac{14}{1}), 14(1,2\frac{4}{1}), \quad {\rm and} \quad
1\frac{1}{1}\frac{18}{1},
1\frac{4}{1}\frac{15}{1},
1\frac{7}{1}\frac{12}{1},
1\frac{9}{1}\frac{10}{1},
$$
with respectively $\chi=16$ and $\chi=-16$. In fact the first
three are known to admit a pseudo-holomorphic realization, but the
$3$ remaining ones are more mysterious in this respect. As another
naive remark if we permute both fractions of the last symbol (to
read $1\frac{10}{1}\frac{9}{1}$), we
%%%contemplate
observe a certain arithmetic progression by $3$ unities.

More precisely among the $104$
%%%%%%%%%%%%[AGAIN: rather $102$!]
logically
possible schemes (after B\'ezout, Gudkov-Rohlin, Fiedler-Viro),
the following scorings of schemes were constructed by:

$\bullet$ Harnack=$2$ (1876), Hilbert=$3/4$ (1891) [WARNING: in
our opinion this is $4$, and there is a misprint in Orevkov 02,
but not in Viro 80], Wiman=$1$ (1923), Gudkov=$2$ (1971),
Korchagin=$1+19=20$ (78/88--89), Viro=$42$ (80), Shustin=$7=6+1$
(85/87/88), Chevallier=$4$ (02), Orevkov=$1$ (02); yielding a
total of $2+4+1+2+20+42+7+4+1=83$ schemes which are effectively
%%%known to be realized.
constructed.

On the other hand post Gudkov-Rohlin, and Fiedler-Viro it remained
$104$
%%%%%%%%%%[NO, rather $102$ as already noted!]
schemes and further ``sporadic'' prohibitions were detected by:

$\bullet$ Viro=$-8$ (1984/86 proof unpublished), Shustin=$-5$
(90/91), Orevkov=$-2$ (02). At this stage it remains exactly $6$
schemes left questionable. Compare our Table below ({\rm
Fig.\,\ref{Degree8-M-curve-TABLE:fig}}) for the exact diagrammatic
showing the state-of-the-art at the time of Orevkov 2002, which
still represents the actual state of affairs. Our table is
essentially Orevkov's 2002 table, less some (minor) misprints and
an improved diagrammatic showing also the prohibitions (entirely
due to Fiedler, Viro, Shustin and Orevkov).
\end{theorem}

Before embarking on the proof which is long (ca. 20 pages and not
yet completed) let us philosophize a bit. As joked in Chevallier
1997 \cite[p.\,4--6]{Chevallier_1997-Secteur-et-def-loc}, it is
widely accepted that the case of $M$-curves should govern the
pyramid and mark the completion of the full Hilbert's 16th without
discriminating non-maximal curves (by a simple combinatorial
descent along the pyramid that could be implemented by the
conjectural Itenberg-Viro contraction of empty ovals). It can
however be remembered that (like in degree 6) there are more
curves which are say $(M-2)$-curves, and actually Polotovskii 1983
\cite{Polotovskii_1983---327-(M-2)-schemes-deg-8} exhibited 327
many $(M-2)$-schemes of degree 8 using Viro's method. This should
be compared with our
main-table=Fig.\,\ref{Degree8-(M-i)-curve-TABLE:fig} where from we
totalize 419 many logically possible $(M-2)$-schemes. Indeed
$(M-2)$-elements of the 1st pyramid may be ranged as
$1+3+5+\dots+15+17+19(+1)=(\frac{19+1}{2})^2(+1)=100(+1)$ many
schemes. In the 2nd pyramid we have a complicated count involving
the 6 layers. In the first layer we count
$1+3+5+\dots+15=(16/2)^2=64$ schemes; in the 2nd layer we count
$2+4+6+\dots+12=2(1+2+3+\dots+6)=2[(7\cdot 6)/2]=42$ schemes; in
the 3rd layer we have $1+3+5+\dots +9=(\frac{9+1}{2})^2=5^2=25$;
in the 4th layer we have $2+4+6=12$ schemes; in the 5th layer we
see $3+1=4$ schemes, and finally in the 6th layer contributes for
nothing. Totalizing gives $64+42+25+12+4=147$ schemes. Finally the
3rd pyramid gives $1+2+3+\dots+18=\frac{19\cdot 18}{2}=19\cdot
9=171$ many schemes. Therefore the total of (logically possible)
$(M-2)$-schemes is exactly $101+147+171=419$ and a good portion
thereof is constructed by Polotovskii. As shown on the main table
(Fig.\,\ref{Degree8-(M-i)-curve-TABLE:fig} and zoomed as
Fig.\,\ref{Degree8-(M-i)-curve-TABLE_I:fig}), there are $6+3=9$
many $(M-2)$-schemes which are prohibited by an obstruction of
Viro. It would be an interesting task to know exactly which
$(M-2)$-schemes are realized. This question albeit more massive by
a factor of about 4, is perhaps easier to complete than for
$M$-curves.

Albeit the title of Polotovskii's article contains the (ambitious)
word classification, it is probably not a {\it complete\/}
%%%% classification.
census. So the full Hilbert problem in degree 8 can probably
occupies several decades until being completed. Actually even the
$M$-case is far from settled and progressing very slowly---not to
say stagnating---since 2002 (last advance due to Orevkov). Then
one would like as well a classification according to the types (as
did Rohlin 1978 for $m=6$), this can perhaps take another couple
of decades (or follows instantly as it was the case for $m=6$
where the $(M-2)$-scheme $\frac{5}{1}3$ was actually not easier to
construct than Gudkov's $M$-scheme $\frac{5}{1}5$). More naively
one may wonder about the altitude  $r=M-i$ (number of ovals) which
produces the largest number (bio-diversity) of schemes. For $m=6$,
a look at Gudkov's table
%%%%(Fig.\,\ref{Gudkov-Table3:fig})
prompts
that the record is scored by $(M-2)$-curves (with 9 schemes so
thrice as many as $M$-schemes). This proportion looks nearly
respected when $m=8$. Probably this abundance of $(M-2)$-schemes
is a general feature (at least when $m$ is even) since there is
not reigning the congruences mod 8 of Gudkov-Rohlin resp.
Gudkov-Krakhnov-Kharlamov. On the other hand, 171 realizations of
$(M-1)$-schemes were worked out by Goryacheva-Polotovskii 1985
\cite{Goryacheva-Polotovskii_1985}, cf. also Polotovskii 1988
\cite{Polotovskii_1988---classif-deg-8} for a general survey. We
shall describe some of them just by adapting Viro's construction
to non-maximal dissipation (yet understanding precisely which
values of the accessory parameters of the dissipation is a bit
heuristic in our treatment.)

\smallskip
\begin{proof}[Proof of theorem
\ref{Hilbert-16th-deg-8-M-curve(Viro-Orevkov):thm}] First we
learned from Orevkov 19XX
\cite{Orevkov_19XX/2001-cplx-orient-M-curves-deg-7-BOZONS} that 2
of these six schemes  (resisting to the settlement of Hilbert's
16th) are explicitly listed as $4(1,2\frac{14}{1})$ and
$14(1,2\frac{4}{1})$.
%(compare Fig.\,\ref{Degree8:fig} for their
%depiction and their geography on the pyramid).
Actually at the
time of the cited Orevkov's article 9 bosons were not yet known to
be either realized or prohibited. A naive idea would be to
prohibit them by the method of total reality via a pencil of
sextics as described in Gabard 2013B
\cite{Gabard_2013B-Riemann's-flirt}. Of course this is
presumptuous as we are not yet even able to reassess the
Hilbert-Rohn prohibitions in degree 6 by this method, but the
technique seems worth exploring further.

Then in Orevkov 2002
\cite[p.\,725]{Orevkov_2001/02-classif-flexible-M-curves-degree-8}
a very detailed survey is given involving the following methods of
prohibitions:

$\bullet$ the Gudkov-Rohlin congruence $\chi\equiv_8 k^2$ for
$M$-curves;

$\bullet$ the obvious B\'ezout obstructions (\`a la
Zeuthen-Hilbert);

$\bullet$ Viro's theorem 1980 (published 1983/84
\cite{Viro_1983/84-new-prohibitions}) forcing oddity of the
contents of a trinested $M$-curve
%%%% with 3 nests
(generalizing an
earlier weaker result of Fiedler);

$\bullet$ another subsequent (unpublished) result of Viro (1984),
yet to be found in Korchagin-Shustin 1989
\cite{Korchagin-Shustin_1988/89}. This is first reported in print
in Viro 1986 \cite[p.\,67]{Viro_1986/86-Progress} where we read:
``In 1984 I found a new possibility of obtaining restrictions of
non-topological origin. It is based on the construction of
membranes in $\CC P^2$ with boundaries in $S_P A$. Here I announce
one special restriction obtained by this method.---(4.12) If
$\frac{\al}{1}\frac{\be}{1}\frac{\ga}{1}$ is the real scheme of an
$M$-curve of degree 8 then the triple $(\al, \be ,\ga)$ cannot be
$(1,3,15)$, $(1,5,11)$ [WARNING: here $\al+\be+\ga\neq 19$ so this
should probably be $(1,7,11)$ as I thought first, but in fact it
is rather $(1,5,13)$, DOUBLE-WARNING: the same misprint is
reproduced in Viro 1989/90 \cite[p.\,1126,
5.3.G.]{Viro_1989/90-Construction}, where Polotovskii 1988
\cite{Polotovskii_1988---classif-deg-8} is cited who probably
gives no proof], $(1,9,9)$, $(3,3,13)$, $(3,5,11)$, $(3,7,9)$, or
$(5,5,9)$. [Further] There does not exist a curve of degree 8 with
the real scheme $4\frac{3}{1}\frac{3}{1}\frac{9}{1}$.'' So 8
schemes are prohibited by this Viro method (compare again
Fig.\,\ref{Degree8-M-curve-TABLE:fig} to appreciate their location
or even better look at
Fig.\,\ref{Degree8-M-curve-TABLE-FIEDLER_VIRO:fig}). Alas, Viro
has so many methods on his active, that the term ``method'' looks
very unappropriate, and we shall speak of Viro's (membranoid or
2nd/sporadic) obstruction, to distinguish it from the 1st
obstruction (partially due to Fiedler). While the 1st Viro
obstruction has a clear-cut statement (and published proof), the
2nd Viro obstruction looks more mystical, less available in print,
and looks at first sight
%in Viro 1986 (where it is first formulated)
%alike a
fairly random. However on using the right diagrammatic
(Fig.\,\ref{Degree8-M-curve-TABLE-FIEDLER_VIRO:fig}), we observe
some biperiodic pattern emerging on the right-part where
$\chi=-16$. Alas some 3 schemes constructed by Shustin interrupt
slightly the symmetric reproduction of Viro's 2nd obstruction.
Additionally, Viro's 2nd obstruction also includes one scheme with
$\chi=-8$, namely $4\frac{3}{1}\frac{3}{1}\frac{9}{1}$ and the
latter also causes a {\it brisure\/} of bi-periodicity on the
sub-plate where $\chi=-8$ (cf. again
Fig.\,\ref{Degree8-M-curve-TABLE-FIEDLER_VIRO:fig}). So it is
slightly puzzling to decipher the exact harmony of the geometry
%
%%%({\it Harmonia Mundum\/} \`a la Kepler, CHECK the LATIN)
({\it Harmonices Mundi\/}
%%%Libri\/}
\`a la Kepler).
%%%CHECKED in Buecher die die Welt veraendern
It may be observed that in the Fiedler-Viro/Viro era all
prohibitions were concentrated on 3-nested schemes. The situation
will change slightly with the next contributor, Shustin and also
more recently with Orevkov (compare
Fig.\,\ref{Degree8-M-curve-TABLE:fig}).

$\bullet$ the result of Shustin 1990/91
\cite{Shustin_1990/91-New-restrictions} excluding
$(1,(20-a)\frac{a}{1})$ with $a>0$; a priori this prohibits circa
20 schemes but by virtue of Gudkov's hypothesis imposing
periodicity modulo 4 on  the parameter $a$ this prohibits only 5
new schemes. For their exact geography we refer again to our
Fig.\,\ref{Degree8-M-curve-TABLE:fig}.

$\bullet$ Orevkov's (2002) pseudo-holomorphic prohibition of the 2
schemes $1\frac{3}{1}\frac{16}{1}$ and $1\frac{6}{1}\frac{13}{1}$.
Prior to that work of Orevkov (2002) it remained
%%%% $9$(=nine)
nine $M$-schemes whose realizability was in doubt. This was the
result of centennial efforts involving the studies of:

$\bullet$ Harnack 1876 \cite{Harnack_1876} (construction of 2
schemes all with $\chi=16$, via the so-called Harnack method);
those are the schemes $18\frac{3}{1}$ and $17(1,2\frac{1}{1})$
(make a figure at the occasion). This is just a matter of
extending the classical picture of Harnack that we traced only up
to degree 6. Alternatively one might try to find a realization \`a
la Hilbert which is a quicker method. In degree $m=6$, Harnack's
original method can be dispensed at least for $M$-curves, yet in
degree 8 it is not clear to me if the same ubiquity of Hilbert's
method is also valid. Of course, the schemes in questions are not
explicitly listed in Harnack's paper (1876), but Russian scholars
are generous enough to ascribe schemes to the inventor of the
method. As we shall see both of  Harnack's two schemes will be
phagocytized as very special of Viro's method that we will expose
subsequently.
%so that it is not really worth wasting time with
%this archaic method.

%for the latter compare Fig.\,\ref{HilbGab2:fig} when constructed
%\`a la Hilbert.
%[It seemed to me that at this place both Viro 1980
%\cite[p.\,568]{Viro_1980-degree-7-8-and-Ragsdale} and Orevkov 2002
%\cite[p.\,726]{Orevkov_2001/02-classif-flexible-M-curves-degree-8}
%contain a misprint when presenting this scheme as $\la 17 \sqcup 1
%\la 1\ra \sqcup 1 \la 2 \ra\ra$, but they are certainly right.]

$\bullet$ Hilbert 1891 \cite{Hilbert_1891_U-die-rellen-Zuege}
(construction of 4 schemes, via a variant of Harnack's vibration
known as Hilbert method); this includes the scheme
$1 (1,2\frac{17}{1})$
%(cf. Fig.\,\ref{HilbGab1:fig} or
%Fig.\,\ref{HilbGab4:fig})
and  $17(1,2\frac{1}{1})$
%(cf.
%Fig.\,\ref{HilbGab2:fig})
plus 2 other schemes not constructed in
this text (but which we will recover along Viro's method). Warning
at this place Orevkov's table of 2002
\cite[p.\,726]{Orevkov_2001/02-classif-flexible-M-curves-degree-8}
contains a slight mistake by accrediting
%ascribing
to Viro instead of Hilbert the last mentioned scheme. Compare our
table Fig.\,\ref{Degree8-M-curve-TABLE:fig} or rather Viro's
original table in 1980 \cite{Viro_1980-degree-7-8-and-Ragsdale}
(which alas is far from complete, and contains another little bug,
namely the 8 last schemes of the series with $\chi=16$ or
$p=19,n=3$ are misplaced and should be in the series $p=11,n=11$).

$\bullet$ Wiman 1923 \cite{Wiman_1923} (construction of one
scheme, via a method of his own, namely
$16\frac{1}{1}\frac{1}{1}\frac{1}{1}$).
%%Added 09.05.13
Actually this scheme belongs to a series of $M$-schemes of even
degree $m=2k$ with an especially pleasant distribution of ovals
involving a square of outer unnested $k^2$ ovals, plus a
``triangular number'' of $\frac{(k-1)(k-2)}{2}=1+2+\dots+(k-2)$
many nests of depth 2 (Fig.\,\ref{Wiman:fig}). Of course this
series of curves like Hilbert's series
%(cf. e.g.
%Fig.\,\ref{HilbGab3:fig})
also adds some slight evidence toward
the truth of Ragsdale's conjecture $\chi\le k^2$ for $M$-curves.
Although Wiman's method is surely a jewel it is not worth exposing
here as this scheme will be subsumed to the much more powerful
method of  Viro, who like his predecessor also worked in Upsala
(after leaving Leningrad). Apparently as pointed in Polotovskii
1988 \cite[p.\,459]{Polotovskii_1988---classif-deg-8}: ``Speaking
about classical methods we mean the methods by Harnack, Hilbert,
Brusotti, and Wiman of construction of $M$-curves. All these
methods are based on smoothing of non-degenerate double points by
small perturbations. These methods except for Wiman's method are
organized recurrently, so that they give series of $M$-curves of
degrees increasing as arithmetic progression. [\dots]''.

\begin{figure}[h]\Figskip
%\vskip-1.2cm\penalty0
\centering
%\hskip-0.7cm\penalty0
\epsfig{figure=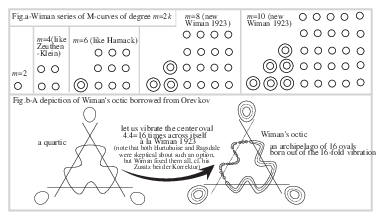,width=122mm} \captionskipAG
  \caption{\label{Wiman:fig}%
  Wiman's series of new $M$-schemes as soon as $m\ge 8$.}
\captionskipAG
\end{figure}

For a less schematic depiction of Wiman's curve compare Orevkov's
picture above (Fig.\,\ref{Wiman:fig}b) borrowed from Orevkov 200X
\cite{Orevkov_2003OR-MORE-Some-examples-Wiman}. On the next page
(p.\,4) Orevkov depicts the next step of a Wiman iteration which
however deviate from our interpretation of Wiman's scheme, so that
perhaps our Fig.\,a is faulty for large $m\ge 10$.

$\bullet$ Gudkov 1971 \cite{Gudkov_1971-const-new-ser-M-curv}
(construction of 2 schemes, via his own method involving Cremona
transformations); namely $14\frac{7}{1}$ and
$13\frac{3}{1}\frac{4}{1}$. (Both those schemes will later be
recovered via Viro's method so that it is not necessary to pay
special attention at Gudkov's realizations.)

$\bullet$ Rohlin 1972 (proof of the Gudkov congruence, with
correction in Marin 1979 \cite{Marin_1979}) ruling out
(statistically) one-quarter of the schemes;
%(compare the pyramid
%Fig.\,\ref{Degree8:fig});

$\bullet$ Korchagin 1978 \cite{Korchagin_1978} construction of one
scheme (namely
$9\frac{1}{1}\frac{10}{1}$, cf. e.g. Viro 1980
\cite{Viro_1980-degree-7-8-and-Ragsdale} table, p.\,568) by a
variant of Brusotti. This will also be subsumed to Viro's method.

$\bullet$ Fiedler ca. 1979 (published 1982/83
\cite{Fiedler_1982/83-Pencil}) (prohibition of 4 schemes, cf. e.g.
Viro 1983/84 \cite[p.\,416]{Viro_1983/84-new-prohibitions})
(published later and englobed in:

$\bullet$ Viro 1980 who establishes an imparity law for 3-nested
$M$-octics, which prohibits 36 additional schemes (proof in Viro
1983/84 \cite{Viro_1983/84-new-prohibitions}). As observed in
\loccit (p.\,416) it seems that Korchagin had some
%pivotal
decisive influence in conjecturing on the basis of specimens
generated by Viro's method the right extension of Fiedler's
obstruction.
%%%that became the Viro's ``odd content'' obstruction.
To visualize the 4 schemes prohibited by Fiedler cf. the hexagons
on our Fig.\,\ref{Degree8-M-curve-TABLE:fig} (or better
Fig.\,\ref{Degree8-M-curve-TABLE-FIEDLER_VIRO:fig} right-below).
To visualize the 36 new obstructions of Viro cf. again on the same
table the squares.
%
%
%
%
%
\iffalse This is alas a bit more tricky on our tabulation. So to
understand this point more geometrically we need to build another
table. For this see rather
Fig.\,\ref{Degree8-M-curve-TABLE-FIEDLER_VIRO:fig} right below: on
it we count only 35 schemes prohibited by Viro and not 36 as
asserted by himself. (Usually in this sort of game Viro is always
right [I often thought at some other places that his count were
sloppy but they always turned out to be right after depicting the
appropriate table], but here I sincerely think that Viro is wrong.
Perhaps his mistake was to count the scheme
$8\frac{3}{1}\frac{4}{1}\frac{4}{1}$ twice by writing it once as
$8\frac{4}{1}\frac{4}{1}\frac{3}{1}$). Finally I realized that
Viro is right when we found one more scheme that we overlooked
namely $12\frac{2}{1}\frac{2}{1}\frac{3}{1}$. \fi

Understanding the Fiedler-Viro proof is fairly tricky. The best
proof is in principle that of Viro 1983/84
\cite{Viro_1983/84-new-prohibitions}, which is however fairly
undigest. If our philosophy of total reality is the right
viewpoint (to obtain even obstruction on some of the six unsettled
cases),  it may perhaps recover the Fiedler-Viro obstruction. In
particular, it may be a wrong idea wasting energy in trying
understanding their proof. As a last remark it seems that Viro's
proof differs from Fiedler's [alternating orientations] in
 using an idea due Rohlin (cf. p.\,414, footnote in
Viro 83/84 and p.\,66 in Viro 86 \cite{Viro_1986/86-Progress}
where Rohlin's formula for a pair of curves is employed). So all
this is a bit tricky and deserves a separate treatment or
approached upon differently via the total reality of $M$-curves.

\begin{figure}[h]\Figskip
%\vskip-1.2cm\penalty0
\centering
%\hskip-0.7cm\penalty0
\epsfig{figure=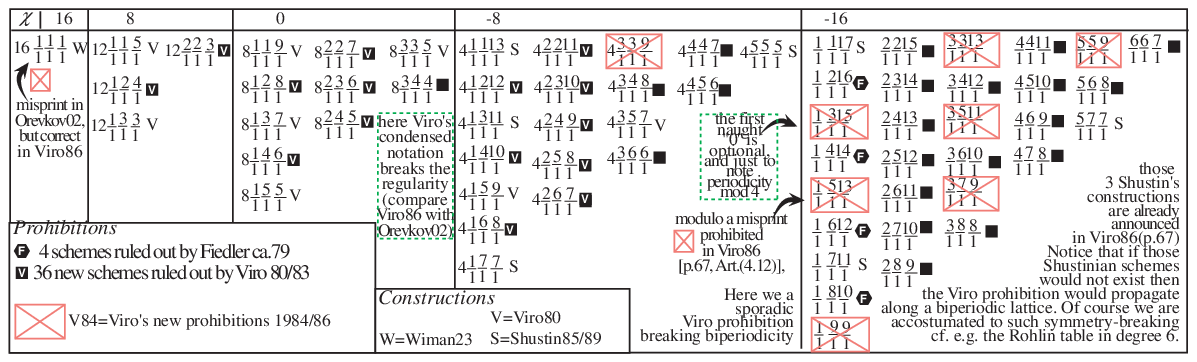,width=122mm}
\captionskipAG
  \caption{\label{Degree8-M-curve-TABLE-FIEDLER_VIRO:fig}%
Diagrammatic of the Fiedler-Viro prohibition for $M$-curves with 3
nests (4 schemes are obstructed by Fiedler 1982/83 and 36 by Viro
1980} \figskip
\end{figure}

$\bullet$ Viro 1980 \cite{Viro_1980-degree-7-8-and-Ragsdale} where
52 types remained open among 104 logically possible (compare our
Table=Fig.\,\ref{Degree8-M-curve-TABLE:fig} to check this count
which is not detailed in Viro's original 1980 paper).
%%%%%%%%%%MY MISTAKE
%[Warning it
%seems that Viro's universe of 104 contains actually only 102
%schemes compare our Table=Fig.\,\ref{Degree8-M-curve-TABLE:fig}].
In this article both the Fiedler obstruction is boosted so as to
exclude 36 types not previously excluded and a revolutionary
method of construction is employed to realize 42 new schemes
(probably encompassing all the previously known constructions but
Viro is modest enough to count just the newcomers); the impressive
list of 42 schemes obtained by Viro is resumed either on our Table
\ref{Degree8-M-curve-TABLE:fig} or on Viro's original table in
1980 \cite{Viro_1980-degree-7-8-and-Ragsdale} (which as we said
contains only a little mistake of 8 schemes which are misplaced in
the series $p=19,n=3$). Viro's table contains additionally the
information of which singularities $X_{21}, J_{10}, N_{16}$ are
dissipated. Some details of Viro's construction are presented in
Viro 89/90 \cite{Viro_1989/90-Construction}. Here we learn the
issue that in contrast to degree 6 where nearly all schemes could
be obtained by perturbing a triplets of 3 coaxial ellipses, the
quadri-axial configuration of 4 ellipses ({\it Viro's earrings\/}
for short) leads only to a special class of $M$-curves and create
``only'' 47 of them. Once we are given the dissipation of the
fourfold $X_{21}$ singularity (4 branches with 2nd order
tangency), compare Fig.\,55 in Viro 89/90, or our figure below
(Fig.\,\ref{ViroDEGREE8:fig}), it is merely a matter of
patchworking to construct the corresponding curves, yet this is so
pleasant that it is worth being published once in full details. As
far as we know this was never exposed in full
%%%details
for
typographical reasons (apart probably in Polotovskii 1988) and so
let us for convenience work out the relevant picture leading to
Viro's breakthrough.

\subsection{Viro's method for $M$-octics}

[04.05.13] First, one of Viro's pivotal idea is to smooth a
configuration of 4 coaxial ellipses
(Fig.\,\ref{ViroDEGREE8:fig}b). Here appears Hawaiian-earing
singularities (or of type $X_{21}$ in Arnold's catalogue of 1975
\cite{Arnold_1975/75-symbols-of-SING-used-in-Viro80}). (Remind
Thom's influence upon Arnold in 1965 and probably also some
overall influence of Klein upon Arnold.) We confess not being
acquainted with Arnold's classification (for the moment), yet this
is no obstacle for understanding the sequel (i.e. Viro's fable).

The next duty is to understand (the maximal) smoothing of this
singularity, those being depicted on Fig.\,\ref{ViroDEGREE8:fig}a.
We see clearly how the 4 branches are reconnected
%%%between
among themselves, while the Greek letters in bracket denotes
pullulation of newly created ovals emerging through dissipation of
the singularity in quantity specified by the table below each
picture. On the right dissipation (V3=Viro3), there is
additionally a little circle which is an enveloping oval. We
personally do not checked (nor do we understand this result), but
again this is no obstacle to understand the sequel. (Full details
seems to be given in Viro 89/90 \cite{Viro_1989/90-Construction}
but the method dates back from the announcement in 1980
\cite{Viro_1980-degree-7-8-and-Ragsdale}.)

The details we give now, albeit elementary, require some tedious
combinatorics that is (never?) user-friendly presented. The idea
is to glue (or patchwork) independently (like in Brusotti 1921)
both singularities so as to create (many) global curves of degree
8. First we can choose any one of the 3 smoothings V1,V2,V3, and
glue it with itself after rotating (and translating) the pattern
by 180 degree. So we get Fig.\,\ref{ViroDEGREE8:fig}c, which are
only $(M-2)$-curves. The reason is fairly simple, namely that 2
large ovals are created while the number of micro-ovals
$\al+\be+\ga$ is always 9, even in the case of
$\delta+\varepsilon=8$ (but keep in mind the micro-circle). To get
better curves we must somehow twist (say) the upper dissipation by
a reflection (symmetry along the vertical axis) to get ``starred''
Viro's smoothing $V1^\ast$, etc., those being depicted on Fig.d.
Then we see 4 large ovals reminding Viro's funny-face (to everyone
knowing  him personally)
%%%or from his web page)
to which are added $2\cdot 9$ micro-ovals reaching therefore
Harnack's bound at $4+18=22=M$. It is then only a matter of
combining all possible smoothings, creating thereby the table
below each picture of Fig.\,d. Here symbols are written along the
Gudkov-Polotovskii as opposed to Viro's symbolism which contains
too much symbols without real significance (like squarecups and
angled-brackets). With Viro's cumbersome symbolism it would never
have been possible to represent everything on a single page as
compactly as we do on Fig.\,\ref{ViroDEGREE8:fig}d.

A minority of those schemes were already obtained by forerunners
of Viro (Harnack76=Ha, Hilbert91=Hi, Wiman23=W, Gudkov72=G,
Korchagin78=K78), but now there is a plethora of new schemes
(marked by ``V'' on the tables).
%when they are really new contributions of Viro 1980).
Those are apparently completely inaccessible to the classical
methods, even when twisted by Cremona transformations like in
Gudkov's trick (fixing Hilbert's 16th in degree 6), or in
Korchagin's variant of Brusotti. Of course we cannot exclude  that
a clever variant of Brusotti being able to create one or two
sporadic schemes, yet Viro's method affords a whole
%%%menagerie
series
%%peloton
of them with
%very
comparatively little efforts.

\begin{figure}[h]\Figskip
%\vskip-1.2cm\penalty0
%\centering
\hskip-1.7cm\penalty0 \epsfig{figure=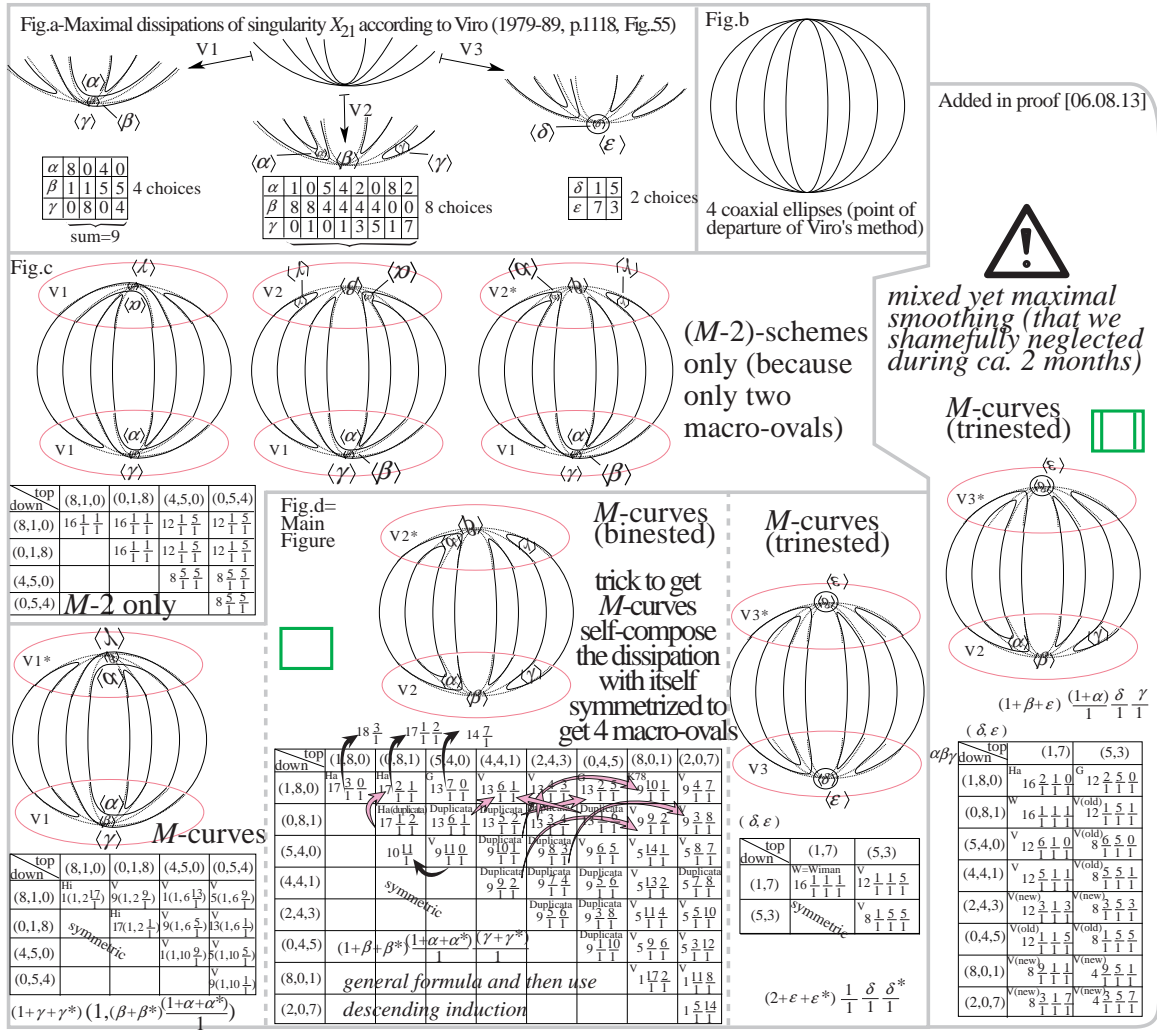,width=152mm}
\captionskipAG
  \caption{\label{ViroDEGREE8:fig}%
  Viro's dissipation (=patchwork) in degree 8}
\figskip
\end{figure}

Working out this table requires some few minutes of concentration.
One trick is to find a general formula for the resulting Gudkov
symbol by contemplating the curve traced on Fig.\,d. After some
few items are calculated, one can propagate the symbols by looking
at increments undergone by the parameters $(\al,\be,\ga)$ and the
story reduces to pure arithmetics without having to refer back to
the picture. This is extremely pleasant to work out and one can
hardly underestimate the level of ecstasy in which Viro must have
been when discovering this ca. 1979/80.

A further pleasant duty is to report all schemes so generated upon
the table (Fig.\,\ref{Degree8-M-curve-TABLE:fig}) by marking with
little green-squared letters ``V'' the schemes so obtained by
Viro. The first curves involving $V1/V1^\ast$ fills $10$ schemes
in the lowest row of the table, the 2nd curves involving
$V2/V2^\ast$ runs through 22 distinct schemes (duplicata being
ignored) in the highest row of the table, while the 3rd curves
involving $V3/V3^\ast$ create only 3 schemes in the
%median
middle row of the table.

At this stage we have realized many schemes (precisely
%%%%% $6
$10+22+3
%%%%%%%%%% =31$
=35$ maximal schemes) arising through perturbation of 4 coaxial
ellipses. This is perhaps worth stating as a separate statement
regardless of the fact that a minority of those
%%schemes
were
obtained by Ha=2, Hi=2, W=1, G=2, K78=1 (yielding a total of $5$
schemes ante-Viro).

\begin{lemma} (Viro 1980).---By dissipating a quadruplet
%of
%$4$
of coaxial ellipses, one can create exactly
%%% 31
$35$ many $M$-curves of degree $8$. This makes precise a bit the
prose in Viro 1989/90 (p.\,1127), ``A very large number of schemes
are realized by our means of small perturbations of the curve in
Figure 72, which is a union of four ellipses having second order
tangency at two points. This curves has two types $X_{21}$
singularities. If we dissipate them using all of the known methods
(see 4.7.A), we can realize $47$ real schemes with $22$[=$M$]
ovals, etc. see Polotovskii [43].''
\end{lemma}

Needless to say we have not yet obtained so many schemes, but only
%%% $31$
$35$ instead of the $47$ many asserted by Viro (on semi-behalf of
Polotovskii 1988). How to explain this gap? Maybe Viro has a
liberal interpretation of the $X_{21}$-singularity in the sense
that it is combined with other tricks (\`a la Gudkov/Newton, i.e.
Cremona or hyperbolism as Viro calls Newton's
%%%% trick).
device).

\section{Flexible exotic patchworking}

\subsection{Bosonic smoothing of the mandarine}

{\it Added\/} [29.07.13] The bosonic strip of doubly nested
schemes with just one outer oval is probably the most mysterious
part of Hilbert's problem for $m=8$. This remains fairly obscure
even after the brilliant interventions of Viro and his disciples,
companions.

It seems of interest---from a naive standpoint at least---to look
at what Viro's method generates when allowing exotic parameters
$(\al,\be, \ga)$ of smoothing. First, in order to land in the
bosonic strip we choose twice $\be=0$ on Fig.\,d. Then we
reproduce the above table yet with extended (unrestricted)
parameters of bubbling $\al, \be=0, \ga$. (The bold faced
characters are the permissible parameters). This gives the large
table of Fig.\,\ref{ViroDEGREE8_BOSONIC:fig} below with obvious
regularity (i.e. each horizontal row is self-reproduced via a
diagonal translation along the South-West direction). All the
other six tables are just replicas of the upper table safe for the
position of the red-crosses. Interestingly the two prohibitions of
Orevkov (namely $1\frac{3}{1}\frac{16}{1}$ and
$1\frac{6}{1}\frac{13}{1}$) forbid conjointly all the crossed
dissipations, i.e. $(7,0,2)$, $(6,0,3)$, $(4,0,5)$, $(3,0,6)$,
$(1,0,8)$, $(0,0,9)$. In some more details we can inspect for each
bosonic curve, which dissipation are killed by an oracle (an
Orevkoracle say) proclaiming the nonexistence of the boson in
question. Actually, the boson $1\frac{1}{1}\frac{18}{1}$ kills
only two dissipations, $(9,0,0)$ and $(0,0,9)$. The pseudo-boson
$1\frac{3}{1}\frac{16}{1}$ (whose nonexistence is due to Orevkov)
kills the 3 dissipations $(7,0,2)$, $(1,0,8)$, $(0,0,9)$. Next,
the boson $1\frac{4}{1}\frac{15}{1}$ kills the 3 dissipations:
$(7,0,2)$, $(6,0,3)$, and $(1,0,8)$. Then, the pseudo-boson
$1\frac{6}{1}\frac{13}{1}$ (ruled out by Orevkov) kills 3
dissipations $(6,0,3)$, $(4,0,5)$ and $(3,0,6)$. Thereafter, the
boson $1\frac{7}{1}\frac{12}{1}$ (positing its nonexistence) kills
the 3 dissipations $(9,0,0)$, $(4,0,5)$ and $(3,0,6)$. Finally,
the boson $1\frac{9}{1}\frac{10}{1}$ (positing its nonexistence)
kills  5 dissipations, namely $(7,0,2)$, $(6,0,3)$, $(4,0,5)$,
$(1,0,8)$ and $(0,0,9)$. It may be observed that nobody succeeds
to kill the dissipation $(5,0,4)$, whose existence would however
not produce new schemes. Also it is quite puzzling to try a
measurement of the eccentricity of a boson via the number of crime
it effects via Viro's method of gluing. Naively they more criminal
a boson is, the more mysterious and difficult to catch it should
be. So perhaps if the criminality index is $\ge 4$, then the boson
does not exist, whereas if it is low $\le 2$ the boson exist.
Perhaps, when this index is 3 then both cases could occur, but we
are rambling into pure speculations due to a lack of geometric
understanding.

It is also tempting to speculate about a soft-universe in which
the dissipation obstructions implied by Orevkov's 2 prohibitions
are the sole ones. In this world, the permissible dissipation
would be $(9,0,0)$, $(8,0,1)$, $(5,0,4)$, $(2,0,7)$. Then, Viro's
method would materialize the 3 bosons $1\frac{1}{1}\frac{18}{1}$,
$1\frac{4}{1}\frac{15}{1}$, $1\frac{7}{1}\frac{12}{1}$, but not
the last one $1\frac{9}{1}\frac{10}{1}$ as shown by a quick
inspection of the table. Of course in principle Viro's theory is
complete for the singularity $X_{21}$ and thus we have just the
two bold-faced dissipation available, and accordingly only the
solitary scheme listed by Viro. Here a solitary scheme means a
doubly nested scheme with only only one outer oval so as to be in
the bosonic strip.

\begin{figure}[h]\Figskip
%\vskip-1.2cm\penalty0
%\centering
\hskip-2.7cm\penalty0
\epsfig{figure=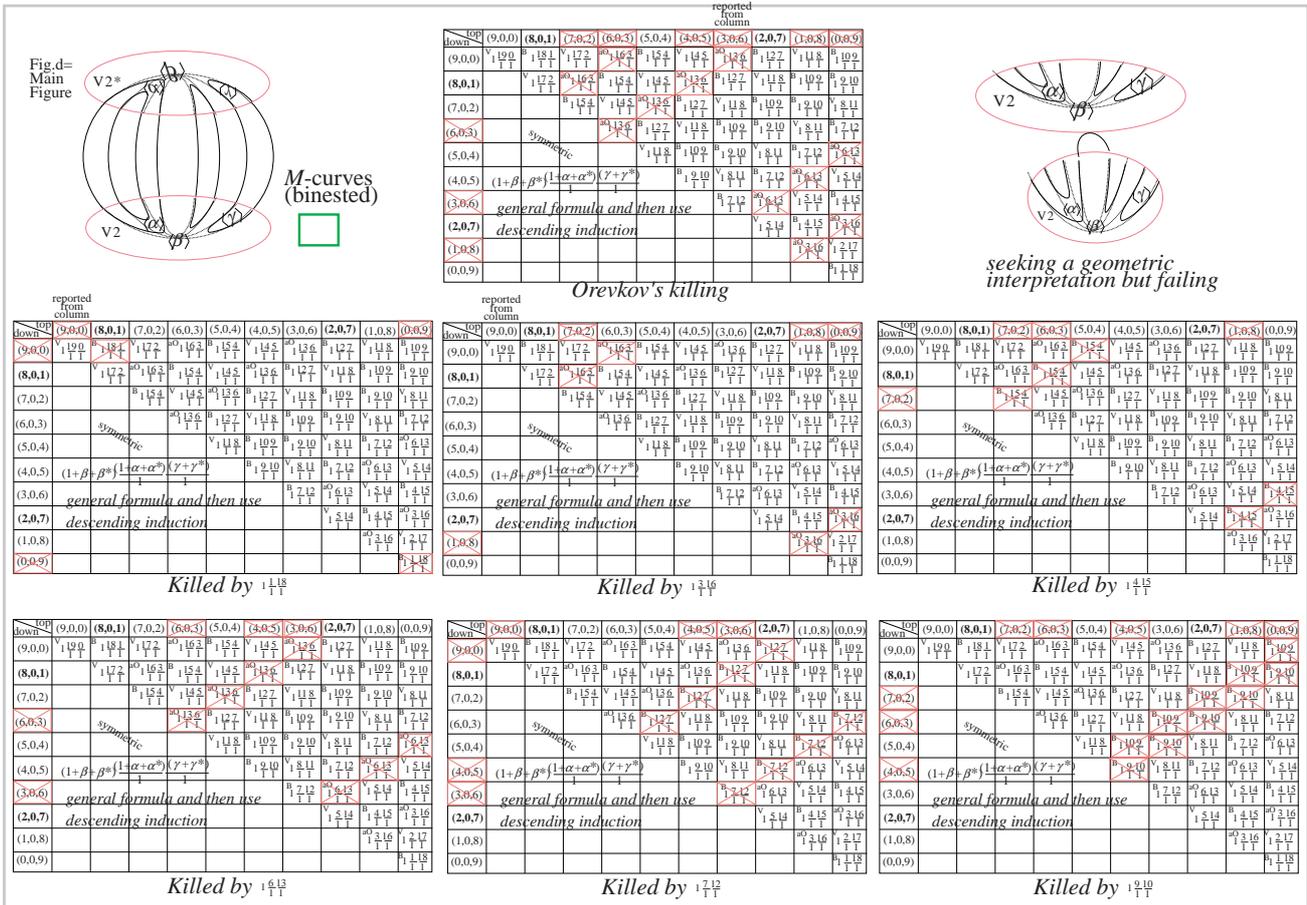,width=172mm} \captionskipAG
  \caption{\label{ViroDEGREE8_BOSONIC:fig}%
  Viro's extended bosonic patchwork, yet no serious foundations}
\figskip
\end{figure}

Alas, we lack a direct geometric interpretation of the dissipation
of $X_{21}$ as a global geometric object, akin to the smoothing of
an ordinary multiple point (say $m$-fold point) directly
interpretable as an affine curve of degree $m$. With such an
analogy available, probably that most restriction could reduce to
Hilbert-Gudkov's classification of $M$-sextics.

%%%%% SHORT PASSAGE KILLED BECAUSE NOT COMPLETED

{\it Added\/}[06.08.13].---Actually it seems worthwhile to
tabulate as well the mixed dissipation V2/V3 which still leads to
$M$-schemes. We had not the time to complete this, but probably we
do  it more systematically elsewhere in this text.

\iffalse

\begin{figure}[h]\Figskip
%\vskip-1.2cm\penalty0
%\centering
\hskip-2.7cm\penalty0
\epsfig{figure=ViroDEGREE8_BOSONIC_MIX.eps,width=172mm}
\captionskipAG
  \caption{\label{ViroDEGREE8_BOSONIC_MIX:fig}%
  Viro's extended mixed bosonic patchwork, yet no serious foundations}
\figskip
\end{figure}

\fi

\subsection{Extended Viro's composition table}

[06.08.13] It seems important to work out extended tables of
compositions as to understand which dissipations are forbidden for
global reasons. We shall construct four large tables extending
Viro's composition (Fig.\,\ref{ViroDEGREE8:fig}) by allowing all
logically permissible values of the parameters. As a result, we
shall either get new schemes (in case Viro's patches list was too
confined), or, patch censorship  whenever the resulting global
scheme is prohibited.

Those extended tables appeal the following comments.

[V1/V1] Schemes of this table V1/V1  land in the sub-nested realm
while realizing both bosons. The structure of this table is very
simple, e.g. invariant along all anti-diagonals. In particular,
the quadri-ellipse could create all schemes of the 3rd pyramid in
a very continuous fashion safe those of Shustin's strip (zero
outer ovals). This could trivialize all $M$-species cooked by
Korchagin (19 many), Chevallier (4 many), Orevkov (one),
constructed along more tricky procedures. Alas,
%(or fortunately?)
 it seems that
there is more rigidity in the dissipation of singularity $X_{21}$
(alias quadruple flat point). Since none of the subnested schemes
is presently known to be prohibited, we cannot exclude any
 of the 15 logically possible dissipations of
$X_{21}$ compatible with Gudkov periodicity. It can be remarked
that the number of big eggs is quantified as 2,6,10,14,18
(fourfold periodicity), and this forces $\be$ to be 1 or its
%periodicals
companions modulo 4 (i.e., 5, 9).

[V2/V2] On this 2nd table we land in the binested realm, where
reigns
%in principle
deep braid-theoretic obstructions of Orevkov. This
%results
%obstructions of
%cribles
%%%%%%%%%%prohibits
%sifts away
rules out six patches: $(7,0,2)$, $(6,0,3)$, $(4,0,5)$, $(3,0,6)$,
$(1,0,8)$ and $(0,0,9)$. However two patches with $\be=0$ are left
intact (namely $(9,0,0)$ and $(5,0,4)$). Those are enough to
create all binested bosons safe one, within Viro's simplest method
of the quadri-ellipse.

\begin{figure}[h]\Figskip
%\vskip-1.2cm\penalty0
%\centering
\hskip-2.7cm\penalty0
\epsfig{figure=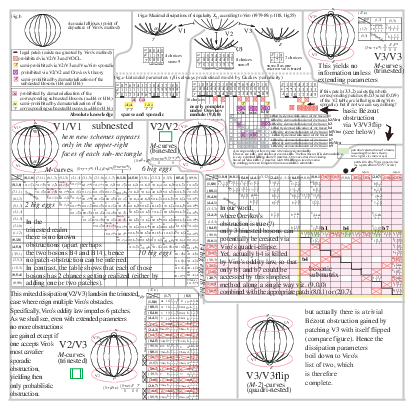,width=172mm}
\captionskipAG
  \caption{\label{ViroDEGREE8_extended:fig}%
  %Viro's extended mixed bosonic patchwork,
  %yet no serious foundations
  Bosonic table of elements (after Viro et cie)}
\figskip
\end{figure}

{\it Added\/} [07.08.13].---On the central table of crosses and
circles (sembling the famous game of life/go), crosses indicate
patches killed by  dematerialization of a boson (scheme). Circles
indicates pseudo-kills of a pair of patches (Heisenberg
incertitude). One remarks quickly that pseudo-bosons known to
exist (three of them thanks to Viro) kill relatively few patches
if they would dematerialize: $b2$ kills only one patch and so does
$b8$ (yet killing twice his victim). As to $b5$ it kills only 2
patches. In contrast, the 2 anti-bosons of Orevkov (known to
dematerialize!) kills both 3 patches. Positing that nature
dislikes criminals, we may expect the following moral akin to {\it
Kant's imperative (moral) law\/}:

\begin{Scholium} As soon as a (bosonic) scheme kills $3$ or more
patches, then it is highly criminal.
%(not to say a serial killer) and
Nature cannot tolerate such serial killers. From this standpoint,
both bosons $b4$ and $b9$ are criminal (with resp. $3$ and $5$(!)
murders) hence judged impossible. Instead, the bosons $b1$ and
$b7$ kill only two patches hence tolerated by society and more
likely to exist. Concretely it suffices for the patch $(9,0,0)$ to
exist for both those bosons to materialize. This would be a big
advance on Hilbert's 16th without requesting more imagination than
Viro's basic method, i.e. without having to resort to artistic
curves like  Viro's beaver, horse, Shustin's medusa, or our
embryos (compare  the sequel of this text and
Fig.\,\ref{SIMPLIFIED-TABLE_gurus:fig}).
\end{Scholium}

{\it Added\/} [08.08.13] (but fairly stupid).---In fact the
V2-patch $(9,0,0)$ plays a pivotal role. Can we disprove its
existence by another method? As yet we only composed the patches
for $X_{21}$ with themselves, yet we could try to glue them with a
double point to get sextics ($4+2=6$). We made some few pictures
below (Fig.\,\ref{ViroDEGREE8_patch:fig}), yet often contradicting
B\'ezout frontally. One configuration yields permissible schemes
yet its singular model is anti-B\'ezout. At any rate, even
patching $(9,0,0)$ yields a permissible Hilbert' sextic, so that
no obstruction is recorded against the $(9,0,0)$ patch. Actually,
one may wonder if a sextic can tolerate the singularity $X_{21}$
at all. One  argument is that an inner perturbation of the tangent
at the quadruple point will intercept the curve in at least 8
points,
%%(due to the local pattern of the singularity),
%
preventing realizability in degree 6. The argument simplifies by
just counting intersections by multiplicity (without perturbing).
So our idea is full rubbish. (As a matter of fuck, if the double
point is dissipated while creating a micro-oval then Harnack would
be foiled by the way.) We can still imagine and trace with our
heuristic embryo method complicated singular octics with a
singular point $X_{21}$ plus another distribution of triple points
while expecting to get curves obstructed  by Viro's imparity law,
say.

\begin{figure}[h]\Figskip
%\vskip-1.2cm\penalty0
%\centering
\hskip-2.7cm\penalty0
\epsfig{figure=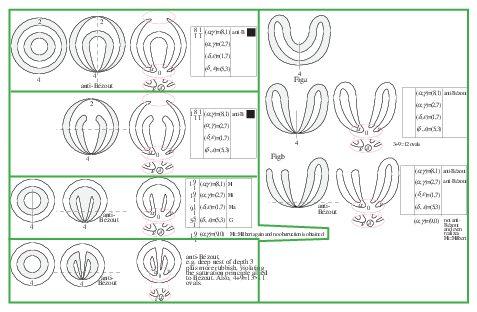,width=172mm} \captionskipAG
  \caption{\label{ViroDEGREE8_patch:fig}%
  Testing the patch $V3(9,0,0)$ in degree 6 (but completely
  erroneous)}
\figskip
\end{figure}

\begin{Scholium}
Even under Orevkov's  prohibitions (only $2$ schemes yet severe
damages over $6$ patches), it is still
%%% envisageable
possible for
Viro's simplest method to create $3$ among the $4$ binested octic
$M$-schemes, namely all but $1\frac{9}{1}\frac{10}{1}$, which
%%merits
deserves perhaps the name of Higgs boson.
%%%%Higgs CHECKED IN BOOK
Yet, a closer look to the table {\rm V2/V3} shows that Viro's
imparity law kills the patch {\rm V2}$(5,0,4)$, and thus the hope
to get $b4:=1\frac{4}{1}\frac{15}{1}$.
\end{Scholium}

\begin{proof}
Look at the grid, and see that all red-colored rows are killed by
Orevkov. Propagating anti-diagonally the Gudkov symbols we see
that each of our bosons can only be realized once on the upper row
of the bosonic submatrix (yellow framed). This holds true with the
exception of Higgs boson which lacks any such realization.
\end{proof}

More information  comes from the combined smoothing V2/V3, landing
in the trinested realm interspersed by a myriad of Fiedler/Viro
prohibitions. This  supplies additional information. Perhaps we
should stay critical about those highbrow prohibitions, which
 eventually falsify the definitive solution of Hilbert's 16th
problem in degree 8.

[07.08.13] Now let us tabulate the V2/V3 composition of patches.
From the scratch we observe that V3-list is very short, yet
perhaps it could be enlarged by introducing a 3rd parameter
$\varphi$ counting ovals at other places, e.g. in the strip of the
two nested arcs, i.e. like the $\al$-position of the patch V2. Yet
let us skip this difficulty for the moment.

On doing this table (Fig.\,\ref{ViroDEGREE8_extended:fig}) we see
that the patches $(3,4,2)$ and (Tupolev) $(1,4,4)$ are prohibited
by Viro's imparity law (VIL).
%%%%% A contrario,
Contrarily, the patch $(9,0,0)$ is left intact as it creates
admissible schemes. Then, $(7,0,2)$ is killed by VIL, but was
already by Orevkov. Next, $(6,0,3)$ is not killed by Viro, but
was by Orevkov. Then, $(5,0,4)$ is killed by Viro but was not by
Orevkov. The sequel is perfectly regular obeying an evident
periodicity of two. So, $(4,0,5)$ is not prohibited by Viro, but
was by Orevkov. As to $(3,0,6)$ it is prohibited by both Viro and
Orevkov, and idem for  $(1,0,8)$. Finally, agent $(0,0,9)$ is not
killed by Viro but was by Orevkov.

This may be summarized

\begin{Scholium} Many,
but by far not all, restrictions on Viro's diagram of dissipations
for {\rm V2\/} are explained by Viro's imparity law. Compare the
red crosses on {\rm Fig.\,\ref{ViroDEGREE8_extended:fig}V2/V3\/}
leaving open  the patches $(9,0,0)$, $(6,0,3)$, $(4,0,5)$ and
$(0,0,9)$. After finer sieving under
%%%believes in
Orevkov's behalf, from those only survives $(9,0,0)$ (apart of
course the patches declared existing
%%% in boldface on the table).
because constructed by Viro).
\end{Scholium}

Reporting those Viro obstructions on the  former table V2/V2 via
black bullets we see that only $(5,0,4)$ is additionally killed in
the bosonic range ($\be=0$), yet the hypothetic smoothing
$(9,0,0)$ still leaves open the realizability of two bosons via
Viro's simplest method (quadri-ellipse=quadri-lips), namely
$b1:=1\frac{1}{1}\frac{18}{1}$ and $b7:=1\frac{7}{1}\frac{12}{1}$.
So:

\begin{Scholium} The simplest bosons are perhaps $b1$ and $b7$.
(Compare eventually with
Scholium~\ref{bosons-elementary-(1_1_18+1_7_12)} which gave
exactly the same conclusion.)
\end{Scholium}

Could it be that experts missed to mention this miraculous patch
$(9,0,0)$? If not, what is the reason prohibiting it? As far as we
know this question is still open today. For cross-reference, let
us formulate this separately:

\begin{ques} {\rm (Patch mirabilis)} \label{patch-mirabilis-V2(9,0,0):ques} Consider the patch {\rm V2\/}$(9,0,0)$ as defined by
Fig.\,\ref{ViroDEGREE8_extended:fig}. If someone succeeds
constructing this patch, then two new bosons are materialized, and
so Hilbert's 16th is advanced. So our question is whether anybody
on the planet knows about a technique to forbid this patch.
(Remind that neither Viro's imparity law nor Orevkov imped this
patch.)
\end{ques}

Another question, is how would the world looks alike if Orevkov's
link theoretic obstruction(s) collapse(s), yet Viro's imparity law
persists true. In this scenario, the periodic table of elements
becomes Fig.\,\ref{ViroDEGREE8_extended_BIS:fig}. Let us call for
simplicity a {\it boson\/} just a binested $M$-scheme with the
minimum possible of one outer oval (as imposed by Gudkov
periodicity, i.e. Rohlin's signature theorem modulo 16). On
reporting on the table V2/V2 the obstruction coming from Viro's
oddity law on the table V2/V3, we get the red colored strip
obstructed. Then for each boson (each appearing twice
palindromically in the upper-right strip of the bosonic
sub-matrix) we may extend its symbol anti-diagonally while
avoiding the red-strips forbidden by Viro. Encoding by b$n$ the
boson with symbol $1\frac{n}{1}\frac{19-n}{1}$ we see that the
zeroth boson $b0=2\frac{19}{1}$ admits one realization nearby the
quadri-ellipse. Remind at this stage that there is only one
(known) construction of this scheme via Viro's horse. For the 1st
boson $b1=1\frac{1}{1}\frac{18}{1}$ it admits 2 realizations, etc.

In the hypothetical scenario that all Viro admissible patches are
practicable we see that {\it all\/} bosons would be constructible,
and Hilbert's 16th puzzle would be  settled (modulo an
understanding of the 2 bosons in the subnested case).

\begin{figure}[h]\Figskip
%\vskip-1.2cm\penalty0
%\centering
\hskip-2.7cm\penalty0
\epsfig{figure=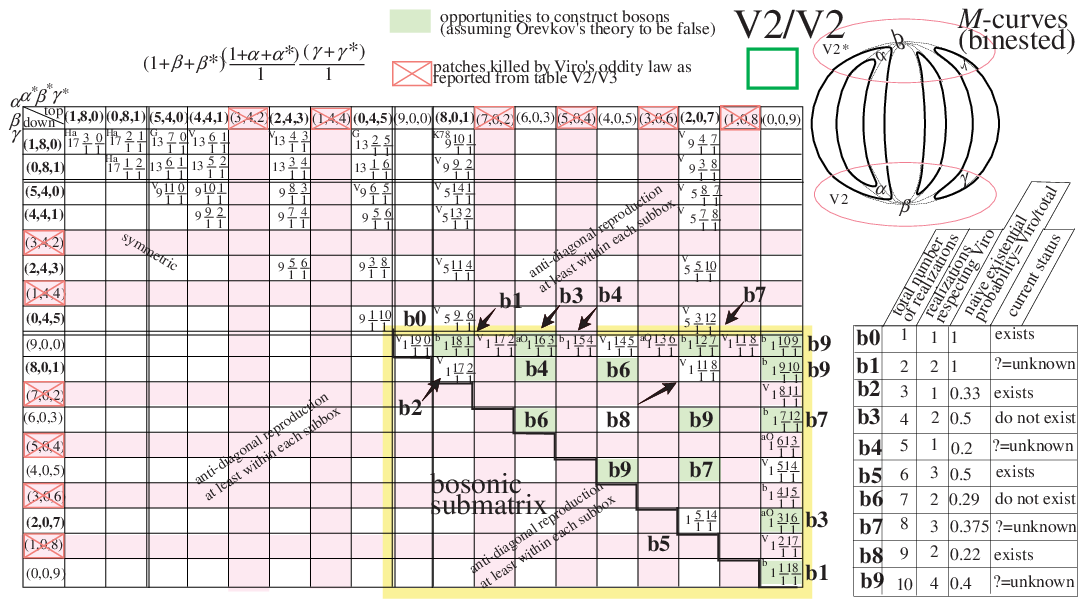,width=172mm}
\captionskipAG
  \caption{\label{ViroDEGREE8_extended_BIS:fig}%
  %Viro's extended mixed bosonic patchwork,
  %under brutal postulation that Orevkov is wrong
  In a world where Orevkov's obstructions are wrong, all binested
bosons could be (spontaneously) created out of Viro's
quadri-ellipse} \figskip
\end{figure}

It remains next to analyze what happens if we introduces an extra
parameter $\lam$ in the patch V3. Here the relevant table is given
in Fig.\,\ref{ViroDEGREE8_extended_TRIS:fig}. First, when
introducing this extra parameter it is not clear  where to place
it: either in the doubly nested lune on the left part of the patch
or in the two simple lunes on the right side. Of course, it cannot
be nested inside the double lune without troubleshooting B\'ezout.
Further it could be placed in the inner simple lune, yet by
analogy with V2, micro-ovals tend perhaps to be
%%%spreading sprayed
spread along the main tangential direction so that they are rather
distributed in the lateral simple lune on the extreme-right of
Fig.\,V3. Apart from this, we could a priori imagine that
micro-ovals can appear in both the double and simple lunes on the
extreme left- and right parts of the patch V3. If so is the case,
%(if not in reality at least in fiction),
we need a 4th parameter (say $\rho$) to describe the generic patch
V3. So $(\lambda, \rho)$ stand for left and  right.

One important remark is that if glue V3 with itself symmetrically,
then the scheme will be quadruply-nested (quadri-nested) at least
for generic values of $\de>0$ and $\rho>0$. Accordingly, it seems
realist to set $\rho=0$ throughout, safe when $\de=0$ but then
maybe Gudkov periodicity is not fulfilled, or alternatively the
patch V3 degenerates to the patch V2.

Hoping  not to miss something essential, let us first work out a
table with only 3 parameters $\de,\ep,\lam$. Now the basic idea is
that the $\de$ many micro-ovals are Swiss-cheese holes white
colored in the Ragsdale membrane and those holes may be transmuted
to the position $\lam$ without changing the Euler characteristic.
So of each of the God-Viro's given parameters generate a little
cascade of new parameters, where $\de$ is successively diminished
by one unit. In contrast it seems that $\ep$ is predestined by
Gudkov periodicity, and cannot do small fluctuations without
changing $\chi$. For $\ep$-ovals, the main option would be to jump
in the inside of the doubled lune, yet this causes troubles with
B\'ezout unless $\de=0$.

On composing V3 with itself we see that Viro's imparity law kills
the dissipation $(4,3,1)$ and $(2,3,3)$, but not all the others.
However it must be remarked that the schemes so created are all
fairly standard,  in particular never conflict with Viro's
sporadic obstructions. It may be hoped that the true secrets will
be revealed in the composition table of V2/V3. Here when filling
along vertical lines the 1st surprise comes when composing (6,0,3)
with $(3,3,2)$ as we meet then Viro's sporadic obstruction (those
being represented by orange crosses, as they are more likely to be
false than the red crosses materializing Viro's imparity law).
Completing this table
%%%in the obvious fashion,
shows that Viro's imparity law prohibits only the V3-patches
$(4,3,1)$ and $(2,3,3)$, hence actually exactly the same as those
ruled out by the V3/V3 table. Of course, it should be no surprise
that we never visited the bosonic strip, since looking at the
scheme of V2/V3, for it to be binested requires $\ga=0$ but then
we see 2 outer ovals (and not just one as in the bosonic strip).
It seems  strange that we only met one of Viro's sporadic
obstruction; yet, contemplating once more the table of elements
(e.g. Fig.\,\ref{SIMPLIFIED-TABLE_gurus:fig}) one sees that all of
Viro's sporadic obstructions (safe one) concerns trinested schemes
without outer ovals, while our gluing V2/V3 exhibits at least one
such oval materialized by the small right lune. Hence
$4\frac{3}{1}\frac{3}{1}\frac{9}{1}$ appears as the most cavalier
of Viro's sporadic obstructions, and one could speculate Viro
being wrong when claiming it.

\begin{figure}[h]\Figskip
%\vskip-1.2cm\penalty0
%\centering
\hskip-2.7cm\penalty0
\epsfig{figure=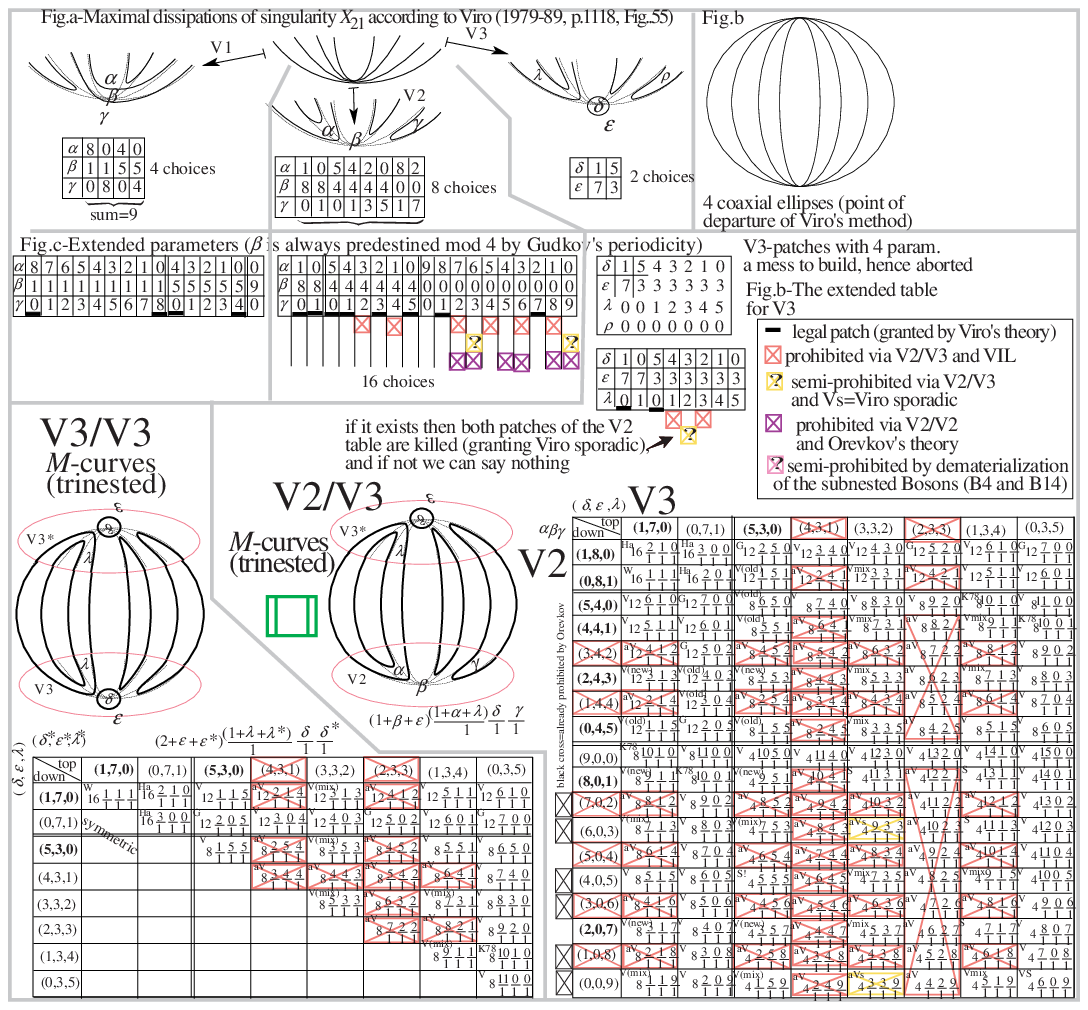,width=172mm}
\captionskipAG
  \caption{\label{ViroDEGREE8_extended_TRIS:fig}%
  Viro's extended mixed patchwork} \figskip
\end{figure}

It is time to synthesize the result coming from those composition
grids:

\begin{lemma}

$\bullet$ Via V2/V3, Viro's oddity law kills six $V2$-patches,
namely those with parameters $(3,4,2)$, $(1,4,4)$, $(7,0,2)$,
$(5,0,4)$, $(3,0,6)$, $(1,0,8)$. Via V3/V3 or V2/V3, Viro's oddity
law obstructs the two V3-patches with parameters $(4,3,1)$ and
$(2,3,3)$.

$\bullet$ Via V2/V3, Viro's sporadic obstruction of
$4\frac{3}{1}\frac{3}{1}\frac{9}{1}$ rules out at least one member
in the pairs $(6,0,3);(3,3,2)$ and $(0,0,9);(3,3,2)$. Put more
concretely, this means that if the $V3$-patch $(3,3,2)$ does exist
(and under the assumption that Viro's sporadic obstruction is
TRUE), then two $V2$-patches $(6,0,3)$ and (secret-agent)
$(0,0,9)$ are killed simultaneously.

Fairly concomitantly with this scenario we have finally:

$\bullet$ Via V2/V2, Orevkov's 2 obstructions rules out six
$V2$-patches among which 3 were already killed by Viro. Precisely
Orevkov's dematerialization of the boson
$b3=1\frac{3}{1}\frac{16}{1}$ kills 3 patches of whose 2 were
already killed by Viro, whereas the evaporation of the boson
$b6=1\frac{6}{1}\frac{13}{1}$ kills two new patches not ruled out
by Viro (at least formally, i.e. via  basic patchwork and his
oddity law). Hence one could speculate that at least the half of
Orevkov's obstruction pertaining to $b6$ is wrong.
\end{lemma}

\begin{proof}
Just look at the tables (especially
Fig.\,\ref{ViroDEGREE8_extended:fig} and
\ref{ViroDEGREE8_extended_TRIS:fig}).
\end{proof}

\subsection{What if Viro's imparity law is false: the big
decongestion}

[08.08.13] Speculating that even Viro's oddity law is wrong, the
dissipation theory of $X_{21}$ could be much richer and the world
would be a completely different smooth porridge, with Hilbert's
16th in degree 8 (potentially) much more trivial. Even if this
scepticism about Viro's oddity law may look retrograde, we remark
two interesting points. First,  it seems important to appreciate
exactly the shape of this simplified world with an abundance of
patches. Second, it may be remembered that without Viro's oddity
law, Shustin's disproof of Klein's (pseudo)-Ansatz (Klein vache)
as well as his disproof of one half of Rohlin's maximality
principle would be ruined.

We may fix as ground postulate that the dissipation theory of
$X_{21}$ is unobstructed, i.e. all values permissible with Gudkov
hypothesis are realized.

First,  notice that the upper right corners of the table V1/V1
(Fig.\,\ref{ViroDEGREE8_extended:fig}) fills out with perfection
the five rows of the 3rd subnested pyramid, safe for the schemes
with zero outer ovals where reigns a Shustin obstruction. Could
Shustin's obstruction be false as well? This looks a serious
challenge because as yet we never succeeded to reach this zone
even under dubious flexible pseudo-construction. At least it is
noteworthy that since the mandarine's range fails exploring
Shustin's strip, the latter fails inducing patches obstructions.
(Added in proof: We shall see later that this is not perfectly
true, if we work more liberally by allowing all patches). So let
us state this as follows:

\begin{lemma}
If the dissipation theory of {\rm V1\/} is unobstructed then all
schemes of the 3rd pyramid are created in a very continuous
fashion via table {\rm V1/V1\/} of
Fig.\,\ref{ViroDEGREE8_extended:fig}. In particular many tricky
constructions of Viro (horse and beaver), Korchagin, Chevallier,
Orevkov could be
%%%trivialized
relegated and everything could be accessed from the mandarine
(alias quadri-ellipse). In particular the 2 subnested bosons
$B4:=4(1,2\frac{14}{1})$ and $B{14}:=14(1,2\frac{4}{1})$ would be
 created.
\end{lemma}

It is important to notice that the patches V2 and V3 are somehow
coupled, i.e.  can be married in the joy of Harnack maximality,
whereas V1 is isolated. Of course we can imagine an avatar of V1
(say V0) with a micro-nest outside, but when gluing V0 with itself
 produce a nest of depth 3 plus one of depth 2 creating
thereby 10 intersections with the line through their centers.

One of the most important paradigm of the theory is independency
of smoothings. Perhaps this has to be revised as well, or requests
hardwork \`a la Viro/Shustin, etc.

Finally one may wonder why Viro does not mention the option of a
4th patch V4 with 2 nested lunes, or even 4 unnested lunes.
Actually one can also have an external branch and three lunes
inside. This latter patch looks especially important as it lands
in Shustin's range (subnested but no outer ovals). This leads us
to the next section of exotic patches, not listed by Viro but
logically possible at least a priori. This will sidetrack us into
the combinatorial study of all those patches.

\subsection{Exotic patches}

[14.08.13] The goal of this section is to prove:

\begin{theorem}
Singularity $X_{21}$ (flat quadruple point with $4$ branches lying
in the same half-plane and having 2nd order contacts between
themselves) can be maximally dissipated (with nine micro-ovals
bubbling out) along any one of the $14$ ways described by Viro
ranked into three classes {\rm V1}, {\rm V2}, {\rm V3}. Together
with the symmetric patches yields a total of $28$ smoothings of
$X_{21}$ (as many as those of the $7$-sphere according to
Milnor-Kervaire).
%%%or vice-versa).

Yet, there is maybe more ``exotic'' dissipations, potentially as
many as the white circles of {\rm
Fig.\,\ref{ViroDEGREE8_exotic_patches0:fig}}.
Along their hypothetic existence, one could primarily construct
some few new bosons (perhaps all of them!), and secondarily
recreate old schemes (especially those of Korchagin, Chevallier)
via the most basic Viro method, thereby gaining a discriminantal
%%%jouxtancy
%%parenty
kinship (Verwandschaft) limitropheness between such schemes and
the quadri-ellipse.

Finally, in principle all patches crossed off
%%%CHECKED IN DICO for BARRER
on {\rm Fig.\,\ref{ViroDEGREE8_exotic_patches0:fig}} are
prohibited.
\end{theorem}

\begin{figure}[h]\Figskip
%\vskip-1.2cm\penalty0
%\centering
\hskip-2.7cm\penalty0
\epsfig{figure=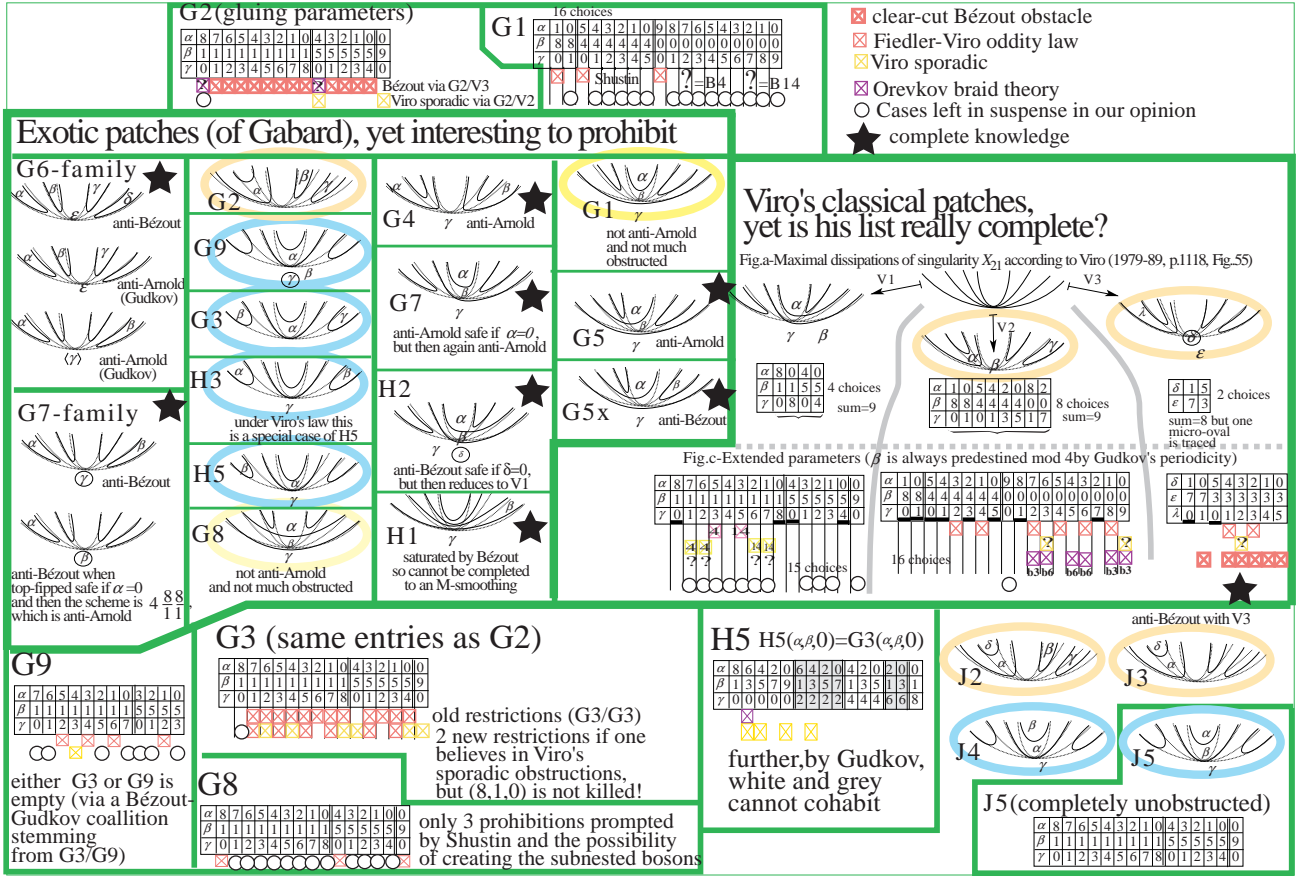,width=172mm}
\captionskipAG
  \caption{\label{ViroDEGREE8_exotic_patches0:fig}%
  A mess of exotic patches (cf. also
  Fig.\,\ref{ViroDEGREE8_exotic_patches0_SYS:fig}
  for a rationalization)}
\figskip
\end{figure}

[08.08.13] We look first at the patch G1 (for Gabard) on
Fig.\,\ref{ViroDEGREE8_exotic_patches0:fig} albeit Viro certainly
thought about it yet without listing it as he was probably not
able to construct it (or maybe knew obstructions  as early as
1980). Then we glue the patch with a symmetric replica and
contemplate the resulting patchwork. On it we see two big eggs are
traced. The table of elements
(Fig.\,\ref{SIMPLIFIED-TABLE_gurus:fig}) reminds us that Gudkov
periodicity predestines this number to be precisely 2 modulo
fourfold periodicity ($6,10,14,18$). Therefore we choose $\be=0$
modulo a periodicity of 4. Then, we can build the table of
parameters $(\al,\be,\ga)$ for the number of micro-ovals bubbling
in the G1-patch, whose (extended) parameters  are actually the
same as for V2. Next, we  build the composition table by
dissipating independently both singularities. The schemes so
obtained are essentially the same as for V1/V1, modulo the crucial
difference that we now obtain the schemes prohibited by Shustin.
Few basic remarks on the table: the evolution is same horizontally
as vertically, with quantum jumps across double bars. Therefore
the table propagates anti-diagonally  and we need just writing the
symbols occurring in the upper-right corners of each sub-boxes.
Remarkably, each of those corner-strips fills precisely one row of
the 3rd pyramid so that the last item is actually on the top of
the 1st pyramid. Also, the first element of each corner-strip is
prohibited by Shustin. Of course globally the whole is diagonally
symmetric.

\begin{figure}[h]\Figskip
%\vskip-1.2cm\penalty0
%\centering
\hskip-2.7cm\penalty0
\epsfig{figure=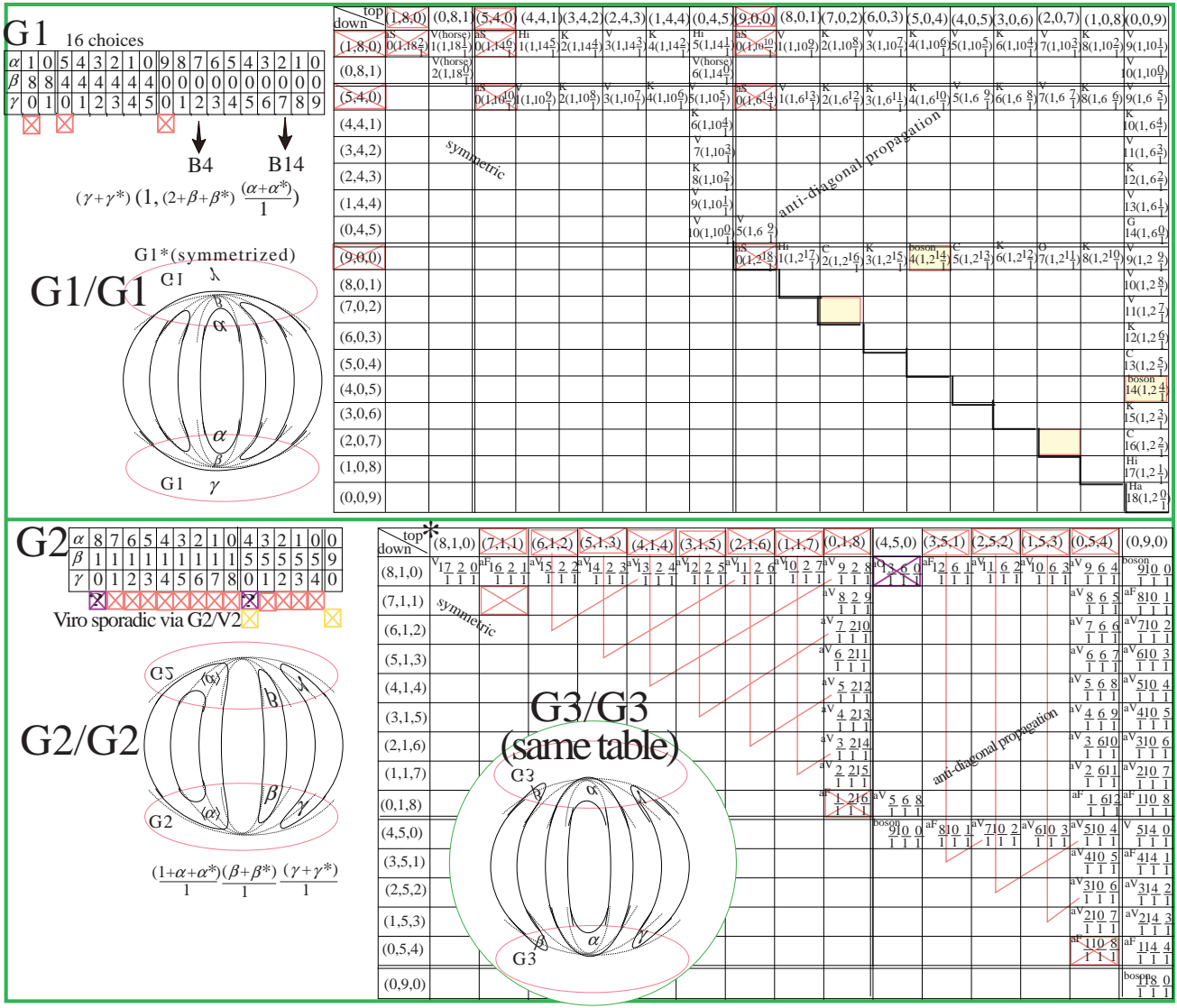,width=172mm}
\captionskipAG
  \caption{\label{ViroDEGREE8_exotic_patchesA:fig}%
  Patchworking exotic patches}
\figskip
\end{figure}

Granting Shustin's obstructions as correct,  each anti-Shustinian
scheme situated on the main diagonal kills the corresponding
patch-parameter. So, Shustin's prohibition of $(1,18\frac{2}{1})$
kills $(1,8,0)$,
$(1,10\frac{10}{1})$ kills $(5,4,0)$,
$(1,2\frac{18}{1})$ kills $(9,0,0)$. The two other obstructions by
Shustin not directly situated on the diagonal kills no definite
patches but a pair of patches with quantum incertitude about who
is exactly killed. For instance the anti-scheme
$S_{14}:=(1,14\frac{6}{1})$ kills either $(1,8,0)$ or $(5,4,0)$,
yet without precising which one. Actually due to the alinement of
all this table, it turns out that both patches are killed by
schemes situated on the diagonal. Yet we could imagine a world
where the diagonal Shustin obstructions are true but not the
others. Further it can be speculated about the dematerialization
of the two bosons and taking the diagonal representative it
results a destruction of a patch, namely $(7,0,2)$ and $(2,0,7)$
respectively.

Starting from zero knowledge (e.g. ignoring Shustin) we do not
know even which patches actually exist. If so then we could have
more destruction of patches than those merely coming from item on
the diagonal.

In reality, it seems that all the G1-series of patch is empty as
it is not listed in Viro 89. It is not clear if Viro just dresses
a list of patch (he is able to construct) or if he is claiming
completeness.
%
%From a naive viewpoint
A priori,  G1 could  lead to monsters (e.g. corrupting B\'ezout)
when glued with other V-patches. Yet, we doubt this to be the
case. So the scholium is:

\begin{Scholium} Could it be that Russian scholars
missed some patches
%or parameters thereof
so that Hilbert's 16th problem is actually trivial to solve in
degree $8$. Of course,  our world of continuous parameters
presupposes that  obstructions \`a la Viro, Shustin, Orevkov, are
false.
\end{Scholium}

We examine now how the patch G1 interact with those of Viro. Each
combination  G1/V1, G1/V2 and G1/V3 has to be envisaged and it
results
%%%fairly boring
(cumbersome) tables of compositions. In the patch G1 we put the
$\al$ micro-ovals in the central lune for otherwise suitably
reflecting the patch gives two subnests and so corrupt B\'ezout.

\begin{figure}[h]\Figskip
%\vskip-1.2cm\penalty0
%\centering
\hskip-2.7cm\penalty0
\epsfig{figure=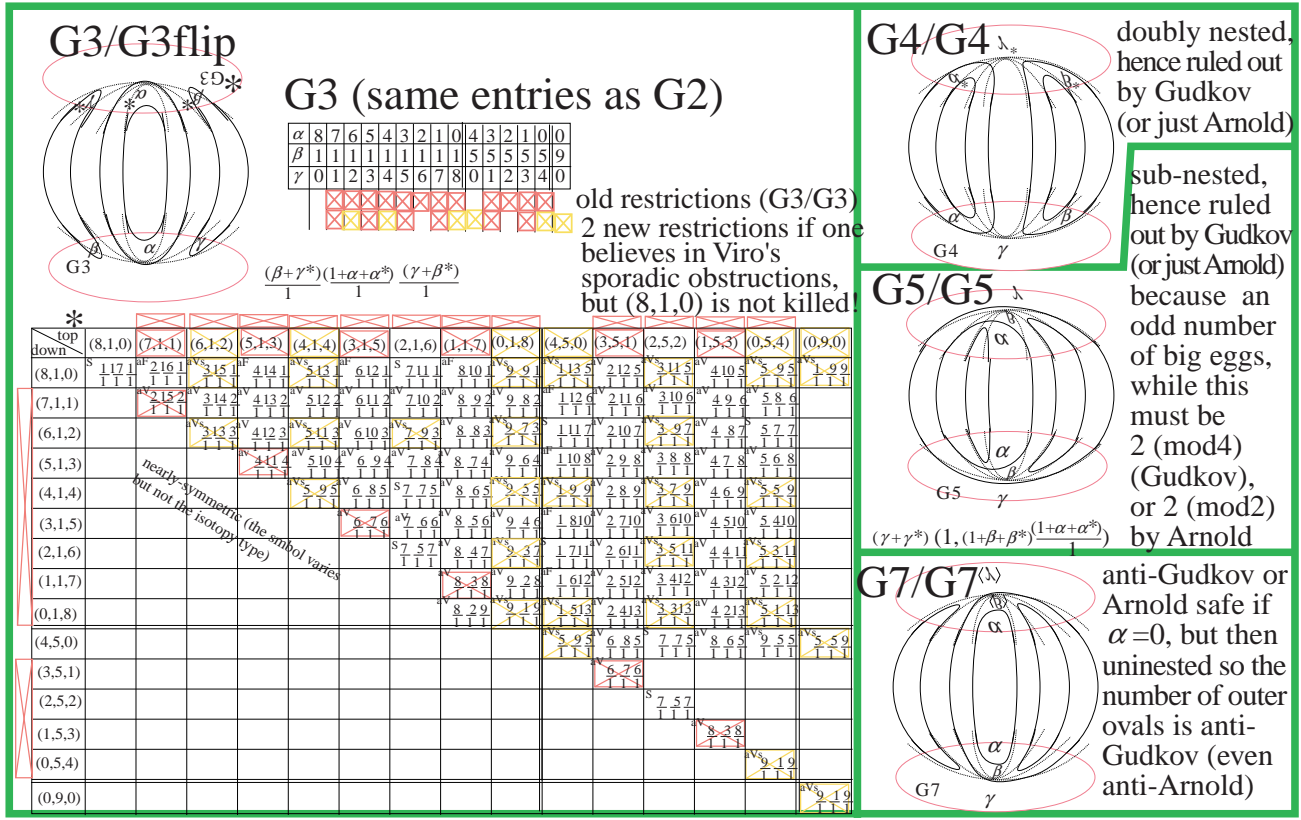,width=172mm}
\captionskipAG
  \caption{\label{ViroDEGREE8_exotic_patchesB:fig}%
  Exotic patchwork}
\figskip
\end{figure}

Patching G1 with V1 produces only $(M-2)$-schemes, yet those are
also subsumed to certain Viro prohibitions and things becomes
fairly tricky. Actually this idea pertains also to V1-patches
 glued with themselves in a non-maximal fashion.

As yet, we never reached really the realm of Fiedler and Viro's
sporadic obstructions mostly concentrated in the trinested case
without outer ovals (safe one exception
$4\frac{3}{1}\frac{3}{1}\frac{9}{1}$). The patch ideally suited to
explore this zone is G2 which is an exact copy of V2 safe that the
parameter $\be$ has been dragged inside the inner lune. On gluing
G2 with itself we get when $\al=\ga=0$ and $\be=max=9$ (so as to
arrange 22 ovals) the scheme $1\frac{1}{1}\frac{18}{1}$ which
although bosonic is at least Gudkov permissible. So we must choose
$\be=9$ and its companions modulo 4. Thus the parameter table for
G2 is actually the same as that of V1 involving 15 values. The
table G2/G2 will strongly conflicts with obstructions and the game
is to see if all patches are killed. When filling the table, one
observes interesting motions along the pyramid with  pleasant
foldings and the phenomenon of palindromic pathes. Granting the
obstructions of Fiedler and Viro, we see---by propagating
anti-diagonally up to the diagonal---that nearly all patches are
killed. For instance Viro's anti-scheme
$\frac{15}{1}\frac{2}{1}\frac{2}{1}$ kills $(7,1,1)$ and so on. At
first, Fiedler's earlier anti-scheme looks unused but will be at
the end via the palindromic effect. Our red broken-lines show how
an anti-scheme is propagated on the diagonal as to kill the patch
above it. The same
%token can be
discourse repeats in the 2nd diagonal block. However this argument
does not rule out the first patch of each
%%quantum
series, namely $(8,1,0)$, $(4,5,0)$ and $(0,9,0)$. Postulating
their (collective) existence would create a simple scheme by Viro,
a corruption of Orevkov and two bosons including Higgs's one
$1\frac{9}{1}\frac{10}{1}$, which is perhaps the most elusive of
all. Of course to be politically correct (joke of Viro, in Geneva
ca. 2010) w.r.t. Orevkov we could only activate sub-collection of
patches. As $(8,1,0)$ self-combined with itself yields a known
scheme of Viro it is the most likely to exist, yet we may imagine
that the 2 others exists individually as well, and this would
create new bosons. Now assume that (8,1,0) exists, and thrusting
in Orevkov's anti-scheme then $(4,5,0)$ is killed yet we could
still posit existence of $(0,9,0)$ and thereby materialize Higgs
boson $b9$ and $b1$. All this without conflicting  with factual
knowledge.

{\it Update} [01.10.13].---In fact the opportunities to get those
bosons via G3, are killed if one considers the symmetry of bending
discussed on Fig.\,\ref{ViroDEGREE8_exotic_patches0_BEND:fig}.
Indeed bending the patch yields one with a double couche (class I
in the catalogue) and those corrupt B\'ezout even when $\ga=0$
since flipping the patch creates two subnests. The sole case
actually is when both lateral parameters ($\be, \ga$) vanish, yet
this is not in line with Gudkov periodicity.

Maybe composing G1 with another patch we can really visit Viro's
sporadic obstructions and so rule out more G2-patches. Yet this
deserves to be studied tomorrow.

\begin{figure}[h]\Figskip
%\vskip-1.2cm\penalty0
%\centering
\hskip-2.7cm\penalty0
\epsfig{figure=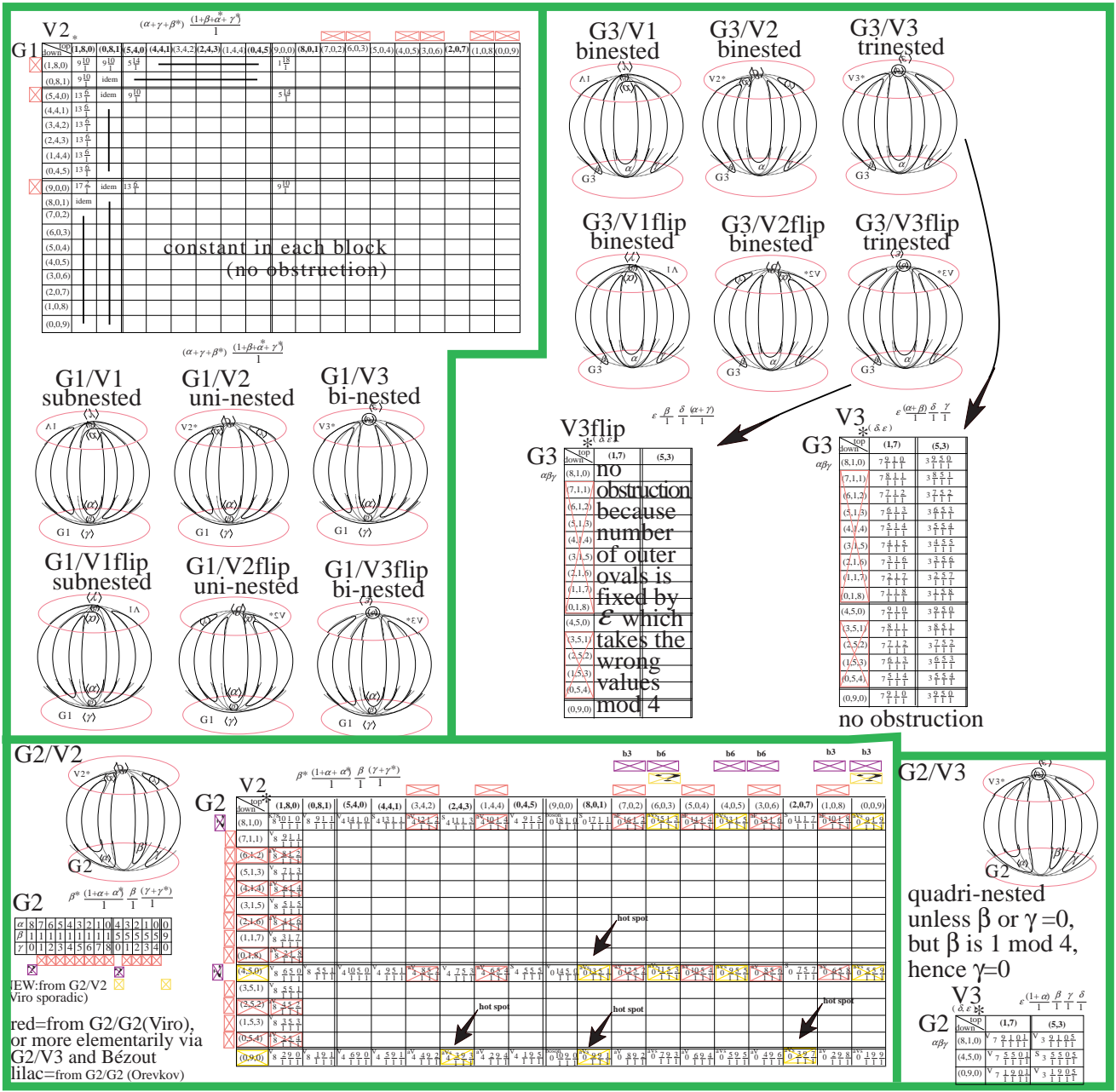,width=172mm}
\captionskipAG
  \caption{\label{ViroDEGREE8_exotic_patches2:fig}%
  Patchworking exotic patches (continued)}
\figskip
\end{figure}

Added [14.08.13].---It seems that we missed to combine G2/V2 and
G2/V3. From the latter it is inferred that G2 is nearly empty safe
3 places without recourse to Viro's oddity law but just B\'ezout.
Precisely, G2/V3 is quadri-nested (because $\de>0$) unless $\be$
or $\ga$ is zero. Yet $\be\equiv1 \pmod 4$, so that $\ga=0$. This
leaves only the three values $(8,1,0)$, $(4,5,0)$, and $(0,9,0)$.
Next, we may dress the G2/V2-table hoping to get more obstacles,
yet this
%hope
expectation is probably not borne out. Actually it suffices to
fill the table along  rows not yet prohibited. Along the first
horizontal row we meet Viro's imparity law yet at a V2-entry not
existing hence no  destruction of G2(8,1,0) is inferrable. On
filling more the table a 1st surprise occurs with entry
G2/V2=(4,5,0)/(8,0,1) where a sporadic Viro obstruction kills
G2(4,5,0). The next surprise occurs when Viro's most sporadic
obstruction on the apocalyptic symbol (year)
$4\frac{3}{1}\frac{3}{1}\frac{9}{1}$ kills the patch G2(0,9,0). A
3rd surprise comes when Viro's anti-scheme
$\frac{1}{1}\frac{9}{1}\frac{9}{1}$ kills a second time G2(0,9,0),
yet still through a sporadic obstruction. A 4th surprise occurs
when Viro's anti-scheme $\frac{3}{1}\frac{7}{1}\frac{9}{1}$ kills
once more again our patch G2(0,9,0), hence
%%%victim of
accusing triple mutilation. At the end, it seems still worth
looking explicitly at the table G2/V3, which is microscopic if we
restrict attention to the sole entry in doubts. This mini-table
yields no (supplementary) obstruction, as all tabulated schemes
are either due to Viro or Shustin. Let us resume this with a:

\begin{lemma}
Provided a suitable sub-collection of Viro's sporadic obstructions
is true, the patch family {\rm G2\/} collapses to a single
representative {\rm G2(8,1,0)} whose existence is not even
granted. It could be imagined that even if this patch existed it
would result no news on Hilbert's problem. This is in part true,
since through the table {\rm G2/G2\/},  the patch in question
$(8,1,0)$ only produces the boring scheme
$1\frac{2}{1}\frac{17}{1}$ due to Viro, but also in part false
because via the table {\rm  G2/V2\/}, {\rm G2(8,1,0)} glued with
the hypothetical {\rm  V2(9,0,0)} creates the boson
$1\frac{1}{1}\frac{18}{1}$ not yet known.
%%%to exist.
\end{lemma}

Added [14.08.13] It seems that we missed a discussion of the patch
H3 of Fig.\,\ref{ViroDEGREE8_exotic_patches0:fig}. First it seems
advisable to consider directly the more general patch H4. On
gluing H4 with a symmetric copy (our notation H4/H4) we see that
Viro's law forces at least one of the parameter $\be$ or $\ga$
being zero. By symmetry of the patch we may w.l.o.g. assume
$\ga=0$, and so we arrive at the new patch H5, where we just
relabelled according to alphabetic order. We may allow a parameter
quantified by 2 and not 4 yet still respecting Gudkov periodicity.
Note that we already normed $\be$ odd not to conflict with Viro's
law. The table H5/H5 shows that we find no {\it pure\/}
obstruction (along the diagonal), safe one implied by Orevkov's
(anti)-boson $b6:=1\frac{6}{1}\frac{13}{1}$. In contrast the other
Orevkov's boson $b3:=1\frac{3}{1}\frac{16}{1}$ misses his chance
to hit properly the diagonal (under anti-diagonal propagation),
thereby failing to induce a direct patch prohibition. Still, it
implies that at least one of (2,7,0) or (0,9,0) is killed.

On filling the table H5/H5flip---which really lands in the
trinested realm---we see that Viro's sporadic obstruction kills
with certitude $(6,3,0)$. We merely try to complete the table
along the diagonal where we have a simple crescendo dynamics of
both
%%wings
lateral coefficients. We see a hecatomb of patches killed by
sporadic obstructions, namely (4,5,0) and (0,9,0). In the second
diagonal box, Viro's most sporadic obstruction kills $(4,3,2)$.
After this 2nd diagonal sub-box passed, we never meet again
sporadic obstructions and nothing more is killed.
%, without having
%to fill more the table.
So what is proved? Answer is given by the:

\begin{lemma} Among all patches in the {\rm H5}-class, only
$(6,3,0)$, $(4,5,0)$, $(0,9,0)$, $(4,3,2)$ are directly prohibited
by Viro's sporadic obstructions (and Orevkov does not give
additional prohibitions).
\end{lemma}

\begin{figure}[h]\Figskip
%\vskip-1.2cm\penalty0
%\centering
\hskip-2.7cm\penalty0
\epsfig{figure=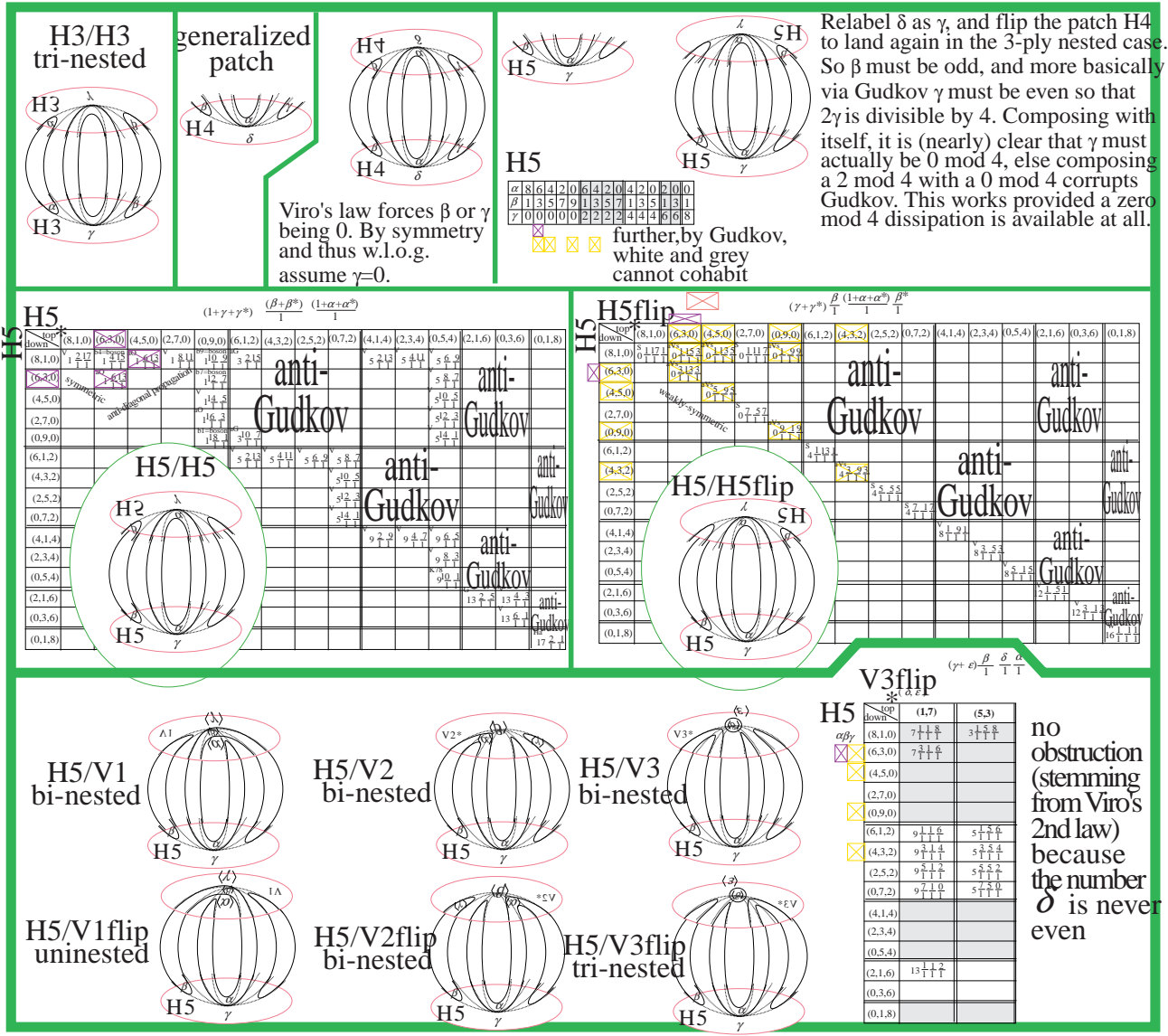,width=172mm}
\captionskipAG
  \caption{\label{ViroDEGREE8_exotic_patches_H3:fig}%
  Patchworking exotic patches (continued): H3/H3, etc.}
\figskip
\end{figure}

Next it is tempting to see if Viro's V$i$-patches induce
restrictions on H5. A little table of patchworks prompts that only
H5/V3flip lands in the trinested realm where there is some chance
to interact with Viro's 2nd law (of thermodynamics). Alas, the
schemes landing in the grey-shaded subregions  have the wrong
number of outer ovals to interact with Viro's law. So we tabulate
only the others, but we see quickly that we never meet Viro's
$(M-2)$-obstruction impeding an RKM-scheme to have only even
numerators. This is just because $\de$ is one of the numerator and
is odd. Of course, roughly speaking it seems that the reason of
this is that we already arranged Viro's law at the $M$-level, so
it is not much surprising that it is likewise respected at level
$(M-2)$. Unfortunately we are not much advanced on our problem.

Further it is perhaps pleasant to note that the first box of table
H5/H5 produce all four unknown (binested) bosons $b1=090\times
090$, $b4=810\times 630$, $b7=630\times 090=450\times 270$ (with
obvious abridged notations), yet all patches used in those bosonic
constructions are destructed by Viro's sporadic obstructions. So
the scholium seems to be that Viro's sporadic obstructions freeze
the
%ebullative
boiling formation of bosons and keeps cold the temperature of the
algebraic universe. Notwithstanding the creationism of any boson
is still not intrinsically prohibited by Viro's sporadic
amendments. Perhaps there is a statement of the shape:

If all of Viro's sporadic laws are true then no boson can be
created out of the quadri-ellipse\footnote{[01.10.13].---This is
not quite true, because even within Viro's patches family C2=V2,
there is via the patch mirabilis $C2(9,0,0)$ a chance to get both
bosons $b1$ and $b7$.}. At least this phenomenon was true for the
patch G2. In contrast, for the subnested bosons B4 and B14
accessible e.g. via G2 (or even V1), it seems that Viro's sporadic
rules have no prohibitive impact upon the bosonic formation out of
the quadri-ellipse.

From this basic composition method (of which Viro surely carefully
examined the combinatorics ca. 3 decades ago) it emerges the
following (seemingly paradoxical) principle: The more a patch
resembles those of Viro, the more he will interact with them at
the (maximum) $M$-level and so more the patch will be prohibited.
This is best exemplified by the patch G2. In contrast if the patch
exploits a totally different geometry then it will not much be
attacked by Viro's prohibitions (typical examples G1 and G8, yet
then slightly attacked by Shustin's rules). Perhaps a patch like
H5 is the medium range (liquid phase) where there is not too much
prohibitions (=cold regulated world).

Next we imagined further the patch J2
(Fig.\,\ref{ViroDEGREE8_exotic_patches0:fig}). The latter cannot
exist with $\de>0$ unless $\be=\ga=0$ (just by B\'ezout). So the
patch becomes J3 which is still anti-B\'ezout when patched with
V3.

\subsection{Working more systematically}

[14.08.13] It is clear that our random exploration must be
rationalized by doing a more proper census of all patches.
Basically this involves first the ground {\it dessin\/} (involving
four arcs), and then the labels counting micro-ovals (topological
circles bubbling out of the blue). Many configurations are just
ruled out by B\'ezout and then sometimes Arnold, or
Gudkov=Rohlin's periodicity. It seems now a vital task to get a
more lucid classification of all logically possible patches than
what we presented before. Each patch has 4 ground branches
traversing the disc (local neighborhood of the patching-surgery).
So there is by Jordan separation, five distinct zones where to
assign the labels. Location which are too deep tend being
prohibited by B\'ezout. For each dessin we shall list all its
incarnations, ideally interpreted as a morphogenetic process akin
to  bifurcations of species (under ``Darwinistic'' evolution).
Having done this properly we should arrive at a more rational way
to enumerate first the dessins, and then all the patches. In
addition each dessin can acquire an island like the patch V3 of
Viro. Evidently the number of island is $\le 1$, because if 2,
then B\'ezout for conics is foiled having already the doubled
quadrifolium when doubling the patch.

\begin{Scholium}
To enumerate properly (\`a la Newton, Hilbert, Polotovskii, etc.)
one can only be guided by a true understanding of the inherent
geometry of the world, which is often akin to morphogenetic rules
of evolutions transcending Darwinism, and hopefully still
available in our aged brains.
\end{Scholium}

All this requests boring hygienical work, yet necessary as it
seems that  we as yet missed the patch J4 of
Fig.\,\ref{ViroDEGREE8_exotic_patches0:fig}, which is B\'ezout
admissible and not corrupted by a gluing with V3.

[15.08.13] Gluing J4 with itself yields a subnested scheme with
$2\al$ big eggs (=oval at depth one). By Gudkov periodicity this
number has to be 2 modulo 4, and so we choose $\al$ as being one
mod 4. (Warning: maybe there is more choices like in our study of
H5.) Finally as in our patch parameters it is $\be$ which is
quantified by fourfold periodicity, it looks desirable to relabel
$\al, \be $ by permuting them. (This yields the patch J5.) The
composition J5/J5 will land in the subnested realm
%where there is
%no known obstruction
%%%freed
exempt of all obstructions apart those of Shustin pertaining to
the absence of outer ovals. Those will not concern us as J5/J5 has
at least 2 outer ovals. Hence filling the table of
Fig.\,\ref{ViroDEGREE8_exotic_patches_J5:fig} is almost
unnecessary, as we know a priori not getting any obstruction
except perhaps hypothetical destructions of the patches 711 and
216 in case of a hypothetical destruction of the bosons
$B4:=4(1,2\frac{14}{1})$ and $B14:=14(1,2\frac{4}{1})$. Further,
even when composing with Viro's patches V1, V2, V3 no obstruction
results as the produced $(M-2)$-scheme are never trinested. So as
far as we can see:

\begin{lemma}
The patch {\rm J5} in completely unobstructed.
\end{lemma}

\begin{figure}[h]\Figskip
%\vskip-1.2cm\penalty0
%\centering
\hskip-2.7cm\penalty0
\epsfig{figure=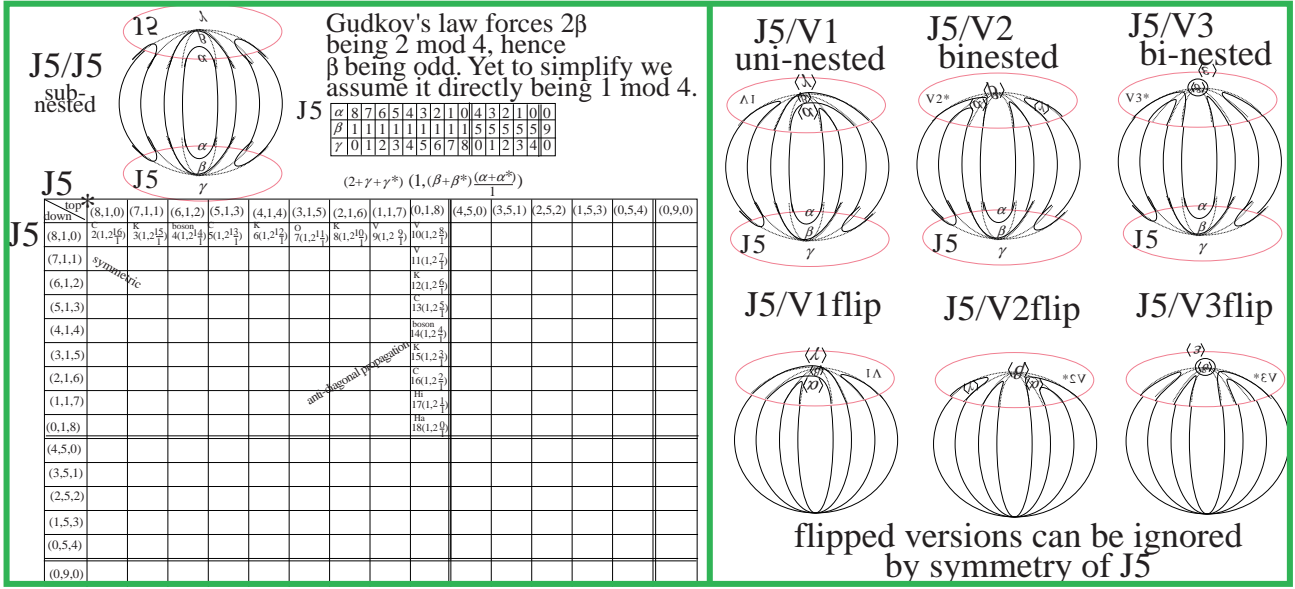,width=172mm}
\captionskipAG
  \caption{\label{ViroDEGREE8_exotic_patches_J5:fig}%
  Patchworking exotic patches (continued)}
\figskip
\end{figure}

One radical obstruction arises if the patch G9 admits a
representant with $\ga>0$, in which case the family J5 collapses
to its 3 items with $\al=0$.

Now we decided to work out a more systematic table of patches. The
idea is to start with a list A, B, C,\dots, J of ground dessins
(ordered from unnested to much nested) which is obviously
exhaustive. For each of them, we have then to imagine the
different places where to put (Greek) label measuring the number
of (bubbling) micro-ovals. On doing this plate we discovered a new
species namely C1, admittedly much akin to Viro's patch C2=V2, and
therefore creating the same schemes. Yet it seems still of
independent  interest to know exactly which patches are realized.
Further the type B1 seems to depend on 4 parameters and we
probably only studied it unsystematically as yet. Of course by
symmetry one could normalize so that $\al \le \ga$. However, we
remind that under Viro's law at least one of $\al$ or $\ga$ must
vanish. As another recompense of our more systematic work, we
discovered another new patch namely G2 of
Fig.\,\ref{ViroDEGREE8_exotic_patches0_SYS:fig}. Of course there
will be frictions between G1 and G2, yet maybe at least one of the
patch could support many representatives.

\begin{figure}[h]\Figskip
%\vskip-1.2cm\penalty0
%\centering
\hskip-3.7cm\penalty0
\epsfig{figure=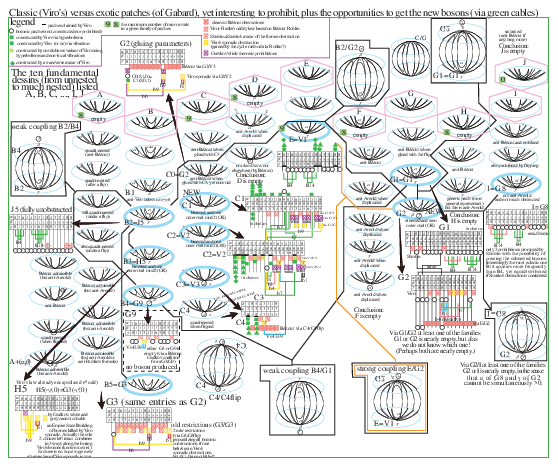,width=192mm}
\captionskipAG
  \caption{\label{ViroDEGREE8_exotic_patches0_SYS:fig}%
  Systematic patches (continued and hopefully finished):
  extended kit
  of patches aiding acute cases of nicotinism
  %tabacosis
  (e.g. Heinz Hopf in World War~I)}
  %to get
  %%%%sevred %%CHECKED IN DICO
  %weaned/tranquilized for a while}
\figskip
\end{figure}

As usual, we compose this (new) patch G2 with itself. Of course we
can by Viro's law rules out those values of $\ga$ which are even
(and positive) without having to work out the full tabulation.
Still, on the diagonal we meet two sporadic obstructions (due to
Viro). One could hope to get more obstructions via G2/V3, but
those schemes have the wrong periodicity on $\chi$ to interact
with Viro's 2nd law. More simply, as $\de$ is odd in our schemes
(so no interaction with Viro's 2nd law which forbids all
numerators being even for an RKM-scheme). Paraphrasing in more
geometric fashion, the $(M-2)$-schemes generated by our table
lands below the $M$-peaks of
Fig.\,\ref{Degree8-(M-i)-curve-TABLE:fig} as opposed to the
depressions (valleys) where Viro's 2nd law is reigning.

Hence, the dissipation theory of the patch G2 is not so much
obstructed as being de facto empty (like in Viro's census). To
know precisely which schemes could result from G2-patches we must
fill more the table away from the diagonal. Our interest was to
measure the degradations effected by Viro's very sporadic
obstruction ($4\frac{3}{1}\frac{3}{1}\frac{9}{1}$). The answer is
disappointingly that this adds no new obstructions as those
already read out from the diagonal. Green frames shows
opportunities to construct {\it new\/} bosons (not yet obtained by
Russian scholars). It should also be remarked that Orevkov's
destructions of $b3$ and $b6$ add no obstruction not already
covered by Viro's sporadic obstructions. Insisting again, a
careful inspection of the table shows that no more obstructions
are deducible than those already offered by the diagonal. So we
are not much advanced on the problem of deciding which patches are
algebro-geometrically realized. This is both annoying and
stimulating as our gap of knowledge raises the hope of new
constructions, namely 800 and 701 combines to the boson $b1$,
while 800 and 107 combines to the boson $b7$.

{\it Added} [01.10.13].---Looking at the bending table
(Fig.\,\ref{ViroDEGREE8_exotic_patches0_BEND:fig}), one notes
however that both those opportunities are killed if one is able to
propagate a Shustin prohibition under bending. Actually, by this
(heuristic) method the patch G2 is completely killed and dead.

\begin{figure}[h]\Figskip
%\vskip-1.2cm\penalty0
%\centering
\hskip-2.7cm\penalty0
\epsfig{figure=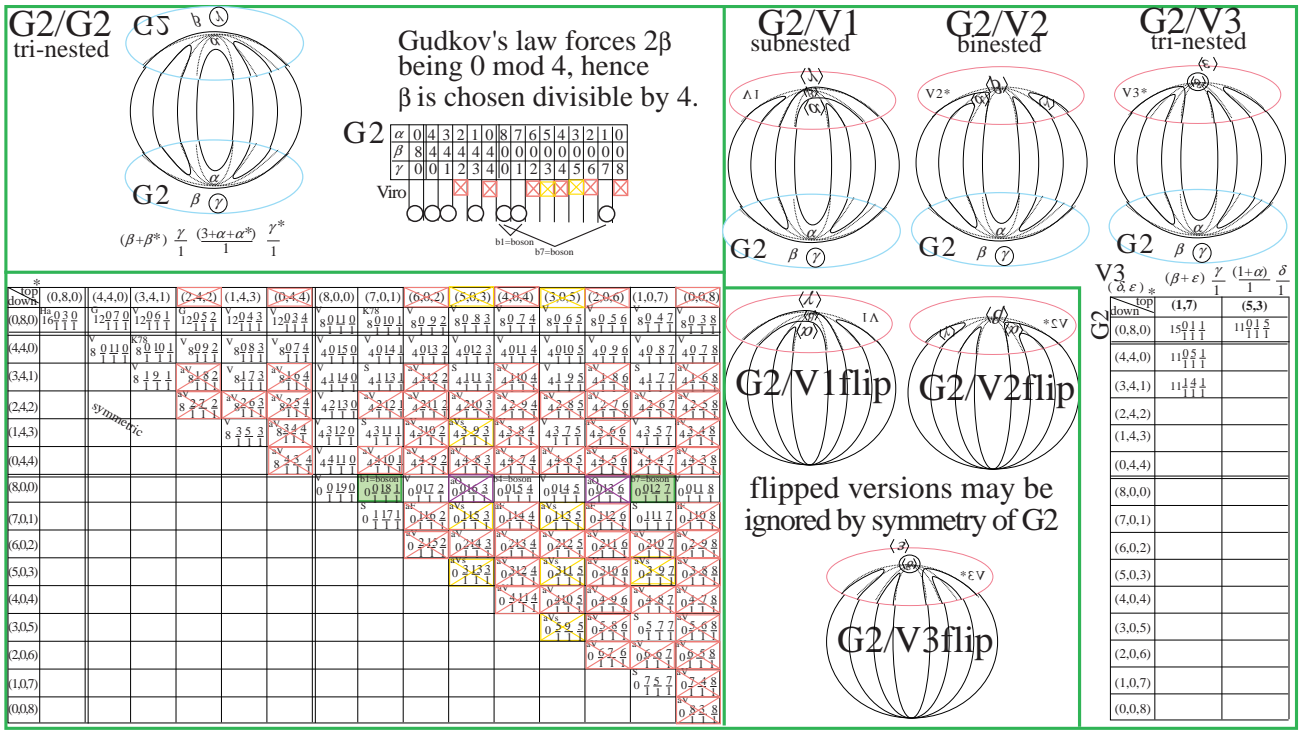,width=172mm}
\captionskipAG
  \caption{\label{ViroDEGREE8_exotic_patches_G2:fig}%
  Exotic patchwork: G2/G2}
\figskip
\end{figure}

Personally, we would found quite a pity if the patch G2 is
completely obstructed (as tacit in Viro if we interpret properly
his text).  Noteworthily,  some of our exotic patches are more
symmetric than those of Viro, and further the G2/G2 gluing looks
very much like a mitosis in
%basic
%mitochondriac living
%energy.
cellular biology. Geometric intuition may suggest that the
resulting picture is too beautiful to be omitted by nature. If so,
we could dream that even the boson $b1$ admits (via
G2(8,0,0)/G2(7,0,1))) a very symmetric realization (at least under
vertical-axis symmetry, as opposed to the  horizontal one impeded
by the asymmetries of the patch parameters employed).

As said earlier, while working out the (novel) patch-table more
carefully we found the new patch C1 of
Fig.\,\ref{ViroDEGREE8_exotic_patches0_SYS:fig}. Evidently this
has the same parameters as C2 at least as far as obstructions are
concerned. It would be of interest to look at the composition
C1/C2 to see if it affords new information.

[16.08.13] On contemplating more seriously the new patch tables,
it should be remarked that it is quite common and easy to get
realized the (capital letters) Bosons B4 and B14 (subnested)
without offending Viro sporadic or Orevkov as those obstructions
really pertains to trinested or binested schemes. In contrast the
small binested bosons $b1$, $b4$, $b7$, $b9$ are harder to
construct without corrupting Viro sporadic. Yet, table G2 is quite
remarkable for supplying legal (i.e. Viro/Orevkov licit)
constructions (hypothetical of course) of the bosons $b1$ and
$b7$. Likewise Viro's table C2=V2 (with extended parameters)
supplies logically permissible (yet still immaterialized
geometrically) constructions of the (same) bosons $b1$ and $b7$.
Hence:

\begin{Scholium}
Even if all of Viro and Orevkov sporadic obstructions are true
there is still some hope that four among the six bosons (namely
{\rm B4}, {\rm  B14} and {\rm b1}, {\rm b7}) are realized
algebro-geometrically via the most basic incarnation of Viro's
method based on the quadri-ellipse. Those bosons would then appear
as basic Kunstformen der Natur (compare optionally Ernst Haeckel's
book of drawings 1899/1904
\cite{Haeckel_1899-1904-Kunsformen-der-Natur}). Roughly speaking
any viable species must have developed a reasonable geometry of
his body-mass-index.
\end{Scholium}

Of course, the patch C1 offers the same opportunities as C2, so
that one more realization of the bosons $b1$, $b7$ is gained. As
to C1, it is evident that C1/C1 yields the same table as C2/C2,
whence the same prohibitions (and the same hypothetical
constructions). Our hope was that new information comes from
C1/C2, yet as this is a {\it mixed\/} table (as opposed to
self-composition table X/X) information is gained only if one
entry contains real patches and this is the case thanks to Viro's
theory (constructions). So obstructions will probably be induced
on C1, yet not on C2. On filling this table C1/C2, it is observed
that Viro's imparity law does not add new restrictions yet a first
one arises with Viro's sporadic anti-scheme
$4\frac{3}{1}\frac{3}{1}\frac{9}{1}$ which kills the C1-patches
$(6,0,3)$ and $(0,0,9)$. On the vertical row indexed by entry
$(9,0,0)$ we get opportunities for bosons $b1$ and $b7$, and even
$b9$ if not obstructed by the side-effect of Viro's sporadic
obstruction on $4339:=4\frac{3}{1}\frac{3}{1}\frac{9}{1}$ (the
year of the apocalypse according to
%Gregorian
%%ultra-
orthodox calendars?). So quite interestingly the boson $b9$
appears as more accessible than $b4$. Yet if Viro is true, only
$b1$ and $b7$ have real chances getting materialized via the
quadri-ellipse. From the next vertical row $(8,0,1)$, we see that
Viro's (sporadic) prohibition $\frac{1}{1}\frac{5}{1}\frac{13}{1}$
kills $(4,0,5)$. In the  603 row we meet again obstruction 4339,
yet in such a fashion that it induces only a probabilistic murder
of either C1(2,4,3) or C2(6,0,3). Next,  row 207 gives a real
chance to materialize the boson $b7$, but otherwise not more
murders of patches (essentially thanks to the presence of
Shustin's schemes). Finally the last 2 rows affords no principally
new information safe for a possible boson $b9$ (granting Orevkov
to be false), and semi-obstructions induced by sporadic
obstructions.

\begin{figure}[h]\Figskip
%\vskip-1.2cm\penalty0
%\centering
\hskip-2.7cm\penalty0
\epsfig{figure=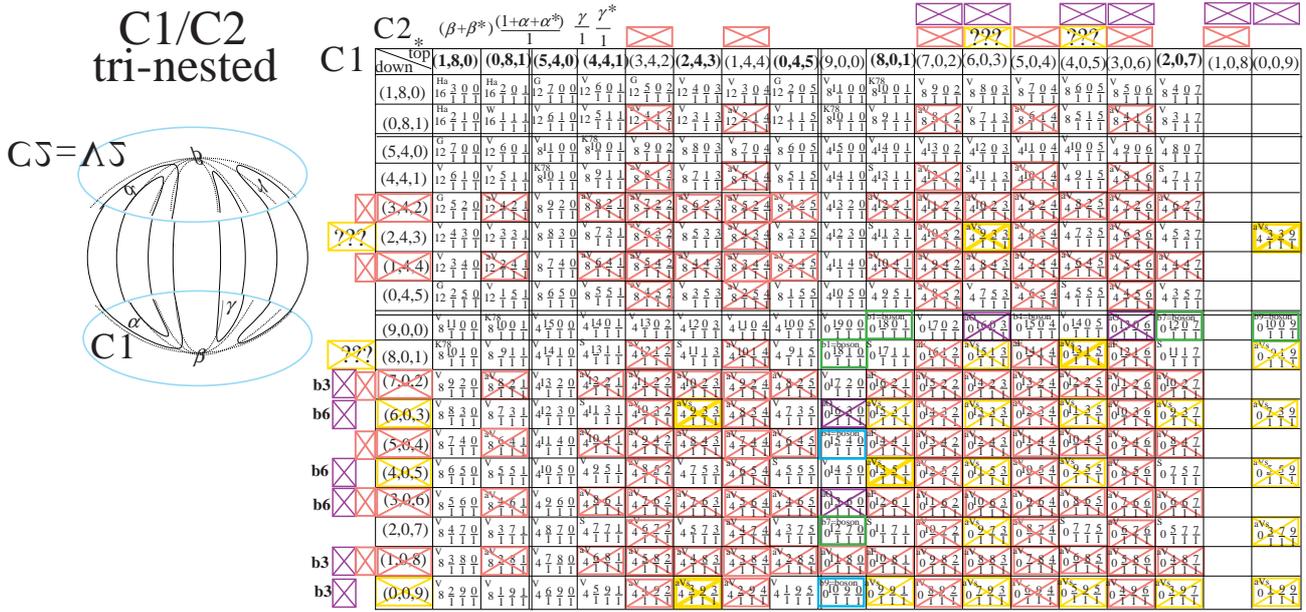,width=172mm}
\captionskipAG
  \caption{\label{ViroDEGREE8_exotic_patches_C1-C2:fig}%
  Patchworking exotic patches (continued)}
\figskip
\end{figure}

In summary the C1/C2 table affords only moderate obstructions on
patches of the class C1 and offers pseudo-construction of the
bosons $b1$ and $b7$ (but not the others). When reporting the
patches obstructions on the C1-table
(Fig.\,\ref{ViroDEGREE8_exotic_patches0_SYS:fig}) we see that the
cumulative effect of  Viro's law plus his sporadic obstructions
covers all of Orevkov's prohibitions (while forbidding actually
one more patch, namely C1(5,0,4)). So even if Orevkov is false but
all of Viro is true then the patch C1 is at least as restricted as
indicated by crosses on that table.

Then again we wondered if V2/V3flip(=C2/C3flip) produces
obstructions on V2, but apparently not according to the
corresponding table.

One could try to use the C1/C2-table to complete our knowledge of
C2=V2. Remind that only the patch C2(9,0,0) is in in doubt
(accepting the Viro+Orevkov theories).

{\it Added} [02.10.13].---As a little detail, we could dispense
Orevkov by using Viro sporadic, provided we can construct the same
C1-patches as Viro constructs in the C2-class. This holds true in
our opinion quite trivially by a simple variation of Viro's
construction  (as we shall see later in
Sec.\,\ref{Viro:patches-construction:sec}).

Hence constructing a patch by C1 could de-construct(=destroy) one
by C2. Alas the vertical row (9,0,0)  contains few schemes known
to be prohibited actually solely the 2 Orevkov schemes yet at
heights which are not constructible. Indeed the horizontal line
$(3,0,6)$ is not constructible by Fiedler-Viro, whereas $(6,0,3)$
is only killed by sporadic obstructions. Positing those to be
wrong, while declaring that C1(6,0,3) exists we deduce (granting
Orevkov's $b3$ to be an anti-scheme) that C2(9,0,0) do {\it not\/}
exist, completing thereby the dissipation theory of C2. Remember
that C1(6,0,3) is killed by disintegration of b6, so in our
scenario (existence of C1(6,0,3)) the boson b6 must materialize
and so one half of Orevkov is wrong.

All this is just sterile logical speculation yet it seems still
puzzling to treat in a systematic way the problem of deciding
which constellation of patches are logically compatible granting
the many interferences between the varied tables.

To summarize, an aspect of the game is to know precisely which
opportunities of creating new bosons arise  within the context of
the most basic Viro method (quadri-ellipse). Answering this should
by now be nearly complete via our new table of patches.

[17.08.13] Next we observe that C0 is forced to have $\ga=0$ and
therefore may be seen as a subclass of the type C1. Actually, in
the patch C0 (upon gluing with V3 and using B\'ezout) at least one
of $\be$ or $\ga$ must vanish (whatever the value of $\al$ is),
and therefore the patch family C bifurcates into the two
subspecies C1 and C2.

It seems also pleasant to combine the ground dessins in all
possible ways to visualize the ground shape of octics in the
vicinity of the quadri-ellipse (see
Fig.\,\ref{ViroDEGREE8_exotic_patches0_SYS2:fig}). Of course D is
empty safe perhaps if $\al=1$, in which case the doubled patch is
the quadri-bifolium $\frac{1}{1}\frac{1}{1}\frac{1}{1}\frac{1}{1}$
but with only 8 ovals. The hope is that some combination yields
$(M-1)$-curves where we could exploit Gudkov-Krakhnov-Kharlamov
periodicity. We find then many interesting combinations that could
induce additional obstructions. For instance C/G can be trinested
and then we have $(M-2)$-curves possibly interacting with Viro's
2nd law (impeding a trinested RKM-scheme to have all its
numerators even).

\begin{figure}[h]\Figskip
%\vskip-1.2cm\penalty0
%\centering
\hskip-2.7cm\penalty0
\epsfig{figure=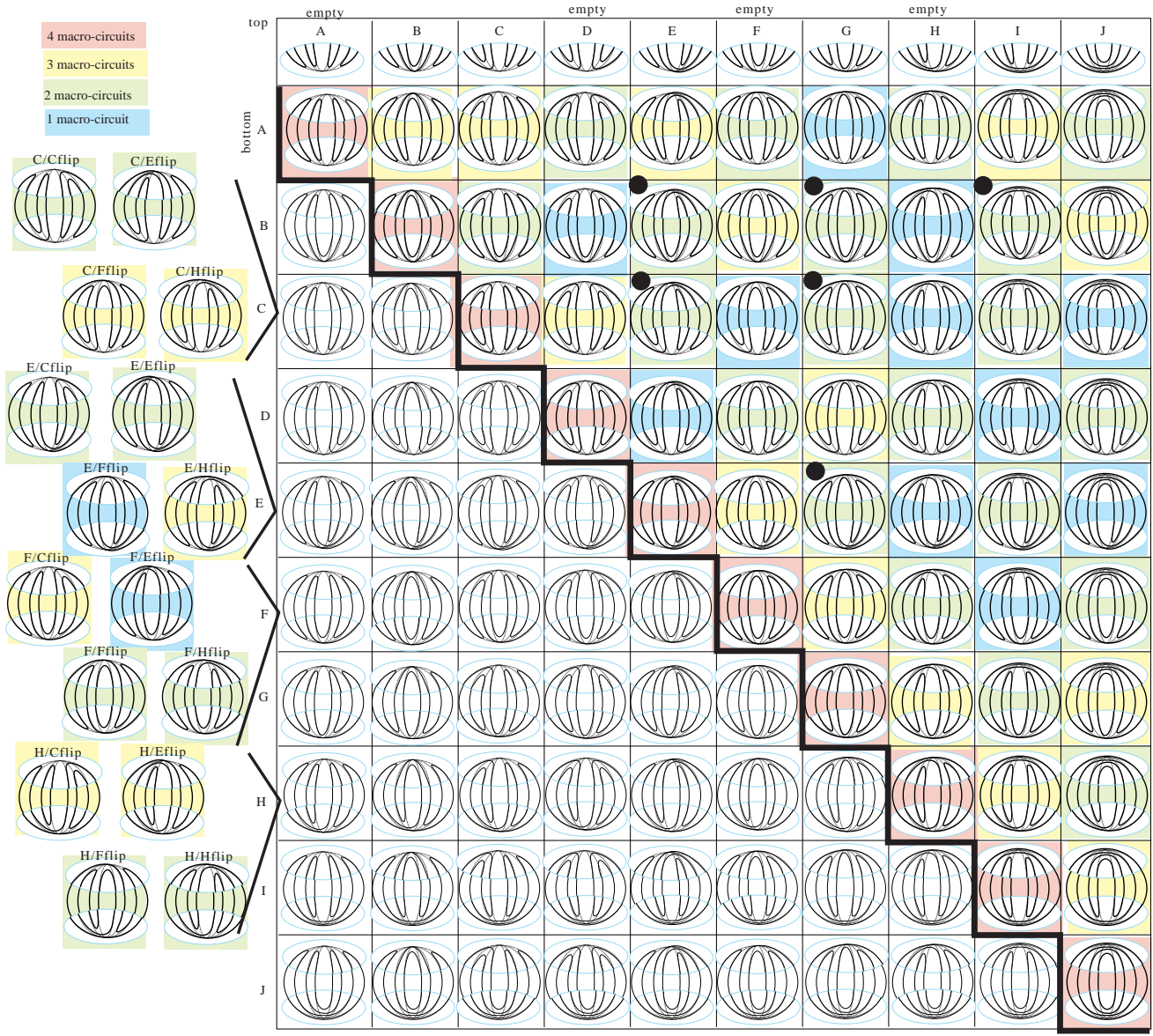,width=172mm}
\captionskipAG
  \caption{\label{ViroDEGREE8_exotic_patches0_SYS2:fig}%
  Patchworking classic (Viro's=CE) and exotic patches}
\figskip
\end{figure}

In fact if only interested in $M$-curves, we may  restrict
attention to patches (potentially) admitting a maximal dissipation
(namely B, C, E, G, I, as we saw earlier as a simple consequence
of Arnold weak-version of Gudkov periodicity). For each
combination we may study the resulting patchworks and
obstructions. This restricted table there is either 4 or 2
macro-circuits (with  3 not occurring anymore). This is a bit
disappointing yet there is still some chance to get obstructions
form Viro's 2nd law. Of course the table is symmetric about the
diagonal and so only the upper-half triangle deserves attention.

\begin{figure}[h]\Figskip
%\vskip-1.2cm\penalty0
%\centering
\hskip-2.7cm\penalty0
\epsfig{figure=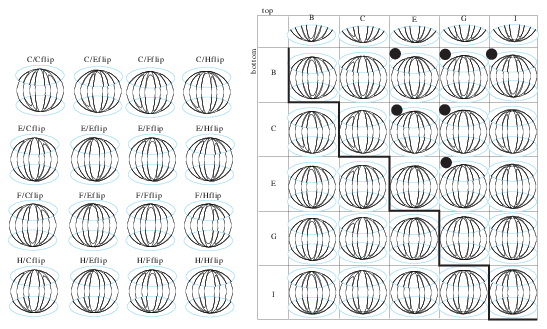,width=172mm}
\captionskipAG
  \caption{\label{ViroDEGREE8_exotic_patches0_SYS3:fig}%
  (Potentially) maximal patches}
\figskip
\end{figure}

[18.08.13] The combination B/B has already been studied, yet
perhaps not in all declinations B2/B2, B2/B3, B2/B4, B2/B5, B3/B3,
B3/B4, B3/B5, B4/B4, B4/B5, B5/B5; and additionally there is
B3/B3flip (which was already considered). However it must be
reminded that from mixed composition of virtual patches no precise
obstructions can be derived. Hence our analysis of cases looks
already complete.

For B/C, it can be trinested in the declination B4/C3 which was
already considered as table G9/V3 which offered no obstruction
(viz. interaction with Viro's 2nd law). But also B3flip/C3 is
trinested. However even without filling the table the presence of
$\de$ odd excludes any interaction with Viro's 2nd law (forbidding
all numerators of the Gudkov symbol being even).

Next we have B/E which is never trinested, even not in the
declination B4/E.

For B/G it is trinested in the declination B4/G2, but as either
$\ga$ is odd or $\de$ congruent to 1 mod 4 no interaction with
Viro's 2nd law can be expected.

For B/I, it is trinested in (and only in) the declination B4/I,
but as $\ga$ is already restricted to being odd (at least when not
zero), no interaction with Viro's 2nd law can be expected.

Next we analyze the C-row. First we have C/C which declines into
C1/C1, C1/C2, C1/C3, C2/C2, C2/C3, C3/C3. All diagonal
combinations were studied a long time ago. C1/C2 was recently
tabulated as Fig.\,\ref{ViroDEGREE8_exotic_patches_C1-C2:fig}
offering 3 new (sporadic) obstructions on C1. C1/C3 offers the
same obstructions on C1 than those coming from C2/C3, studied
earlier. C3/C3 was also taken into account but gave little
information. Next there is also in C/C flipped versions C/Cflip
which are $(M-2)$-curves. This becomes trinested in the case C1/C3
or C2/C3, but $\de$'s oddity (even 1 mod 4) impedes any direct
interaction with Viro's 2nd law. Finally, the case C3/C3flip was
already considered and yielded the strong (and complete) B\'ezout
restrictions upon C3.

Next we have to worry about C/E which is trinested in the
specialization C3/E, but as usual we cannot expect an additional
interaction with Viro's law. Further we must not forget the
flipped version C/Eflip, but this never becomes trinested even
when specialized as C3/Eflip (as there is no micro-oval population
the right lune of C3/Eflip).

For C/G it becomes trinested only in the declination C3/G2, but as
usual as both $\de$ and $\ga$ were already calibrated to odd there
is no hope to get further obstruction form Viro's 2nd law.

As to C/I it never becomes trinested.

Next we have E/E. Remind that E admits a sole declination namely
E=E(=V1) Viro's 1st dissipation mode. This was studied long ago,
and offers no obstruction as we land in the subnested real free of
all obstructions according to present knowledge. Of course we must
not miss the flipped variant E/Eflip, yet as E cannot form a
micro-island the scheme E/Eflip stays binested, and thus no
interaction with Viro's (trinested) 2nd law is expectable.

As to E/G it is at most binested in the coloration E/G2, and thus
nothing can be inferred from Viro's law. As G is symmetric (in any
coloration G1 or G2), there is no need to consider the flipped
version.

Finally (as long as E is concerned), we have E/I which is also at
most binested since the I-class cannot produce an island.

Next we have G/G where G admits 2 colorations G1 and G2. Each pure
table G1/G1 and G2/G2 were already considered, yet one must not
miss the mixed table G1/G2. Denote as usual with stars the top
parameters corresponding to G2. We meet a B\'ezout obstruction
whenever $\ga^{\ast}>0$ and at least one of both $\al$'s is
positive. Can we deduce that G2 has no patch with $\ga$ positive,
and that G1 has no patch with $\al>0$? Nearly but actually, it is
only a simultaneous realization of both conditions that would
violates B\'ezout. So the impact of G1/G2 is hard to quantify yet
it means roughly speaking that at least one of the family is empty
safe perhaps for degenerate parameters (i.e. $\al=0$ in G1 and
$\ga=0$ in G2). So we can say that G1 and G2 are coupled against
B\'ezout, but alas we do not know how to extract concrete
information  from this. Of course a construction in any one of the
category G1 vs. G2 should be the real opportunity to
%%%tranch
fix the question of deciding if the patch is rather subnested (G1
with $\al$ positive) or insulated (G2 with $\ga$ positive), but
alas it could of course be that G is none of both. Then the patch
G would reduce to G1($\al,\be \ga$) for $(\al,\be \ga)$ equal to
$(0,8,1)$, $(0,4,5)$ and $(0,0,9)$ or eventually be even smaller
(perhaps even empty).

Next we have G/I which is at most binested in the coloration G2/I.
No, actually, G2/I is anti-B\'ezout provided $\ga>0$ and $\al$ or
$\al^{\ast}$ is positive. Hence we have (warning skip this lemma
where there is a mistake, but look at the next version):

{\footnotesize

\begin{lemma}
It suffices the patch family {\rm I} containing a single
representant to force the family G2 being nearly empty. Precisely
$\ga$ and $\al$ should be both zero and therefore {\rm G2} reduces
to the single patch $(0,8,0)$. Caution: this is a mistake as we
misplaced the parameters $\al, \be, \ga$ of G2.

Likewise we can reformulate the previous token involving G1/G2, as
follows:

It suffices the family {\rm G1} being nonempty (or rather to
contain a representant with $\al>0$) to force a collapse of {\rm
G2} to the patch $(0,8,0)$. (Of course $G2(0,8,0)=G1(0,0,9)$ so
that $G2$ can be considered as empty.)

But now the ``clou'' of the argument is that the hypothetical
G2-patch can
%%%be considered
serve as the G1-patch effecting the closing and therefore it is
deduced that {\rm G2} contains at most $(0,8,0)$. Sorry, it seems
that this is rather a misconception (due to the fact that I
misplaced the micro parameters).
\end{lemma}

}

Now the corrected lemma is as follows and should be interpreted as
a reciprocity law between patches:

\begin{lemma}
It suffices for the family {\rm G1}  to contain a single
representant with $\al>0$ to force a collapse of {\rm G2} to the
patches with $\ga=0$, i.e. $(0,8,0)$, $(4,4,0)$ or $(8,0,0)$. (Of
course $G2(0,8,0)=G1(0,0,9)$, $G2(4,4,0)=G1(0,4,5)$ and
$G2(8,0,0)=G1(0,8,1)$ so that $G2$ can be considered as empty.)

Reciprocally, it suffices for the family {\rm G2} to contain a
single patch with $\ga>0$ to force a collapse of {\rm G1} to the
patches with $\al=0$, namely $(0,8,1)$, $(0,4,5)$, and $(0,0,9)$.
All those patches can be considered as element of G2 via the
formula $G1(0,\be,\ga)=G2(\be,\ga-1,0)$.

Hence as the intersection $G1 \cap G2$ reduces to the three
patches listed above ($G2(0,8,0)=G1(0,0,9)$, $G2(4,4,0)=G1(0,4,5)$
and $G2(8,0,0)=G1(0,8,1)$) we deduce that at least one of both
families G1 and G2 can be considered as empty.
\end{lemma}

 Likewise there is a reciprocity law between G2 and I.

\begin{lemma}
It suffices for the family {\rm I}  to contain a single
representant with $\al>0$ to force a collapse of {\rm G2} to the
patches with $\ga=0$, i.e. $(0,8,0)$, $(4,4,0)$ or $(8,0,0)$. (Of
course $G2(0,8,0)=G1(0,0,9)$, $G2(4,4,0)=G1(0,4,5)$ and
$G2(8,0,0)=G1(0,8,1)$ so that $G2$ can be considered as empty.)

Reciprocally, it suffices for the family {\rm G2} to contain a
single patch with $\ga>0$ to force a collapse of {\rm I} to the
patches with $\al=0$, namely $(0,1,8)$, $(0,5,4)$, and $(0,9,0)$.

Hence  at least one of both families G2 and I can be considered as
nearly empty.
\end{lemma}

In summary, either G2 contains a non trivial patch (with $\ga>0$),
in which case both G1 and I collapse to their 3 representatives
with $\al=0$, or alternatively G2 reduces to a subcollection of G1
and then there is no strong coupling and in both families G1 and I
could flourish many patches of potential subnested bosonic
interest. The first scenario (G2 non trivial) seems to favor the
materialization of the binested bosons $b1$ and $b7$, while the
second scenario (G2 trivial) corroborates rather existence of the
subnested bosons $B4$ and $B14$.

Very finally, we have I/I which only produces prohibitions already
analyzed.

It seems further that strong coupling occurs at other places like
with B2/G2, C1/G2 and C2/G2, or also E/G2. But now by Viro's
theory we know that E has patches with positive $\al$'s and thus
the patchwork E/G2 shows the:

\begin{theorem}
The patch G2 does not admit representatives with $\ga>0$, i.e.
insulated patches.
\end{theorem}

We have also a weak coupling B4/G1 from which however no tangible
obstruction can be drawn due to a sterile lack of construction \`a
la Viro in those patch classes. More generally similar couplings
arise whenever the ground curve of
Fig.\,\ref{ViroDEGREE8_exotic_patches0_SYS4:fig} contains a nest
and the corresponding patch can inject ovals in the nest while the
other patch creating an island. So we have a coupling B2/B4 from
which no concrete information can be extracted due to the merely
formal stature of our patches. Next B/G suggests two couplings,
namely B2/G2 and B4/G1. But both are only weak couplings, from
which no concrete information can be drawn. Naively, the recent
collapse of G2 (in the above theorem), suggests that in the most
plastic world there is a series of patches flourishing in G1 which
being coupled with B4 will kill many patches there (those with
$\ga>0$). In turn plasticity and the coupling B4/B2 suggest many
patches in B2 and via the coupling B2/G2 many patches are killed
in G2, in accordance with the theorem above. Of course the circle
is now complete.

Let us continue our analysis of all couplings. The next case of
nesting occurs with C/C and so we have  couplings C1/C3 and C2/C3
which are precisely forcing extinction of the micro-ovals in the
%%%%inner %%%%critical change made the 02.10.13
nested lune of both C1 and C2. Further we have the coupling
C4/C4flip (Fig.\,\ref{ViroDEGREE8_exotic_patches0_SYS:fig}) which
produces the collapse from C4 to the restricted family C3.

{\it Added} [02.10.13].---This looks true provided both $\de,
\lam$ are positive. So we erroneously ruled out the patches
C4(0,7,1) and C4(0,3,5). But, those patches are respectively equal
to C2(1,8,0) and C2(5,4,0), which are both constructed by Viro.
Thus we had to make a little correction on our catalogue (as yet
only refreshed on
Fig.\,\ref{ViroDEGREE8_exotic_patches0_BEND:fig}). Of course this
mistake has  little impact
%%{\it per se\/}
since those patches
where already catalogued in the C2-family.

Next we have C/G, but alas this does not create couplings because
neither C nor G
%able to
injects ovals in the nest (compare sub-figure C/G e.g. on
Fig.\,\ref{ViroDEGREE8_exotic_patches0_SYS:fig}).

The next case of nesting concerns E/E, and here the coupling
basically prevents the patch E to form an island, so that actually
E reduces to the single incarnation E=E=V1 (i.e. Viro's 1st
family). Note also that the flipped version E/Eflip kills the
nesting, and therefore nothing tangible can be inferred from it.

The case to come next is E/G, where we have the strong coupling
already discussed yielding all the severe obstructions on G2
incarnated by the above theorem.

Then we have G/G whose couplings implies the evident restrictions
that there cannot be micro-ovals in the lateral lunes, and also
that if there is an island growing then the central lune becomes
void too. Paraphrasing this is  the basic splitting in classes G1
and G2. Further we have the coupling G1/G2 from which we above
failed to derive decent information, which we have now gained via
the coupling with Viro's E-family. Actually, now we have the
perfect circular circuit of couplings. First, G1/G2 (with G2
nearly empty via E), then G2/B2 (with plastically B2 nearly full),
and then via the coupling B2/B4, B4 would be nearly empty, and
finally, the closing coupling B4/G1 plus plasticity would make G1
nearly full.

Next we have G/I which induces a coupling G2/I, which as G2 is
nearly empty could make I nearly full (by hypothetical
plasticity).

Lastly,  I/I  produces no coupling due the inability of I to form
an island.

\begin{figure}[h]\Figskip
%\vskip-1.2cm\penalty0
%\centering
\hskip-2.7cm\penalty0
\epsfig{figure=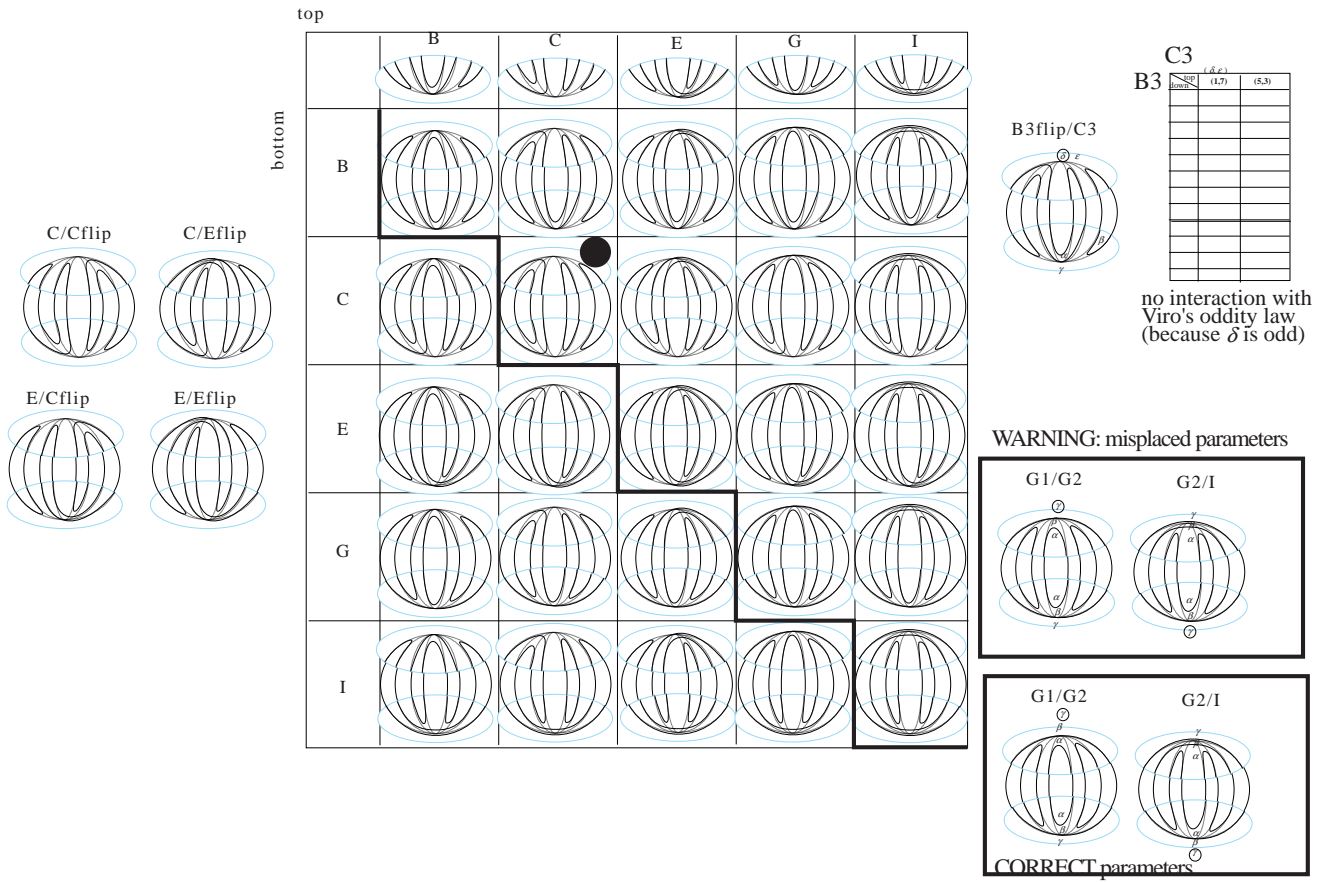,width=172mm}
\captionskipAG
  \caption{\label{ViroDEGREE8_exotic_patches0_SYS4:fig}%
  Admissible patches}
\figskip
\end{figure}

At this stage it seems that we have analyzed all logically
possible couplings, and thus established an exhaustive list of
patch restrictions derivable by the naive composition method (i.e.
formal patchwork based on the superstition of independency of
smoothing of both singularities of the quadri-ellipse). From this
analysis, it results one major prohibition on G2, as well a cyclic
structure on couplings along which the density of patches could
alternate with periodicity two. This is just to say, one is nearly
empty and the successor nearly full, and so on. In particular G2
is nearly empty and so would be B4 (yet which produced no bosons).
In contrast we expect that G1, B2, and I are nearly full thereby
contributing to the materialization of the boson B4 and B14 in the
subnested realm where many constructions post-Viro were supplied
by Korchagin, Chevallier, Orevkov, and where up-to-date nobody
could find any obstruction except Shustin in the case of zero
outer ovals.

So despite the recent basic B\'ezoutian destruction of bosons via
G2/G2 (notably $b1$ and $b7$), there is still some hope
%borne out
for their realizations via C1/C1 (whose composition table is the
same as C2/C2 via an evident isotopy), or C2/C2, or even via
C1/C2.

Further there is also a coupling C0/C3 forcing $\ga=0$ on C0=G2
yielding some other opportunities for bosons (yet violating some
of Viro's sporadic obstructions). Actually the few surviving
C0-patch are readily covered by family C1 so that we cannot speak
of a principally new realization.

[19.08.13] Interestingly,  simple B\'ezout obstructions  often
recover all the deep Fiedler-Viro obstructions, as for instance
with C4$\approx$C3  or G2-patches (compare
Fig.\,\ref{ViroDEGREE8_exotic_patches0_SYS:fig}). One may wonder
if all patches obstructions (about $X_{21}$) can be subsumed to
B\'ezout directly, yet this looks quite unlikely as we proposed
(we believe) an exhaustive search of all couplings relation
entertained by patches.

Further we remind that there is another anti-B\'ezout coupling
with B4/B5 which is generically quadri-nested, and when not it
turns anti-Gudkov as we saw in an earlier table (G3/G9 in older
notation). Unfortunately this coupling relates exotic patches and
therefore there is not enough grip to infer any concrete
information from it.

Next, there is also the coupling B2/B5 and even B2/B3, yet all
those are weak couplings and it looks hard to infer any concrete
information. So actually B2 is coupled  with all others B's, i.e.
B3, B4, and B5. But it seems that the loose (sparse) information
we have is caused by the fact (compare
Fig.\,\ref{ViroDEGREE8_exotic_patches0_SYS4:fig}) that the patch B
does not interact in a nested way with the real patches of Viro
(letters C and E in our catalogue). This absence of nest on the
ground figure explains why we fail meeting a B\'ezout obstruction
involving a nest of depth 3 plus one of depth 2. In the case of no
nest on the ground naked figure (prior to the  addition of islands
and micro-ovals, cf. again
Fig.\,\ref{ViroDEGREE8_exotic_patches0_SYS4:fig}) we can still
expect to find a B\'ezout obstruction involving saturation of the
bi-quadrifolium $\frac{1}{1}\frac{1}{1}\frac{1}{1}\frac{1}{1}$
(i.e. 4 nests of depth 2). So for instance there could be a
coupling between B with an island and C insulated as well, that is
between B4 and C3, but as B4 has empty lateral lunes (B\'ezout
applied to B4/B4flip) it is seen that the right lune of B/C stays
empty and the full scheme becomes generically trinested, yet not
quadri-nested as initially expected.

A similar discussion shows that there is no conical (equivalently
quadri-nested) coupling when pairing B with E whose ground figure
is a lune plus a snail (cf. still
Fig.\,\ref{ViroDEGREE8_exotic_patches0_SYS4:fig}). Here E cannot
produce an island nor fill the unnested lune with micro-ovals, and
so fails to do
%%%also
the insulated incarnation of the patch B.
Hence B/E is at most binested.

So it seems that our poor understanding of the B-patches is caused
by a lamentable geometric interaction of the B-patch with both C
and E the (fundamentally European) collections of Viro.
(Memnotechnic trick: CE=communaut\'e europ\'eenne.)

Also at this stage we enumerated B\'ezout obstruction for lines
and conics, and one may wonder if there are obstructions induced
by
%%ancillary
auxiliary cubics.

\section{Toward sophisticated B\'ezout obstructions}

\subsection{Speculating about obstructions}

[19.08.13] Perhaps one could so (via cubics) obstruct some of the
bosons yet we believe that the general method of total reality
(alias the Riemann-Schottky-Bieberbach-Grunsky theorem, perhaps in
the synthetic variant of Gabard 2013B
\cite{Gabard_2013B-Riemann's-flirt}) should be the true weapon to
detect additional prohibitions if there is any (in degree $m=8$).
Of course this has be mixed perhaps with standard homological
methods (construction of $2$-cycles$\approx$membranes) designed by
Arnold, Rohlin, Viro, Fiedler, etc. or link theory \`a la Gilmer,
Orevkov. For instance beside Arnold's surface which can fail to be
orientable there is a myriad of natural membranes, e.g. the Rohlin
surface obtained by filling all ovals by their bounding disc. This
is in the $M$-curve case represented by a singular sphere, whose
self intersection yields Rohlin's complex orientation formula. As
we said often, it seems that total reality offers a sort of
transverse structure (a bit like Haefliger, etc.) and so perhaps a
good deal of Hilbert's puzzle can be tractable in the context of
holomorphic foliation theory \`a la Painlev\'e et ali (i.e. Brazil
and Dijon, like Cerveau, Camacho, etc.).

An idea (that flashed us ca. 1 month ago) is that given an
$M$-curve of degree $m$ (say $m=8$) we can look at all pencils of
$(m-2)$-tics, which have degree 2 less and which we shall call of
co-degree 2, hence co-conics (so coconuts or just nics). By the
elementary argument in Gabard 2013B
\cite{Gabard_2013B-Riemann's-flirt}, we know that there is always
such a pencil which is totally real. So we could look in the
Grassmannian parametrizing all those pencils at the sub-body
consisting of total pencils. This must have a marvellous geometry,
especially when it comes to look at the boundary of the body.

{\it Added} [02.10.13].---Consider the trivial case ($m=3$) of an
$M$-cubic ($r=2$ circuits) swept out by a pencil of lines. Then
total reality holds iff the center of perspective is located
inside the oval of the cubic, or its boundary. In the latter case
the degree of the total map lowers to 2 (instead of 3 when looking
from inside). In general, it seems evident by analogy that low
degree total maps will emerge along the boundary of the body of
all total pencils. Okay but actually the recipe of total reality
on plane $M$-curves (as in Gabard 2013B
\cite{Gabard_2013B-Riemann's-flirt}) readily gives such maps of
lowest possible degree. So its seems natural to expect that those
maps describes explicitly the boundary of the body of all total
maps.

Studying all this very precisely should perhaps advance the
resolution of Hilbert's problem in degree $m=8$ and higher. Alas,
it is also evident that  several tour de forces are requested.

It is only now that a  basic aspect came transparent to us. It is
clear that real algebraic curves like nesting but not
%%%too much
excessive nesting as there is evident B\'ezout bounds upon the
nesting complexity. So an octic cannot be quadri-nested, and when
it is it reduces to the quadri-nest (alias bi-quadrifolium).
Likewise when trinested the schemes suffer severe Fiedler-Viro
prohibitions. In the binested case there is only for the moment
two (striking) prohibitions of Orevkov upon $b3$ and $b6$, which
less surprisingly pertains to curves with the minimum number of
outer ovals. It is clear that one could hope to attack the problem
via B\'ezout for cubics (a bit like along the strategies of Le
Touz\'e).

The basic idea is that if the octic is much nested (e.g. in the
binested bosonic range $1\frac{x}{1}\frac{y}{1}$ with $x+y=19$ in
order to have a total number of $M=22$ ovals) then one could
select in both nests a quadruplets of points and let pass through
the resulting 8 points a connected rational cubics. Remember that
through 8 points there is a pencil of cubics containing among the
12 complex singular curve at least 8 which are real, and so we
find a curve of the desired type. The trick would be then to
control somehow a salesman travelling between those group of 4
points so as to create excessive intersection. For instance if we
could arrange 8 transitions from white to black (being the colors
of both nests) then it would result $8\cdot 4=32$ intersections
overwhelming B\'ezout's $3\dot 8=24$. Seven transitions would be
enough, but looks hard to have an odd number of transition. Six
instead would not be enough to corrupt B\'ezout.
%Of course
It would be of paramount importance to understand if some
obstruction in the bosonic range can be drawn from this simple
device, especially if it recover the Orevkov obstructions on $b3$
and $b6$ (i.e. $x=3$ and 6 respectively).

The point is that the rational cubic creates 4 intersections when
it salesman travels from black to white, but still create 2
intersections when linking monochromatic points (i.e. in the same
nest). Thus six transitions is enough, affording
$4.6+2.2=24+4=28$ intersections, while 4 transitions produce only
$4.4+2.4=16+8=24$ intersections without corrupting B\'ezout.

Then one can imagine the 8 basepoints colored black and white
according to the splitting $8=3+5$ instead of $4+4$. Then there is
still enough room to get 6 transitions and we could so perhaps
reprove Orevkov's obstruction of $1\frac{3}{1}\frac{16}{1}$.
Indeed the technique would be to choose 3 points inside the ovals
of the small nest and to choose 5 of them inside the big nest
containing 16 eggs(=empty ovals), and to pass a connected cubic
through those 8 points in such a way that there is 6 transitions
from black to white. Then as before $4.6+2.2=28$ many
intersections are granted in $C_3\cap C_8$ and B\'ezout is
overwhelmed.

But of course in this sort of games the crude theory is very easy
yet the practice is very hard to polish. We mean of course that
the argument should not rules out schemes constructed by Viro.

The method could equally apply to the trinested case,
 in principle more elementary just for historical
reasons (Fiedler-Viro came prior to Orevkov). Here we have three
colors say black, white and red. The 8 points are distributed
along the cubic circuit and since more color are available there
should also be more transitions, while having 6 of them would suit
our desire. Of course in the worst case it is possible for this
distribution of colors to be very monotonic and forcing only 3
color changes.

It would be fascinating if it is so possible  to reprove the
Fiedler-Viro oddity  law, or some of Viro's sporadic rules (or
decree) (e.g. the famous
4339:=$4\frac{3}{1}\frac{3}{1}\frac{9}{1}$).

\begin{figure}[h]\Figskip
%\vskip-1.2cm\penalty0
%\centering
\hskip-2.7cm\penalty0
\epsfig{figure=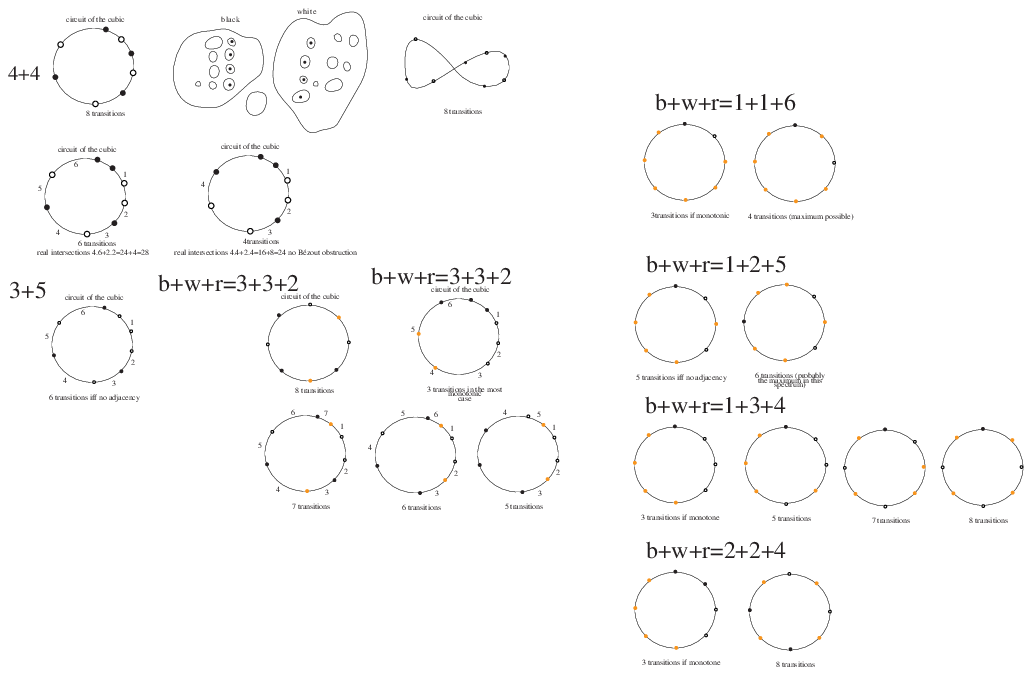,width=172mm} \captionskipAG
  \caption{\label{ViroDEGREE8_Cubics:fig}%
  Salesman travelling with cubics}
\figskip
\end{figure}

First the 8 points may be partitioned into the 3 colors as
$1+1+6$, $1+2+5$, $1+3+4$, $2+2+4$, $2+3+3$ (alias color
spectrum). In each case the monotonic distribution of colors
implies only 3 transitions, hence only $4.3+2.5=12+10=22<24=3.8$
intersections. What about 5 transitions? Then there is
$4.5+2.3=26$ intersections and this is the minimal number
corrupting B\'ezout (as 25 is not possible being odd while the
octic circuit is even by the M\"obius-von Staudt law). So in the
case of 3 colors we can observe an odd number of transitions the
minimal of them corrupting B\'ezout being 5.

Now given an $M$-scheme of degree $m=8$ which is trinested one can
choose 8 basepoints among the deepest ovals, and in such a clever
way as to maximize the number of transition. Of course by B\'ezout
the number of transition is bounded by 4, yet we may hope that
topological reasons makes this sometimes higher.

Let us assume the $M$-scheme $t\frac{x}{1}\frac{y}{1}\frac{z}{1}$
to have one even numerator (number of eggs in a nonempty oval).
Then one other numerator has to be even (and  nonzero), thus we
can employ the painting $2+2+4$, except on the 1st layer of the
pyramid (Fig.\,\ref{SIMPLIFIED-TABLE_gurus:fig}). In fact $1+2+5$
can be injected in all trinested schemes safe
$12\frac{1}{1}\frac{2}{1}\frac{4}{1}$.

Of course the (color) spectrum $1+1+6$ admits at most 4
transitions and so is useless to corrupt B\'ezout. The spectrum
$1+2+5$ can have 5 and even 6 transitions (which is probably the
maximum possible). Similarly the spectrum $1+3+4$ can reach 8
transitions, and so do its companions $2+2+4$ and $2+3+3$.

So as soon as we can employ the color spectrum $1+2+5$ or the
higher more colorful avatars, we can expect that for a suitable
rational member of the pencil of cubics assigned to visit the 8
basepoints the allied colorimetry (induced by the 3 nests of the
$M$-octic) will assume high chromodynamic level, and so B\'ezout
will be corrupted.

Basically this is rather plausible just for statistical reason
that the distribution of colors (on the ground circle) will be
fairly random and so likely to exhibit 5 or more transitions.

One possible scenario could be that for the  first spectrum
admitting 5 transitions, namely $1+2+5$ one can always find a
singular cubic with (at least) 5 transitions through any 8 points
colored along this spectrum. If true, this would explain Viro's
oddity law safe apparently for
$12\frac{1}{1}\frac{2}{1}\frac{4}{1}$. Alas, our hypothesis would
also destroy the scheme $8\frac{1}{1}\frac{3}{1}\frac{7}{1}$ (the
first constructible into which $1+2+5$ may be injected). This
being constructed by Viro we see that our hypothesis is just
superstition.

Let us look at the next spectra $1+3+4$, $2+2+4$, $2+3+3$. Assume
again for this first one $1+3+4$, a universal law (involving
merely configuration of 8 points and cubics) of chromodynamics
telling that there is, for any distribution of 8 basepoints
colored in this fashion, a rational connected cubic interpolating
the eight points with at least 5 transitions (of colors). If true
universally this would also kill Viro's scheme
$8(1,3,7):=8\frac{1}{1}\frac{3}{1}\frac{7}{1}$.

Positing the same law for $2+2+4$ would kill Viro's scheme
$8(3,3,5)$ or $4(3,5,7)$. Finally, this bad state-of-affairs is
not arranged when looking at the last spectrum available $2+3+3$.

So it seems that there is no universal law of chromodynamical
excitation, at least in the vacuum, i.e. regardless of the
distribution of 8 basepoints.

So we seems blocked. One reaction could be that cubics are not
flexible enough with their 8 basepoint assignable to visit
maximally the 19 empty ovals of a trinested $M$-octic. Especially
important is the case where there is no outer ovals, because it is
maximally obstructed apart from 3 constructions due to Shustin.
For such schemes it seems natural to impose 19 basepoints. To
lines we can impose 2 points, to conics 5=2+3 points, to cubics
9=5+4 points, to quartics 9+5=14 points and to quintics 14+6=20
point. However for a pencil we have precisely 19 basepoints and so
we can expect a singular quintic with less ovals than $M_5=7$.

Then we apply the same methodology. We split first the 19 points
in three colors. Here there is plenty of such partitions. First
$19=1+1+17$ up to $19=6+6+7$ which the most energetical one, i.e.
anti-capitalist and best distributed. Now the problem is that the
singular quintic of the pencil has 6 circuits, and a priori the 19
colored points can land monochromatically into those circuits. In
this case there is no transition and only 19.2=38<5.8=40 real
intersections are granted. Yet, we see that if we manage to gain a
bit more chromatism we are going to violate B\'ezout.

Again we assume the scheme trinested without outer ovals (i.e.
Gudkov symbol $\frac{x}{1}\frac{y}{1}\frac{z}{1}$ with
$x+y+z=19$). We consider the corresponding spectrum $x+y+z$, and
distribute the 19 basepoints among the 19 empty ovals of the
$C_8$. We have a corresponding pencil of quintics interpolating
those 19 points, and $3(5-1)^2=3.16=48$ members of it will hit the
discriminant (over the complexes at least). It seems evident that
there is at least one real singular member in the pencil, and let
us assume that there is even one member which is dichromatic in
the sense that two different colors lands in the same oval of the
$C_5$. Then we have two transitions at least and so
$2.4+17.2=8+34=42>40$ real intersections and B\'ezout is
corrupted.

Of course this scenario in abstracto would conflict with Shustin's
3 constructions which are perhaps wrong albeit this is quite
unlikely.

The methodology employed here is just the classical trick of
``interpolation through the deep nests'' and was used
systematically by Zeuthen, Harnack, Hilbert. We are just now
trying to see if it can explain most of Viro's sporadic
prohibitions. So our strategy is quintic as tool to interpolate
deep nest of M-octics.

First it should be noted that the number of transition of 19
points distributed on the at most 7 (6 if singular) circuits of
the quintics will be equal to the number of points regardless of
the fact that $C_5$ is not anymore a single circle like our
previous cubic. A transition can of course just be interpreted as
the arc resulting from cutting along the points; and in a
(triangulated) circle (or more generally a compact Hausdorff
one-manifold) there is always a bijection between edges and
vertices given e.g. by an orientation.

Actually, it seems therefore not even indispensable to lower the
number of ovals of the interpolating quintic by acquisition of a
singular point. What is crucial is rather the dichromatism
effecting that 2 points belonging to different nests of  the $C_8$
lands in the same component of a suitable $C_5$. Actually if this
is not the case then all quintics of the pencil (considered as a
dynamical \"Olfleck) would effect (Morse) juncture only between
themselves. Imagine so 3 groups of $x$, $y$ and $z$ many points
summing to 19($=x+y+z$) and the corresponding pencil of quintics
through them (which is unique under harmless genericity
assumptions). Then as time evolves the initial (say smooth)
quintic $C_5$ is deformed and at some stage it seems forced that
there is a conjunction of 2 ovals coalescing together in which
case we would have right after the critical level a quintic curve
with dichromatism.

Of course this looks the psychologically simplest phenomenon, yet
perhaps it is not a necessity. (Actually, the existence of
Shustin's three maximally trinested curves, i.e. no outer ovals,
incarnates an obstacle along our scenario of forced collision
between ovals belonging to different colors).

Finally the method adapts (nearly mutatis mutandis) to the
binested case where the game is still open (as we are in the
bosonic strip of Fig.\,\ref{SIMPLIFIED-TABLE_gurus:fig}). Here we
just have to split the 19 basepoints in two colors instead of 3,
and the same dichromatism phenomenon would an obstruction of the
corresponding schemes. Yet all the art is to do this without
conflicting with the 3 schemes constructed by Viro in this bosonic
range. So as before the method need to be refined, or perhaps it
is true in brute force generality in which case few of Viro's and
Shustin's constructions would be erroneous. As to Viro it is even
more unlikely as his construction relies on the quadri-ellipse,
yet for patches coming from the C/C combination (as E/E is
subnested as seen on
Fig.\,\ref{ViroDEGREE8_exotic_patches0_SYS3:fig}), especially as
C2/C2 for the 2 admissible parameters with $\be=0$. Can we imagine
that this parameters are killed (i.e. that Viro is wrong when
claiming their existence)? Probably, not yet our understanding is
so weak that we cannot exclude this  option for the moment.

The philosophical principle that could make some of the above
argument work, is that {\it any pencil is color-mixing\/} (like in
fine arts).

\subsection{The method of the deepest penetration}

[20.08.13] It is clear that the previous method can be declined in
several contexts depending on the degree of the interpolating
curves. Basically, we can fix our attention on the case of
$M$-octics, and look at varied interpolating curves of degree
either 1, 2, 3, 4, 5, 6. Degree 1, and 2 yields basic B\'ezout
obstructions on the depth of nest and the maximality of the
bi-quadrifolium. Degree 3 was as far as we know never successfully
exploited to draw an obstruction on octic.

In the former section, we explained how this could be used by
imposing 8 basepoints while still having a pencil. Now one can
also impose directly 9 basepoints inside the deepest ovals. One
get then (generically) a smooth cubic with 2 circuits. If our
octic is purely trinested (i.e. no outer ovals) then we have 3
colors corresponding to the three nests. By the pigeonhole
principle two distinct colors must land in the same circuit of the
cubic. Therefore we have two color-transitions, at least, and
therefore $7.2+2.4=14+8=22<24=3.8$ intersection granted. So we
need more chromatism.

So assume given 9 points. We suppose given a color-spectrum of 3
colors B, W, R (black, white, red, say). Those are in
correspondence with partitions of nine of length three, i.e.:
$9=1+1+7$, $9=1+2+6$, $9=1+3+5$, $9=1+4+4$, $9=2+2+5$, $9=2+3+4$,
$9=3+3+3$. One must imagine those nine colors  falling into the 2
circuits of the $C_3$. As three of them are distinct, at least two
must land in the same circuit creating two color-shift
(transition). However we need more than that to corrupt B\'ezout,
namely four shifts as then we have $5.2+4.4=10+16=26>24=3.8$ many
intersections.

We could posit that any purely-trinested $M$-octic admits a
distribution of 9 points among the deepest 19 ovals which has at
least 4 color-transitions. Even more than that we could suppose
that the four-color principle holds true universally for cubics
without reference to any octic, but this looks hazardous as we may
choose the 9 points on a given cubic while choosing the coloration
very monotonically e.g. one black and one white point on the oval
of $C_3$, and all remaining 7 red points on the pseudoline: then
there is only 2 color-changes, instead of the 4 desired.

By the way even the version conditioned by an octic cannot be true
universally at least without conflicting with 3 of Shustin's
constructions.

Further as an additional technical difficulty, it seems that if
one of the (nine) points lands alone on a circuit then this means
roughly that the cubic has a micro-oval visiting only one inner
basepoint without having to cross necessarily the oval to salesman
travel in another oval. So we can even loose one of the
weighted-by-2 intersection, and all our count can be jeopardized.
This phenomenon is a traditional difficulty in the field (which we
call the {\it traquenard of mini-ovals\/}).

Next, the method of the cubics can be adapted to the case of 2
colors, and then we may expect to derive old (Orevkov's) or new
obstructions in the bosonic strip
(Fig.\,\ref{SIMPLIFIED-TABLE_gurus:fig}). Again we need now 4
color-shift to corrupt B\'ezout, and this cannot hold universally
without corrupting three constructions by Viro.

Next we can augment the degree of the interpolating curve, first
to 4 and the 5 or 6. The gain is that we can penetrate through
more basepoint as we have with increasing degree more freedom in
assigning basepoints, or should perhaps rather say anchor points
when we choose so many as to have a single curve. (This is the
statical penetration method as opposed to the dynamical one using
a whole pencil.) Of course the price to pay is that with
increasing degree the interpolating curve has a priori more ovals
(potentially as many as Harnack's bound), and it becomes harder to
ensure chromatism, i.e. color changes can be vacant in case all
three colors lands in different circuits.

For an interpolating quartic, we may impose $2+3+4+5=14$ points
distributed among the 19 deep ovals of a trinested $M$-octic. If
the latter is purely trinested (i.e. no outer ovals) then we get
so, assuming $\tau$ many transitions, $(14-\tau).2+\tau.4$
intersections which exceeds B\'ezout's $32=4.8$ as soon as
$(14-\tau).2+\tau.4=34$, i.e. $2\tau=34-28=6$, that is $\tau=3$
transitions. Alas, a priori the 3 colors can fall apart in the 4
ovals of the $C_4$ without any chromatic interaction, and then
$\tau$ is as low as zero.

In contrast we can posit some higher intelligence able to show the
existence of at least 3 color-changes. Yet, as before the trick of
choosing the 14 points on a given $M$-quartic in a very
monochromatic fashion (say 1 black point on one oval and one white
point on another oval, plus the 12 remaining ones on the same
oval), yields a distribution for which the interpolating quartic
has zero color-change.

Then we can move to quintics. There the space of coefficients has
dimension $\binom{5+2}{2}=\frac{7.6}{2}=7.3=21$, and we can impose
20 anchor points. In first approximation, we may choose them
inside the 19 deep ovals, but there is now one more point
available for which there is no preferred position. Here it may
seem that the dynamical variant involving a pencil was better
suited as we used the idea that any pencil is color-mixing. Maybe
the 20-th (twentieth) point should be chosen externally of the 3
nests of the $C_8$, and it may thus be considered as belonging to
a 4th color (say G=green). A color-transition to green imposes
only 2 (instead of 4) real intersections, and so out twenty points
grants only $20.2=40=5.8$ intersections without provoking
B\'ezout.

Actually this count ignore the traquenard, and can be cleaned by
imposing the anchor points on the $C_8$ as opposed to inside their
ovals. Then as all circuits of an octic are even (null-homotopic)
we gain one more intersection on each twenty marked ovals and so a
total of 40. Now this holds true for any choice of 20 points
(injectively) distributed on the 22 ovals regardless of marking
primordially the 19 deep ovals. One could hope that there is a
special turbo-injection of 20 such points such that more
intersection are gained.

Maybe first note that the interpolating quintic cannot intercept
the two unmarked ovals, because the maximum number of 40 is
already reached by the boni intersections  given M\"obius-von
Staudt. So if we suppose given a purely trinested $M$-octic
(visualizable as 3 Swiss cheeses of the type Emmenthal) and if we
fix the marking of 20 on all 19 empty ovals plus one nonvoid oval,
we see that the 2 remaining nonvoid ovals are trapping the quintic
ovals which cannot intersect them. We call them therefore
fundamental barriers to the proliferation of ovals. We can now
imagine that the 20 points are dragged inside their respective
ovals moving therefore in a 20-dimensional torus. It can be
imagined that for some special position the constellation of 20
points fails to impose independent conditions upon quintics, or
paraphrasing that we have suddenly a pencil of quintics through
our twenty points. (Observationally this is somehow reminiscent of
the solar magneto-hydrodynamics with flow lines becoming so
distorted that a finally violent global rupture is then necessary
causing explosions, which forms the famous auroras borealis when
reaching the terrestrial atmosphere.) If so is the case, we can
further impose our quintic to visit one of the two barriers ovals
and we get a contradiction with B\'ezout.

If this argument works universally then we get a proof of most of
Viro's sporadic obstructions but alas also a disproof of 3 of
Shustin's constructions. A little look at
Fig.\,\ref{SIMPLIFIED-TABLE_gurus:fig} shows that there is exactly
$9+7+6+4+3+1=30$ purely trinested $M$-schemes, of
which---according to the Germano-Russian pact
(Fiedler-Viro-Shustin)---only 3 of them are constructible: all the
others being prohibited. So in some probabilistic sense our
argument is true with probability of 90 percents, or universally
true in case Shustin's constructions are wrong (albeit the seem
rather plausible, compare our
Fig.\,\ref{ViroDEGREE8_SHUSTIN:fig}).

If we assume that each twenty-points $g_{20}$ on the $C_8$
determines a unique quintic we get a continuous mapping to $\vert
5H \vert$ the hyperspace of all quintics isomorphic to $\PP^{20}$.
The source of the mapping is essentially a torus once we restrict
the location of the marking on some definite 20 ovals among the 22
available. Naively one could expect the mapping $T^{20}\to \RR
P^{20}$ being surjective, but then all quintics are swept out: so
in particular one visiting the barrier-ovals and we corrupt
B\'ezout.

Of course if there is at least two quintics interpolating the
group of 20 points $g=g_{20}$, then the spanned pencil is also
interpolating the same data.

So again, given a purely trinested $M$-octic we mark 20 ovals on
it (e.g. by omitting two nonempty ovals, which we call the
barriers). We consider for each 20-tuple distribution on those
ovals, a quintic interpolating them (which exists by basic linear
algebra). We would like to show that this $C_5$ is not always
unique. If so, then we can impose---additionally to the already 40
granted intersections situated on the 20 marked ovals---one more
intersection by forcing to visit one of the two available
barriers. B\'ezout is then contradicted.

Assume the contrary (i.e. perpetual uniqueness). Then we can
define a mapping $T^{20}\to \vert 5 H \vert \approx \RR P^{20}$,
which cannot be surjective (else we could impose a visit of the
barrier). On the other hand it could be hoped that a homological
mapping degree argument \`a la Brouwer prompts  surjectivity of
the mapping in case the top-dimensional representation $H_{20}$ on
homology is non-zero. Of course as $\RR P^{20}$ is nonorientable
we must confine on homology modulo 2. Alas, it seems hard to tell
anything on this mapping degree without penetrating better into
the geometry of the map. We can still try to imagine each
interpolating quintic through the $g_{20}$ as this group of 20
points varies along the 20 marked ovals. One could argue that the
subset of $\RR P^2$ swept out by the collection of all those
quintics is open (by a balayage argument) and compact (by general
topology, plus a simple fibering argument) and therefore a nonvoid
clopen (=closed open set) in the connected set $\RR P^2$.
Therefore the sweeping set is full yet this contradicts the fact
that the barriers cannot be visited (by interpolating quintics).

This contradiction would prove the:

\begin{Scholium} The interpolating quintic cannot be perpetually
unique for any location of the group $g_{20}$ of twenty markers
distributed on the ovals of a purely trinested $M$-octics. But
then we can impose a visit through the barrier and B\'ezout is
foiled. In conclusion there would be not a single trinested
$M$-octics, jeopardizing thereby construction by Shustin.
\end{Scholium}

In fact, the argument would nearly work as well regardless of the
ovals distribution. It is only essential to be able to mark 20
ovals, and thus our scholium seems to kill all $M$-octics which
seems a bit too apocalyptic if one believes in the elementary
construction of Harnack/Hilbert.

Of course it  could be that there is a basic mistake in what we
called above the balayage argument. A priori as parameters varies
the curve $C_5$ moves, but it can move like a wave front coming
back and forth and thereby not sweeping an open domain, but doing
rather what Whitney calls a fold.

Actually, it seems that the stable portion of our reasoning gives
the:

\begin{lemma} Any group of 20 points injectively distributed
on twenty ovals of an $M$-octic imposes independent conditions on
quintics, and therefore determines unambiguously a unique quintic
interpolating those points. Actually the assertion holds as well
for $(M-1)$-curves.
\end{lemma}

\begin{proof} Choose a $C_5$ through the 20 points of the $C_8$.
By M\"obius-von Staudt there is one more intersection on each of
the 20 ovals reaching thereby already the maximum permissible of
40. But there is one more (so-called barrier) oval on the $C_8$
(provided it is at least an $(M-1)$-curve), and so we can---in
case of non-uniqueness---impose additionally to the interpolating
quintics to visit one of the barrier ovals, but then B\'ezout is
contradicted.
\end{proof}

Intuitively this means that everything is very stable and there is
no solar irruption causing sudden jumps in the dimension of the
space of interpolating quintic.

Note this being valid universally independently of the oval
distributions. Yet, through the work of Fiedler-Viro-Orevkov we
expect severe prohibitions in the maximally nested cases (i.e.
binested with one outer oval and trinested without outer oval).
Can we detect them by our naive B\'ezout style approach, i.e. via
an elaboration of the above lemma?

By the lemma (uniqueness of the quintic interpolating a
distribution) we have a continuous mapping from varied tori
(amounting to the $22.21/2=11.21=231$ markings of twenty ovals
among the 22) to $\vert 5 H\vert$ the hyperspace of quintics.
Further each of the interpolating quintics avoids two
barrier-ovals.

It may be observed that given any group $g_{20}$ of twenty point
injectively distributed on the ovals (i.e no two lands on the same
oval), we have a unique quintic through them and therefore also a
dual group of 20 points cut by the same quintic, and which is
still an injective distribution among the same twenty ovals as
those where $g_{20}$ lives. So we obtain an involution on $G_{20}$
the variety of all distributions of 20 points preserving its
(toric) components, and whose operation leaves the interpolating
map invariable. Actually if one imagine a 20-tuple and the
corresponding (interpolating) quintic each points has a unique
companion of the same oval of the $C_8$, and we can flip each of
them independently to gain $2^{20}$ many 20-tuple inducing the
same quintic.

It is fairly puzzling to imagine that as the 20 points moves along
their respective ovals the corresponding quintic can never cross
the 2 barriers. So the variability of the $C_5$ is much
%%%% entraved
hindered by the 2 barrier ovals. So for instance if we imagine
Fiedler's prohibited curve $\frac{1}{1}\frac{2}{1}\frac{16}{1}$
and we choose as marked 20 ovals all but the two containing an
even number of ovals we get qualitatively the following picture.
There we switched to dashed the 2 ovals which we decreed as being
barriers (i.e. just those ovals where we choose no marking). On
the remaining 20 ovals we choose one point on each and trace then
(in red) the unique quintic interpolating them. (By the lemma it
is unique, otherwise we can violate B\'ezout, et ``\c{c}a baise
tout'' as we say in French). Now whatever the position of the 20
points the corresponding quintic will never sweep across the two
barriers into and outside of which it stays confined perpetually.
This seems a rather strong property, but alas we do not know if
one can derive from this the  Fiedler-Viro regular obstructions
(oddity law) and perhaps the sporadic avatars too, along a purely
elementary B\'ezout line of thoughts.

\begin{figure}[h]\Figskip
%\vskip-1.2cm\penalty0
%\centering
\hskip-2.7cm\penalty0
\epsfig{figure=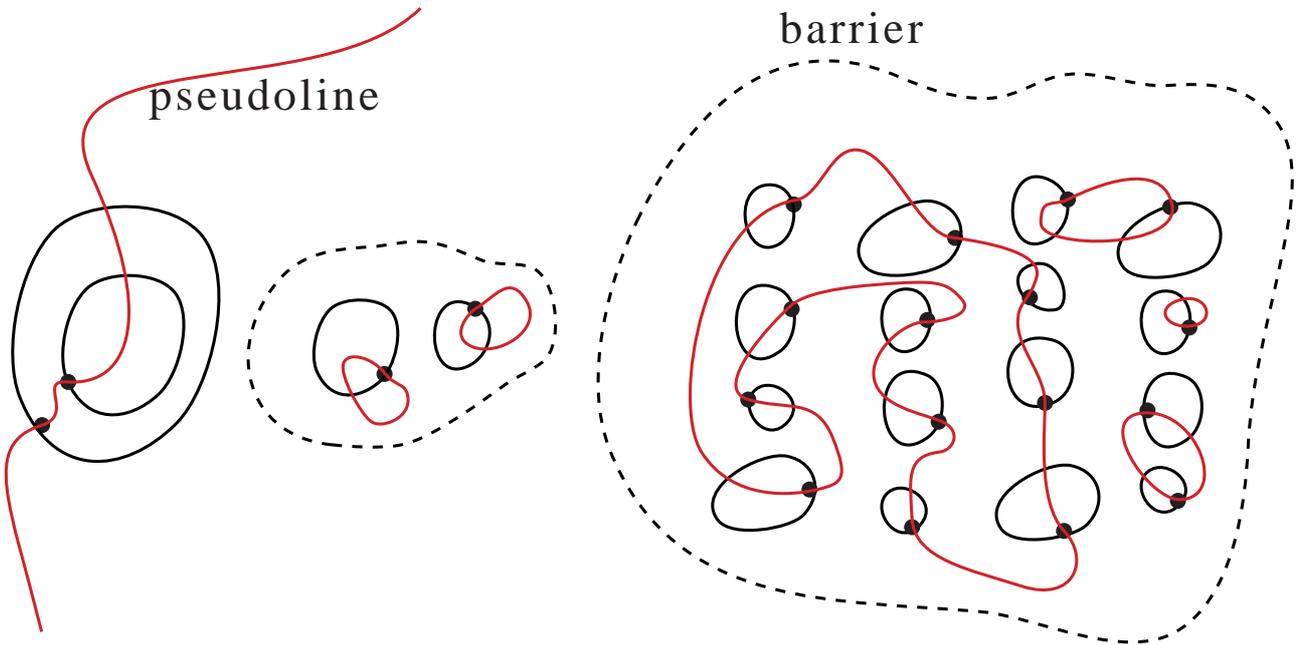,width=172mm} \captionskipAG
  \caption{\label{ViroDEGREE8_Fiedler:fig}%
  Confinements by the barriers (on the Fiedler
  anti-curve $\frac{1}{1}\frac{2}{1}\frac{16}{1}$)}
\figskip
\end{figure}

\subsection{More total reality: new higher superconductivity cases }

[21.08.13] Albeit unfinished and promising, we leave now the
quintics to move to interpolating sextic. Here we have the basic
phenomenon of total reality \`a la Riemann made synthetical as in
Gabard 2013B \cite{Gabard_2013B-Riemann's-flirt}. We recall
briefly how it works on the case at hand of octics. We look in
codegree 2, here to sextics with $\binom{6+2}{2}=4.7=28$ free
coefficients so that we may impose 26 basepoints to a pencil. On
an $M$-octic we therefore impose 22 basepoints (injectively)
distributed on the 22 ovals and we choose the remaining 4 points
either as a tower of 4 concentrated on one oval or as two little
towers of height 2 dispatched on two different ovals. In both case
by the M\"obius-von Staudt principle of intersection we have one
bonus intersection gained by continuity on each oval bringing the
total number to $26+22=48$, which is miraculously equal to $6.8$,
whence the total reality of the considered pencils.

It may perhaps be argued that total reality can merely reassess
Harnack's bound (think e.g. with the total reality of an
$M$-quartic via a pencil of conics, also of co-degree 2, i.e. 2
units less than the curve under inspection). If so, maybe this
total reality is not a serious weapon toward inspecting the
distribution of ovals of curve. Yet we believe that via the pencil
one should be able to draw by the dextrogyration principle some
valuable information upon complex orientations, and this should in
turn make possible further advances on Hilbert's problem.

Next why stopping at degree 6? In higher degrees one can imagine
again imposing (more) basepoints and eventually several layers
(couche in French) of basepoints imagines as eggs ranged in the
ovals conceived as pigeonholes. Perhaps provided we take care
imposing an odd number of them on each oval, M\"obius-von Staudt
will still create for us one boni intersection on each oval and it
remains merely to count at what happens, i.e. the net profit of
hanseatic capitalism.

So assume auxiliary degree $k=7$ to study (absolute) degree $m=8$.
We have then $\binom{7+2}{2}-2=9.4-2=34$ basepoints assignable. We
distribute them on the 22 ovals, while dispatching the 12
remaining ones as 6 pairs either horizontally or as vertical
towers. All combinatorial possibilities are accepted a priori. By
evenness of the ovals, we gain one more bonus intersection on each
oval and thus the total number of intersection is $34+22=56$ which
is again equal to $7.8$, and total reality is obtained anew. It is
clear that this miracle must reproduce in (all) higher degrees
$m\ge 8$.

Besides, for $k=8$ the miracle probably holds as well. We have
then $\binom{8+2}{2}-2=5.9-2=43$ basepoints assignable. We
distribute them on the 22 ovals, while dispatching the 21
remaining ones as 10 pairs interpreted as towers, but there is one
extra point left alone.  By evenness of the ovals, we gain one
more bonus intersection on each oval but one and thus the total
number of intersection is $43+21=64$ which is again equal to
$8.8$, and total reality is granted anew.

All this looks quite formidable but it remains to inspect if this
can be employed as a tool to investigate distributions of ovals.

Let us consider finally $k=9$. We have then
$\binom{9+2}{2}-2=11.5-2=53$ basepoints assignable. We distribute
them on the 22 ovals, while dispatching the 31 remaining ones as
15 pairs interpreted as towers, but there is one extra point left
alone.  By evenness of the ovals, we gain one more bonus
intersection on each oval but one, and thus the total number of
intersection is $53+21=74$ which is again equal to $9.8$, and
total reality is granted anew.

So we have (modulo a trivial arithmetical check) an infinity of
ways to exhibit total reality of $M$-curves, which are perhaps
relevant to the problem of the distribution of ovals.

For concreteness, we must concentrate on degree $m=8$. The
challenge would be to recover the Fiedler, Viro, Shustin and
Orevkov obstructions (assuming them to be all correct) while also
possibly discovering new obstructions on the six bosons not yet
realized. A priori the game can be dangerous as even the most
basic looking obstruction of Fiedler-Viro could be completely
erroneous. For instance remember from our earlier composition
table that patchwork with extended patches could easily produce
all the schemes prohibited by the Fiedler-Viro oddity law.
Notwithstanding, let us hope that this law is true and then the
question becomes: how to prove it in the most elementary way and
ideally in such a fashion that the new bosonic obstructions (yet
unknown) appear as likewise trivial consequences of B\'ezout. This
we call the principle of the Grande Nation, i.e. French
post-revolutionaries annexing all Prussia and Russia with a single
pseudo-hero, Napol\'eon.

As we see there is many experiments that can be imagined by
studying varied auxiliary curves of possibly very high degree. Of
course if we look say at (absolute) quartics (i.e. $m=4$) then
auxiliary curves of degree 1 gives the classical obstructions (no
nesting for $M$-quartics), and those of degree 2 produces
Harnack's bound. So it seems that all information is obtained by
looking at adjoint curve of co-degree $k=m-2$. Whether this is a
general principle is not clear to us but perhaps quite likely even
for $m=8$. Maybe the situation is just opposite and one can infer
information from higher order $k>m-2$ curves. We shall loosely
refer to them as cases of superconductivity.

Alas, apart from the just observed extension of total reality to
all higher degrees we have not yet a single concrete manifestation
of the principle that superconductivity should afford new
information. Yet , this seems quite likely.

Basically, all those superconductions incarnates total reality
hence implies Harnack's bound, but perhaps with increasing
energetic levels so has to contain additional information upon the
distribution of ovals themselves. This basic idea looks plausible
yet needs to be substantiate with more tangible evidence.

Concretely we look again at $k=6$, then we have 4 extra
basepoints, and it is not clear how to choose them. We suppose
given a trinested $M$-curve with an even number of ovals in one
nest. Then there is actually  two such even nests and one which is
even (compare the pyramid Fig.\,\ref{SIMPLIFIED-TABLE_gurus:fig},
which is merely a combinatorial traduction of Gudkov periodicity).
Maybe we should distribute the 2 extra pairs of basepoints on the
two even nests.

For $k=7$, we had $\binom{7+2}{2}-2=9.4-2=34$ basepoints
assignable, and distributed them on the 22 ovals, while
dispatching the 12 remaining ones as 6 pairs either horizontally
or as vertical towers. All combinatorial possibilities are
accepted a priori. By evenness of the ovals, we gain one more
bonus intersection on each oval and thus the total number of
intersections is $34+22=56$ which is again equal to $7.8$, and
total reality is obtained anew. Maybe here we can infer Viro's
sporadic obstructions for judicious choices of the 6 extra pairs.

Next we have $k=8$, where we had $\binom{9+2}{2}-2=11.5-2=53$
basepoints assignable. We distribute them on the 22 ovals, while
dispatching the 31 remaining ones as 15 pairs interpreted as
towers, but there is one extra point left alone.  By evenness of
the ovals, we gain one more bonus intersection on each oval but
one, and thus the total number of intersections is $53+21=74$
which is again equal to $9.8$ (NO SORRY this is 72), and total
reality is granted anew. Since this is the maximum possible, we
see that the one oval where we assigned only one extra basepoint
is so-to-speak just a double couche, and on it no new
intersections can be created, because the boni intersections
gained outside of this oval already saturate B\'ezout's
%graciousness.
%%%%generosity
hospitality. Hence it is quite puzzling to see a total pencil
where there is no mobile point circulating on one oval. This seems
to contradict all what we knew about the Riemann-Ahlfors map,
since Riemann, Schottky, Bieberbach, Grunsky, etc. Of course one
trivial explanation could be that the pencil degenerate somehow by
splitting off the ground octic as a subfactor, yet then our pencil
would be a very statical object. Of course the same phenomenon
occurs when $k=9$, then perhaps in a less statical incarnation.

Then with $k=10$ evenness of the excess is restored again.
Precisely, we have then $\binom{10+2}{2}-2=6.11-2=64$ basepoints
assignable. We distribute them on the 22 ovals, while dispatching
the 42 remaining ones as 21 pairs interpreted as towers.  By
evenness of the ovals, we gain one more bonus intersection on each
oval, and thus the total number of intersections is $64+22=86$
which is larger than $10.8$, and total reality is lost. Presumably
due to the excess intersection the decaics (degree 10) have to
split off the ground octic. Yet it is curious that as we have a
triple couche on all ovals but one (in simple couche), it seems
that we have a decent circulation \`a la Riemann-Bieberbach, but
apparently B\'ezout is unhappy. Naively it seems nearly that this
reasoning shows that there in not a single $M$-octic in degree 8,
which is blatantly false (in principle).

At least, it seems that our expectation of the phenomenon of total
reality as admitting  infinite repetition in all higher degrees is
foiled, but seems to appear only at degrees $k=m-2,m-1$ and
perhaps $m, m+1$.

Next with $k=11$ evenness of the excess is still conserved. We
have then $\binom{11+2}{2}-2=13.6-2=76$ basepoints assignable. We
distribute one of them on each of the 22 ovals, while dispatching
the 54 remaining ones as 27 pairs interpreted as towers.  By
evenness of the ovals, we gain one more bonus intersection on each
oval, and thus the total number of intersections is $76+22=98$
which is larger than $11.8$ (as $k=10$ the excess $80<86$ was of
six and now it is larger by 4), and total reality is lost.

Of course in the cases where we get excess intersections (over
B\'ezout), we could change the distributions of basepoints as to
produce less boni intersections.

Next we noted that due to a sordid arithmetical mistake of us
($9.8=74$ instead of 72) already the case $k=8$ presents excess
intersection, and geometrically this should mean that the pencil
must split off the ground octic. Yet by linear algebra the points
are interpolated by a pencil (at least or some linear series  of
higher dimension), yet this sounds quite paradoxical because as no
residual curve is available, we really seems to face a paradox of
the sort quite common to algebraic geometers (compare
Enriques-Chisini's discussion of MacLaurin, or so).

So it seems important to settle the paradox, as there is some hope
to derive from it a tension potentially valuable to Hilbert's
16th. Alternatively we can ignore this and hope that the real
information on Hilbert's problem is
%condensed
stocked in the cases $k=6$ and 7 where total reality works well
without overheating B\'ezout. Of course the genuine information is
perhaps also stocked by using lower order curves with $k=3,4,5$ as
we tried unsuccessfully to sketch some few days (or pages) ago.

For $k=7$, we had $\binom{7+2}{2}-2=9.4-2=34$ basepoints
assignable, and distributed them on the 22 ovals, while
dispatching the 12 remaining ones as 6 pairs either horizontally
or as vertical towers. All combinatorial possibilities are
accepted a priori. By evenness of the ovals, we gain one more
bonus intersection on each oval and thus the total number of
intersections is $34+22=56$ which is again equal to $7.8$, and
total reality is obtained anew. Maybe here we can infer Viro's
sporadic obstructions for judicious choices of the 6 extra pairs.
As we reach B\'ezout's bound no more intersection exist than those
imposed and the single one created by continuity on each oval.
This seems to be a strong constraint, and maybe there is a clever
way to deduce something out of it. Where to choose the 6 extra
pairs for a given $M$-scheme (which is prohibited)?  Perhaps the
simplest prohibition (at least historically) seems to be Fiedler's
prohibition of 4 purely trinested $M$-schemes, starting say with
$\frac{1}{1}\frac{2}{1}\frac{16}{1}$ which as an anecdote contains
also Hilbert's beloved number 16 in the 3rd numerator.

So imagine Fiedler's (anti) curve as on
Fig.\,\ref{ViroDEGREE8_Fiedler2:fig} below where we already
imposed a basepoint on each ovals but we are still free to place 6
pairs of basepoint. Naively it could seem that the best choice is
to impose them on the deep ovals so that the septics being
assigned to visit those deep ovals more frequently there will be a
higher probability that a certain curve of the pencil will cross
the separating nests too often. (Recall that each curve of the
pencil can cross only once each oval outside of the imposed
basepoints.) So if we can force more intersections B\'ezout is
corrupted and the scheme prohibited. However on doing a
qualitative picture we see that even with this deep assignment
there is no difficulty to trace by free hand a curve visiting all
basepoints while crossing each oval at most once. Of course as we
in reality a pencil we should trace not just a single curve but
worry about the whole induced foliation (with mild singularities
at the base point or at the critical curves). The point is that
the pseudo line will intersect itself and so contributes to the
creation of unassigned basepoint. But of course it can also be the
case that the pseudo-line visits assigned basepoints.

\begin{figure}[h]\Figskip
%\vskip-1.2cm\penalty0
%\centering
\hskip-2.7cm\penalty0
\epsfig{figure=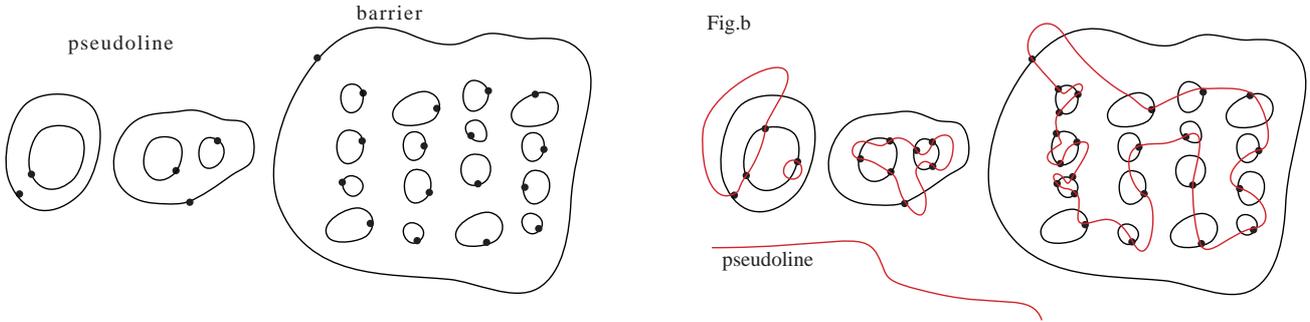,width=172mm}
\captionskipAG
  \caption{\label{ViroDEGREE8_Fiedler2:fig}%
  Confinements by the barriers (on the Fiedler
  anti-curve $\frac{1}{1}\frac{2}{1}\frac{16}{1}$)}
\figskip
\end{figure}

So the problem looks fairly difficult, yet perhaps trivial once
one has the right idea. Actually, unclever peoples just needs a
lot of time to test all ideas (=possibilities) until finding the
true argument (reminds Poincar\'e's story about {\it les
combinaisons stables\/}).

At any rate, it is clear that  we are in realm quite fascinating
(Hilbert, Rohlin, Fiedler, Viro) becoming even more so once the
connection with conformal mapping is noticed (Riemann, Bieberbach,
Grunsky, Ahlfors, etc.), yet it is still very hard to understand
properly the Fiedler-Viro obstruction in the most synthetic way.
Of course their proofs might be the correct one yet we may expect
a simpler argument leading perhaps to new insights (i.e. new
restrictions that as far as we can judge where not yet derived by
the methods of Fiedler, Viro, Shustin, Orevkov) which are
methodologically a bit disparate and loosely unified. One could
dream that very simple B\'ezout style arguments do obstructs some
of the octic schemes. The difficulty is that there are several
fronts where to attack the problem $k=3,4,5,6,7$ and perhaps
higher so that once energy is much dissipated by the variety of
situation to analyze. Further each situation request long hours of
concentrations to get mentally familiar with and looks like a big
Eiger Nordwand face hard-to-climb upon. So a human intelligence is
quickly desperate by the duty,  and the sole consolation is that
the proof must be trivially beautiful once found.

So one can hope to have a reliable geometric flair of where to
find the argument without wasting to much energy. On the one hand
using say the case of degree $m=4$ as prototype, we see that once
Harnack's bound is known curve of degree 1 (lines) suffice to
settle the isotopy classification of (in particular) $M$-curves.
Those having co-degree $m-k=4-1=3$, we may expect that in degree 8
similar information (isotopic classification) is gained by curves
of degree 5. Another sloppy reason would be the role of the
canonical class which for plane curves of degree $m$ is cut out by
adjoint of degree $m-3$. However when it comes to $m=5$, the
crucial role should be played by conics, yet it seems to be still
played by lines. Even when $m=6$, curves of co-degree 3 i.e.
cubics plays a little role safe in the fundamental Rohlin-Le
Touz\'e phenomenon which explains all of Gudkov's prohibitions in
some ad hoc way.

Maybe this gives the following idea: to get obstruction on
$M$-curves one must look lower at $(M-2)$-curves of the RKM type
(hence universally orthosymmetric) and tabulating upon Rohlin's
maximality conjecture this would kill new schemes. One problem
with this approach is that Rohlin's maximality principle  looks
severely foiled in degree $m=8$ (compare our counter-examples
\ref{RMC:cter-example-via-Viro-1st-curve=BEAVER} via Viro's 1st
curve (beaver) or \ref{RMC:cter-example-via-Shustin} via Shustin's
medusa).

Of course, it would be perhaps of interest to see if our naive
methods (deepest penetration=DEPP, or total reality=TOR) do work
in degree 6 already to get the classical prohibitions of Hilbert,
Rohn, Gudkov. This is already a (ingrate) nontrivial exercise, yet
perhaps necessary to prepare further progresses, or gain
confidence in the method: notably, to decide the question if it is
better to work with codegree 3 (DEPP) or codegree 2 (TOR).

Through 8 points we can pass a pencil of cubics in particular a
connected cubics. Through 9 points we can pass a cubic. So
imposing 9 points on the 11 ovals of an $M$-sextic we get twice so
many intersection granted, i.e. 18 by M\"obius-von Staudt
principle of evenness, which is the maximum permissible. So if we
impose those 9 points on the 10 deep ovals of the $M$-sextic we
see that the nonempty oval of the $C_6$ will act as a barrier upon
the situation of the cubic (since the intersection $C_3\cap C_6$
is already saturated to 18), and will actually fall it apart

Fig.\,\ref{ViroDEGREE8_Fiedler3:fig}. One sees in particular that
any such distribution of 9 points impose independent conditions on
cubics. Else if there would be a pencil freedom we could impose
visiting the barrier and B\'ezout is corrupted. We recognize here
a perfect analogy with the situation $(m,k)=(8,5)$ studied
earlier. In particular, we see again the phenomenon of confined
thermo-excitation, namely: whatsoever the position given to the
nine anchors the unique cubic through them avoids the barrier. Now
we can additionally infer that it must be smooth, otherwise it is
connected and so hit the barrier or possess a solitary node but
then cannot fulfill its interpolating duties inside the $C_6$'s
nonempty oval (except perhaps for Harnack's scheme
$9\frac{1}{1}$).

Philosophically, we see that an $M$-curve of degree $m$
constitutes a strong trap for curves of codegree 3, i.e. degree
$m-3$. From this principle one would lie to deduce the classical
prohibitions (Hilbert, Rohn, Gudkov as done without the
simplifying assistance of Arnold-Rohlin using systematically
homology).

Let us assume that we have Rohn's scheme $\frac{10}{1}$. Then we
may impose 9 basepoints inside the nonempty oval and the resulting
cubic has to be trapped in it (because the intersection is already
saturated to 18 by topology), but this is impossible as the
cubical circuit cannot be null-homotopic. This contradiction
supplies a very elementary proof of Rohn's prohibition. (Discovery
of us at [14h24, 21.08.13]). So:

\begin{theorem}
There is a completely trivial proof of Rohn's prohibition using
just M\"obius-von Staudt principle of intersection (you cannot
penetrate in an oval along a recurrent motion without escaping
once of it). This proof is one just given!
\end{theorem}

Yet, Fig.\,c kills our pseudo-proof, since our argument overlooked
the option that it is just the oval of the cubic which visits the
9 basepoints.

\begin{figure}[h]\Figskip
%\vskip-1.2cm\penalty0
%\centering
\hskip-2.7cm\penalty0
\epsfig{figure=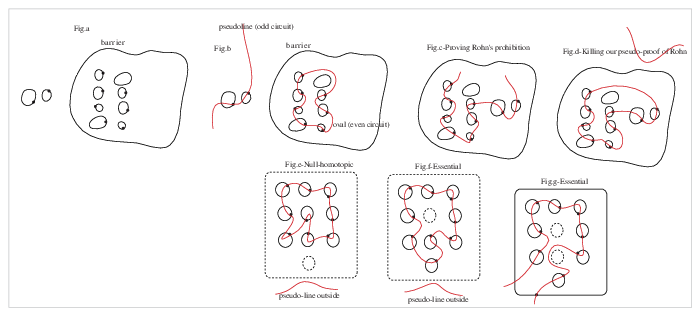,width=172mm}
\captionskipAG
  \caption{\label{ViroDEGREE8_Fiedler3:fig}%
  XXX}
\figskip
\end{figure}

[22.08.13] Nonetheless, we see that if we impose 9 points on the
ovals (injectively), then the cubic has by M\"obius-von Staudt
already 18 intersection with the $C_6$ and therefore we deduce
that the two remaining ovals are not intercepted. In case the
sextic has a non-empty oval, and we distribute the 9 points
outside of it, then the interpolating cubic cannot meet the
nonempty oval nor the other oval left unmarked. Thus provided
points are marked inside of the {\it nonempty oval\/} (barrier for
short), the latter acts as a separator splitting the cubic into 2
pieces. If the cubic would be connected, then it would be trapped
inside the barrier,  forcing it to be null-homotopic, which is
impossible. Incidentally, we see also that the cubic must be
smooth and forced to have two components. Moreover we see that the
9 points imposes independent conditions on the cubic, since there
were 2 of them we could impose to visit a point on the barrier
(and B\'ezout is corrupted).

So we really get  a stringent context, where we can drag all 9
points on their respective ovals while seeing always a unique
corresponding (split) cubic whose oval stays confined inside the
barrier. Alas, we do not know how to draw a contradiction under
the standard hypotheses (i.e. outside  of the Harnack, Hilbert,
and Gudkov curves).

Several ideas are as follows.

First, one could expect the mapping assigning to each
9-distribution the unique interpolating cubic to be \'etale or
just open, in which case the image would be a clopen hence full.
(Remind, indeed, that the source is a compact finite union of
tori, whence closed-ness of the image.) Yet this violates the
issue that the interpolating cubics are always split, i.e. smooth
with 2 components (hence confined in one chamber past the
discriminant). Hence it seems that there is no chance for the
canonical map being open.

Next, one can observe that the cubic's oval will belong to a
unique homotopy class past the 2 barriers ovals which are
unmarked, just because we can deform along the 9-torus of all
distributions. To be more concrete, if we assume Rohn's type
$\frac{10}{1}$ and the 9 marked ovals to be inside then the
resulting cubics cannot meet the nonempty oval nor the little
empty oval left unmarked. Thus the cubic's oval is either
null-homotopic or winding once around the hole formed by the
little unmarked oval. Probably both cases are possible upon
varying the marking (cf. Figs.\,e, f). Further one can of course
also mark the nonempty oval (Fig.\,g) but then it cannot be
anymore guaranteed that the interpolating cubic is split.

It is clear that the strategy looks hard to complete.

Again, in degree 6 most obstructions can be inferred from the
Rohlin-Le Touz\'e total reality of the two $(M-2)$-schemes which
are RKM hence universally orthosymmetric. One can wonder if this
method applies as well in degree 8, as to re-explain the
obstructions of Fiedler, Viro, Shustin, Orevkov lying beyond
Gudkov periodicity. A loose counter indication from the scratch is
that in degree 6, the Rohlin-Le Touz\'e phenomenon explains only
(as it should since there are no more obstructions) Gudkov
periodicity. Still one could expect that a total reality at level
$(M-2)$ could kill $M$-schemes. So $(M-2)$-total reality acts as a
cosmic censorship over $M$-schemes, exactly as the G\"urtelkurve
of degree 4 explains the classification of $M$-quartics as being
necessarily unnested.

Now looking at the main pyramid
(Fig.\,\ref{Degree8-(M-i)-curve-TABLE:fig} or its enlargement
Fig.\,\ref{Degree8-(M-i)-curve-TABLE_I:fig}) we can try to
speculate of where to locate totally real $(M-2)$-schemes so as to
induce the known obstructions.

For instance to explain Orevkov's obstructions, we could imagine
that the scheme $T_1:=\frac{3}{1}\frac{15}{1}$ is totally real
(under a pencil of quintics, i.e. co-degree 3 like by Rohlin-Le
Touz\'e). As good news our map remembers us that this scheme
easily exist via Viro's simplest method. Albeit not RKM, we posit
that this scheme is totally real. Quintics have 19 basepoint
assignable for a pencil, and so 38 intersections are granted and
we need just a miracle of 2 to reach total reality at $40=5.8$.
Now what are the enlargements of this scheme. First there is
$1\frac{3}{1}\frac{15}{1}$ which is dominated by Viro's
anti-scheme $\frac{1}{1}\frac{3}{1}\frac{15}{1}$, and this is
another good news as our supposition of total reality would kill
this scheme too. Next on the right wing there is the scheme
$\frac{3}{1}\frac{16}{1}$ that would be killed as well, and this
does not contradict experimental data available to us. As another
enlargement of $T_1$ we have $\frac{4}{1}\frac{15}{1}$ and
$1\frac{4}{1}\frac{15}{1}$ and this would explain the
disintegration of the boson $b4$ too.

Likewise Orevkov's 2nd obstruction on
$b6:=1\frac{6}{1}\frac{13}{1}$ could be explained by a total
reality of the $(M-2)$-scheme lying below it, namely
$\frac{6}{1}\frac{12}{1}=:T_2$. This would kill the boson $b7$.
Alas as yet our table does not report a realization of $T_2$.
Additionally, our scheme $T_2$ is dominated by Fiedler's scheme
$\frac{1}{1}\frac{6}{1}\frac{12}{1}$ whose prohibition would be
derived anew.

Likewise we can imagine a total reality killing the boson
$b1:=1\frac{1}{1}\frac{18}{1}$, hence concerning the scheme
$\frac{1}{1}\frac{17}{1}=:T_3$. Alas, this would kill also Viro's
scheme $1\frac{2}{1}\frac{17}{1}$ which exist however, and also
Shustin's scheme $\frac{1}{1}\frac{1}{1}\frac{17}{1}$. Hence total
reality at $T_3$ is unlikely.

Finally, we can imagine a total reality killing the boson
$b9:=1\frac{9}{1}\frac{10}{1}$, hence concerning the scheme
$\frac{9}{1}\frac{9}{1}=:T_4$. A first good news is that $T_4$
exists by a simple Viro method. In view of its extremal position
our scheme as no enlargement  in the first pyramid safe those
readily visualized above it. However an enlargement occur in the
2nd pyramid where we find $\frac{1}{1}\frac{9}{1}\frac{9}{1}$,
which is prohibited by a Viro sporadic obstruction. In conclusion
total reality at $T_4$ is fairly likely, and compatible with
Russian knowledge/folklore.

In contrast, we expect now total reality below the boson
$b7=1\frac{7}{1}\frac{12}{1}$, i.e. at
$T_5:=\frac{7}{1}\frac{11}{1}$ as this would kill Viro's scheme
$b8=1\frac{8}{1}\frac{11}{1}$.

\subsection{Total reality as a cosmic censorship}

[22.08.13] Can we continue this game (of total reality as a cosmic
censorship) as to explain also the other Fiedler, Viro
obstructions. Yes, we can; it seems at least worth trying to
elucidate  this.

First, below Fiedler's obstruction of
$\frac{1}{1}\frac{2}{1}\frac{16}{1}$ we find
$\frac{1}{1}\frac{2}{1}\frac{14}{1}$ which we posit totally real.
The impact would be to kill $\frac{1}{1}\frac{3}{1}\frac{14}{1}$
and $\frac{1}{1}\frac{3}{1}\frac{15}{1}$, in accordance with
prohibitions by Shustin and Viro respectively. A third enlargement
involves $\frac{2}{1}\frac{2}{1}\frac{14}{1}$ and
$\frac{2}{1}\frac{2}{1}\frac{15}{1}$ also prohibited by Shustin
and Viro respectively.

Next, below Viro's obstruction of
$\frac{1}{1}\frac{3}{1}\frac{15}{1}$ we find
$\frac{1}{1}\frac{3}{1}\frac{13}{1}$ which kills the way right
above it (Shustin, Viro), also that one in front of it (Shustin,
Fiedler), and that in the shifted layer (Shustin, Viro).

Next, below Fiedler's obstruction of
$\frac{1}{1}\frac{4}{1}\frac{14}{1}$ we find
$\frac{1}{1}\frac{4}{1}\frac{12}{1}$ which kills the way right
above it (Shustin, Fiedler), also that one in front of it
(Shustin, Viro sporadic), and that in the shifted layer (Shustin,
Viro regular).

Next, below Viro's obstruction of
$\frac{1}{1}\frac{5}{1}\frac{13}{1}$ we find
$\frac{1}{1}\frac{5}{1}\frac{11}{1}$ which kills the way right
above it (Shustin, Viro sporadic), also that one in front of it
but this time conflicting with a construction claimed by
Polotovskii namely $\frac{1}{1}\frac{6}{1}\frac{11}{1}$. So either
Polotovskii is wrong or so is our censorship principle. Eventually
our censorship could still be true, but our scheme
($\frac{1}{1}\frac{5}{1}\frac{11}{1}$) would lack total reality.
This is a possible scenario, yet would be annoying because then
total reality would not explain all prohibitions. However
$\frac{1}{1}\frac{5}{1}\frac{13}{1}$ could be prohibited by
$\frac{1}{1}\frac{4}{1}\frac{12}{1}$ as we saw, and this restores
the hope to explain everything via total reality. On the shifted
layer our scheme kills dully schemes prohibited by Shustin and
Viro.

Next, below Fiedler's obstruction of
$\frac{1}{1}\frac{6}{1}\frac{12}{1}$ we find
$\frac{1}{1}\frac{6}{1}\frac{10}{1}$ which overkills the way right
above it (Polotovskii, Fiedler), and also makes conflicting
damages in front of it (Polotovskii, Shustin), but in the shifted
layer we recover known obstructions (Shustin, Viro). Yet, in
summary it seems that we cannot expect total reality for this
$(M-2)$-scheme.

At this stage no more comments should be necessary, and it
suffices to mark on the table
(Fig.\,\ref{Degree8-(M-i)-curve-TABLE:fig}) by TOR schemes
susceptible of total reality and by NIET those for for which there
there is ``No Instinctive Evidence for Total reality''. It should
be noted that contrary to what we said there is no conflict
between our principle of censorship and Polotovskii's
constructions.

Next we move to the 2nd layer of the 2nd pyramid. First we meet
the scheme $12\frac{2}{1}\frac{2}{1}\frac{3}{1}$. If censorship is
a universal reason for $M$-prohibitions, one would naively posit
total reality for the scheme below it
$11\frac{2}{1}\frac{2}{1}\frac{2}{1}$. Actually there is other
$(M-2)$-schemes in the first layer dominated by
$12\frac{2}{1}\frac{2}{1}\frac{3}{1}$, yet positing their total
reality (TOR) would corrupt a construction by Polotovskii. Still
we can in the 1st layer posit TOR for
$10\frac{1}{1}\frac{2}{1}\frac{4}{1}$, and several other schemes
marked by TOR on the main table (all this being consistent with
Polotovskii and re-explaining Viro's law). Now back to the2nd
layer, our TOR-postulation on
$11\frac{2}{1}\frac{2}{1}\frac{2}{1}$ would via censorship
prohibit $12\frac{2}{1}\frac{2}{1}\frac{2}{1}$ (a question left in
suspense in Shustin 90/91
\cite{Shustin_1990/91-New-restrictions}).

Next we have $8\frac{2}{1}\frac{2}{1}\frac{7}{1}$ which as before
cannot be prohibited from the 1st layer (without conflicting with
Polotovskii), and so we put a TOR-tag right below it. Note at this
stage that in this 2nd layer we have $(M-2)$-schemes prohibited by
Viro, and to rules them out via our method we should posit certain
total realities at the $(M-4)$-level running thereby out of our
tabulation. Maybe we should content to explain only $M$ and
$(M-1)$-prohibitions.

The next interesting case is $\frac{2}{1}\frac{2}{1}\frac{15}{1}$
which as we saw can be ruled out by a TOR-prescription below the
appropriate Fiedler's scheme
($\frac{1}{1}\frac{2}{1}\frac{14}{1}$), hence there is no need to
impose a TOR-tag on $\frac{2}{1}\frac{2}{1}\frac{13}{1}$. Beside
an obvious principle of economy the net bonus is that we get so a
very regular distribution of TOR's on the first layer (on the
first row at least as examined up to now). It should be remarked
however that other more frequent distributions of TOR's could
explain the same prohibitions: so we could instead of TORing the
central item of each monticulus we could TOR the two lateral items
at the basis of the triangle. In that case Shustin's uncertain
scheme $12\frac{2}{1}\frac{2}{1}\frac{2}{1}$ would not be killed
except if we TORize the left basis
$12\frac{2}{1}\frac{2}{1}\frac{1}{1}$, which is however a ghost
copy of a scheme in the 1st layer below a Polotovskii trademark
(construction). Hence Shustin's incertitude remains very vivid,
and not easy to settle even after acceptance of our censorship
principle.

Next we arrive at $8\frac{2}{1}\frac{3}{1}\frac{6}{1}$. This
prohibition cannot be explained by TOR in the 1st layer, and so we
are forced to put a TOR at $6\frac{2}{1}\frac{3}{1}\frac{6}{1}$.
Alas, this breaks our central positioning of TOR's in first row
(of the 2nd layer). One checks quickly that this causes no undue
damage in the 3rd layer. (This is because the right basis of a
monticulus=triangle has no superior in the upper layer, as a
consequence of GKK-periodicity=Gudkov, Krakhnov, Kharlamov).

Next we have $4\frac{2}{1}\frac{3}{1}\frac{10}{1}$. One can of
course impose a TOR on  $3\frac{2}{1}\frac{3}{1}\frac{9}{1}$ to
explain the 3 prohibitions right above, and in the 3rd layer this
implies 2 additional prohibitions namely
$3\frac{3}{1}\frac{3}{1}\frac{9}{1}$ (Shustin) and
$4\frac{3}{1}\frac{3}{1}\frac{9}{1}$ (Viro's most sporadic).
Alternatively one can fix TOR's at the 2 corners of the triangle
($4\frac{2}{1}\frac{3}{1}\frac{8}{1}$ and
$2\frac{2}{1}\frac{3}{1}\frac{10}{1}$) to get the same censorship
on the 2nd layer, yet now causing different damages on the 3rd
layer namely killing rather the left versant of the monticulus.
Thi is to say that $4\frac{3}{1}\frac{3}{1}\frac{8}{1}$ and
$4\frac{3}{1}\frac{3}{1}\frac{9}{1}$ are now killed.

Next we examine $8\frac{2}{1}\frac{4}{1}\frac{5}{1}$. There is no
way to eliminate it by a TOR in the 1st layer (due to a
construction of Polotovskii, that we checked via Viro), and so we
are forced to place the TOR at the natural location
$7\frac{2}{1}\frac{4}{1}\frac{4}{1}$. This kill besides 2 schemes
in the 3rd layer (Shustin and Viro).

Next we examine $4\frac{2}{1}\frac{4}{1}\frac{9}{1}$. Again there
is no way to eliminate it by a TOR in the 1st layer (due to a
construction of Polotovskii, that we did not checked but probably
just as a consequence of the fact that we initially missed all
combinations of Viro's method). So we are forced to place the TOR
at the natural location $3\frac{2}{1}\frac{4}{1}\frac{8}{1}$. This
kill besides 2 schemes in the 3rd layer (Shustin and Viro).

The next case of interest arises with
$\frac{2}{1}\frac{5}{1}\frac{12}{1}$ where there is not anymore a
reduction to the 1st layer, and thus it seems now necessary to
impose a TOR on the 2nd layer at the natural place, namely
$\frac{2}{1}\frac{5}{1}\frac{10}{1}$.

Next we can move to the 3rd layer. First, albeit not completely
necessary (depending on what we did on the 2nd layer) we are
invited to put a TOR at $3\frac{3}{1}\frac{3}{1}\frac{8}{1}$.

Next we inclined to put a TOR on
$\frac{3}{1}\frac{3}{1}\frac{11}{1}$ albeit this could be
dispensed if we had introduced one at
$\frac{2}{1}\frac{3}{1}\frac{12}{1}$. Our idea of looking at
totally real $(M-2)$-schemes as a trick to find the true reason
behind the seemingly chaotic distribution of prohibitions is only
half efficient. Yet our Rohlin-style philosophy that total reality
should regulate the distribution of ovals looks to us an extremely
appealing law  that the divine nature is probably following. Of
course it can be that all our discussion is biased by Viro's
theory in case he used anomalous patching parameters. Of course
the censorship principle is merely a strong form of total reality
permitting one to sweep out the curve by a totally real pencil,
hence forbidding the presence of any additional ovals. We think it
is easy to justify theoretically. Philosophically, it seems that
total reality is fairly ubiquitous and therefore explaining the
many prohibitions of Fiedler, Viro et cie (Shustin, Orevkov). So
we imagined that prohibition are bad, but in reality they are the
reverberation of the goodness of total reality. Furthermore the
presence of many prohibition in low degree will permit (via the
satellite principle) the presence of more schemes in higher
degrees multiple of 8 (since an $M$-scheme is totally real and so
will kill all extensions of its satellites).

Next a TOR at $6\frac{3}{1}\frac{4}{1}\frac{4}{1}$ is not even
requested. For $4\frac{3}{1}\frac{4}{1}\frac{8}{1}$ and its
$(M-1)$-companion $3\frac{3}{1}\frac{4}{1}\frac{8}{1}$ we do not
need to place a TOR at $2\frac{3}{1}\frac{4}{1}\frac{8}{1}$ since
$3\frac{2}{1}\frac{4}{1}\frac{8}{1}$ does already the killing-job.

For $\frac{3}{1}\frac{4}{1}\frac{12}{1}$ we need a TOR, in case we
did not placed one in the 2nd layer. So the discussio can be
continued in a quite tricky way. When reaching
$\frac{3}{1}\frac{8}{1}\frac{8}{1}$ it seems necessary to
introduce a TOR at $\frac{3}{1}\frac{7}{1}\frac{7}{1}$, yet this
not even needed as we ``TORed''
$\frac{2}{1}\frac{7}{1}\frac{8}{1}$. We can remark that it must be
nearly possible to avoid any TOR in the 3rd layer if we distribute
them suitably on the 2nd layer. This request to be better
analyzed.

Despite tired and poorly organized, our troupes can now enchain
with an attack of the 4th layer. A construction of Polotovskii in
the 3rd layer forces us to tag (by a TOR) the scheme
$3\frac{4}{1}\frac{4}{1}\frac{6}{1}$.

Actually we remark that in the 3rd layer
$\frac{3}{1}\frac{7}{1}\frac{7}{1}$ cannot be TOR as it is
dominated by schemes of Polotovskii and Shustin in the 4th resp.
5th layer. Our thesis is really that all prohibitions of Hilbert's
16th (in degree 8 those being due to Fiedler, Viro, Shustin,
Orevkov) can in reality be reduced and uniformized through a
Rohlin-Le Touz\'e series of phenomena of total reality. This is an
unifying theme as to comprehend whole of them in as single soup
incarnating the telluric plasm behind each little volcanic
irruptions (seemingly completely random), yet governed reality by
a deep flow at the $(M-2)$-level and not just by the little summit
of volcanos at the $M$-level. So we a perfect metaphor with the
famous geological story about the foss Marianes where much
terrestrial crust goes absorbed by tectonic translation, resulting
thereby in the chain of volcanoes.

It is at this stage that we discovered the curvy arrows of
domination moving two rows upwards as depicted on
Fig.\,\ref{Degree8-(M-i)-curve-TABLE:fig}. With those even less
TORs are required to explain all prohibitions. For instance
$\frac{2}{1}\frac{7}{1}\frac{8}{1}$ does not seem to request
anymore a TOR, since all its prohibited entourage can be reduced
to other TORs. Of course it is then requested that
$\frac{3}{1}\frac{6}{1}\frac{8}{1}$ is a TOR.

So it seems that to approach the problem more systematically we
need first to look at construction (e.g. Shustin), accept them as
legal, and then look at all arrows (magma motions of the telluric
flow) and put NIETs whenever we are the endpoint of an arrow
issuing from a constructible $M$-scheme. We see then that
Hilbert's problem is an organic whole not concerning isolated 6
schemes but there is global coherence and potentially everything
Viro included must revised from the very beginning. Then by
wondering where those curvy arrows lands when approaching the top
we discovered the green curvy-arrow, of which there is a menagerie
not all traced on our diagram.

At this stage it seems wise to rationalize the diagrammatic by new
figure with less curvy arrows. With the new algoritm in mind
namely by flowing from the constructible  $M$-curves we can
sharpen information for instance via the curvy arrow from Viro's
$4\frac{3}{1}\frac{5}{1}\frac{7}{1}$ we move down to
$4\frac{3}{1}\frac{3}{1}\frac{7}{1}$ which is therefore NIET, i.e.
no total reality (or as we say German ``eine Niete''). Alas we
must chancge NIET to NOT to save room (abridging still no total
reality).

Next we realized (aided by this better diagrammatic) that the
scheme 7.2.3.6 (abridged notation for
$7\frac{2}{1}\frac{3}{1}\frac{6}{1}$) via 7.2.2.6 without that it
it necessary to impose a TOR on 6.2.3.6. As a consequence we can
impose a more regular distribution of TOR's on the 2nd layer. It
seemed also at this stage advisable to change the whole
diagrammatic by interpreting the 1st pyramid as the ground-zero
layer (0th layer). With this better diagrammatic it also apparent
that Orevkov's obstructions when interpreted via an underlying
total reality one sporadic obstruction of Viro (on 1.3.15) and one
obstruction by Fiedler (namely 1.6.12).

\begin{figure}[h]\Figskip
%\vskip-1.2cm\penalty0
%\centering
\hskip-2.7cm\penalty0
\epsfig{figure=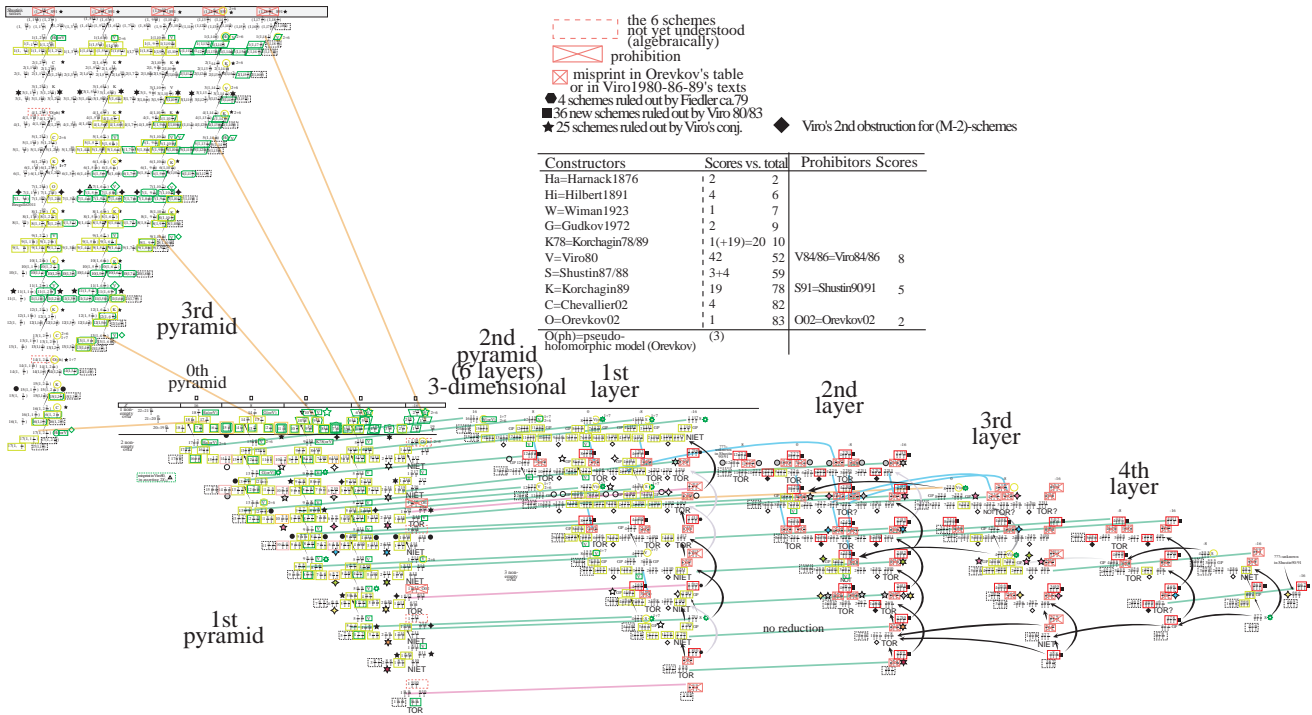,width=172mm}
\captionskipAG
  \caption{\label{Degree8-(M-i)-curve-TABLEBIS:fig}%
  Yet another view of the pyramid}
\figskip
\end{figure}

\subsection{Some conjectures \`a la Rohlin-Le Touz\'e
implying a completion of Hilbert-Viro's 16th problem}

[22.08.13] At this stage we started to get some understanding in
Orevkov and Viro's sporadic obstruction by gaining an
understanding of their magmatic coherence (despite apparent
randomness when one does does not take care sufficiently about the
architecture/combinatorics of the pyramid). In particular it seems
that there is a natural distribution of TOR at the $(M-2)$-level
explaining via total reality and the allied Rohlinian phenomonon
of censorship (i.e. a scheme totally flashed by a pencil of curves
cannot be augmented without corrupting B\'ezout) all prohibition
in an uniform fashion. In particular extrapolating a bit and
hoping that the distribution of TOR is uniquely determined on the
basis of (already) available constructional knowledge (i.e. Viro's
theory and his many companions), we shall get a complete
resolution of Hilbert's 16th. Of course the solution remains
heuristic unless we are able: first to establish the requested TOR
by synthetical algebraic geometry (i.e. Rohlin-Le Touz\'e type
theorems) and second to establish rigorously the principle of
censorship (that must be easy).

At this stage it seemed advisable to improve the overall
architecture by shifting the layers diagonally. Perhaps the whole
exercise should from the top of the telescope of the pyramid where
Shustinian information is reigning. Actually at the very summit of
the telescope we have Viro's anti-scheme 6.6.7. Looking at what is
below we find 5.6.7 (constructed by Polotovskii via Viro) and
4.6.7 which must therefore be TOR. Notice that 5.5.7 is also
constructed by PV(=Polotovskii-Viro) or can be interpreted loosely
as a double contraction of Shustin's 5.7.7. Hence it seems that we
are really forced to put a TOR on 4.6.7. Yet this has the
disastrous effect of killing 5.6.7 which is constructed by PV. So:

\begin{Scholium}
There is noway to explain all prohibtion of the actual census in a
fashion respecting Rohlin-Gabard's desideratum of total reality
and censorship.
\end{Scholium}

By the way the domination of Shustin's 5.7.7 over 4.6.7 forces the
latter to be NIET. In our opinion this sad issue incarnates a
severe anomaly in the architecture of the pyramid. Did God
constructed such a disgraceful world? Maybe Shustin is again
responsible of the turmoil of the business?? Of course we know
that there is a certain friction between the thesis of Rohlin and
Shustin's discovery,  yet maybe Shustin's curve are too much
free-hand traced and a microscopic mistake went unnoticed through
the decades.

\begin{figure}[h]\Figskip
%\vskip-1.2cm\penalty0
%\centering
\hskip-2.7cm\penalty0
\epsfig{figure=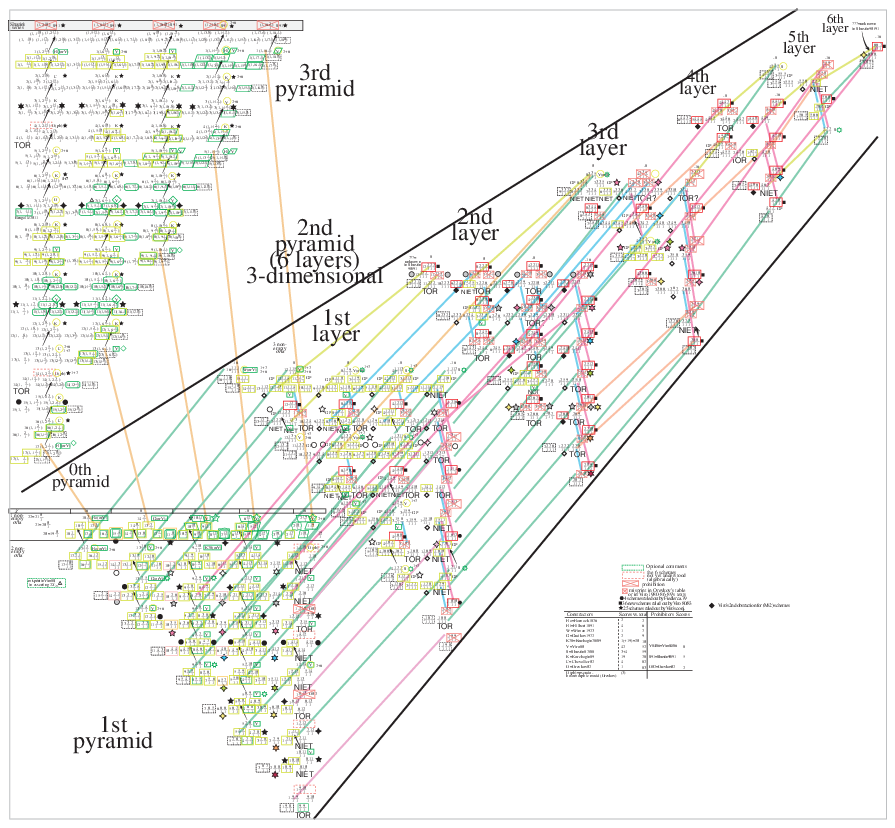,width=172mm}
\captionskipAG
  \caption{\label{Degree8-(M-i)-curve-TABLETRIS:fig}%
  Still another (more comfortable?) view}
\figskip
\end{figure}

Finally, let us look in the subnested (3rd) pyramid. Assuming that
there is a disintegration of the bosons B4 and B14, it is likely
that degenerating Chevallier we can get $4(1,2 \frac{13}{1})$ and
also moving above there is constructions by Korchagin. So the
safest way to kill the boson $4(1,2 \frac{14}{1})$ is to put a TOR
on $4(1,0 \frac{14}{1})$. Likewise to kill the boson $14(1,2
\frac{4}{1})$ the sole reasonable choice seems to put a TOR on
$14(1,0 \frac{4}{1})$. As to  Shustin's $M$-obstructions on the
top of the subnested pyramid they can be either explained by
putting central TORs doing severe damages on the subordinated
$(M-1)$-schemes, or by just placing left corner TOR's damaging
only the left versant of each hills climbing to Shustin's
(anti)-schemes. Curiously, enough we did not as yet gathered
enough data to decide which option is more likely.

[23.08.13] It would be interesting to see if there is a coherent
distribution of TOR if we put in discredit Shustin's
constructions.

Next it is tempting to look at the walls crossing past the
discriminant as yielding the magmatic dynamic of the whole
pyramid. Especially we had the idea of looking at eversions of the
nonempty oval. Of course empty ovals can in principle be shrunk,
but not so with nonempty ovals which can instead be eversed. Under
this operation,
%%%described in Sec.\,\ref{Eversion:sec}
%
all what what inside the oval appears outside of it and viceversa
the outside of the oval is now captured in the inside of the
eversed oval. Of course this is much akin to Steiner's
Wiedergeburt und Neuauferstehung, when it comes to inversion. If
we imagine a trinested scheme then its eversion will be binested.
If we evert an non empty oval of a trinested scheme then the
result will have 2 subnests and so violates B\'ezout. Hence
eversion cannot be performed on trinested curves. Also, we can
imagine to evert uninested curves, naturally turning to
themselves. Yet, starting with say Harnack's scheme
$18\frac{3}{1}$ we would get $3\frac{18}{1}$ which does not verify
Gudkov periodicity (even in the simple version of Arnold). So
unlike the case $m=6$, the case $m=8$ has some asymmetry in its
Gudkovian structure, yet this is probably restored in degree
$m=10$, etc.

Notwithstanding eversions could give some secret passage to travel
in the pyramid and so perhaps aid to guess the isotopic
classification, especially that of the still open bosons. For
instance let us consider the boson $14(1,2\frac{4}{1})$. On
everting the biggest nest we get the scheme
$2\frac{4}{1}\frac{14}{1}$ which does not satisfy Gudkov
periodicity. Okay, but maybe it is possible for the oval to evert
as to capture only one outer oval, but the
Fig.\,\ref{ViroDEGREE8_eversion:fig} prompts rather the contrary
intuition. Hence it seems that in degree 8 Gudkov periodicity
forbids all form of eversions.

\begin{figure}[h]\Figskip
%\vskip-1.2cm\penalty0
%\centering
\hskip-2.7cm\penalty0
\epsfig{figure=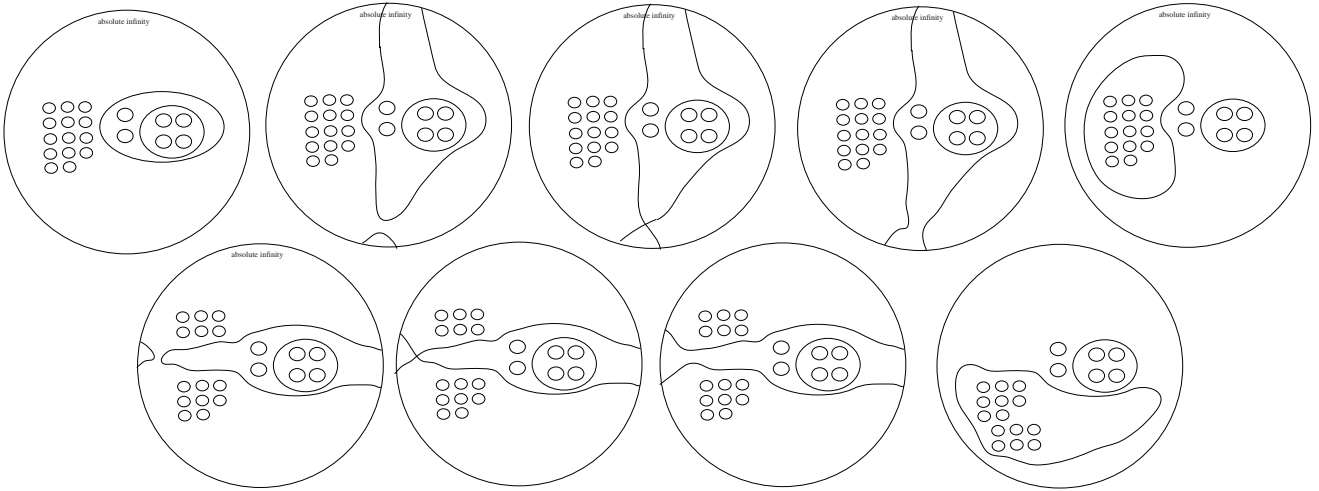,width=172mm}
\captionskipAG
  \caption{\label{ViroDEGREE8_eversion:fig}%
  Eversion}
\figskip
\end{figure}

Now it remains to tackle the problem of speculating about the
falsity of Shustin's construction as to see if---using only Viro's
potentially stabler construction (as using only the quadri-ellipse
as ground curve)---it is possible to distribute TOR in such a way
that censorship explains all prohibitions of Fiedler, Viro,
Shustin, Orevkov and perhaps some few more not yet known to exist.

Of course another option would be that the Fiedler-Viro oddity law
is wrong and then we would have less obstructions and virtually no
obstructions, as Gudkov somehow conjectured at the end of his 1974
survey. (Curiously, there, Gudkov seems to count 102 schemes
instead of the ca. 104+40=144 that exists prior to applying the
Fiedler-Viro amendment (which kill forty guys).

Now we first try to mistrust Shustin's constructions involving the
medusa (Fig.\,\ref{ViroDEGREE8_SHUSTIN:fig}), plus the other
construction of Shustin ($4\frac{5}{1}\frac{5}{1}\frac{5}{1}$)
which we failed as yet to understand. Killing all those Shustin's
scheme or perhaps keeping the last one, causes a little trouble
because right below Shustin's schemes there are schemes
constructed by Polotovskii, and so it looks hard to explain our
postulated prohibition of Shustin via censorship, without killing
simultaneously Polotovskii's scheme. In fact we  constructed most
of Polotovskii's schemes just by extrapolating Viro's
$M$-parameters to $(M-1)$-parameters, i.e. with 8 micro-ovals
nascent instead of 9.

\begin{figure}[h]\Figskip
%\vskip-1.2cm\penalty0
%\centering
\hskip-2.7cm\penalty0
\epsfig{figure=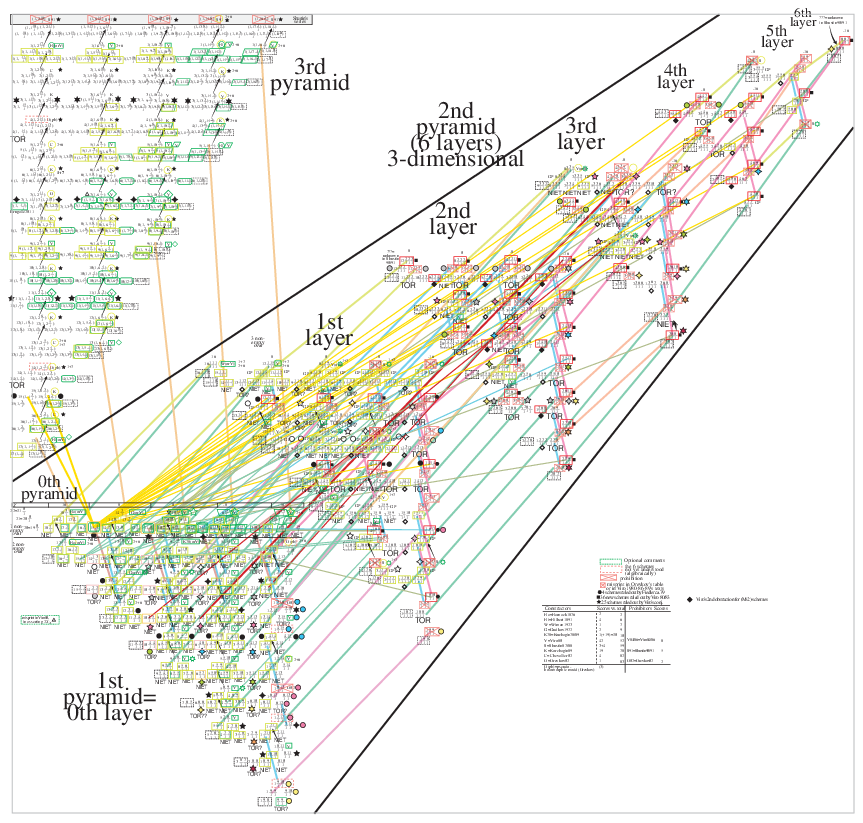,width=172mm}
\captionskipAG
  \caption{\label{Degree8-(M-i)-curve-TABLEShustin:fig}%
  Trying to depict all enlargements (a mess!)}
\figskip
\end{figure}

[25.08.13] Crudely put, we may accept as valid only those of
Viro's construction attained by the most basic recipe of the
quadri-ellipse, considering e.g. Viro's more exotic constructions
as fallacious too (like Shustin's).

[27.08.13] After one day of vacation with Markus Jura Suisse and
Traugott Schneider, we noted that the scheme
$10\frac{2}{1}\frac{6}{1}$ seems to respects censorship, because
its enlargements (lying in the 2nd layer) are killed by Viro's law
and a Shustin $(M-1)$-obstruction. So it seems relevant to check
if there is any concomitance between Viro's law and Rohlin's
maximality principle for RKM-schemes.

Let us do the search systematically, starting  with the 1st
RKM-scheme namely $15\frac{4}{1}$. As we already remarked a long
time ago (ca. June 2013, yet after posting v2 of this text) we
have the black bullets showing extensions of this $(M-2)$-scheme
and those are often dominated by Viro's scheme so that Rohlin's
maximality principle is trivially corrupted. However according to
a refined interpretation of the latter we can just infer that
$15\frac{4}{1}$ lacks TOR (total reality). So we cannot expect
(provided Viro's method is as reliable as it fame is) that any RKM
scheme ($\chi\equiv_8 k^2+4$) prompts the phenomenon TOR. Hence
practically, this means that we may tag the label NIET on
$15\frac{4}{1}$.

Further the black circle on
Fig.\,\ref{Degree8-(M-i)-curve-TABLEShustin:fig} shows  that there
are interesting extensions from the simply-nested to the subnested
nested realm by adding a separating oval in the nest enclosing a
certain number of eggs while leaving the other outside. Also the
same separation can be applied to the outer ovals leading thus to
the black-circled schemes in the binested and trinested realms.
Especially we reinterpreted 3 of Viro's oddity law via Rohlin's
maximality principle. However this cannot be done in an uniform
fashion because of the extensions $13\frac{2}{1}\frac{4}{1}$,
hence $15\frac{4}{1}$ despite being RKM lacks total reality, and
must be ascribed the (poor-quality) label ``NIET''.

The same rating (ranking) ``NIET'' must be given (for the same
reason) to $11\frac{8}{1}$, $7\frac{12}{1}$ and  $3\frac{16}{1}$.

A ``NIET'' must also be given to $14\frac{1}{1}\frac{3}{1}$,
$10\frac{1}{1}\frac{7}{1}$. When it comes to
$6\frac{1}{1}\frac{11}{1}$ we find 3 extensions in the 1st layer
yet all dominated by Shustin's $M$-schemes, so that if the latter
are erroneously constructed we could posit TOR for
$6\frac{1}{1}\frac{11}{1}$. But of course if Shustin's
construction are solid this scheme should be assigned a NIET. Alas
as we said Shustin's construction cannot hold true in a way
compatible with TOR and censorship.

For $2\frac{1}{1}\frac{15}{1}$ we have only one extension in the
1st layer which correctly prohibited by Shustin and Fiedler. Due
to the binesting there is no extensions in the subnested realm.
Hence it seems that the scheme $2\frac{1}{1}\frac{15}{1}$ is the
first deserving the label ``TOR''.

Then $14\frac{2}{1}\frac{2}{1}$ has several extensions in the 2nd
layer, all prohibited by either Shustin or Viro's law, safe for
the $(M-1)$-scheme $12 \frac{2}{1}\frac{2}{1}\frac{2}{1}$ (of
Shustinian uncertainty). Accordingly it seems plausible to ascribe
a TOR to scheme under inspection.

The same verdict---with more certitude even---can be applied to
$10\frac{2}{1}\frac{6}{1}$.

For $6\frac{2}{1}\frac{10}{1}$ we have likewise a TOR. For
$2\frac{2}{1}\frac{14}{1}$ we have only one extension, which is
prohibited by Shustin and so we can expect a TOR.

In contrast for $10\frac{3}{1}\frac{5}{1}$ we have extensions
below a Viro's $M$-scheme and so we must impose a NIET.

For $6\frac{3}{1}\frac{9}{1}$ we have 2 (pure) extensions in the
3rd layer, which are (correctly?) prohibited by Viro-sporadic and
Shustin, and Viro-regular. So it seems advisable to ascribe a TOR
to the scheme under scrutinity.

For $2\frac{3}{1}\frac{13}{1}$ we have one pure extension in the
2nd layer (Shustin prohibited) and one more (impure) in the 3rd
layer (Viro's sporadic). Hence the label TOR seems conceivable for
the scheme investigated.

Then it comes to $\frac{3}{1}\frac{15}{1}$---which albeit not
RKM---has only one extensions in the 1st layer (Viro sporadic),
and thus we may expect the  TOR label. {\it Warning:} one should
not miss two extensions in the 0th layer, including the boson
$1\frac{4}{1}\frac{15}{1}$ which is therefore virtually
prohibited. Of course there are also 2 direct extensions abutting
to Orevkov's anti-scheme $1\frac{3}{1}\frac{16}{1}$ via the
$(M-1)$-scheme  $\frac{3}{1}\frac{16}{1}$ which to the best of our
knowledge is not yet realized nor prohibited (i.e. bosonic for
short).

The story continues with $10\frac{4}{1}\frac{4}{1}$ where we have
extensions in the 4th layer all prohibited by Viro's law or
Shustin's $(M-1)$-avatars. Actually one must not miss extensions
in the 2nd and 3rd layers, yet all prohibited, and thus we may
ascribe TOR to the scheme inspected.

For $6\frac{4}{1}\frac{8}{1}$ we have 3 extensions (in the 2nd,
3rd and 4th layers) all prohibited by Shustin or Viro's regular
law, so that the TOR label is  expectable.

The scheme $2\frac{4}{1}\frac{12}{1}$ permits only one pure
extension in the 2nd layer and one impure one in the 3rd layer.
Those being prohibited either by Shustin or Viro regular, the TOR
label at the position inspected is expectable.

For $6\frac{5}{1}\frac{7}{1}$ it has an extension in the 3rd layer
subsumed to a Viro $M$-construction, and this suffices to prevent
a TOR at the inspected place. So the rating agency ascribes a NIET
to the scheme inspected.

For $2\frac{5}{1}\frac{11}{1}$ we find extensions in the 2nd and
3rd layers (all prohibited) sot that a TOR label is expectable.

For $6\frac{6}{1}\frac{6}{1}$, we find extensions in the 2nd, 3rd
and 6th layers. The first two are prohibited while the 3rd is
still-open in Shustin 90/91
\cite{Shustin_1990/91-New-restrictions}. Accordingly, it is not
easy to make a decision yet the TOR label is likely.

For $2\frac{6}{1}\frac{10}{1}$ we have extensions in the 2nd and
3rd layers all prohibited by either Viro or Shustin. Thus a TOR
label is quite likely.

Then we have $\frac{6}{1}\frac{12}{1}$, which admits only one
extension in the 1st layer ($\frac{1}{1}\frac{6}{1}\frac{12}{1}$),
yet prohibited by Fiedler's special case of Viro. {\it Warning:\/}
actually, one should not miss two extensions in the 0th layer,
including a possible prohibition of the boson
$1\frac{7}{1}\frac{12}{1}$.

Next we must examine $2\frac{7}{1}\frac{9}{1}$, which admits
extensions (only) in the 2nd and 3rd layers (all prohibited by
Viro sporadic or Shustin). Hence, the TOR-label is expectable.

As to $2\frac{8}{1}\frac{8}{1}$, we detect 3 extensions in the 2nd
and 3rd layers (all prohibited by Shustin or Viro regular), so
that the TOR-label  is probable.

Finally, for $\frac{9}{1}\frac{9}{1}$ we have extensions in the
0th layer (all bosonic, i.e. unknown) or in the 1st layer (Viro
sporadic). Thus, it is plausible (despite sparse experimental and
theoretical knowledge) that a TOR label is present on this scheme.

Curiously enough, it seems that there is no way to rule out the
boson $b1:=1\frac{1}{1}\frac{18}{1}$ via censorship unless one
puts a TOR on $1\frac{18}{1}$ but this would (by censorship)
jeopardize Viro's horse construction of $2\frac{19}{1}$.

%%SAXOOO:
DO-NOT-FORGET: No\'emie Combe showed me today a set of
notes by Viro where in ca. Chap. 4 by page ca. 60 Fiedler's
alternation rule is explained in (seemingly) full details. I
cannot remember if I cited this or analyzed properly this source
which is perhaps superior than the original papers (Fiedler 83,
Viro 83).

[28.08.13] {\it Little side-remark:\/} As we know since a long
time (cf. e.g. \ref{RMC:cter-example-via-Viro-1st-curve=BEAVER})
Rohlin's maximality conjecture is trivially foiled by Viro's
construction. To be specific the RKM-scheme $15\frac{4}{1}$ admits
an extension as the $M$-scheme $13 \frac{2}{1}\frac{5}{1}$ which
is actually even constructed by Gudkov (and arguably in a more
elementary/versatile way by Viro). This explains why Viro was in
his letter (cf. Sec.\,\ref{e-mail-Viro:sec}) fairly sure that
Ahlfors' theorem cannot prove Rohlin's maximality conjecture, just
because the latter false. So:

\begin{Scholium} Rohlin's maximality conjecture is trivially
false because of Viro's method, yet the latter is too polite to
attack frontally his teacher leaving thus in literature a little
cloud of unclearness (apart from the little critiques expressed by
Shustin, yet not explicit enough regarding the direct sense of
Rohlin's conjecture: type~I scheme is maximal). Actually, it
suffices even to know Gudkov's construction of the $M$-scheme
$13\frac{2}{1}\frac{5}{1}$ to corrupt Rohlin's maximality
conjecture in view of the enlargements $15\frac{4}{1}<
13\frac{2}{1}\frac{4}{1}<13\frac{2}{1}\frac{5}{1}$.
\end{Scholium}

{\it Back to the main diagrammatic.}---Yesterday, we traced the
faisceaux of all strokes emanating from the uninested RKM-scheme
(e.g. $15\frac{4}{1}$). Here there is a bunch of extensions
radiating in higher layers but all those rays abut always (as the
comprehend the sub-symbol $\frac{4}{1}$) to schemes prohibited by
Viro's oddity law. As a moral imposing a TOR on such schemes
($15\frac{4}{1}$, $11\frac{8}{1}$, $7\frac{12}{1}$,
$3\frac{16}{1}$) would explain a substantial part (all?) of Viro's
law but alas in a way not defendable since by Viro's method there
is also extensions in the 1st pyramid=0th layer.

Let us now continue our game of rating-agency for the labels NIET
and TOR granting Viro's model of the theory as being the reliable
one. Yesterday we rambled the whole 1st layer and it remains now
to elevate to the higher layers.

The 1st case of interest is the RKM-scheme
$13\frac{1}{1}\frac{1}{1}\frac{2}{1}$. Surprisingly, the latter
lacks any extension in the pyramid (just because of Gudkov and GKK
essentially). Therefore the TOR-label looks amply merited.

The same applies to the other RKM-schemes in the 1st row of the
1st layer (i.e. $9\frac{1}{1}\frac{1}{1}\frac{6}{1}$,
$5\frac{1}{1}\frac{1}{1}\frac{10}{1}$,
$1\frac{1}{1}\frac{1}{1}\frac{14}{1}$). The basic raison seems to
be that the double occurrence of ones forces any extension to
stays in the 1st layer, except if we increase both one numerators
simultaneously, e.g. as $9\frac{2}{1}\frac{2}{1}\frac{6}{1}$ yet
this lands outside the range specified by Gudkov periodicity. In
conclusion, all those RKM schemes can be ascribed the label TOR
yet it does not result any prohibition beyond Gudkov periodicity.

The next case of interest is
$10\frac{1}{1}\frac{2}{1}\frac{4}{1}$. Augmenting successively
each coefficient of this symbol we always get schemes out of
GKK-periodicity, so that the scheme admits only those extensions
immediately visible above it, but prohibited by Shustin and Viro's
oddity law. Hence the TOR-label is likely. Additionally, it may be
noted that those 2 prohibitions (above this scheme) cannot bwe
induced from one of the 3 schemes of the 1st layer below
$11\frac{1}{1}\frac{2}{1}\frac{4}{1}$, obtained by conserving two
subfraction of the symbol.

For $9\frac{1}{1}\frac{2}{1}\frac{5}{1}$ we find again no
extension and thus the TOR-label is likely.

The same holds for $6\frac{1}{1}\frac{2}{1}\frac{8}{1}$ and
$5\frac{1}{1}\frac{2}{1}\frac{9}{1}$.

Next it comes to $4\frac{1}{1}\frac{2}{1}\frac{10}{1}$ whose
status becomes ambiguous as  we reject Shustin's construction.
Nonetheless it can be that a Viro construction from the
quadri-ellipse produces the $(M-1)$-scheme
$4\frac{1}{1}\frac{2}{1}\frac{11}{1}$ as asserted by Polotovskii.
In that case we must assign the NIET-label to both
$4\frac{1}{1}\frac{2}{1}\frac{10}{1}$ and
$3\frac{1}{1}\frac{2}{1}\frac{11}{1}$. Actually, we realized
$4\frac{1}{1}\frac{2}{1}\frac{11}{1}$ via a damped version of
Shustin's construction (see
Fig.\,\ref{ViroDEGREE8_SHUSTIN_NEW:fig}b2). However one could
expect that there is a cleaner version via Viro's quadri-ellipse.
At any rate it seems worth expanding the case of Viro's
construction that we missed. However it seems that all extensions
of the scheme $4\frac{1}{1}\frac{2}{1}\frac{11}{1}$ never
interacts with a Viro $M$-scheme, so that it seems unlikely that
this $(M-1)$-scheme can be constructed by Viro's method involving
the quadri-ellipse (at least for the small collection of patches
available to Viro). Hence it seems plausible that discrediting
Shustin's construction discredits as well Polotovskii's
construction of $4\frac{1}{1}\frac{2}{1}\frac{11}{1}$ although it
is not clear whether Polotovskii's work is logically subsumed to
 Shustin's. Using our exotic variant of Viro's method there is soe
 hope to get $(M-1)$-scheme by using maximal dissipation of an
 exotic sort. For a catalogue of such opportunities cf. the yellow
 shaded combination of patches shown on
 Fig.\,\ref{ViroDEGREE8_exotic_patches0_SYS2:fig}, but alas none
 of which is readily accessible to Viro's confined toolkit of
 patches.

 Nonetheless, it seems clear that the resulting $(M-1)$-patchwork
 will often
 conflicts with Shustin's $(M-1)$-obstructions so that we should
 discover new patching obstructions. Looking at
 Fig.\,\ref{ViroDEGREE8_exotic_patches0_SYS2:fig} we see four
 tri-ovals curves on the 1st A-row, but as A is an empty patch
 family we will not be able to extract any obstruction. Likewise
 the B-row will not be much instructive say just because we lack
 any concrete patch in the B-family. More promising is the C-row,
 which interact in the yellow-way (3 macro-ovals) with the D-patch
 (which is already known to be empty via B\'ezout and Arnold). In
the C-row one must not forget the flipped combination (in the
margin of Fig.\,\ref{ViroDEGREE8_exotic_patches0_SYS2:fig}), but
those interacts with patches F and H, already known to be void
families. Next for the D-row we have a yellow interaction with
column G, but as D is empty this will not supply any bit of
information. Next we arrive at Viro's inhabited row E, where we
get yellow $(M-1)$-interaction with column F (which is already
known to be empty). As to the flipped versions we get an
interaction with H, which is actually already known to be empty.
Hence:

\begin{Scholium}
The many Shustin's $(M-1)$-obstructions do not cause any damage
upon  (exotic) patches that could not already have been drawn from
B\'ezout or Arnold.
\end{Scholium}

This is quite disappointing, but in reality is not as it leaves
room for more patches than in Viro's theory, and so perhaps the
materialization of some bosons not yet known to exist.

Let us return to $4\frac{1}{1}\frac{2}{1}\frac{11}{1}$. Can we
construct it by Viro's quadri-ellipse? A priori there are two
ways. Either one uses maximal dissipation yet not optimally
reglued so as to have only 3 macro-ovals like on the yellow cases
of Fig.\,\ref{ViroDEGREE8_exotic_patches0_SYS2:fig}. (As we just
saw this variant leads nowhere because Viro's patches only
interacts in the appropriate $(M-1)$-way with patches families
known to be empty). Alternatively one can use an optimal gluing
(red on the same Fig.), but with a damped dissipation with only 8
instead of 9 micro-ovals. This is a sort of social dumping that
could maybe create new-or-old  $(M-1)$-schemes. For instance one
may start with the $B$-patches and glue it with themselves modulo
a social dumping. A priori it could seem that if the patches are
already restricted according to Viro's law, it will appear as
difficult to reach Shustin-prohibited $(M-1)$-scheme as the latter
often lies below those excludes by Viro. However a more systematic
search deserves to be done. So the philosophy is always the
alienating ``patchworker de toutes les mani\`eres possibles et
imaginables''. So for instance we may try to consider B3/B3flip
with one factor damped and the resulting table will be essentially
the same as that of H5/H5flip modulo the damping. Evidently a
priori we may draw little information from the experiments unless
we suppose a strong correlation principle akin to the
Viro-Itenberg principle of contraction of empty ovals (in a
germinal version for patches). So we shall assume that whenever a
patch exists, its versions where empty ovals are imploded also
exist, and then we expect additional prohibitions coming from
Shustin's $(M-1)$-prohibitions (thought the latter look merely a
compromise between Viro-Shustin constructions and Viro's
prohibitions).

On working the table, we can comment as follows. On composing
(8,1,0) with (5,3,0) we find an anti-Shustin scheme. So if one
could prove existence of the patch B3(5,3,0), this would rules out
the $M$-patch B3(8,1,0) (a task we are as yet unable to complete).

[29.08.13] It is evident that Hilbert's 16th problem is completely
determinist, but nobody understand the determinism. Now continuing
our boring task of filling that table it seems evident at an early
stage already that as Shustin's prohibitions are much dominated by
those of Viro's that no new patch obstructions will be derived.

However as a little surprise the entry $(2,3,4)\times (5,3,0)$
creates an anti-Shustin scheme so that one can infer:

\begin{lemma}
Either the $M$-dissipation $B3(2,3,4)$ or the damped
$(M-1)$-dissipation $B3(5,3,0)$ is prohibited provided Shustin's
$(M-1)$-avatar of Viro's most sporadic obstruction is true.
Similarly, via $(2,3,4)\times (6,2,0)$, either the $M$-dissipation
$B3(2,3,4)$ or the damped $(M-1)$-dissipation $B3(6,2,0)$ is
prohibited provided a Shustin obstruction (below a Viro regular
one) is right.
\end{lemma}

Alas, it must be confessed that unless one is able to construct
the damped patches this supplies no additional information on the
$M$-dissipation of type B3. Further, our initial hope to meet
Shustin's obstruction along the diagonal (or pseudo-diagonal we
should say in view of the damping) is not borne out, or rather
when it occurs it is already covered by Viro's $M$-obstructions.
So in conclusion it seems that nothing tremendously new follows
from Shustin.

\begin{figure}[h]\Figskip
%\vskip-1.2cm\penalty0
%\centering
\hskip-2.7cm\penalty0
\epsfig{figure=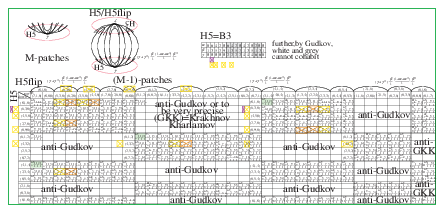,width=172mm}
\captionskipAG
  \caption{\label{ViroDEGREE8_exotic_patches_B3:fig}%
  Patchworking exotic patches at the $(M-1)$-level: B3/B3flip}
\figskip
\end{figure}

On the other hand we are not much advanced on the question of
deciding if there is a construction of the $(M-1)$-scheme
$4\frac{1}{1}\frac{1}{1}\frac{12}{1}$ \`a la Viro (quadri-ellipse)
independent of Shustin's tricky medusa. Of course this scheme
appears on our table
(Fig.\,\ref{ViroDEGREE8_exotic_patches_B3:fig}) but is by no mean
a regular construction unless one is able to prove existence of
the involved patches.

Now let us take again a global look at the pyramid
(Fig.\,\ref{Degree8-(M-i)-curve-TABLEShustin:fig}). On it we
already experimented that not all RKM-schemes seems to be capable
of the phenomenon of total reality at least in the strong way
implying censorship of all enlargements (i.e. maximality of the
scheme in the sense of Rohlin). For instance
$10\frac{3}{1}\frac{5}{1}$ cannot be TOR (so is NIET) because
there exist extension by Viro's method (in the 2nd and 3rd layers,
subordinated to Viro's $M$-schemes 8.3.3.5 and 4.3.5.7). So in a
very radical world where much importance is given to our
Rohlin-style maximality principle it could be that even Viro's
constructions must be revised.

Loosely out it seems that the architecture of the Gudkov-Viro
pyramid in degree 8 (alias the Grand-Pyramid of Gizeh) and more
generally any Gudkov pyramid of higher degree is governed by two
principle a bit akin to the well-known inclusion-exclusion of
rudimentary combinatorics. One the one side there is a principle
of inclusion saying roughly that whenever we have a scheme
algebraically realized all smaller schemes will also be present.
This means roughly that a stone of the pyramid (or cathedral if
you prefer Judeo-Christianism over Egyptian civilization) cannot
stay in the air in levitation without supporting elements. This
intuition is sustained by the Viro-Itenberg contraction principle
for empty oval.

On the other hand there is---Rohlin's intuition, and to some
extend also Wiman 1923 \cite{Wiman_1923}---a dual principle of
exclusion that whenever a scheme is subsumed to total reality it
should kill all its enlargements. This we call the censorship
principle. It imposes roughly the architecture of the cathedral
being pure, i.e. without too much ``fioritures''=embelishments \`a
la Gaudi (in Barcelona). Alas it seems that this principle cannot
hold at all RKM-positions (i.e. $(M-2)$-schemes with $\chi\equiv_8
k^2+4$) in view essentially of Viro's rich method of construction.
Actually, even more basically the RKM-scheme $15\frac{4}{1}$
violates already the censorship principle (interpreted as a sort
of denuded Roman architecture) since it accepts the
$M$-enlargement $13\frac{2}{1}\frac{5}{1}$ due to Gudkov (showing
incidentally that the ``direct sense'' of Rohlin's maximality
conjecture is foiled). Could it be that this construction of
Gudkov is erroneous (and so are many of Viro construction) in
order to restore Rohlin's exclusion principle (censorship).

\subsection{The dream/nightmare of a world with censorship}

[29.08.13] We can dream of a world where Rohlin's maximality
principle is true in the strong form of its original formulation,
any RKM-scheme is of type I (true theorem of
Rohlin-Kharlamov-Marin), and therefore maximal=saturated in the
sense of Rohlin. (This last clause is in principle false, as we
just noted and actually was first disproved by Gudkov prior than
Rohlin formulated his conjecture). So rating-agency conclusion:
Rohlin did not perfectly his homework in 1978 but is completely
excused as he already endured an hearth attack ca. 1976 (compare
Vershik's obituary ca. 1986). Notwithstanding we may adhere to the
radical position that Rohlin's principle is universally true, in
which case there would be much less $M$-curves than in the
politically correct world (where all published results are
considered as true). Of course in reality, we may have a more
nuance landscape where the censorship principle applies only to a
subcollection of RKM-schemes and perhaps other $(M-2)$-schemes,
typically if the are subsumed to the phenomenon of total reality
of a pencil which via B\'ezout should prevent any extension of the
scheme, whence the censorship.

In this section we adopt the most radical attitude, in order to
see which $M$-schemes are really resistant.

First, the uni-nested RKM-schemes admits many extensions in the
3rd=subnested pyramid killing thereby many $M$-schemes claimed by
Korchagin, Viro, Orevkov-Viro, Korchagin-Viro (compare the red
faisceaux on Fig.\,\ref{Degree8-(M-i)-curve-TABLERohlin:fig}).
Diagrammatically, it seems more efficient to just use our circles,
stars, etc. propagation of symbols along the extension of a given
RKM-scheme. Positing Rohlin's maximality principle causes then a
hecatomb of $M$-scheme also in the 1st pyramid, all marked by a
red rectangle.

Alone the RKM-scheme $15\frac{4}{1}$ would kill a plethora of
$M$-schemes (marked by black circles on the figure), and so on for
all other RKM-schemes. One could expect that this hecatomb forced
by a strong Rohlin's maximality principle is still compatible with
Viro's theory in the sense that the last survivals $M$-schemes
(which are easy to list explicitly) are all constructible via
Viro's method for a subcollection of his patches (taking for
granted that he may have had a too liberal acceptation of patches
based upon fraudulent constructions).

\begin{figure}[h]\Figskip
%\vskip-1.2cm\penalty0
%\centering
\hskip-2.7cm\penalty0
\epsfig{figure=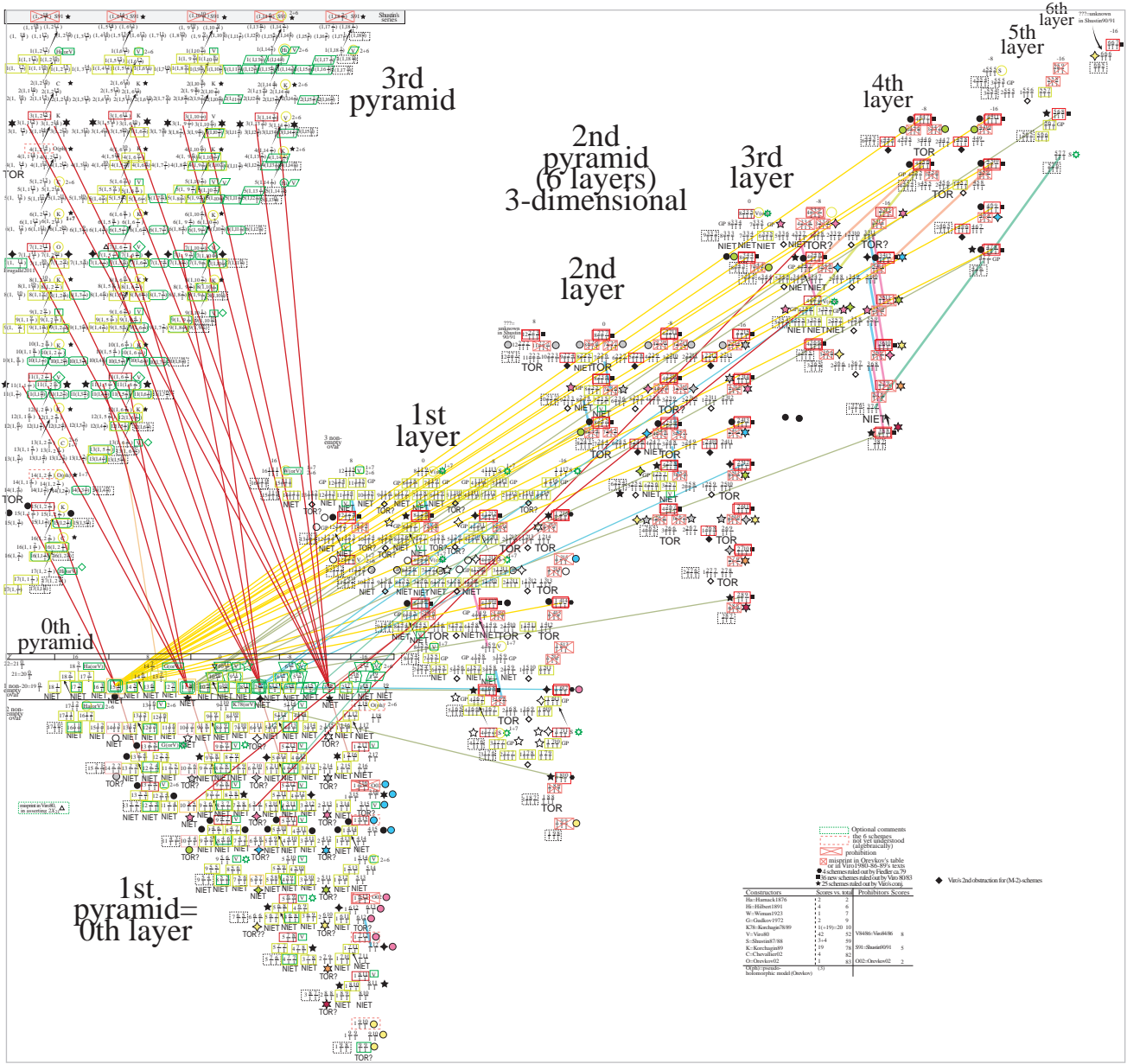,width=172mm}
\captionskipAG
  \caption{\label{Degree8-(M-i)-curve-TABLERohlin:fig}%
  Assuming the truth of Rohlin's maximality principle:
  an epured architecture killing Gudkov, Viro, etc.}
\figskip
\end{figure}

At this stage one loose a lot of energy due to spatial
distantness, hence it seemed advisable to manufacture a colimasson
depiction to save space. This representation permits one to save
much energy dispensed in scrolling the windows in the computer.

\begin{figure}[h]\Figskip
%\vskip-1.2cm\penalty0
%\centering
\hskip-2.7cm\penalty0
\epsfig{figure=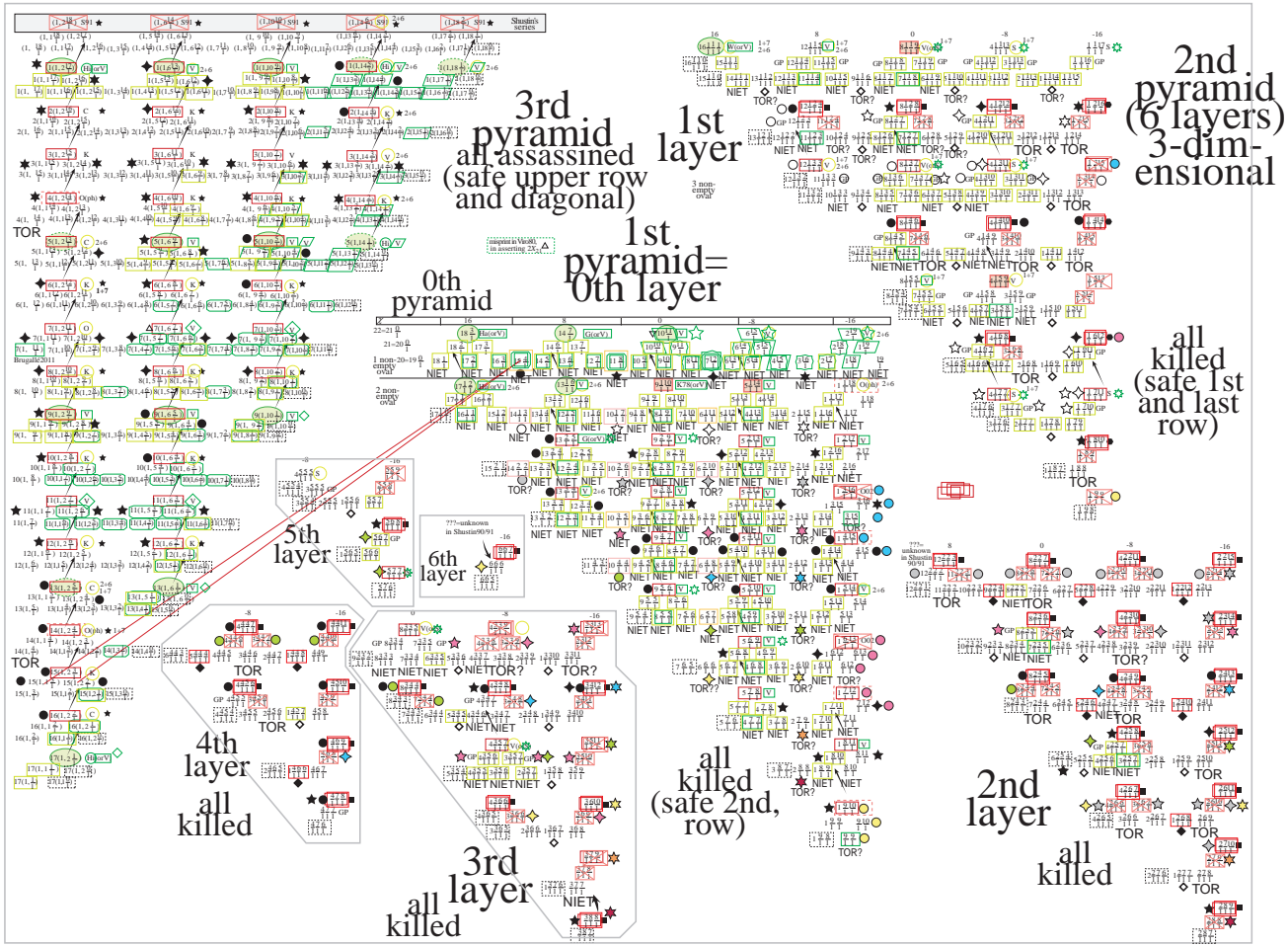,width=172mm}
\captionskipAG
  \caption{\label{Degree8-(M-i)-curve-TABLERohlinBIS:fig}%
  Assuming the truth of Rohlin's maximality principle:
  an epured architecture killing Gudkov, Viro, etc.}
\figskip
\end{figure}

Of course albeit our scenario of Rohlin maximality is quite
apocalyptic for Viro (even Gudkov) it is much in line with our
paradigm of total reality and the parsimony of $M$-schemes prior
to Viro's intervention.

[30.08.13] So now we hope to have detected all extensions of
RKM-schemes and it is time to look at a restricted Viro theory
that would corroborate this. For this it suffices to take Viro's
original table of compositions (patchwork) and censure what is
forbidden by Rohlin's maximality principle. We find that the first
frank obstruction (along the diagonal) occurs in the V2/V2-table
with the entry $(4,4,1)$ which is forbidden. It results a
destruction of all shaded schemes, which all admits a replica so
that the are definitely killed (or rather lost). The next serious
damage occurs along the diagonal entry (8,0,1), which is therefore
killed, and it results severe lost of schemes on the corresponding
row. Namely it seems that we loose the K78-scheme of Korchagin
1978, and the scheme $5\frac{1}{1}\frac{14}{1}$ of Viro. Of course
crosses not  situated on the diagonal also causes collateral
damages yet of an indefinite nature. Notwithstanding there is
latent patch damages impeding that the patch list stabilizes to a
sole censorship of $V2(4,4,1)$ and $V2(8,0,1)$.

Then the table V2/V3 only gives uncertain patch damages and also
jeopardizes some of Viro's constructions, namely of
$8\frac{1}{1}\frac{1}{1}\frac{9}{1}$ and
$4\frac{1}{1}\frac{5}{1}\frac{9}{1}$.

To explain the collateral damages it seems likely to abandon the
patch V2(2,4,3) which has well on V2/V2 as on V2/V3 monopolize
many Rohlin-forbidden schemes. The sole negative side-effect is a
lost of Viro's realization of $5\frac{5}{1}\frac{10}{1}$; not in
contrast that $9\frac{5}{1}\frac{6}{1}$ is not lost as it admits
another realization as V2(5,4,0)$\times$V2(0,4,5).

For the same reasons (concentration of anti-schemes) so as to
minimize the number of patches killed it seems wise to abandon the
patch $V2(2,0,7)$, and the sole side-effect is a inhibition of
Viro's construction of $1\frac{5}{1}\frac{14}{1}$. At this stage
all of Rohlin-style prohibitions have been taken into account safe
for Gudkov's scheme $13\frac{2}{1}\frac{5}{1}$ occurring in both
tables. So we must either abandon $V2(1,8,0)$ or $V2(0,4,5)$.
Killing (1,8,0) would have the effect that Viro's method would not
any more cover Harnack's elementary construction of
$18\frac{3}{1}$, and accordingly it may seem more advisable to
abandon $V2(0,4,5)$. It seems then that $9\frac{5}{1}\frac{6}{1}$
is lost, and likewise for $9\frac{1}{1}\frac{10}{1}$. Via the
table V2/V3, a serious damage is the lost of the patch $V3(5,3)$.
So at this stage Viro's method shrinks dramatically and is only
able to produce the green-colored schemes. So in this Rohlin-style
scenario Hilbert's 16th problem would be nearly settled yet still
mysterious. An exact statement is cumbersome but as follows:

\begin{Scholium}
If Rohlin's maximality principle is true, then many of Gudkov and
Viro's construction are foiled and in Viro's method the list of
patches must be severely shrunk. All schemes with a thin red frame
on Fig.\,\ref{Degree8-(M-i)-curve-TABLERohlinBIS:fig} are
prohibited, while deciding which one are realized depends upon the
exact stabilization of the theory, and cannot be decided on ground
of sole combinatorics. Yet, positing that Viro's method should
englobe Harnack's and Hilbert's, it seems that the most likely
frozen Viro's theory is only capable to produce the schemes marked
by green-ellipses on the table
(Fig.\,\ref{Degree8-(M-i)-curve-TABLERohlinBIS:fig}).
\end{Scholium}

\begin{figure}[h]\Figskip
%\vskip-1.2cm\penalty0
%\centering
\hskip-2.7cm\penalty0
\epsfig{figure=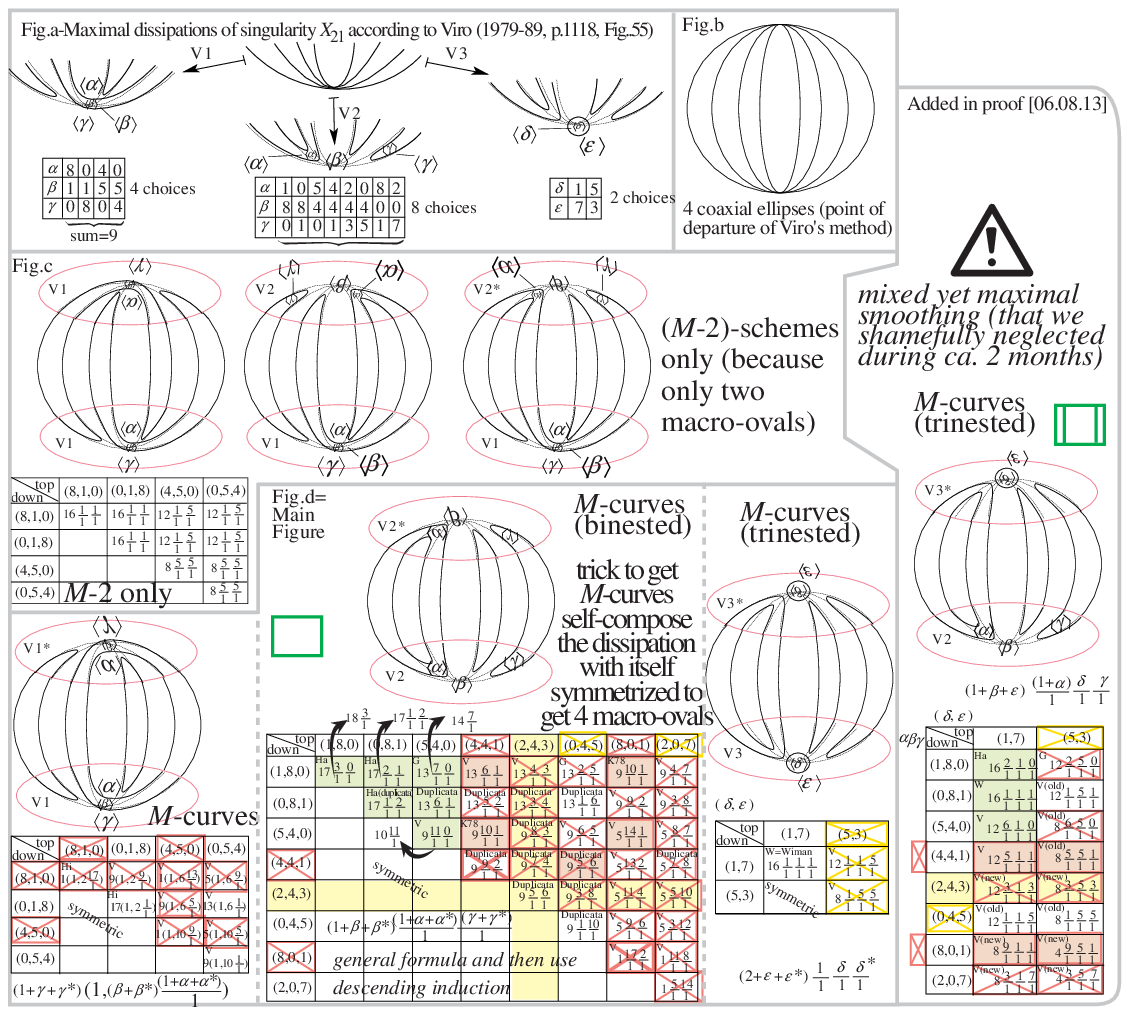,width=172mm}
\captionskipAG
  \caption{\label{ViroDEGREE8_censorship:fig}%
  Assuming the truth of Rohlin's maximality principle:
  an epured architecture killing Gudkov, Viro, etc.}
\figskip
\end{figure}

It is always puzzling that Viro's purest method (via the
quadri-ellipse) is not even able to reproduce two of Hilbert's
schemes, namely those with 14 big eggs (=ovals at depth 1 which
are empty). Naively put one could imagine that Viro's patching
parameters for $V1$ are too restricted. Copying the earlier
extended table V1/V1 we get the following table
(Fig.\,\ref{ViroDEGREE8_extendedROHLIN:fig}). This employs the
usual method of destruction along the diagonal, or diagonal
destruction along constructed rows. We see that Rohlin's
maximality principle of type~I schemes rules out all patches not
constructed by Viro, yet leaves moreover intact the patches
$(4,1,4)$ and $(0,9,0)$, which could exist. In that case both
Hilbert's schemes with 14 big eggs would be accessible to Viro's
method. Actually for this to be the case it suffices to have the
patch (0,9,0). Adding this create no more schemes except Viro's
horse-type  $1(1,18\frac{1}{1})$. Adding instead the patch (4,1,4)
create two schemes by Chevallier but otherwise nothing new to
Viro. So perhaps it quite likely that Viro's list of patches for
V1 is a bit overcautious and frozen (even under the stringent
Rohlin's maximality principle). On reporting the so-constructed
schemes on Fig.\,\ref{Degree8-(M-i)-curve-TABLERohlinBIS:fig}
(with dashed ellipses) we get a very regular pattern of fourfold
periodicity, yet with half of the rows not constructed nor
prohibited.

\begin{figure}[h]\Figskip
%\vskip-1.2cm\penalty0
%\centering
\hskip-2.7cm\penalty0
\epsfig{figure=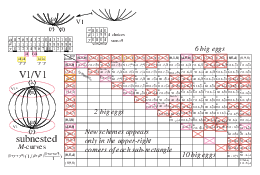,width=172mm}
\captionskipAG
  \caption{\label{ViroDEGREE8_extendedROHLIN:fig}%
  Extended table V1/V1 and patches selection dictated by
  Rohlin's maximality principle}
\figskip
\end{figure}

Apart from detail in the finishing it is already clear that
positing Rohlin's maximality principle leads to a strong deviation
from the actual state of knowledge (assuming the latter to be
true). Thus it is probably more likely that a refined version of
Rohlin's maximality principle which we shall call censorship
induced by total reality holds true. Actually as we saw even that
is not really compatible with Shustin's construction.

Then we noted that we forgot to notice that the boson
$b9=1\frac{9}{1}\frac{10}{1}$ is killed by a black star (extension
of $11\frac{8}{1}$) and likewise we omitted many black-circles on
the row of $9\frac{5}{1}\frac{6}{1}$. So we had to correct much of
our token. We also forgot initially all those extensions indicated
by the thin arrows in the 3rd pyramid (subnested case), and so we
get considerably more prohibitions in the 3rd pyramid. Yet,
upgrading those on the table V1/V1
(Fig.\,\ref{ViroDEGREE8_extendedROHLIN:fig}) does not cause any
additional patch damages (quite surprisingly).

Further it seems puzzling that Rohlin's maximality principle says
nothing on Shustin's scheme $4\frac{5}{1}\frac{5}{1}\frac{5}{1}$,
and more generally tells us little about the 5th layer (apart of
two prohibitions including one aggressing Shustin's medusa
construction of $\frac{5}{1}\frac{7}{1}\frac{7}{1}$). Basically
the reason for little information on Shustin's 4.5.5.5 seems to be
the fact its primitive antecedent in the 0th layer is
$8\frac{5}{1}\frac{5}{1}$ which fails being RKM.

Further starting from the scheme $3\frac{16}{1}$, one may add an
all enclosing oval phagocytizing the whole configuration. Alas
this violates Gudkov periodicity, but we may correct the situation
by leaving one oval outside. So we find the extension
$1(1,2\frac{16}{1})$, which was not reported as yet on our
tabulation. Alas the impact is rather dramatic, as Rohlin's
principle would then kill a scheme due to Hilbert. Hence:

\begin{Scholium}
Rohlin's maximality principle is not even compatible with
Hilbert's method of construction.
\end{Scholium}

Continuing we find analogous extensions for the other RKM-schemes,
e.g. $7\frac{12}{1}$ extends to $1(1,6\frac{12}{1})$, but also as
$5(1,2\frac{12}{1})$. Likewise $11\frac{8}{1}$ admits the 3
extensions $1(1,10\frac{8}{1})$, $5(1,6\frac{8}{1})$ and
$9(1,2\frac{8}{1})$. This continue along the obvious way, and it
results a generalized hecatomb of $M$-schemes in the 3rd layer,
apart few schemes resisting along the diagonal of this pyramid.
Actually two schemes of Hilbert are killed by Rohlin's maximality
principle which is therefore (granting the evident absence of
mistakes in Hilbert) trivially false. Of course in this scenario,
Viro's patches collection for V3 shrinks considerably, yet could
stabilize to the pair $(0,1,8)$ and $(0,5,4)$, in which case still
some 3 subnested $M$-schemes would be accessible to Viro's method
(namely the three on the bottom part of the diagonal): that is
Hilbert's $17(1,2\frac{1}{1})$, and Viro's $13(1,6\frac{1}{1})$
and $9(1,10\frac{1}{1})$.

All this looks very claustrophobic and apocalyptic, yet it is just
the combinatorial/rude consequence of Rohlin's (otherwise lovely)
principle pushed to its ultimate retrenchment.

It now of course time to better understand the substance of Viro's
method to check its truth, and potential extensions of its
flexibility. This is another topic that we shall try to analyze in
Sec.\,\ref{Viro:patches-construction:sec}.

Of course in the previous token about Rohlin's maximality as
corrupting even Hilbert's method, it is tacitly supposed by us
that the RKM-schemes are realized algebraically. This is not
completely obvious from zero knowledge, but we had somewhere in
this text vague realizations of those schemes via suitable
variants of Viro's method (hopefully hygienical). Compare for this
Lemma~\ref{RMC:cter-example-via-Viro-1st-curve=BEAVER} which is a
construction via Viro's beaver of $15\frac{4}{1}$, $11\frac{8}{1}$
and $7\frac{12}{1}$. It would be interesting to know if there is a
simpler construction say via Viro's quadri-ellipse.This seems to
be the case according to the green-rectangle on our table. Alas
browsing through our old tables (especially
Fig.\,\ref{ViroDEGREE8-(M-2):fig}) we did not found the requested
schemes.

In fact looking at
Fig.\,\ref{ViroDEGREE8_exotic_patches0_SYS2:fig}, we see that
using Viro's patches there is several ways to get $(M-2)$-schemes
via the gluing C/E, C/Cflip, C/Eflip (or equivalently E/Cflip) or
finally E/Eflip. It remains to check that we tabulated all those
patchworks. The figure below gathering all relevant information
seems to show that we missed as yet to tabulate the combination
C/Eflip, actually not since it corresponds to Fig.\,d. Inspecting
carefully the tables it seems however that we never get the
uni-nested RKM-schemes.

\begin{figure}[h]\Figskip
%\vskip-1.2cm\penalty0
%\centering
\hskip-2.7cm\penalty0
\epsfig{figure=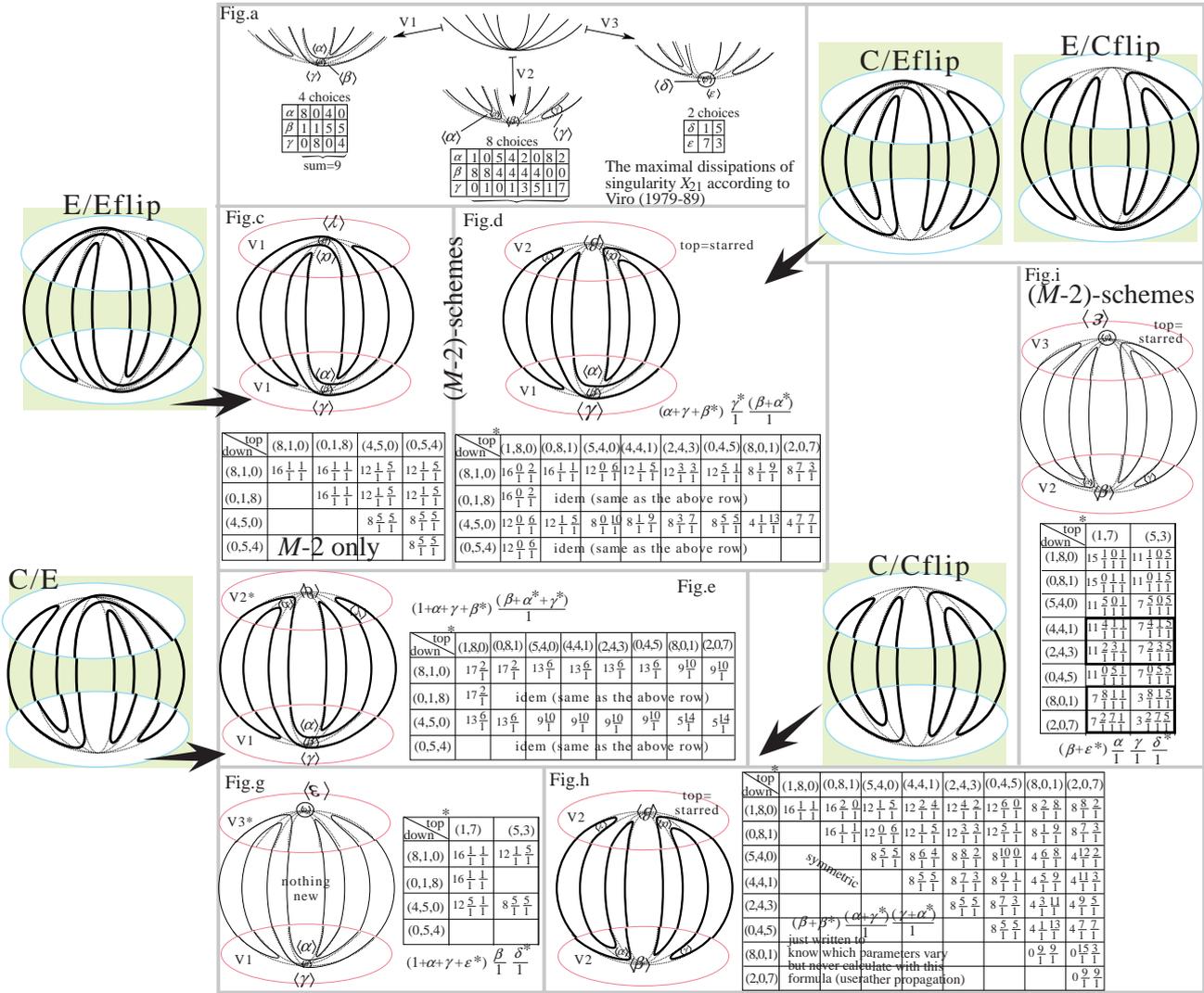,width=172mm}
\captionskipAG
  \caption{\label{ViroDEGREE8_exotic_patches0_SYS2BIS:fig}%
  Trying to get the RKM-schemes via Viro standard}
\figskip
\end{figure}

[31.08.13] In fact it seems that we may also interpret Rohlin's
maximality principle in a more restricted sense, namely by
considering only extensions by a circle bounding a disc disjoint
from the rest of the curve so that the initial scheme can be
thought of as resulting by shrinking an empty oval. This
interpretation is more in line with the contraction principle for
empty ovals (Klein-Viro-Itenberg), but apparently so trivially
true that it seems of little predictive significance when it comes
to explain prohibitions beyond Gudkov periodicity and its
Krakhnov-Kharlamov $(M-1)$-avatar.

Also puzzling are Viro's $(M-2)$-obstructions on RKM-schemes in
the 2nd and 4rt layers (black rhombs on the main-table). If those
are correct and interpretable in terms of a censorship allied to
total reality it would involve (at least) a scheme lying deep at
the $(M-4)$-level (recall Klein's obvious congruence $r\equiv_2
g+1=:M$), and somewhat along Arnold's philosophy of the mushroom
we could expect these as being very deep mushrooms explaining all
other prohibitions by resurfacing along some unexpected (but
completely deterministic) pathes of the pyramid. To speculate
about this issue it seems necessary to scroll out more the
diagrammatic up the $(M-4)$-level. This requests a new depiction.
As shown below there is actually two scenarios depending on
whether we fix the phenomenon of total reality. Perhaps the most
likely variant involves fixing total reality (TOR) of the RKM
$(M-4)$-schemes, i.e. those with $\chi\equiv_8 k^2+4$. However if
we put a TOR on $8\frac{2}{1}\frac{2}{1}\frac{3}{1}$, we kill
$8\frac{2}{1}\frac{2}{1}\frac{5}{1}$, which is dominated by Viro's
$M$-scheme $8\frac{3}{1}\frac{3}{1}\frac{5}{1}$. Hence Scenario A
with a TOR on $9\frac{2}{1}\frac{2}{1}\frac{2}{1}$ is more likely
to be correct, and by periodicity we have a TOR on
$5\frac{2}{1}\frac{2}{1}\frac{6}{1}$, and then the sawtooth
climbing up to $4\frac{2}{1}\frac{2}{1}\frac{9}{1}$ somewhat
uselessly as the scheme is not derived from an $M$-scheme in view
of Viro's sporadic obstruction. Much more importantly, it should
be noted that Scenario A (somehow forced by Viro's construction of
$8\frac{3}{1}\frac{3}{1}\frac{5}{1}$) looks much incompatible with
Arnold's periodicity for dividing curves forcing them to live in
the 2-by-2 lattice centered at $M$-schemes and propagating
downwards through RKM-schemes. So Scenario B is much forced by
Arnold's will (congruence), and therefore the censorship principle
looks once more jeopardized, except if one of Viro's construction
is foiled due to over-liberal patchworking. Of course another
option is that the principle of contraction is false so that
$8\frac{2}{1}\frac{2}{1}\frac{5}{1}$ does not exist, yet we would
then still have $8\frac{3}{1}\frac{3}{1}\frac{5}{1}$ in the
censorship cone above $8\frac{2}{1}\frac{2}{1}\frac{3}{1}$.

We emphasize once more that by virtue of Arnold's periodicity for
dividing curves it is quite likely that there is a phenomenon of
total reality at the (non RKM) scheme $\frac{6}{1}\frac{12}{1}$
that would sustain Orevkov's obstruction of
$1\frac{6}{1}\frac{13}{1}$ and simultaneously kill the boson
$1\frac{7}{1}\frac{12}{1}$. Likewise it is likely that a TOR
phenomenon is active on the scheme $\frac{3}{1}\frac{15}{1}$.

When we arrive to the 4thlayer it seems natural (to explain all of
Viro's prohibitions) to put a TOR on
$\frac{4}{1}\frac{4}{1}\frac{7}{1}$, and by periodicity likewise
on $4\frac{4}{1}\frac{4}{1}\frac{3}{1}$. However the normalized
representative of this scheme lives in the 3rd layer under the
guise $4\frac{3}{1}\frac{4}{1}\frac{4}{1}$ and there would kill a
scheme by Polotovskii subsumed to Viro's $M$-scheme
$4\frac{3}{1}\frac{5}{1}\frac{7}{1}$.

At this stage it seems also of some relevance to see what Viro's
method produce at the $(M-3)$-level. However a look at
Fig.\,\ref{ViroDEGREE8_exotic_patches0_SYS2:fig} shows that
configuration with only one macro-circuit occurs only when
combining patches one of whose member is empty, so that no
$(M-3)$-scheme can be created out of Viro's maximal dissipations.

For instance it is puzzling whether the $(M-1)$-scheme exist. In
principle it is via a contraction of an empty oval on Viro's
$M$-scheme $1\frac{2}{1}\frac{17}{1}$. We could expect getting it
by a Viro patchwork with 3 macro-circuit, yet referring back once
more to Fig.\,\ref{ViroDEGREE8_exotic_patches0_SYS2:fig} we see
that Gods leave no room for such a construction as all those
tri-macro-oval combination involve a foolish (faul) empty patch
collection.

Now as to the 1st layer it could be posited that the bosonic
obstructions are caused by a mushroom of total reality at level
$(M-4)$. The most plausible for this to occur would be to place a
TOR on $1\frac{1}{1}\frac{14}{1}$, yet destroying thereby also
Viro's $M$-scheme $5\frac{1}{1}\frac{14}{1}$, however this need
perhaps just to amputate suitably Viro's list of patches. Of
course it is more politically correct (for Viro) to place the TOR
only at $\frac{1}{1}\frac{17}{1}$, yet doing so would also ill
Viro's scheme $1\frac{2}{1}\frac{17}{1}$, except if we imagine
that $\frac{1}{1}\frac{17}{1}$ causes only a degenerate censorship
cone, killing only what is lying immediately above it (in our
depiction mode which contains admittedly some arbitrariness). In
this scenario of course we would have counterexample to the
Viro-Itenberg contraction principle, or speak heuristically Viro's
$M$-scheme would be an instable clef-de-voute of the cathedral,
i.e. sitting firmly but only supported through a thin collection
of stones.

\begin{figure}[h]\Figskip
%\vskip-1.2cm\penalty0
%\centering
\hskip-2.7cm\penalty0
\epsfig{figure=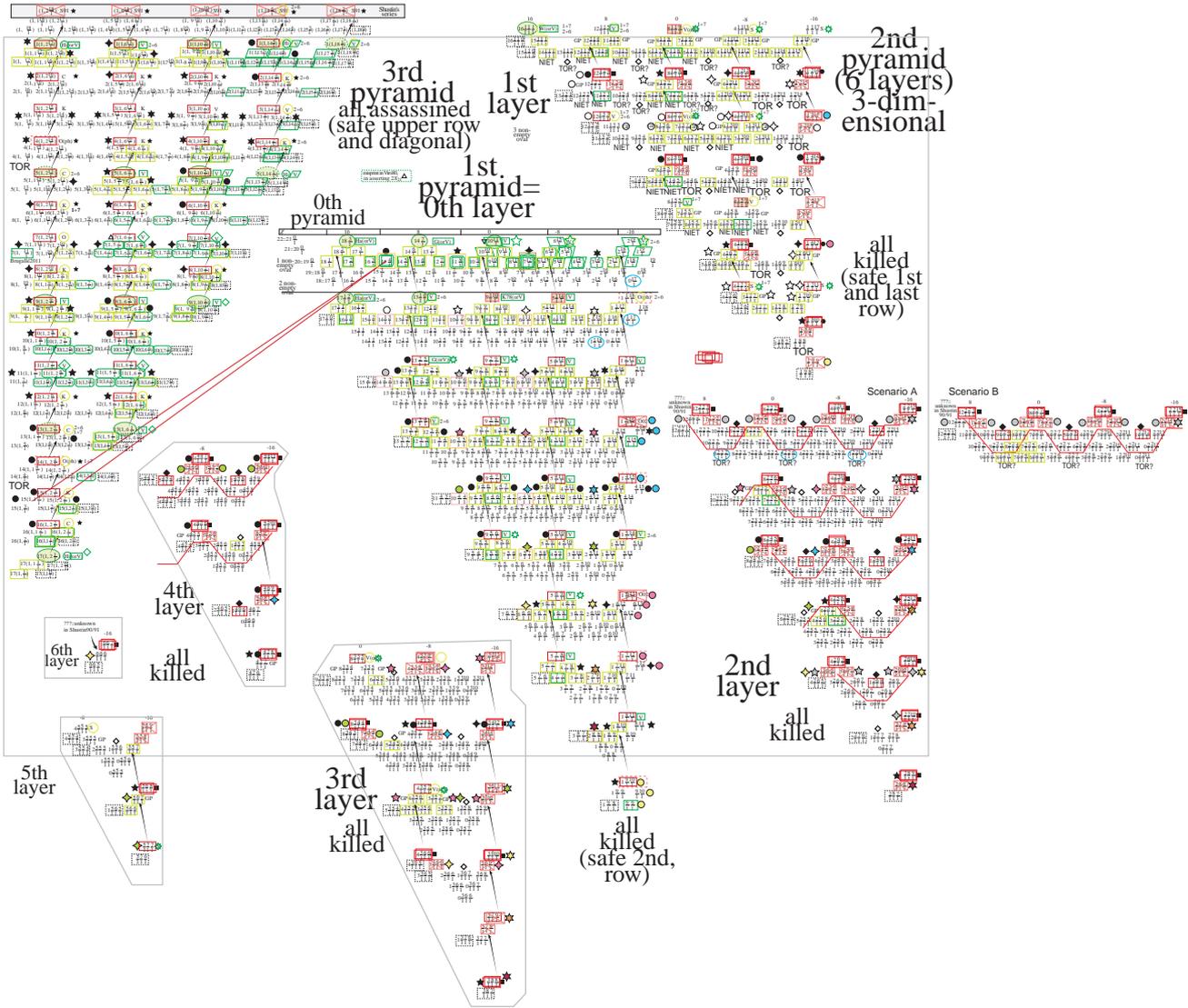,width=172mm}
\captionskipAG
  \caption{\label{Degree8-(M-i)-curve-TABLERohlinTRIS}%
  Trying to interpret obstructions via censorship above
  mushrooms situated at the
  $(M-4)$-level}
\figskip
\end{figure}

As we saw one important question (in order to refute Rohlin's
maximality conjecture) is the constructibility of the uni-nested
RKM-schemes. are not readily accessed through asymmetric gluings
of $M$-dissipation ($M$=maximal number of nine micro-ovals) it is
reasonable to expect that Viro's method yields all of them when
using non-maximal dissipation where the patches collection is
anymore severely restricted by Gudkov periodicity. Nonetheless
Viro's $M$-patches are sufficiently restricted to fail getting all
uni-nested $M$-schemes (only $18\frac{3}{1}$, $14\frac{7}{1}$ and
$10\frac{11}{1}$ being realized by Viro's purest method of the
quadri-ellipse). Therefore it is not completely clear that all
uni-nested RKM-schemes should succumb to Viro's quadri-ellipse
(VQE). It seems yet that the first three such schemes
($15\frac{4}{1}$, $11\frac{8}{1}$ and $7\frac{12}{1}$) are
coverable by VQE modulo sloppy extrapolation of Viro's
$M$-parameters (compare Fig.\,\ref{ViroDEGREE8_2:fig}f).

In contrast it is perhaps even more obscure if VQE can reach the
scheme $3\frac{16}{1}$. Of course, this is quite likely if Viro's
$M$-parameters can attain the $M$-schemes $6\frac{15}{1}$ and
$2\frac{19}{1}$. As shown by the composition table with extended
parameters (Fig.\,\ref{ViroDEGREE8_extended:fig}) $2\frac{19}{1}$
can be obtained as $(9,0,0)\times(9,0,0)$, and $6\frac{15}{1}$ is
obtained as $(5,4,0)\times (9,0,0)$. Moreover we remind (see
Fig.\,\ref{ViroDEGREE8_exotic_patches0_SYS:fig}) that to the best
of our knowledge the patch $V2(9,0,0)$ is still non-obstructed so
that there is a little chance that Viro's method (especially VQE)
has more swing than presently known. Remind also that this patch
$V2(9,0,0)$ implies directly the materialization of the boson $b1$
and $b7$ (see the cited composition table).

Thus with the patch (9,0,0) we may reach the extreme right of the
1st pyramid and thus for damped parameters we may also reach the
yet unrealized scheme $\frac{19}{1}$. However the picture of the
patchwork V2/V2 shows that there is at least one outer ovals so
that the scheme $\frac{19}{1}$ remains
%%% unattacked
unattained by Viro's method albeit diagrammatically dominated by
Viro's $M$-scheme $2\frac{19}{1}$. Hence speculating about an
intrinsic obstruction against existence of $\frac{19}{1}$
(Petrovskii's bound fails doing this), then this scheme could be a
possible  counterexample to the Viro-Itenberg contraction
principle of empty ovals.

[01.09.13] It is maybe worth noting that the purely nested
uni-nested $(M-2)$-scheme $\frac{19}{1}$ admits many
$(M-1)$-extensions below Shustin series (subnested case without
outer ovals). Alas $\frac{19}{1}$ has wrong Euler-Ragsdale
characteristic for being of type~I, and also it is for all those
sub-Shustinian scheme to be prohibited since half of them are
deducible as shadows of the subnested $M$-scheme with one outer
oval constructed by either Hilbert or Viro.

Of course, we can imagine that Shustin's $M$-obstruction are
interpretable as a censorship to a total reality reality
concentrated on the $(M-2)$-schemes below it (e.g.,
$(1,1\frac{17}{1})$, $(1,5\frac{13}{1})$, etc.), but then both
$(M-1)$-schemes below Shustin's scheme (which admits only two ways
of contraction in view of the absence of outer ovals) are killed,
and so the contraction principle is jeopardized again. As a
conclusion, it seems that there is no way to explain Shustin's
$M$-prohibitions via censorship (at least in a way respecting the
Viro-Itenberg contraction principle).

Otherwise we can imagine a deep mushroom of total reality reigning
for the Arnold-congruent $\chi\equiv_4 k^2$ scheme $\frac{17}{1}$
lying at level $(M-4)$. Of course this deep mushroom would cause
serious damages in the uni-nested pyramid (more than Viro's beaver
construction permits, yet this is perhaps erroneous), and
collaterally it would explain Shustin's obstructions.

The question about $\frac{19}{1}$ is rather puzzling. As an
attempt to get it via Viro's method (with extended parameters) we
may first take a look on the gluing-table
(=Fig.\,\ref{ViroDEGREE8_exotic_patches0_SYS2:fig}). As the scheme
as $r=M-2$ many ovals we confine our attention (using Viro's
maximal dissipation) to the gluings with 2 macro-ovals
(green-colored on the gluing table). The amendment of pure nesting
inherent in $\frac{19}{1}$ forces limiting attention to the cases
where the 2 macro-ovals are nested, and a short inspection of the
table leaves leaves only the few possibilities: B/G, C/G, D/F
(alas this is a no-mans-land because F is empty), D/H (also
empty), D/J (cannot reach sufficiently many ovals because J is
saturated), E/G, F/H (empty), F/J (too small because J is
saturated), G/I, H/J (not good as H is empty and J is saturated).
So it remains a good deal of candidates. Additionally, one must
not forget the gluings with flipped version in the margin of
Fig.\,\ref{ViroDEGREE8_exotic_patches0_SYS2:fig} (specifically
F/Fflip, F/Hflip, H/Fflip symmetric, and H/Hflip), but all those
of interest involve an empty family of patches (specifically F or
H) and so those gluings can  ultimately be discarded.

First for B/G, we see that we may specialize to B2/G1 and to
control the topology we choose $(\al,\be,\ga)=(0,9,0)$ at B2 and
$(\al^\ast,\be^\ast,\ga^\ast)=(0,9,0)$ above at G1. Controlling in
the table of patches
(Fig.\,\ref{ViroDEGREE8_exotic_patches0_SYS:fig}) we see that
B2(0,9,0) is unobstructed, but alas G1(0,9,0) is not adjusted to
Gudkov periodicity (hence an illegal patch).

It remains to browse the other options. The next opportunity is
C/G, but again whatsoever we choose below (C1 or C2, where $\ga$
is just transplanted in the more internal lune), we are forced to
select the same upper patch as before G1(0,9,0), which is an
illegal move.

Our next chance involves E/G, but again the prescribed
target-topology of $\frac{19}{1}$ forces us to employ the illegal
patch G1(0,9,0), while on the bottom we use the potential patch
E(0,9,0) not given by Viro's theory (yet not prohibited too). Of
course this no consolation and there is no opportunity to get the
requested scheme in this flavor E/G.

The next (and last) opportunity is G/I, and once more it is clear
that we are forced to appeal to the same illegal patch G1(0,9,0)to
get the desired curve.

We have proven the following:

\begin{lemma}
Even with extended parameters Viro's method of  the quadri-ellipse
will never succeed in constructing the extreme right uninested
scheme $\frac{19}{1}$ at least by using $M$-dissipations. However
it remains the hope that it is accessed via non-maximal
dissipation.
\end{lemma}

\begin{figure}[h]\Figskip
%\vskip-1.2cm\penalty0
%\centering
\hskip-2.7cm\penalty0
\epsfig{figure=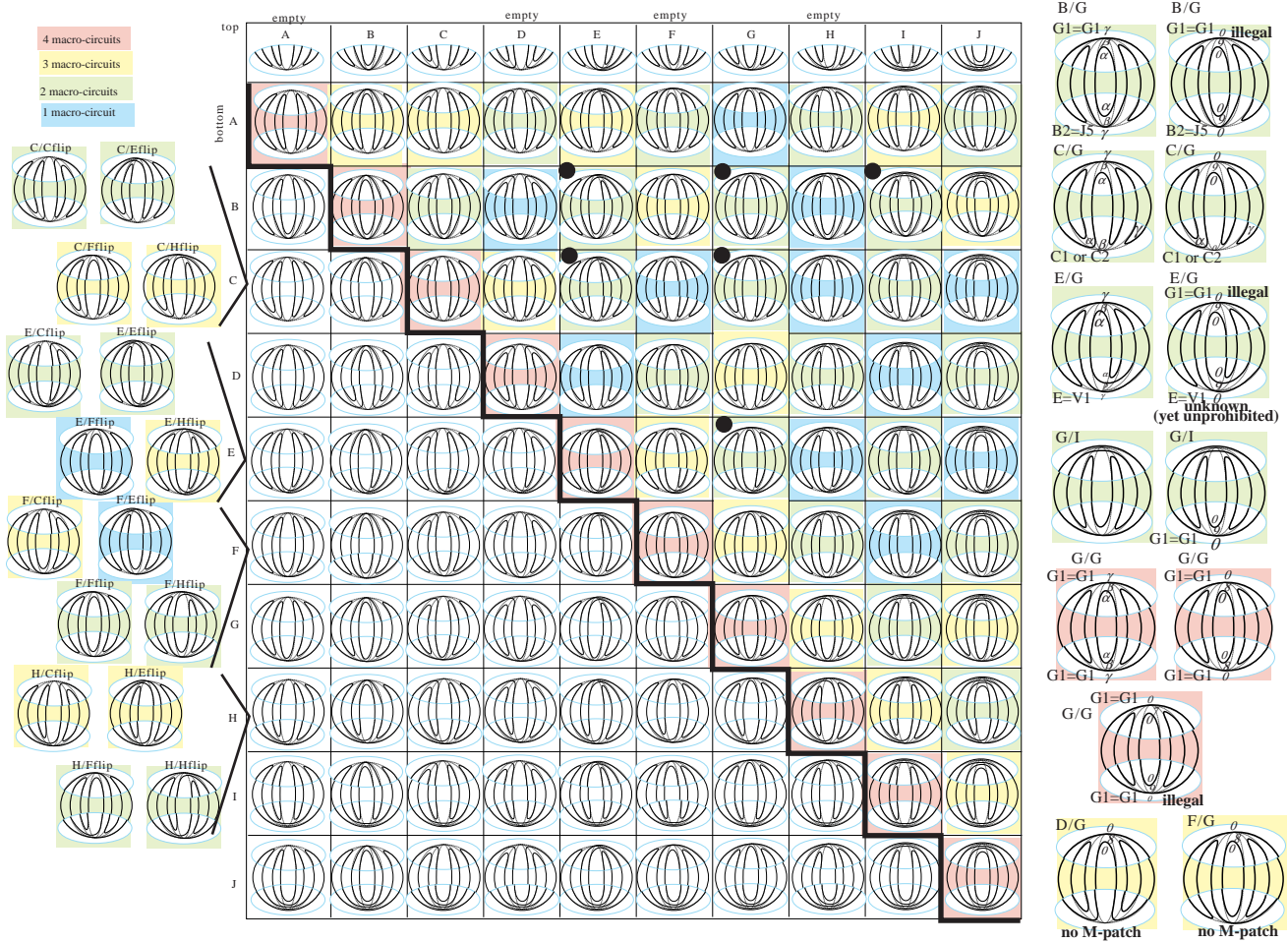,width=172mm}
\captionskipAG
  \caption{\label{ViroDEGREE8_exotic_patches0_SYS2_19-1}%
  Trying to get $\frac{19}{1}$ (via $M$-patches)}
\figskip
\end{figure}

It remains to analyze if it is possible to realize the scheme by
using non-maximal dissipation. Again we look at the ground-shapes
of gluings, and we may first consider those with 4 macro-circuits.
The absence of outer ovals and presence of a single nest restrict
our attention to the sole option of G/G. Then our sole chance  to
realize the given scheme is via G1(0,8,0)$\times$G1(0,8,0).

Of course the next reflex is to compose this patch with those of
Viro (C,E), but alas the morphogenetic table shows that those
patchworks (C/G and E/G) have only 2 macro-ovals, so that the
resulting schemes have $M-3$ ovals.

A further possibility is to exploit the ground curves with 3
macro-ovals. We can then focus attention to the yellow cases of
the table, especially the combinations D/G (like David Gauld!),
F/G and G/H. First, as to D/G we can only take G=G1(0,8,0)
(granting existence of this patch) but then we are forced to take
an $M$-patch for D and those are known to be prohibited (Arnold's
weak version of Gudkov periodicity). For F/G we can tell the same
story, and likewise for G/H. This proves the following:

\begin{lemma}
There is one and only one chance to realize the extreme-right
uni-nested  scheme $\frac{19}{1}$ via Viro's method of the small
perturbation of a quadri-ellipse, und zwar als
G1(0,8,0)$\times$G1(0,8,0). Whether this opportunity is actually
realized is unknown to us. (If it is not we get maybe a
counter-example to the contraction conjecture of Viro-Itenberg.)
\end{lemma}

\subsection{Patchwork vs. Artwork}

[09.08.13] Next we see that the defect of the patch G2 is that to
glue it maximally with itself we necessarily create doubloons of
$\be$ and $\ga$ in the same lunes. if take G3 the avatar of G2 yet
with bifolium at the center then we will dispose of one more way
to self glue the patch with itself by adding a vertical symmetry.
So we will also be able to exploit Viro's sporadic law and more
patching obstacles will be gained than for G2.

Initially we wanted to compose as shown on Fig.\,pre G3/G3 where
$\al$, $\be$, $\ga$ appears along the ``natural'' left-right sense
of reading. Yet, as a trick, it seems advisable to permute the
parameters in order that G3/G3 involves the same table as G2/G2,
and so we changed the labelling of ovals as $\be$, $\al$ and $\ga$
when read from left-to-right (compare Fig.\,G3/G3). Hence as
before the same patches as in the G2/G2-table are killed. Now the
quick is that G3 can be composed with itself plus an additional
vertical symmetry so that both parameters are not anymore in the
same lune creating thereby evenness of the content of an oval.
Alas, we realize only now that even before (e.g. for G2), it is
actually $\be$ and $\be^{\ast}$ which can vary independently so
that no parity of the content is forced. The real reason is
actually that $\be$ is fixed to 1 mod 4 and therefore all schemes
of the  table have even content on the central nest. Yet, it is
precisely this consanguinity that we can now jeopardize by
symmetrizing the patch G3. This gives us the table G3/G3B where B
stands for bis. Actually as most patches are already prohibited we
do not even need to feel carefully that table, but can restrict at
the nodes of where there are quantum jump which are the white
patches in suspense (doubt). Then we see that Viro's sporadic
obstructions (it is enough to consider them along the diagonal)
kills the $G3$-patches $(4,5,0)$ and $(0,9,0)$, but however
$(8,1,0)$ is not killed since it creates a Shustin's scheme
accessible via the medusa. Even if modest, this would be an
interesting result, for it would provide a more elementary
construction of Shustin's scheme. At any rate, we see that if the
three patches not forbidden by Fiedler-Viro's oddity law are
available then we could refute some of Viro's sporadic
obstructions, and we would so-to-speak raise the death (i.e.
schemes prematurely killed by Viro). So:

\begin{Scholium}
Could it be that Viro's presentation of the dissipation theory of
$X_{21}=F4$---(in Arnold vs. Gabard's notation), i.e. quadruple
flat point with 4 branches having a 2nd order tangency (with all
branches curved in the same half-plane (by the way it would be
interesting to study the other cases sembling chromosomes, where
some of the 4 branches lies in different half-planes, as those
could admit realization in degree 6 as well at least from the
naive real viewpoint)---is too rigid or frigid?? That is to say
did Viro missed certain patches? And if yes this can be done in
varied level of dequantization sometimes by attacking Viro's
sporadic obstructions. In contrast if all of Viro's sporadic
obstruction are correct (actually two of them suffices, namely
prohibition of $\frac{5}{1}\frac{5}{1}\frac{9}{1}$ and
$\frac{1}{1}\frac{9}{1}\frac{9}{1}$), then the patch G3 still
admits the dissipation $(8,1,0)$, which is not listed on Viro's
table (=Figure~55 on p.\,1118 of Viro 89
\cite{Viro_1989/90-Construction}).
\end{Scholium}

Again at this stage it seems worth reading once more Viro's remark
(p.\,1119 of \loccit):

``We do not yet have a complete topological classification of the
dissipations of $X_{21}$ singularities. Shustin [32] proved that
all dissipations of type $X_{21}$ singularities with a given
number of real branches have the same topological type; however
there is still a big gap between what is given by the
constructions and the prohibitions. Curiously, the problem has
been completely solved for dissipations that can occur in the
construction of nonsingular $M$-curves. These dissipations are
considered in the next theorem. It can be shown that any
dissipation of an $X_{21}$ singularity with four real branches in
the course of which nine new small ovals appear (this is the
maximum possible number) is topologically equivalent to one of the
dissipation in Theorem~4.5.A.'' [{\it Gabard'a addition}:
Incidentally it seems that there is a misprint here and that one
should read Theorem~4.7.A.]

GABARD's comment on this prose. Personally, we do not see how to
prove this and believe that the ``it can be shown'' of Viro is
hazardous or presumptuous/arrogant (as we are unable to reprove
it). As we said the point is that even when admitting Viro's
sporadic obstructions it does not result a prohibition of the
patch G3(8,1,0). Question: how can Viro rules out this patch?

Next it is evident that this sort of argument via composition
table can be multiplied to other exotic patches (exotic in the
sense that Viro claims their nonexistence).

So for instance the patch G4 when self-glued with itself (with or
without an additional vertical flip) will produces a scheme with
invariably $\ga+\ga^\ast$ outer ovals which will therefore be even
when evaluated along the diagonal of the composition table.
However Gudkov periodicity tells us this number of outer ovals
being 1 modulo 4 in the doubly-nested case (compare e.g. the
pyramid). Hence all the G4-patches are prohibited merely by Gudkov
periodicity (even in the weak version of Arnold).

The patch G5 will conflict with Shustin's restrictions, yet
probably for only few values of the parameters. Hence, we do not
see how Viro is able to prohibit all those patches. On gluing G5
with itself, we note that the number of big eggs (i.e. craw's eggs
at depth 1 is $1+\be+\be^\ast$ hence always odd along the diagonal
of the composition table). However Gudkov periodicity imposes this
number as being congruent to 2 mod 4. Actually Arnold's weaker
version would actually suffices for our purpose of inferring
therefore, that all patches of the family G5 are killed. So in
principle it not even worth tabulating the G5/G5, yet let us do it
just for the fun of being a ``angst Hase''. So let us choose $\be
=0$ as initial value, albeit this is not satisfactory when
$\be^\ast=0$ too, and so it is really not worth doing any
tabulation.

Next we have a menagerie of patches involving 4 unnested lunes
(the G6-family on Fig.\,\ref{ViroDEGREE8_exotic_patches0:fig}),
and some similar thinking based Gudkov periodicity (often in the
simple variant of Arnold) suffices to kill all those patches.
Indeed the most universal G6-type involves parameters $\al,
\be,\ga,\de,\ep$. If all first four of those are positive then
B\'ezout for conics is violated via the saturation of doubled
quadrifolium. So we may assume one of them zero. Without much loss
of generality, assume $\de=0$ (compare picture right below), then
when gluing with a symmetric replica we land in the trinested
realm where the number of outer ovals is 0 mod 4 and not
$1+\ep+\ep^{\ast}$ which is odd when the patch is self-composed
with itself, hence violating even Arnold's version of Gudkov. So
we may assume one more parameter zero and this yields e.g. the
picture just right below the former one. But again the same
argument shows that it is anti-Gudkov in the version of Arnold,
since in the doubly nested case there is one outer oval mod 4 (or
mod 2 suffice). Finally we must analyze variant were there is
additionally a micro-nest (see G7-family on the figure). Of course
if there is too much parameters $\al, \be$ positive, we can
corrupt the bound of 4 nests, and so we are led to patches with
only one of the 4 ``semi-lune'' inhabited, yet even in this case
an additional vertical flip of the upper (replicated) patch will
corrupt the saturation of the doubled quadrifolium, and so
independently of $\al$'s position which cannot be central!
Accordingly the sole issue is to have all semi-lunes empty, but
then the resulting curve has the scheme $4\frac{8}{1}\frac{8}{1}$,
with number of outer ovals improper for Gudkov (even Arnold's)
periodicity. So we have more or less proved:

\begin{lemma}
Arnold's periodicity implies that there is no dissipation of
$X_{21}$ such that the perturbed germ exhibits 4 unnested
semi-lunes (compare G6-family on
Fig.\,\ref{ViroDEGREE8_exotic_patches0:fig}).
\end{lemma}

Philosophically it is quite remarkable that all those patches are
killed albeit they look decent traffic transition from the real
viewpoint (of traffic circulation). Yet, the reason behind this
rigidity is Arnold's deep insights on the complexification and
4-manifold. Therefore those prohibition reduce essentially to the
basic algebraic fact that an even unimodular integral quadratic
form has signature divisible by 8 (due to old arithmeticians like
Zolotarev, etc.), while Rohlin's deepens this by a factor 2 if the
form arise as the intersection form of a smooth spin
$4$-manifolds.

Yet, in contrast to all this our patch G1 is nearly unprohibited
even under Shustin's obstructions. So how can Viro claim that his
list of patches is exhaustive!?? Maybe we should try to compose G1
with other patches and exploit other obstructions of Viro relative
to $(M-2)$-curves. We recall therefore Viro's corresponding
statement (Theorem~\ref{Viro-Fiedler-prohibition:thm}) but which
we reproduce now for convenience, and of which we shall try to
exploit the 2nd clause.

\begin{theorem}
{\rm (Viro 1983 \cite{Viro_1983/84-new-prohibitions})}
\label{Viro-Fiedler-prohibition-COPY:thm}

$\bullet$ ($M$)---If $
\frac{\alpha}{1}\frac{\beta}{1}\frac{\gamma}{1} \delta$ is the
real scheme of an $M$-curve of degree $8$ with $\alpha, \beta$ and
$\gamma$ nonzero, then $\alpha, \beta$ and $\gamma$ are odd.

$\bullet$ ($M-2$)---If $
\frac{\alpha}{1}\frac{\beta}{1}\frac{\gamma}{1} \delta$ is the
real scheme of an $(M-2)$-curve of degree $8$ with $\alpha, \beta$
and $\gamma$ nonzero and with $\alpha+\beta+\gamma\equiv 0 \pmod
4$, then two of the numbers  $\alpha, \beta, \gamma$ are odd and
one is even.
\end{theorem}

(Remind that albeit a bit undigest, this 2nd clause has a
clear-cut diagrammatic impact, compare
Fig.\,\ref{Degree8-(M-i)-curve-TABLE:fig}).

The idea is that if  we compose G1 with V1 (or some other patch)
then we hope to land in the trinested realm where there are
serious prohibitions. Actually this is not readily the case, yet
we can try first V1/V2 etc, or V2/V2fipped, etc, and hope thereby
to get the restriction that we were as yet not able to explain by
looking merely at $M$-curves. After drawing some combination of
patches we realize that combination V2/V3flip yields trinested
(M-2)-schemes.

[10.08.13] However upon working out the table it is clear that we
never interact with Viro's $(M-2)$-obstruction. Indeed the latter
concerns schemes with number of outer ovals equal to $1 \pmod 4$,
whereas those just created have this number equal to $3 \pmod 4$.
Therefore, we are a bit disappointed to get no new obstructions
but perhaps should not as it is precisely our secret dream that
Viro's dissipation theory was over-atrophied, i.e. handicapped by
too many restrictions, exacerbating
%%%and therefore exaggerating
thereby the difficulty of $Hi(8)$, Hilbert's in degree 8.

\begin{figure}[h]\Figskip
%\vskip-1.2cm\penalty0
%\centering
\hskip-2.7cm\penalty0
\epsfig{figure=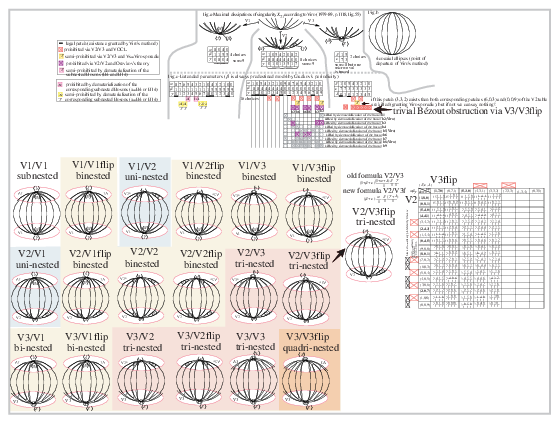,width=172mm} \captionskipAG
  \caption{\label{ViroDEGREE8_V1-V2:fig}%
  Nonmaximal $(M-2)$-gluings of V1/V2, etc. yet no direct patchwork
  obstruction are gained from Viro's $(M-2)$-law.}
\figskip
\end{figure}

What to do next? After some sleep, we realized that one must
always think in a structured fashion, and so arranged the
super-table of composition as an array with double entry. Doing
this we discovered the type V3/V2flip also leading to the
trinested realm. Yet it seems clear that this is just the earlier
one rotated. Finally, working systematically we find the
combination V3/V3flip which is even quadri-nested, hence B\'ezout
(or saturation allied to total reality) implies that $\lam=0$, and
so:

\begin{lemma} The dissipations
of type $V3$ really reduces to Viro's list of two.
\end{lemma}

Alas, it seems then that we have exhausted all combinations,
because V3/V2 is isotopic to V2/V3 via a reflection (recall
optionally that the mapping class group of $\RR P^2$ is trivial).
In conclusion, we still do not understand why some of Viro's
patches (especially V1 and V2) list are so restricted.

[11.08.13] So we have to continue this boring game of tabulation.
For instance we tabulated the G3/V3 and G3/V3flip combination on
Fig.\,\ref{ViroDEGREE8_exotic_patches2:fig},  yet it resulted
schemes with the wrong number of outer ovals.

Next we analyzed G7/G7 and it seems clear that there is an evident
Gudkov-Arnold obstruction killing the full G7-family. Indeed when
the patch G7 is glued with itself (modulo horizontal reflection),
we get Fig.\,\ref{ViroDEGREE8_exotic_patchesB:fig}(G7/G7), where
we see $1+2\be$ big eggs thereby violating the Gudkov-Arnold
periodicity since the scheme belongs to the subnested case. The
only way to resolve the contradiction is to impose $\al=0 $, but
then the scheme is simply-nested with an odd number $1+2\ga$ of
outer ovals, jeopardizing again  Gudkov-Arnold periodicity.

Eventually one can expect that some combination yields
$(M-1)$-curves where there is the GKK-prohibition, but as yet we
could never arrange this.

Next, we worked out more carefully the G3/G3flip-table of
Fig.\,\ref{ViroDEGREE8_exotic_patchesB:fig}, but this is a bit
anecdotic because we can surely rules out two more patches if we
believe in Viro's sporadic obstruction, but still cannot exclude
the G3-patch $(8,1,0)$.

Eventually, we imagined also the patch G8 (of
Fig.\,\ref{ViroDEGREE8_exotic_patches0:fig}). On gluing G8 with
itself under a horizontal symmetry we get the combination G8/G8
and a corresponding table.

\begin{figure}[h]\Figskip
%\vskip-1.2cm\penalty0
%\centering
\hskip-2.7cm\penalty0
\epsfig{figure=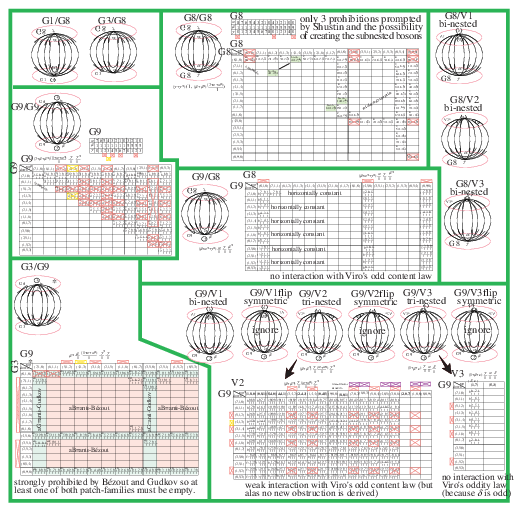,width=172mm}
\captionskipAG
  \caption{\label{ViroDEGREE8_exotic_patchesBIS:fig}%
  Patchworking exotic patches (continued)}
\figskip
\end{figure}

By Gudkov periodicity the parameter $\be$ can be chosen as $\be=1$
plus its companions mod 4, i.e. 5, and 9. It may be observed that
the general formula of the Gudkov symbol differs only slightly
from that arising via G1/G1, yet as we choose the parameters
differently it seems that we must restart anew the whole
calculation. As the evolution rule is so simple it is a simple
matter to cut-and-paste the symbols while adjusting a bit the
geometry. Once the table dressed, we see that G8/G8 sweeps out
slightly less schemes than G1/G1, but again as Shustin's
prohibition concerns only the top of the 3rd pyramid only 3 direct
obstructions are gained , namely $(8,1,0)$, $(4,5,0)$ and
$(0,9,0)$. Again it is not clear how Viro can prohibit the
majority of  those surviving patches. In case of a realization of
them, the net bonus would be the creation of new bosons, as well
as trivialization of many schemes due to Viro (ad-hoc's
construction with horses and beavers), Korchagin, Chevallier, and
Orevkov. The philosophy would be that the simplest animal (namely
the protozoan of the quadri-ellipse), incarnating something like a
primitive biological shape is sufficient to access most curves of
the world. Hence, the embryology of the most complex real octic
could be readily traced back to its most primitive ancestor by a
direct morphogenesis (dissipation).

Next, we can glue G1/G8 or G3/G8, yet the slight defect of this
scheme is that it is only doubly-nested. So we imagined a variant
of G3 called G9 where we introduced $\be$ subsequently as to make
possible Gudkov periodicity upon choosing $\be=1$ and its
companion mod 4. While filling this table we wondered if did not
made basic mistakes in the earlier tables??? Further a simple
trick is to discover the evolution rule along the diagonal without
having to fill the full table. At any rate, by the usual method we
can kill via Viro's law 4 patches, and one more is killed by
Viro's sporadic obstruction.

For instance we decided to check the table of G2/G2, or maybe even
V2/V2. This is a boring task and we reserve it for later as we
need an optical pause.

[12.08.13] Next, it remains to compose G3 with G9 which leads to
an $M$-curve since both patches have the same underlying graphics.
The resulting patchwork G3/G9 looks quadrinested so there should
be severe obstructions. More precisely, the gluing G3/G9 is
quadri-nested as soon as $\be, \ga, \ga^{\ast}$ are positive. On
working out the table we remark that all schemes which are not
frankly anti-B\'ezout are anti-Gudkov in the fine form of the
latter not covered by Arnold's weaker version thereof. Therefore
it must be recognized that all schemes of the G3/G9-table are
nonexistent, and therefore that at least one of dissipation mode
G3 or G9 is empty (i.e. contains not a single patch). Okay but can
we really conclude that both types G3 and G9 are empty as Viro
seems to claim?

From the table G3/G3 (which is the same as that of G2/G2 on
Fig.\,\ref{ViroDEGREE8_exotic_patchesA:fig}) we know that only 3
patches at most can survive to the Fiedler-Viro oddity law, namely
$(8,1,0)$, $(4,5,0)$, and $(0,9,0)$. Granting moreover Orevkov's
obstruction(=desintegration) of the (pseudo)-boson
$1\frac{6}{1}\frac{13}{1}$ we see that that both first patch
cannot coexist simultaneously. Which one is not clear a priori
albeit the patch $(8,1,0)$ looks more down-to-earth as its
self-gluing yields a simplest Viro's scheme, namely
$1\frac{2}{1}\frac{17}{1}$. So accepting Orevkov as true, three
scenarios are possible, according to the exhaustive list of
G3-patches:

(1) $G3$ contains $(8,1,0)$ and $(0,9,0)$ (and nothing more); [in
this case the bosons $b_1$ and $b_9$ are materialized]

(2) $G3$ contains $(4,5,0)$ and $(0,9,0)$ (and nothing else); [in
this case the same bosons $b_1$ and $b_9$ are materialized]

(3) $G3$ contains only $(8,1,0)$; [in this case no boson is
materialized]

(4) $G3$ contains only $(4,5,0)$; [in this case the bosons $b_9$
is materialized]

(5) $G3$ contains only $(0,9,0)$; [in which case the bosons $b_1$
is materialized]

(6) $G3$ is the empty  set [no new boson is materialized].

This looks really meagre and sloppy piece of knowledge, and one
would like to know more! A priori we could hope to draw sharper
information by the interaction G3/G9. For instance if we knew that
the G9-family of patches is nonempty, we could deduce that G3 is
empty, and we would be in the arid scenario (6). Alas, it is not
clear which of G3 and G9 is more likely to be empty. From a naive
probabilistic perspective, the diagonal of G9/G9 contains fairly
many realized schemes due to Gudkov (once), Viro (thrice), Shustin
(thrice). In contrast the G3/G3-table hits more systematically
against the Fiedler-Viro obstruction so that one could suspect
that emptiness of G3 is more likely. Further, G3/G3 can create
bosons whilst G9/G9 only creates standard schemes (at least when
not violating Viro's highly sporadic obstruction). Accordingly, it
could be that emptiness of G3 is more likely than that of G3, of
course not impeding an extinction of both species (as seems
implicit in Viro), yet we know no theoretical justification for
this collapse.

However remind from the table G3/G3flip that it suffices to
believe in Viro's sporadic obstructions to kill the patches
$(4,5,0)$ and $(0,9,0)$, so that in principle only scenarios (3)
and (6) are logically possible. Alas, we do not know how to
produce more information, as one of the pity is that since G9 is
symmetric we cannot infer more information by composing the patch
with itself flipped.

Of course the general feature to be remembered from our
experiments with V2/V3 and G3/G9 is that both patches are coupled
with one differing from the other merely by the addition of an
outer island (and slight restriction in the parameters). Each pair
of such partner patches will create an additional table of
$M$-schemes from which one may infer additional
patching-prohibitions. Alas, one does not receive complete
information by this method, at least not as complete as tacitely
claimed by Viro.

One of the other problem is that as yet we could not exploit the
Viro's prohibitions on $(M-2)$-schemes of RKM-type in the
trinested realm. Those will perhaps enter into action when
composing G9/G8, which by construction will land in the
3ply-nested zone. From the scratch, i.e. by filling the first
entry of the table and comparing with the main pyramid
(Fig.\,\ref{Degree8-(M-i)-curve-TABLE:fig}) we see that we land in
the RKM-domain where $\chi\equiv k^2+4 \pmod 8$, so there is a
good chance to meet Viro's $(M-2)$-obstruction. Yet on the same
moment since the table G9/G8 is not a self correspondence it will
only results hypothetical patches-obstruction unless we know a
specific construction in one of both patch family. So it is not
clear that the device will lead to any tangible piece of
information. Yet having no better idea in the pockets let us
peacefully fill the complete G9/G8-table. It is quickly observed
that the table is horizontally constant, and during the tabulation
we see that the vertical rows do the usual palindromic motions in
the pyramid yet {\it without\/} interacting with Viro's
$(M-2)$-obstruction.

Now what about composing G9 with V1, or V2, etc. Of course G9/V1
does {\it not\/} lead to the 3-ply nested realm, but so does G9/V2
and G9/V3. On tabulating the patch-works, we expect so to infer
valuable information from G9/V2 and G9/V3. On filling G9/V2 the
very first item is actually the scheme $11\frac{8}{1}$ which is
RKM. Then we propagate the table, e.g. first along the vertical
row, where we observe the usual palindromic path. Then we copy and
adapt each raw along the evolution rule. A Viro obstruction
appears first in the 5th vertical row, but alas as the vertical
$V2$-entry $(3,4,2)$ is not existing we cannot infer damages on
the corresponding $G9$-entries. We hope that things will improve
in the next row, but alas there we meet no Viro obstruction. On
the next row (0,4,5), we have the 3rd numerator equal to 5 so
there is no Viro obstruction, as the latter amounts to forbid all
three numerators being even integers. Of course, prior to filling
the whole it is already evident that we will gain no (new)
obstruction along the way, because the schematic obstruction of
Viro always fails to land in bold-faced rows corresponding to
materialized patches. As a result we see

\begin{Scholium} Viro's oddity law for $(M-2)$-schemes of type RKM does
not implies more patches obstructions than those already gathered
from the oddity law for $M$-schemes.
\end{Scholium}

It remains to hope that the situation is changed when using G3/V3
which despite naive expectation is only triply- (and not
quadruply) nested. As V3 is inherently restricted to two patches
by B\'ezout the tabulation-efforts are more user friendly. The
first entry of the table is reassuring (yielding an RKM-scheme).
However as $\de$ is odd (even 1 mod 4) it is clear that the
schemes so created never interact with Viro's obstruction and no
additional information on G9 is gained.

At this stage the methodology is fairly clear, nearly
algorithmizable as follows:

Step 1: Imagine all the patches that you want by tracing them as
if you where an artist like Saint-Exup\'ery's Petit Prince,

Step 2: combine all the patches between themselves (especially
with the V-family of Viro which is known to contain explicit
dissipations), and report the elementary obstruction of B\'ezout
and Gudkov, and optionally the highbrow obstruction of Fiedler,
Viro, Shustin, Orevkov. Keep also in mind the option of looking at
non-maximal (hybrid) smoothing which could potentially interact
with Viro's 2nd oddity law (not all 3 numerators of a trinested
RKM-scheme can be even). Further the important trick is to combine
a patch with its partner which has the same look plus an extra
island, so that it is a fake hybrid, leading thus to the
$M$-realm, where perhaps the most stringent obstruction are known.

Step 3: Self-criticism keep in mind the option that some of those
Russian obstructions are potentially erroneous, and that the
actual state-of-the-art on $H(8)$(=Hilbert in degree 8) is over
rigidified by too many prohibitions. Keep in minds that Russian
citizens prefer life with not excessive prohibitions, especially
when it comes to distilled vegetables.

Step 4: Notice, that apart from being a combinatorial looser, a
certain discrepancy between the obstruction so generated and those
claimed (semi-tacitely) by Viro 89 is quite big, leaving the topic
in an unsatisfactory state of affair.

To be concrete again it seems that we missed as yet to compose G8
with V1, etc. However it seems (see
Fig.\,\ref{ViroDEGREE8_exotic_patchesBIS:fig}) that all those
schemes G8/V1, G8/V2 and G8/V3 are never trinested and so we can
interact with Viro's $(M-2)$-prohibition. Of course flipping the
top patch does not help as the bottom one i.e. G8 is symmetrical.

\section{Viro's construction of the patches for $X_{21}$(=quadruple rainbow)}
\label{Viro:patches-construction:sec}

\subsection{The method employing hyperbolisms (Huyghens,
Newton,\dots, Viro)}

[13.08.13] After the energy spent in prohibiting patches it seems
advisable (and even logically requested in the prohibitive aspect)
to construct patches. We follow the geometric approach in Viro
1989/90 \cite{Viro_1989/90-Construction}, along a methodology
using Gudkov, Polotovskii, etc, and hyperbolism \`a la Huyghens,
Newton.

[30.08.13] Quite strangely Viro's exposition is fairly geometric
and not brute-force combinatorial patchwork. Also puzzling from
the entrance is the issue that the patch is something local but
constructed by excision out of a global object.

[01.09.13] Now we give the details. Viro starts with affine
quintics due to Harnack, Gudkov and Polotovskii and deduce some
quintics with special position w.r.t. the 3 coordinate axes. The
patch will be deduced via a hyperbolism (Huygens-Newton-Cremona).

It may be noted that Polotovskii's affine quintics yields when
smoothed (along with the line ``at infinity'') the RKM-schemes of
degree 6 studied by Rohlin and Le Touz\'e ($2\frac{6}{1}$ and
$6\frac{2}{1}$). Also pleasant is the issue of interpreting affine
curves as waves profiles, and then Polotovskii species correspond
to giant waves forming a rouleau, while the additional oval of the
quintic may be imagined as droplets of water resp. bubbles of air
injected in the sea profile. This metaphor looses
%some
of its
pertinence once it is remembered  that the pseudoline fails to
divide the plane $\RR P^2$. Hence, there is no distinguishable
mediums of water and air.

From all those 4 fundamental species of affine quintics, it is
argued that by dragging the line at $\infty$ to a tangent we get
six types of projective quintics
%%%flirting
behaving as depicted w.r.t. the fundamental triangle allied to
projective coordinates. At this step, one may wonder if there is
not also a configuration as shown below Polotovskii's where the 6
air bubbles are converted to droplet surfing below the rouleau of
the wave front. (One should perhaps keep this or other possibility
as an attempt to expand Viro's list of patches, with the possible
net bonus of getting the boson of $M$-octics not  yet known to be
realized, at least some of them.) Of course the classification of
affine $M$-quintics is in principle a closed chapter of geometry
certainly going back to Polotovskii himself. For simplicity, let
us leave aside this difficulty for the moment to continue Viro's
argument.

The next step involves applying a hyperbolism to all those six
quintics. This is nothing else than a Cremona transformation
involving the net/web of all conics through the 3 fundamental
points of the triangle.

[02.09.13] The resulting curve resembles a Bugs-Bunny (rabbit with
two big ears), and can be recognized as a margarita with 3 petals.
Actually, we fail to understand exactly why the resulting curve
has degree 8 as perhaps a hyperbolism is slightly different from
the standard Cremona transformation. When applied to Polotovskii's
curve we get a variant with both ears invaginated into the head of
the rabbit. Of course one may wonder if there is not a variant
where one ear is invaginated while the other lies outside (and
this naturally comes up when tracing the curve on an electronic
tracer like Illustrator, see Fig.\,a). It is also puzzling that
both Harnack's and Gudkov's curve creates two projective curves
while those of Polotovskii just a single one. So it is  a hard
duty to get convinced that Viro's exposition is exhaustive. For
instance we could wonder of what happens if the tangency is
arranged along the extremity of the wave in Polotovskii's curve
(Fig.\,b).

The next step involves dissipating the triple point while creating
a new little oval (as is evident from the theory of cubics) and
this in a somewhat erotical way so that one of the connecting
branch of the patch
%(i.e. reaching the boundary of the gluing
%disc)
performs a meander of 3 crossings along one of the fundamental
lines. This step we call a vibratory dissipation and in substance
goes back to Harnack himself, yet in the present twist this seems
to be due to the Russian scholars (Gudkov, Polotovskii, Viro,
etc.) Remember  that a similar trick was used in Gudkov's 2nd
construction of his novel sextic $5\frac{5}{1}$.
%%%(see Fig.\,\ref{GudkovCampo-5-15:fig}).
Actually as the curve is of
degree 8 one  could try achieving more oscillation than four, yet
it should probably be kept in mind that there is dormant a
supermassive black-hole of an $X_{21}$-singularity  with
invisible/imaginary branches located at the upper vertex of the
fundamental triangle and this absorbs 4 intersections with the
oscillated-about line.

\begin{figure}[h]\Figskip
%\vskip-1.2cm\penalty0
%\centering
\hskip-2.7cm\penalty0
\epsfig{figure=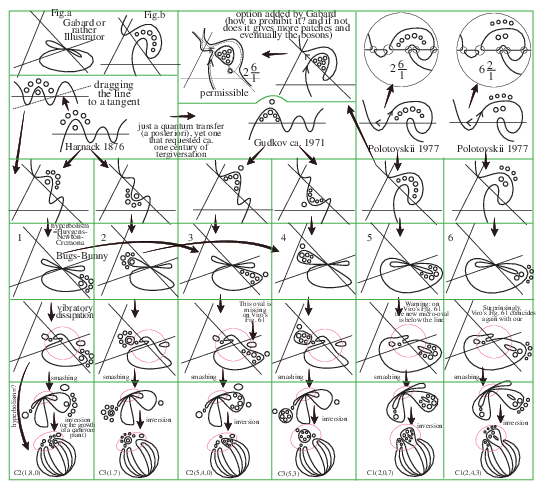,width=172mm} \captionskipAG
  \caption{\label{ViroDEGREE8_PATCH1:fig}%
  Viro's construction of the patches (via hyperbolisms): nearly
  complete modulo the mirrors to be found on the next plate}
\figskip
\end{figure}

The last step involves another hyperbolism, but alas here Viro's
paper lacks any depiction and becomes fairly incomprehensible. We
were thus blocked for a while at this place. Our heuristic idea is
that the hyperbolism might be interpretable as a smashing followed
by an inversion, and then indeed we recognize the patch C2(1,8,0)
in the notation of
Fig.\,\ref{ViroDEGREE8_exotic_patches0_SYS:fig}. Again the trick
of the inversion gives us a ``Wiedergeburt und a Neuauferstehung''
in Jakob Steiner's picturesque language. Of course, this inversion
is a naive vision-de-l'esprit and can be directly interpreted as
an isotopy akin to a mitosis in cell-biology, or better, as the
growth of a carnivore plant with several tentacles merging at the
opposite pole. The other cases are then self-explanatory and we
get successively the patches C3(1,7), C2(5,4,0), C3(5,3). Then the
first Polotovskii specimen yields the patch C1(2,0,7) not readily
listed by Viro, and which may thus be considered as
relatively-new. Remind that when combined with the hypothetical
patch C2(9,0,0), this yields the boson
$b7=1\frac{7}{1}\frac{12}{1}$.

This divergence from Viro is a bit puzzling and may suggest that
our interpretation of the hyperbolism is not entirely adequate.
Alternatively, it could be that our patches are correct, but those
of Viro are {\it not\/} since everything depends just upon the
location of the letter $\ga$ on his picture, and one can easily
imagine  this can be misprinted. This hypothesis, although poorly
founded, has to be envisaged with non-zero probability in view of
the density of misprints detected in Viro's paper. Of course it
can also be that both collection of patches (C1 and C2) both
exist, so that both Viro and ``Gabard'' are right without
excluding themselves.

Note incidentally that our earlier incertitude about the exact
location of the micro-oval arising from the triple point
dissipation (as being either above or below the oblique line)
seems anyway irrelevant in view of the merely isotopic nature of
the problem.

So at this stage we have the important:

\begin{Scholium}
Looking at the entraille of Viro's method we see that his list of
patches is possibly not exhaustive. Maybe Viro can be excused
because the new patches are producing isotopic schemes than those
generated by the old list, or because a slight misprint infested
his picture. Alternatively, it can be that our visual
interpretation of the hyperbolism is foiled, and then Viro is
correct.
\end{Scholium}

Now, Viro's argument has still not yet been exploited in full,
because he also proposes working with symmetrized replicas of the
ground octics numbered 1,3,5,6 (on
Fig.\,\ref{ViroDEGREE8_PATCH1:fig}); compare our
Fig.\,\ref{ViroDEGREE8_PATCH2:fig} especially the ``mirror''
curves M1, M3, M5, M6 near the bottom of that plate. Those produce
the patches C1(0,8,1), C1(0,4,5), C1(8,0,1), C1(4,4,1). It is
slightly puzzling that those are initially the palindrome of the
original patches (prior to the mirroring), but this palindromic
law is not always respected.

At this stage, we really exhausted Viro's discourse, yet getting
rather the C1-patches instead of the C2-versions. Further one can
wonder if there is the smoothing M52 where the micro-oval born out
of the triple point is located more on the right. This would give
the patch C1(7,0,2), which is however ruled out by Viro's oddity
law or by Orevkov's disintegration of the boson b3. Hence it seems
that the dissipation M52 has to be excluded.

\begin{figure}[h]\Figskip
%\vskip-1.2cm\penalty0
%\centering
\hskip-2.7cm\penalty0
\epsfig{figure=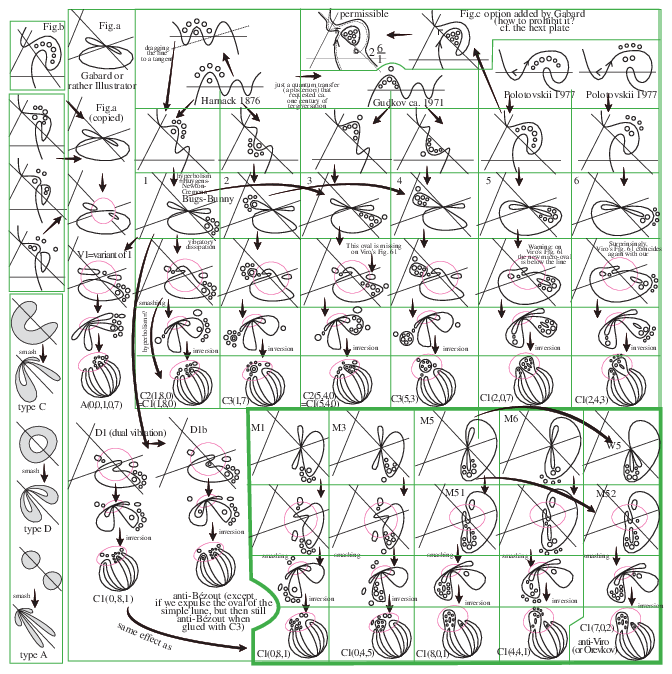,width=172mm} \captionskipAG
  \caption{\label{ViroDEGREE8_PATCH2:fig}%
  Viro's construction of the patches (continued, i.e. the mirrors)}
\figskip
\end{figure}

Besides going back to the very first construction (curve 1), we
may also a priori imagine a dual dissipating vibration (shown as
D1), but this yields again the patch C1(0,8,1) already obtained
via M1. Of course we could imagine that the micro-oval as on
Fig.\,D1b by being roughly speaking inside of the meander, yet the
resulting patch just violates B\'ezout.

Of course whenever $\ga=0$,then C1=C2, and so all the patches
constructed can be interpreted as belonging to class C1. Hence in
our interpretation all the patches constructed belong to C1 and
not C2 as asserted by Viro's Fig.\,55 (p.\,1118 of Viro 89/90
\cite{Viro_1989/90-Construction}). Of course it could be that
there is another constructions yielding the class C2, or if not it
could be that Viro's Figure contains a perfidious misplacement of
the symbol $\ga$ in the outer lune instead of the intern lune.

Incidentally all those constructions are based on the tricky issue
of the possibility of a vibratory dissipation and so it could be
that all constructions are actually foiled in case this crucial
step is fallacious.

At any rate we see that much variants have to be explored (for
instance our Fig.\,a and our Fig.\,c). More philosophically we see
that Viro's method is very much in continuity with Gudkov's
technique.

We may note that our Fig.\,c leads to a path violating B\'ezout
(when it is doubled). Still, we may hope that the method can be
more varied to produce new patches (so perhaps new $M$-schemes),
especially if we use also the rabbit with one ear invaginated
(Fig.\,a).

\begin{figure}[h]\Figskip
%\vskip-1.2cm\penalty0
\centering
%%%%\hskip-2.7cm\penalty0
\epsfig{figure=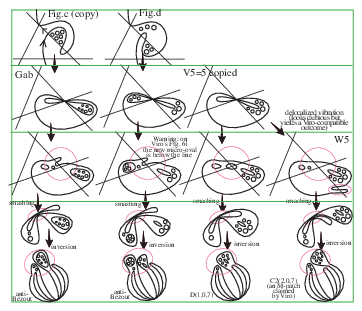,width=122mm} \captionskipAG
  \caption{\label{ViroDEGREE8_PATCH2B:fig}%
  Prohibiting Fig.\,c and other variants of Viro's hyperbolic method} \figskip
\end{figure}

[02.09.13, 23h39]
%%%At this stage
Disinhibited by some alcohol (cf. Ahlfors or the Greco-Roman
tradition in ``Vino veritas''), and according to our previous
dubitativeness (especially that based on Rohlin's maximality
principle corrupted by---despite its intrinsic beauty---Viro,
Gudkov, and even Hilbert, cf. one of the previous sections), we
must confess that most of those Viro's constructions of patches
seems to lack (severely?) in rigourousness to be
algebraically %%realizable.
christianizable. Of course, by social-ethnical affinity with Viro
(Leningrad, Kronstad[t]) the writer wishes Viro all the best,
i.e., all his constructions being correct (and perhaps more
importantly exhaustive). It seems that whatever happens (i.e.
Viro's death vs. Viro's triumph)
%%%and immortality)
there is---in the present state-of-affairs, at least---a serious
lack of didactic value in the presentation, admittedly imputable
to the intrinsic difficulty of the subject. Actually, the subject
itself looks not so much intrinsically difficult, yet sufficiently
boring in requesting high-memory faculty from the
%%%(alcoholized)
worker that the task turns quickly to an existential stress.

{\footnotesize  In comparison, some natural philosophers (like
Markus Schneider, oder aber, his young brother Traugott Schneider
(=two of the uncles of the writer on the maternal side of the
genetical tree) who are epistemologically satisfied with
phenomenologically more trivial apparitions of lesser logical
structuring (like Buddhas appearing in the clouds), yet still
apparently  able to enjoy life at some more primitive level of an
irrational Weltanschauung. This is ungef\"ahr wie bei dem Otto
Walkes (aus Ostfriesland), der den Markus so gut immitieren kann.

}

[03.09.13] An obvious challenge would be to realize the
wonder-patch C1(9,0,0) in which case we could  progress the
Hilbert-Viro 16th problem by constructing two new bosons namely
b1$=1\frac{1}{1}\frac{18}{1}$ and b7$=1\frac{7}{1}\frac{12}{1}$.
As yet the closest patch is C1(8,0,1) obtained via M5 and we could
try to contort this construction as to get the wonder-patch.
Basically, this would involve turning the upper loop of M5
upsidedown like on Fig.\,W5 for wonder. Alas it is clear that such
a contortion cannot be effected, just because then the number of
semi-branches emanating from the singularity would be split w.r.t.
the line smashed under the hyperbolism as $1+5$ instead of $3+3$
which is the sole reasonable splitting. So it seems that the
wonder patch---if it exists at all---requests a principally new
mode of construction.

We remind that Viro's construction as not yet been exposed in full
since it remains to construct the 4 patches of the family E=V1
(i.e., those with a trinested lune). Those are obtained by another
device of Viro (exposed in the next section), which is actually
somewhat more elementary in circumventing the use of hyperbolisms.

On the other hand we could suspect an affine quintic like our
Fig.\,d which differs from Polotovskii's just by dragging a
certain number of ovals in the wave extremity. However the
resulting patch is clearly anti-B\'ezout  because by gluing with a
symmetric copy we get 2 subnests.

Our next random idea is to wonder about the case where the
fundamental oval of the pre-smashed octic vibrates rather like a
horse-shoe. Yet, on more mature thinking we realize that the
example considered by Viro are already horse-shoes. Actually it
seems that the horse shoe is the only admissible shape for a cell
to intercept four-times a line. Another option is to have a ring
(annulus) where of course the cell-shape is changed. This suggests
looking at a variant of Fig.\,5, our Fig.\,V5 which gives the
patch D(1,0,7) referring to the notation of
Fig.\,\ref{ViroDEGREE8_exotic_patches0_SYS:fig}. Of course this is
not an $M$-patch, yet maybe still a valuable instrument to study
$(M-1)$-curves when patched with Viro's patches. When glued with
itself it gives the $(M-2)$-scheme $\frac{3}{1}\frac{15}{1}$ which
is fairly familiar (below Orevkov's obstruction and accessible via
a twisted gluing of Viro's $M$-patches).

Our next crazy idea is materialized by Fig.\,W5. Here we assume
that dissipating the triple point the oscillation does not take
place  in the vicinity of the smoothed point but faraway like in a
non-local phenomenon of quantum chromodynamics. The net effect of
this crazy modification is just that the ovals are transferred in
the outer lune and so we get the patch C2(2,0,7), referring as
usual to the notation of our
Fig.\,\ref{ViroDEGREE8_exotic_patches0_SYS:fig}. It is evident at
this moment, that we can (assuming that this delocalized quantum
vibration is always possible) reach all the patches of the C2
family as claimed  in Viro 89/90 (granting that there is no
misprint on his fundamental figure, i.e. no misplacement of the
$\ga$ parameter).

{\it Added\/} [15.09.13].---Let us check this assertion more
pedestrianly. So we copy the Viro table, and rationalize a bit our
depiction of it by killing the ``carnivore-plant stage'' which is
quite unnecessary by the way. Then, for each of Viro's
constructions, we consider the avatar with a vibration acting at
long distance (see Fig.\,\ref{ViroDEGREE8_PATCH2_QUANTUM:fig}).
Although counter-intuitive this is still B\'ezout compatible so
that there is perhaps an algebraization of such pictures. By this
recipe we get indeed patches claimed by Viro starting with
C2(0,8,1). We report the patches so realized in the catalogue
(Fig.\,\ref{ViroDEGREE8_exotic_patches0_SYS:fig})  by yellow-green
rectangles to emphasize the issue that those patches are  more
dubious than the evergreen patches gained by a localized
vibration. Next, we get successively the patches C3(1,7) (as
above), C2(0,4,5) (relatively new and the palindrome of the
above), C3(5,3) (the same as above), C2(2,0,7) (new and the
C2-avatar of the above), C2(2,4,3) (new and the C2-avatar of the
above).

At this moment, a look on the patch-catalogue
(Fig.\,\ref{ViroDEGREE8_exotic_patches0_SYS:fig}) shows that we
still miss two patches claimed by Viro, namely C2(4,4,1) and
C2(8,0,1).

Another idea is to invert the sense of the meander as shown on
series~C. Here we get first a monotonic repetition of the first
four patches of the original Viro's series (A), but when it comes
to 5, we find the patch C2(7,0,2) violating Viro's oddity law as
well as Orevkov's (link-theoretic) prohibition of b3. Likewise the
sixth configuration, yields the patch C2(3,4,2) corrupting Viro's
law of oddity.

\begin{figure}[h]\Figskip
%\vskip-1.2cm\penalty0
%\centering
\hskip-2.7cm\penalty0
\epsfig{figure=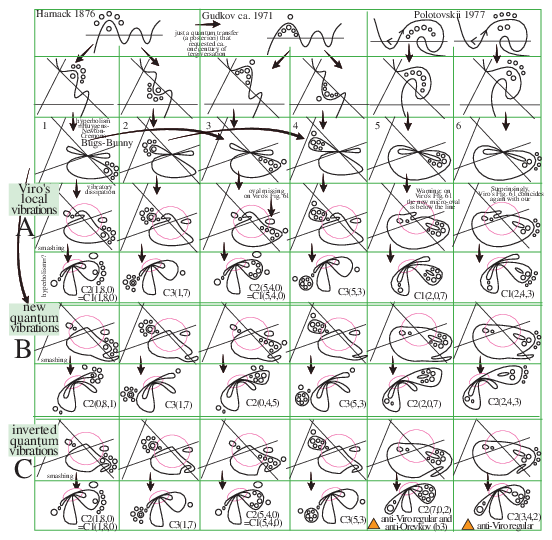,width=172mm}
\captionskipAG
  \caption{\label{ViroDEGREE8_PATCH2_QUANTUM:fig}%
Quantum vibration acting at long distance} \figskip
\end{figure}

Alas we still do not have obtained all patches claimed by Viro,
and still miss the two mentioned ones in the C2-class.

At this moment we had the idea (especially when looking at the
entry C-5 of the previous table) that we may get something
different (hence new) when smashing instead of up as we did till
now. The next figure gives a systematic tabulation of those
smashing down by giving them in the little right window of each
configuration of the previous table. It is observed that the down
smashing is isotopic to to up version for the first four entries
1,2,3,4, but leads to a different patch for fifth entry, namely we
get C2(2,0,7) at entry A5-down, instead of C1(2,0,7) at A5-up. As
far as our presentation is concerned this corroborates (with
delocalization) Viro's claim of existence of this patch. Likewise
along the sixth entry 6, we get the C2-version or the earlier
C1-patch. All this concerns Viro's A-series, but little impede us
repeating the story of downwards smashings for the other more
exotic quantum vibrations (B- and C-series). For the B-series we
get nothing new, and likewise for the C-series (apart a repeated
conflict with the Viro/Orevkov obstructions).

\begin{figure}[h]\Figskip
%\vskip-1.2cm\penalty0
%%%\centering
\hskip-2.7cm\penalty0
\epsfig{figure=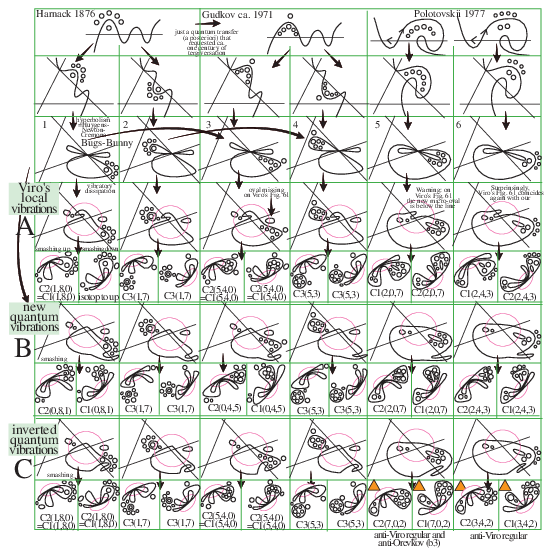,width=172mm}
\captionskipAG
  \caption{\label{ViroDEGREE8_PATCH2_QUANTUM_DOWN:fig}%
Smashing down} \figskip
\end{figure}

The abstract mechanism of the smashup-smashdown option is that it
seems to put in duality the C1- and C2-classes. Hence whenever we
have a C1- or C2-patch it suffices altering the up/down option to
get the ``same'' patch in other class that is with same symbols
$\al,\be,\ga$. This being understood it is now a simple game to
get the two C2-patches claimed by Viro, of which we were only able
to get the C1-decoration. To be specific to get C2(4,4,1), we look
at C1(4,4,1) and browsing back through our figures we identify the
construction as being via the mirror M6, and it suffices altering
this by smashdown to get the desired patch C1(4,4,1), as shown on
Fig.\,\ref{ViroDEGREE8_PATCH2_M6_DOWN:fig}. Likewise to get
C2(8,0,1), we look at its $C1$-companion C1(8,0,1), and then
browse back through earlier table to find its mode of generation
via the mirror M5, and just alter this by a smashing down.

\begin{figure}[h]\Figskip
%\vskip-1.2cm\penalty0
\centering
%%%%\hskip-2.7cm\penalty0
\epsfig{figure=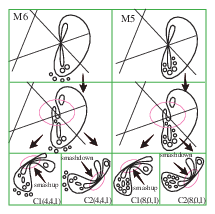,width=122mm}
\captionskipAG
  \caption{\label{ViroDEGREE8_PATCH2_M6_DOWN:fig}%
Getting the last two patches claimed by Viro} \figskip
\end{figure}

At this stage we have completed our understanding of Viro's
theory, and believe that our C1-patches also exist albeit this is
not explicitly mentioned in Viro 89/90.

\subsection{On a bending principle for patches}

[16.09.13] It seems that one could supply a formal definition of a
patch as being the trace of a real algebraic curve on the unit
ball $B^4\subset \CC^2$ and bounding the link of the singularity
$X_{21}$, which has 4 components (in general as many as the
singularity has branches). Now it seems that given a patch one can
bend its content by changing the curvature of all its branches
simultaneously. Referring to our patch catalogue, we see then that
C1 turns to C2 and viceversa. The class A is self-dual. The class
B transmutes to the class I. Class C is self-dual, with its
subspecies C1 and C2 permuted, while C3 stays invariant under
bending. The class D is also invariant under bending.

The class E is likewise invariant under bending. So at least at
the topological level, but it seems interesting to speculate about
a higher form of algebraic invariance, that is of algebraic
patches. As noted yesterday, this is in part motivated by the
issue that smashing-up or -down hyperbolisms gave within the
C-class the dual patch. Alas the E-patches are not as yet realized
via hyperbolisms, and some heuristic thinking about the shape of
embryos under smashing inclines one to believe in the
impossibility of a such a realization.

On the other hand, the bending duality transforms the patch
E($\al,\be,\ga$) into E($\ga,\be,\al$) (compare the
catalogue=Fig.\,\ref{ViroDEGREE8_exotic_patches0_SYS:fig}). As we
shall see in the next section Viro's construction produces four
patches in the E-class, which are in stable equilibrium under the
bending symmetry. Hence the hypothesis of invariance under bending
appears as experimentally verified, but it would be nice to find a
theoretical justification (if any). One can imagine that a simple
hyperbolism (Cremona transformation) transmutes a patch into its
bending, and if this exists it must easy to write down. So:

\begin{Scholium}
The whole theory of Viro's patches is probably invariant under the
bending involution transforming an $X_{21}$-patch to its companion
with branches of inverted curvatures.
\end{Scholium}

Albeit only heuristic, this principle suggests new results or at
least novel predictions as side effects. For instance in the
I-class we have 3 restrictions caused by Shustin (hoping this to
be true). But as the patch I($\al,\be,\ga$) bends to
B2($\ga,\be,\al$), it would result three new obstructions in the
B2-patches as reported by red rhombic crosses on the catalogue
(Fig.\,\ref{ViroDEGREE8_exotic_patches0_SYS:fig}).

Next the patch F bends to the patch H, compatibly with the fact
that both families are empty as inferred from Arnold's weak
version of Gudkov periodicity.

The patch G is in contrast self-dual under bending, and again
bending acts palindromically over the symbols $\al,\be,\ga$. So at
least for the class G1 where bending G1($\al,\be,\ga$) yields
G1($\ga,\be,\al$). Thus again by the posited principle of
invariance one can reflect three of Shustin's prohibitions to get
three new (but hypothetical) restrictions marked by rhombic
crosses in the catalogue. Bending is somehow akin to an inversion
of the magnetic poles of planet Earth. In particular it causes
serious damages on the patch family G2. This bends indeed to a
configuration which is anti-B\'ezout [at least if $\ga>0$]. As a
novel consequence bending would permit to rule out the three
remaining patches of the G2-family as to make it completely empty.
[Warning: in fact in view of the proviso in bracket, I am not
completely sure about this conclusion.] In fact  bending G2(0,8,0)
gives G1(9,0,0), which is prohibited by Shustin, hence so is
G2(0,8,0). For G2(4,4,0) the bending is G1(5,4,0) (because in
general G2($\al,\be,0$) bends to G1($\be+1, \al,0$)). However the
 latter is prohibited by Shustin. Finally the last standing man
G2(8,0,0) bends to G1(1,8,0) also prohibited by Shustin. Thus:

\begin{lemma}
Shustin's obstruction diffuses from the G1-family into the
G2-family via bending and par-achieves killing completely this
family, whose generic member G2($\al,\be, \ga>0$) was already
killed by B\'ezout via the gluing E/G2, or via bending of G2 to a
scheme enlarging the nest of depth 4.
\end{lemma}

Then our discussion is essentially complete since H is in duality
with F, I with B, and J is self-dual (and trivially understood via
B\'ezout).

\begin{figure}[h]\Figskip
%\vskip-1.2cm\penalty0
%%%\centering
\hskip-2.7cm\penalty0
\epsfig{figure=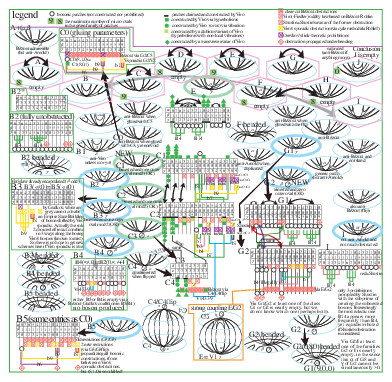,width=172mm}
\captionskipAG
  \caption{\label{ViroDEGREE8_exotic_patches0_BEND:fig}%
Catalogue of all patches under bending duality} \figskip
\end{figure}

The moral of all this is that invariance bending looks an
extremely versatile tool. On the one hand it is like a hidden
symmetry explaining why the $M$-patches C and E are the most
populated (even the sole populated by algebraic representatives
according to the present state of knowledge), and this  despite
lacking symmetry along the vertical axis.

Further it seems that bending affords a powerful method of
prohibition. For instance, bending the patch B3 corrupts B\'ezout
for line provided $\be >0$. Thus we can conclude that the family
B3 is nearly empty. Actually, it may be considered as empty since
if $\be=0$ the family B3 degenerates inside B2. Similarly bending
B4 is anti-B\'ezout when $\ga>0$, and if not then B4 collapses to
B2. Finally, on bending B5 we get a configuration anti-B\'ezout as
soon as either $\be$ or $\ga$ is positive (glue with the flip).
Hence strikingly, this re-explains all the first series of
red-crosses (prohibitions derived from Viro's oddity law) in a
much more elementary fashion. Even more, even the remaining three
B5-patches with $\ga=0$ are killed, e.g. B5(8,1,0), since---as
noted---gluing with the flip corrupts B\'ezout. Hence it would
follow (from bending invariance) that the family B5 is completely
extinct.

Looking closer to the remaining B4-patches with $\ga=0$ we
identify them as B4($\al,\be,0$)=B2($0,\al,\be+1$), and so the
specific element B4($7,1,0$) is B2($0,7,2$) which is not even
present in the list as a consequence of the periodicity modulo
four imposed on $\be$. Yet, doubling B2($0,7,2$) gives the scheme
$6\frac{15}{1}$ which respects Gudkov and even exists since Viro.
Thus we are never entirely sure that we did not from the scratch
overrestricted the parameters. Actually, in the B2-class $\be$ was
pre-calibrated as $1 \pmod 4$, but it could be just odd. The
crucial point is just that when doubled B2 yields a subnested
$M$-scheme with $2\be$ big eggs (=ovals at depth 1). But this
quantity has to be $2 \pmod 4$ by Gudkov periodicity (compare the
periodic table of elements), and thus $\be $ has to be odd. The
point however is that the case  $\be \equiv 1 \pmod 4$ and
$\be^\ast \equiv 3 \pmod 4$ cannot cohabit because  the gluing of
both patches would have $\be+\be^\ast\equiv 0 \pmod 4 $ big eggs,
violating Gudkov periodicity. So we have, alas, to enrich the
catalogue by the family B2$^\ast$ where $\be$ is $3 \pmod 4$, and
we lack unfortunately a recipe to prohibit them. Of course all of
B2$^\ast$ would be killed in one stroke if there were any element
in B2, yet  Viro's theory seems rather to tell that B2 too is
empty.

[17.09.13] In contrast, one can entertain the dream that the B2
family is non-void. As we saw, the patch B is not readily seen as
the product of a smashing hyperbolism acting upon a ground shape
(topological cell, binion=annulus, or 2 disks).
Fig.\,\ref{ViroDEGREE8_exotic_patches_B:fig} shows indeed that
under a topological smashing only patches of type C, D and A do
occur, from respectively a cell of the horseshoe type, a ring or a
double cell. This basic experiment adumbrates also why only type-C
patches achieve Harnack-maximality of 9 micro-ovals, since in the
other two cases (ring or bi-disc) one oval is already wasted as
contour of the ground shape. Albeit heuristic, this justification
stands in perfect accordance with Viro's theory.

\begin{figure}[h]\Figskip
%\vskip-1.2cm\penalty0
\centering
%%%\hskip-2.7cm\penalty0
\epsfig{figure=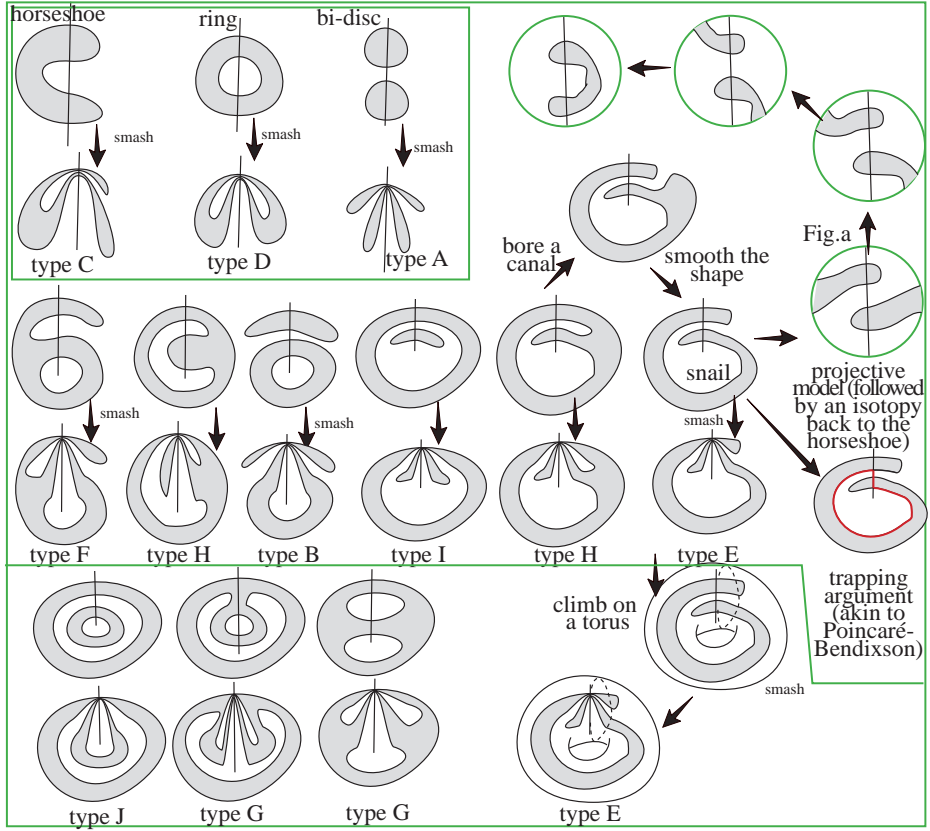,width=122mm}
\captionskipAG
  \caption{\label{ViroDEGREE8_exotic_patches_B:fig}%
Morphogenesis of patches under smashing hyperbolisms} \figskip
\end{figure}

It is evident that the embryology of
Fig.\,\ref{ViroDEGREE8_exotic_patches_B:fig} is exhaustive, and so
we get:

\begin{lemma}
Viro's method of the hyperbolism is only capable producing
$M$-patches of type {\rm C\/}, and eventually $(M-1)$-patches of
type {\rm D} and {\rm A}. In particular there is no
%%hope
chance to get new $M$-patches of type {\rm B} via the hyperbolism
method.
\end{lemma}

Actually, we can obtain the B-type through smashing  a ground
shape involving 3 contours. As a slight surprise we can also
access the E-type through smashing a snail (with only one
contour). So this raises some hope that the method of hyperbolism
could as well produce patches of type E, and eventually new ones,
as it must be confessed that Viro's knowledge of the E-class
remains fairly
%%%sparse
lacunary (compare the catalogue).

Our snail-model leading to type E involves the smashing of a line
somehow trapped inside the bag of the snail. So by Jordan
separation applied to the Bendixson-style bag formed by the spiral
arc closed by the transverse segment (cf. red contour on the
figure), we get a trapping of the one-end of the line. Hence the
line should actually have more than the four visible intersections
with the snail, but after smashing this yields a singularity with
more than four branches (hence not a dissipation of $X_{21}$).

One can try to find a refuge by  imagining a model of the snail in
the projective plane, but ours (Fig.\,a) seems isotopic to the
horseshoe. So it seems that the Poincar\'e-Bendixson-style
trapping argument prevails, preventing thereby the creationism of
E-patches via hyperbolisms (in accordance with Viro's praxis, yet
disappointing from an avantgardist viewpoint). More generally all
other configurations outside of the sub-frame of
Fig.\,\ref{ViroDEGREE8_exotic_patches_B:fig} involve a trapping
and therefore inadmissible for creating  the corresponding patches
(all types safe C,D,A) by the method of hyperbolisms.

Despite the trapping obstruction, it is amazing that from this
naive embryological viewpoint only the types C and E arise by
smashing a single simple cell,  those being precisely the
patch-families populated by algebraic representatives according to
Viro's theory. This is a striking coincidence, but alas
(apparently) not vivid enough to bring the E-class within the
range of hyperbolisms.

One way to get
%%rid off
around the trapping obstruction would be to replace the plane by a
torus, and imagine the figure of E-type traced on a such
(materialized as usual by a quadric in 3-space abstractly
isomorphic to $\PP^1 \times \PP^1$). Then we can smash along a
meridian of the torus without encountering a trapping obstruction.
It maybe speculated that the whole construction projects down to
the plane, as to get patches of type-E, hopefully of a new sort
not yet listed by Viro. Of course, all this  looks fairly tricky
and hard-to-implement geometrically, yet perhaps there is some
chance to concretize this idea.

Of course, it would be overall simplifying if all patches of the
E-class existed. In this scenario the two subnested bosons would
be created and most of Korchagin, and Chevallier's octic schemes
trivialized to Viro's simplest method (as shown by the relevant
composition-table, i.e., V1/V1 on
Fig.\,\ref{ViroDEGREE8_extended:fig}). Alas, for the moment it is
extremely hard to decide which scenario corresponds to  reality.

[16.09.13] The overall philosophy is that under  bending
invariance we can rule out huge family of patches just by basic
reliance upon B\'ezout, hence trivializing most of the highbrow
Viro-style obstructions.

Another consequence of the bending hypothesis is that upon looking
at the E-patches the obvious creationism of the boson B4 (via the
patch E(5,1,3) would automatically create the palindromic patch
E(3,1,5), and consequently materializes the ``dual'' boson B14.
Unfortunately, since Viro's method is just a method, and  not an
intrinsic feature of the universe, one cannot deduce the existence
of one boson being coupled to that of its dual, i.e., B4 exists
iff B14 does exist. However this becomes true as soon as one of
the boson is realized as perturbation of the quadri-ellipse,
provided our bending principle holds true.

{\it Stupid remark.}---Of course, bending translates to
palindromic symmetry in the E-class but not in the C-class.
Otherwise, applying palindromic symmetry in the C1-class leads to
serious contradictions in mathematics, like Viro's construction of
C1(4,4,1) conflicting with Viro's prohibition of the palindrome
C1(1,4,4). This is just mentioned in order to avoid the reader
doing the same basic mistake as the tired writer.

Also, note that the class C3 is self-dual under bending, and
actually pointwise invariant. Hence we cannot infer any new
information, but the census in this family was complete just under
basic B\'ezout obstructions and Viro's constructions (either in
the hyperbolism setting or via the more elementary vibratory
method, as we shall see later).

Let us summarize the situation as follows:

\begin{lemma}
If the bending principle is true, then most families of our
catalogue of exotic patches are empty by virtue of a trivial
reduction to B\'ezout. More precisely the $M$-patches families
{\rm B3}, {\rm B4}, {\rm B5} are empty, and so is {\rm G2}.
However the dual pair {\rm B2} and {\rm I}, and also self-dual
class {\rm G1} are potentially non-empty, as to contain algebraic
patches favoring the creation of the bosons {\rm  B4} and {\rm
B14}. Alas, presently nobody ever succeeded constructing any such
patch. Besides, the classes {\rm  C} and {\rm E} are self-dual yet
not completely elucidated, when looking very deeply into the glass
(all white-circles in the catalogue are not yet known). As a last
remark, the method of hyperbolism amounting to a smash, seems only
able to create the {\rm C}-patches (when smashing a horse-shoe),
or {\rm D}-patches (when smashing a ring), or finally {\rm
A}-patches when smashing a pair of discs. Accordingly, it should
be no surprise that the {\rm E}-patches are not accessed by
hyperbolisms, which in the realm of maximality produce only {\rm
C}-patches.
\end{lemma}

In sum this means that under bending we are fairly close to
getting a complete census of all patches modulo the ambiguity left
in the B2/I-classes, the G1-class,  the completion of the C1/C2
pair which is nearly settled if Viro-Orevkov are true (sole
exception C1(9,0,0)=C2(9,0,0), and the E-class which is still
elusive (no known prohibitions).

Of course to substantiate all this bending hypothesis one should
take the pain  to write down a simple (algebraic) transformation
doing the requested deformation of bending. We imagine this must
be a trivial task. Besides, bending invariance seems to fit with
all factual data available up to now.

As we saw the idea of bending trivializes at the basic B\'ezout
level several obstructions first derived via Viro's oddity law.
This reduction truly concerns patches, yet it would be of interest
to wonder about a global avatar of bending acting on projective
curves and reducing the Viro obstruction to B\'ezout.

As a last philosophical touch, we always found Viro's census of
patches shocking, as it seems to corrupt the ``Didon principle''
of extremality: {\it most solutions to  extremal problems are
  inhabited by deep character of symmetry\/}, like in the
isoperimetric problem, the Bloch constant, etc.) However, in Viro
census the vertically-symmetric patches (typically classes B or I)
lack apparently any representative, while the asymmetric families
C and D extremalize the number of ovals. As already noted, this
violation of the Didon principle is relaxed if we imagine bending
as a hidden symmetry of the problem.

\subsection{Speculating about a global bending: alias the NAS=
Neuauferstehung}

[17.09.13] If we look the table of $M$-symbols
(Fig.\,\ref{SIMPLIFIED-TABLE:fig}), and in it on the 3rd pyramid
of subnested schemes one is struck by a certain symmetry of the
symbols putting in duality the symbol $x(1,n \frac{y}{1})$ with
$y(1,n \frac{x}{1})$. Geometrically, this is merely a reflection
of each vertical rows about its midpoint. Strikingly,  this
duality almost everywhere respects the constructor of the scheme.
The sole exceptions are:

$\bullet$ in the 1st row,  the pair $7(1,2\frac{11}{1})$
(O=Orevkov) and $11(1,2\frac{7}{1})$ (V=Viro);

$\bullet$ in the 2nd row, the pair $3(1,6\frac{11}{1})$
(K=Korchagin) and $11(1,6\frac{3}{1})$ (V=Viro);

$\bullet$ in the 3rd, 4th and 5th rows, there is no exception to
the symmetry of the constructors.

It seems natural to speculate about a global duality of a
geometric nature explaining directly this symmetry. Of course one
way would be the symmetry on patches, but this would really work
if knew how to construct all E-patches.

The sole objection against this duality comes from Shustin's
obstruction of subnested schemes without outer ovals. Indeed under
our duality, those schemes  correspond to the relevant
simply-nested schemes, all of which exist since Harnack, Gudkov
and Viro. Thus, there is a violent break of symmetry, only
remediable by rejecting Shustin's obstruction. More politically
correct, is to imagine that our duality reigns only over a
restricted range not going  as far as Shustin's series.

Is there a direct geometric interpretation of this hypothetical
duality? At the (soft) topological level via
Fig.\,\ref{ViroDEGREE8_duality:fig} below, we may just interpret
the duality as being merely a human face with an indigestion of
smarties (chocolates) in the mouth, and exchanging this out and
in. So we can call this the digestive python-duality, or the
smarties-duality. Of course, this move or rather exchange can be
imagined as the combination of two Morse surgeries: first the
mouth content is liberated by a fission of the buccal cavity, and
in turn new labial expansions phagocytose the outside. Of course
this digestion process looks magic and delicate to ape in the
algebraic category. Besides, the whole process transits through a
single intermediate uni-nested $(M-1)$-curve.

\begin{figure}[h]\Figskip
%\vskip-1.2cm\penalty0
%%%\centering
\hskip-2.7cm\penalty0
\epsfig{figure=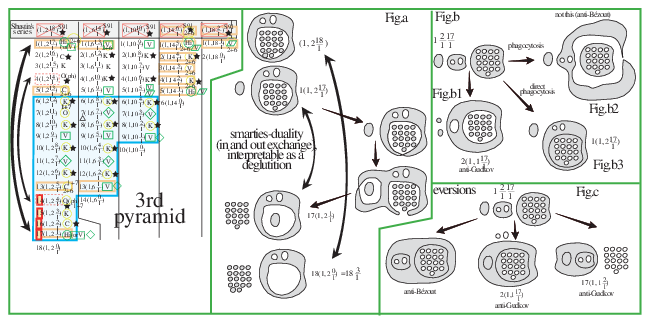,width=172mm} \captionskipAG
  \caption{\label{ViroDEGREE8_duality:fig}%
Subnested Duality} \figskip
\end{figure}

 Imaginatively, and without referring to a rigid isotopy,
 this duality could be
realized instantly like an inversion with respect to the ring
formed by the nonempty oval at depth 0 and the subnest (i.e.
nonempty oval at depth 1). In case of Shustin's series the duality
appears as broken since there is no canonical way distinguishing
the mouth from the eyes in the figure for $18\frac{3}{1}$. One may
dream about a God-given Cremona transformation effecting such an
inversion directly. By any reasonable Nullstellensatz (algebraic
rigidity), we cannot expect the transformation to fix (pointwise)
the ring, nor any of its contour. In parallelism, one could search
 a synthetic rule doing the inversion like via a transformation
of reciprocal radii. \`A la Zeuthen, we would  like to infer some
convexity properties of the deepest oval, and exploit this for a
synthetical inversion keeping
 algebraicity and the degree constant to 8. If implementable
both bosons B4 and B14 would be
%%married
coupled in strong-duality, with joint destiny of life or death,
like
%faithful
inseparable, fusional partners.

\subsection{Trying (but failing) to extend  duality outside
the subnested realm}

[18.09.13] It seems of even greater importance and predictive
power to extend the duality to the other pyramids (binested and
trinested case), since in those cases the available theory looks
more lacunary and of lesser symmetry that in the 3rd pyramid where
everything just depends on the existence of the two bosons. As
above we may first like to guess the symmetry at the combinatorial
level of symbols, as some intrinsic symmetry of the table of
periodic elements.

From the very beginning, one feels subconsciously the duality
linking say both of Orevkov's schemes $1\frac{3}{1}\frac{16}{1}$
and $1\frac{6}{1}\frac{13}{1}$ (both prohibited by link theory).
More generally this involves the involutary symmetry
$k\frac{x}{1}\frac{y}{1}\mapsto k\frac{x+10}{1}\frac{y-10}{1}$,
where a packet of 10 is traded between both nests. A priori this
looks special astrological numerology involving the number 10 (of
our fingers), and there is of course many variants involving other
integers. So let us abort this viewpoint in the hope to find some
more intrinsic symmetry based on the idea of a topological
inversion.

For instance starting from Viro's binested scheme
$1\frac{2}{1}\frac{17}{1}$, one may imagine to invert  with
respect to the non-empty oval containing the least number of ovals
(Fig.\,b1). It results the scheme $2(1,1\frac{17}{1})$ which does
not even respect Gudkov periodicity (as the number of big eggs
must be $2 \pmod 4$). Of course the same aberration occurs when
inverting w.r.t. the nonempty oval of greatest content. The
problem is that the outer ovals of a bi-nest (congruent to $1\pmod
4$) becomes under inversion the number of big eggs of the subnest,
but the latter has to be $2 \pmod 4$. Maybe, an ad hoc  correcting
intervention can repair this, but looks unnatural.

Another piste of longstanding is some natural symbolic duality
between symbols like $1\frac{2}{1}\frac{17}{1}$ (V=Viro) and
$1(1,2\frac{17}{1})$ (Hi=Hilbert), etc, given generally by the
formula $k\frac{x}{1}\frac{y}{1} \mapsto k(1,x\frac{y}{1})$. For
Orevkov's $1\frac{3}{1}\frac{16}{1}$, this leads to
$1(1,3\frac{16}{1})$ which is outside of Gudkov's range. Viro's
existing scheme $1\frac{5}{1}\frac{14}{1}$ is carried to an
anti-Gudkov scheme, hence foiling the invariant character of the
postulated symmetry. Next Orevkov's anti-scheme
$1\frac{6}{1}\frac{13}{1}$ dualizes to $1(1,6\frac{13}{1})$ which
exists by Viro's basic theory of the quadri-ellipse. Hence our
symbolic duality does not seem to respect the intrinsic nature of
algebraic-geometry. Notwithstanding we may seek a direct geometric
interpretation of it. Alas  several phagocytosis attempts failed
miserably. A first such, is a cannibalistic attempt of the small
nest to annex (eat) the large one, as shown on  Fig.\,b2.
Unfortunately, the resulting curve violates B\'ezout (extension of
the deep quadri-nest). So the desired phagocytosis involves a
direct absorption of the big nest without forming a buccal cavity
so-to-speak (Fig.\,b3). Alas it seems that there is no Morse
theoretical historiography for such a move. Further, in contrast
to the earlier smarties-duality, the present one seems
chromatically anomalous w.r.t. to the natural black-and-white
coloration of the Ragsdale membrane.

In conclusion, both experimental and theoretical evidence seem to
fight against the symbolic duality $k\frac{x}{1}\frac{y}{1}
\mapsto k(1,x\frac{y}{1})$, where it is assumed $x\le y$.

Besides, there are several eversions shown on Fig.\,c. Yet, as we
already knew, those fails dramatically  to respect Gudkov
periodicity, and therefore are prohibited despite the lack of
evident topological obstructions.

At this stage our quest of hidden symmetries in the periodic table
elements seems already exhausted, and to be in ``panne''.

In conclusion, the smarties-duality (or {\it deglutition\/}) seems
to be the sole global tangible symmetry  we could detect in the
periodic table of elements.

In fact, we can also speculate more about the dubious symmetry by
10. In the binested realm, especially in the bosonic strip, this
symmetry takes the boson $1\frac{1}{1}\frac{18}{1}=:b1$ to
$1\frac{11}{1}\frac{8}{1}$ (constructed by Viro), so giving some
evidence for the materialization of the boson $b1$. Next Viro's
scheme $1\frac{2}{1}\frac{17}{1}=:b2$ dualizes to the boson
$1\frac{7}{1}\frac{12}{1}=:b7$, which get so some existential
probability. Next we have Orevkov's pair, both prohibited, so that
our postulated symmetry by ten-trading is still plausible. Finally
the boson $1\frac{4}{1}\frac{15}{1}=:b4$ is allied to Viro's
scheme and therefore likely to materialize.

Next we may also try to extrapolate the ten-trading symmetry in
the realm of trinested schemes. Looking in the 4th row of the 2nd
pyramid, we see first $4\frac{1}{1}\frac{1}{1}\frac{13}{1}$
(S=Shustin) in duality with $4\frac{1}{1}\frac{11}{1}\frac{3}{1}$,
also due to Shustin. Next we have
$4\frac{1}{1}\frac{2}{1}\frac{12}{1}$ which is self-dual, or
dualizes to $4\frac{11}{1}\frac{2}{1}\frac{2}{1}$, all being
anti-Viro regular. Next $4\frac{1}{1}\frac{3}{1}\frac{11}{1}$ is
also in duality (rather trinity) with
$4\frac{11}{1}\frac{3}{1}\frac{1}{1}$, that is itself, so that
everything is right. However when its comes to the next scheme
$4\frac{1}{1}\frac{4}{1}\frac{10}{1}$, the trading-by-ten brings
it to
$4\frac{1}{1}\frac{14}{1}\frac{0}{1}=
5\frac{1}{1}\frac{14}{1}$ (due to Viro). So we get an existential
conflict w.r.t. to our dubious symmetry. Of course the scheme in
question dualizes also to $4\frac{11}{1}\frac{4}{1}\frac{0}{1}=
5\frac{4}{1}\frac{11}{1}$ (also constructed by Viro).

In conclusion, it seems that our naive trading by 10 is not
compatible with actual knowledge, and of course looks very dubious
numerology, without strong geometric support.

[19.09.13] As a last attempt it can be imagined that one of the
bi-nest is directly glued  inside one of the oval of the other
nest, while another oval is created outside to compensate the
loss. But basically this gluing yields again either Fig.\,b1 or
Fig.\,b3 depending on how  ``outside'' is interpreted. The first
option is anti-Gudkov, while the second interpretation amounts to
a direct phagocytosis of Fig.\,b3. As we showed this runs against
troubles when it comes to Viro's scheme $1\frac{5}{1}\frac{14}{1}$
whose dual $1(1,5\frac{14}{1})$ violates Gudkov periodicity (even
in the weak sense of Arnold). Likewise our duality is disrupted
for Orevkov's anti-scheme $1\frac{6}{1}\frac{13}{1}$ whose dual
$1(1,6\frac{13}{1})$ exists by Viro's simplest method
(quadri-ellipse).

\subsection{Overview}

[19.09.13] Most of the mathematicians are not discovering new
fruits but just distilling old ones, to high-condensed beverages,
that nobody is virtually able to drink, without serious
%damage on his brain.
intoxication (brain damages). Viro's theory is
%certainly not an exception to this pitiful state-of-affairs.
no exception to the rule. Somehow, we need to present its
``deploiement universel'', at the level of the primitive fruits
%to
so that everybody can understand (and check) the whole
distillation process. In particular the theory of patches looks to
us still unachieved, and so is the global Hilbert's 16th problem
in degree 8. It seems that what remains left are just peanuts (a
negligible proportion of 6 bosons over the 104 logically possible
cases), yet this can safely occupy several generation of workers
unless one finds the correct ideas, perhaps possible revisions,
and rationalizations of the existing theory.

Roughly, it seems that several strategies could intermingle to
complete our understanding of degree 8.

1. Construction of new patches for $X_{21}$ via an adaptation or
renovation of Viro's methods. Note the plural since the
%%core
epicenter of Viro's method splits apart into hyperbolisms and a
basic vibrational method with tangency (a sort of non-transverse
avatar of Harnack-Hilbert-Brusotti).

2. Prohibition either via the method of total reality (Riemann et
ali) or just via tracing a curve  interpolating the deepest nests.
This is what we call the {\it method of  deepest penetration\/},
which we are quite incapable to implement seriously. Basically,
one should imagine that several regions of the periodic table of
$M$-elements (Fig.\,\ref{SIMPLIFIED-TABLE:fig}) are frozen because
they are over-nested. Algebraic curves like nesting, but they
cannot be too nested as evidenced by B\'ezout for lines. Higher
order curves than lines or conics should prompt novel (or
semi-novel) obstructions of the Rohlin-Fiedler-Viro era, with
recent ramification in Orevkov's link theory.

One can employ the usual metaphor about phase-changes between
solid, liquid and gaseous states, with the bosons unambiguously
identified to the intermediate liquid-state with ultimate destiny
 yet undecided. What freezes a scheme to the algebraic crystal or
 in contrast evaporate it as a nebulous gas unobservable with
 naked eyes is the mystery of those bosons.

3. The method of hidden symmetry in order to detect symmetry
patterns in the table of elements (especially $M$-elements), while
guessing the underlying geometrical motives. Alas, presently we
failed to
%%% deceal
disclose %%%% deceler in DICO OR REVEAL
any such symmetry beyond the deglutition-symmetry of subnested
schemes exchanging the inside of the subnest with the outside of
the primary nest (see Fig.\,\ref{ViroDEGREE8_duality:fig}a). Alas,
this hypothetical duality does not readily afford new concrete
information on Hilbert's 16th, safe for a coupling of the
existential destiny of both subnested bosons.

4. This method of duality admits apparently a semi-local avatar at
the level of patches, where it seems experimentally
%%realist
%%tangible
sound to expect a duality of bending
%yielding
prompting a global symmetry
%on
over all patches for $X_{21}$. This could be an important tool to
complete the classification of patches.

5. Finally, the whole philosophy of the method of small
perturbation up to its ultimate era of glory reached in Viro's
method seems to use a infinitesimal gluing principle of patches
inside algebraic objects much akin to surgeries feasible usually
in the smooth category. That the surgeries works algebraically
while keeping the degree controlled is much miraculous, and in the
case of the simplest nodal singularity already amounts to the
Severi-Brusotti transcription of Riemann-Roch. The general case of
independence of smoothing  is the credit of Gudkov-Viro-Shustin
under varied decorations. Of course course conceptually it seems
that the principle of gluing as local surgeries encompass the
principle of independence of smoothing: just patch locally and
contemplate globally.

6. In principle, it is expected that Hilbert's problem at least in
degree 8 is reducible to this method (Viro's patchwork), either
from the sole quadri-ellipse or via more elaborate ground curves,
themselves generated by ad hoc recipes (hyperbolisms, Cremona,
etc.). Remind here the biotope of curves imagined by Russian
scholars: Viro's beaver, horse, Shustin's medusa, etc (see our
Fig.\,\ref{SIMPLIFIED-TABLE_gurus:fig}). It is here that
patchworking  degenerates, or rather ramifies, to an artwork
 difficult to implement systematically, yet offering highly
 arborescent possibility for the artistically inclined worker.

7. If Viro's method fails to construct all curves, this means
there are smooth curves in remote mysterious chambers past the
discriminant. This is somehow akin to dark energy/matter hard to
interact with (within the standard model). Yet, in down-to-Earth
reality each smooth curve can degenerate toward the discriminant,
and a priori along several faces of it, so that we get a highly
singular curve with controlled singularities, and to which Viro's
method applies. Admittedly, all this is somewhat ill-posed, but
perhaps there is a reasonable way to claim-and-prove that Viro's
method is omnipotent, i.e. able to construct all smooth curves as
perturbation of a suitable curve with controlled singularities
(the {\it protozoan\/} so-to-speak). The latter has not to be same
throughout the hyperspace of all curves, and rather one expects
the presence of several gurus (=prototypical curves) required to
explore the full universe. The situation is somewhat akin to a
universe with several big-bangs with overlapping zone of
influences.

At least, the following quantitative problem seems senseful:

\begin{prob}
For a fixed degree $m$, what is the least number $\pi(m)$ of
protozoans (i.e., points of the discriminant) requested, so that
arbitrarily small neighborhoods of those
%%%intersect
overlap all chambers past the discriminant (or at least all
isotopy type of curves). In particular, each smooth curve is
rigid-isotopic to a small perturbation of one of the protozoan.
\end{prob}

This magnitude $\pi(m)$ admits alas several variants depending on
whether we restrict attention to
%the creationism of
$M$-curves, or admit all curves in the competition. In the
$M$-context we denote it $\Pi(m)$. For instance by Viro's
revisiting of the Harnack-Hilbert-Rohn-Gudkov theory we have
$\Pi(6)=1$, since all $M$-sextics arise as small perturbation of
the tri-ellipse. We guess (cf. also one of Viro's text) that all
sextics are small perturbation of the tri-ellipse safe the empty
sextic (this must follow rather easily from Nikulin's theory). In
that case $\pi(6)=2$, using as other protozoan the curve
$x^6+y^6=0$ with an isolated real point, but six linear branches
over $\CC$.

Intuitively one could expect that the geometry past the
discriminant is so intermingled, or  that Viro's method is so
versatile, that those (protozoan) numbers grow quite slowly in
function of $m$, and constitute so to speak black-holes governing
a whole galactic ama. Maybe $\pi(m)$ and $\Pi(m)$ even grows only
linearly in $m$. At the opposite extreme, one may speculate that
Hilbert's problem is so messy that even under this condensed
viewpoint there is an exponential growth of protozoans when the
degree $m$ increases.

All this is interesting yet one would like in a more narrow-minded
and stubborn fashion first fix the case of Hilbert's 16th in
degree $8$. The work for this is still immense, and decomposable
in  the following great lines:

1. Ensure (or perhaps refute?) that Viro's dissipation of $X_{21}$
is complete so as to rule out the option of creating new bosons by
the most rudimentary protozoan (namely the quadri-ellipse).

2. Prove that the Fiedler-Viro oddity obstruction is right, and
then prove (or disprove) Viro's sporadic obstructions as well as
Shustin's obstruction of subnested $M$-schemes with outer ovals.

3. Understand the link theory of Orevkov and the two resulting
obstructions.

4. Hope to use total reality or the method of the deepest
penetration as a way to unify and ideally to discover new
prohibitions.

5. Try to implement the deglutition-duality in order to link the
destiny of both subnested bosons.

6. If there is still some hope to construct new $M$-schemes and if
Viro's simplest method seems to have reached its limits, then try
\`a la Viro-Shustin to flexibilize the whole method by free-hand
tracing some protozoans creating new curves. Here there are
several ramifications: either via decomposing curves  of all
possible orders splitting 8 ($4+4$, $3+5$, $2+6$, $1+7$), or via
prescribed singularities of multiplicity splitting 8 too, like
$4+4$ (Viro's quadri-ellipse, Shustin's medusa), $3+5$ (Viro's
beaver and horse), $2+6$, etc. As we already saw, it is usually an
easy matter to discover qualitative configurations leading to
certain bosons, yet it is another ``paire de manche'' to ensure
algebraicity of the construction.

\subsection{Some new artwork hybridizing Viro and Shustin
(gorillas, yetis, etc)}

[19.09.13] It is evident that the number of ideas susceptible to
make progress the problem is fairly enormous, and one needs some
clairvoyance to find the right path to the goal. Fortunately, bad
geometers like to waste their time in this mess. For instance we
may wonder if there is any curve hybridizing Viro's mandarine with
Shustin's medusa, or to speak more concretely with one singularity
$X_{21}$ (quadri-contact) and one of type $Z_{15}$ (tri-contact
plus one crossing, alias the candelabrum). We cannot remember to
have tried this idea already, so we explore it anew. Of course a
reasonable attitude to have in this Hilbert problem, is to not
fear to repeat oneself since one can easily make mistakes leading
to erroneous conclusions or miss a combinatorial possibility in
the arborescences of the method.

After some trials on how to
%recombine
combine $X_{21}$ with $J_{15}$ (under
 the obvious constraint that nothing more must traverse the line
joining both singularities), we arrive at the gorilla curve
depicted below. This results from a search guided by the
desideratum of landing in the bosonic strip (one outer oval),
hence also preferring the $X_{21}$-dissipation of type E, as those
leave precisely hanging out one lune only (instead of the two
involved in type C smoothings). For a suitable quantization of the
gorilla (i.e. materialization of the 2 quantum ovals), we nearly
get the boson $1\frac{1}{1}\frac{18}{1}$ provided we could employ
the patch E(0,9,0). This is alas not available in Viro's theory
(but as far as we know not prohibited too).  Another smoothing
(said to be {\it external\/}, as the trunk of the candelabrum
connects with the extern branch) of our quantized gorilla leads to
a curve corrupting B\'ezout (saturation of the deep nest). In
conclusion,  the quantum oval cannot bubble out in the inner loop
as we did. In contrast it looks permissible to let it grow in the
outer loop, but even this violates B\'ezout as shown by the
suitable (external) smoothing with both a nest of depth 3 and one
of depth 2.

\begin{figure}[h]\Figskip
%\vskip-1.2cm\penalty0
%\centering
\hskip-2.7cm\penalty0
\epsfig{figure=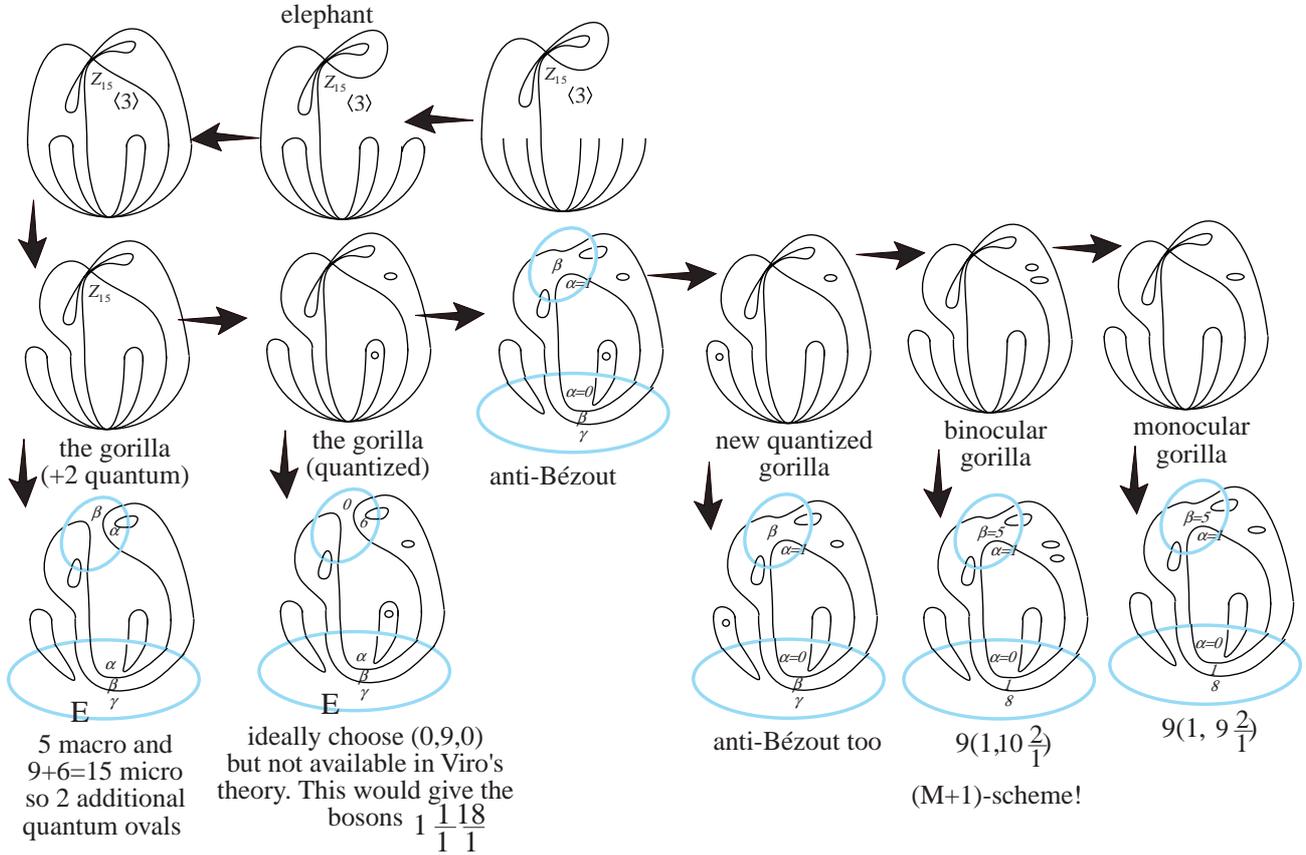,width=172mm} \captionskipAG
  \caption{\label{ViroDEGREE8_hybrid:fig}%
  Hybrid of Viro-Shustin}
\figskip
\end{figure}

So we finally opt for the {\it binocular gorilla\/}. Its depicted
smoothing is an $(M+1)$-curve violating Harnack! So in reality the
gorilla can accept only one eye ({\it monocular gorilla\/}), and
we made a mistake at the beginning when evaluating the number of
quanta by using a suboptimal smoothing.  Alas, the resulting
(monocular) smoothing yields the symbol $9(1,9\frac{2}{1})$, which
corrupts Gudkov periodicity (forcing the number of big eggs to be
$2 \pmod 4$). It is clear that wherever the quantum oval is
%%materialized
localized the number of big eggs will  always be 9 or eventually 8
if the quantum oval is ejected outside or injected inside, but
never $2 \pmod 4$. It follows the:

\begin{lemma}
Gudkov periodicity forbids any singular octic to be modelled over
the ground shape of the gorilla curve. Roughly put, the gorilla is
not sufficiently civilized to be algebraizable.
\end{lemma}

In fact the gorilla curve may also be imagined as a
carnivore-plant two of which protuberances are French-kissing in
the crepuscule. If preferring ornithology, you can also imagine a
mother bird feeding %%% CHECKED DICO ALIMENTER
its progeniture by direct transfer in the gullet (throat).
%%% CHECKED DICO gosier = throat
  With this view it is simple to imagine
variant of the gorilla. First we got the external kiss, but this
seems involving as optimal smoothing of the bottom
$X_{21}$-singularity the A-patch, which fails being an $M$-patch
(essentially by Arnold). Hence we considered rather the curve we
call the {\it yeti\/} (a nordic avatar of the gorilla).
Considering the dissipation depicted below it we see that there is
place for one quantum oval. This can be materialized either on the
right or on the left getting so the right- and left-yeti. From the
right versions we get two bosons while the left avatar yields
schemes due to Viro. Hence yeti-right is a good candidate to
create bosons, but before getting too excited we shall impose it
some harder resistance tests by evaluating on this architecture
all possible dissipation known in the Viro/Korchagin catalogues.
Besides, it may be of interest to test also our (Gabard's)
dissipation of type-C1, not mentioned by Viro. However using those
Gabard's patches yields corruption of Fiedler and Viro-regular.
More frankly, using Viro's class C3 we get likewise corruptions of
Viro's oddity law. So if the Fiedler-Viro theorem is right, our
right-yeti is jeopardized and so is the corresponding stratagem to
get the new bosons $b4$ and $b9$.

\begin{figure}[h]\Figskip
%\vskip-1.2cm\penalty0
%\centering
\hskip-2.7cm\penalty0
\epsfig{figure=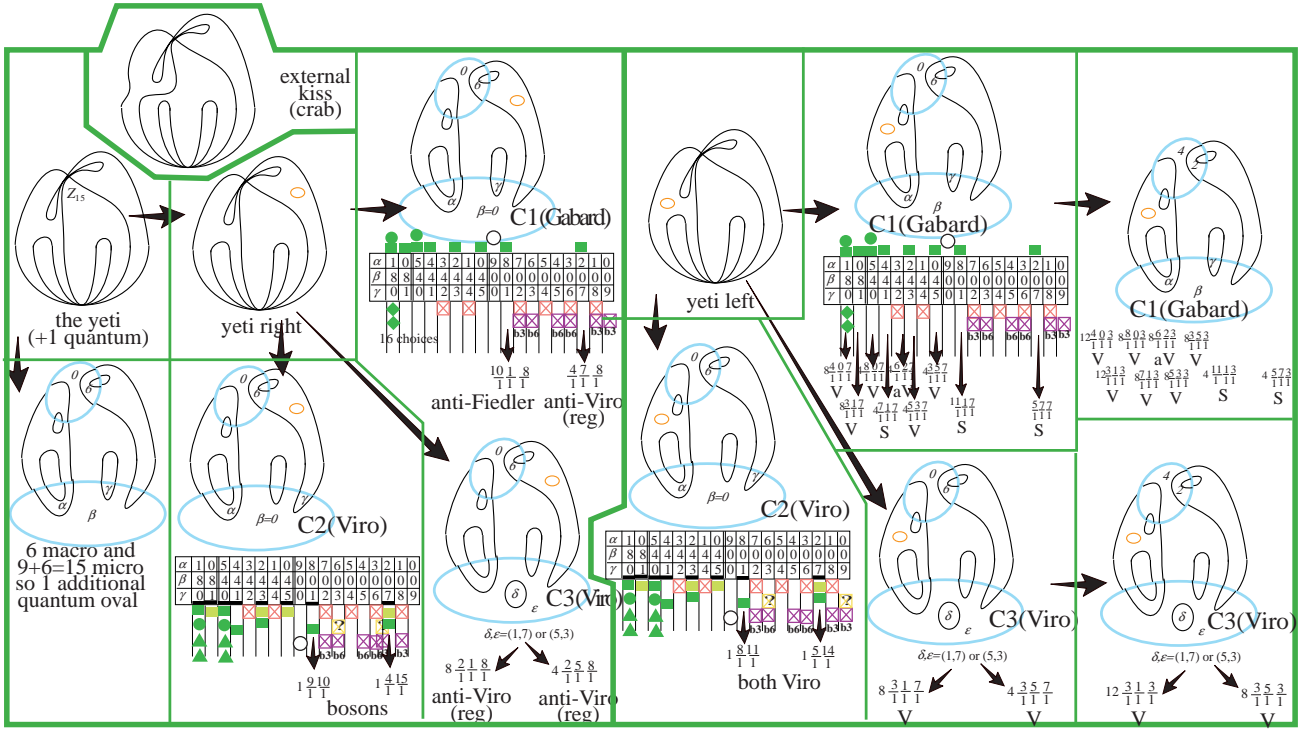,width=172mm} \captionskipAG
  \caption{\label{ViroDEGREE8_hybrid2:fig}%
  The yeti: another hybrid of Viro-Shustin}
\figskip
\end{figure}

It remains now to investigate whether the left-yeti is a more
respectable species supporting the patchworking
test-de-resistance.

{\tiny

(Skip this messy paragraph of dubious philosophy.)---[19.09.13,
23h51] Principles of relativity (\`a la Einstein) are only
required in physics where the basic concepts (forces, mass,
acceleration, speed,   matter, energy, time, physical space) are
ill-defined, perhaps because the are Gods creation instead of our
owns. More dramatically, since teenaged, we believe that reality
is only an appearance, without sound conceptual foundation except
the divinity. In contradistinction, the mathematical world (be it
human or divine Sch\"opfung) is perfectly well-defined and
therefore absolutist. No relativism is required to arrange the
intrinsic misconceptions.
 Each well-formulated problem reduces to a
yes or no answer\footnote{Compare optionally Steve Smale's list of
problem as a palish avatar of Hilbert's own, modulo  the
Gottschalk conjecture ca. 1958.}, apart maybe the so-called
undecidable problem \`a la G\"odel, hopefully reducible to a
matter of semantical misconception. In the real geometrical world
(like Hilbert's 16 th, etc.) it is evident that any reasonable
question receives a reasonable answer in finite time. Hilbert'
16th itself, after a small combinatorial cleaning, receives itself
a yes or no binary treatment, since there a re only finitely many
possible schemes by virtue of Harnack's bound, and evident
combinatorics. The puzzling issue, however, is that the shortness
of the question involves usually a very alembicated answer.
However there is then in principle algorithms of rationalization
allowing one to trivialize his long quest to a short explanation,
yet of a violent nature since it ignores the whole random
exploration process requested to find the solution by lucky
stroke, or by natural selection which is nearly synonym as
inefficient as it is. Otherwise we would since  the Jurassic era
already be immortal!

}

[20.09.13] Indeed when plugging in the left-yeti the C1-patches of
Gabard, or the C3-patches of Viro we get throughout respectable
schemes due either to Viro or Shustin. Hence the left-yeti seems
to pass the exam-test of resistance under patching. Of course,
this does not alone imply algebraicity of the left-yeti, but may
give supporting evidence for this. Unfortunately, the left-yeti is
conservative in the sense that it does not produce {\it new\/}
schemes.

Maybe there is still other variants of gorillas, yetis, etc.
Further there exist maybe also variants of Shustin's medusa, as
say the curve depicted below where both arms of the medusa are
merged together. Call it the {\it octopus\/}. This curve is
composed of 2 circuits and we believed a long time ago this being
an obstacle toward Harnack-maximality. Let us abort this prejudice
and explore the setting more liberally. So we smooth the
configuration in the optimal way and count the number of extra
quantum ovals required to reach Harnack-maximality. Here we find 4
quantum ovals. By B\'ezout those cannot emerge inside of the 4
loops emanating from both singularities. Using the smoothings
along the exterior branches, we see that Shustin's obstruction
forces at least one quantum oval to be outside. Also the interior
smoothing and Gudkov periodicity shows that exactly one quantum
oval must be outside. Then looking again at the interior smoothing
we see that the quantum ovals cannot be in the ventricle of the
octopus, and so have to be essentially as on our picture
(binocular octopus). When smoothed (interiorly) this produce the
boson $1\frac{4}{1}\frac{15}{1}$. When smoothed exteriorly we get
17 big eggs, violating thereby Gudkov periodicity. Hence we are
faced a serious dilemma as we like to arrange Gudkov periodicity
on both panels of interior and exterior smoothing yielding
respectively binested schemes (forced to have one outer oval mod
4) and subnested schemes (forced to have 2 big eggs mod 4). It
seems that there is no solution of compromise arranging both
relations in one stroke. In fact even without appealing to
Shustin's obstructions, it seems that there is no positioning of
the quantum ovals on the octopus so that both the interior and
exterior smoothings verify Gudkov periodicity. For instance we may
drag the outer quantum inside the ventricle of the octopus. When
smoothing interiorly then the situation is unchanged, but under
the exterior smoothing the situation has not been improved.

\begin{figure}[h]\Figskip
%\vskip-1.2cm\penalty0
%\centering
\hskip-2.7cm\penalty0
\epsfig{figure=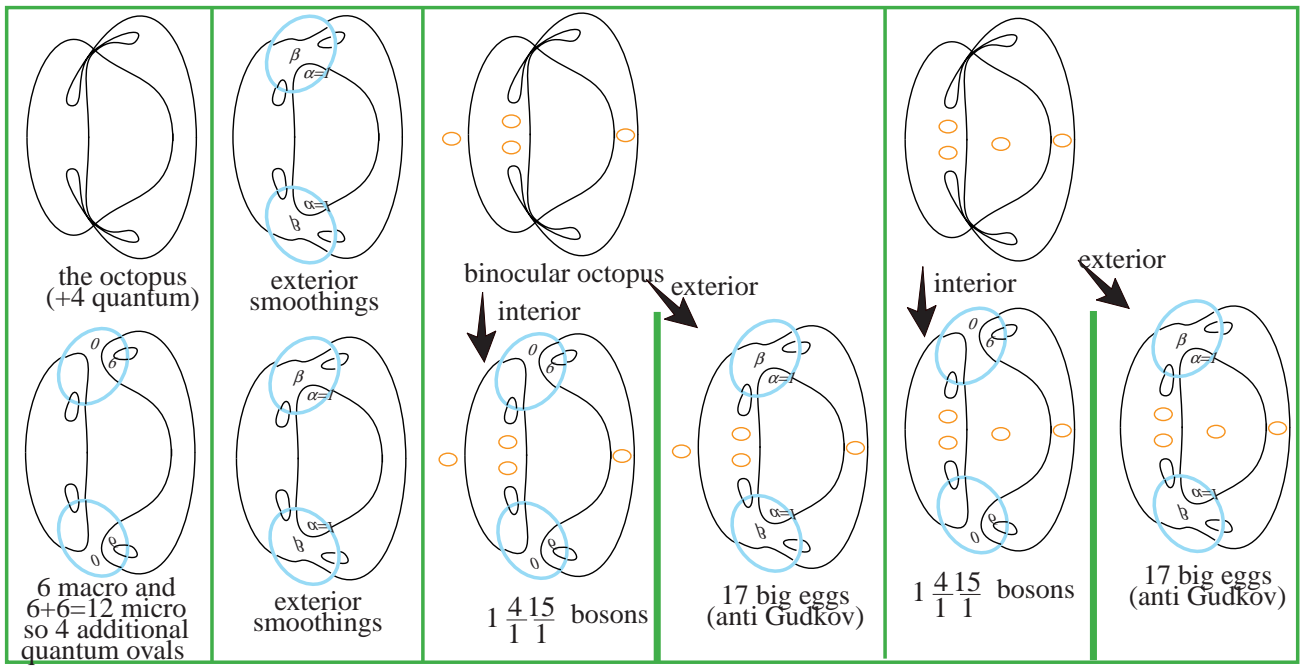,width=172mm} \captionskipAG
  \caption{\label{ViroDEGREE8_hybrid3:fig}%
  Octopus as a variant of Shustin's medusa: one boson
  is created, with anti-Gudkov dual particles}
\figskip
\end{figure}

Hence we hope to have proven:

\begin{lemma}
There is no algebraic curve whose topology is that of the octopus
plus $4$ quantum ovals whatever their location.
\end{lemma}

Next at the very beginning of our idea to recombine Viro's
mandarine (quadri-ellipse) with Shustin's medusa we traced a curve
(the {\it elephant}) which we neglected to consider more seriously
as it seemed to have two outer ovals  driving us outside the
bosonic strip (=binested with one outer oval). This could be
remedied by injecting a quantum oval into the loop, yet there are
2 objections against this. First, our curve would be trinested.
Second, applying Viro's C3-dissipation gives a quadri-nested curve
violating the saturation principle of the quadri-bifolium
$\frac{1}{1}\frac{1}{1}\frac{1}{1}\frac{1}{1}$. Next we find a
variant of the yeti, which we call the cobra. Probably its status
is essentially the same as the yeti. Next we can trace the full
zoo of animals encountered in a safari-tour: rhinoceros,
hippopotamus, scorpion, crab.
%%%CHECKED IN DICO= crabe
%%

\begin{figure}[h]\Figskip
%\vskip-1.2cm\penalty0
%\centering
\hskip-1.7cm\penalty0
\epsfig{figure=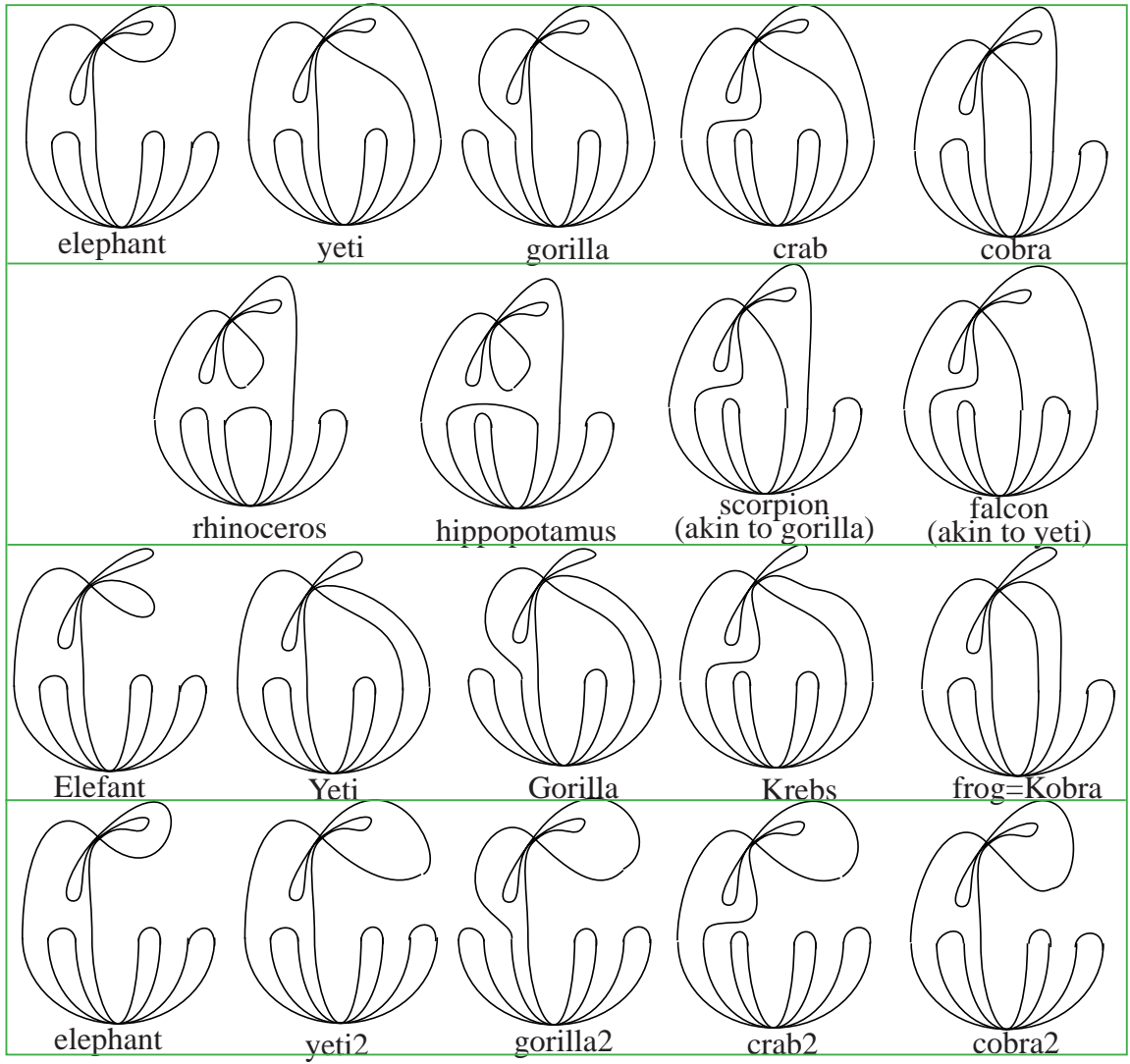,width=152mm} \captionskipAG
  \caption{\label{ViroDEGREE8_hybrid4:fig}%
  A zoo of animalistic curves}
\figskip
\end{figure}

[21.09.13] From all those animals which one is best suited to
revolutionize Hilbert's 16th problem in degree $m=8$? As we saw
the yeti was a good candidate but alas its most interesting
decoration creating new bosons  turned out to conflict with
Fiedler-Viro's oddity law. Of course we could reject the latter,
but this looks a bit cavalier. Still, we must confess that as yet
our brain never had the patience to study carefully the
Fiedler-Viro theorem. As to the cobra it is when smoothed
essentially isotopic to the yeti, especially if the collection of
patches C1 and C2 are symmetric as we could infer from our
interpretation of Viro's method. This would be more conceptually
explainable via the hypothesis of invariance under bending.

Let us be more systematic. First we see the elephant, but its
natural smoothing with a C-patch at the bottom create 2 outer
ovals so that we miss the bosonic strip. The natural parade is to
inject  a micro (quantum) oval in one of the loop emanating from
$X_{21}$, but then either by using the appropriate dual patch C1
or C2 as to fill the other loop,or just by taking C3, we arrive at
a quadri-nested scheme violating B\'ezout. Hence:

\begin{lemma}
There no chance for the elephant, even if it exists algebraically,
to realize the four binested bosons.
\end{lemma}

The destiny of the yeti was already discussed and there is no
chance for him to create boson unless the Fiedler-Viro obstruction
is wrong.

The case of the gorilla was already analyzed and conflicts with
Gudkov periodicity.

The case of the crab seems ruled out from entrance, because the
optimal smoothing is the A-type where only 8 micro-ovals can
appear. Despite this defect, one could hope still reaching
Harnack-maximality provided there are sufficiently many quantum
ovals (probably 2). But even if possible, one can argue that there
will be 2 outer ovals (at least), so failing to land in the
bosonic strip. It can then be counter-argued that one of the
quantum oval could appear in one of the loop.

So a more thorough analysis of the crab seems necessary. First, as
we said a smoothing of type A shows that to reach an $M$-curve
with 22 ovals the crab must be capable of receiving 2 additional
quantum ovals. (As usual quantum ovals are just ovals whose exact
location is not yet determined.) Next to drive the crab in the
bosonic region we force it to accept a quantum oval in one of the
loop. Now choosing on the top the external branch dissipation of
the candelabrum, and one the bottom a type~I smoothing (granting
its existence) we get would get a configuration violating
B\'ezout. Hence we suspect that the crab cannot be quantized as to
reach the bosonic strip, and therefore even if it existed it would
be useless to Hilbert's problem, except maybe in the subnested
case.

\begin{figure}[h]\Figskip
%\vskip-1.2cm\penalty0
%\centering
\hskip-2.7cm\penalty0
\epsfig{figure=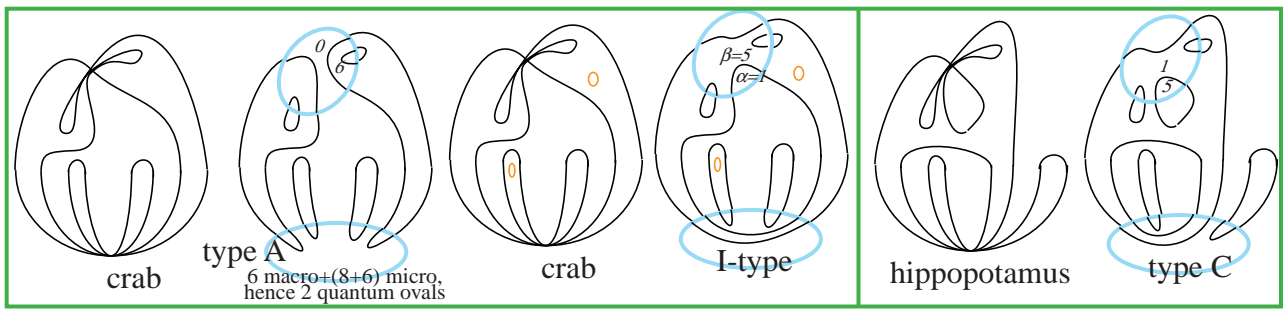,width=172mm} \captionskipAG
  \caption{\label{ViroDEGREE8_crab:fig}%
  The crab attempting to reach Harnack-maximality despite
  bad predispositions (followed by the hippopotamus, who
frankly corrupts B\'ezout)} \figskip
\end{figure}

After the crab we have the cobra, but the latter is much akin to
the yeti, and by virtue of the structural symmetry between C1- and
C2 patches (Gabard's belief, but check once by e-mail if Viro
agrees), both curves produce isotopic schemes. Hence nothing new
can be expected from the cobra, that the yeti not already
revealed. To remind the latter only produced boring schemes of
Viro-Shustin, when not conflicting with the Fiedler-Viro law.

Next it comes to the rhinoceros: this may be ruled out from the
scratch as the requested patch is of type F and only capable of 8
micro-ovals (spermatozoid droplets).

The hippopotamus deserves more respect as there is $M$-patches to
smooth the bottom singularity, namely those of type E. As to the
candelabrum it can be smoothed in the optimal way (compare Viro's
Fig.\,39, p.\,1112 in Viro 89/90
\cite{Viro_1989/90-Construction}), but the end-result involves two
subnest and so foils B\'ezout. Hence:

\begin{lemma}
There is no octic whose morphology is that of the hippopotamus. In
particular the latter will not aid us to advance Hilbert's  16th
problem on the qualitative theory of algebraic curves.
\end{lemma}

Our story continues with the scorpion, which is symmetric to the
gorilla, hence of no use; and idem for the falcon which is
symmetric to the yeti.

One question arises: was our safari tour in the zoo exhaustive, or
did we missed a species of special noteworthy-ness? Besides, some
principle of graphical elegance of algebraic curves often allied
to a principle of minimization (least effort law) seems to
corroborate the fact that the hippopotamus is ruled out. One may
first classify species along the number of connection linking both
singularities. For instance the yeti has 4 such connections, while
the hippopotamus only two. At this stage only our brain noted that
our picture of the hippopotamus is ruled out from entrance by
tracing the line through both quadruple points.

Trying to answer the above question one can create new animals by
surgery. For instance starting from the cobra and reconnecting
some braids we get the frog. Of course this alteration can be
operated on all animals listed, and so we get the 3rd row of
animals whose name is just translated in German (often just
amounting to a capitalization of the word). In particular we
rebaptize the frog as Kobra. As a first remark in those germanic
version of the animals there is always a smashed loop at the
candelabrum, so that the curve contains (at least) two circuits.
From earlier experience, we think this being a defect as somehow
Harnack's bound cannot then be optimized, but maybe we were a bit
prejudiced by a misconception. So let us start a naive browse
through the German bestiary.

Before doing this we see that the primary bestiary can undergo
another surgery amounting to cut the right arm of the animal
viewed as a carnivore plant, and this gives the series 2:
involving yeti2, gorilla2, etc. It may be observed that yeti2 and
cobra2 are nothing but the elephant. The other species looks
suboptimal as their best smoothing does not involve an $M$-patch.
As this stage it seems important to investigate more thoroughly
the elephant without prejudice about the two outer ovals.

Starting from the elephant and doing the optimal smoothing
depicted, involving the type C patch plus an internal bifurcation
on the candelabrum, we get a curve with 7 macro ovals, and 6+9
micro ovals so that Harnack's bound is alread attained with having
to introduce quantum ovals. In other word the elephant
configuration is already saturated (i.e. peasant enough that
nothing more can appear in the horizon). Somehow this is a bad new
as we hoped to use a quantum to kill one outer oval by injecting
stuff inside of it. (Keep in mind our intention to reach the
bosonic strip.) Though a lesser hot-spot, it seems still of
interest to investigate what schemes arise as progeniture of the
elephant. (One of our hope would be to get Shustin's last scheme
that we as yet never succeeded to construct.) On patchworking
``Gabard's'' series of patches C1 we get schemes violating Gudkov
periodicity. Of course the same outcome would result from using
Viro's more respectable patches of type C2. Hence:

\begin{lemma}
Gudkov periodicity (and probably nothing more elementary) impedes
any (singular) octic to acquire the morphology of an elephant.
\end{lemma}

\begin{figure}[h]\Figskip
%\vskip-1.2cm\penalty0
\centering
%%%\hskip-2.7cm\penalty0
\epsfig{figure=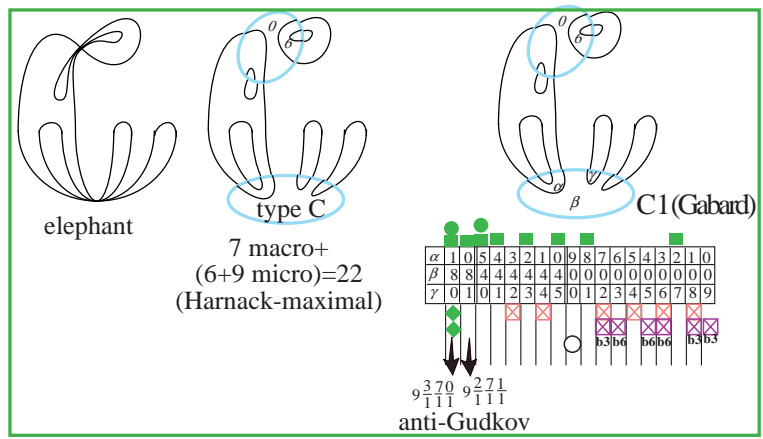,width=102mm}
\captionskipAG
  \caption{\label{ViroDEGREE8_elephant:fig}%
  Elephant trying to reach the bosonic strip despite
  bad predispositions} \figskip
\end{figure}

[22.09.13] {\it Optional digression}.---As a side remark to clean
at the occasion, one can start with any $M$-octic and a line
cutting it four times, and perform a karate move (smashing)
generated by a hyperbolism. This just amounts picturesquely to
hang on the curve like a dead medusa over a nail. Then a
$X_{21}$-singularity is created and one can dissipate it along the
usual Viro patches. It seems evident then, especially in case when
the line hits only one oval like a horseshoe pattern, that the new
curve will have 9 additional micro oval coming from the patch,
while the smashing deteriorating at most one oval. Hence it seems
clear that the new curve will have circa $22+9=31$ ovals,
overwhelming seriously Harnack. Presumably this paradox is
explained by the issue that the hyperbolism does not conserve the
degree to 8, but might increase it. Sorry for this loose idea, but
we just wrote it to not forget it, and in the hope to clarify it
at the occasion.

{\it Spruch der gut klingelt}.---Si on commence \`a s'enliser dans
les d\'etails arith\-m\'etiques, on ne comprend
%%%pas
plus la structure g\'eom\'etrique du cosmos.

Next, we were
%%%%%%%%%struck
sidetracked by the following idea. All animals depicted as yet
(models of qualitative octics) have the special feature of not
intersecting the line at infinity over the reals. We may thus
imagine more general (projective) curves drawn in the projective
plane (disc with boundary antipodically identified). Perhaps we
get then new animals susceptible of producing the new bosons. The
unfinished picture below tries to explore this idea, and must be
completed at the occasion.

\begin{figure}[h]\Figskip
%\vskip-1.2cm\penalty0
%\centering
\hskip-2.7cm\penalty0
\epsfig{figure=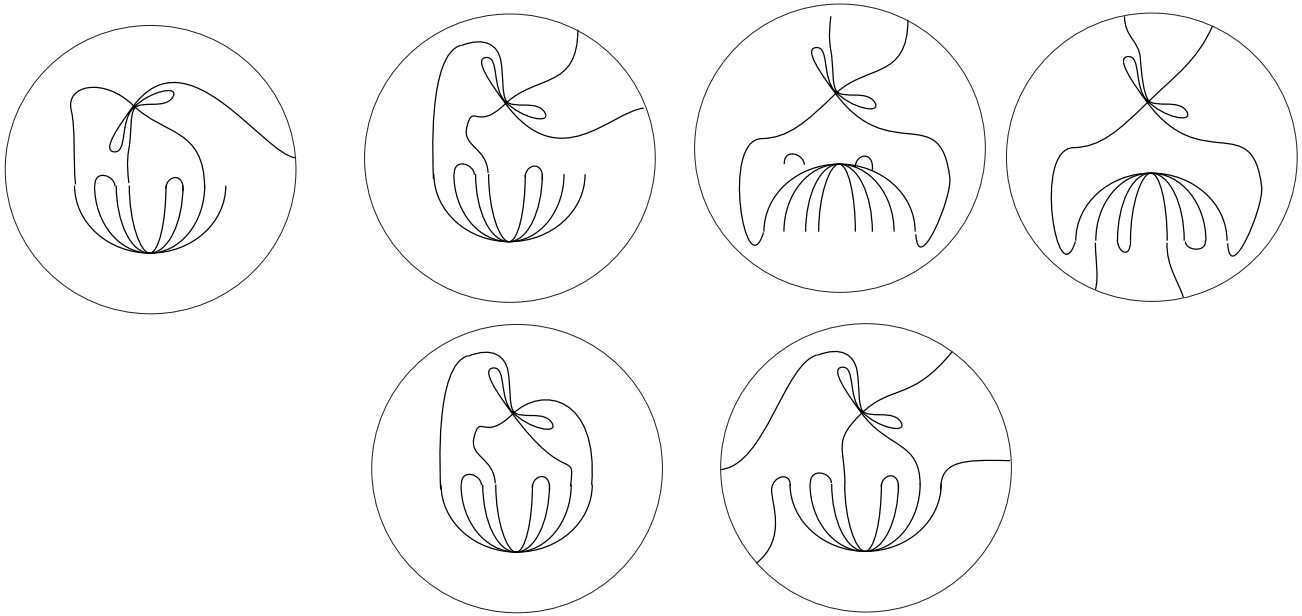,width=172mm} \captionskipAG
  \caption{\label{ViroDEGREE8_proj:fig}%
  Projective connections: a new menagerie of animals} \figskip
\end{figure}

Besides, now that we have understood the epicenter of Viro's patch
method (in its two decorations hyperbolism versus vibrational) it
seems also realist to adventure into new singularities type like a
quintuple flat point plus a triple point. By a {\it $k$-tuple flat
point} F$k$ we mean $k$ branches with 2nd order tangency; so for
instance $X_{21}=$F4.

The work then decomposes in two steps.

1. First, exploration of the dissipation theory of F5 by an
extension of Viro's two methods. In one decoration this merely
involves looking at one more step in Viro's vibrational process.

2. Second, enumeration of singular octics with a F5+F3 pair of
singularities.

\begin{figure}[h]\Figskip
%\vskip-1.2cm\penalty0
%\centering
\hskip-2.7cm\penalty0
\epsfig{figure=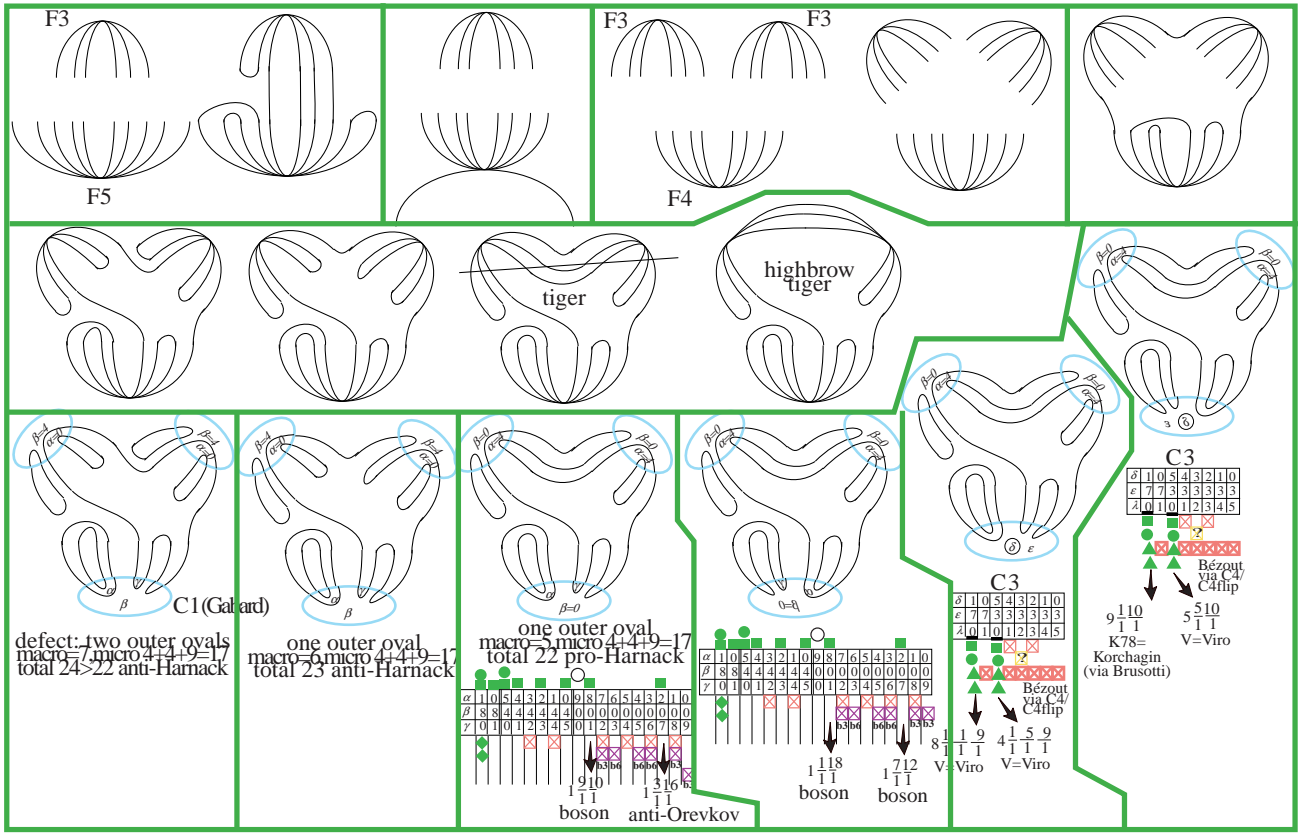,width=172mm} \captionskipAG
  \caption{\label{ViroDEGREE8_flat3+5:fig}%
  Singularities F3+F5, F3+F3+F4 including the tiger, as a powerful
  detergent of bosons, yet anti-B\'ezout after more mature thinking} \figskip
\end{figure}

Unfortunately, one sees quickly that the singularities F5 seems to
involve 10 intersections with a perturbed tangent to the ``south''
pole. Hence the singularity F5 seems ruled out for octics. Things
could be salvaged if the outer branch had inverted curvature, but
in reality on taking the tangent to the south pole (F5) we get
again a multiplicity intersection of at least 10, violating
B\'ezout. Hence:

\begin{lemma}
In the construction of octics (more genrally curve of order $m=2
\ell$), flat singularities {\rm F$k$} can have at most
multiplicity $k=4$ (resp. $k=\ell$).
\end{lemma}

Maybe one can investigate F4+F3+F3, i.e. curves with one flat
quadruple point and two triple points. After few attempts we found
a curve called the {\it tiger\/}, which respects Harnack's bound.
When smoothed in the clever way as to land in the bosonic strip
this yields the boson $1\frac{9}{1}\frac{10}{1}$, plus the
anti-Orevkov scheme $1\frac{3}{1}\frac{16}{1}$. Hence:

\begin{lemma}
Orevkov's link theoretic obstruction (if true at all\footnote{We
have no specific objection against Stepan, but only confess to
have not yet found the energy to check his proof.}) kills the
tiger, and so our
%%%dessein
plan
%%%CHECKED IN DICO
to get a new boson.
\end{lemma}

Besides, if we still believe in the tiger (price-to-pay=misthrust
Orevkov), then opting for the  vertically symmetrized patch we get
two other bosons. Hence the tiger looks very puissant at the
bosonic level, i.e. a good particles detector in CERN's jargon.
Further, gluing Viro's C3-patches yields no obstruction but
recovering standard schemes due to Viro's quadri-ellipse method.
In summary:

\begin{Scholium}
The tiger offers a real opportunity to win three new bosons in the
binested realm (all safe $1\frac{4}{1}\frac{15}{1}$), provided
Orevkov's obstruction of $b3=1\frac{3}{1}\frac{16}{1}$ is wrong.
However as we shall see,  there is a basic B\'ezout obstruction
working against the tiger.
\end{Scholium}

{\it Note}.---The reader may have noticed sooner than us, that the
tiger as it stands is intercepted ten times by the depicted line.
This is a violation against B\'ezout, yet a soft one since by
inflating the tiger's brain we may get the {\it highbrow tiger\/}
were this defect is remedied upon.

Further, flipping C3 along vertical-axis symmetry yields common
schemes due to K78=Korchagin 1978 (via Brusotti), and one scheme
due to Viro (quadri-ellipse).

The tiger seems to respect B\'ezout for lines, but what about
conics? The test involves passing a conic through the
singularities. From the five points through which one may pass a
conic we may amalgamate 2 pairs in the infinitely small so as to
offer tangency conditions. Imposing tangency at $X_{21}$ (the
quadruple point), and one more tangency at  the triple point, plus
a simple passage through the remaining triple point, we get a
multiplicity intersection of $8+6+3=14+3=17>16=2.8$, overwhelming
 B\'ezout. This seems to kill the tiger and comforts thereby
Orevkov's obstruction.

Of course the tiger is probably not an isolate species in the
class F3+F3+F4, and probably there is another curve leading to the
fourth (binested) boson. Yet this hypothetical curve will be
subsumed of course to the same B\'ezout obstruction just sketched.

In fact to lower this ``singular multiplicity'' one could replace
one triple point by a double one, or alternatively,  trade the
quadruple point $X_{21}$ for a candelabrum (where one of the four
branches is transverse). Then the multiplicity intersection of the
conic interpolating the singularities is only $7+6+3=16$ and
B\'ezout is respected. This brings us to our next picture
(Fig.\,\ref{ViroDEGREE8_candel3+3+4:fig}).

\subsection{C4+F3+F3: one candelabrum and two flat triple points}

{\tiny [23.09.13] Usually, mathematicians exposes their results,
but not the methods. Presumably, the reverse-engineering would be
at least as useful. A typical example is Viro's sporadic
obstruction.

}

{\tiny {\it Side-remark}.---Actually, Ahlfors extremal problem, is
probably more a Mittel zum Zweck (biased by the
Koebe-Carath\'eodory tradition) than an intrinsic feature of the
problem.

}

[23.09.13] On this picture
(Fig.\,\ref{ViroDEGREE8_candel3+3+4:fig}) we get several animals
and corresponding curves. Alas, apart from recovering Shustin's
last scheme $4\frac{5}{1}\frac{5}{1}\frac{5}{1}$, we were not able
(after a boring tedious  search) to reach the bosonic strip. We do
not know if this is caused by our incompetence, or an intrinsic
feature of this distribution of singularities.

Specifically, we found first the {\it lion\/} producing Shustin's
last scheme. Then we have several curves with anti-Gudkov
smoothings, hence not worth paying attention at (those include the
{\it panther\/}, {\it pig\/}, {\it cingallo\/}, etc.). Several
others do not attain Harnack-maximality (at least without
injecting extra ``quantum'' ovals); so for instance the {\it
gu\'epard\/} (=cheetah), {\it pork\/}, {\it cow\/}, etc.
%%%%CHECKED IN DICO
Browsing through the whole figure there is---apart from the
lion---only the {\it cat\/} which is Gudkov compatible, hence
susceptible to admit an algebraic model. Alas, it does not produce
any new boson but still Shustin's last scheme. We do not know alas
if there is an octic isotopic to the lion or the cat. But even if,
this does not impact tremendously upon Hilbert's 16th (except if
it may help to clarify Shustin's construction that we were as yet
unable to digest). So it seems reasonable to leave the question
asides for the moment.

\begin{figure}[h]\Figskip
%\vskip-1.2cm\penalty0
%\centering
\hskip-2.7cm\penalty0
\epsfig{figure=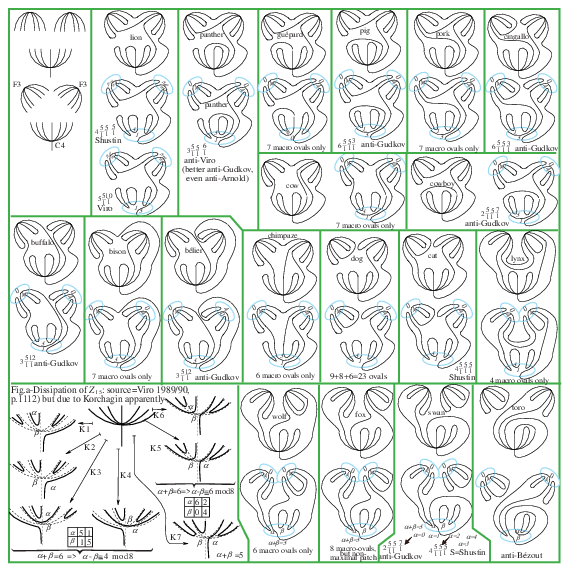,width=172mm}
\captionskipAG
  \caption{\label{ViroDEGREE8_candel3+3+4:fig}%
  Singularities  F3+F3+C4 (candelabrum) including the lion, etc.} \figskip
\end{figure}

And we continued our search, on a 2nd figure
(Fig.\,\ref{ViroDEGREE8_candel3+3+4B:fig}) yet not finding any
species worth of commentary. So our boring pictures are given
without any comments. So far so good but alas TeX is not happy
with little comments because then the
%%%mass
overflow of pictures overwhelms
%what he seems able to page-make on
%a single page.
its page-making aptitudes. So let us comment against our will.
First we have a {\it ch\`evre\/} (=coat), which produces 9
macro-ovals, but this stays okay because the suited patch permits
only 5 micro-ovals (compare the Viro-Korchagin catalogue
reproduced on our previous figure). But then there are apparently
no restriction on $\al,\be$, and thus we frequently collide
against Gudkov periodicity. In conclusion there should be no octic
taking the form of the goat.

Next we have the {\it sheep\/}, but this looses one oval over the
goat, and thus should not be able to reach Harnack-maximality.
Next, we imagined a {\it zebra\/}, which is however anti-B\'ezout
as it contains two nests of depth 3 and 2 respectively. The {\it
borsuk\/} (=blaireau=badger?) is for the same reason
anti-B\'ezout. The {\it dolphin\/} has only 5 macro-oval, and so
its smoothing severely fails to be an $M$-curve, except if we
would add quantum ovals. Next we have a {\it shark\/} with 7
macro-ovals, but this is still not enough. One can increase by
going to the {\it orc\/} (=orque in French=grampus?),
%%%CHECKED in DICO
but  its progeniture under smoothing is alas anti-Gudkov (even in
the simple form of Arnold). Next we have the {\it wale\/} whose
smoothing produces only 7 macro-ovals. This can improved by
choosing a better patch splitting apart the ``huge'' contorted
oval, and we get so an $M$-scheme, alas anti Gudkov periodicity.
(The latter forces in the binested case the number of outer oval
being $1 \pmod 4$.) From the wale there is an obvious
morphogenesis to the {\it Stier\/} (=German for bull), but its
smoothing is anti-Gudkov. By the way the more suited
split-smoothing would violate Harnack's bound.

\begin{figure}[h]\Figskip
%\vskip-1.2cm\penalty0
%\centering
\hskip-2.7cm\penalty0
\epsfig{figure=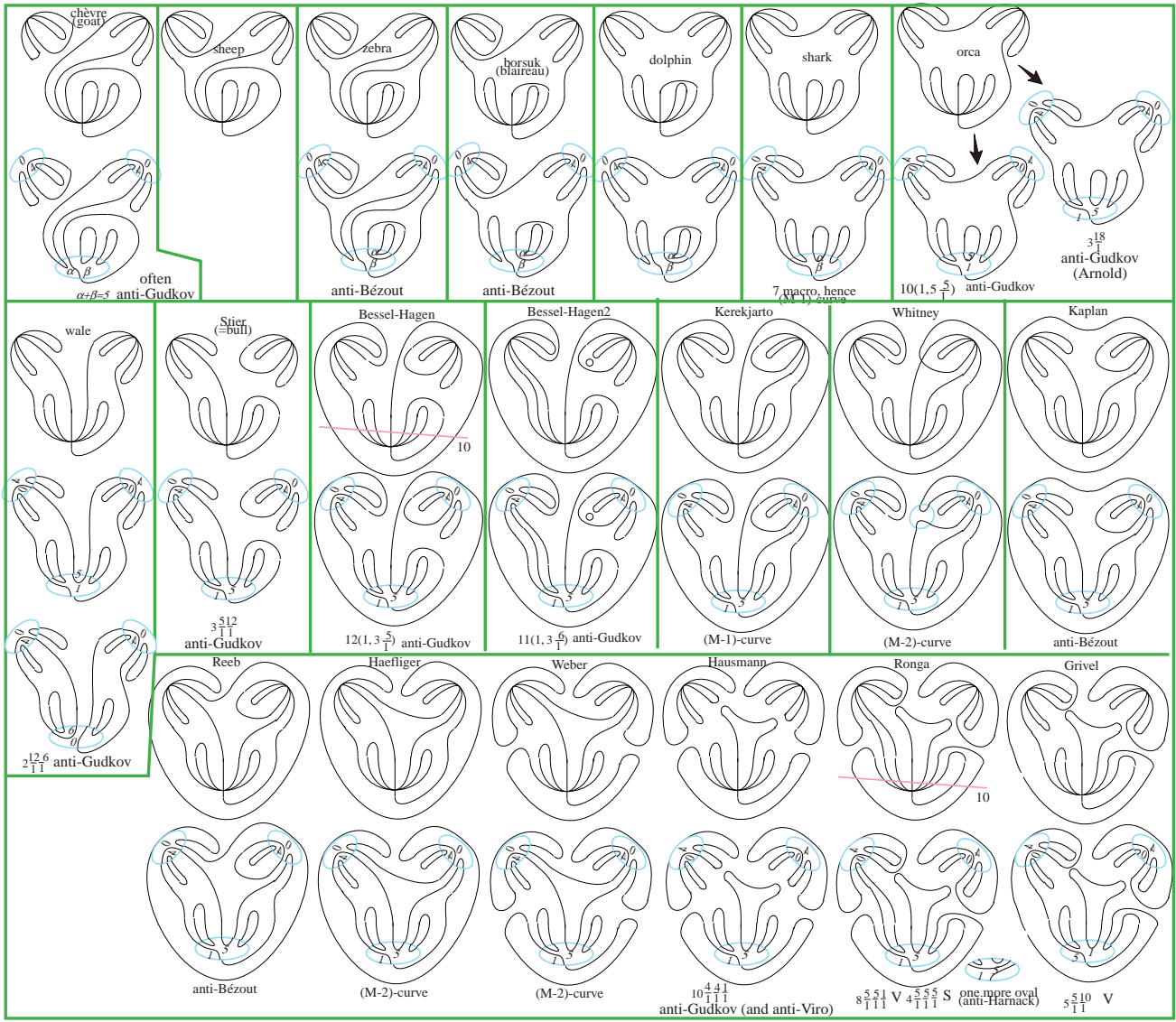,width=172mm}
\captionskipAG
  \caption{\label{ViroDEGREE8_candel3+3+4B:fig}%
  Singularities  F3+F3+C4 (candelabrum): from the goat to
  Grivel} \figskip
\end{figure}

Next being \`a cours de vocabulaire, we decided to opt for names
of famous geometers instead of animals. The basic idea is to
consider a curve like {\it  Bessel-Hagen\/} which is more
``claustrophobic'' or squat (trappu in French) with a branch
winding around the whole configuration before closing back to the
singularity. However its production is anti-Gudkov, and of course
there is also a line cutting the curve along 10 points. We
explored so a long list of curves termed after Kerekjarto,
Whitney, Kaplan, Reeb, etc. The sole interesting species is that
called Ronga, which produces Shustin's last scheme
$4\frac{5}{1}\frac{5}{1}\frac{5}{1}$. Of course our model of Ronga
is still intersectable in 10 points by a line. Further on choosing
another best suited smoothing (depicted in the margin) we can gain
one more oval, and this brings Ronga's curve outside the realm of
Harnack-maximality. Next we have Grivel's curve to which the same
token applies mutatis mutandis.

[24.09.13] Of course the canonical idea to work
%%%this out
more systematically is to start from the $M$-smoothing of the
candelabrum, and then connect the branch so as to reach
maximality. Hence, one may start from the Viro-Korchagin catalogue
of dissipations and then complete the curve. Doing so while
connecting the trunk of the candelabrum to the other side, we
found first a curve (called {\it Garfield}) violating two of
Shustin's prohibitions ($(1,14\frac{6}{1})$, $(1,18\frac{2}{1})$).
So:

\begin{lemma}
Either two of Shustin's prohibitions are wrong or there is no
octic curve isotopic to the Garfield. However it seems to us that
the singular version of Harnack's bound (due to either Harnack,
Klein or Hurwitz, who else?) easily expels the Garfield outside
the algebraic realm. Hence, Shustin is probably safe.
\end{lemma}

[Added 27.09.13.---One may also wonder about avatars of the
Garfield attacking the other three Shustin's obstructions.
Besides, one may also imagine another Garfield with ovals
quantized elsewhere as to get  the boson
$1\frac{1}{1}\frac{18}{1}$. Of course there is still the critique
of the singular Harnack bound, yet its seems still worth tracing
that curve, as {\it goret\/} on
Fig.\,\ref{ViroDEGREE8_candel3+3+4C:fig}. Of course its smoothing
turns to be anti-B\'ezout. So sorry for that stupid example.]

Another smoothing of the Garfield yields a more respectable
Hilbert's $M$-curve. Yet, this is certainly not enough evidence to
fight against Shustin. We shall soon give  an argument based on
%the genus of the complex curve
 the Harnack-Klein
bound preventing  the Garfield's existence.

{\footnotesize

[Added 27.09.13.---Actually, if we pass a conic tangent to both
fat branches of the Garfield at the point C4 and F3 while passing
simply through the other F3, we get $3.2+1+3.2+3=16=8.2$, and
B\'ezout is still happy. Sorry, this is a stupid remark, since by
construction we know that this distribution of singularity is
B\'ezout permissible.]

}

\begin{figure}[h]\Figskip
%\vskip-1.2cm\penalty0
%\centering
\hskip-2.7cm\penalty0
\epsfig{figure=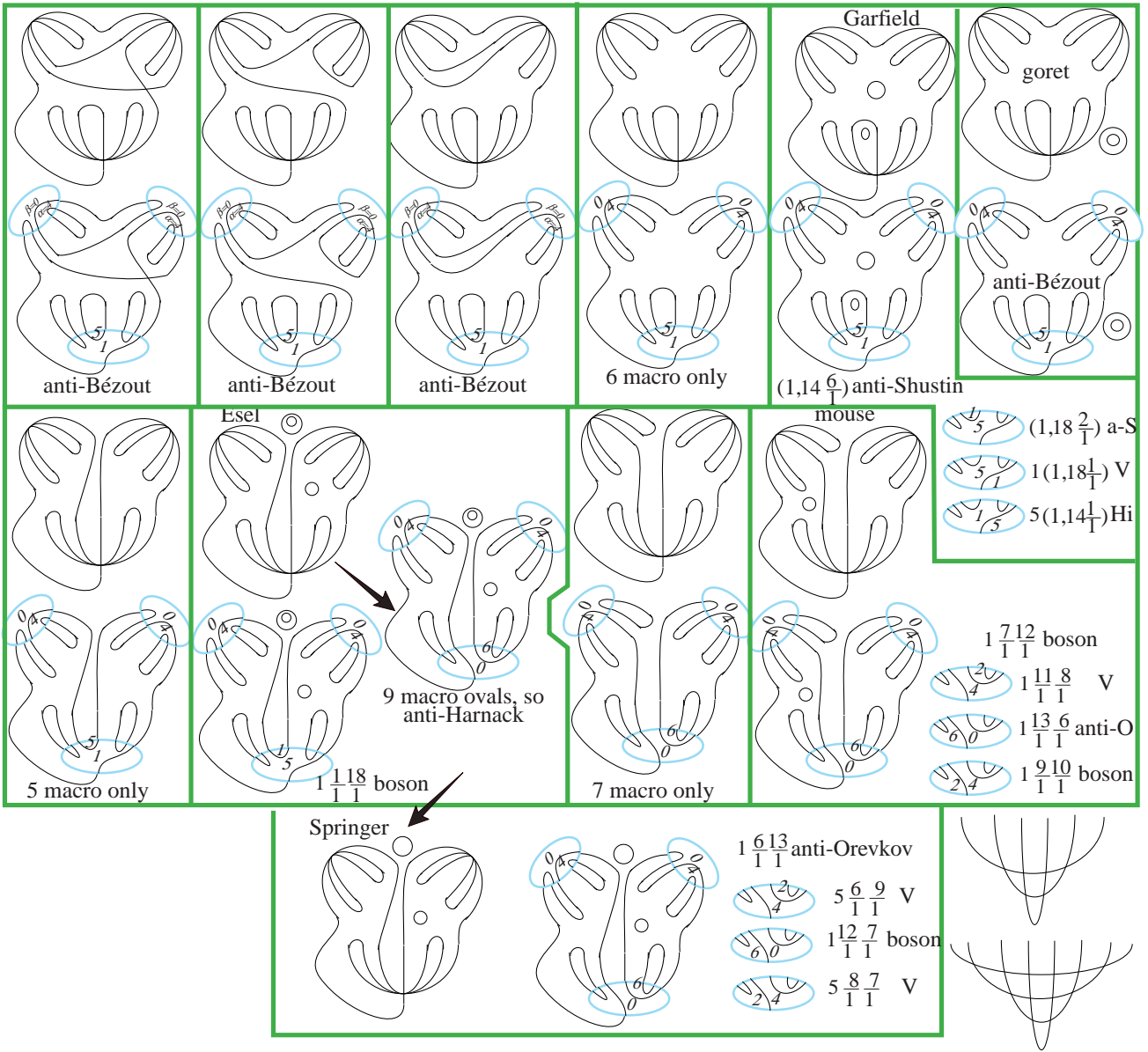,width=172mm}
\captionskipAG
  \caption{\label{ViroDEGREE8_candel3+3+4C:fig}%
  Singularities  F3+F3+C4 (candelabrum): Garfield, Esel, etc.} \figskip
\end{figure}

Next, stubborn as a mule,
%%%
%%FROM DICO: tetu comme un ane
%
we found a curve called the {\it Esel\/} (=\^ane=ass, or donkey)
which produces the boson $b1=1\frac{1}{1}\frac{18}{1}$. However
smoothing it differently violates Harnack's bound. Of course we
can just kill one quantum oval, to get the {\it Springer\/}. Its
most virulent smoothing  produces the scheme
$1\frac{6}{1}\frac{13}{1}$ violating Orevkov. Other patches give
two Viro schemes, plus the boson $1\frac{7}{1}\frac{12}{1}$. In
conclusion the Springer is only executed (killed) by Orevkov, but
we think that the singular Harnack bound also prohibits the
Springer.

Further,   the {\it mouse\/} creates two (new) bosons ($b7$ and
$b9$) but conflicts once more with Orevkov's $b6$. Hence:

\begin{lemma}
Orevkov is either false, or prohibits  the mouse.
\end{lemma}

Of course our mouse is just found by successive trials, especially
introduction of additional (quantum) ovals as to force
Harnack-maximality. Perhaps the mouse is readily ruled out by
Harnack's bound in the singular realm.

\subsection{The singular Harnack bound argument}

For this one must predict the (salaries) dumping effected on the
genus by singularities F3 and C4 (so-called {\it Mindestlohn\/} in
Germany). We would have preferred to skip this issue for the
moment, but let us improvise despite our unculture, since the
method seems a powerful tool of censorship against our
pseudo-counterexamples to Shustin, Orevkov.

If we imagine a sextic $C_6$ with two F3-points it will split
toward a tri-ellipse of genus $-2$ (imagine spherical
modifications, each lowering the genus by one). As  a smooth
sextic has genus 10,  {\it each {\rm F3} must drop the genus by
$6$}.

[Added 27.09.13].---Another method consists in perturbing the
singularity into an arrangement with normal crossings while
counting the number of double-points so created. For F3 we get 3
elliptical branches with a total of 2+4=6 nodes (cf.
Fig.\,\ref{ViroDEGREE8_candel3+3+4C:fig}). The same method for
F4=$X_{21}$ gives 2+4+6=12 nodes in accordance with the result
that one may derive by the first ``genus'' method.

\begin{figure}[h]\Figskip
%\vskip-1.2cm\penalty0
%\centering
\hskip-0.7cm\penalty0
\epsfig{figure=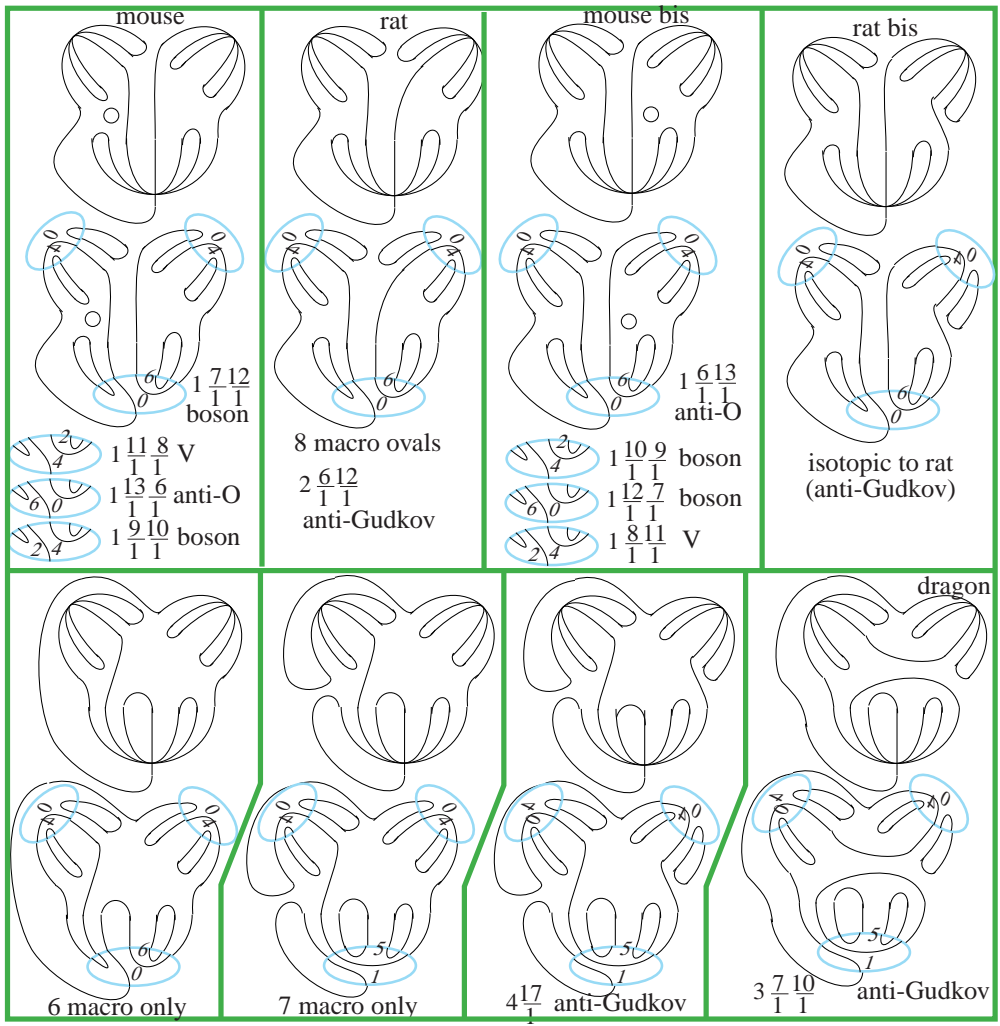,width=132mm}
\captionskipAG
  \caption{\label{ViroDEGREE8_candel3+3+4D:fig}%
  Singularities  F3+F3+C4 (candelabrum): rat, etc.} \figskip
\end{figure}

As to the candelabrum C4, it may be conceived as F3 plus a
transverse branch. Hence we have  3 additional crossings, each
eating one unity to the genus. So C4 decreases the genus by 6+3=9.
This count is compatible with Shustin's medusa
(Fig.\,\ref{SIMPLIFIED-TABLE_gurus:fig}) which has two
candelabrums C4, hence genus $21-18=3$, while totalizing precisely
4 circuits in accordance with Harnack's bound. In sum, our
distribution of singularities 2.F3+C4 drops the genus of a smooth
octic ($g=21$) to $21-12-9=0$. Thus our mouse has too many
circuits. Using curvature conventions along branches,  the mouse
consists actually of 3 circuits; yet, already without them
Harnack's bound is violated.

The same argument applies {\it a fortiori\/} to the Garfield (2
quantum ovals), which by curvature conventions even raises to 5
circuits. It applies also to the Springer (2 quantum ovals).

Then we transform the mouse to a {\it rat\/}
%%%
%%%%CHECKED IN DICO rat in French = the same in English
%
producing more ovals when smoothed. Alas the resulting $M$-scheme
violates Gudkov periodicity even in the simple version of Arnold.
Note en passant that the quantum oval of the mouse could as well
have appeared in the other half of the curve (cf. mouse bis). This
produces exactly the same collection of four $M$-schemes modulo a
shuffle. If our above genus dropping count is correct, we are not
allowed to add quantum ovals, and the game becomes fairly rigid,
in the sense that it becomes hard to land in the bosonic strip.

So the problem seems to be: is it possible to interconnect the
branches of two singularities of type F3  plus one of type C4 as
to get a single circuit (Harnack's bound for the singular curve)
while simultaneously arranging an $M$-scheme in the bosonic strip,
i.e. binested with one outer oval. This is tantamount the Gudkov
symbol being $1\frac{x}{1}\frac{y}{1}$ with $x+y=19$ (w.l.o.g.
$x\le y$).

\subsection{Sidetracked to C4+C4: two candelabrums}

[24.09.13] Albeit our search on 2.F3+C4 was far from systematic,
there is maybe more freedom when composing two candelabrums
(C4+C4) like in Shustin's medusa. This gives us the following
picture (Fig.\,\ref{ViroDEGREE8_candel4+4:fig}).

\begin{figure}[h]\Figskip
%\vskip-1.2cm\penalty0
%\centering
\hskip-1.7cm\penalty0
\epsfig{figure=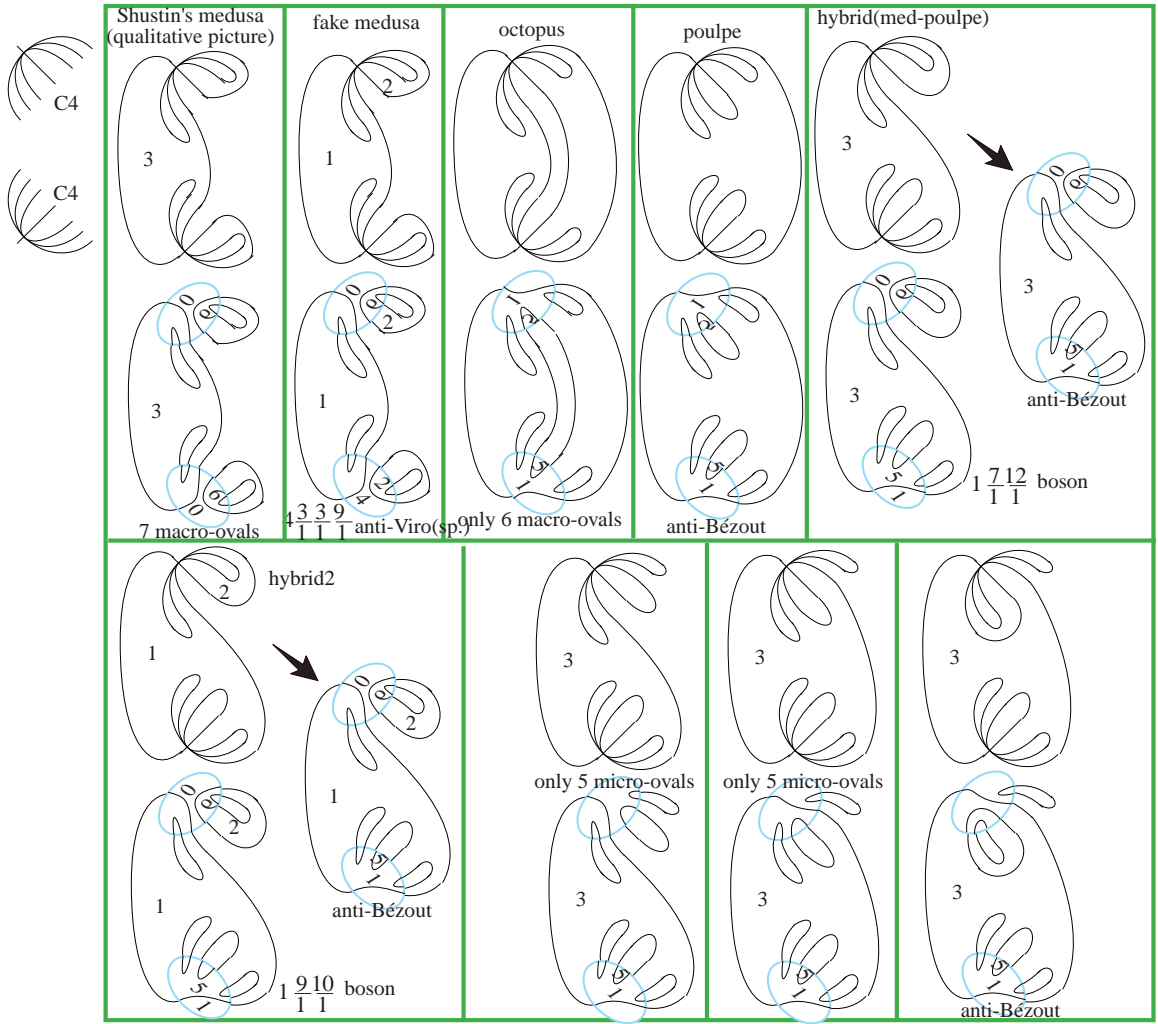,width=152mm}
\captionskipAG
  \caption{\label{ViroDEGREE8_candel4+4:fig}%
    C4+C4 (two candelabrums): bosons  created but
    anti-B\'ezout cousins} \figskip
\end{figure}

Here we start with a distribution of two candelabrums  with four
branches (notation C4). A first natural way to connect them is
like in Shustin's medusa, whose natural smoothing has 7
macro-ovals. As $6+6=12$ micro-ovals are given by dissipation
theory, we can add 3 quantum ovals. One adds them traditionally in
the {\it core\/} like on the medusa picture (i.e., the region
limitrophe to both candelabrums). As a variant one may imagine a
{\it fake-medusa\/} where only one quantum oval is centrally
placed, but two delocalized in the double-loop (see the {\it
fake-medusa\/} picture). Then it is a simple matter to arrange the
gluing patches as to
%obtain
offend Viro's most cavalier sporadic obstruction, namely
$4\frac{3}{1}\frac{3}{1}\frac{9}{1}$. Hence:

\begin{lemma}
Either Viro's sporadic obstruction
($4\frac{3}{1}\frac{3}{1}\frac{9}{1}$) is false or it kills the
fake-medusa (of Fig.\,\ref{ViroDEGREE8_candel4+4:fig}). As we
shall see later a simple trapping argument \`a la
Poincar\'e-Bendixson combined with B\'ezout rather corroborates
this 2nd alternative (alas without proving Viro's obstruction).
\end{lemma}

It seems of interest to
%%%look at
inspect the full d\'eploiment of this fake-medusa as it seems to
corrupt other Viro sporadic obstructions. This deserves a separate
plate (Fig.\,\ref{ViroDEGREE8_fake-medusa:fig}).

\begin{figure}[h]\Figskip
%\vskip-1.2cm\penalty0
%\centering
\hskip-1.7cm\penalty0
\epsfig{figure=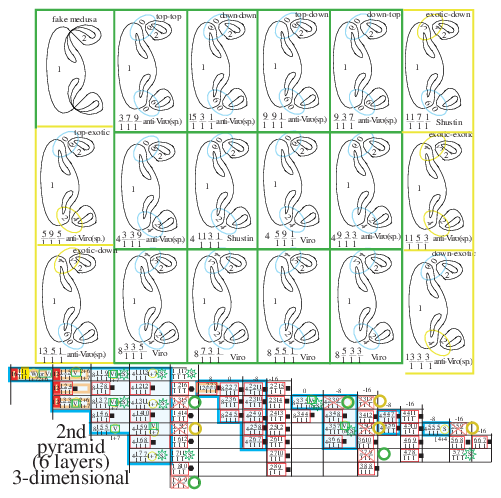,width=152mm}
\captionskipAG
  \caption{\label{ViroDEGREE8_fake-medusa:fig}%
  Fake medusa and its full d\'eploiment: killing half of Viro sporadic} \figskip
\end{figure}

After the fake-medusa we have the {\it octopus\/}, arising by
conjunction of the tentacles of the medusa, and so look
structurally incapable to reach Harnack-maximality. The sole
deliverance could come from more spontaneous quantum ovals, but as
discussed earlier (Fig.\,\ref{ViroDEGREE8_hybrid3:fig}) it seems
impossible to produce bosons without conflicting (radioactively)
with Gudkov periodicity. Next we have a {\it poulpe} (=French for
Devil-fish),
%%%in my dictionary),
but this mutates anti-B\'ezout on gluing the appropriate patch
(two subnests). Next, we imagined a hybrid of the medusa and the
poulpe (called {\it med-poulpe\/}). This even produces an exciting
boson, but alas runs against B\'ezout when
 smoothed differently. Then we can imagine {\it hybrid2\/} which
produces another boson, but again along an illegal way foiling
B\'ezout. The next two curves involves the non-maximal patch (K7),
and so Harnack-maximality seems out of reach, except if more
quantum ovals (viz. four) are created.

\subsection{On the fake-medusa}

[25.09.13] Now we turn back to the project of enumerating all
smoothings of the fake-medusa.
Fig.\,\ref{ViroDEGREE8_fake-medusa:fig} shows than we can attack
several of Viro's sporadic obstructions from this single position.
%%(of the discriminant in the hyperspace of octics).
Precisely:

\begin{lemma}
If there is a fake-medusa, then all the following four sporadic
Viro obstructions are wrong:
%$$
$\frac{3}{1}\frac{7}{1}\frac{9}{1},\;
4\frac{3}{1}\frac{3}{1}\frac{9}{1},\;
\frac{1}{1}\frac{9}{1}\frac{9}{1},\;
\frac{1}{1}\frac{3}{1}\frac{15}{1}.
$
%
%
%$$
Conversely, it suffices one of those obstructions being true to
rule out the fake-medusa.
\end{lemma}

%Furthermore,
Further, on using dubious (yellow-colored) patches we can also
construct the four remaining Viro obstructions. Hence supposing
that the dissipation of the candelabrum was not fully explored (by
Korchagin-Viro), it could be that all the eight sporadic Viro
obstructions in Hilbert's 16th for $m=8$ are wrong.

So the point is twofold. First we see a  splitting of Viro's
sporadic obstructions in two classes of four, one
%being
more
suspect than the other. Besides, it
%%seems to us
is
%%quite
pretty remarkable that the fake-medusa seems  refuted only by
sporadic obstructions and nothing more tangible, like B\'ezout,
Gudkov, Viro's law of oddity, etc.

One naive idea to get a contradiction is to opt for the interior
smoothing without nesting (see Viro's Fig.\,39, p.\,1112 or
equivalently our K7 on Fig.\,\ref{ViroDEGREE8_candel3+3+4:fig}).
Apparently here Viro claims that for $\al+\be=5$ each values
$\al=0,1,2,3,4,5$ can be realized. One could imagine this
%causing
%troubles with
falsifying GKK-periodicity (GKK=Gudkov-Krakhnov-Kharlamov). In
reality, when gluing K7, we loose both a macro-oval and one
micro-oval and so land with an $(M-2)$-scheme only, where
periodicity is abolished. Using instead the patch K1 we may get an
$(M-1)$-scheme but GKK-periodicity seems respected. Maybe there is
a theological reason explaining that as the $M$-production of the
curve respects Gudkov, so must its $(M-1)$-descendance respects
GKK-periodicity. More experimentally, a quick browse---hopefully
exhaustive---through all patches listed by Viro(-Korchagin) does
not conflicted with GKK when glued inside the fake-medusa (compare
Fig.\,\ref{ViroDEGREE8_fake-medusaGKK:fig}).

\begin{figure}[h]\Figskip
%\vskip-1.2cm\penalty0
%\centering
\hskip-1.7cm\penalty0
\epsfig{figure=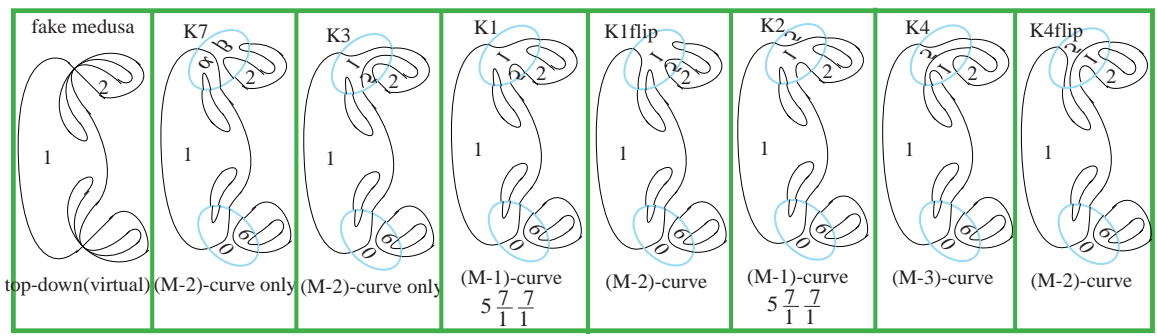,width=152mm}
\captionskipAG
  \caption{\label{ViroDEGREE8_fake-medusaGKK:fig}%
  Fake medusa versus GKK} \figskip
\end{figure}

Of course, it is a serious challenge to construct
algebro-geometrically the fake-medusa. A loose essay would be to
ape Shustin's construction of the original medusa, while trying to
deviate from it as soon as the occasion presents itself. This is
of course pure opportunism without tangible knowledge of the
terrain.

Another idea would be to search a direct B\'ezout obstruction on
the singular model prior to smoothing. The method,
by-now-standard, is to pass a conic
%%%curve  interpolating
through the singularities. First impose 2 tangencies along the fat
branches of both candelabrums. Besides,  impose another
anchor-point, typically inside one of the quantum oval. Counting
intersection we find $2.7+2=16=2.8$ so that B\'ezout's bound is
already attained, forbidding any further unassigned intersections.
Applying this recipe to Shustin's medusa, it seems plausible that
the interpolating conic through one of the 3 quantum ovals will
not intercept anymore the singular octic $C_8$ apart in the
assigned loci (cf.
Fig.\,\ref{ViroDEGREE8_fake-medusa-Bezout:fig}\,left). In
contrast, applying it to the fake-medusa (with one central oval
and two peripheral ones), while asking the conic to visit
%%%the inside  or just
a point situated on one of the two
%%noncentral
peripheral ovals then  the loop defines a trap \`a la
Jordan-Poincar\'e-Bendixson where the conic stays confined without
possible issue (Fluchtweg). So the conic-circuit is actually
forced to revisit the singular point of the $C_8$, but this is
intolerable (either for a conic to become a figure eight or the
intersection $C_2\cap C_8$ would then exceed B\'ezout). At any
rate it seems evident that:

\begin{lemma} \label{fake-medusa:lem}
A Poincar\'e-Bendixson trapping argument combined with B\'ezout
%prevents the
%existence of
excludes
%%%expels
the fake-medusa
%in
from the temple of algebraic-geometry.
\end{lemma}

\begin{figure}[h]\Figskip
%\vskip-1.2cm\penalty0
\centering
%%%%%\hskip-1.7cm\penalty0
\epsfig{figure=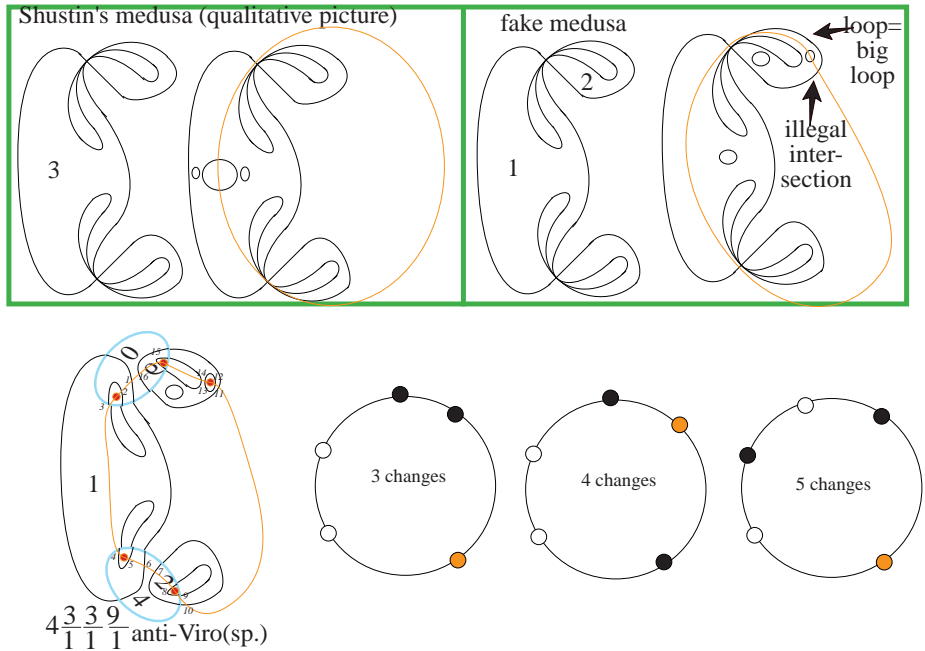,width=122mm}
\captionskipAG
  \caption{\label{ViroDEGREE8_fake-medusa-Bezout:fig}%
  %%%Direct anti-B\'ezout
  Trapping
  argument on the singular curve: genuine vs. fake
  medusas} \figskip
\end{figure}

This is
%to be
a fatal jeopardy of our essay to corrupt Viro's sporadic laws. If
the latter are true, one may
%now
wonder if an
elaboration of
%the lovely
this Poincar\'e-style argument
%just given
could assess those Viro obstructions. One method could be to smash
the smooth curve to a fake medusa by a suitable ``hyperbolism'',
broadly interpreted as a geometrical recipe preserving degree and
contracting certain lines. Yet this would probably involve deep
enumerative properties valid in some universal sense. This looks
of course completely out of reach.

[Added 27.09.13.---Further it is clear that the above lemma
(\ref{fake-medusa:lem}) applies as well to  two other types of
fake-medusas where the three quantum ovals are distributed either

(FM2) as one in the core (=central region limitrophe to both
singularities), and one in each legs (periphery) of the medusa;

(FM3) as two in the core, and one in one of both legs.

Of course, existence of such medusas is also ruled out by Viro's
oddity law.]

\subsection{Falling back to the method of deep penetration
with salesman travelling}

[26.09.13] Alternatively to the method of hyperbolizing to a
fake-medusa, one may try to ape the infinitesimal tangency
conditions imposed on the osculating conic at the more global
level of the smooth curve. For concreteness, starting with Viro's
anti-scheme $4\frac{3}{1}\frac{3}{1}\frac{9}{1}$, we may impose 5
basepoints (red-points on
Fig.\,\ref{ViroDEGREE8_fake-medusa-Bezout:fig}, which are located
in the deepest lunes akin to fossil residues of the candelabrum
dissipation. Of course, our picture holds only in the vicinity of
the dissipated
%singular
curve, yet we expect afterwards extending the argument to a
general curve in the fixed isotopy class
($4\frac{3}{1}\frac{3}{1}\frac{9}{1}$). The corresponding conic
will cut the $C_8$ in at least 16 points, still in accordance with
B\'ezout. Can we be more clever?

In summary,  for any trinested octic holds already a phenomenon of
B\'ezout saturation when passing a conic through five deep ovals
provided those are dispatched in all three nests. To contradict
B\'ezout it suffices to arrange  a salesman travelling, i.e. one
more color change than the three granted ones. Each color change
forces 2 intersections (out from the old nest to get in the new
one). Thus a conic with 4 color-changes intercepts the $C_8$ along
$5.2+4.2=10+8=18>16$ points, violating B\'ezout.

Hence the whole game of prohibiting curves reduces to that of
ensuring dichromatism, as opposed to the monochromatism of an
uniform color distribution with only 3 changes. First, it is
advisable to impose the 5 basepoints on the ovals themselves
instead of their insides. (This avoids  a minor technical
%%%point
worry.) Next we may let vary the location of the five basepoints,
inside a five-dimensional tori, where each element is assigned the
interpolating conic (generically unique), and in turn a color
distribution. On varying the position this color distribution
stays constant by continuity unless a catastrophe happens. What
are the catastrophes of the problem? One is certainly the
degeneration of the conic to a pair of lines. If this happens---by
the pigeonhole principle---at least three of the five points land
in the same line, and we get a three-in-line condition which
violates B\'ezout. Indeed, recall that we have two white and two
black points plus a red one (coloring being by appurtenance to a
given nest). Hence our line with 3 points intercepts the $C_8$ in
$3.2+2.2=10>8$ points, since it involves 2 color changes.

From hereon, it seems possible to infer that the conic is always
uniquely defined. Otherwise, there would be  a pencil of such
conics (interpolating the five points=pentagon), but then there
would be also a singular member in this pencil (e.g. by a crude
dimension count of the discriminant=hypersurface), yet this
violates B\'ezout for lines as just observed.

\begin{lemma}
Given any trinested octic, $C_8$, plus a trichromatic pentagon
(injectively) inscribed in the deep ovals of the $C_8$ such that
each color appears at most twice. (As usual, the three colors are
assigned in reference to the three nests.) Then, there is a unique
conic interpolating the pentagon, and it is smooth. Hence we have
a canonical mapping from a five-torus to the hyperspace of conics
$\PP^5$ (which avoids the discriminant). The number of chromatic
changes in the pentagon coursed along the conic is thus constant
through continuous variation of the pentad, and is at least 3
%%%(discolored).
%%%%(unichrome).
(unicolor), in which case the conic cut the $C_8$ already in
$5.2+3.2=16$ points. If
%%%more colorful
multicolor (i.e. at least 4 changes),
%then
B\'ezout is violated
and the octic
%curve
prohibited.
\end{lemma}

So far so good, but can one implement this basic method on any
concrete octic curve prohibited by Fiedler, Viro, Shustin,
Orevkov, or maybe even in the more select realm of bosons not yet
prohibited? (Of course this involves primarily the four binested
bosons, and so the discourse has to be slightly adapted; i.e. only
two colors instead of three.) One capable doing this, gets
probably his name graved on the
%glamorous
obelisk of Hilbert's
problems solvers (Dehn, Bieberbach, Arnold, Gudkov ($m=6$), Viro
($m=7$), etc.)

Let us formulate a (loose) avatar for {\it  binested\/} octics:

\begin{lemma}
Given any binested octic, $C_8$, plus a dichromatic pentagon
inscribed in the deepest ovals of the $C_8$ such that each color
appears at least twice. If a color appears only once, then we
cannot expect more than two color changes, and B\'ezout is
respected ($5.2+2.2=14<16$). Then there is a unique conic
interpolating the pentagon, and it is smooth. (Maybe not true
because it may split off in two lines each monochromatic, hence
respecting B\'ezout). Hence we have a canonical mapping from a
five torus to the hyperspace of conics $\PP^5$ which avoids the
discriminant. The number of chromatic changes in the pentagon
coursed along the conic is thus constant through continuous
variation of the pentad, and can be at most equal to 3
(discolored). If more colorful, then B\'ezout is violated and the
octic curve prohibited.
\end{lemma}

So a direct avatar looks dubious yet we may perhaps impose a 3rd
color by looking at the outer oval (existence granted by Gudkov
periodicity). This we declare as defining the red color. A shift
through it forces only one (bonus) intersection (instead of the
two gained by changing of nest). So we distinguish  strong
(black-to-white) color-changes from weak ones (black-to-red or
white-to-red).

Assume given a binested $M$-octic, and suppose it to have one
outer oval. Choose a pentagon with say 2 pairs of vertices in each
nest and one vertex on the outer oval. The conic interpolating
this pentad has $2.5+2+2.1=14$ real intersections granted in case
of the worst possible colorimetry. This is not enough to foil
B\'ezout, but the distribution with 5 color-changes is enough to
attack B\'ezout. By the way, the interpolating conic could still
split off a line without violating B\'ezout.

Hence the binested case looks intrinsically harder, yet perhaps
subsumable to the same basic method. This is in accordance with
the factual knowledge assembled by Russian scholars, especially
Fiedler, Viro, Orevkov. So it seems wise modesty to first
understand the trinested case (while recovering only old truths),
hoping that no revisionism of Viro sporadic is necessary.

Of course in the trinested case certain curves do exist while
other do not apparently. It is a very subtle matter of deciding
under which circumstance a
%colorful
multicolor pentad can be arranged, and thus the corresponding
scheme prohibited. A first condition for our method to apply is
the presence of at least two nests containing at least two ovals.
Diagrammatically, this merely amounts to rule out the first line
in the 1st layer of the trinested pyramid. This represents no loss
of generality as all those (five) schemes (containing
$\frac{1}{1}\frac{1}{1}$ as sub-symbol) are resp. constructed by
Wiman, Viro, and Shustin.

This being said, we can fix a pentad in our trinested scheme which
is tricolored in such a way that each color is represented at most
twice. Under this condition we can grant no degeneration of the
conic interpolating the pentad, and therefore its uniqueness too.

To each pentad inscribed inside the deep ovals, we may assign a
coloration according to the nest of appurtenance, and count the
number $C$ of color-changes  when circulating along the
interpolating conic (which is unambiguously defined). Of course
$3\le C\le 5$, and it suffices to have $C\ge 4$ in order to
corrupt B\'ezout, since each color-change forces 2 intersections
with $C_8$. We say then that the pentad is {\it multicolor\/}, and
one needs a trick ensuring a multicolor pentad on certain
hypothetical curves as a weapon for their prohibitions.

Alas, it is here that things start becoming difficult. We would
like to show that certain curves prohibited by Viro always contain
a multicolor pentad.

To fix the idea we would like to solve this problem along three
levels of successive difficulties:

(1) the Fiedler oddity law for trinested $M$-schemes without outer
ovals.

(2) Viro's oddity law extending Fiedler's to an arbitrary number
of outer ovals.

(3) Viro's sporadic obstructions (mostly concerned with the case
of naught outer ovals safe one exception).

\begin{figure}[h]\Figskip
%\vskip-1.2cm\penalty0
\centering
%%%%%\hskip-1.7cm\penalty0
\epsfig{figure=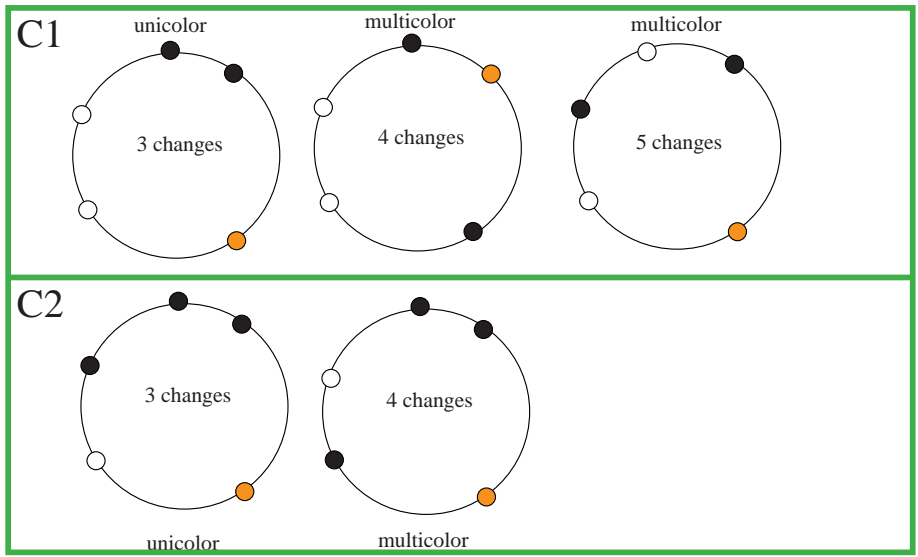,width=122mm} \captionskipAG
  \caption{\label{ViroDEGREE8_DP:fig}%
  %%%Direct anti-B\'ezout
  Trapping
  argument on the singular curve: genuine vs. fake
  medusas} \figskip
\end{figure}

So let us consider Fiedler's setting first in the hope that
history was right in finding it first. We have then basically
22-3=19 ovals ranged in 3 nests.

We must first distinguish between two tricolor spectra:

(C1) two black, two white and one red (5=2+2+1); or

(C2) three black, one white and one red (5=3+1+1).

Perhaps one can do all the work solely with the color spectrum C1,
as the 2nd one really pertains to the 1st line of the 1st layer
which is already completely elucidated by the constructions of
Wiman, Viro and Shustin.

Further, as we said it seems that the color-palette (C1) has the
definitive advantage that the interpolating conic cannot split off
a line, because then the five points migrate apart in the two
lines, one of which containing three of them (pigeonhole), and in
the colorimetry C1 this forces a dichromatism, hence
$3.2+2.2=10>8$ intersections with a line (against B\'ezout).

So in the color (C1), it is ensured that the interpolating conic
is smooth and therefore unique.

To fix better  ideas, we examine Fiedler's scheme
$\frac{1}{1}\frac{2}{1}\frac{16}{1}$. Here we distribute the 5
basepoints on the deep ovals along the coloration C1. So we choose
a pentad with one point on $\frac{1}{1}$, two points on the deep
ovals of $\frac{2}{1}$, and 2 points on the 16 ovals of
$\frac{16}{1}$. Here there is $\binom{16}{2}=8.15=120$ ways to
proceed up to continuous deformation.

Again by continuous variation of the pentad, the conic varies
continuously, and as it stays smooth (C1 hypothesis) the
distribution of colors stays constant during the deformation.
Accordingly,  each distribution of 5 points on the deep ovals of
colorimetry C1 defines unambiguously the number of color-changes
which is either 3 (unicolor), or 4,  5 (multicolor). So from a
brute statistic viewpoint, among the 120 possible distributions
(in case of Fiedler's curve), it would be pure miracle if all 120
chromatic numbers would be 3. The probability for this event would
be ca. $(\frac{1}{3})^{120}$.

From this perspective, it looks very miraculous that trinested
$M$-schemes exist at all, e.g. Shustin's scheme
$4\frac{1}{1}\frac{3}{1}\frac{11}{1}$ where in technicolor C1,
there is $3.\binom{11}{2}=3.55=165$ possible distributions of
pentads. All of them have to be unicolor, and this is a
probabilistic miracle, yet made real by geometry.

So the big problem is to find a {\it technique d'existence\/} of a
multicolor pentad on a Fiedler curve of type
$\frac{1}{1}\frac{2}{1}\frac{16}{1}$ (more generally on any curve
declared prohibited by Fiedler, and especially Viro)

[25.09.13] Besides, the whole theory of the dissipation of the
candelabrum can obviously be correlated (via a tri-ellipse plus a
line)  to Viro's census of septics (especially $M$-septics). In
particular one of Viro's global obstruction should prohibit
certain candelabrum patches. It must be of primary interest to
work this out in full detail at the occasion.

\subsection{F4+F4}

[28.09.13] To construct curves one natural method is that of
dissipation theory, small perturbation of common objects. The
point is that singular (in particular decomposed curves) are
better known and anchor-bases toward the exploration of new
continents. This is the philosophical substance of Viro or older
methods. Alas even in degree 8 the method seems in panne toward
solving the isotopy classification (Hilbert's problem) unless one
is able to prohibit what has not yet been constructed.

One basic idea would be to smooth two F4=$X_{21}$, yet not
incarnated as a quadri-ellipse, but as  a more complicated curve
interconnecting those germs.

However we had already this idea a long time ago, and actually any
octic with this prescribed configuration is a quadri-ellipse.
Indeed the osculating conics with prescribed contacts of tangency
along the 2 singular points form a pencil. Imposing to visit any
additional point of the $C_8$ gives $2.8+1=17>16$ intersections so
that the octic has to split off the conic of the pencil through
that point.

\begin{figure}[h]\Figskip
%\vskip-1.2cm\penalty0
\centering
%%%%%\hskip-1.7cm\penalty0
\epsfig{figure=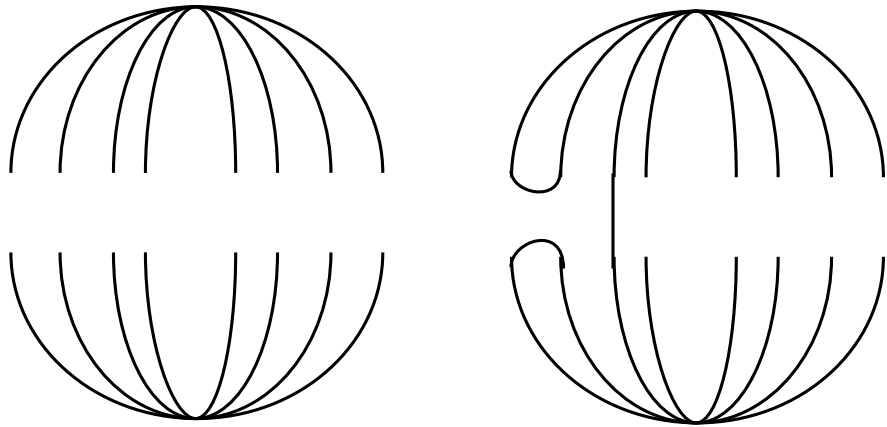,width=122mm} \captionskipAG
  \caption{\label{ViroDEGREE8_F4+F4:fig}%
  %%%Direct anti-B\'ezout
  F4+F4: two rainbows} \figskip
\end{figure}

\subsection{Some general ideas: Stonehenge
%lineation
alinement of all Gudkov pyramids}

[27.09.13] (Written down but based on a older idea ca. April 2013,
when publishing v2 of this text).

As we wrote in the introduction (of v.2) it seems  that there is
phenomenon of stability under satellites that was anticipated by
Wiman 1923, and Rohlin 1978. Here it is understood that if a
scheme of a certain degree is saturated (usually via B\'ezout)
then all its satellites
%multiples
are likewise saturated.

The special corollary is that the Gudkov pyramids in degrees an
integer which has a rich decomposition into primes will be more
lacunary than those in degree a prime where  no censorship is
induced by satellites. So Gudkov pyramids in primes degree will
appear as dense crystal with a minimal number of prohibitions,
while those of compound degrees will have a much more lacunary
architecture, with several flaps and wings of the edifice
completely missing.

Of course the drama is that censorship under satellites (even
combined with all versions of Gudkov periodicity) does not explain
all prohibitions as best exemplified by degree $m=8$. This constat
follows from the Fiedler, Viro, etc. prohibitions, at least
granting them as being correct.

So the general Hilbert's 16th problem splits into two parts:

(1) the regular prohibitions explainable by periodicity, and
satellites censorship,

(2) the irregular part formed by several sorts of prohibitions,
whose raison d'\^etre is poorly understood (at least by the
writer). Of course it could be that all those Fiedler-Viro-Orevkov
style prohibitions are subsumed to a basic B\'ezout-M\"obius style
of prohibitions; or to the method of total reality.

Another way to emphasize our ignorance is along the following
quantitative idea. As we know, of the about 144 logically possible
$M$-schemes in degree $m=8$ satisfying Gudkov periodicity only
about 83 are constructed   at most 89 of them are constructible
(if the present state-of-knowledge is reliable). Hence, for each
integer $m$ we may list all logically possible schemes (under
Gudkov periodicity when $m$ is even) and denote their numbers by
$G(m)$. Of course we also take into account the basic
B\'ezout-Hilbert style prohibitions prompted by intersection with
lines and conics, etc. Here already it becomes a bit messy to
dissociate trivial from nontrivial obstructions. Yet let us assume
that there is  a well-defined $G(m)$ taking into account all
trivial obstructions plus Gudkov periodicity. In contrast one
defines $R(m)$ the number of schemes which are effectively
realized. For instance $R(6)=G(6)=3$. Then $83\le R(8)\le 89 \le
G(8)=144$. The ratio $R(m)/G(m)$ measures essentially the
existential probability for an $M$-scheme to be algebraic.
Naively, see especially the note by Kharlamov-Orevkov 2003
\cite{Orevkov-Kharlamov_2003-asymptotic-growth}, it seems clear
that even $R(m)$ grows exponentially with $m$. However algebraic
curves may become a rarety as $m$ increases and we could imagine
that $R(m)/G(m)$ tends quickly to 0 as $m\to \infty$.

Of course here Gudkov periodicity only intervenes for a censorship
factor of four and so can be actually ignored without altering the
qualitative behavior of the asymptotic ratio. As we said above we
expect that when $m=p$ is prime the score $R(m)$ of algebraic
schemes is high and viceversa it is low when $m$ is much
compounded. Thus perhaps $R(p)$is sufficiently high that the ratio
starts an oscillating behavior without tending to a definite
limit. Of course all this a very naive speculations and just
supply a vertiginous feeling of imagining which sorts of
combinatorial tour-de-force is requested to get some intuition of
how high order algebraic curves looks alike.

\subsection{Dissipation of the candelabrum  via the
theory of septics}

[25.09.13] [not yet written, but a straightforward adaptation of
what we did in degree 8, with Viro's $X_{21}$]. The basic idea
here is to attempt to get as many curves as possible from a
tri-ellipse plus a transverse line. As a reasonable competitor,
one can consider a basic septic consisting of a tri-ellipse plus
the line tangent to the triple branch. Then we have again an
$X_{21}$-singularity (F4 in our more naive notation), and perhaps
the theory of septics affords prohibition on the $X_{21}$-patches.

It is clear that the whole topic is so much ramified that the
researcher quickly looses his strength and moral along the
menagerie of pathes to be explored.

\subsection{Dissipation of $X_{21}$ via septics}

[25.09.13] In this section we focus on the {\it sunset\/} septic
depicted below, consisting of a tri-ellipse plus the line tangent
to one of the singularity. Then one may apply the usual patchwork
method, while hoping to infer ``relatively new'' (i.e. new for our
own personal understanding of the topic) obstructions on patches.
By the way, the universal (absolute=Russian) knowledge is not
complete in our opinion.

First, we tabulate the table induced by Gabard's patch C1 for
$X_{21}=F4$. Here the scheme depends only upon the value of $\ga$
and we get effective constructions of five $M$-septics (marked by
little green squares on the table below). Using Viro's patches C2
we get very monotonically the sole and same scheme $15$ (unnest)
due to Harnack first. Using the patch C3 (where $\al=1,5$) gives
two schemes already obtained via C1. Finally, employing the patch
E, we get four $M$-schemes marked by green circles, two of which
being ``new''. However our expectation to deduce new prohibitions
is {\it not\/} borne out because the (sole) prohibited $M$-septic
is never encountered.

\begin{figure}[h]\Figskip
%\vskip-1.2cm\penalty0
%\centering
\hskip-1.7cm\penalty0
\epsfig{figure=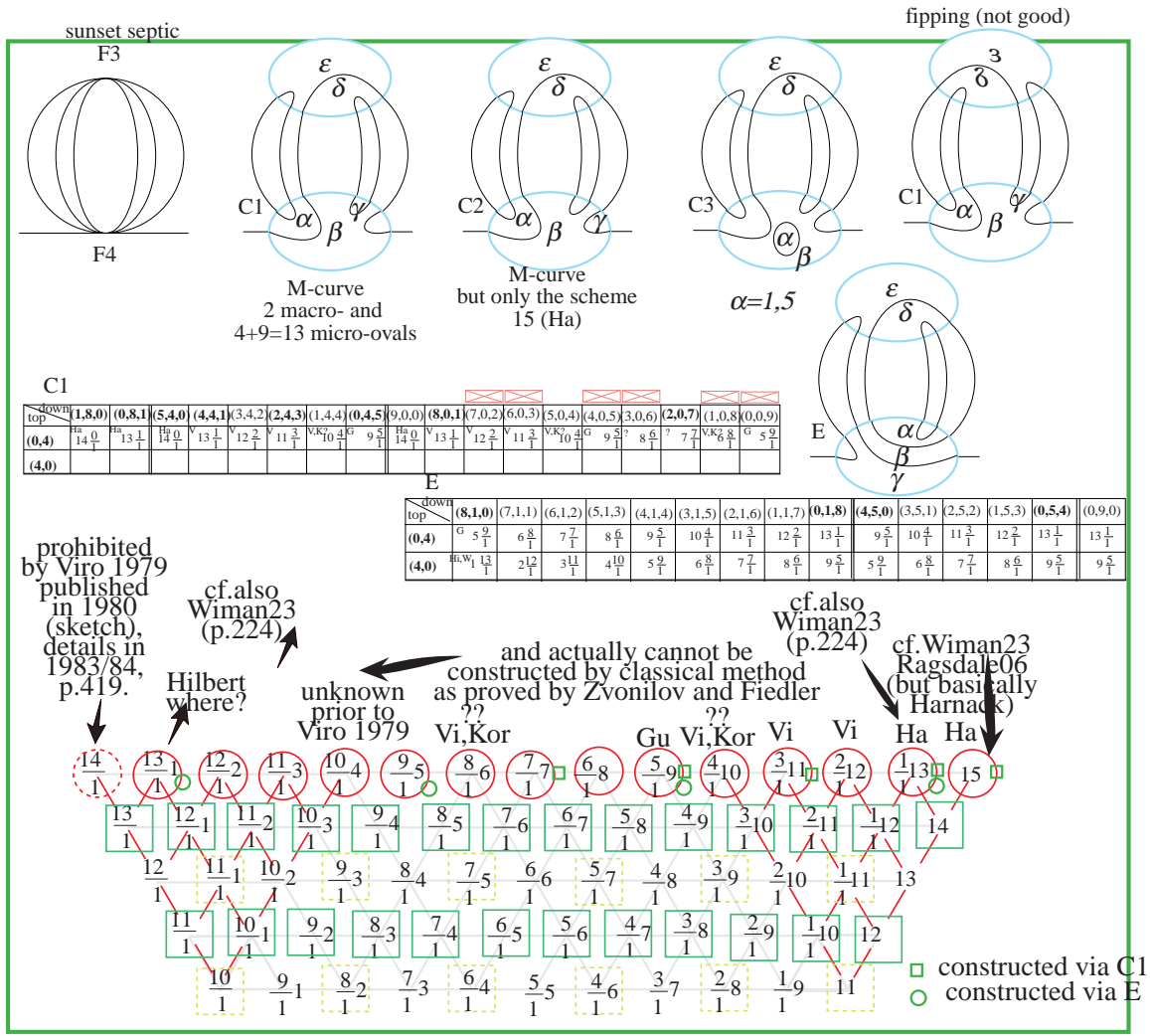,width=152mm} \captionskipAG
  \caption{\label{ViroDEGREE8_septics:fig}%
  Septics and $X_{21}$} \figskip
\end{figure}

Further one must also analyze the other types of patches (A, B, D,
F, G, H, I, J) according to our catalogue
(Fig.\,\ref{ViroDEGREE8_exotic_patches0_SYS:fig}). Working this
out, we note that first the G-patch gives interesting $M$-curves,
yet with one outer oval visible as macro-oval hence there is no
chance to draw a prohibition via Viro's obstruction of
$\frac{14}{1}$ (maximally nested scheme). Then, interestingly, the
patch H creates 3 macro-ovals, and so we get an $M$-scheme despite
non-maximality of the patch employed. It would be interesting to
work out exactly which $M$-schemes are so obtained, but the
presence of one outer lune will not produce any prohibition. Then
we may flip the H-patch, but this forms a snakelike oval wasting
much of the energy in meanders. Finally, the patch I looks the
most
%interesting
promising as there is no outer ovals. Indeed, the patch I(9,0,0)
would create the $M$-septic prohibited by Viro, but the former
(patch) was already prohibited by Gudkov periodicity applied to
the doubled patch (with 0 big eggs, hence not $2\mod 4$). So:

\begin{Scholium}
Quite disappointingly the theory of septics does not prompt any
prohibition upon the patch for $X_{21}={\rm F}4$ the flat point of
multiplicity four.
\end{Scholium}

\begin{figure}[h]\Figskip
%\vskip-1.2cm\penalty0
%\centering
\hskip-1.7cm\penalty0
\epsfig{figure=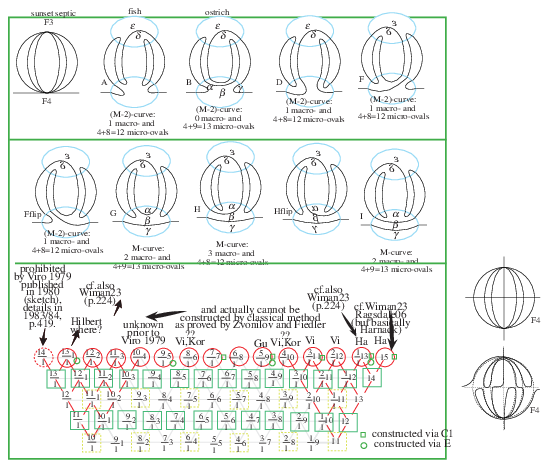,width=152mm}
\captionskipAG
  \caption{\label{ViroDEGREE8_septicsB:fig}%
  Septics and $X_{21}$ (continued)} \figskip
\end{figure}

Of course when less tired one can do the same game for the
candelabrum using the septics consisting of a tri-ellipse, plus
the line through both singularities.

As a guess, it seems that---since the theory of $M$-septics is so
little obstructed---we may not be able to draw any serious
prohibitions on patches by this method. So, one may wonder what is
the avatar in degree 8 of the quadri-ellipse allied to
$X_{21}=F4$, when it comes to the candelabrum C4. Perhaps the
natural candidate  is Shustin's medusa.

\subsection{Sequel of old text}

{\it Sequel of the old text.}---Now at some more fundamental level
the ubiquity of the horse-shoe---as a fundamental shape crossing
four times a line---becomes when smashed Viro's pattern of
dissipation of type C (i.e. lateral double-lune plus two simple
lunes). Actually as shown by Fig.\,e it seems that the method only
yields the types C, D, and A. Of course the operation of smashing
is akin to chocolate and cream decoration in French gastronomy,
namely the experience of taking a knife and dragging through black
chocolate and white-colored cream so as to created the depicted
patterns.

As we note yesterday already it seems that Viro does not exploit
the rabbit with one invaginated ear (as depicted on our Fig.\,a),
and this can be interpreted as a so-called ``angst-Haase'', i.e.
an {\it anxious rabbit}. Of course, it is hard to imagine which
affine quintic could be the antecedent (primitive) of this
angst-rabbit, since Polotovskii's curves seems to be the only
possible alternative to standard undulations. Let us yet,
cavalier, inspect which sort of patch could result from such a
possibility. In fact, noting that the angst-Haase has in fact two
circuits we can better imagine which sort of quintics is the
primitive of the configuration. This amounts just to breakdown of
the wave from Polotovskii's model, and this by B\'ezout can only
occur if the separating mass of water contains no bubbles of
oxygen (otherwise B\'ezout for quintic is violated unless the
configuration reduces to the deep nest $\frac{1}{1} J$), which has
however a too ridiculous number of oval to merit our attention.
Notwithstanding if the breaking mass of water is empty (of
oxygenation) then the configuration is permissible, but of course
we loose one micro-oval having only five of them on the quintic
(since one is consumed by the breaking mass of water). Hence
before doing any specific depiction we look wrong engaged to reach
any $M$-patch from the angst-Haase.

Now, albeit this fails miserably, we get the idea that starting
say from Fig.\,1 we can do the dissipation in such a way that the
ear intercept the smashed line, and so we get the variant V1. This
produces the patch A(0,0,1,0,7) using more-or-less
self-explanatory notation. Of course, this only an $(M-1)$-patch
because to get the type A, we essentially wasted one oval just to
create the singularity $X_{21}$, and thus we cannot expect an
$M$-patch. Our patch when doubled gives the scheme
$17\frac{2}{1}$, which is well-known (i.e. accessible to Viro's
purest method via the quadri-ellipse).

Now we confess being a bit a ``cours-d'imagination'', yet we can
still explore the result of opting always for the dual vibration.
It seems that this trick essentially amounts to the symmetrization
device used by Viro as shown by our Figs.\,D1 and M1. Still a
systematic search looks desirable. To explore this properly we
need a new figure (Fig.\,\ref{ViroDEGREE8_PATCH3:fig}). For the
next dualization D2 we get the same patch C3(1,7) due to an
evident symmetry. Working out D3, we get as expectable from D1
just the symmetrized patch C1(0,4,5) where over the original
construction O3 the lateral symbols are just switched. Alas this
patch is not new as it was already cooked by Viro's M3 (on the
former plate=Fig.\,\ref{ViroDEGREE8_PATCH2:fig}). At this stage it
seems that the work is automatic, i.e. the dual vibration should
just produce the patch with palindromic parameters, i.e. $(\al,
\be, \ga)$ changes to $(\ga, \be, \al)$. Yet some surprise (cf.
the mirrors of the earlier plate) still encourage us to tabulate
naively the dualized patches. (Philosophy: Nothing is more
concrete than mathematics, especially geometry.)

Then we arrive at D4 and this is certainly not even worth
depicting by same symmetry as that encountered by D2. Next we
arrive at D5(dual), but then it seems
%%%%that we must
necessary to distinguish two cases depending on the location of
the newly formed micro-oval arising through dissipation of the
triple-point. In the first version D5=D5A, the new oval is on the
left, while in the 2nd version D5B, it sits on the right of the
line smashed under the hyperbolism. The 1st cast yields the patch
C1(7,0,2),  heavily prohibited by Viro's oddity law or by
Orevkov's dematerialization of the boson b3. At this stage the
philosophy is two-fold: first Viro's method (liberally
interpreted) seems nearly to violate Viro and Orevkov's
obstructions as we saw via D5=D5A. However dissipation D5B leads
to an admissible patch, and is by the way kinematically more
likely, since the new micro-oval is located in the prolongation of
the branch performing the vibration. Still, one could
counter-argue that even on Fig.\,D5 we could arrange this property
by increasing the curvature so as to form an isthmus nearly
connecting the branch to the micro-oval (cf. detail D5D). Perhaps
one can counter-counter-argue that the infinitesimal cubical patch
then seems to violate B\'ezout by tracing the orange-line which
seems to intercept five times the cubic. So we gain perhaps here
some insight of why it is not so easy to corrupt Viro's oddity law
nor Orevkov's obstruction.

Then we have D6 (dual) with again a surrealist micro-oval formed
on the left and the resulting patch again violates Viro's oddity
law. However, the more realist version D6B produces the patch
C1(4,4,1) already found by Viro (at least provided that there is a
misplacement of the symbol $\ga$ on his Fig.\,55). Of course it is
also more likely that our visualization of the hyperbolism is
slightly incorrect leading to a twist of all the results. In any
event this is merely a psychological difficulty that should be
easy to fix once more time is available.

Another point is to wonder if there is also a bifurcation of the
dual dissipation in the earlier cases 1,2, 3, 4 depending on the
location of the micro-oval. Of course we can drag on D1 the
micro-oval on the right of the smashed line, yet this will not
affect the isotopy type unless we drag this oval to its ultimate
confinement, namely inside the meander, back again to the left
side of the smashed-line. However the resulting patch will frankly
corrupt B\'ezout (for lines) as the duplicated patch will exhibit
a nest of depth 3 plus an outer nest of depth 2, forcing 10
intersections with a line through their centers.

At this stage we must (rather disappointingly) confess that Viro's
search looks exhaustive unless one can imagine a really new twist
of the construction.

\begin{figure}[h]\Figskip
%\vskip-1.2cm\penalty0
%\centering
\hskip-2.7cm\penalty0
\epsfig{figure=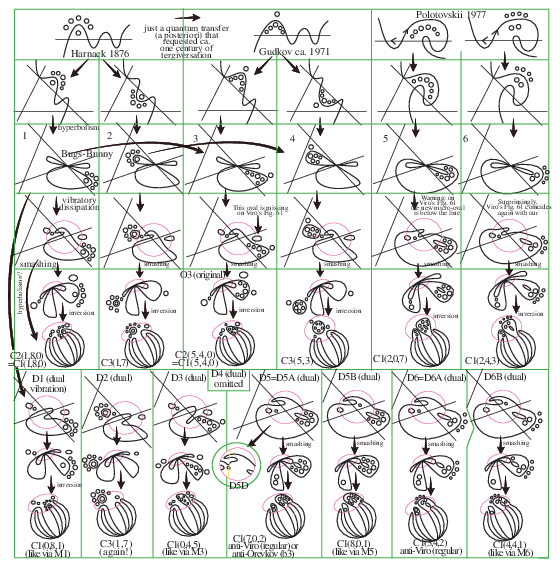,width=172mm} \captionskipAG
  \caption{\label{ViroDEGREE8_PATCH3:fig}%
  Dual vibrations as those of Viro}
\figskip
\end{figure}

So far so good, and we have modulo the C1-versus-C2 ambiguity a
complete understanding of Viro's theorem regarding patches of type
C. It remains now to understand those of type E (trinested lune),
which are explained in Viro, p.\,1119 (especially Fig.\,56). Alas,
this figure is
%%%a bit
awkwardly depicted in the Bible (Viro 89/90), but seems to involve
another genius stroke of Oleg Yanovich's imagination.

\section{Viro's vibratory method}

This section explains another fundamental construction by Viro, in
some sense even more elementary than the one involving hyperbolism
presented earlier.

\subsection{Viro's trick for patches of type E (trinested lune)}

[03.09.13] Viro starts with a pair of conics tangent at one point
and transverse at the remaining two points (cf.
Fig.\,\ref{ViroDEGREE8_PATCH4_Eclass:fig}). A suitable
perturbation of their union  offers a quartic $C_4$ oscillating as
depicted across $C_2$. Note that the intersection $C_2\cap C_4$ is
totally real involving 4 transverse and two 2nd order contacts at
the ``north pole''. Actually, as the sequel of Viro's Figure~56
involves an $A_3^-$-singularity we wondered if the perturbation
$C_4$ is not rather  involving a tangency. Recall that $A_k^-$ is
the germ of $y^2-x^{k+1}=0$. Then we managed finally understanding
Viro's picture despite
%the fact the original
%is
being really poorly traced (at least on my small sized Xerox copy
of the article). Notwithstanding Viro's construction is genial,
and the crucial step is to count properly the contacts to get the
right perturbation (Fig.\,a). At the $A_3$-point we have two
contacts of order two between $C_2$ and $C_6$, while at the
$J_{10}$-point we have 3 contacts of order 2. Hence the
intersection $C_2\cap C_6$ consists already of $2.2+3.2=4+6=10$
intersections (counted by multiplicity), whence the possibility to
impose two additional intersections as shown by the bump on
Fig.\,a. The sequel of the construction should be self-explanatory
from the figure. It is perhaps still puzzling that at some stage
of the argument we thought that the line through both
singularities of the $C_6$ would corrupt B\'ezout, but apparently
not so. Another slightly puzzling aspect is that on Fig.\,b the
curvature of the branches of the tripod-singularity $J_{10}^-$
does not seem respected: maybe there is a topologico-metrical
parade identifying this as mere optical illusion).

\begin{figure}[h]\Figskip
%\vskip-1.2cm\penalty0
%\centering
\hskip-2.7cm\penalty0
\epsfig{figure=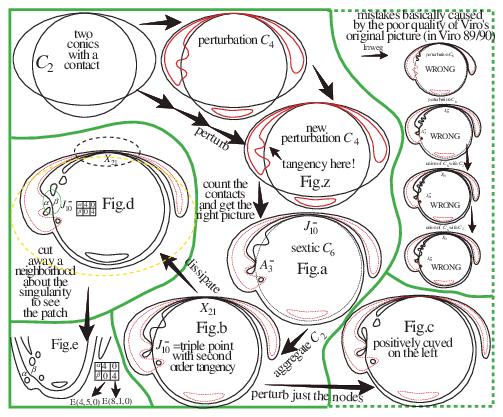,width=172mm}
\captionskipAG
  \caption{\label{ViroDEGREE8_PATCH4_Eclass:fig}%
  Viro's vibrational method leading to the E-class (trinested lune)}
\figskip
\end{figure}

{\it Added\/} [14.09.13].---Actually there is at last two parades.
A first involves  Fig.\,c consisting in first dissipating the
nodes of the $C_8$ (Fig.\,b) and one may expect  the resulting
curve $C_8$ having three branches positively curved inside the
same half-plane. Then we are in a position to apply the usual
dissipation theory of this triple point ($J_{10}$). The other
parade is that our Fig.\,b is actually much distorted. In reality,
the two branches at $J_{10}$ which looks curved to the left are in
reality much closer to the circle $C_2$ to such a point that those
branches are in fact curved to the right. Concomitant to this,
remark that on our picture of the quartic (Fig.\,z), the line
tangent $T_p C_4$ to the bicontact of $C_4\cap C_2$  seems to
intersect 6 times the quartic. This aberration is dissolved if the
curve $C_4$ is imagined much closer to the circle. Getting a
metrically accurate vision is a challenging task, compare
optionally our free-hand Fig.\,\ref{ViroDEGREE8_PATCH4_ZOOM:fig}.

Next, dissipate the triple-point $J_{10}^-$ to get
Fig.\,\ref{ViroDEGREE8_PATCH4_Eclass:fig}d. Finally,  cut away a
neighborhood of the $X_{21}$-singularity to find with the
complement a patch for the same singularity $X_{21}$. This is a
trivial, yet somewhat miraculous step, reminiscent of Steiner's
Wiedergeburt und Neuauferstehung (when it came to philosophize
about inversions). From Fig.\,e we easily recognize the patches
E(4,5,0) and E(8,1,0), in the notation of our catalogue
(=Fig.\,\ref{ViroDEGREE8_exotic_patches0_SYS:fig}). Keep  maybe in
mind the following slight objection:  the excised object is an
$\RR P^2$ less a disc so a M\"obius band, hence not so much a
topological disc, as one imagine the patch substratum. This defect
can be resolved by choosing instead the yellow-colored ellipse and
by keeping its inside instead (see again Fig.\,d).

\begin{figure}[h]\Figskip
%\vskip-1.2cm\penalty0
%\centering
\hskip-2.7cm\penalty0
\epsfig{figure=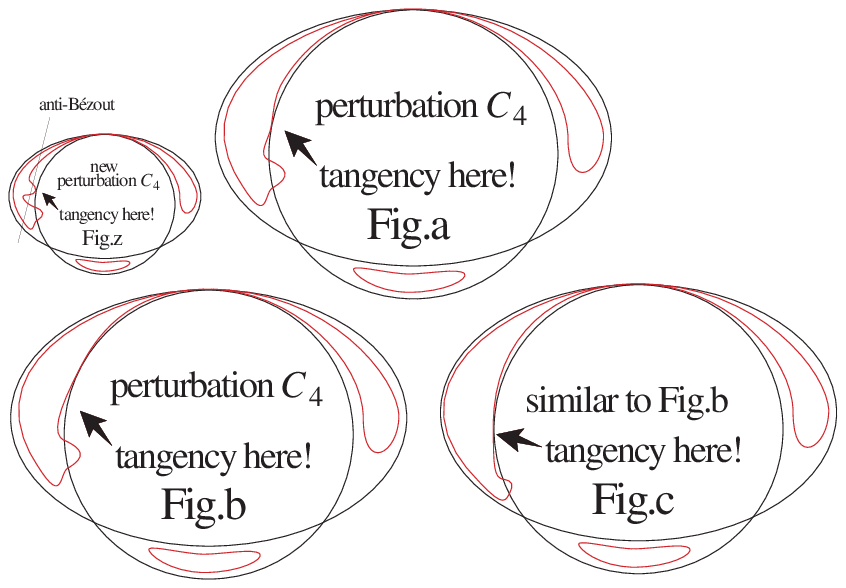,width=172mm}
\captionskipAG
  \caption{\label{ViroDEGREE8_PATCH4_ZOOM:fig}%
  Viro's E-class (trinested lune)}
\figskip
\end{figure}

Note: If we could choose other values of $(\al,\be)$, e.g.
$(2,2)$, we could get more patches, but unfortunately $(2,2)$
though realistic as being involved in subdivision of Gudkov's
sextic $5\frac{5}{1}$---imagined as a patchwork of
$J_{10}$-singularities---are not realized since when glued with
Viro's patches (4,0)/(0,4) yields sextics corrupting Gudkov
periodicity.

Comparing with Viro's catalogue, we still miss two of his patches.
In his article (Viro 89/90) he proposes a conceptual argument we
were not able to follow. Surely, there is alternatively a variant
of the  construction in view of the palindromic reversion of the
missing parameters. (The sequel will indeed supply an alternative
of Viro's construction doing this job.)

More importantly, the class E is presently not much obstructed and
one is naively expecting that a variant of Viro's trick should be
capable producing more patches, especially those materializing the
two remaining subnested bosons B4 and B14.

One naive idea is what happens in the above construction if the
oscillation (bump) of Fig.\,a is is effected on the other side of
the circle, or eventually below the tacnode singularity
($A_{3}^-$).

A more elaborate idea that we had later on that day (22h31), is
wonder about the case where two tangencies are arranged. However
it seems then that we get already $(2+2+3)=7$ second order
contacts between the $C_6$ and $C_2$ (even prior to introducing
any bump). This yields a multiplicity intersection of $7.2=14$
corrupting severely B\'ezout. So it seems that we cannot arrange
such a double bicontact (as on Fig.\,A). Perhaps it would still be
of interest to see which sort  of patches results from
transgressing this B\'ezout obstruction (probably one which
overwhelms violently Harnack's bound). To our little surprise
Harnack is respected (9 micro-ovals) despite the double production
of ovals allied with the pair of triple points, but the patches so
obtained E(1,8,0), E(9,0,0) and E(5,4,0), when doubled certainly
violates Gudkov periodicity (e.g. 2.E(1,8,0)$=$). For instance,
E(9,0,0) doubles to the scheme $1(1, 0 \frac{19}{1})$, which is
not to be found on the periodic table of elements.

\begin{figure}[h]\Figskip
%\vskip-1.2cm\penalty0
%\centering
\hskip-2.7cm\penalty0
\epsfig{figure=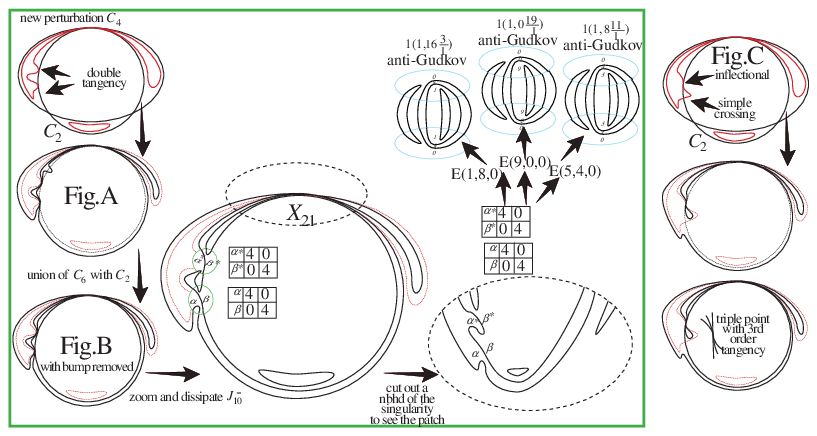,width=172mm}
\captionskipAG
  \caption{\label{ViroDEGREE8_PATCH4B_Eclass:fig}%
  Double bicontact (fails miserably against
  B\'ezout and then Gudkov)}
\figskip
\end{figure}

Despite failing miserably, our attempt may adumbrate other
combinations of singularities and contacts permissible for
B\'ezout procuring more patches (potentially leading to new
breakthroughs of more flexibility in Hilbert's 16th).

[05.09.13] More modestly, we may wonder if this elementary
vibrational method of Viro (yielding the E-patches) modifies as to
offer as well the C-patches constructed by the somewhat different
technology of hyperbolisms. (We shall find a positive answer at
least upon admitting some common elasticity of algebraic
geometry.)

Turning again to our last miserable construction (double
bicontact), one can wonder if the topos is improved if we
amend on
the outer side of the ellipse (Fig.\,X). Of course, this will
change nothing to the numerology of bicontacts, and B\'ezout is
still jeopardized.

One another possibility is to have a triple contact followed by a
transverse crossing (as shown on Fig.\,C). It remains then to
improvise the dissipation theory of triple point with 3rd order
tangency.

Another option is that of introducing a finger-move on the quartic
like on Fig.\,D., on the variant of Fig.\,E. Of course B\'ezout
looks foiled when tracing a suitable line centered through the
Hohlraum (=trap) formed by the finger-move. Still, on resorbing
progressively the ``Falaise-pocket'' through the isotopy suggested
by  Figs.\,E,F,G we may rehabilitate B\'ezout, and so perhaps
there is some quartic perturbation realizing the qualitative
picture of Fig.\,G. Alas, the latter as another B\'ezout defect
with respect to the dashed line, but this can be remedied by
deflating the inner bump below the ``horizon'' as shown on
Fig.\,H. As a conundrum,
%%% enigme devinette CHECKED IN DICO CASSELL's French Dictionnary
it seems that in the sinuous $S$-shaped tube there will be a
bitangent line which when escaping from the circuit has to create
at least 6 intersections with the $C_4$.

Despite all those defects, let us apply Viro's algorithm to
Fig.\,G (as being a respectable isotopic model). On perturbing
$C_2\cup C_4$ we get Fig.\,g0 with an flex on the left of $C_2$,
which we perturbed transversally  on Fig.\,g1 so as to avoid
referring to an obscure dissipation theory. On the latter figure,
we count in the intersection $C_2\cap C_6$ as many points as
$2.2+3.1+3.2=4+3+6=13>2.6=12$, overwhelming B\'ezout. Of course a
stupid parade is to lower this
%%Pech
mischance-number 13 to 11, by dissipating the undulation yet there
is an anomaly with B\'ezout-Galois (i.e. the B\'ezout count modulo
two over the reals). Another option, would be that during the
perturbation we do not have anymore 3 bi-contacts at the north
pole of the circle $C_2$. Actually, this is the forced scenario as
soon as we take notice that the most in-curved branches through
the north pole actually crosses the fundamental circle. Hence the
true multiplicity count for $C_2\cap C_6$ as materialized on
Fig.\,g1 is $2.2+3+5=12$ and B\'ezout is intact. Working out the
resulting patch gives C1(3,5,0) and C1(7,1,0). Alas, those have
only 8 (micro) ovals, so not $M$-patches. Still, it seems of
interest to stress that so Viro's vibratory method reaches the
C-class of patches, though it looks apparently difficult to gain
maximal patches. (We shall soon see that we can arrange maximality
in the C-class as well!)

\begin{figure}[h]\Figskip
%\vskip-1.2cm\penalty0
%\centering
\hskip-2.7cm\penalty0
\epsfig{figure=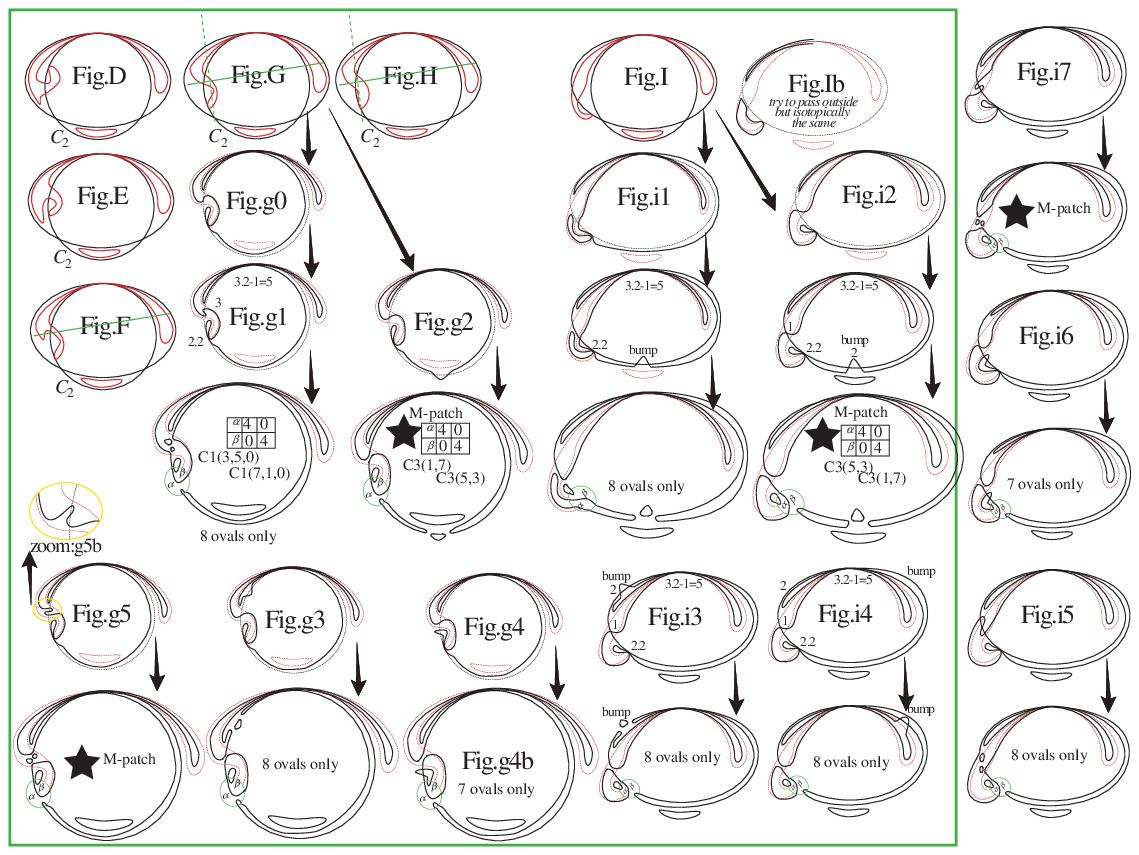,width=172mm}
\captionskipAG
  \caption{\label{ViroDEGREE8_PATCH4C_Eclass:fig}%
  Finger moves: elementary construction of the C3-patches}
\figskip
\end{figure}

Of course it could be that we misplaced the bump, imagined as
cached in the oscillation. So there is perhaps more clever bumps
leading to $M$-patches. So, considering Fig.\,g2, yields
indeed---somewhat miraculously---$M$-patches, und zwar (=and
actually, in German) those with symbols C3(1,7) and C3(5,3). Of
course, those are {\it not\/} new, but now obtained via a perhaps
more elementary method avoiding hyperbolisms. Of course, one
challenge could be to obtain all of Viro's patches (and more if
divinity agrees) by this uniform method (due to Viro, but perhaps
twistable).

Evidently, we may the alter the finger-move trick by oscillating
instead across the eccentric ellipse. This idea materializes to
Fig.\,I.

\begin{Scholium} Generally speaking, especially in the fingers of Gudkov,
Polotovskii, Orevkov, this suggests that most of the
constructional aspect of Hilbert's 16th must reduce to a catalogue
of erotical position
%%%%(so-called kamasutra?)
%%%%YES saw a German film on this recently
%
(kamasutra like) adopted upon by algebraic curves (especially
those of decomposing type where both components can interlace
along
%%%nearly all imaginable ways
fairly complicated patterns).
\end{Scholium}

Back to Fig.\,I, we derived only  $(M-1)$-patches via Fig.\,i1,
but perhaps we missed a more strategic option. Indeed, adhering to
the more clever smoothing of Fig.\,i2 we arrive at $M$-patches,
but alas the same as those already obtained (via Fig.\,g2). Did we
exploited all possibilities? As usual it is here that the brain
starts blocking, as the problem requires both memory and
combinatorial skills (creativeness), which are somehow
incompatible hemispheres of the brain (like the dead and vive
memory in computing machines).

A naive idea is to vary the bump location. So from Fig.\,g2 we
manufacture Fig.\,g3 (with a bump on the left-fringe of Viro's
hairs which in reality are short-cut). The resulting octic patch
has 8 ovals only (and belongs to type A).

Next Fig.\,g4 shows the case where the bump is placed inside of
the ventricle. This variant looks quite erotical, yet hopefully
still algebro-geometrizable. We get so Fig.\,g4b where there is
two maximizing options of smoothing (nested or not), but alas
there is only 7 ovals created (either way).

Of course it remains now to work this out more systematically
(i.e. all bump locations also the case of Fig.\,I), but we wanted
prior to this to investigate the fairly contorted case, akin to
Fig.\,g1, yet where the undulation sense is ``reversed''. By this
we mean Fig.\,g5, which admittedly does not look very natural, but
it looks wise exploring all the options as to get a better grasp
of the phenomenology. On the zoom:g5b we try to show how the real
oscillation looks alike, but it is quite puzzling to know if this
works algebraically. Working out the next perturbation, gives us
indeed $M$-patches, but unfortunately the same ones as those
cooked by Fig.\,g2.

[06.09.13] Next we tried Fig.\,i3 where the bump is placed on the
left fringe, but it results only an $(M-1)$-patch with 8 ovals.
This suggested also placing the bump on the right fringe
(Fig.\,i4), which creates only 8 ovals. Further there is 2 options
for smoothing the bumpy part, and one yielding a patch of type D,
which as far as we can remember was not yet realized.

Then Figs.\,i5, i6, i7 are obvious variants, the latter of which
giving $M$-patches, but alas still the same two of type C3.

Next deforming Fig.\,I we get Fig.\,J, yet with oscillatory
pattern across the ground ellipse isotopic to that of the previous
configuration, so that there is no
%%%hope
chance to get something new. Fig.\,K shows a more contorted
position for the $C_4$, yet one violating B\'ezout (intersect with
the circle), except if one branch is actually transverse to the
north pole (the most incurved branch actually traverse the circle,
so the contact must be odd and it would arrange us being just
one). The result is an $M$-patch of type A, which for $\al >0
(=4)$ violates
 B\'ezout when glued with its flip. (Then we get an extension of the
biquadrifolium $\frac{1}{1}\frac{1}{1}\frac{1}{1}\frac{1}{1}$,
which is saturated in degree 8).

As a moral too much erotical contortion foils B\'ezout and kills
any viable progeniture. Moreover Fig.\,k2 yields a patch with
supernumerary 10 ovals (so violating Harnack, aber Hallo!). Hence
Fig.\,K looks definitively too violent.

Then we have of course Fig.\,L with an upward finger-move, but by
a trapping argument (based on Jordan separation like in
Poincar\'e-Bendixson), it seems impossible that this will ever
satisfy B\'ezout. Despite, it is perhaps still informative to
inspect which sort of patches results from this unlikely specie.
The perturbation of Fig.\,l1 yields only a patch with 8 ovals,
while that of Fig.\,l2 seems to corrupt B\'ezout if we keep the
bump. Transgressing this, the resulting patch has nine ovals and
is of type A, but of the sort violating Arnold's weak version of
Gudkov-Rohlin periodicity. The version of Fig.\,l2 without bump
(at the south pole of the ellipse) probably exists, and yields
A-patches with 7 ovals, namely $A_{+}(5,1)$ and $A_{+}(1,5)$. When
doubled those produce the $(M-4)$-schemes
$6\frac{5}{1}\frac{5}{1}$ and $14\frac{1}{1}\frac{1}{1}$ (whose
geography can be checked on
Fig.\,\ref{Degree8-(M-i)-curve-TABLERohlinTRIS}).

Next, the earlier Fig.\,G suggests that the left-fringe of the
$C_4$ may travel as far as to touch the ellipse, and it remains
then the option to get a bump as on Fig.\,M. This gives two
$M$-patches each interpretable as belonging to either class-C1 or
C2. Although those patches are structurally new, it is a quite
spectacular methodological success because Viro accessed to these
patches via the somewhat more elaborate technique involving
hyperbolism. Further, this is the first intrusion of Viro's
elementary method in the realm of C1- and C2-patches.

Then, we realized that Fig.\,M may be changed to Fig.\,N, which
may look more natural (hence algebraizable), yet still leading to
the same patches as before, since the construction is independent
of the circular ellipse.

The next natural variant is Fig.\,O, but that was already analyzed
via Fig.\,\ref{ViroDEGREE8_PATCH6_Eclass:fig}z (perhaps later in
this text but earlier in our historiography).

At this moment, we see that we can let oscillate the left-fringe
(w.l.o.g.) either first tangentially and then transversally
(abridged tantra like on Fig.\,O) or vice-versa first
transversally and then tangentially (tratan as on Fig.\,N).
Additionally we can tell the same story by vibrating across the
circular ellipse (so-call internal vibration). In principle, we
have already followed all those options. However, on comparing
Fig.\,N with Fig.\,I we see that it is not enough specifying first
transverse and then tangent, but there is two options depending on
whether the appendicitis makes it contact from inside or outside
(again compare Fig.\,N with Fig.\,I).

\begin{figure}[h]\Figskip
%\vskip-1.2cm\penalty0
%\centering
\hskip-2.7cm\penalty0
\epsfig{figure=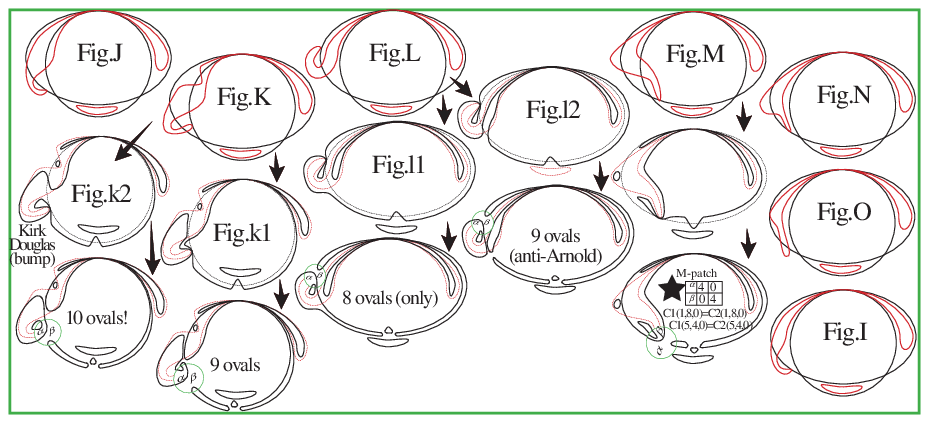,width=172mm}
\captionskipAG
  \caption{\label{ViroDEGREE8_PATCH4D_Eclass:fig}%
  Finger moves: elementary construction of the C3-patches}
\figskip
\end{figure}

So at this stage it seems that we can work more systematically
with Fig.\,\ref{ViroDEGREE8_PATCH4SYS_Eclass:fig} that should be
self-explanatory. Note yet from the eight possibilities two are
isotopic for tautological reasons. Next we apply Viro's algorithm
(of perturbation with bump) to all possibilities and mark by green
triangles the resulting patches on the main-catalogue
(=Fig.\,\ref{ViroDEGREE8_exotic_patches0_SYS:fig}). As we already
experimented it seems that to get an $M$-patch we are forced to
bump between the north pole ($X_{21}$) and the point of tangency
of $C_4\cap C_2$ (which becomes a $J_{10}$-singularity). After
completing this genealogical tree (where it is not necessary to
work out the last specimen of each series of four), we recover
with triangles all the patches gained erratically by circles via
the same method. However we do not get all the C1- and C2-patches
gained by the hyperbolisms method. Unfortunately, it seems that we
have exhausted the power of the method, but we may still hope that
suitable twist of this basic method will give more patches.

\begin{figure}[h]\Figskip
%\vskip-1.2cm\penalty0
%\centering
\hskip-2.7cm\penalty0
\epsfig{figure=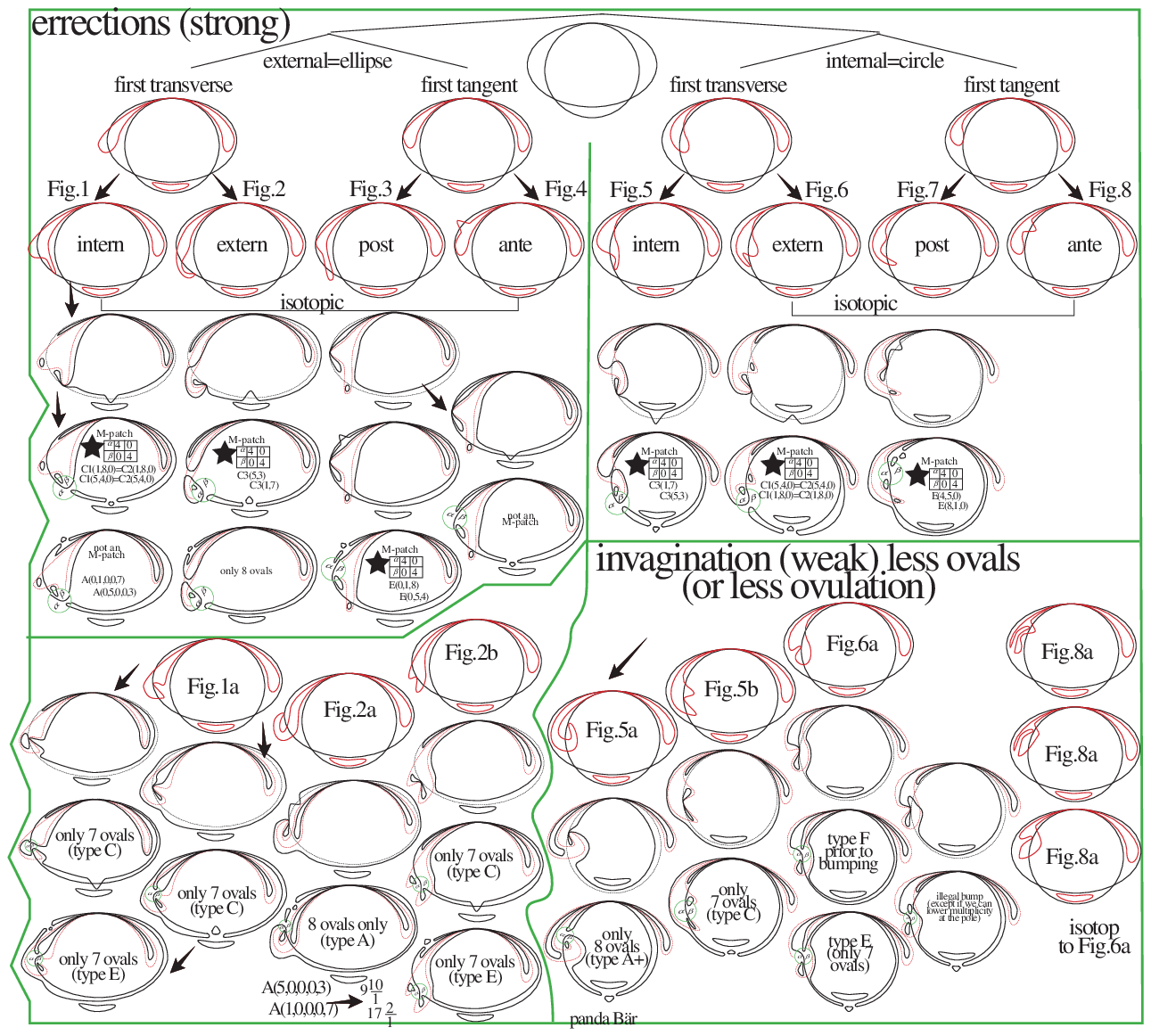,width=172mm}
\captionskipAG
  \caption{\label{ViroDEGREE8_PATCH4SYS_Eclass:fig}%
  Finger moves: systemic construction of the patches via
  erections (which turn out to be more energetic than
  the dual invaginations)}
\figskip
\end{figure}

A perhaps promising route is to explore situations where the 2nd
cytoplasmic expansion involves a curly protuberance (like on
Fig.\,5a), but alas often those versions enter in conflict with
B\'ezout. Moreover the first $C_4$ in the first census of eight
admit an invaginated version depicted below (as Fig.\,1a).
Although sembling erotical, we reached only 7 ovals so that some
vibratory energy is lost due to the invagination.

Then we may look at the dual erection of the 2nd $C_4$, getting so
Fig.\,2a, which probably leads nowhere due to a basic B\'ezout
corruption. Somewhat more exciting looks Fig.\,2b, but alas we
could only press 7 ovals out of it (one of the reason being that
the triple-point lacks a maximal dissipation with the external
branch as leaf). This is essentially a corollary of Rohn's
prohibition of the $M$-sextic $\frac{10}{1}$.

[07.09.13] Then along the idea of the invaginated protuberance we
have also (dually) Fig.\,5b, but here again the maximal
$J_{10}$-patches are not ideally suited to perturb the
corresponding $C_6$.

{\footnotesize ({\it Natur-sozial parenthesis.}---Everybody must
know precisely why his work is useful, and not be the slave of a
capitalism. The goal of all  work---since Neanderthal and
earlier---is to reach free immortality of the individuum, and this
must be conscious in every mind as the true motor of life. No
money is required in such a system, and its usage even impedes the
system being perfectly rentable.)

}

{\footnotesize

Of course, one of the little difficulty of Hilbert's 16th is its
Warenhaus catalogue nature (compare optionally a commentary upon
Edmund Landau). Besides, the position adopted by algebraic curves
as rigid configuration yet able doing the most
erotical/acrobatical position provided B\'ezout restriction are
respected emphasize a sort of elasticity of the algebro-geometric
crystal, which for some deep reason is both allied to gravitation
(ellipses by Kepler-Newton), and the role of higher order curve in
optics by Newton (etc.) By Gabard 2012 (v1 of this text), it was
also clear that there is some connection with dynamics of the
electrons around any massive nucleus, at least so can one
interpret the fantastic dance of points allied to a totally real
map \`a la Ahlfors.

}

Alas, as to our concrete problem it seems (being in a bad day)
that our above organigram
(Fig.\,\ref{ViroDEGREE8_PATCH4SYS_Eclass:fig}) have exhausted all
the swing of this Viro (vibratory) method.

Of course we can still imagine a replica of the catalogue with all
other positioning of the bump of the $C_6$, although we think to
have always exploited the maximizing option. Then it remains also
to investigate all the variants where there is a bicontact (i.e.
two points of tangency between $C_4$ and $C_2$). After that it
remains also to study the configuration $C_4\cup C_2$ with two
pairs of transverse points.

 At any rate, the core of Viro's method is much akin to the
vibratory methods of Harnack-Hilbert (ovals=ovules, etc.). In the
case of bi-contacts (as we already once experimented), but
re-experiment again (Fig.\,\ref{ViroDEGREE8_PATCH4SYS_TAN:fig}) we
can even reach a patch with 10 ovals thereby  corrupting Harnack's
bound. Of course, our intersection $C_6 \cap C_2$ involves
$4.2+3.2=7.2=14>12$ supernumerary intersection, yet we can imagine
that multiplicities at the north pole can be lowered (from 6 to 4)
by transverse behavior. Notwithstanding there is some
psychological frustration that the involved lovely picture
(reminiscent of a galvanic current) does not produce a reasonable
patch.

As a loose idea, we have not yet the energy to follow, one could
imagine that the fringe of the $C_4$ has contact with both the
ellipse and the circle and that during the vibratory process
always accompanied by the aggregation of the conics, one alters
between the circle and the ellipse. This would be a sort of
alternating Viro method, but  probably this idea leads to no
serious result.

\begin{figure}[h]\Figskip
%\vskip-1.2cm\penalty0
%\centering
\hskip-2.7cm\penalty0
\epsfig{figure=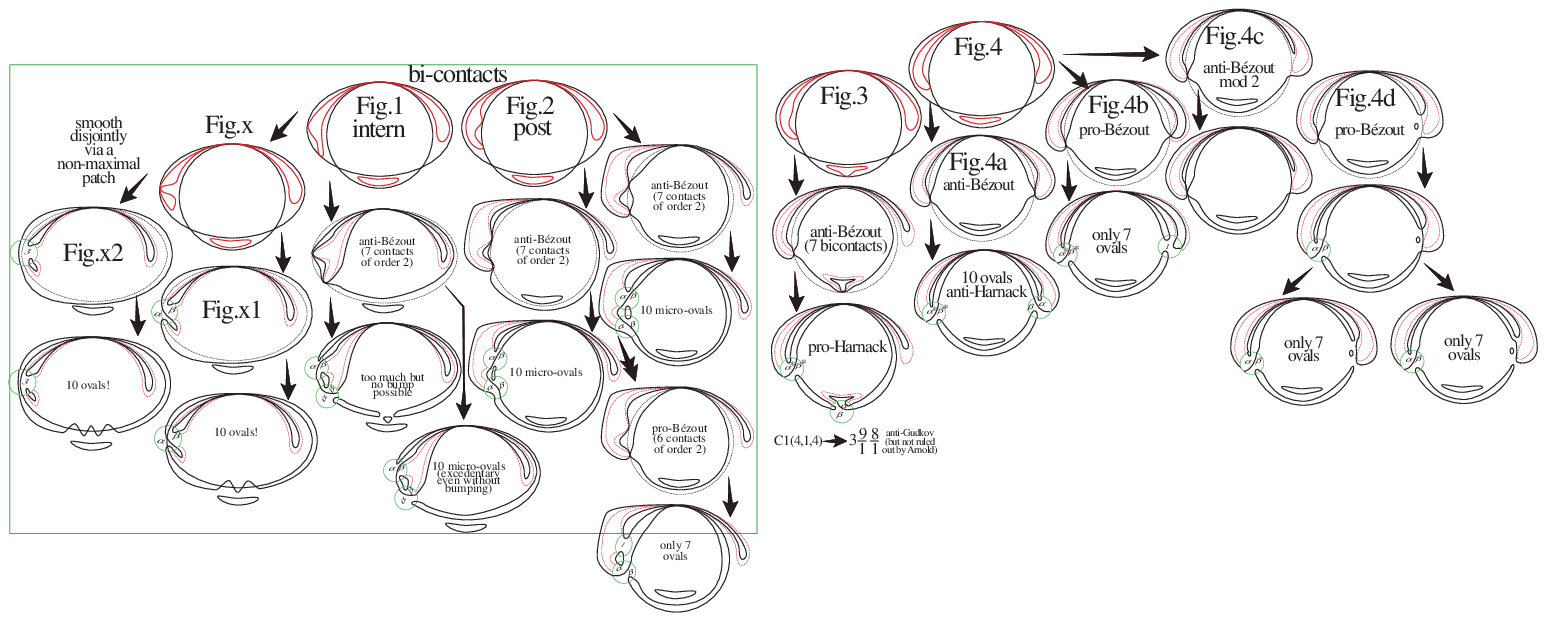,width=172mm}
\captionskipAG
  \caption{\label{ViroDEGREE8_PATCH4SYS_TAN:fig}%
  Bicontacts}
\figskip
\end{figure}

Now we adapt the table to transverse behaviorism
(Fig.\,\ref{ViroDEGREE8_PATCH4SYS_TRAN:fig}). It seems evident
that  transversality  will {\it not} aid attaining $M$-patches.
However the little surprise is that the first Fig.\,1 yields the
patch V(0,4) for $J_{10}$, i.e. three branches with 2nd order
tangency (i.e. what is fundamental to degree 6, and so Gudkov's
solution to Hilbert's problem can be derived along Viro's method,
compare his letter
%%%in Sec.\,\ref{e-mail-Viro:sec}).
in v.2 of Gabard 2013). On the next step of the iteration we
recover the $M$-patch C1(1,8,0) (for $X_{21}$), and also E(0,1,8)
if vibrating the left fringe. Albeit, not novel this is a slight
methodological breakthrough, since it trivializes Viro's method at
the Harnack-Hilbert-Brusotti level involving only
 dissipation of  ordinary double points. Dually, via the
internal vibration (e.g. the model of Fig.\,7) we get first the
$M$-patch of degree six V(4,0) (4 ovals in the bi-lune), and at
the next step of the iteration the patches C1(1,8,0) and E(8,1,0)
depending on the location of the vibration (in the beard or in the
hairs). The corresponding patches are catalogued by green-rhombs
on Fig.\,\ref{ViroDEGREE8_exotic_patches0_SYS:fig}. Alas this
transverse case leads only to a minim proportion of all Viro's
patches and therefore the force of Viro's seems to reside in its
inherent tangential-ness (complicated singularities) as a more
versatile angle of attack upon the optical phenomena allied to
Newton-Hilbert et cie. Hence the power of Viro's method is
unantastbar, contrary to the loose opinion expressed some few line
above.

\begin{figure}[h]\Figskip
%\vskip-1.2cm\penalty0
%\centering
\hskip-2.7cm\penalty0
\epsfig{figure=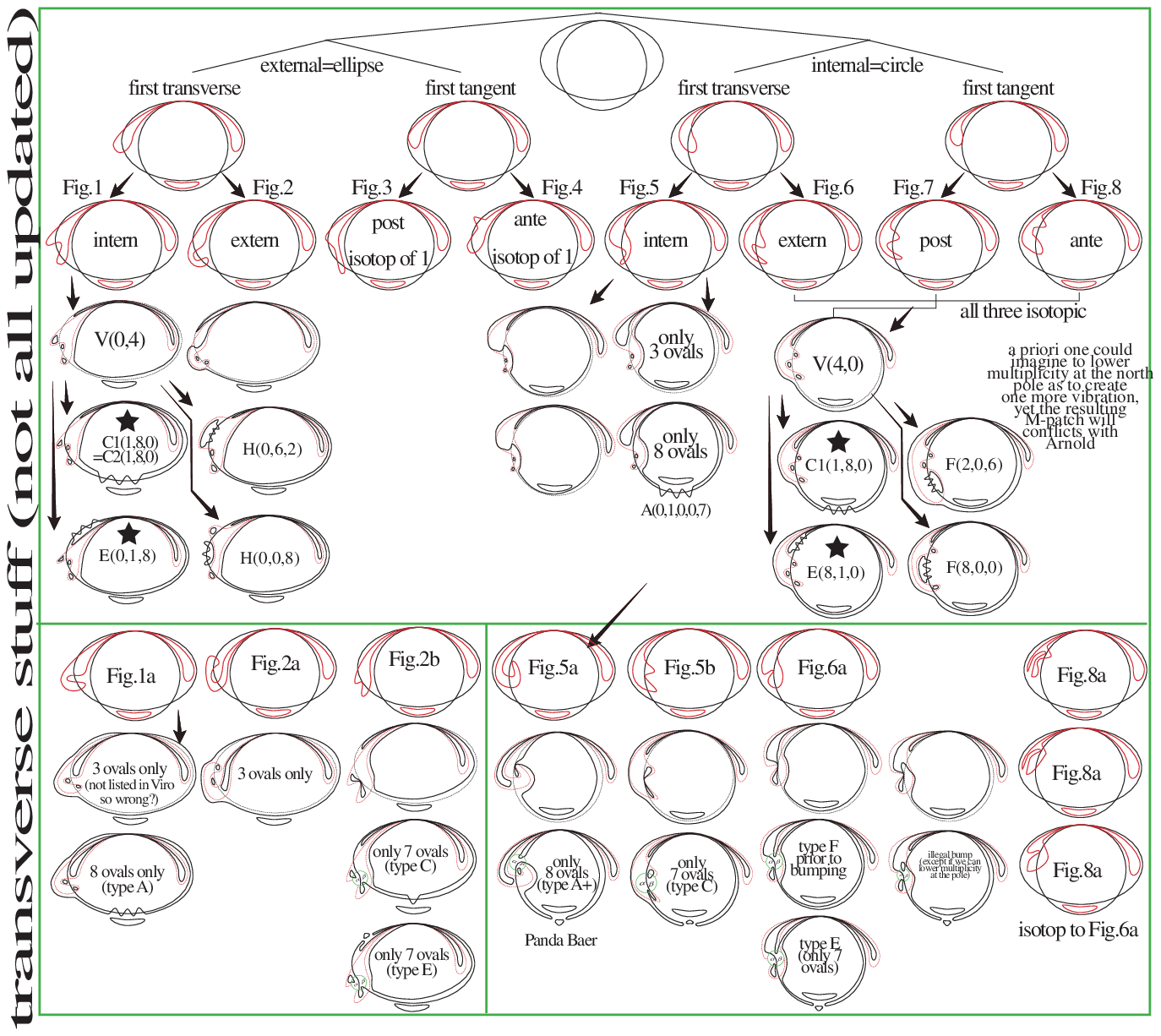,width=172mm}
\captionskipAG
  \caption{\label{ViroDEGREE8_PATCH4SYS_TRAN:fig}%
  Transverse behaviorism}
\figskip
\end{figure}

[08.09.13] Those patches are so-to-speak the most elementary one,
and when doubled produce the schemes $2.C_1(1,8,0)=18\frac{3}{1}$
(Harnack), $2.E(0,1,8)=17(1,2\frac{1}{1})$ (Hilbert), and
$2.E(8,1,0)=1(1,2\frac{17}{1})$ (Hilbert).

Further from Fig.\,1 we can make vibrate one of the oval to get
the patch H(0,6,2) with alas only 8 ovals. If instead the inner
oval is vibrated we get the $(M-1)$-patch H(0,0,8). Doubling those
gives $2.H(0,6,2)=4(1,13\frac{1}{1})$, which is below a Korchagin
scheme, and $2.H(0,0,8)=16(1,1\frac{1}{1})$, which is below a
Chevallier $M$-scheme.

In the overall we see that transversality (code $1+1+1+1$) amounts
to only few patches, while Viro's tangentiality $1+1+2$ leads to
more patches. Hence it seems that the more tangentiality is
reigning, the more flexible is the method (of small perturbation).
In this optic, we have then to analyze the case $1+3$ (one contact
of order three) and $4=4$ (one contact of order four). Besides we
have the case of $2+2$ (double bi-contact), which as we saw seems
to violate Harnack.

Perhaps one can even imagine the situation where both bi-contacts
are located at the same place, the so-called {\it place-to-be\/}
(compare Fig.\,\ref{ViroDEGREE8_PATCH4SYS_TAN:fig}x). Of course,
this looks bizarre as Fig.\,x does not look a small perturbation
of both ellipses, since the red-circuit close to the circle
deviates violently to reach the ellipse. Looking at the
corresponding patch (Fig.\,x1 and below) we get one violating
Harnack's bound with 10 micro-ovals. The variant Fig.\,x2 may look
more clever as it employs a patch disjoint from the ellipse so
that more vibrations can be forced on the $C_6$. Alas, the
ultimate result also involves too many micro-ovals.

Then we were struck by the idea of using a bitangent initial
configuration of ellipses. According to the Viro-style philosophy
that more tangentiality leads to stronger patches, this idea
should not be completely stupid, but alas led us nowhere.

\begin{figure}[h]\Figskip
%\vskip-1.2cm\penalty0
%\centering
\hskip-2.7cm\penalty0
\epsfig{figure=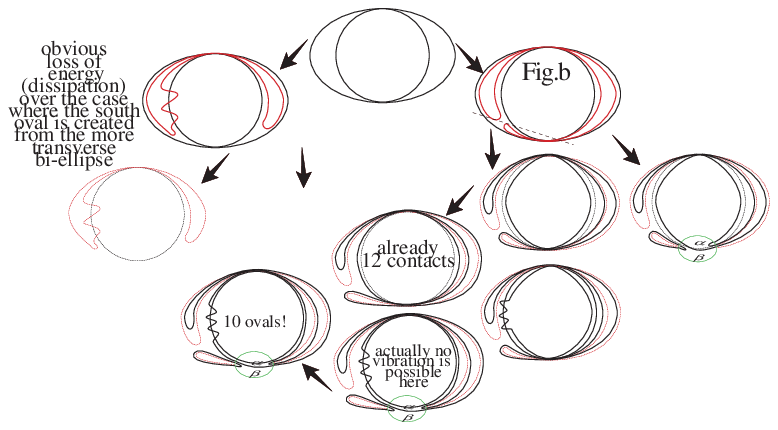,width=172mm}
\captionskipAG
  \caption{\label{ViroDEGREE8_PATCH4SYS_TANG:fig}%
  Transverse behaviorism}
\figskip
\end{figure}

{\footnotesize [22h36, vor dem Einschlafen]:
Sozial-Philosophischer Spruch auf Franz\"oszich:
%(musste sogar den eigenen Vater davon
%\"uberzeugen am vorigen Tag: la vermine du pripet wie der Adolf
%sagte):
La qu\^ete de l'immortalit\'e est un motif suffisant pour
que chacun travaille librement, sans \^etre exploit\'e ni
exclavagis\'e, au projet d'une vie meilleure, infinie, et
affranchie du joug du capitalisme.

[09.09.13] A lost day due to capitalistic duties. (Christa's Konto
nicht f\"ur Ruthli zust\"andig geleistet, Sozialschmarotzern bei
der Bank, usw.)

}

[10.09.13] Next, we noted that the pseudo-quartic $C_4$ of
Fig.\,\ref{ViroDEGREE8_PATCH4SYS_TANG:fig}b probably violates
B\'ezout due to the bitangent along both fringes (cf. the dashed
line).

Next our brain came back again to
Fig.\,\ref{ViroDEGREE8_PATCH4SYS_TAN:fig} where there are too much
(seven) contacts of order two in $C_6 \cap C_2$. It
seems puzzling
 that Figs. 1 or 2 of that plate must exist but they do not lead
to reasonable perturbation in degree 6. Of course the algebraic
realm is flexible in the sense that any (reducible) curve can be
deformed just by perturbing the coefficients. It is puzzling
therefore that we lack as yet any realist perturbation for say
Fig.\,2. Eventually, we discovered Fig.\,2d where the correct
number of 12 contacts is totalized. Here we meet candelabrum-type
singularity with 3 branches. Hence it suffices to know the
dissipation theory of the former to get a patch. A priori, when
there is 3 branches there should be 4 micro-ovals (compare the
patches for $J_{10}$=three tangential branches). But on comparing
with the higher candelabrum with 4 branches (Fig.\,39 in Viro
89/90, p.\,1112), we note that the former accepts at most
$\al+\be=6$ micro-ovals, and not nine like for $X_{21}$ which is
also four branched. By analogy, it seems that the three-branch
candelabrum lacks a smoothing with 4 ovals, as this maximum is
preferably achieved in the purely tangential setting. For a more
intrinsic reason if the 3-candelabrum had a smoothing with 4
micro-ovals, then (whatever their location) applying this patch on
a quintic union of two ellipses plus a transverse line would
create a quintic with 3+4+4=11 circuits (cf.
Fig.\,\ref{ViroDEGREE8_PATCH4SYS_CANDEL:fig} if necessary),
violating frankly Harnack's bound.

So again we have a Viro's style philosophy: tangentiality is the
motor of Harnack-maximality. The contrary  would more readily
serves our purpose, as then Fig.\,2d could perturb to an
$M$-patch.

\begin{figure}[h]\Figskip
%\vskip-1.2cm\penalty0
%\centering
\hskip-2.7cm\penalty0
\epsfig{figure=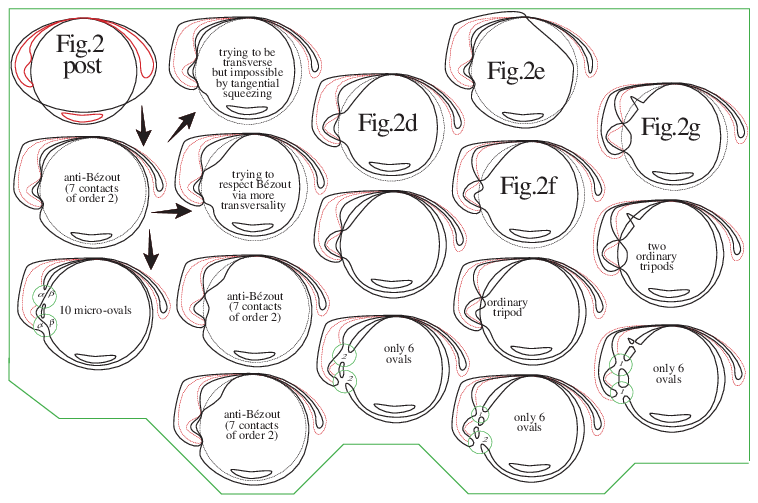,width=172mm}
\captionskipAG
  \caption{\label{ViroDEGREE8_PATCH4SYS_TAN2:fig}%
  Transverse taming of the excedentary contacts}
\figskip
\end{figure}

So we need to understand the  dissipation of the tri-branched
candelabrum. An imperfect attempt is done below
(Fig.\,\ref{ViroDEGREE8_PATCH4SYS_CANDEL:fig}), but we still find
the patch of Fig.\,xx, which is a plausible model in view of the
knowledge  a-priori of $M$-quintics. Alas, when glued in our
earlier Fig.\,2d,  yields only a patch with 6 ovals.

Fig.\,2e exploits the idea of lowering the multiplicity of
intersection by coining transversality at the north pole. Alas,
this destroys the very basic desideratum of getting an
$X_{21}$-patch. Accordingly, this seems to be a cul-de-sac.

\begin{figure}[h]\Figskip
%\vskip-1.2cm\penalty0
%\centering
\hskip-2.7cm\penalty0
\epsfig{figure=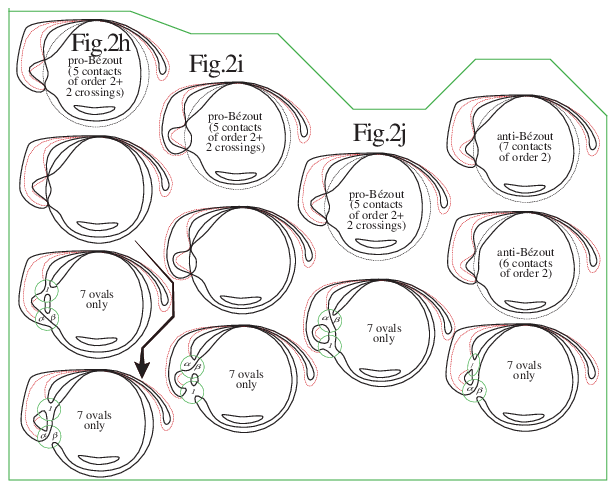,width=172mm}
\captionskipAG
  \caption{\label{ViroDEGREE8_PATCH4SYS_TAN3:fig}%
  Transverse taming of the excedentary contacts (continued)}
\figskip
\end{figure}

Philosophy: The true core of any mathematical truth is a
geometrical picture, hence so-to-speak a physical reality.

Next, we found Fig.\,2h by concealing some transversality as to
respect B\'ezout. So this is the first  perturbation of a
potentially algebraic character, but alas the resulting patches
%%%do have
%%%ceil  at
%%%%% plafonne
reach  only 7 ovals.

Fig.\,2i is an obvious variant,  also procuring only 7 ovals.

\begin{figure}[h]\Figskip
%\vskip-1.2cm\penalty0
%\centering
\hskip-2.7cm\penalty0
\epsfig{figure=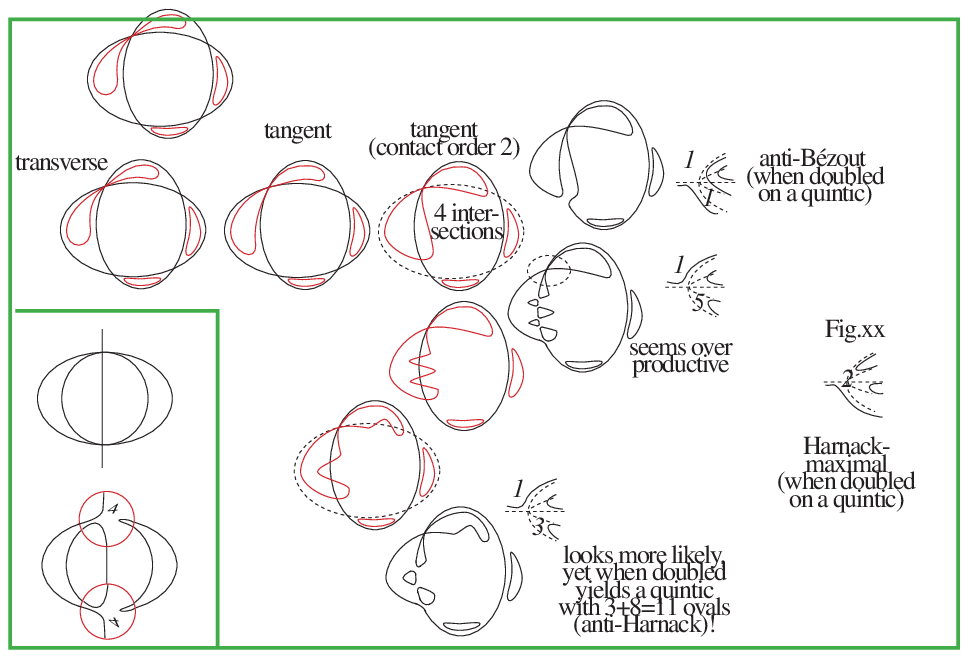,width=172mm}
\captionskipAG
  \caption{\label{ViroDEGREE8_PATCH4SYS_CANDEL:fig}%
  Transverse behaviorism}
\figskip
\end{figure}

At this stage we look blocked and it is somehow disappointing that
we are unable from the $C_4$-configuration  with a bi-contact  to
produce $M$-patches for $X_{21}$.

[11.09.13] So our problem is still the same: can we produce more
patches with Viro's elementary method of perturbation? In
particular can we get all of Viro's patches including those
obtained via the more tricky methodology of hyperbolism (Huyghens,
Newton, Cremona, Gudkov, Viro).

If yes, which sort of initial geometric configuration for $C_4\cap
C_2$ has to be employed? Apparently the case of bicontact $4=2+2$
leads nowhere.

It is only at this moment, that we were flashed by the simple idea
that if the bicontacts of
Figs.\ref{ViroDEGREE8_PATCH4SYS_TAN:fig}(1,2) are over-productive
(10 micro-ovals) one can just consider the situation of Fig.\,3
where some of the energy is lost by splitting the bicontacts over
two different circuits. It results indeed then $M$-patches based
on the geometry of a smiling-face. Precisely, we get when both
$\al$ and $\al^\ast$ are 4 the patch C1(4,1,4). This alas seems to
violate Gudkov periodicity. Another objection is that the involved
pair $C_6\cap C_2$ does not respect B\'ezout, as we see 7
bicontacts between the sextic and the conic (circle). For $\al=4$
and $\al^\ast=0$ we find, etc... but by experience those will also
certainly corrupt Gudkov.

Next, we have also Fig.\,\ref{ViroDEGREE8_PATCH4SYS_TAN:fig}(4)
where both fringes have a bicontact (contact of order 2). Again
our perturbation $C_6$ cannot be a genuine algebraic one, as the
intersection with $C_2$ involves seven bicontacts violating
B\'ezout. Hence in Fig.\,4 exists in the algebraic category (and
there is no B\'ezout obstruction to this), we infer that the a
perturbation of $C_4\cup C_2$ must have another look. For instance
we have the perturbation of Fig.\,4b which respects B\'ezout and
involves a transverse behavior on the right fringe, but still a
tangential one along the left fringe. Alas, making this
concession, we get only a patch with 7 ovals instead of the nine
ones requested to reach Harnack-maximality.

Next, we can imagine the perturbation of Fig.\,4c, but this
appears to corrupt B\'ezout modulo 2. But this little defect may
be corrected by Fig.\,4d where we better into account the
non-traversing issue about the branch of
highest curvature through
the north pole.

It seems clear at this moment that we have exhausted the power of
the method, or to speak frankly, we missed any single $M$-patch
from the bicontact trick. Alas we do not know if this is due to
the incompetence of the writer (Gabard), or an intrinsic
state-of-affairs.

\subsection{Study of the contact $3+1$}

[11.09.13] Let us know look at the case of a contact of type $3+1$
as shown on Fig.\,\ref{ViroDEGREE8_PATCH_3+1:fig}. Fig.\,1 and
Fig.\,2 shows two ways of having a contact of order 3 between a
quartic $C_4$ and a conic $C_2$. By perturbing $C_4\cup C_2$ we
get the sextic of Fig.\,1b. Ten there is some conceptual
difficulties of how to interpret the picture in terms of Puiseux
(?) branches. We mean basically that either Fig.\,1b1 or 1b2 could
occur where the numbers 2,3 label a given branch, while its
specific value measure the contact with the ground ellipse. Then
there are other several conceptual obstacle but a rapid run
through the lead us to the conclusion that the local singularity
involved (say $X_{31}$ improvising notation) should accept
dissipations with 6 micro-ovals at most (as suggested by the count
on Fig.\,x). Unfortunately, applying such a dissipation to
Fig.\,1b1 seems to create only 8 micro-ovals.

\begin{figure}[h]\Figskip
%\vskip-1.2cm\penalty0
%\centering
\hskip-2.7cm\penalty0
\epsfig{figure=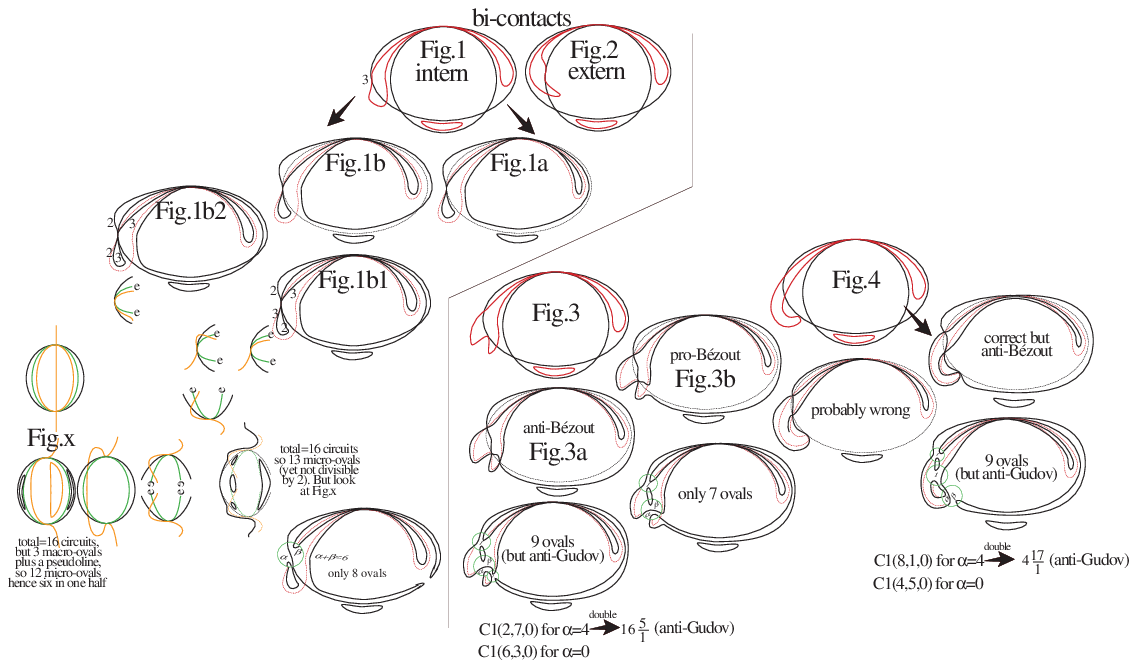,width=172mm}
\captionskipAG
  \caption{\label{ViroDEGREE8_PATCH_3+1:fig}%
  Contact $3+1$}
\figskip
\end{figure}

[12.09.13] Our next idea is materialized by
Fig.\,\ref{ViroDEGREE8_PATCH_3+1:fig}(3), where  like Viro we
exploit the idea of keeping all singularities unsmoothed until the
end in the hope that accumulated tension  when ultimately
liberated will act as a devastating flood offering a real
breakthrough on Hilbert's 16th.

To be concrete we get first Fig.\,3a, but the intersection
$C_2\cap C_6$ involves 5 bicontacts plus 4 crossings, hence a
total of 14 points (counted by multiplicity), violating B\'ezout.
Incidentally, the resulting patches corrupt Gudkov periodicity.
Fig.\,3b shows a variant where B\'ezout is respected, but the
resulting patch exhibit only 7 ovals,and so does not arise much
interest. Next we tried Fig.\,4, but again the intermediate
product violates B\'ezout, and the final patch corrupt Gudkov
(periodicity).

Our next idea is to consider an avatar of the previous systematic
figure implementing Viro's vibrational method (acronym VVM) by
passing the less curved branches outside instead of inside as in
the original method. This gives us
Fig.\,\ref{ViroDEGREE8_PATCH4SYST:fig} where we were only able to
reach 8 ovals in the erectile case, and only 7 in the invaginated
case. So it seems evident that there is a loss of energy by using
those inflated version of the earlier main-figure.

\begin{figure}[h]\Figskip
%\vskip-1.2cm\penalty0
%\centering
\hskip-2.7cm\penalty0
\epsfig{figure=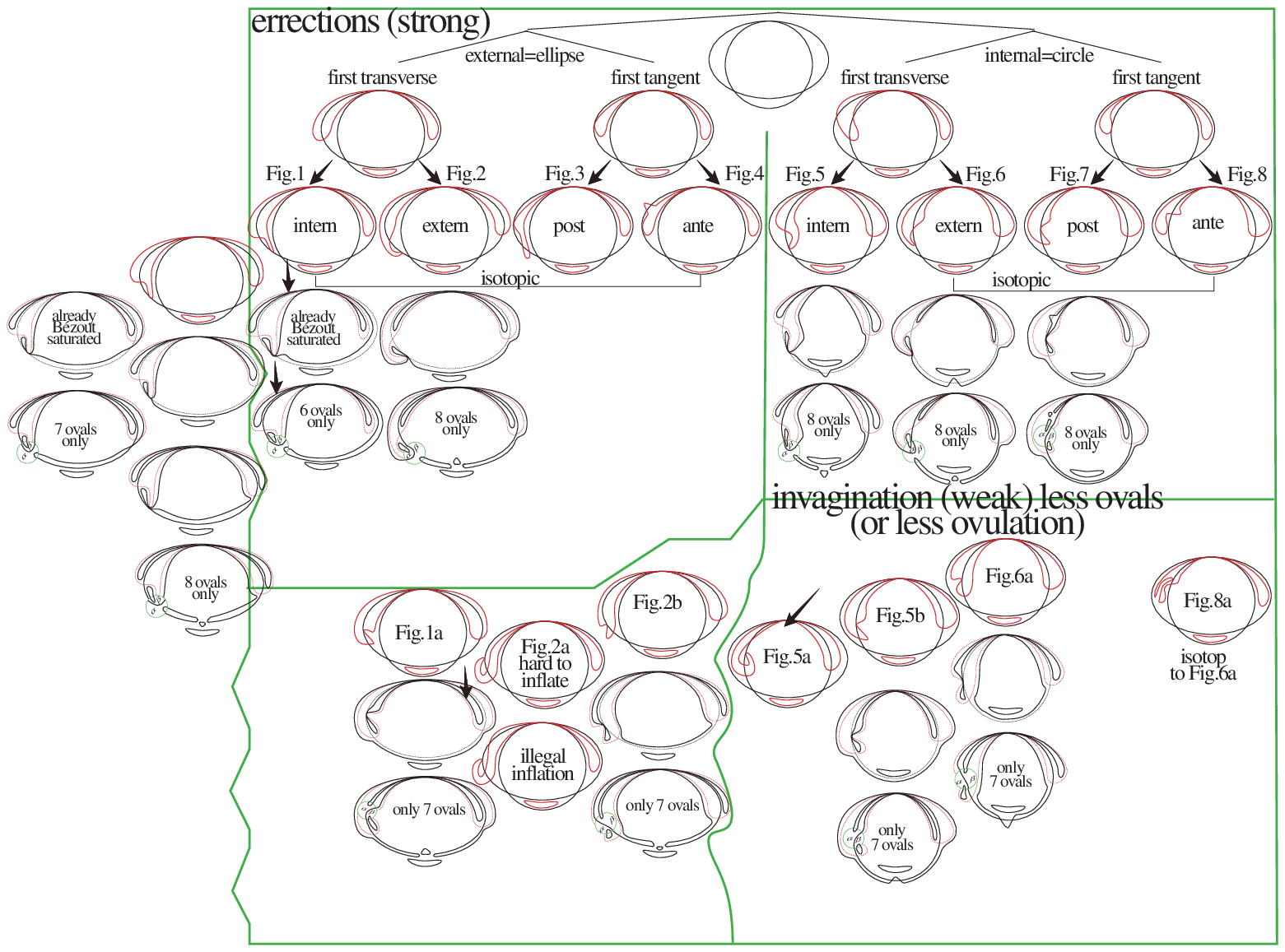,width=172mm}
\captionskipAG
  \caption{\label{ViroDEGREE8_PATCH4SYST:fig}%
  Inflated (or deflated) variants of the main-table, yet
  with a loss of one oval over the primordial construction}
\figskip
\end{figure}

At this stage one is clearly lost in a sterile labyrinth, and some
clairvoyance is requested to come out of it alive.

[21h42, 12.09.13] Vor dem Einschlafen, we were flashed by the
following modest idea which we only report due to our  deep level
of depression. The idea is to consider a configuration of ellipses
somewhat more transverse than Viro's namely that depicted below
(Fig.\,\ref{ViroDEGREE8_ellipse:fig}). After more lucid
experimental thinking it seems evident a priori to us that this
lack of tangentiality will be incapable reaching
Harnack-maximality. Nonetheless it would be nice to know the
maximum number of ovals accessible by mean of a perturbation of
such a configuration. A naive drawing suggests the answer being
16, because we see on the picture below 8 ``macro'' ovals and 8
micro-ovals coming from the  $\al, \be $
 parameters. Needless to say our deception is great, but was
 somehow anticipated by our experimental knowledge. This is
 another typical illustration of the fact that too much
 transversality (between low degree objects)
 %%%%nuit (de nuire)
 impedes Harnack-maximality.

\begin{figure}[h]\Figskip
%\vskip-1.2cm\penalty0
%\centering
\hskip-2.7cm\penalty0
\epsfig{figure=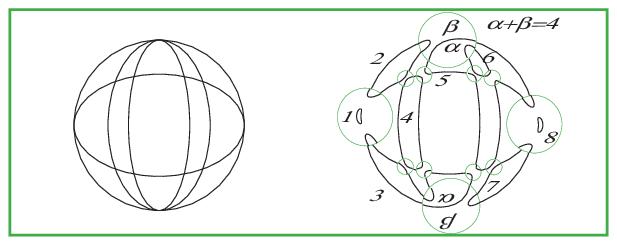,width=172mm} \captionskipAG
  \caption{\label{ViroDEGREE8_ellipse:fig}%
  Another variant of Viro}
\figskip
\end{figure}

\subsection{Falling in love with Oleg Yanovich}

[04.09.13] Of course in the above construction of O.\,Ya. Viro it
should be noted that it exploited the only two possible
dissipations of $J_{10}^-$ (which are precisely those involved in
Hilbert's 16th in degree $m=6$, and the only possible by virtue of
the Hilbert-Rohn-Gudkov obstructions: compare the list of
$M$-patches for $J_{10}^-$). Further, it is pleasant to
contemplate the ``gigogne'' (telescopic) nature of the dissipation
theories of all those singularities as a sort of big inductive
process. (Hilbert would say a Einschachtelung: i.e. to understand
the dissipation of the quadruple point with 2nd order tangency, we
rest on same knowledge for the triple point.) Notwithstanding,
Viro's method appears as split into the method for class $C$ based
split sextics of bidegree $1+5$ (so-called affine $M$-quintic) and
now a method of tangential vibration (somewhat reminiscent of
Hilbert's method), which furthermore avoids the usage of
hyperbolism. Hence a basic task could be to find a more unified
treatment of all Viro's patch by one and the same method.

Again we insist that on Fig.\,\ref{ViroDEGREE8_PATCH4_Eclass:fig}
there is no ore freedom in the dissipation of the triple-point
with 2nd order tangency than those tabulated because otherwise we
would get sextics violating the Hilbert-Rohn-Gudkov census (HRG).

Now let us explore some variant of Viro's (last) E-method, as we
shall call it, since it produce patches of type E (i.e. with a
triple lune).

One of our idea was one the previous Fig.\,a to move the location
of the oscillation below the tacnodal singularity $A_3^{-}$.
Another idea (simpler to depict) is just to permute the location
of the oscillation and that of the tangency on the preliminary
quartic $C_4$, see Fig.\,x below. Disappointingly, the resulting
patch has only $8$ micro-oval (not an $M$-patch), yet still of the
interesting class A which is thus non-empty outside of the maximal
realm. (Note: one of this patch namely $A(1,0,0,0,7)$ is nearly
equal to the one we speculated about yesterday).

\begin{figure}[h]\Figskip
%\vskip-1.2cm\penalty0
%\centering
\hskip-2.7cm\penalty0
\epsfig{figure=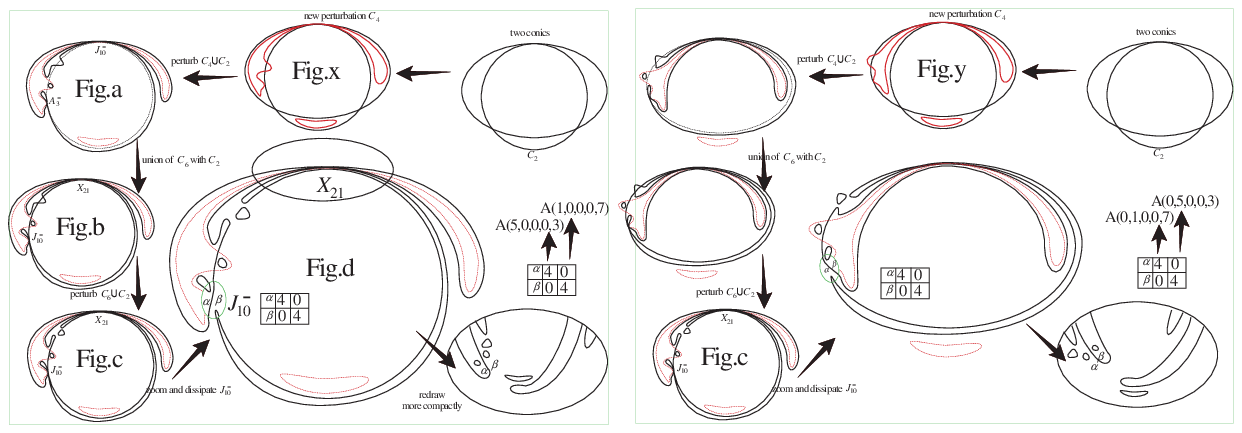,width=172mm}
\captionskipAG
  \caption{\label{ViroDEGREE8_PATCH5_Eclass:fig}%
  Viro's E-class (trinested lune)}
\figskip
\end{figure}

Next one can imagine the same game yet oscillating against the
other ellipse, cf. Fig.\,\ref{ViroDEGREE8_PATCH5_Eclass:fig}y.
This gave us the patch A(0,1,0,0,7) which is by virtue of an
obvious mirror-symmetry on the first 4 parameters
$\al,\be,\ga,\de$ equivalent to A(0,0,1,0,7), which is exactly the
patch on which we speculated yesterday. Of course the deception is
that this is still not an $M$-patch, but there is certainly yet
another variant where we get maximality. The vague experimental
qualification (based on month of experimentation with the Harnack,
Hilbert and Viro methods) prompts that usually damped dissipations
comes from a division of the oscillating energy, and therefore a
loss of newly created ovals (hence ``dissipation of energy''). Put
more concretely, this suggests looking at
Fig.\,\ref{ViroDEGREE8_PATCH6_Eclass:fig}z where the tangency
occurs more internally. Repeating Viro's trick (algorithm adapted
to this situation) we get the two $M$-patches $E(0,1,8)$, and
$E(0,5,4)$ that where precisely the last two patches claimed by
Viro that missed us as yet. Of course all this harmony and duality
(Wechselbeziehung) is highly reminiscent of Hilbert's method where
depending upon an internal or external vibration we get either
Harnack's curves or Hilbert's. In the present Viro context it
seems that is a contiguity between the tacnode and the oscillation
that permits to maximize the vibratory energy and therefore the
number of springing ovals.

\begin{figure}[h]\Figskip
%\vskip-1.2cm\penalty0
%\centering
\hskip-2.7cm\penalty0
\epsfig{figure=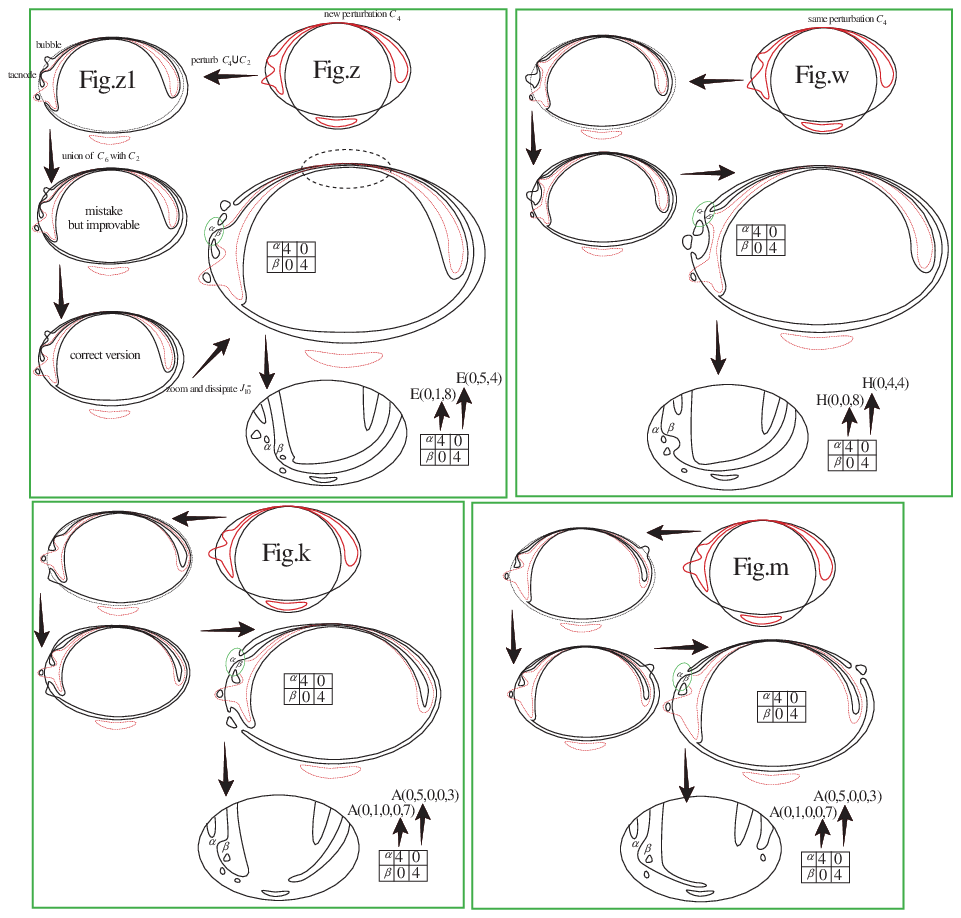,width=172mm}
\captionskipAG
  \caption{\label{ViroDEGREE8_PATCH6_Eclass:fig}%
  Viro's E-class (trinested lune)}
\figskip
\end{figure}

At this stage we have a complete knowledge of Viro's theory (as
far as $X_{21}$ is concerned), yet there is still the hope that
his assertions are not exhaustive that his method can be further
varied as to give more patches. At least this the naive hope and
we shall try to explore some additional possibilities (requesting
poor level of imaginativeness).

One of the idea to explore is to vary the position of the
``bubble'' of Fig.\,z1, by placing it rather in the loop below the
tacnode. Of course this is somewhat hard to depict, and so let us
retrace Fig.\,z as Fig.\,w to get more free-room at the critical
place. (We may so expect to find: a variant of Viro's method
leading to the materialization of some bosons in degree 8.) Alas
the resulting patch is not maximal and belongs to type H (in the
notation of Fig.\,\ref{ViroDEGREE8_exotic_patches0_SYS:fig}). More
specifically it realizes the patches H(0,0,8) and H(0,4,4) (where
as usual we count the ovals' number $\al, \be, \ga$ from inside to
outside).

Another idea is to place the vibratory energy of the bottom oval
of the quartic $C_4$ as shown on
Fig.\,\ref{ViroDEGREE8_PATCH7_Eclass:fig}a. Alas the resulting
patches though belonging to the interesting class I, have only 7
``micro-ovals'' and so are $(M-2)$-patches (i.e. 2 units below the
maximum possible). Remind beside, that the interest of class I is
that it is not much obstructed (a priori at least), and it could
lead to the creation of bosons as summarized on
Fig.\,\ref{ViroDEGREE8_exotic_patches0_SYS:fig}. {\it Optional
note:} If we double the patch I(7,0,0) we get the $(M-4)$-scheme
$(1, \frac{16}{1})$) which 4 steps below Shustin's
(anti)-$M$-scheme $(1,2 \frac{18}{1})$.

Perhaps for a more clever placement of the bump we may at least
create an $(M-1)$-patch of type~I. This suggests the variant
Fig.\,b1 where the bump is placed on the left of the tacnode.
First we trace Fig.\,b which is the same as Fig.\,a modulo leaving
some more room at the place where intend to bump the curve.
(General comment: as usual in those iterative methods \`a la
Harnack, Hilbert, Gudkov, Viro just to name the Gods) we are
permitted at each step the depiction to alter slightly the
metrical proportions as to emphasize the topological properties of
the curve under construction. This may be interpreted as an
elasticity of the underlying ether.) For Fig.\,b we find then
$(M-1)$-patches indeed, yet of type A, precisely A(0,0,0,5,3) and
A(0,0,0,1,7). Doubling the first gives the $(M-2)$-scheme
$9\frac{10}{1}$ (which is not extremely exciting).

Another option involves placing the bump on the very right of the
scene, and this yields Fig.\,c, which creates ultimately the
patches of type A with an island, hence denoted A+. Specifically,
we get A+(5,2) and A+(1,6), but those are alas only
$(M-2)$-patches.

\begin{figure}[h]\Figskip
%\vskip-1.2cm\penalty0
%\centering
\hskip-2.7cm\penalty0
\epsfig{figure=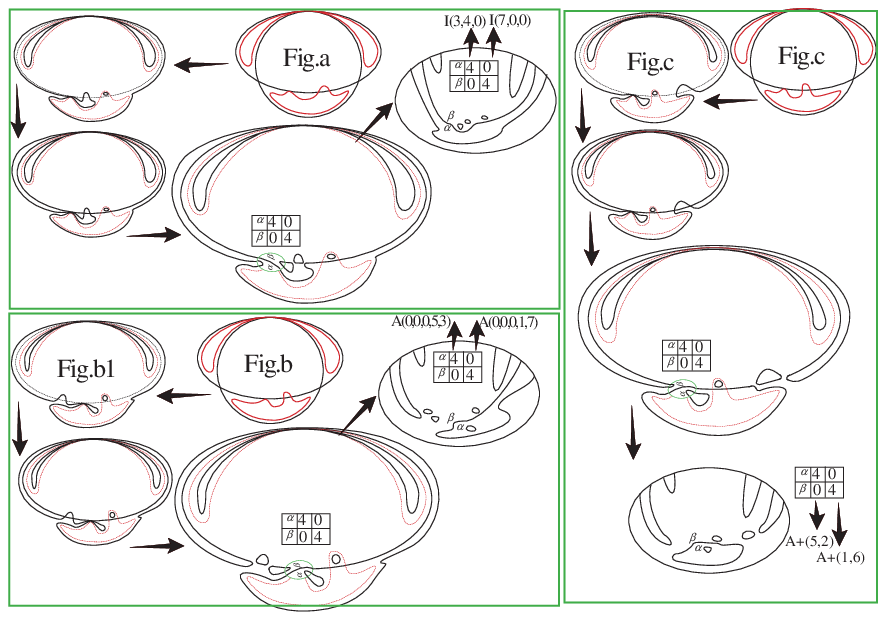,width=172mm}
\captionskipAG
  \caption{\label{ViroDEGREE8_PATCH7_Eclass:fig}%
  Viro's E-class (trinested lune)}
\figskip
\end{figure}

Besides it remains thereafter the option of doing an outer
vibration across the ellipse which is circle-like. This brings us
to Fig.\d below. Of course all this series of pictures where the
central bottom oval is vibrated never yield $M$-patches, in view
of the fact that if left tranquil this oval  contributed to one of
the nine micro-ovals. Specifically Fig.\,d only gives an
$(M-2)$-patch of type A (vut for which we have alas no canonical
naming as yet).

Next we can as before continue the game by varying the location of
the bump, and this produces Fig.\,e and Fig.\,f. First, Fig.\,e
gives only $(M-2)$-patches, yet of the not yet found type B2.
Finally Fig.\,f, where the bump will be placed on the left. Again
we only get $(M-2)$-patches, actually the same as those already
derived via (the previous) Fig.\,c.

\begin{figure}[h]\Figskip
%\vskip-1.2cm\penalty0
%\centering
\hskip-2.7cm\penalty0
\epsfig{figure=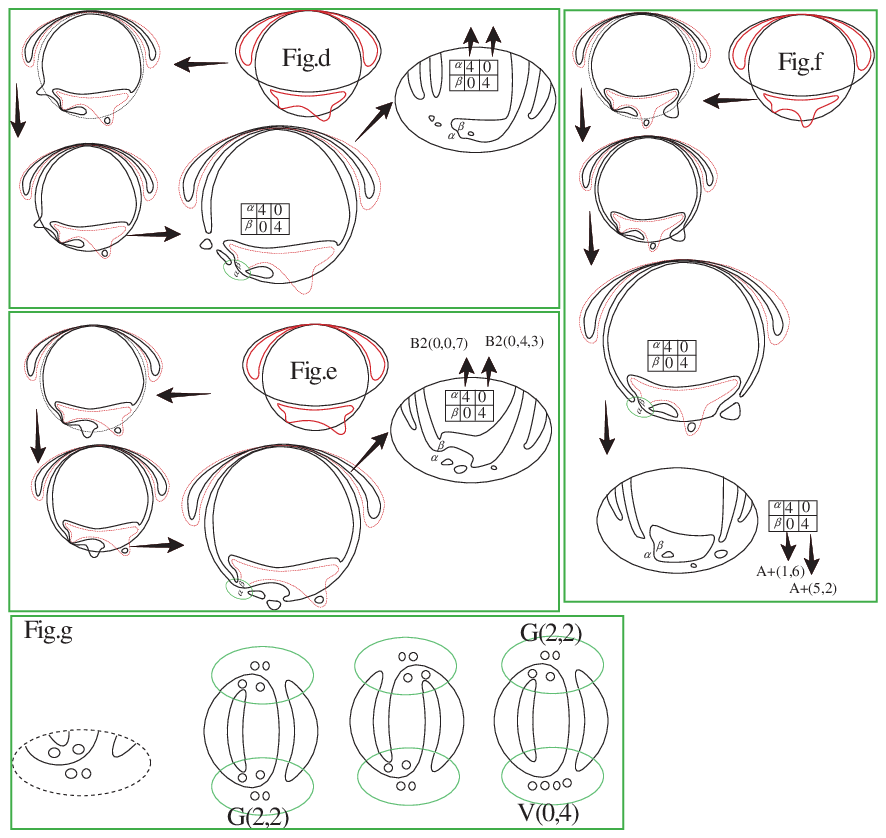,width=172mm}
\captionskipAG
  \caption{\label{ViroDEGREE8_PATCH8_Eclass:fig}%
  Viro's E-class (trinested lune)}
\figskip
\end{figure}

At this stage we are slightly disappointed, yet we should
definitively study all the possibilities on a more
charitable/peaceful day, as today some financial stress is
apparent. As said Hassler Whitney let us do research naturally.

Then we can imagine variants of the very first constructions where
the bump is placed on the other side of the circle, yet  this will
certainly leads to a loss of the vibratory energy.

Besides, after a sleep we were flashed by the dubious idea that
perhaps we have not yet exploited the most general dissipation of
the triple points with 2nd order tangencies (double contacts for
short), alias $J_{10}^{-}$, in  abstruse Arnold's notation. Our
idea is just to imagine Gudkov's sextic curve $5\frac{5}{1}$ as
split into two patches of the type depicted on
Fig.\,\ref{ViroDEGREE8_PATCH8_Eclass:fig}g. A priori we knew no
obstruction against the existence of such a patch, yet after some
minute of thinking resurfaced in our memory Fiedler's enhancement
by a factor 2 of Gudkov periodicity in the case of symmetric
$M$-curves. This Fiedler's result ($\chi\equiv k^2 \pmod 16$)
obstructs the existence of the posited patch for F3=$J_{10}^-$
(flat singularity with 3 branches),and is as far as we know the
sole obstruction against this patch. Still, it seems of interest
to look at which sort of patches results for F4=$X_{21}$ assuming
that this patch exist for F3=$J_{10}$ (we omit the minus from
Viro-Arnold's notation as there is no ambiguity in principle).
However some more basic thinking (that just escaped from our
memory) shows that when gluing Viro's patch V(0,4) with our patch
G(2,2), we get the $M$-scheme of degree 6 with symbol
$7\frac{3}{1}$ corrupting Gudkov periodicity. So our patch cannot
coexist with Viro's, and since the latter exists, our does not.
(Hence Fiedler's periodicity can be circumvented.)

Next coming back to the previous Fig.\,z, it seems that there is a
3rd option consisting in placing the bump on the bottom as shown
on Fig.\,\ref{ViroDEGREE8_PATCH6_Eclass:fig}k, yet as expectable
we get only an $(M-1)$-patch along this route. So it seems evident
that to reach an $M$-patch the bump must really sit between the
tacnode and the $X_{21}$-singularity (sorry for this vagueness).

For the sake of exhaustiveness, let us study the case of a bump
placed on the right circuit of Fig.\,z1. This is shown on
Fig.\,\ref{ViroDEGREE8_PATCH6_Eclass:fig}m, and gives actually the
same patches as the just studied Fig.\,k, yet still no additional
$M$-patches at the horizon.

Of course the study  can be continued along any more systematic
way, yet it seems also clear that we (or rather Viro) have
extracted all the juice from his method.

\subsection{Toward a classification of all patches for
$X_{21}$ via knot theory}

[14.09.13] Today only, we were flashed by the idea that
potentially much obstruction on patches could come from knot
theory alone. The idea is the simple trick of what is known as the
Milnor fibre, or rather just the link of  a singularity, a trick
going back to Wirtinger and earlier. For $X_{21}$, if we intersect
the singularity by a small $3D$-sphere centered at the singularity
we see four branches and thus a link with 4 components all
unknotted, but maybe collectively linked.

Remind that the basic normal crossing $x.y=0$ determines by this
recipe Hopf's link, and likewise it would of interest to know
which link of $S^3$ is determined by $X_{21}$. Evidently it must
be like a catenary yet with all items interlaced.

Next the fascinating issue must be to interpret the complexified
patch as  a bordered Riemann surface filling smoothly this link
inside of $B^4$. So this bordered surface has 4 contours, and if
we imagine the 4 ellipses as spheres creating 3 holes when pieced
together, it remains a genus of $(21-3)/2=18/2=9$, i.e. nine, for
the membrane of the patch.

Further the patch is invariant under complex conjugation and this
shall give the number of ovals in case of an $M$-patch. All this
context looks sufficiently fine and explicit that via some
knowledge of topology the usual combinatorial tricks \`a la
Poincar\'e, Rohlin, etc. should yields obstruction on the patches
directly, without having to refer to global obstructions in $\CC
P^2$ coming from the Rohlin-Fiedler-Viro theory.

Naively (or abstractly) the symmetric patch may be visualized as a
bordered surface with symmetry.

\begin{figure}[h]\Figskip
%\vskip-1.2cm\penalty0
%\centering
\hskip-2.7cm\penalty0
\epsfig{figure=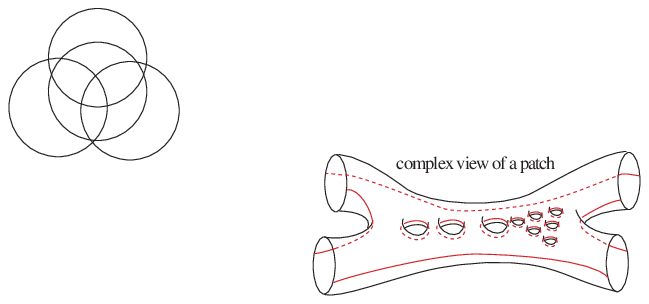,width=172mm}
\captionskipAG
  \caption{\label{ViroDEGREE8_PATCH_LINK:fig}%
  Link theoretical view of the patches}
\figskip
\end{figure}

Apart from possible purely topological obstructions as we said in
the realm of combinatorial topology that could reproduce those of
Arnold, Gudkov that we already derived via doubling the patch, it
could be that highbrow differential geometric obstruction arise
when considering the patch as  a minimal surface filling the link.

So Viro's theory of the patches ought to be directly connected
with soap film experiments in 4D-space, namely the ball $B^4$.

The real challenge would be to prohibit all (or at least some) of
the $M$-patches not constructed by Viro, and not already
prohibited by classical or sporadic obstructions (compare for the
present state-of-the-art our
Fig.\,\ref{ViroDEGREE8_exotic_patches0_SYS:fig}).

Of course even if we could prohibit a patch (by this or another
method) then it would not directly result a prohibition of the
corresponding scheme deduced by doubling, because it can be
imagined that the scheme in question can be constructed by another
method than via Viro's quadri-ellipse. Still, one may expect that
at least if the curve is symmetric under a mirror involution then
there is more global splitting of the curve into two patches. This
reminds slightly Fiedler's ideas on symmetric curve, but we are
too tired to make any direct connection.

Further more even if we could prohibit a patch then it could be
that its double (which is  an $M$-scheme) is realizable by another
combination of patches, so that nothing could be inferred.

Still, we believe that obtaining a complete classification of all
$M$-patches is a problem of independent interest, certainly much
intermingled with Hilbert's 16th (perhaps the key to its ultimate
solution in degree 8) in case Viro's method has not yet been
exploited to its full regime. Alternatively it can be the case
omniscience of all patches is not even enough wisdom to fix
Hilbert's 16th (for $m=8$), in case Viro's method is not flexible
enough to reach some of the bosons (assuming there existence). Put
more intrinsically, it could be that bosonic curves (if the
materialize) are far away from the quadri-ellipse. In this
scenario of bosonic chambers far from the quadri-ellipse (let us
speak of ghost curves) even a complete knowledge of all
$X_{21}$-patches would not suffice to fix Hilbert's 16th.

By the way ghost chamber seem to exist as exemplified by Viro,
Shustin, Korchagin, Chevallier, Orevkov exotic construction not
readily based on a perturbation of the quadri-ellipse. Accordingly
ghost schemes may be a reality and impede solving Hilbert's
problem purely in the vicinity of the quadri-ellipse.

In contrast, forcing a bit the passage, it is still conceivable in
a more generous world of patches that all bosons are constructible
via perturbation of the quadri-ellipse.

\iffalse
{\it Micro-slogan}: Think geometrically and prove
algebraically, we read in dubious book (LNM by Kirby on
4-manifolds). We know at least since Thurston, that the right
method is: Think geometrically and prove geometrically.

\fi

\subsection{A naive idea for prohibitions via satellites}

[13.08.13] In principle, if a plane curve is dividing then its
satellite(s) should also be dividing, because it will be totally
swept out by the same total pencil. In abstracto, this thesis is
hard-to-defend because we may lack a concrete pencil realizing the
total map whose existence is abstractly granted by Ahlfors'
theory. Yet in the case of $M$-curves, total reality is trivially
granted (compare Gabard 2013B
\cite{Gabard_2013B-Riemann's-flirt}). So a naive method of
obstruction could be to look at a curve and build its doubled
satellite, which being dividing has to verify Rohlin's complex
orientation formula. Additionally, all circuits which are doubled
are circulated in the same sense due to a naive dextrogyration
occurring within a real tubular neighborhood of the curve, and
therefore form negative pairs of ovals in Rohlin's sense.

We can try to take any of the boson in degree 8, and satellites it
to get a curve of degree 16 with $2\cdot 22=44$ ovals and with
complex orientations partially controlled by the above rule. On
this sixteen-tics $C_{16}$ we may hope to get sometimes trouble
with Rohlin's formula. If the method does not readily apply to the
boson themselves, it may at least perhaps re-explain some of
Viro's obstructions (granting their correctedness of course).

All this request some investigation, yet the methodology looks a
bit overnaive to lead somewhere. Imagine first the boson
$1\frac{1}{1}\frac{18}{1}$. When doubled all ovals can be imagined
as coupled in infinitesimal pairs, each of which is a negative
pair in the sense of Rohlin, because it is swept out in the same
sense by the total pencil of sextic (cf. again Gabard 2013B
\cite{Gabard_2013B-Riemann's-flirt} for the fact that total
reality of an $M$-curve of degree $m$ is always flashed by a
pencil of $(m-2)$-tics).

Denote by $C_{16}$ the curve of degree 16 doubling the octic
$C_8$. It is dividing, with $r=44$ ovals and subsumed to Rohlin's
formula $2(\pi-\eta)=r-k^2=44-8^2=-20$. Further, Hilbert's tree of
the doubled scheme is for the boson $1\frac{1}{1}\frac{18}{1}$ as
depicted below (Fig.\,\ref{Double:fig}b), on which we may count
the total number of pairs $\pi+\eta$ regardless of (complex)
orientations. The count is effected by breaking along the 3
obvious components and according to the length $\ell $ of the
pair, to get
$$
\pi+\eta=
\begin{cases} 1 \quad &\textrm{of } \ell=1 \cr
0  \quad &\textrm{of } \ell=2\cr
0  \quad &\textrm{of } \ell=3
\end{cases}
+
\begin{cases} 3 \quad &\textrm{of } \ell=1 \cr
2  \quad &\textrm{of } \ell=2\cr
1  \quad &\textrm{of } \ell=3
\end{cases}
+
\begin{cases} 1+2.18 \quad &\textrm{of } \ell=1 \cr
2.18  \quad &\textrm{of } \ell=2\cr
18  \quad &\textrm{of } \ell=3
\end{cases}=1+6+1+5.18=98.
$$
Hence, $2\pi=88$, whence $\pi=44$.

On the other hand we may calculate $\pi$ using Hilbert's tree with
signed branches prescribed by complex orientations. For this
purpose we use the evident signs-law
%%%%%(\ref{Signs-law:lem})
%%%H16 mask desactivated
to the effect that mixing the gene is good while consanguinity is
bad. Further it is essential to know that during the doubling
process Hilbert's tree grows but has a sort of trunk incarnating
the original undoubled curve, and the latter complex orientations
can be transplanted, while all other newly created branches of the
doubled tree are just negatively charged by the dextrogyration
argument (sketched right above).

\begin{figure}[h]\Figskip
%\vskip-1.2cm\penalty0
%\centering
\hskip-2.7cm\penalty0 \epsfig{figure=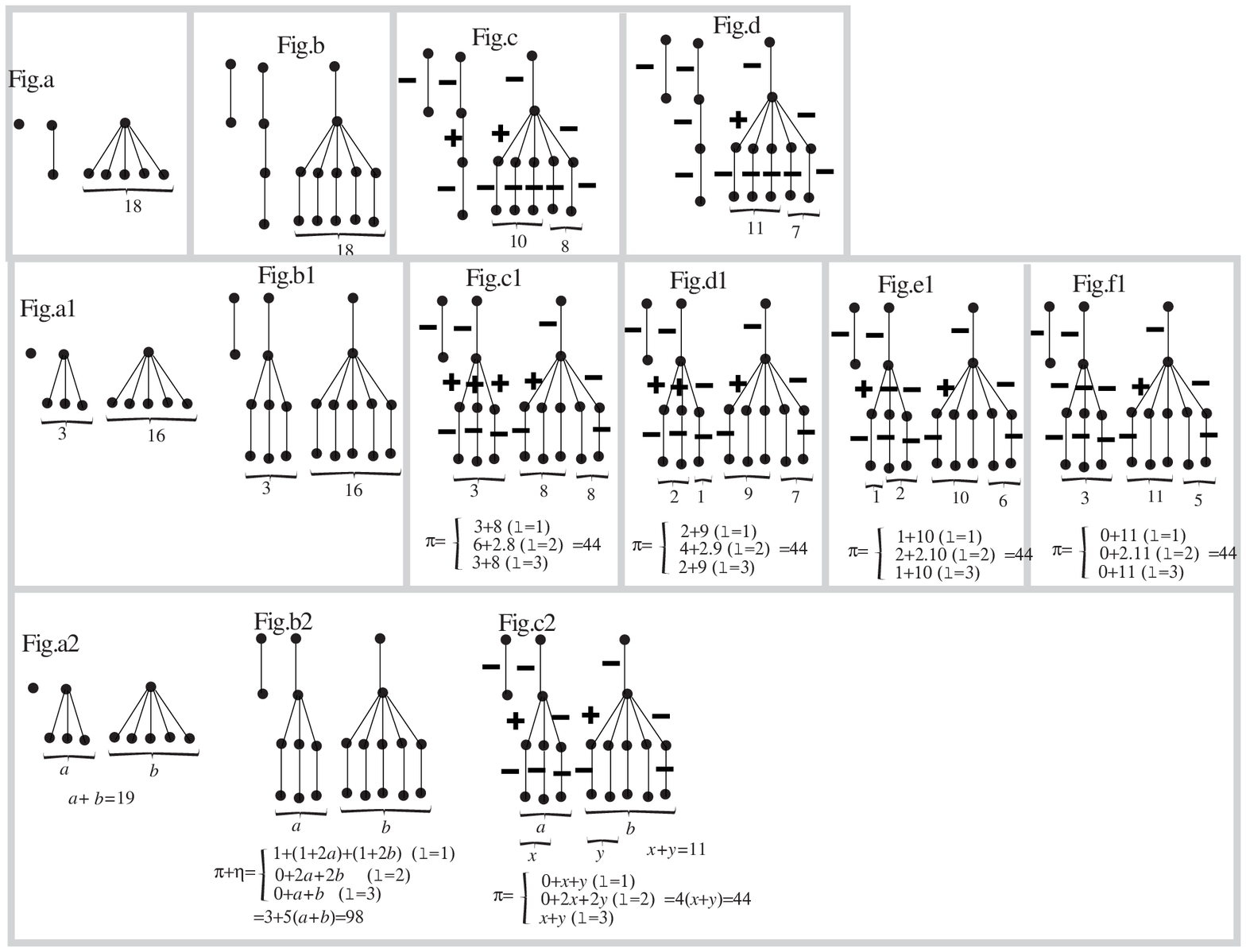,width=172mm}
\captionskipAG
  \caption{\label{Double:fig}%
  A Rohlin tree argument}
\figskip
\end{figure}

Finally, applying Rohlin's formula to the original octic we have
$2(\pi-\eta)=r-k^2=22-16=6$, whence $\pi-\eta=3$ and $\pi+\eta=19$
(looking at Fig.\,a), and therefore $2\pi=22$, whence $\pi=11$. We
find then (assuming the distribution of charges as on Fig.\,c)
$$
\pi=
\begin{cases} 0 \quad &\textrm{of } \ell=1 \cr
0  \quad &\textrm{of } \ell=2\cr
0  \quad &\textrm{of } \ell=3
\end{cases}
+
\begin{cases} 1 \quad &\textrm{of } \ell=1 \cr
2  \quad &\textrm{of } \ell=2\cr
1  \quad &\textrm{of } \ell=3
\end{cases}
+
\begin{cases} 10 \quad &\textrm{of } \ell=1 \cr
10  \quad &\textrm{of } \ell=2\cr
10  \quad &\textrm{of } \ell=3
\end{cases}=0+4+4.10=44,
$$
somewhat miraculously in accordance with the earlier calculation.
(It may be checked that for the other distribution of signs
depicted on Fig.\,d we get still the same result of 44.)

Alas, we would rather have preferred a disagreement so as to gain
an obstruction of this boson, but which perhaps  exists. Of course
one needs then to repeat the procedure for other schemes
(especially the bosons or pseudo-bosons which are doubly-nested
and with one outer oval) yet it turns out that our doubled
inference of Rohlin's formula is always respected, and as a
result, it seems that no obstruction stems from our simple device.
Especially interesting, would be the case of Orevkov's anti-boson
$1\frac{3}{1}\frac{16}{1}$, yet repeating our method still yields
twice $\pi=44$, quite regardless of the changing isotopy type.
Okay, but one has to check this for all logically possible
original complex orientation, and this involves a menagerie of ca.
four cases.

Let us do this more concretely. First Fig.\,a1 shows the nesting
tree of the Orevkov's (anti)-scheme $1\frac{3}{1}\frac{16}{1}$.
Fig.\,b1 shows the tree of the doubled satellites of a
(hypothetical) octic. Of course, for the $C_8$ we still have
$\pi=11$, as the total number $\pi+\eta$ of pairs is still 19, and
the difference prescribed by Rohlin's formula
$2(\pi-\eta)=r-k^2=22-16=6$.

Next we apply Rohlin's formula on the doubled curve of degree 16,
obtaining as before $2(\pi-\eta)=r-k^2=44-8^2=-20$, whence
$\pi-\eta=-10$. Further looking at the tree (Fig.\,b1), we count
the total number of pairs of nested ovals to find
$$
\pi+\eta=
\begin{cases} 1 \quad &\textrm{of } \ell=1 \cr
0  \quad &\textrm{of } \ell=2\cr
0  \quad &\textrm{of } \ell=3
\end{cases}
+
\begin{cases} 7 \quad &\textrm{of } \ell=1 \cr
6  \quad &\textrm{of } \ell=2\cr
3  \quad &\textrm{of } \ell=3
\end{cases}
+
\begin{cases} 1+2.16 \quad &\textrm{of } \ell=1 \cr
2.16  \quad &\textrm{of } \ell=2\cr
16  \quad &\textrm{of } \ell=3
\end{cases}=1+16+1+5.16=98.
$$
So $2\pi=88$, and $\pi=44$. Next we calculate this magnitude $\pi$
by using the charged tree assuming the distribution of signs being
that materialized by Fig.\,c1 (other charges being those of
Figs.\,d1, e1, f1). A plain calculation (mimeographed right below
the corresponding Fig.\,c1, etc.), using the signs law, shows that
$\pi$ is invariably equal to 44 regardless of the charges
distribution. So we cannot expect deriving Orevkov's obstruction
by our naive method, and more generally Figs.\,a2, b2, c2 and the
allied calculation shows that $\pi$ is always equal to 44 so that
no violation of Rohlin's formula can be obtained. In particular no
one of the 4 doubly-nested bosons can be prohibited by our
pseudo-method.

\subsection{The Hawaiian earing of Chevallier}

[30.07.13]
%%%[Old text]
Another naive idea (yet suggested by the
papers of Chevallier and Orevkov) is to consider an Hawaiian
earing alone. But then we fails strong to approach an $M$-curve as
we have only 9 micro ovals coming from the patch. Hence one must
really appeal to another singularity than $X_{21}$, namely one
where the contact between each branches is of order 4 and not just
two.

\begin{figure}[h]\Figskip
%\vskip-1.2cm\penalty0
%\centering
\hskip-2.7cm\penalty0
\epsfig{figure=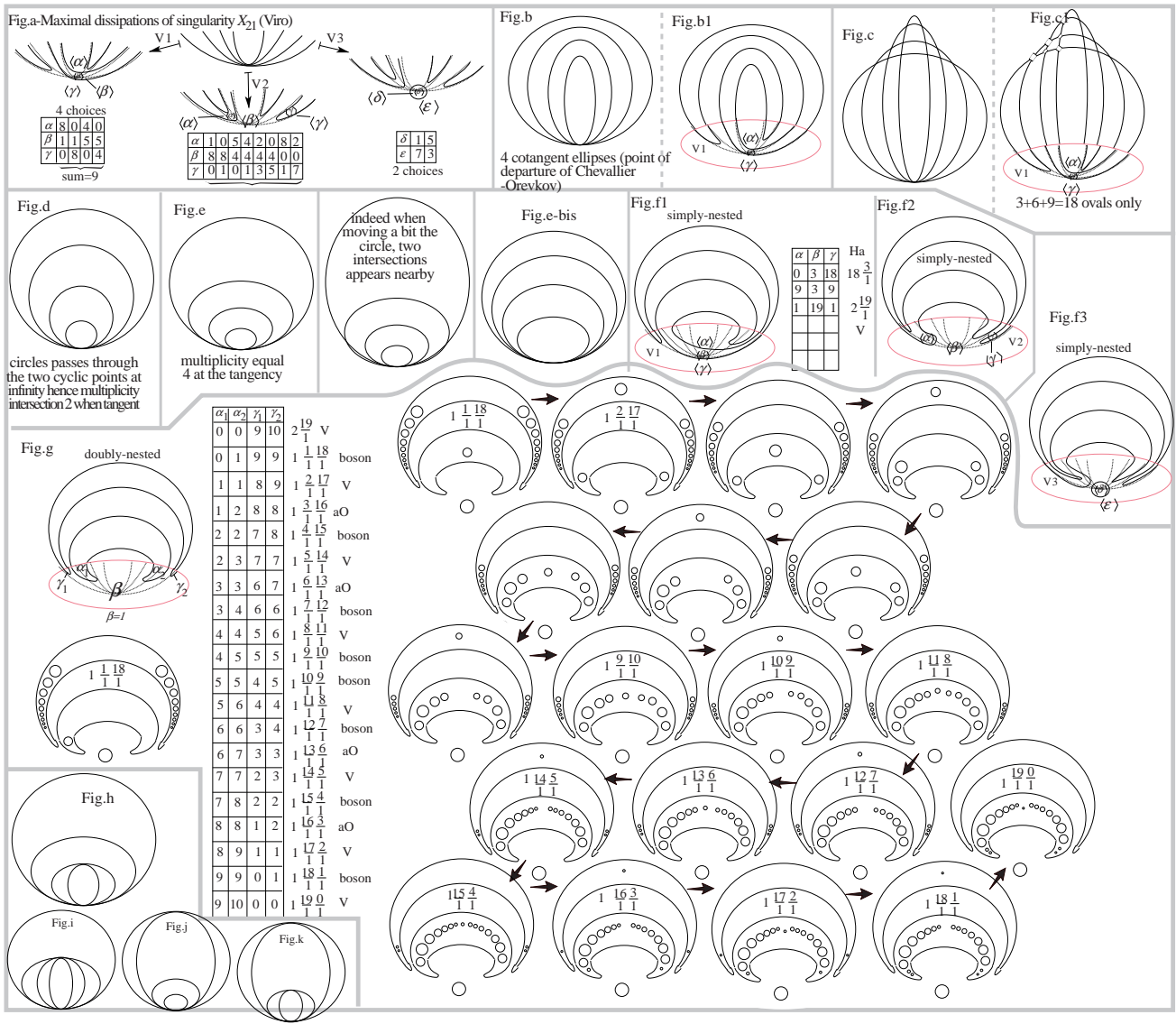,width=172mm} \captionskipAG
  \caption{\label{ViroDEGREE8_HAWAI:fig}%
  Viro's method in the Hawaiian context}
\figskip
\end{figure}

 It seems therefore that the singularity to be used has
a rich dissipation theory, as it creates the 4 schemes of
Chevallier and one due to  Orevkov. Of course it produces 5 new
schemes but probably overlap with many other older construction
and this is actually what is really worth exploring. Indeed the
point is really to decide if there is a best  curve from where all
other (or so many as possible) derive. Alas we lack the technique
to understand dissipations. A first improvisation is to posit that
the dissipating pictures are the same as for $X_{21}$, safe of
course, for the value of parameters $\al, \be, \ga$. If so is the
case then Hawaiian earing produces only simply nested schemes as
shown by the three Figs.\,f (f1,f2,f3). Actually f1 resembles the
fire-fox icon reminiscent of a snail, while f2, f3 involves a sort
of snake. At any rate it seems clear that much controtion is
involved (i.e. one oval is very long and contorted) and there is
subconscious experimental principle telling us that we cannot
thereby reach Harnack maximality under so much energy for creating
a single oval. (Alas we do not know a formal phenomenon behind
this but is the result of many accumulated evidence.) In contrast
we can then dissipate as on Fig.\,g. To land in the bosonic strip
we fix $\be=1$. To simplify we could hope to exploit symmetry as
to choose $\al_1=\al_2$ and $\ga_1=\ga_2$. However then the
content of both nonempty ovals would be even, contradicting a
basic feature of the fundamental table of periodic elements
(Fig.\,\ref{SIMPLIFIED-TABLE:fig}). Of course one could imagine an
inflation of the deformation-zone (red ellipsoid) so that both
$\al_i$ merges together to a single $\al$. This could help finding
symmetric models \`a la Fiedler.

The long table of Fig.\,g shows how all bosonic schemes can be
swept out by a single procedure. Alas we are not able to prove
existence of any such dissipation, so that we can only infer from
Orevkov's octic obstructions the nonexistence of certain patches,
notably $(\al_1,\al_2,\be_1,\be_2)$ equal to $(1,2,8,8)$,
$(3,3,6,7)$, $(6,7,3,3)$ and $(8,8,1,2)$. Of course those are only
the mostly equilibrated patches, yet their shuffles are likewise
prohibited, so for instance $(1,2,8,8)$ can be shuffled to
$(0,3,8,8)$ which is likewise prohibited.

On the positive side this heuristic method does little in proving
existence of any of the bosons, yet maybe there is free room for
progresses along this line until someone claims to have fully
exhausted the dissipation theory of this singularity (which we
call the {\it Hawaiian contact\/}$\approx$the French kiss of
Chevallier).

Further it may be noted that there are other way to land in the
bosonic strip (doubly nested with one outer oval), e.g. by
injecting ovals in the little central island while letting the
upper lune becomes empty. (More about this soon, see the next
picture.)

Another idea is that for the curves of Fig.\,g it seems worth
trying to sweep out the curve via a pencil of line based in the
central island and so perhaps to infer some information on complex
orientation via Fiedler's alternating rule. When combined with
Rohlin's formula  this can perhaps prohibit certain curve, and so
the allied dissipation. Of course this would not prohibit the
scheme itself but maybe only its realizability through our imposed
procedure. Let us explore this method. Suppose the picture to be
done more realistically as Fig.\,g1. Then sweeping from the island
it could be that Fiedler's rule implies that the chain of ovals
have alternating (complex) orientations. So among the 18 ovals
inside the upper lune (croissant) there are as many positive pair
than negative pair of ovals so that their contribution to Rohlin's
formula $2(\pi-\eta)=r-k^2=22-16=6$ cancel out. So $\pi-\eta$ is
equal to $\pm 1$ depending on the orientation of the small lune
and its inner oval. Yet in any event Rohlin's formula cannot be
respected.

Examining the dissipation with injection of ovals in the central
island gives the following table of schemes
(Fig.\,\ref{ViroDEGREE8_HAWAI2:fig}a). Again this says little
because our dissipation are fictional ones (not known by us to
exist). Again all what can be inferred is the negative result
prohibiting all patches converging to Orevkov's two schemes. Of
course we cannot exclude the possibility that some patch leading
to a boson do exist in which case we draw a new existence result.
Alas we do not master the technology granting existence of a
patch. Again there should be a geometric interpretation as a
global (projective) object (probably in a toric manifold), but
alas we do not understand properly the quintessence of Viro's
method.

Finally there is still another possible mode of generation of the
bosonic strip via the patch of
Fig.\,\ref{ViroDEGREE8_HAWAI2:fig}b, where now both croissants are
nested. Again no existence result can be inferred unless one has a
deeper understanding, which apparently nobody has despite the
related work by Chevallier and Orevkov (2002
\cite{Orevkov_2002/XX-New-M-curves-degree-8}) yet landing in
another pyramid, i.e. the sub-nested realm.

\begin{figure}[h]\Figskip
%\vskip-1.2cm\penalty0
%\centering
\hskip-2.7cm\penalty0
\epsfig{figure=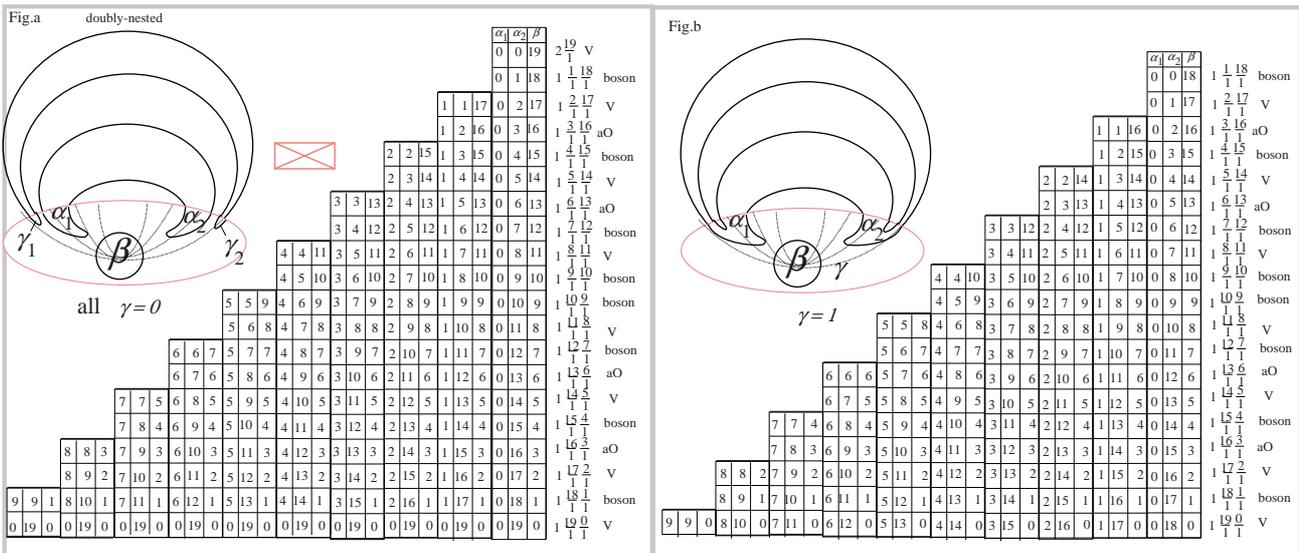,width=172mm} \captionskipAG
  \caption{\label{ViroDEGREE8_HAWAI2:fig}%
  Viro's method in the Hawaiian context}
\figskip
\end{figure}

As before one could expect that via a control of complex
orientations one could prohibit such and such geometric
realization of the abstract scheme via those more explicit model.

Alas it is clear that we technologically limited, and cannot
advance further on the problem due basically to a lack of
understanding of Viro's method.

Nonetheless, it would be of extreme importance to describe all
patches for the Hawaiian singularity, so to know globally which
schemes are adjacent to the Hawaiian configuration of 4 ellipses.
As we said there is sub-nested schemes due to Chevallier and
Orevkov that do the job, but what is demanded is a more
exhaustive
search (even when overlapping with older constructions by Viro).
Of course the subconscious desideratum is that the Hawaiian curve
has not yet delivered all its nectar (i.e., some new bosons could
perhaps derive from it).

We could try at the occasion to understand Chevallier/Orevkov
constructions in the hope that the do adapt to the doubly-nested
bosonic region, yet this seems fairly difficult to implement as
otherwise they could have done it themselves. With our poor
understanding of Viro's method it seems that what is required is a
little miracle of a convex triangulation of an arithmetic
character (integral lattice points) yielding a sort of finitary
crystallography in the Newton polygon. So what is need is before
our eyes, but nobody sees it  merely because there is a myriad of
such crystals (about one billion?),  the patch being precisely
constructed out of such a crystal.

\subsection{Viro vs. le Chevallier du temple solaire,
or monotheism vs. the necessity of several Gurus (Luc de Jouret et
son charisme extraordinaire) as to
%%%cover
account the full morphogenetic freedom of algebraic curves}

[30.07.13, aber sp\"at in der Nacht.]
%
%RICHTIG SO (nicht in die Nacht), ich fragte Christeli
%
 Actually the philosophical
issue is about a sort of absence of monotheism and religion. More
precisely, we have Viro's mandarine versus Chevallier's Hawaiian
earring. What is more powerful? As we know it seems that
Chevallier's earing was able to create 4+1 new schemes (when
combined with Orevkov's variant thereof). Albeit modest those 5
schemes were previously inaccessible through the combined efforts
of the eminent predecessors:

$\bullet$ Harnack 1876, Hilbert 1891, (and their Miss Ragsdale
[1906] who perhaps made more explicit the output of the formers
constructions),

$\bullet$ Wiman 1923 \cite{Wiman_1923} (who discovered a new
scheme in degree 8 through a methodology that apparently both
Hilbert and Ragsdale imagined as sterile and impossible),

$\bullet$ Gudkov ca. 1971 \cite{Gudkov_1971-const-new-ser-M-curv}
(who after his successes in degree $m=6$ made sporadic
contribution in degree $m=8$),

$\bullet$ Korchagin 78 (who using Brusotti, managed to construct
one more scheme),

$\bullet$ Viro 80 (who suddenly learned us that $M$-curves are
%almost
not so rare diamonds, but rather proliferating as fast as the {\it
vermine du Pripet\/}),
%or
%Chinese workers/citizens),

$\bullet$ Shustin (who completed some of the efforts by creating a
certain medusa (Fig.\,\ref{ViroDEGREE8_SHUSTIN:fig}) of an
interesting type basing himself on the dissipation theory of the
candelabrum with 3 branches transverse to the trunk), plus,
subsequently,  a very clever construction using threefold symmetry
to get the scheme $4\frac{5}{1}\frac{5}{1}\frac{5}{1}$. (Not yet
assimilated by us at the time of writing.)

$\bullet$ Korchagin (who created a myriad of 19 new schemes by
merely exploiting a complex dissipation theory and the Newton
polygon).

So Chevallier's Hawaiian earring inspects another region of the
spectrum, but we cannot exclude the crude option that it can
phagocyte completely Viro's method via the mandarine
(quadri-ellipse). Our belief (from pure intuition as we do not
know the precise dissipation of Hawaii earrings) is that both
Viro's mandarine and Chevallier earrings are complementary object,
perhaps even not sufficient to cover all octic distributions of
ovals. Crudely put, it seems that in degree $m\ge 8$ (and higher a
fortiori) there is no monotheism, but rather a pluralist world
requesting several gurus (privileged curves) enabling one to
access all rooms past the discriminant (or at least their isotopic
incarnations).

{\it Added} [04.10.13].---Besides, it seems that one can imagine
configuration of 4 ellipses hybridizing Viro's contact of order 2
with Chevallier's contact of order 4, see e.g.
Fig.\,\ref{ViroDEGREE8_HAWAI:fig}h, i, j, k.

\begin{ques}
Is there any chance that any one of those (hybrid) quadri-ellipses
produce new bosons, as the Chevallier earing managed to produce
new $M$-schemes not formerly accessed by Viro's method(s).
\end{ques}

Despite this philosophical wisdom, it seems of interest to prove
that no monotheism is possible and to measure the power of the
varied gurus. By this we just mean, given a special singular
octic, count the number of ($M$-)schemes, that can be gained by
(infinitesimal) deformation (of the guru curve).

The list of gurus is as follows:

$\bullet$ Viro's quadri-ellipses (coaxial) produces ca. 40
schemes. In principle Viro claims that the dissipation theory of
the allied singularity $X_{21}$ is fully understood (so that no
new boson could emerge through that procedure). One may wonder
where the proof of those assertions are supplied (source to
compare Viro 89/90 \cite{Viro_1989/90-Construction}).

$\bullet$ Chevallier's (Hawaiian) earring produces at least
$4+1=5$ schemes (according to literature, the last one being due
to Orevkov) but probably creates much more (possibly even schemes
not yet known to exist).

$\bullet$ perhaps the sequel of our text, or many lovely curves by
Viro, Shustin, etc, depicted in the sequel of this text shows that
there  is plenty of other gurus-curve allowing potentially to
construct new scheme not yet known to exist.

So  questions worth elucidating  are (focusing on $m=8$ just for
the sake of concreteness):

(1) What is the best  guru,
%the
humanity can dream about? Jesus
Christ? Mohammed? Herbert Gr\"otzsch? Markus Schneider?

Natural candidates: Viro's quadri-ellipse, but check if
Chevallier's earing is not stronger. Alas no complete data is
available to us.

(2) Assuming that there is no monotheism, what is the minimal
number of Gurus required to cover all schemes in degree 8? In
general call this $G(m)$ the Guru constant in degree $m$. For
instance $G(4)=2$, as all schemes are deducible from two
transverse ellipses safe the empty scheme. Likewise $G(6)=2$ as
all non-empty schemes are deducible from Viro's mandarine. This
number has an intrinsic significance and it would be interesting
to know if the despite the bewildering exponential rate of growth
of all schemes as a function of $m$, if the number of Guru's
required can be kept into reasonable growth like a linear growth.
Roughly put, mankind is stupid or ``moutonesque'' enough to
tolerate very few Guru's. Of course in an ideal world everyone
should be his one guru (or at least take the responsibility
thereof).

Of course we do not have any moral message to transmit, but just
used this image to get red of the geometric insigificance of our
search. Assuming that all 6 bosonic octic $M$-schemes are
prohibited it seems that the Guru number $G(8)$ is circa one for
Viro (plus two for its tricky variant), two  for Shustin, one for
Chevallier and Orevkov. So all $M$-schemes known up to date
derives from only $3+2+1=6$ Gurus. Of course assuming taht all
remaining 6 bosons exist we could let degenerate them by brute
force on a wall of the discriminant and thereby get 6 new (adhoc)
gurus of poorly charismatic nature.

Hence it seems that actual knowledge (discussed in more detail in
the sequel) implies the following:

\begin{lemma}
The Guru constant in degree $m=8$ is at most $12$, $G(8)\le 12$.
Of course it can be much lower, say perhaps as low as $G(8)=1$ in
case the Chevallier (du temple solaire de Toulouse) is able to
crack all $M$-schemes (the empty one being discarded incarnating
the opposite spectrum of maximality).
\end{lemma}

Of course there is (at least) two ways to define the Guru constant
depending on whether we want to access all schemes or confine to
the maximal ones. We denote $G(m)$ that for $M$-schemes, and by
$g(m)$ that encompassing all schemes. Evidently $G(M)\le g(m)$.
For instance $G(6)=1$ via Viro's mandarine, while $g(6)=2$ via
Viro's mandarine plus Fermat's solitary curve $x^6+y^6=0$.

Another way to define the Guru constant involves the rigid version
of really spotting out all chambers past the  discriminant and not
crude isotopy type. This we call $\Gamma(m)$ the rigidified Guru
constant. For $m=6$ this is probably computable via the
Rohlin-Nikulin-Kharlamov theory employing K3 surfaces. Probably
$\Gamma(6)=2$ once Viro's method is known.

Of course assuming that there is a linear complexity upon the
number of Gurus then Hilbert's  16th could still be
%so-to-speak
(sozusagen) tractable algorithmically, but probably energy
resources to complete this work will be wasted by private Bankers
and oriental oil merchant prior than  humanity becomes conscious
of its real mission. All the best and good luck for the sequel.

[31.07.13] Also, it would be nice to do a careful map of all the
gurus and their zones of  influence.

As we said we have as gurus the following curves:

(1) Viro's mandarine (Fig.\,\ref{ViroDEGREE8:fig}). Its zone of
influence is marked by green rectangles enclosing the letter ``V''
on Fig.\,\ref{SIMPLIFIED-TABLE_gurus:fig}. This is very compact on
the 1st pyramid, yet with severe lacunas on the bosonic strip. The
mandarine zone attacks only 3 schemes in the trinested region (2nd
pyramid).

(2) Viro's beaver (castor) which realizes 7 schemes
(\ref{ViroDEGREE8_TRICKY:fig}).

(3) Viro's horse which realizes

(4) Shustin's medusa (cf. Fig.\ref{ViroDEGREE8_SHUSTIN:fig})

(5) Shustin's flower with threefold symmetry (alas no depiction
available for the moment)

(6) Korchagin but no understanding

(7) Chevallier's earrings

It seems further that a natural Guru is the following curve

(8)  (Gabard's) margarita: a singular octic with an ordinary
septuple point (alias the margarita flower), cf.
Fig.\ref{GabardDEGREE8_septuple:fig}. The dissipation theory of
this guru seems to reduce to that of affine $M$-septics with
pseudo-line oscillating simply across the line at infinity. Hence
the radiation influence of the margarita includes all $M$-schemes
realized via affine $M$-septics. Albeit this zone is a priori
equal to the blue region of
Fig.\,\ref{SIMPLIFIED-TABLE_gurus:fig}, but alas the dissipation
of the septuple point or equivalently the theory of affine septics
is not as yet sufficiently developed to make an exact description
of the zone.

\begin{figure}[h]\Figskip
%\vskip-1.2cm\penalty0
%\centering
\hskip-2.7cm\penalty0
\epsfig{figure=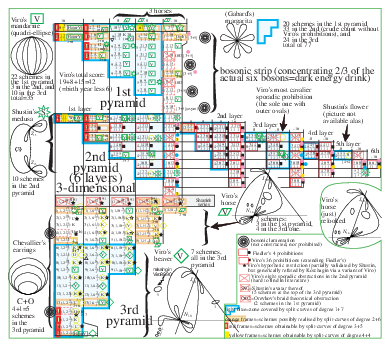,width=172mm}
\captionskipAG
  \caption{\label{SIMPLIFIED-TABLE_gurus:fig}%
  The Gurus and their zone of influence/radiation.}
\figskip
\end{figure}

The next figure (Fig.\,\ref{gurus:fig}) is an attempt to discover
Shustin's flower yet we lost ourselves in random nonsense. In fact
after a long search we came to the idea that if an octic has 3
candelabrum points then the conic through them and tangent to both
both of them (think with 2 infinitely close points) will have an
intersection of $3\cdot 4+6=18> 16=2\cdot 8$, and B\'ezout is
overwhelmed. So it seems that:

\begin{lemma}
There is no singular octic with 3 candelabrums, except perhaps if
it splits off a conic.
\end{lemma}

\begin{figure}[h]\Figskip
%\vskip-1.2cm\penalty0
%\centering
\hskip-2.7cm\penalty0 \epsfig{figure=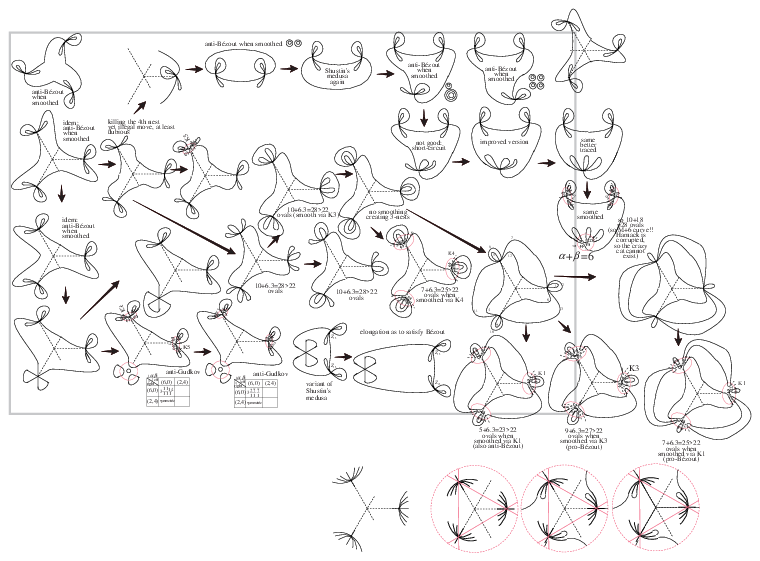,width=172mm}
\captionskipAG
  \caption{\label{gurus:fig}%
  A la recherche of Shustin's curve}
\figskip
\end{figure}

[01.08.13] Further it is natural to wonder if there is any natural
species (animal) between the horse and beaver that would permit
the creation of new (or old) schemes. As we said in this game one
must not strive toward extreme originality but rather try to get a
global understanding through much overlap. Also one problem is to
find the hottest curve permitting to create the greatest number of
smooth curve.

For instance one can imagine a version of Viro's horse where all 4
inner ovals are transferred outside. Then the schemes are just
changed by a fluctuation of 4 ovals so that probably just a left
translation on the table (Fig.\,\ref{SIMPLIFIED-TABLE_gurus:fig})
is effected. Similarly one can imagine a version of the beaver
with 4 outer ovals transplanted inside. Let us work out the scheme
realized by those genetically modified birds. The little surprise
is that we get so new schemes claimed by Viro (like
$3(1,14\frac{3}{1})$ or $3(1,10 \frac{7}{1})$) but which we were
never able to construct as yet.

\begin{figure}[h]\Figskip
%\vskip-1.2cm\penalty0
%\centering
\hskip-2.7cm\penalty0
\epsfig{figure=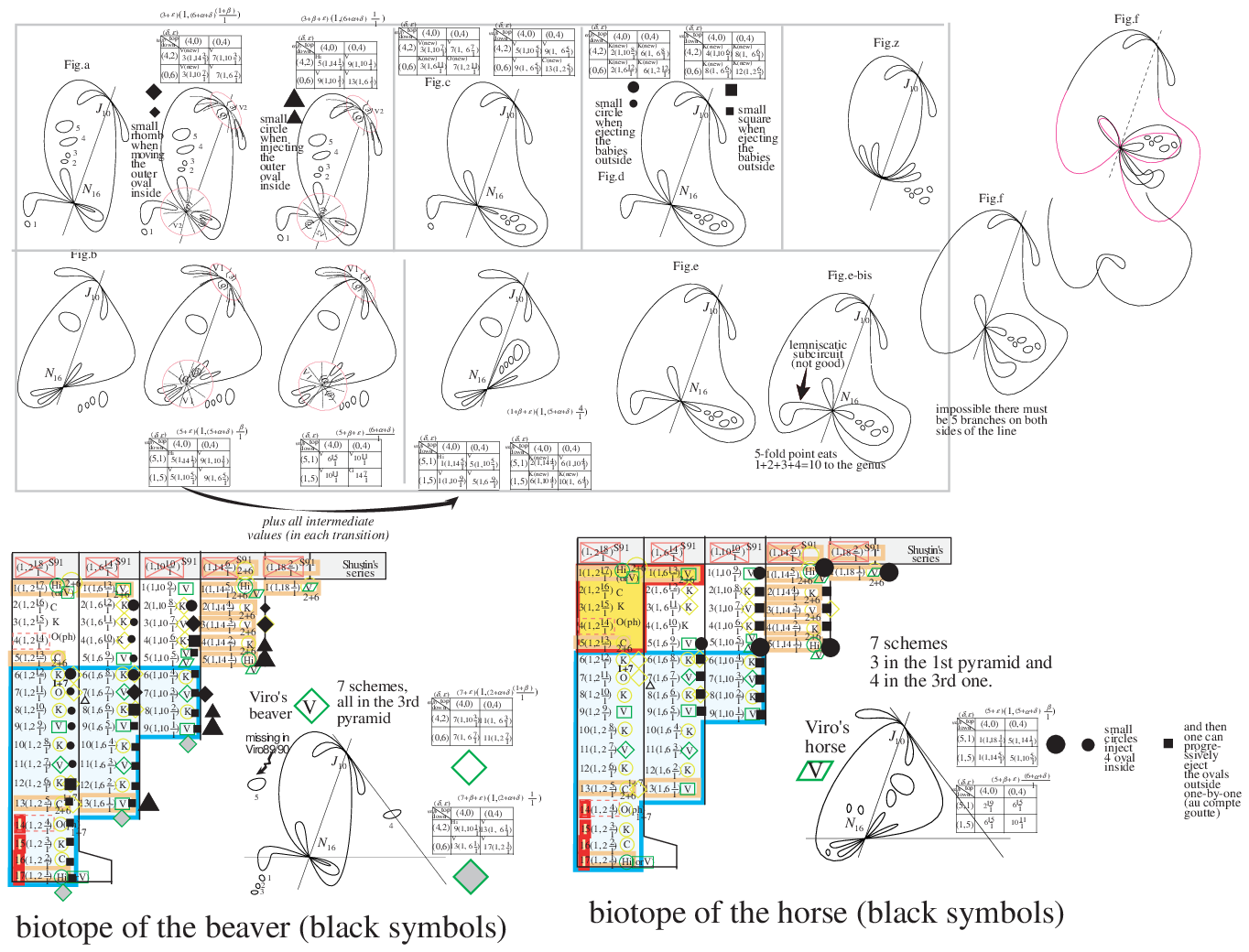,width=172mm} \captionskipAG
  \caption{\label{Beaver-modified:fig}%
  Beaver and horse genetically modified}
\figskip
\end{figure}

Of course one can then do even more radical genetic modification
of beaver and horse, like transferring the beaver's eggs in its
foot (Fig.\,c). The resulting schemes are obtained just by
shifting by 4 the 2 earlier tables. Then we get one scheme due to
Korchagin but that we were never able to understand and which is
thus relatively new, and also the scheme $7(1,2\frac{11}{1})$ due
to Orevkov. Likewise the 2nd table yields the interesting scheme
of Chevallier $13(1,2\frac{5}{1})$.

\begin{Scholium}
There is perhaps a fairly elementary construction of Orevkov's
(unique) octic scheme by a variant of Viro's method using a
genetically modified beaver (Fig.\,\ref{Beaver-modified:fig}c).
\end{Scholium}

At this stage we were striked by the idea of using a quintuple
flat point (i.e. with coincident tangent). This should have more
power, or otherwise in the beaver we could trade the upper (flat)
triple point for an ordinary one which admits a more elementary
dissipation theory.

Fig.\,d transplants the outer egg in the pregnancy bag of the
beaver, and then we get the same collection of schemes modulo a
shift. Actually Fig.\,d is the extremal pregnancy level  for the
beaver and it creates then the schemes $2(1,10\frac{8}{1})$ and
$2(1,6\frac{12}{1})$ plus all five schemes below them (on the
table) amounting to eject the babies outside of the beaver's
belly. From $6(1,2\frac{12}{1})$ we can also run 5 schemes below
on the table visiting in particular Orevkov's scheme. Once this is
understood we have a clear view of the possibility (geographic
aptitude) of the beaver. Roughly it is most interesting when
pregnant, yet it still misses Chevallier's scheme
$5(1,2\frac{13}{1})$ and what is above.

So it seems that we need a completely new animal to explore this
zone, or perhaps the beaver's queue may be invaginated as on
Fig.\,e. Alas on improving the picture it seems that there is a
short-circuit then this represents usually a waste of energy
leading us away from the realm of $M$-curve.

More formally, we know that a fivefold point eats 10 units to the
genus, so as the beaver curve has 6 circuits (hence genus 5) we
infer somewhat indirectly that the triple point eats 6 units to
the genus ($21-10-x=5$, whence $x=6$). At any rate the fact that
Fig.\,e does not represent a perfect circuit diminish our aptitude
to create ovals, and we are lost. It seems that we have completely
exploited the capacity of the beaver.

As a last remark concerning the beaver, when Fig.\c produces
Chevallier's scheme $13(1,2\frac{5}{1})$, then we may additionally
eject the deepest ovals outside of the belly and so we get the 4
schemes below it on the main-table, namely
$14(1,2\frac{4}{1})$~(boson), $15(1,2\frac{3}{1})$~(K=Korchagin),
$16(1,2\frac{2}{1})$~(C=Chevallier),
$17(1,2\frac{1}{1})$~(Hi=Hilbert).

Actually one can represent the biotope of the beaver by black
symbols on the table (see Fig.\ref{Beaver-modified:fig}) and one
sees that actually permitting on Viro's original beaver quantum
jump inside one can raise always one five positions so that from
the fundamental parallelogram (green rhombs with signature
``V=Viro'') plus the 3 aligned same rhombs but shaded one can
sweep out the full biotope of the beaver safe the right row. Then
in a similar way it is as simple matter to delimit the biotope
(habitat) of the horses and their genetic modification. Indeed one
marks first the 4 fundamental schemes coming from Viro's horse
(parallelogram of thickest circles), and then one can inject  4
median ovals in the subnest. Diagrammatically this implies a shift
to the left (smaller circles), and once this is done one can
progressively eject the deep ovals outside without changing the
Euler-Ragsdale characteristic $\chi$ so that ejection one-by-one
is permissible (without conflicting with Gudkov periodicity).
Diagrammatically this has the effect of moving down along 4
positions, and so we get the black region.

It seems therefore that both the beaver and horse are never able
to visit the polar (or desertic) region yellow-colored ont our
map. In particular it seems to us that the corresponding boson
($4(1,2\frac{14}{1})$)---which is in some sense protectionist with
many ovals protected at depth 2---looks somewhat harder to detect
than its open-minded/liberal companion $14(1,2\frac{4}{1})$. Of
course it can be as well as case that we did not as yet explored
all type of animals by relying merely on Viro's beaver and horse.

So some artistic imagination power is required to draw new animals
and here comes a sort of theory des dessins d'enfants into the
game at least to make the task psychologically supportable.
However it seems that there is some obstruction to modify the
beaver in the requested geographical locus.

Hence, we may of course imagine more artistic tracing of curves
representing principally new animals.  Fig.\,z looks to be a
serious candidate using a quintuple point with 5 tangent branches.
Yet alas this requires a new dissipation theory that we are not
familiar with. Of course one can improvise, yet we feel too tired
today to explore this in any systematic fashion.

[02.08.13] Actually one can imagine such curves as coming from a
Jordan domain (splashed disc) with tentacles growing (cytoplasmic
expansion) and self-penetrating an invagination of the domain. So
we get for instance the curve depicted once smashing the vertical
point so to create higher singularities. For instance Fig.\,e
looks promising because it has only one external teats hanging
outside, and therefore the smoothed curve should land in the
bosonic strip (where four $M$-schemes are in suspense/doubt).
Fig.\,f is just a relooking with less B\'ezier points so that the
curve looks better (a Kunstform der Natur). Actually the smoothing
Fig.\,f1 leads potentially to Korchagin's (desertic) scheme
$3(1,2\frac{15}{1})$, and when modified even the (protectionist)
boson $4(1,2\frac{14}{1})$.

\begin{figure}[h]\Figskip
%\vskip-1.2cm\penalty0
%\centering
\hskip-2.7cm\penalty0 \epsfig{figure=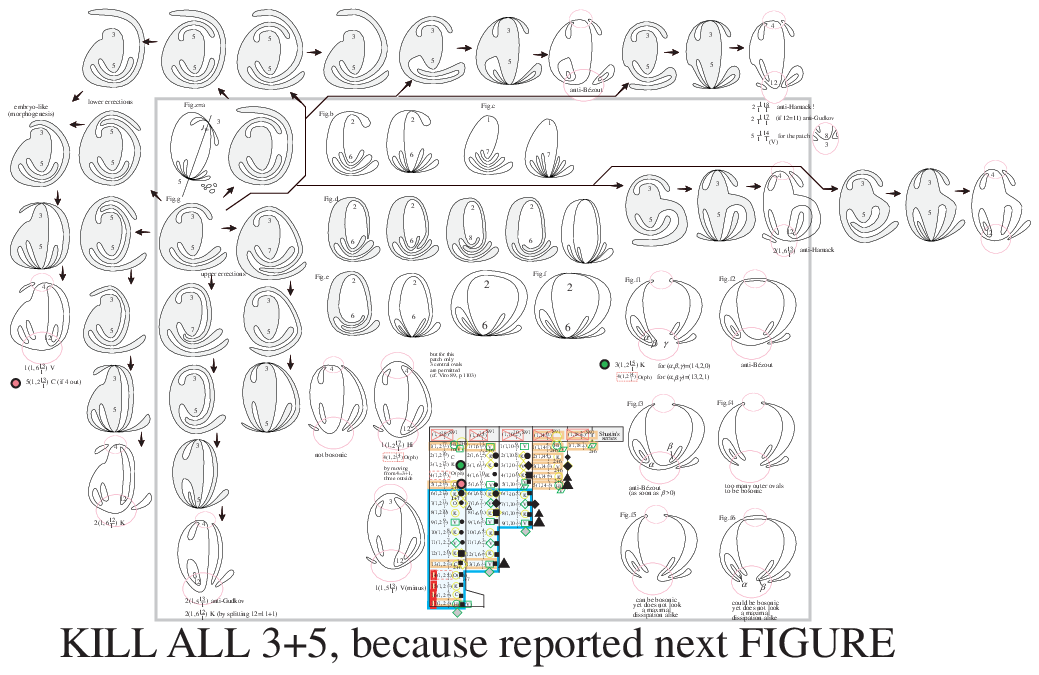,width=172mm}
\captionskipAG
  \caption{\label{gurus_2:fig}%
  More gurus}
\figskip
\end{figure}

Then Fig.\,g shows a myriad of four fashions for this ground form
to evolve through varied erections. The ground idea is to have a
cell-like shape with cytoplasmic expansions interpenetrating
themselves, and then smashing the cell we get two singularities of
multiplicity 3 and 5 respectively. Eventually after a long search
we realized that the snake of Fig.\,h or that of Fig.\,i
incarnating the simplest forma for having one protuberance leads
us potentially to the bosonic strip. Of course Fig.\,h1 may be an
aggression for B\'ezout (look at the dashed line), hover it seems
plain that a contortion of the queue could repair this misfortune.
Of course in both cases (Figs.\,h or i) we just depicted the
extremal value of the parameter, but transplanting one-by-one
ovals through (hypothetical) patches yields in all cases all the
bosons even in some palindromic repeated fashion. From Fig.\,h2 we
get first for $(\al,\be)=(12,1)$ the boson
$1\frac{1}{1}\frac{18}{1}$, and then run through all successors
$1\frac{2}{1}\frac{17}{1}$, $1\frac{3}{1}\frac{16}{1}$ up to
$1\frac{13}{1}\frac{6}{1}$. Of course if Orevkov's obstruction are
correct, then we get a corresponding prohibition of the patch for
$F5$ (flat quintuple point, in Gabard's notation not acquainted
with Arnold's numerology). So again little can be gained but
potentially all those 4 bosons are constructible if one had a
understanding of the smoothing of $F5$, and of course if the curve
of Fig.\,h1 (or better its contortion Fig.\,h3) exists
algebraically. As to Fig.\,i it is noteworthy that despite the
dissipation of F3 used is not maximal from the viewpoint of the
only 3 quantum=small ovals created by the dissipation described in
Viro 89/90 (p.\,1103), globally the smoothed curve has 6 macro
ovals (so one more that the former Fig.\,h2) and therefore we
still obtain an $M$-curve.

Of course our method has to be explored much more systematically,
yet it gives already a sort of algorithm to generate the singular
circuit possibly realized by singular octics, via a pleasant
embryology of the basic cell (disc).

\begin{figure}[h]\Figskip
%\vskip-1.2cm\penalty0
%\centering
\hskip-2.7cm\penalty0 \epsfig{figure=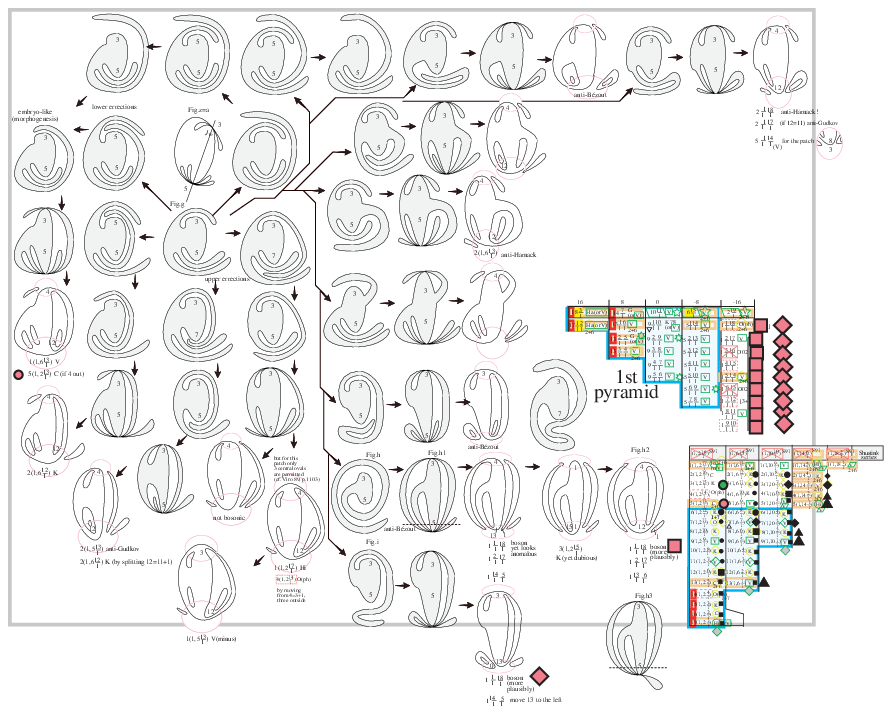,width=172mm}
\captionskipAG
  \caption{\label{gurus_3:fig}%
  More gurus}
\figskip
\end{figure}

Then of course beavers and horses corresponds to the combination
O5+F3, where O5 in an ordinary quintuple point, while F3 is a flat
triple point. One can imagine then to explore all possible
combination like rather F5+O3, where we have a quintuple flat
point and an ordinary triple points.

[03.08.13] Our algorithm of the cytoplasmic expansions of a cell
(embryology) can as well applied as to produce curves with a flat
pair of quadruple points ($8=4+4$) or one double plus a sextuple
point or even one simple combined with a septuple point. In the
case of  even splitting the sole difference is that the ground
cell is outside of the vertical line and it is merely its
cytoplasmic expansions (teats) that intercepts it. So for $2+6$ we
merely have one teat above and 3 teats below in the form of a
comb.

Tracing all those cells and the allied  curves requires as usual
some artistic skill and patience plus the appropriate dissipation
theory. Evidently the case $4+4$ deserves special interest as the
dissipation theory is then completely known according to efforts
of Viro and perhaps some intervention of Korchagin. At least we
may hope that their classification is complete (as stated e.g. in
Viro 89/90, p.\,XXXX). This being said let us explore the
embryology in degree $4+4$. (Of course since at least yesterday it
is clear also that this artwork is closely connected to the
cartoons of an eminent artist known as Ibl al Rabin, alias Mathieu
de Baillif, who is famous for a minimalist bande dessin\'ee,
characterized by purely 2D-black shapes representing humans in all
their positions and social duties). Further our embryology (and
more generally Hilbert's 16th) has some close connection evidently
with H.\,A. Schwarz's \"Olfleck. Yesterday we experimented with
oil in water (or wine) how one can create subnest by injecting a
wine drop in oil, and then again put oil droplet in that  wine to
form arbitrarily complex Einschachtelungen \`a la
Hilbert-Ragsdale. Of course when everything is exited dynamically
one gets beautiful patterns.

Doing the picture of Fig.\,\ref{gurus_4+4:fig} we realized that
most (all?) of our embryos violates severely Harnack's bound.
Indeed Viro's patches of $X_{21}$ permits to create always 9 extra
micro ovals whatever the type of dissipation used $V1,V2, V3$, so
that the total number of ovals is equal to the apparent number
plus $2\cdot 9=18$. Eventually, we realized that for a deformed
embryology Viro's original method can also be interpreted in terms
of two cells degenerating to the topology of an annulus. We can
imagine more (non-connected) embryos like the copyright symbol on
the left of Viro's species. This, when smashed and smoothed, leads
to $5+18=23>22$ ovals, violating again Harnack. So:

\begin{lemma} All the qualitative singular octics depicted on
Fig.\,\ref{gurus_4+4:fig} below, do not exist algebraically simply
because when dissipated \`a la Viro they produce curves violating
Harnack's bound.
\end{lemma}

Probably there is a theoretical justification a priori without
invoking Viro's theory (dissipation of $X_{21}$), like arguing
that if an octic has two $F4$ singularities ($F4=X_{21}$ in
Arnold's notation) then the conic through both points with the
prescribed tangents has intersection number with the $C_8$ of at
least $2\cdot 4=8$ at each point, so a total of 16 intersections.
This hold true for any conics of the pencil, yet when the
curvature is made to coincide with one of the 4 branches then the
contact exceeds 16, and so by B\'ezout the octic is forced to
split off an ellipse (that ellipse). Hence:

\begin{lemma}
Any singular octic with a distribution of two flat quadruple
points is necessarily the union of 4 ellipses. (Disappointing
corollary: there is no rich embryology as that just explored
susceptible to produce new octic schemes, especially the six
unsettled bosons not yet known to exist).
\end{lemma}

\begin{figure}[h]\Figskip
%\vskip-1.2cm\penalty0
%\centering
\hskip-2.7cm\penalty0 \epsfig{figure=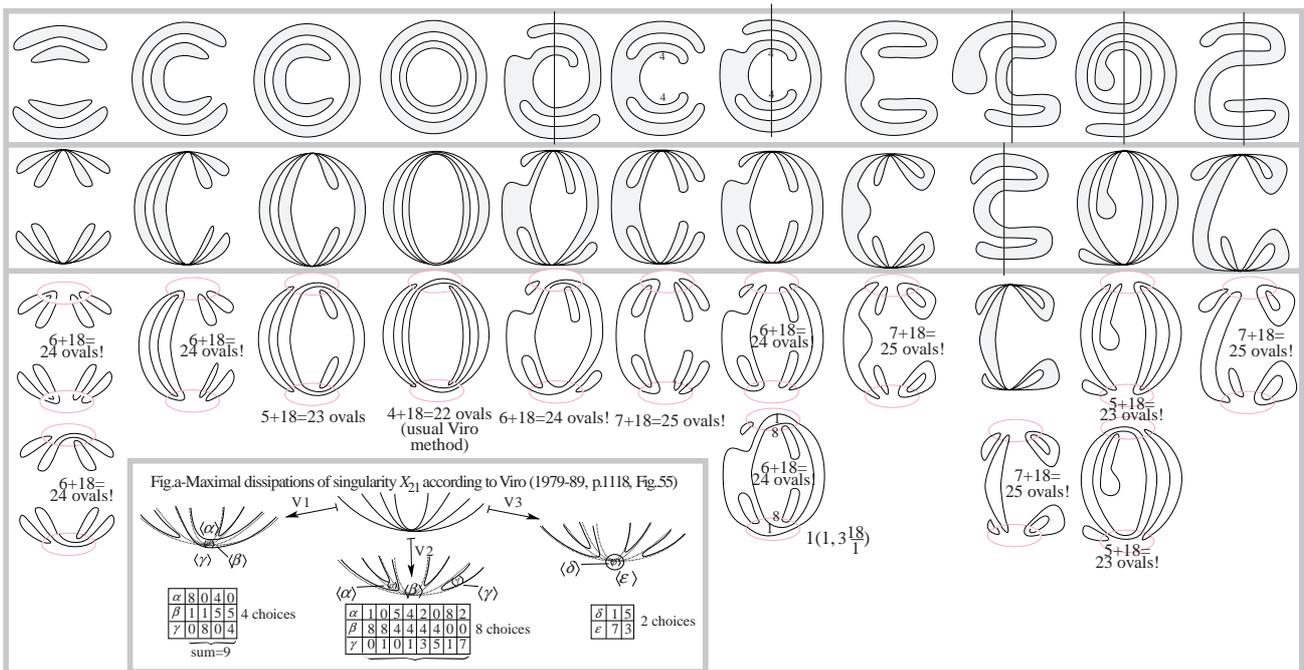,width=172mm}
\captionskipAG
  \caption{\label{gurus_4+4:fig}%
  Embryology of gurus in degree $4+4$}
\figskip
\end{figure}

So what next? Of course one could hope that embryology is still
useful in other context like $8=3+5$ or $8=2+6$, or even $8=1+7$.
Alas we run then again more obscure dissipation theory (not
described in Viro's works to the best of our knowledge). One idea
could be to use more complex pattern (embryos) involving other
distribution of foldings than F4+F4, e.g. F4 plus several F3, or
just several F3.

After some few trials we arrived at a somewhat ad~hoc pseudo-octic
(depicted as Fig.\,\ref{gurus_4+3:fig}a) able to create the 2
bosons $1\frac{1}{1}\frac{18}{1}$ and $1\frac{7}{1}\frac{12}{1}$
via the usual dissipation theory for F3=$J_{10}$ and F4=$X_{21}$
in Arnold's nomenclature. So:

\begin{Scholium}
Possibly Viro's standard dissipation theory has not yet been fully
exploited if one dispose of a global curve tracer able to
substantiate our embryology in the algebraic category.
\end{Scholium}

Of course our method is no mystery, as we started from a $3+4$
configuration of Fig.\,a and then expanded to Fig.\,b by adding a
2nd triple point (of the flat type F3). Then smoothing the curve
(downwards on the Fig.) we got a scheme with 25 many ovals. So we
just have to kill 3 of them as this and this is achieved by our
Fig.\,c. Of course another choice is to ``bore the canal'' like on
Fig.\,d. Further our picture is intercepted 12 times by the red
line (so B\'ezout is severely foiled). So the hard game is to get
more respectable pictures of such pseudo-octic. An improvement
along this way  is given just by rotating the head of the animal
to get Fig.\,e. Alas, there is still a line with 10 interceptions
of our hypothetical octic $C_8$. Of course one could hope that
this is merely a defect of our model and that there is a
diffeotopic model  mimicking better the behavior of an algebraic
octic.

[04.08.13] Next we had the idea of using a more symmetric embryo
(Fig.\,f sembling a smiling face) and after slight morphogenetic
adjustment (lowering the number of ovals within Harnack's range
while also arranging one outer oval) we find Fig.\,f1 producing
both a schemes prohibited by Orevkov and one of the bosonic type
(namely $1\frac{7}{1}\frac{12}{1}$). Again our figure is not ideal
as one can trace a line intercepting the $C_8$ along 10 points.
Fig.\,f2 shows how to corrupt one Viro sporadic obstruction while
also getting Shustin's scheme $\frac{5}{1}\frac{7}{1}\frac{7}{1}$.
Fig.\,f3 yields the same as Fig.\,f1, while Fig.\,f4 yields the
same schemes as Fig.\,f2. Next we realized that appealing to
Viro's dissipation V3, all of our pseudo-octics (Figs.\,f2,f3,f4)
produce other interesting schemes.

More precisely Fig.\,f0 produces under the smoothing V2 some basic
Viro's scheme, but under V3 schemes violating Viro's
oddness(oddity?) law. As a result either the latter is wrong
(unlikely because published), or we deduce that there is no
singular octic isotopic to our Fig.\,f0 (referring of course to
its smashing depicted right below).

As to Fig.\,f1, it creates a scheme anti-Orevkov, one boson, and
two classical schemes (a simple one of Viro plus the most tricky
one of Shustin). Therefore either Orevkov's sporadic obstruction
is wrong or the latter implies that there is no singular octic
like Fig.\,f1 (smashed as depicted right below).

When it comes to Fig.\,f2 (perhaps embryologically even nicer than
the previous pictures), we see the creation one scheme violating a
Viro sporadic obstruction, but besides 3 perfectly standard
schemes of Viro and 2 of Shustin. Therefore, either Viro's
sporadic obstruction is wrong or it implies the non-existence of
our singular octic Fig.\,f2. So, as a scholium, it perhaps quite
probable that some of Viro's sporadic obstructions are wrong

Likewise Fig.\,f3 violates one of Orevkov's two prohibition while
producing a boson plus two standard schemes due respectively to
Viro and Shustin. It may be noted that Fig.\,f3 creates actually
the same schemes as Fig.\,f1, yet perhaps afford perhaps a better
aggression against Orevkov's obstruction, as the picture of the
singular octic is more symmetric. This seems remarkable enough to
deserve a statement:

\begin{lemma}
Either Orevkov's obstruction of the octic $M$-scheme
$1\frac{6}{1}\frac{13}{1}$ is wrong, or if true then there is no
singular octic isotopic to our Fig.\,f3.
\end{lemma}

As to Fig.\,f4 it produces the same schemes as Fig.\,f2, yet
through a more symmetric model, and thus represents perhaps a more
severe offense against the relevant Viro's sporadic obstruction.
Precisely, we have the:

\begin{lemma}
Either Viro's obstruction of the octic $M$-scheme
$\frac{1}{1}\frac{5}{1}\frac{13}{1}$ is wrong, or if true then
there is no singular octic isotopic to our Fig.\,f4, nor to our
Fig.\,f2. Crudely put there are two flexible singular octics
fighting against the truth of the sporadic Viro obstruction, so
that the latter is perhaps wrong. (Of course it would be of
interest to realize the other Viro anti-scheme
$\frac{1}{1}\frac{3}{1}\frac{15}{1}$ via our pseudo-method).
\end{lemma}

\begin{figure}[h]\Figskip
%\vskip-1.2cm\penalty0
%\centering
\hskip-2.7cm\penalty0 \epsfig{figure=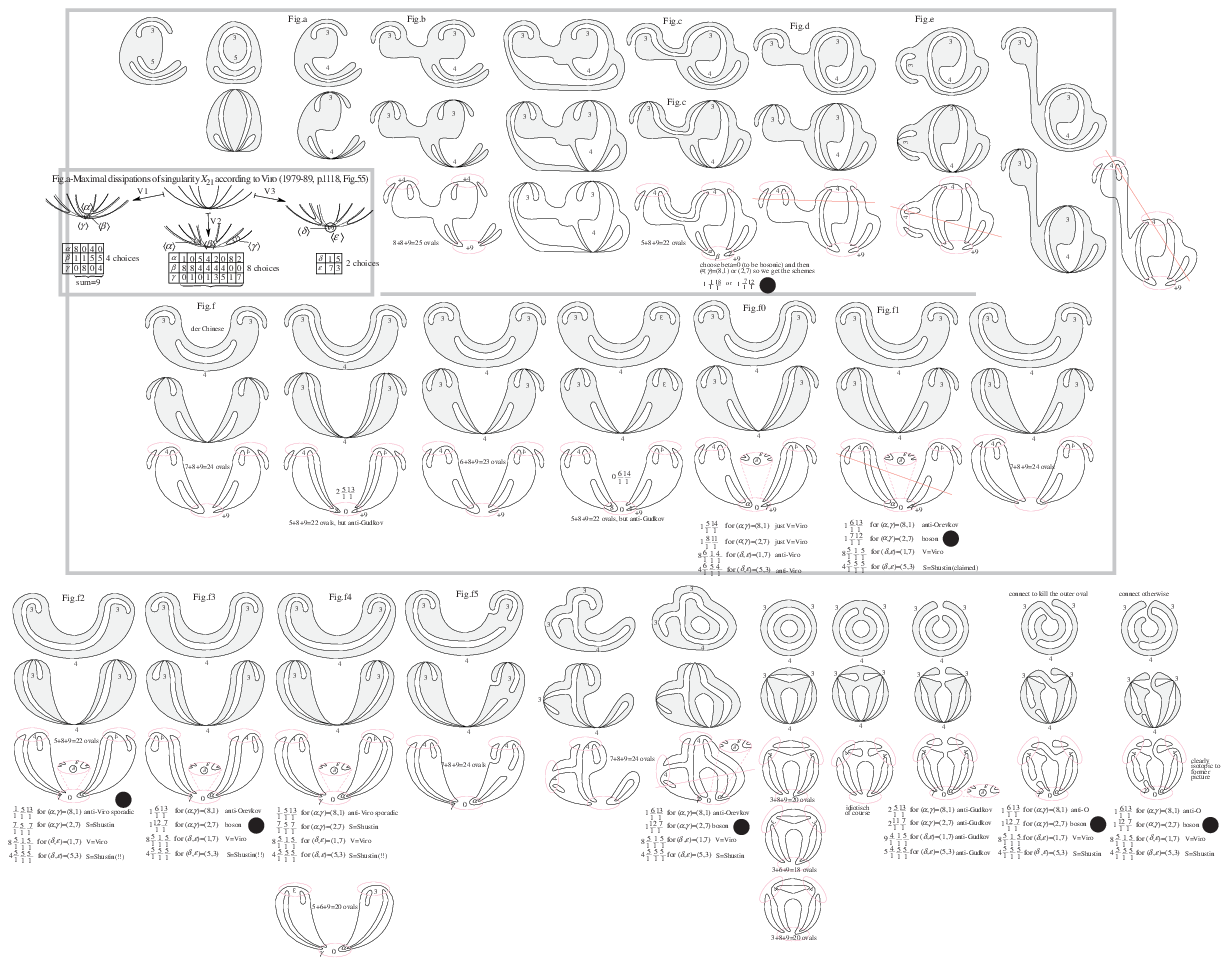,width=172mm}
\captionskipAG
  \caption{\label{gurus_4+3:fig}%
  Embryology of gurus in degree $4+4$}
\figskip
\end{figure}

\subsection{More bosonic embryology as applied
to Hilbert's 16th}

[05.08.13] Then there is more possibilities when quadruple point
occurs through the interpenetration of 2 teats. After few trials
leading to anti-Gudkov configurations we had the idea of taking an
gemelar embryo like Fig.\,\ref{gurus_4+3_BIS:fig}a. Albeit the
latter leads again to an anti-Gudkov scheme, it is easy to imagine
appropriate surgeries  correcting the number of outer ovals (which
in the doubly nested case has to be one mod 4, with special
interest in the case when this is really one in the bosonic
range), and so we get quickly Fig.\,\ref{gurus_4+3_BIS:fig}b. The
latter produces a pair of bosons, namely
$1\frac{1}{1}\frac{18}{1}$ and $1\frac{7}{1}\frac{12}{1}$, plus
two elementary schemes both accessible to Viro's simplest method
(via the quadri-ellipse).So again this gives support that both
bosons do have some chance to exist. Further this time there is no
contradiction with Orevkov's obstruction when producing the boson
$1\frac{7}{1}\frac{12}{1}$.

Is there other embryos leading to other bosons? Fig.\,c is another
option for bringing the number of outer ovals to one, yet the
resulting schemes are isotopic to those of Fig.\,b.

Fig.\d shows another surgery leading to the trinested realm, yet
outside Gudkov's range. Fig.\,e shows a simple correction to zero
outer ovals as it should be in the trinested case. Alas, the
resulting schemes are anti-Viro (imparity law), yet it is simple
to correct by surgering  at different places as to get Fig.\,f.
Its smoothing gives the scheme
$\frac{1}{1}\frac{9}{1}\frac{9}{1}$, which very shamefully is
missing from our table say just at the combinatorial level. So we
added it yet on Fig.\,\ref{SIMPLIFIED-TABLE_gurus:fig} yet beware
that other copies of that table may not be up-to-date.

More intrinsically we have the:

\begin{Scholium} (WRONG, cf. right-below)
There is a potentially new boson that we just overlooked at the
combinatorial level. In particular, it seems to us that nobody
(i.e. neither Viro nor Shustin or whoever else) never constructed
nor prohibited this scheme which is therefore a potentially new
boson in the trinested realm. In particular Shustin's assertion
that Hilbert's 16th (isotopic classification) of $M$-curves is
complete in that case is (slightly) erroneous. Compare e.g.
Shustin 1990/ \cite{Shustin_1990/91-New-restrictions} or some
other work by this author.
\end{Scholium}

Sorry very  much, the problem was that this symbol disappeared
during an electronic cut-and-paste procedure. So the real scholium
is:

\begin{Scholium}
Never thrust
%%%the work of
electronic computers (especially when cutting and pasting images
in Adobe Illustrator, while forgetting to unlock stuff).
\end{Scholium}

In reality the scheme in question is prohibited by a sporadic Viro
obstruction (cf. Viro 86 \cite[p.\,67]{Viro_1986/86-Progress}).
Let us now be more serious, Fig.\,f creates besides another
anti-Viro scheme ($\frac{3}{1}\frac{7}{1}\frac{9}{1}$) plus two
schemes constructed by Viro (albeit we were as yet not able to
digest his construction, yet this is merely a matter of detail).
Hence again we get:

\begin{lemma}
Either both of Viro's sporadic obstructions against
$\frac{1}{1}\frac{9}{1}\frac{9}{1}$ and
$\frac{3}{1}\frac{7}{1}\frac{9}{1}$ are false, or there is no
singular octic like our Fig.\,\ref{gurus_4+3_BIS:fig}f.
\end{lemma}

It is clear that the number of such statement can be multiplied
and it is merely now a combinatorial quest to make an exhaustive
list of such hypothetical statement. Though hypothetical we
believe that they could offer new insights on the ultimate destiny
of Hilbert's 16th in degree $m=8$.

\begin{figure}[h]\Figskip
%\vskip-1.2cm\penalty0
%\centering
\hskip-2.7cm\penalty0
\epsfig{figure=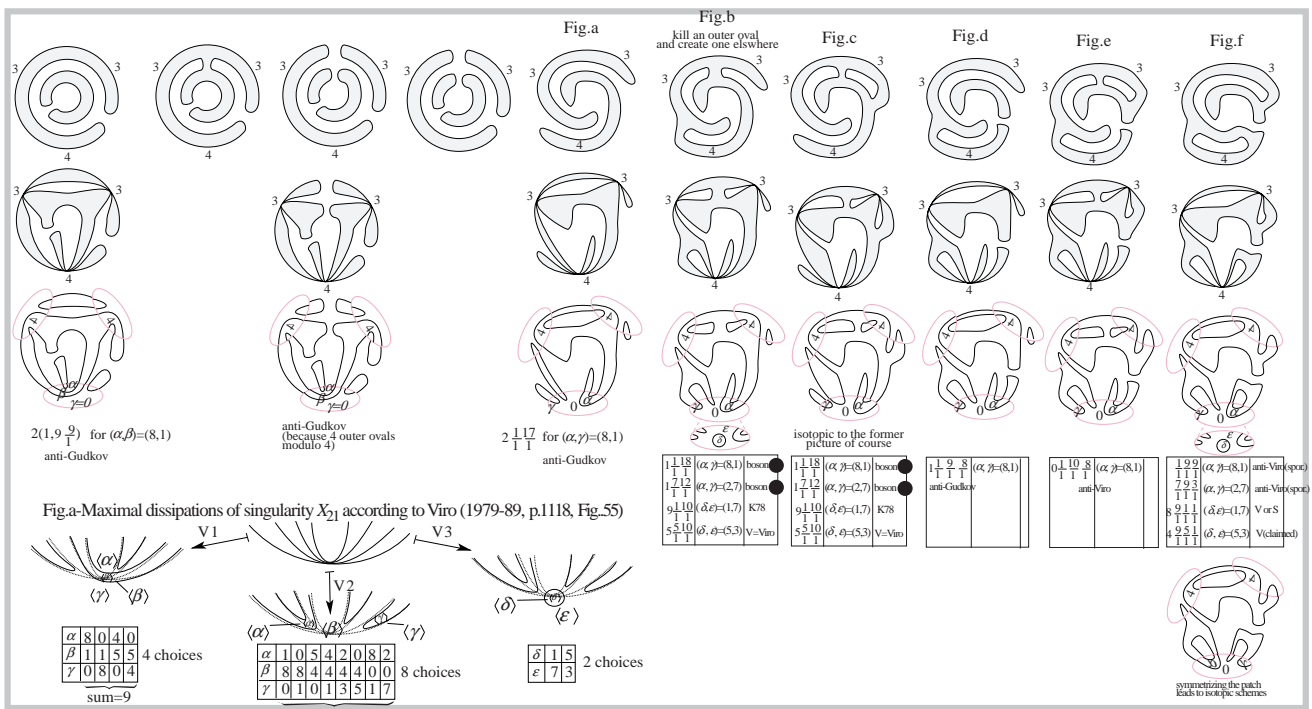,width=172mm} \captionskipAG
  \caption{\label{gurus_4+3_BIS:fig}%
  Gurus 4+3 BIS}
\figskip
\end{figure}

{\it Summary of the method.}---It seems that what is missing is an
organizational way to explore all embryos, i.e. a morphogenesis
(e.g. in R.en\'e Thom's jargon). This is to say, that we could
start from any ground shape and perform adequate smashing to get
fundamental shapes of s-octics, where the prefix ``s'' refers to
singular. Any $s$-octic generates then schemes via Viro's method
of the patch(work), and this permits to explore either new
schemes, or in contradistinction when the scheme is obstructed
either by Harnack, Gudkov, Viro, Orevkov to deduce a corresponding
obstruction on the $s$-octic. Possibly the method can be used to
kill present day obstruction (say of Viro and Orevkov) which are
possibly wrong. If not it is hoped at least that the embryology
method will aid us to detect $s$-octic producing new schemes among
the six bosons not yet known to be realized or prohibited. As a
typical example we believe that Fig.\,\ref{gurus_4+3_BIS:fig}c
could exist thereby producing in one stroke the two bosons
$[1,18]$ and $[7,12]$ without conflicting with Viro nor with
Orevkov's obstructions.

Let us implement our algorithm again. We start with a ring
slightly distorted like a horse shoe (i.e ``U-shaped''). We select
on it three series of anchor points (as the are called in the
theory of B\'ezier curves). (French geometers include too many
baiseurs: B\'ezout, B\'ezier, who else?). Then to create our
distribution of singularity $F4+F3+F3$ we smash respectively 4 or
3 of them to a single point. When we collapse 3 points in a tetrad
there are two options amounting to coalesce either down or up.
Then we merely have to tabulate all options as on
Fig.\,\ref{embryo:fig}.

Fig.\,a enters in conflict once with a sporadic Viro obstruction,
but otherwise create respectable schemes due to Viro/Shustin. So
either Viro's obstruction is false or there is no $s$-octic like
Fig.\,a. Fig.\,b conflicts with Orevkov's obstruction and would
produce one boson $(1,7,12)$. So either Orevkov's obstruction is
false or there is no singular octic like our Fig.\,b. Fig.\,c is
actually isotopic to Fig.\,b. Fig.\,d produces an interesting
boson but conflicts with Orevkov's obstruction. Disappointingly,
the next two figures (Figs.\,e,f) are isotopic to Fig.\,d. Finally
Fig.\,g conflicts with Viro twice and more seriously with B\'ezout
for conics. Hence we can safely claim the:

\begin{lemma}
There is no singular octic with real locus isotopic to
Fig.\,\ref{embryo:fig}g.
\end{lemma}

As yet the method seems confusing yet we have not exhausted all
possibilities as we saw earlier that it is possible to find
embryos creating bosons without conflicting with the Russian
obstructions (which looks to us as random the Russian roulette
game with the pistolet used in transiberian cow-boys movies).

 Then we continued our loose algorithm using as fundamental embryo
 a ``S'' digestive tube like Superman's logo. Fig.\,j gave one
 boson yet conflicting with Orevkov. Varying the position of the
 collapsing segment we arrived eventually at Fig.\,o where 2 bosons
 are created without conflicting with Russian obstructions.
Curiously enough it seems that this curve permit only two
$M$-smoothing (where as usual with Petrovskii $M-$ abridges
Harnack maximality). Then we arrived at Fig.\,p which seems to
lack an $M$-smoothing.

\begin{figure}[h]\Figskip
%\vskip-1.2cm\penalty0
%\centering
\hskip-2.7cm\penalty0 \epsfig{figure=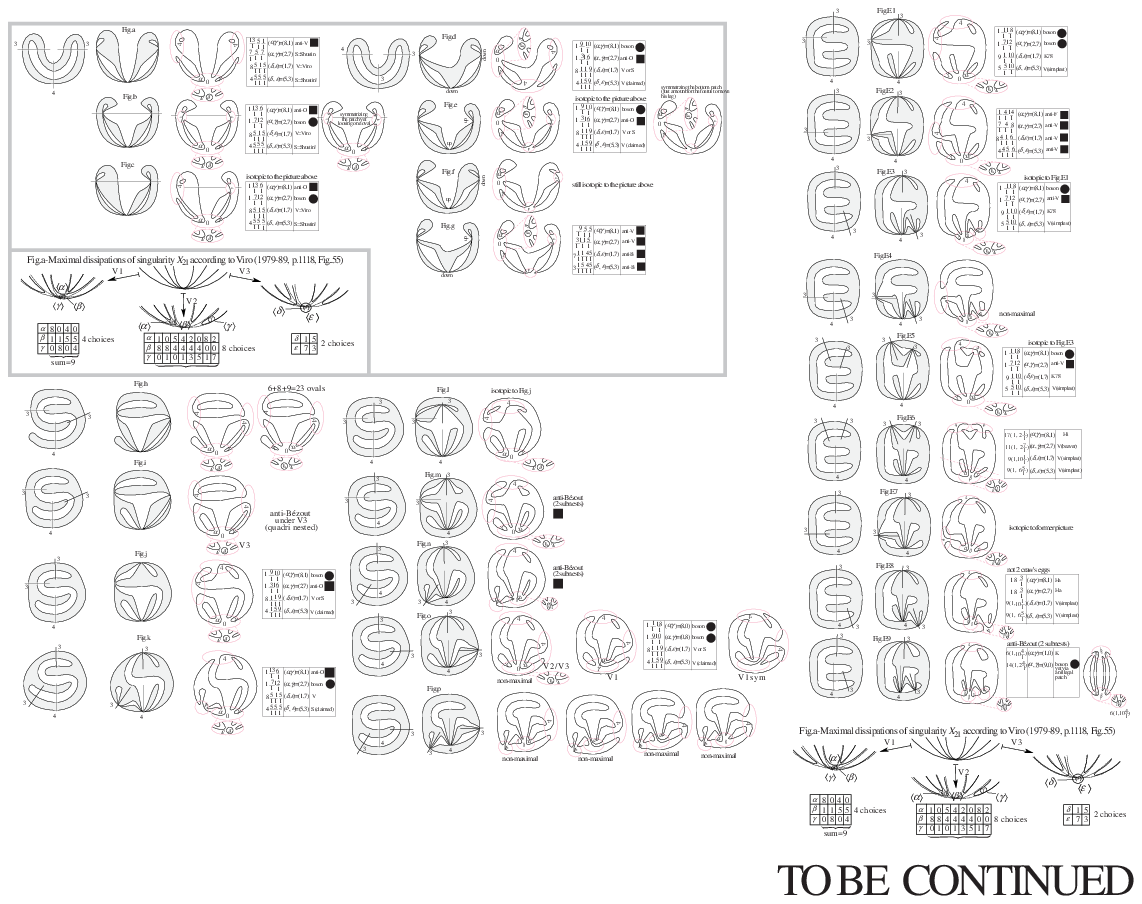,width=172mm}
\captionskipAG
  \caption{\label{embryo:fig}%
  Embryo}
\figskip
\end{figure}

Eventually we decided to change of embryo by tacking an
``E''-letter shaped embryonic substratum for the morphogenesis.
Let us remark from the intrinsic viewpoint that we did not as yet
realized all bosons, so for instance we missed
$1\frac{4}{1}\frac{15}{1}$ if my short-run memory is not failing.
Using this letter we first found Fig.\,\ref{embryo:fig}E1
realizing two bosons without that the other $M$-smoothing
conflicts with Russian scholastic prohibitions (Viro/Orevkov). It
should be remembered that Fig.\,\ref{gurus_4+3_BIS:fig} already
showed this phenomenon of frictionless creationism while involving
exactly the same two bosonic schemes. Therefore let us posit the
following somewhat cavalier:

\begin{Scholium}\label{bosons-elementary-(1_1_18+1_7_12)}
Among all remaining 6 bosons perhaps that
$1\frac{1}{1}\frac{18}{1}$ and $1\frac{7}{1}\frac{12}{1}$ are the
less mysterious one, i.e. they exist algebraically as it possible
to create them out of a plastic curve (singular octic) without
friction against the Russian prohibitions.
\end{Scholium}

Then we got Fig.\,E2 violating twice Fiedler and twice Viro's
extension thereof (imparity law in the trinested case). So safe
for a global mistake in Fiedler-Viro there is no curve like
Fig.\,E2 for fairly deep reasons  beyond immediate B\'ezout
intuition. Fig.\,E4 shows again a phenomenon of non-maximality
apparently allied to the issue that the chosen bridges (to be
collapsed) do not exploit sufficiently the topological contortion
of the ground shape. Intuitively, bridges have to be placed in an
economical fashion so that all isthmus are properly visited. So
for instance a system of efficient bridges is that generating
Fig.\,E5, which alas produces schemes isotopic to those of
Fig.\,E3. Fig.\,E6 shows the result of when the long bridge of
length 4 is pushed inside the embryo. The net results seems to be
that  we land in the subnested realm. Alas the schemes so obtained
are super-classical (Hi, V/beaver or V simplest mandarine method
with the quadri-ellipse). It should be remarked taht on E6 one can
trace a line with 10 intersections, yet perhaps this is only a
defect of our depiction and not an absolute property  of the
isotopy class of our singular octic.

[06.08.13] Let us call a nested oval an egg, as algebraic curves
are like biotopes with several birds constructing the nests and
placing their eggs in them with the possibility that small birds
construct their nests inside one's larger bird nest. (Nid de
moineau dans nid de corbeau.) In the subnested case (i.e. Gudkov
symbol of the form $x(1,y\frac{z}{1})$) the bosonic strip involves
schemes with the minimum number $y=2$ of craw eggs. A priori we
would like to do construction of curves with controlled topology
yet it seems that the piece of information missing are so large
that it is not worth trying to be deterministic, preferring rather
a random search as in each zone there is something to learn, or
de-construct (i.e. mistrust Viro/Orevkov, etc.). So let us look at
Fig.\,E7 which produces a configuration  isotopic to the former
one yet in a more acceptable way w.r.t. B\'ezout.

Fig.\,E8 shows how to access the boson $14(1,2\frac{4}{1})$, yet
by using an illegal patch over a curve violating B\'ezout as soon
as $\ga$ is positive. Yet naive question why is the (maximal)
patch $(\al,\be,\ga)=(9,0,0)$ illegal? We hoped that gluing the
patch with itself along the mandarine yields a contradiction, yet
this produces the scheme $2\frac{19}{1}$. Of course we studied
already this question more systematically on
Fig.\,\ref{ViroDEGREE8_BOSONIC:fig}, and their we noticed that the
existence of the smoothing $V2(9,0,0)$ would only create the 2
bosons $1\frac{1}{1}\frac{18}{1}$ and $1\frac{7}{1}\frac{12}{1}$,
which according to our
Scholium~\ref{bosons-elementary-(1_1_18+1_7_12)} are the most
plausible bosons. Should we therefore mistrust Viro when claiming
that his table of dissipation is complete?

Sorry, but it is at this stage only that we missed to combine the
dissipation V2 and V3 which are perfectly compatible, i.e. married
in  strong Harnack maximality.

[NB 04.08.13: While doing all this picture requires some patience
and skills with Adobe Illustrator. I would like to thanks a
certain artist known as Cathia for giving us the requested
encouragements to look always for the perfect curve. Compare her
illustration for children. Alias I do not remember exactly her
family name, yet googling Cathia, Geneva, illustration for
children is perhaps enough to get something from the web.]

\subsection{Irrational thinking and morphogenesis}

Skip this subsection if you dislike rubbish.

[06.08.13, at 1:21, flashed by some dubious philosophy]. Why do we
like curve? The primate (as says Gromov) is perhaps erotically
attracted by curves. Especially aesthetical are the algebraic
curves requesting only a finitary description and thereby
incarnating a principle  of least action that even nature adhered
on, at least since the vision of Johannes Kepler, ca. 1603. Yet,
why do we
%care
prise nowadays about those rigid objects often overwhelming our
visualization power and seemingly very special and rigid for
modern standards? The Answer is: in die Ruhe des Newton liegt die
(Schwer)kraft! Yes, indeed it seems that it is fairly attractive
to contemplate how an object like an algebraic curve [grooving in
a moduli (better parameter) space of fairly big dimension
$\binom{m+2}{2}-1$], despite its intrinsic rigidness is nearly
able to adopt highly contorted shapes and sees the singularity of
any of its singular representant able to deviate along all
possible dissipation of its singularities in an independent
fashion (alias the
Pl\"ucker-Klein-Harnack-Brusotti-Gudkov-Viro-Shustin principle
incarnating the extreme graphical flexibility of algebraic
objects.) Yet it seems that algebraic geometry is slightly (and
probably violently) more rigid than combinatorics as exemplified
by the case of octics, where we have some Russian obstruction that
probably deserve to be more closely examined in view of our
previous experiments. If not it seems valuable that those highbrow
Russian obstructions (Fidler=Fiedler in RuSSian calligraphy, Viro,
Shustin, Korchagin, Orevkov, nobody else!) receives better
treatment in literature. If not, within the next few month, it is
fairly probable that few of them, especially those given in random
(only semi-published) fashion turns out to be fairly erroneous,
inhibiting thereby the morphogenetic flexibility of algebraischen
Gebilde. The latter to our greatest surprise turns out to be
extremely flexible, and perhaps realize  more schemes that Viro's
original guess suspected. As an historical antecedent of this Viro
``debandade'', we need only to remind the 19 schemes constructed
by Korchagin, to which may be safely added the 4 schemes by
Chevallier (du temple solaire), plus that one exhibited by
Orevkov. It seems now evident that several bosns as well as some
few schemes now believed to be prohibited due to hasty
exhaustivity in Viro's primitive methods (and a lack of
imagination about contortion of algebraic shapes) is responsible
of an actual probably biased state of affairs when it comes to
list octic schemes.

Another idea striking the writer' understanding is the mediocre
trend allied to Newton (cf. e.g. the anti-Riemann viewpoint
expressed in Chevallier 1997
\cite{Chevallier_1997-Secteur-et-def-loc}, that neither the
d\'eploiment universel nor the Riemann surface is requested since
everything is readable from the Newton polygon). This looks to us
severe mysticism or rather lack of poetical continuity method
striving us to the claustrophobic realm of combinatorics or
capitalism. However this suggested us a little flash, when
thinking about Newton's method (the basic one=roots searcher in
one variable) as compared to the Viro-Itenberg polyhedral method
in two variable for tracing curves in the piecewise linear
category. We realize only today that, by virtue of the principle
of linearization in the small (alias calculus or analyse des
infiniment petits) presumably the Viro-Itenberg method is nothing
else than Newton's (parallelogram?) method amounting to detect the
zero loci of algebraic equations trough linearization in the
small.

\subsection{More embryos}

Not yet written, but easy to do-it-yourself.

\subsection{Non-maximal dissipation of the mandarine}

[22.06.13] We shall call the union of 4 coaxial ellipses either
Viro's 0th curve or the {\it mandarine\/}. It would be of interest
to know as well which non-maximal dissipations arise through
smoothings of singularity $X_{21}$. The remarks to be found in
Viro 1989 (p.\,1119) imbue to some vagueness: ``We do not yet have
a complete topological classification of the dissipations of
$X_{21}$ singularities. Shustin [32] proved [\dots]; however there
is still a big gap between what is given by the constructions and
the prohibitions. Curiously, the problem has been completely
solved for dissipations that can occur in the construction of
nonsingular $M$-curves.''

Despite this we can naively postulate (inspired by the
Itenberg-Viro contraction principle of empty ovals) that all
non-maximal dissipations derive from a maximal one via extinction
(=contraction) of an empty (micro) oval. After more mature
thinking this is certainly erroneous (imagine an RKM
$(M-2)$-scheme lying in the depression of Gudkov's sawtooth), but
we really intended to say that that all contractions are realized
algebro-geometrically so as to produce a smoothing. Then we can
naively extend the dissipation pattern to the pre-maximal cases in
order to tabulate the corresponding $(M-1)$- and $(M-2)$-curves.
This requires some tedious tabulation
(Fig.\,\ref{ViroDEGREE8_2:fig}) with the evident drawback is that
there is no absolute warranty about existence of such curves.
Accordingly the corresponding schemes (distribution of ovals) will
%brown
only be yellow-green colored on the main census plate
(Fig.\,\ref{Degree8-(M-i)-curve-TABLE:fig}).
As expectable, working out this tabulation
(Fig.\,\ref{ViroDEGREE8_2:fig}) essentially involves removing the
empty ovals as to get the schemes below Viro's schemes obtained
via the 4 coaxial ellipses (tangent at 2 points). This needs
little commentary apart from a long contemplation of this
tabulation which looks really subsumed to the maximal
dissipations. Perhaps one noteworthy detail is the obtention of
the scheme $18\frac{1}{1}$ (cf. the red case on the 2nd table of
Fig.\,\ref{ViroDEGREE8_2:fig}), which resisted to the other
methods (that will be exposed subsequently). Another little
comment is that at some stage of the tabulation (cf. green-case)
we get a scheme (namely $6\frac{14}{1}$) which is apparently not
dominated by a Viro-style $M$-scheme, but which in reality is via
the $M$-scheme $5\frac{1}{1}\frac{14}{1}$. Eventually, once the
full triangle is filled, it remains a certain rectangle which
combines $M$- with $(M-2)$-smoothings.
%and the latter is depicted
%below to save space.
Rosa-colored case depicts the first occurrence of a scheme (alas
we did not did it systematically from the beginning so our
tabulation is not perfect along this first occurrence option).
Anyway, after completing the tabulation one sees that nearly all
positions dominated by a Viro's $M$-scheme are filled, apart some
few schemes on the right-side of the 1st pyramid like
$\frac{2}{1}\frac{17}{1}$, $\frac{5}{1}\frac{14}{1}$,
$\frac{8}{1}\frac{11}{1}$, etc. It would be
interesting to know if
those schemes are realized.

\begin{figure}[h]\Figskip
%\vskip-1.2cm\penalty0
%\centering
\hskip-3.5cm\penalty0
\epsfig{figure=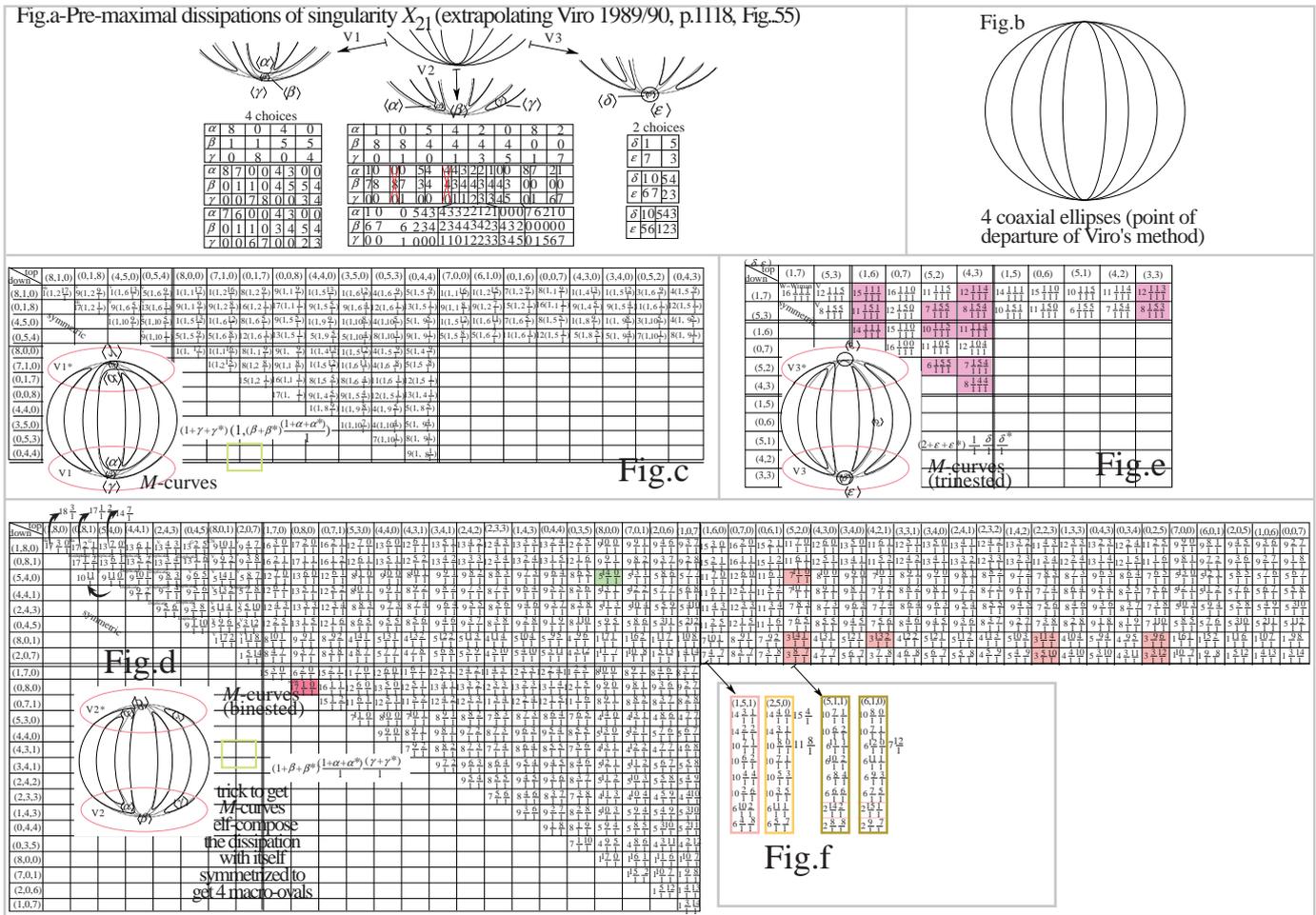,width=182mm} \captionskipAG
  \caption{\label{ViroDEGREE8_2:fig}%
  Non-maximal dissipation of $X_{21}$ and their gluings in degree 8}
\figskip
\end{figure}

\subsection{An error to avoid}

{\it Insertion (of a misconception)} [06.05.13].---The morning
after having told Viro's story to Misha Gabard (born in same year
1948), we wondered why the ``braids'' of
Fig.\,\ref{ViroDEGREE8:fig}a cannot be reconnected in a very
symmetric fashion as to produce a curve $C_8$ with 2 nests of
depth 2, i.e. like the left-half of the right-part of
Fig.\,\ref{ViroDEGREE8:fig}d extended by symmetry to give a fourth
type of ``mandarine'' (with 4 lunes). So we can speculate about a
fourth mode of dissipation in Fig.\,\ref{ViroDEGREE8:fig}a which
is like V3, modulo symmetric reproduction of its left-half. Alas,
it seems that this smoothing was overlooked in Viro 89/90
\cite{Viro_1989/90-Construction}, and that the list of accessory
parameters was not specified. But actually Viro 89/90 writes on
p.\,1119 (4.7.A.) ``and also all of the [quasihomogeneous]
dissipation obtained from them by reflection about the vertical
axis.'' So we get ``new'' smoothing depicted as S2 and S3 on
Fig.\,\ref{ViroDEGREE8_SYM:fig}a, where the accessory parameters
(counting the micro-ovals) are to be chosen as the corresponding
V2 or V3 table. Sorry that I missed this but what is not depicted
is not read (Thurston's philosophy) as opposed to Sullivan's (what
is not written is not read). However when tracing the
corresponding curves using the gluing of S2 with itself we get
schemes violating Gudkov's hypothesis (alias Gudkov-Rohlin
congruence), cf. Fig.\,\ref{ViroDEGREE8_SYM:fig}. Indeed the
general formula for the scheme induced by S2/S2 is
$(\be+\be^\ast)\frac{(1+\al+\ga^\ast)}{1}
\frac{(1+\ga+\al^\ast)}{1}$ whose orientable Ragsdale membrane has
Euler characteristic
$\be+\be^\ast+(1-1-\al-\ga^\ast)+(1-1-\ga-\al^\ast)=
\be+\be^\ast+(-\al-\ga^\ast)+(-\ga-\al^\ast) $, ah sorry too
arithmeticized. Look rather at the table (Fig.\,c), and we see
that we get always schemes violating Gudkov's congruence
$\chi\equiv_8 k^2$. On the left corner of each cases is written
Ga, abridged for Gabard the first constructor of those curves
which do not exist algebraically [sic!], while the right
upper-corner indicates the value of $\chi$ (Euler characteristic
of the Ragsdale membrane bounding the oval from ``inside''.). We
see that $\chi$ runs along a  periodicity mod eight
$14,6,-2,-10,-18$. Similarly, Fig.\,d based on the smoothing using
$S3$ also produces schemes violating Gudkov's hypothesis. In
conclusion, our idea seems a misinterpretation of Viro's prose,
and he really depicted all smoothings $V1,V2,V3$, without
necessity  to consider symmetrized smoothings. His phraseology
just means that we are allowed to take the symmetric of the
asymmetric smoothings.

\begin{figure}[h]\Figskip
%\vskip-1.2cm\penalty0
\centering
%\hskip-0.7cm\penalty0
\epsfig{figure=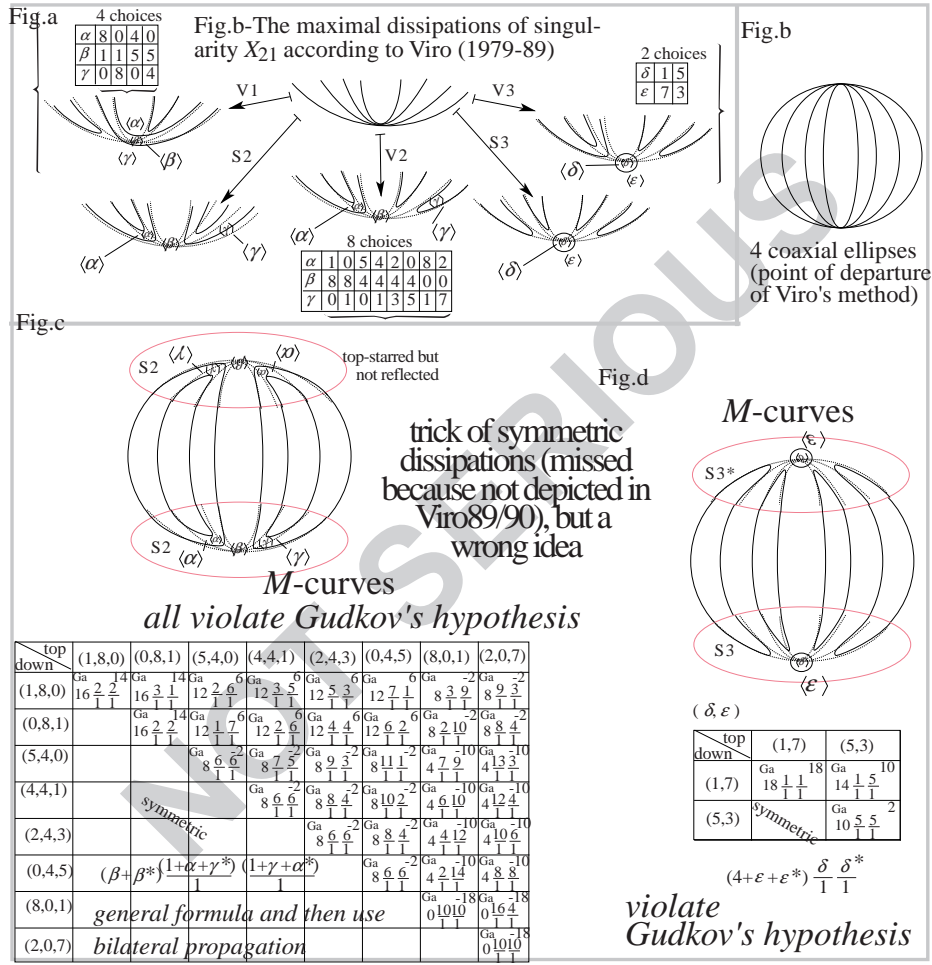,width=122mm} \captionskipAG
  \caption{\label{ViroDEGREE8_SYM:fig}%
  A misinterpretation of Viro 1980 contradicting
  Gudkov's hypothesis. Yet another Irrweg!}
\figskip
\end{figure}

\subsection{Discrepancy between Viro and our
homemade patchwork (our mistake)}

{\it Comment on this subsection}.---It was written before we
realized that Viro's patches C2 and C3 can be glued in a maximal
way. Hence, most of the sequel is based on a basic mistake of us,
and so should be skipped by the reader.

So in reality, the Hawaiian=Leningradian dissipation of $X_{21}$
via a quadruplet of coaxial ellipses creates $35-8=27$ new
schemes. So this is less than the 42 revendicated by Viro 1980,
and therefore the latter had another trick in his pocket.

At this stage we checked Viro's table of 1980, especially
regarding the symbol $2X_{21}$ certificating that the scheme may
be obtained by the above method of dissipation. Safe for our
misunderstanding,  Viro's table contains here another some few
misprints, located by the symbol $\triangle$ on our
Fig.\,\ref{Degree8-M-curve-TABLE:fig}. (Additionally, but less
importantly, for the scheme $10\frac{11}{1}$, Viro's table omits
its realizability via $2X_{21}$, i.e. the above construction.)

How can we explain this discrepancy between our patchwork and
Viro's table? One possible explanation, is that (since we missed
several schemes of the central row) Viro's list of values for
$(\delta,\varepsilon)$ was not complete in the 1989 article.

Explaining the single missing scheme of the last row, namely
$7(1,6\frac{7}{1})$, is more tricky to guess a reason. We checked
once more our first table (Fig.\,\ref{ViroDEGREE8:fig}d, left)
which looks perfectly correct.

Perhaps all these discrepancies may be ascribed to the typographer
of Viro's article 1980, or we missed some detail. A next step is
to understand Viro's other constructions (i.e. not via
$2X_{21}$=dissipation of a quadritangent quadruplet of ellipses).

\section{Artwork (Viro, Shustin)}

\subsection{Viro's twisted constructions (beaver)}

We now explain other constructions due to Viro starting with a
more complicated singular octic than the fundamental union of 4
coaxial ellipses. By these somewhat more complicated procedure
some few other sporadic curves will be gained, yet alas not
covering in full the Hilbert's 16th problem.

Again the key idea is explained in Viro 1989/90
\cite{Viro_1989/90-Construction} (especially his Fig.\,77), where
4 new schemes are obtained (see rhombic ``V'' on
Fig.\,\ref{Degree8-M-curve-TABLE:fig}), by a clever construction
we shall now attempt to summarize. Alas, Viro's original picture
seems to contain minor bugs, and we hope our presentation being
more reliable.

First, Viro constructs a singular octic with 2 singularities which
are respectively  an ordinary (nondegenerate) 5-tuple point
($N_{15}$ in Arnold's census) and a triple branch each having 2nd
order contact ($J_{10}$ in Arnold's census), compare
Fig.\,\ref{ViroDEGREE8_TRICKY:fig}a. We shall detail this
construction in the sequel, but first show its
%application
utility to Hilbert's 16th in degree $m=8$. It is a good exercise
to see how the complicated branch is unicursally travelled by a
particle according to a simple law, namely whenever we cross
$J_{10}$ we stay on the branch having the same ``curvature'',
while when crossing the ordinary 5-fold node $N_{16}$ (locally
like 5 concurrent lines) we always have to count 5 branches
cyclically to find our way out of the singularity. This singular
curve $C_8$ is constructed via a hyperbolism \`a la Newton
(essentially akin to Gudkov's clever use of Cremona transformation
ca. 1972). Let us trace first the singular octic in question
(Fig.\,\ref{ViroDEGREE8_TRICKY:fig}a), and the game will be as
before to dissipate the 2 singularities in all possible fashions
while
%%%listing
tabulating all resulting schemes (distribution of ovals).

\begin{figure}[h]\Figskip
%\vskip-1.2cm\penalty0
\centering
%\hskip-0.7cm\penalty0
\epsfig{figure=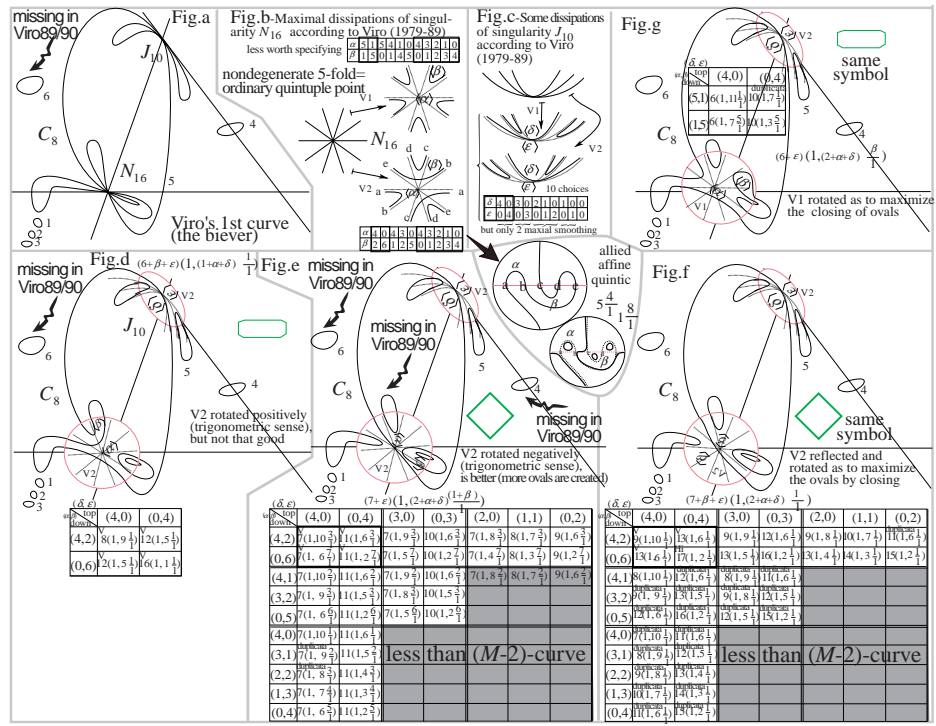,width=122mm} \captionskipAG
  \caption{\label{ViroDEGREE8_TRICKY:fig}%
  Detailing Viro 1980 via Viro89/90 (while correcting
  his picture).}
\figskip
\end{figure}

The next step is to remember the dissipation of $J_{10}$ (cf.
Fig.\,c). Those were already used when exploring  Viro's gluing in
degree 6 hence no further
%%%explanation
ado is required, safe that $(\al,\be)$ went relabelled
$(\delta,\varepsilon)$. For $N_{16}$ the maximal dissipation are
given on Fig.\,b following Viro89/90 \cite[p.\,1109,
Fig.\,34]{Viro_1989/90-Construction}. There we read that when
$\al+\be=6$ on V2, then $\al-\be \equiv 2 \pmod 8$ and this leaves
only the 2 possibilities listed. The case of V1 is not even needed
%%%%in the sequel as we shall see, because
as it does not create $M$-curves, but is of some interest if
attention is paid to $(M-2)$-curves too.

We first traced Fig.\,\ref{ViroDEGREE8_TRICKY:fig}d by gluing the
V2 dissipation of $N_{16}$ (ordinary quintuple point) rotated so
as to close the double petal to a nest of depth 2. Alas doing so
the upper petal is not closing, but rather connecting with the
west petal (creating thereby a non-maximal smoothing). Albeit
``less'' interesting than $M$-curves, this still gives three
$(M-2)$-schemes worth tabulating and  reporting on
Fig.\,\ref{Degree8-(M-i)-curve-TABLE:fig}.

On rotating clockwise the V2-template at $N_{16}$, we get
Fig.\,\ref{ViroDEGREE8_TRICKY:fig}e, which leads to four
$M$-curves,
%%tabulated on the diagram
namely $7(1,10\frac{3}{1})$, $11(1,6\frac{3}{1})$,
$7(1,6\frac{7}{1})$ and $11(1,2\frac{7}{1})$. Those  4 new schemes
(due to Viro 1980) are reported on
%%%the map
Fig.\,\ref{Degree8-M-curve-TABLE:fig} via rhombic squares
enclosing the letter ``V''. It seems evident that no more
$M$-schemes can be
%%%gained
catched by this method, because the smoothing $V1$ of the
quintuple point $N_{16}$ fails creating the maximum number of
ovals. (We warn the reader that some details of Viro's picture in
1989/90 (esp. his Fig.\,77) looks anomalous albeit his final
(symbolical) results are
%%% right (
in accordance with our own checking. The graphical omissions we
detected are marked by squig-arrows on our figures.) For
subsequent
%%%interest,
purposes, it is worth tabulating the non-maximal smoothings too of
Fig.\,e. One can modulate further by examining also non-maximal
smoothing of $N_{16}$, for which we must consult Viro 89/90. There
on p.\,1109, the structural constant of the dissipation are
prescribed by a congruence $\pmod 8$ which leads easily to the
values
%%%we
tabulated above on Fig.\,b. Actually, for an
$(M-2)$-smoothing it seems that there is no obstruction and the
full range of pairs is realized from $(4,0)$, $(3,1)$, etc. up to
$(0,4)$. Interestingly, it may be observed that those smoothings
are just interpretable as derived from the maximal smoothings
through contraction of the empty micro-ovals. This explains how
$(5,1)$ creates $(5,0)$ and $(4,1)$, and so on. Subsequently, only
$(2,2)$ is of a new stock and not created by contraction. This is
surely allied to  Rohlin's maximality conjecture for RKM-scheme.
Further there are probably other smoothings leading to more
schemes.

The V1-smoothing  leads (alas) only to \hbox{$(M-2)$}-curves
depicted on Fig.\,\ref{ViroDEGREE8_TRICKY:fig}g which
%%are
were new (at the time we found them). Then we can progressively
rotate the
%%%gluing-template
patch as shown on
Fig.\,\ref{ViroDEGREE8_TRICKY_SUITE_XXX_15_4-1:fig}h,\allowbreak
i,j,k,l,m,n,o,p,q,r,s, but
%%%the resulting schemes yield
nothing new is obtained.

\begin{figure}[h]\Figskip
%\vskip-1.2cm\penalty0
%\centering
\hskip-1.7cm\penalty0
\epsfig{figure=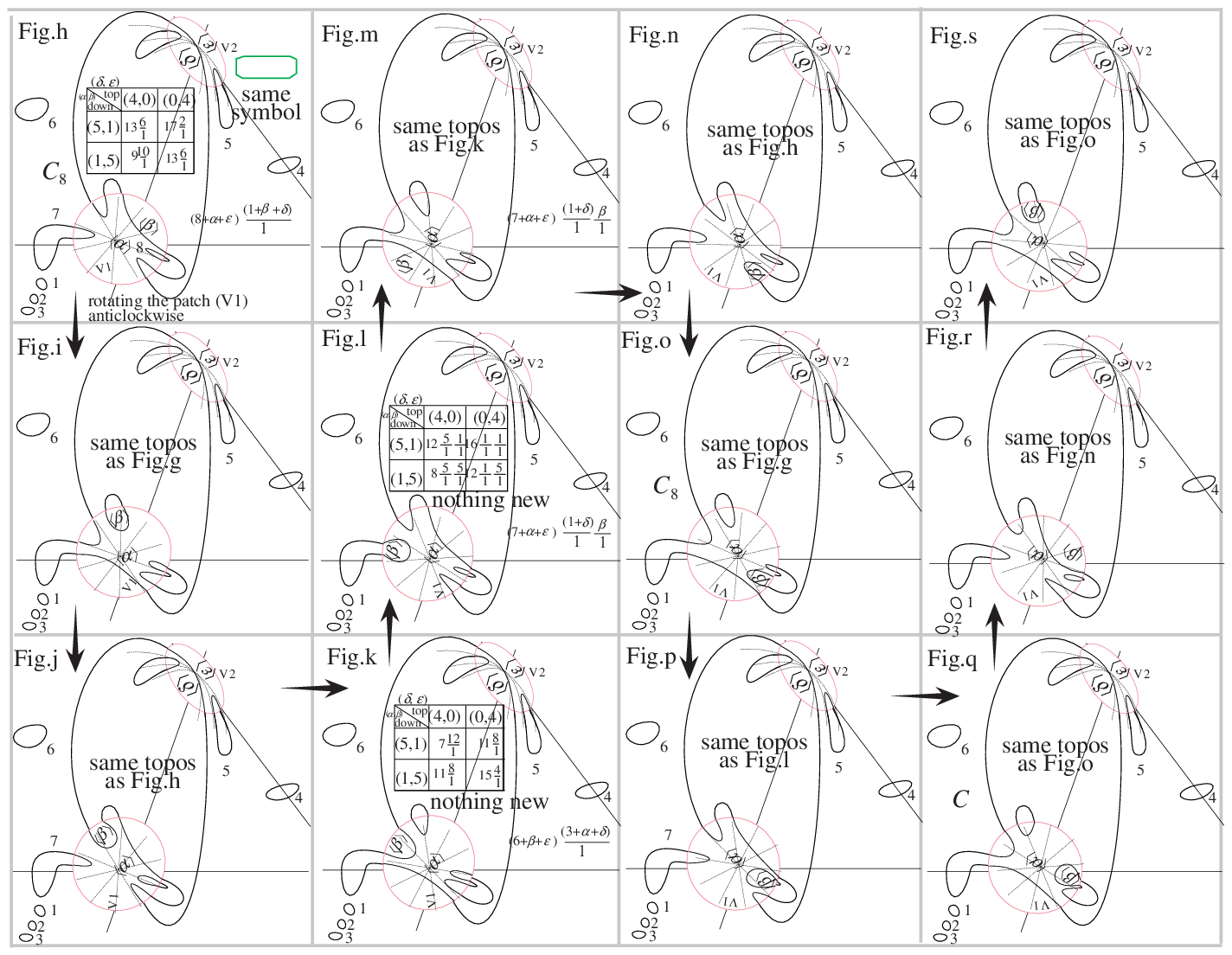,width=152mm}
\captionskipAG
  \caption{\label{ViroDEGREE8_TRICKY_SUITE_XXX_15_4-1:fig}%
  Little twists of Viro's construction giving some $(M-2)$-curves.}
\figskip
\end{figure}

Then it is also reasonable to reflect the patch and rotate it to
get the curves of Fig.\,\ref{ViroDEGREE8_TRICKY_SUITE_B:fig}. The
series of curves so obtained is extremely boring producing not a
single new curve. Actually, the symmetry of the patch V1 affords
probably a metaphysical explanation a priori for this lack of
newness.

\begin{figure}[h]\Figskip
%\vskip-1.2cm\penalty0
%\centering
\hskip-1.7cm\penalty0
\epsfig{figure=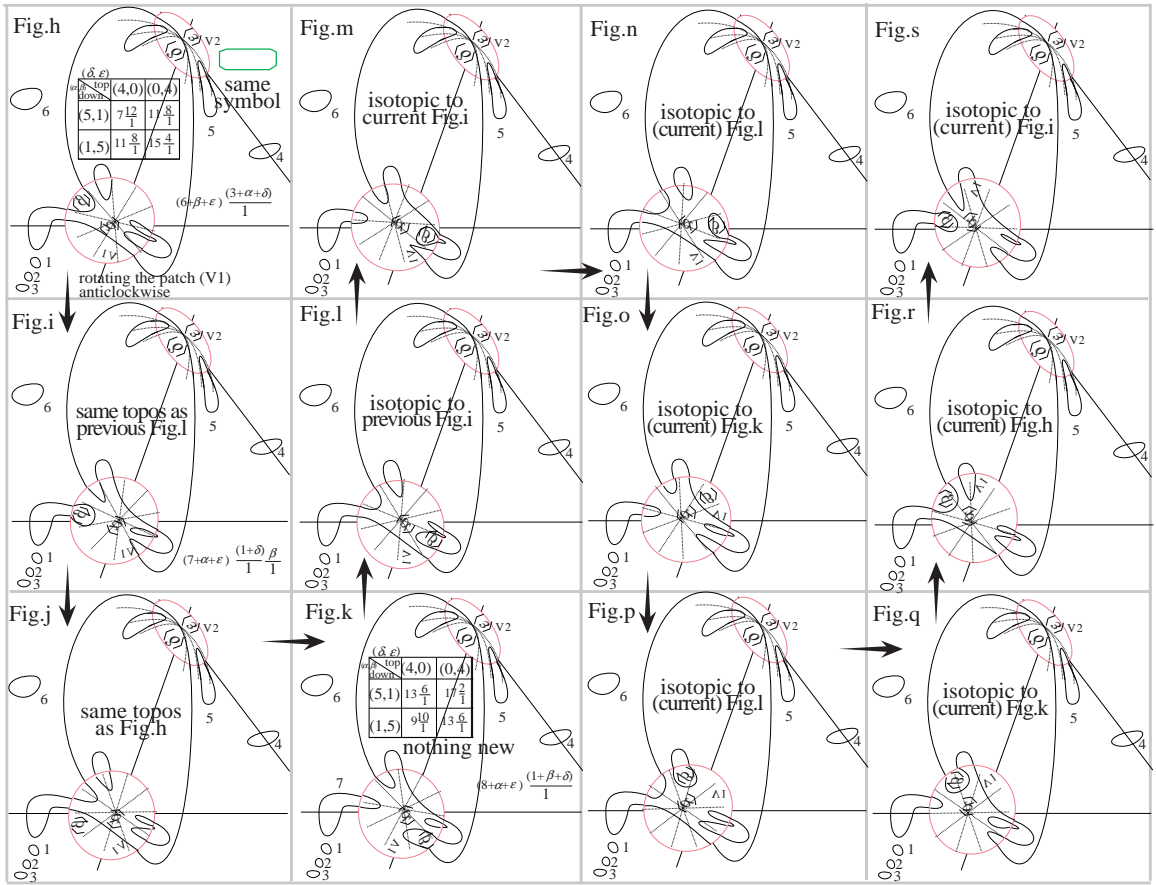,width=152mm}
\captionskipAG
  \caption{\label{ViroDEGREE8_TRICKY_SUITE_B:fig}%
  Reflecting the patch V1 at the bottom}
\figskip
\end{figure}

 So we get:

\begin{lemma}\label{RMC:cter-example-via-Viro-1st-curve=BEAVER}
Considering Fig.\,\ref{ViroDEGREE8_TRICKY_SUITE_XXX_15_4-1:fig}k
or Fig.\,\ref{ViroDEGREE8_TRICKY2:fig}f creates interesting
RKM-schemes. (Those are also obtainable via a Shustin curve, see
the subsequent Fig.\,\ref{ViroDEGREE8_SHUSTIN2:fig}.) Hence Viro's
method applied to Viro's 1st
%%%%2nd
curve (alias the beaver) realizes $3$
%%% one-nested
simply-nested RKM-schemes, namely $15\frac{4}{1}$,
$11\frac{8}{1}$, $7\frac{12}{1}$. Further
%as we shall see later (Theorem~{\rm
%\ref{RMC:cter-example-via-Shustin}})
any one of those (simply-nested) schemes suffices to refute
Rohlin's maximality conjecture.
\end{lemma}

\begin{proof}
Only the second clause deserves
%comments.
commenting upon. For instance we may consider the RKM-scheme
$15\frac{4}{1}$. Aided optionally by the geography of
Fig.\,\ref{Degree8-(M-i)-curve-TABLE:fig} (or better its enlarged
version Fig.\,\ref{Degree8-(M-i)-curve-TABLE_I:fig}), we have an
$(M-1)$-enlargement $13\frac{2}{1}\frac{4}{1}$, where so-to-speak
two among the 15 free ovals went captured
%%%enclosed
by a loop as to become ``nested''. In turn this scheme may be
enlarged to the $M$-scheme $13\frac{2}{1}\frac{5}{1}$ which
(though historically first constructed by Gudkov) may be
constructed \`a la Viro  via dissipation of the quadri-ellipse
(cf. our previous Fig.\,\ref{ViroDEGREE8:fig}). The proof is
complete and Rohlin's conjecture refuted.
\end{proof}

{\it Philosophico-bibliographical comment.}---It is fairly
puzzling that this result holds true because as far as we could
interpret the literature this problem was still open yet simpler
to settle than the refutation of the reverse sense of Rohlin's
maximality conjecture (abridged RMC henceforth). Recall that RMC
posits that a scheme is of type~I (purely orthosymmetric) iff it
is maximal (in the hierarchy of all schemes in the prescribed
degree). The reverse sense $\Leftarrow$ was refuted by Shustin,
while the other sense seemed to stay open. However in Viro's
survey of 1986 \cite{Viro_1986/86-Progress} the text is fairly
confusing on this topic and it seems that there is an interversion
of logical data when it comes to this topic, cf. p.\,?? for the
exact passage. As a vague guess, it may be argued that Viro stayed
confuse on this point because he did not wanted to attack
frontally the conjecture of his beloved teacher (Rohlin). So
either sentimentality or over-modesty contributed to add confusion
to the topic. It is also possible that our above refutation is
new, but since it is merely  based on old techniques of Viro, what
is new is not the geometric quintessence but merely the
combinatorial skills allied to the contemplation of the
%%% combinatorial structure
architecture of higher Gudkov pyramids (notably
Fig.\,\ref{Degree8-(M-i)-curve-TABLE:fig}).

{\it Added\/} [28.06.13].---At this stage it seems fairly
plausible that the top-row RKM-schemes (i.e. the list of
$(M-2)$-symbols $15\frac{4}{1}$, $11\frac{8}{1}$, $7\frac{12}{1}$,
$3\frac{16}{1}$) should also be constructible via dissipation of
the quadri-ellipse if we knew the (extended) accessory parameters
of the smoothings. Indeed by analogy with the smoothings of
$N_{16}$-singularity, those of singularity $X_{21}$ must also have
mutant species not derived by emptifying (=contracting) the
micro-ovals. As a consequence our previous tabulation of the
$X_{21}$-smoothings (Fig.\,\ref{ViroDEGREE8_2:fig}) cannot claim
exhaustiveness, and it could be an interesting duty to make it
complete. At this stage we could choose more or less random values
of the parameters and experiment what global curve is resulting
thereof. Since it is the V2-dissipation which leads to
simply-nested curves, we consider the first listed
$(M-2)$-dissipation, i.e. $(1,6,0)$ and alter it to the nearby
value $(1,5,1)$. Then we get many new $(M-2)$-schemes of type RKM,
that are tabulated on Fig.\,\ref{ViroDEGREE8_2:fig}f. Similarly if
we can twist the parameters to $(2,5,0)$, then we obtain the
simply-nested RKM-schemes $15\frac{4}{1}$ and $11\frac{8}{1}$ that
suffices to corrupt Rohlin's maximality conjecture (RMC). So we
have the following hypothetical lemma:

\begin{lemma}
Viro's most elementary method involving dissipation of the
quadri-ellipse suffices to corrupt Rohlin's maximality conjecture
provided singularity-$X_{21}$ admits the dissipation V2 of
Fig.\,\ref{ViroDEGREE8_2:fig}a with parameters
$(\al,\be,\ga)=(2,5,0)$.
\end{lemma}

Turning back to Viro's 1st curve, we can also rotate the patch to
obtain Figs.\,\ref{ViroDEGREE8_TRICKY2:fig}g,h,i,j,k,l,m, yet what
is so obtained is not
%extremely
strikingly original.

\begin{figure}[h]\Figskip
%\vskip-1.2cm\penalty0
%\centering
\hskip-1.7cm\penalty0
\epsfig{figure=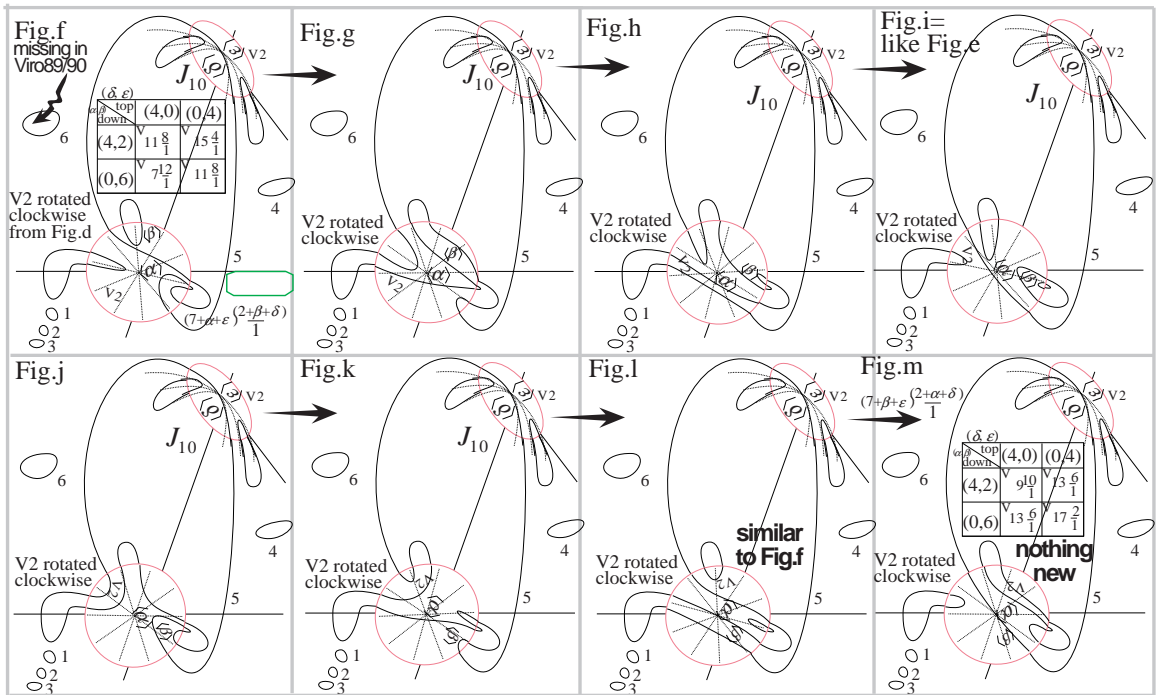,width=152mm} \captionskipAG
  \caption{\label{ViroDEGREE8_TRICKY2:fig}%
  Counterexample to RMC via Viro, plus rotating the patch.}
\figskip
\end{figure}

Finally we can also rotate the reflected patch V2 to get the
following series of Fig.\,\ref{ViroDEGREE8_TRICKY3:fig}. Some few
new species do occur denoted by ``new''. Usually no more comment
would be requested but again TeX is unhappy and is overflowed by
the avalanche of pictures versus the little of text we have to
say. Of course geometry is the art of staying silent in front of
the beauty of the landscape to be contemplated. So let us be
silent and write some anodyne prose to get our text synchronized
with the figures.

\begin{figure}[h]\Figskip
%\vskip-1.2cm\penalty0
%\centering
\hskip-1.7cm\penalty0
\epsfig{figure=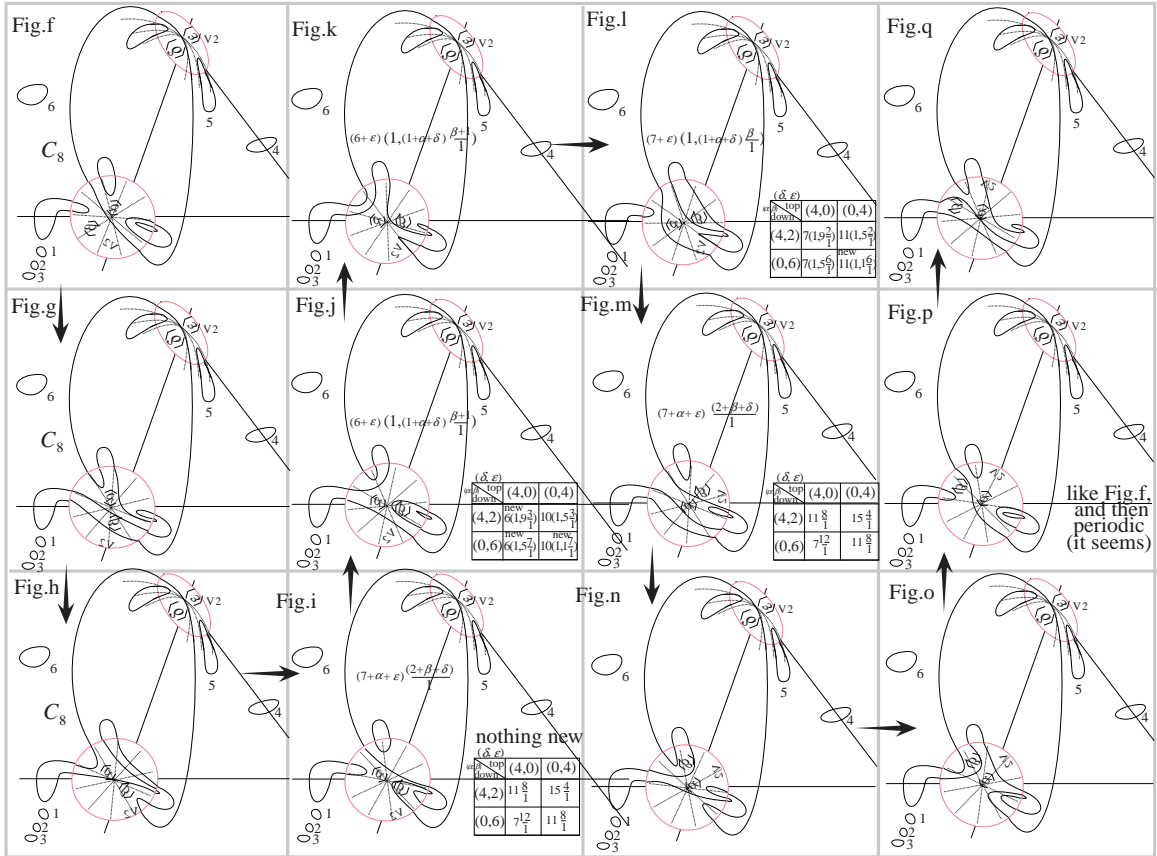,width=152mm} \captionskipAG
  \caption{\label{ViroDEGREE8_TRICKY3:fig}%
  Rotating the reflected patch to get some new $(M-2)$-curves.}
\figskip
\end{figure}

To construct the singular octic (of
Fig.\,\ref{ViroDEGREE8_TRICKY:fig}a) Viro rests on the following
picture (Fig.\,\ref{ViroDEGREE8_SING:fig}) in spirit akin to
Gudkov's trick. The first step is a Harnack-style vibration
creating a cubic $C_3$ oscillating 6 times across the conic $C_2$.
The 2nd step involves a (partial) smoothing of $C_3 \cup C_2$ to
get a quintic $C_5$ with a unique ordinary double point. On Viro's
figure (Fig.\,76 of Viro 89/90) a branch of the 6th circuit is
missing.
%%%(not that grave).
The 3rd step involves a hyperbolism \`a
la Newton (which
%%%is really
reminisces Gudkov's use of Cremona transformation). The strict
transform of the $C_5$ under this map is an octic because the
pull-back of a line is a conic through the 3 basepoints of the
fundamental triangle,
%%and the latter
which cuts the $C_5$ along 8 mobile points since 2 of them are
statically located on the unique singularity of the $C_5$. Our
Fig.\,d shows one additional oval that was overlooked on Viro's
picture (Fig.\,73 of Viro 89/90). Beware that our location of the
6th oval might be exotic but it
%%%is
lies certainly outside of the singular circuit, at least on behalf
of the sequel of Viro's text.

\begin{figure}[h]\Figskip
%\vskip-1.2cm\penalty0
%\centering
\hskip-2.7cm\penalty0
\epsfig{figure=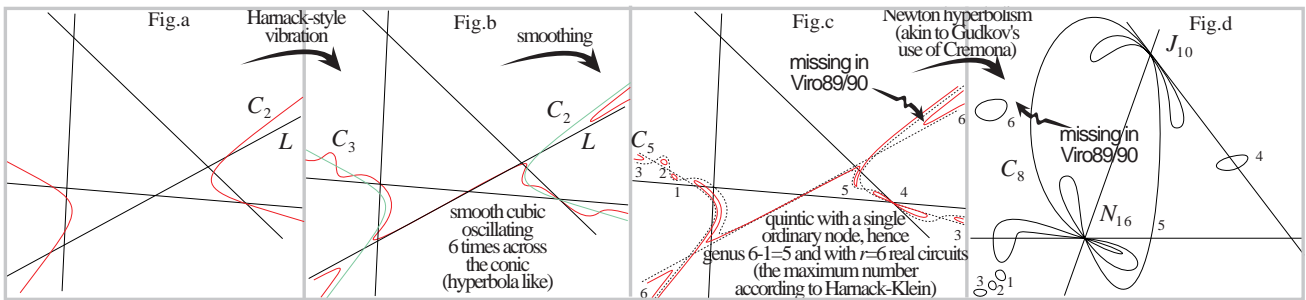,width=172mm} \captionskipAG
  \caption{\label{ViroDEGREE8_SING:fig}%
  Creating Viro's 1st curve (from Viro 89/90 but
  corrected picture)}
\figskip
\end{figure}

\subsection{Viro's 2nd construction (horse)}

[05.05.13] Finally, Viro proposes a 2nd fundamental curve $C_{8}$
(cf. our Fig.\,\ref{ViroDEGREE8_HORSE:fig}a based upon
p.\,1129--30 of Viro 89/90 \cite{Viro_1989/90-Construction}) which
leads to another series of $M$-octics. This curve resembles the
profile face of a horse, hence refer to this curve as Viro's 2nd
curve or the horse.
%This requires to be
%presented in some more details but the idea is essentially the
%same.
%

\begin{figure}[h]\Figskip
%\vskip-1.2cm\penalty0
\centering
%\hskip-0.7cm\penalty0
\epsfig{figure=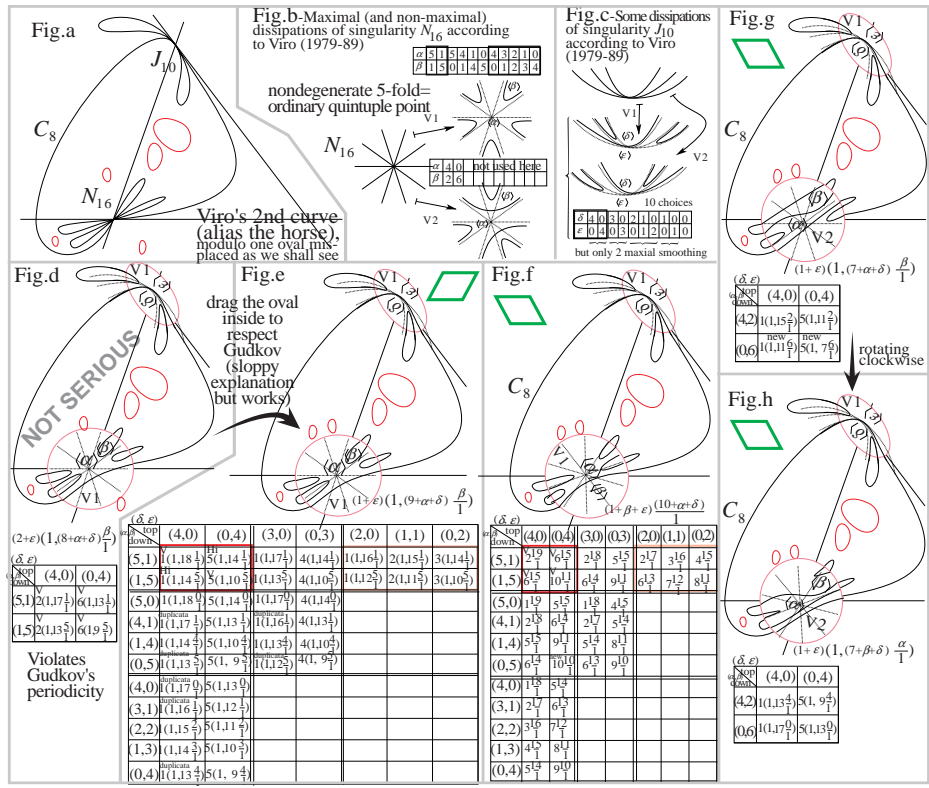,width=122mm} \captionskipAG
  \caption{\label{ViroDEGREE8_HORSE:fig}%
  Following Viro 89/90 yet correcting
  his picture.}
\figskip
\end{figure}

First, Viro constructs another singular octic (compare his
Fig.\,78 specialized to $k=2$, or our
Fig.\,\ref{ViroDEGREE8_HORSE:fig}a). We differ Viro's construction
of this ground curve to later, to first work out the patchwork. On
gluing the dissipation V1 of the quintuple point $N_{16}$ where
the allowed parameters $(\al, \be)$ with $\al+\be =6$ have to
satisfy the congruence $\al-\be\equiv 4 \pmod 8$ (hence restricted
to take the values $(5,1)$ or $(1,5)$ as tabulated on Fig.\,b), we
get Fig.\,d after choosing the appropriate V1 dissipation of the
triple point $J_{10}$. Alas the curves so constructed do not
%enter
fit on the tabulation (Fig.\,\ref{Degree8-M-curve-TABLE:fig}), the
intrinsic reason being that those $M$-schemes violate Gudkov's
hypothesis ($\chi\equiv_8 k^2$). A fairly simple explanation is
that we wrongly located the oval of the fundamental curve on
Fig.\,a. It suffices indeed to transfer the outer oval of $C_8$
inside of the singular circuit to get Fig.\,e, which create 4 new
$M$-schemes, namely $1(1,18\frac{1}{1})$, $5(1,14\frac{1}{1})$,
$1(1,14\frac{5}{1})$, $5(1,10\frac{5}{1})$. Actually, the 2nd and
3rd one were first realized by Hilbert's construction, yet the 1st
and last one are pure creation of Viro 1980. As usual, we report
the geography of those scheme on the main table
(Fig.\,\ref{Degree8-M-curve-TABLE:fig}) by using this time a
green-parallelogram enclosing the letter ``V'' honoring as usual
Viro.

As before, we may rotate the patch, as much as we please, and do
reflections. This is fairly tedious (space-consuming) to depict
and it is not impossible that there is a more expediting way to
construct those schemes maybe via Viro's most elementary
4-ellipses method if we knew the possible non-maximal smoothings.
In particular upon rotating the bottom patch $V2$, we get the
following series of curves (Fig.\,\ref{ViroDEGREE8_HORSE2:fig})
extending the former Fig.\,g, yet yielding nothing tremendously
revolutionary, since the obtained isotopy types seem subsumed to a
law of repetition yielding a poor level of bio-diversity.

\begin{figure}[h]\Figskip
%\vskip-1.2cm\penalty0
%\centering
\hskip-2.7cm\penalty0
\epsfig{figure=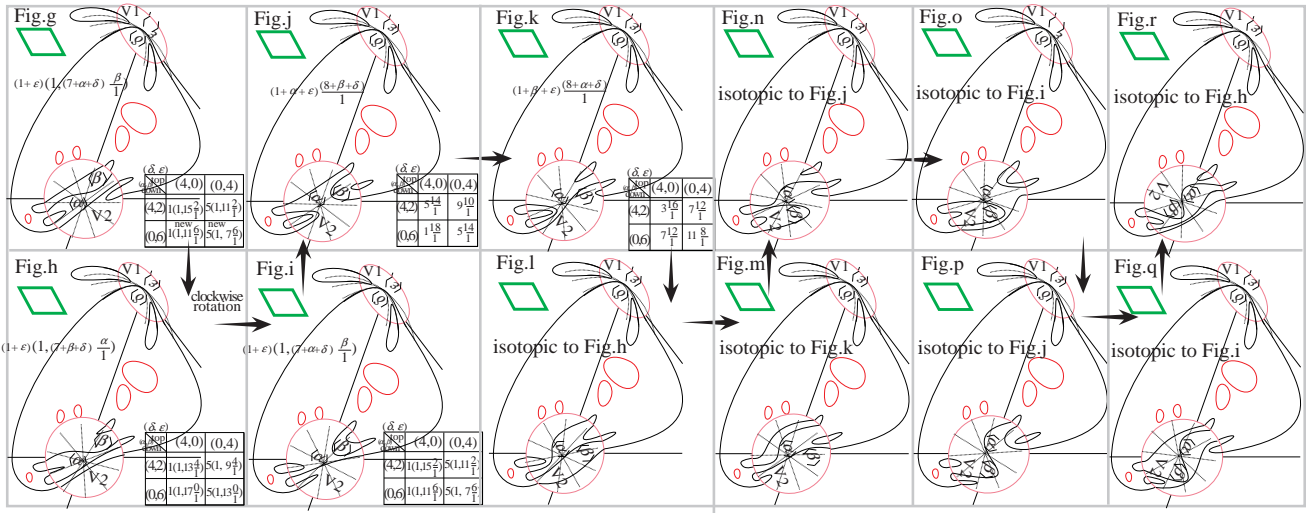,width=172mm} \captionskipAG
  \caption{\label{ViroDEGREE8_HORSE2:fig}%
  Rotating the bottom patch V1, while getting few new species.}
\figskip
\end{figure}

It remains then to symmetrize the bottom patch V2, and also to
explore the rotations of the (symmetric) V1 patch (dispensing us
to consider its reflection). Let us first rotate V1 (say in the
clockwise sense), starting from Fig.\,e of the main former figure
(Fig.\,\ref{ViroDEGREE8_HORSE:fig}), we get the following series
of curves (Fig.\,\ref{ViroDEGREE8_HORSE3:fig}). The first so
obtained (Fig.\,f) is not even worth tabulating as there is zero
(naught) ovals created nearby the bottom-patch, so that only
$(M-4)$-curves are created. Continuing the rotating process we
eventually arrive at Fig.\,k, which violates Gudkov's hypothesis
(conclusion found the [21.06.13]), e.g. because it is not
catalogued on our Table
(Fig.\,\ref{Degree8-(M-i)-curve-TABLE:fig}). This is fairly
puzzling and it is tricky to locate the plague of the reasoning.
Maybe Viro's 2nd curve is a hallucination (so-called
phatamorgana=mirage in German)? It should be remained that the
fundamental curve of Fig.\,\ref{ViroDEGREE8_HORSE:fig}a had
precisely a defect w.r.t. Gudkov periodicity that was remedied
upon dragging one outer oval inside, yet this naive trick turns
out to create another tension with Gudkov at the later level of
rotation of the patch. Oh sorry, it seems rather that we made a
fatal mistake when moving from Fig.\,i to Fig.\,j. Correcting this
defect we get Fig.\,j-star (and so on), yet no tremendous
gold-mine is discovered along this way giving only curves isotopic
to Fig.\,e.

\begin{figure}[h]\Figskip
%\vskip-1.2cm\penalty0
%\centering
\hskip-2.7cm\penalty0
\epsfig{figure=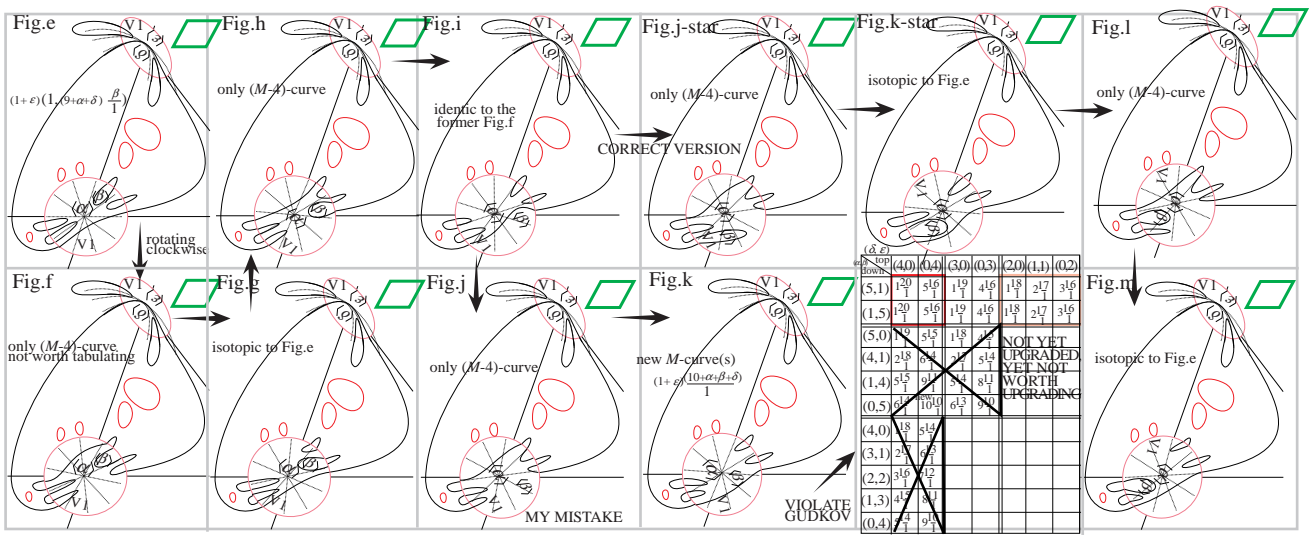,width=172mm} \captionskipAG
  \caption{\label{ViroDEGREE8_HORSE3:fig}%
  Rotating the bottom patch V2, while getting little new species.}
\figskip
\end{figure}

We can then explore the symmetrized bottom patch V2, while
performing a rotation initiating with
Fig.\,\ref{ViroDEGREE8_HORSE4:fig}g below (reflecting the original
Fig.\,\ref{ViroDEGREE8_HORSE:fig}g). This gives the following
series of curves (Fig.\,\ref{ViroDEGREE8_HORSE4:fig}). Seen
dynamically all this may be interpreted as the mastication of a
herbivore, typically a horse whose resemblance with Viro's 2nd
curve is self-explanatory. Again all this sentimental prose has to
be introduced for otherwise there is a boring shift (d\'ecalage
between pictures and text), due to TeX's rigidity. At any rate the
conclusion is that reflecting the patch leads apparently always to
configurations isotopic to those listed on the previous
tabulation. Hence nothing original is created. After tracing
Fig.\,l it is already evident that nothing tremendous will appear
in the firmament, yet as we were
%%hung-over and
fairly tired and bored by the game we decided to continue to be
sure to miss nothing. Upon continuing up to Fig.\,r it seems clear
that we explored everything albeit we did not closed completely
the ``loop'', or that some periodicity (in the sense of a boring
repetition) start to predestine the story. Hence it is intuitively
clear that no new schemes are created along this procedure.

\begin{figure}[h]\Figskip
%\vskip-1.2cm\penalty0
%\centering
\hskip-2.7cm\penalty0
\epsfig{figure=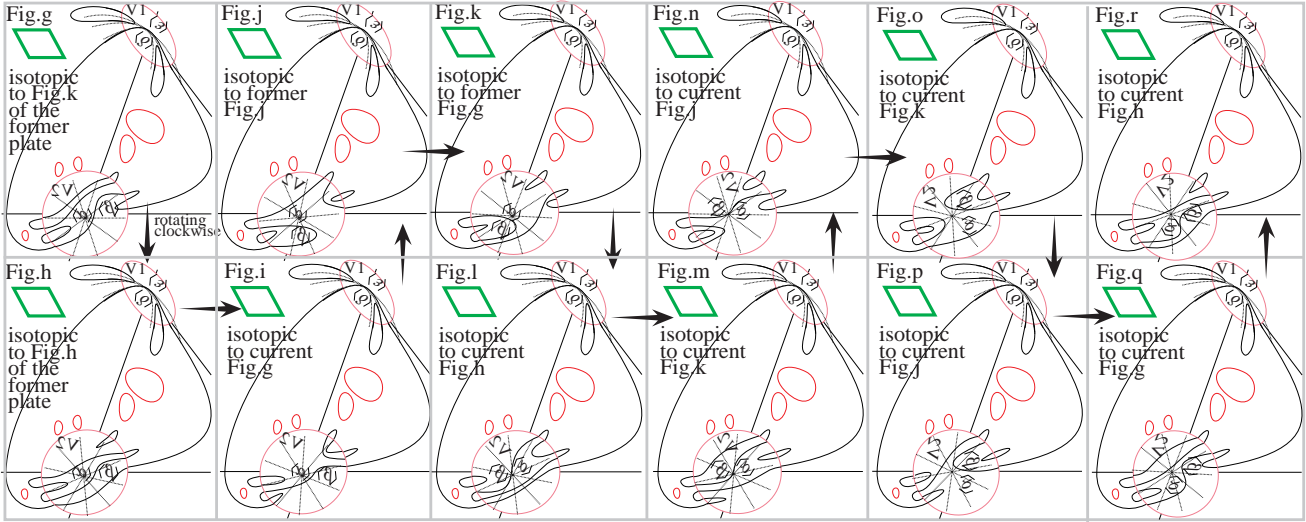,width=172mm} \captionskipAG
  \caption{\label{ViroDEGREE8_HORSE4:fig}%
  Reflecting the bottom patch, and then rotating it.}
\figskip
\end{figure}

At this stage it seems that we have exploited all possible
smoothings leading to $(M-2)$-curves (or better) of Viro's (3rd)
``horse'' curve.

%%%%%   --HIER WWWWWWWWWWWWWWWWWWWWWEITER

It remains to explain Viro's construction of this auxiliary
singular octic curve. Again the methodology is nearly the same and
seems to owe some inspiration from both Brusotti and Gudkov. Again
we follow Viro 89/90. We get started with a cubic $C_3$
oscillating across the axes $L$ and $L_1$ (lines) as on
Fig.\,\ref{ViroDEGREE8_SING_BIS:fig}a. Then $C_3\cup L$ is
smoothed to a $C_4$ of Fig.\,b. Here Viro's figure seems to miss
an oval. Then $C_4\cup L$ is vibro-smoothed \`a la Harnack to the
quintic of Fig.\,c. Further it seems that Viro proposes to
contract the oval 1 to a solitary node (isolated double point so
as to ensure that the strict transform under the Newton-Cremona
transformation will have degree 8 when one of the 3 basepoints of
the pencil of conics is chosen on the isolated node). This gives
therefore the octic $C_8$ of Fig.\,d. Here it is useful to
introduce letters a,b,c,d,e,f in order to understand a bit the
highbrow distortion of such a map. Warning: it seems that actually
our picture created curves violating Gudkov, so that actually the
oval 2 of Fig.\,d should be inside the complicated circuit of
$C_8$. At any rate, it seems that this construction of Viro is the
most tricky as it uses as well a semi-regional (large) deformation
principle. It is incidentally quite puzzling to wonder if this
construction contributed to Viro's general formulation of the
Itenberg-Viro contraction principle.
%%%(\ref{Itenberg-Viro-contraction:conj}).
Of course in degree 5 it
could be that the latter conjecture is true essentially as in
degree 6 (Itenberg 1994) and by virtue of Nikulin-Kharlamov's
theory, but this is probably not really required.

\begin{figure}[h]\Figskip
%\vskip-1.2cm\penalty0
%\centering
\hskip-1.7cm\penalty0
\epsfig{figure=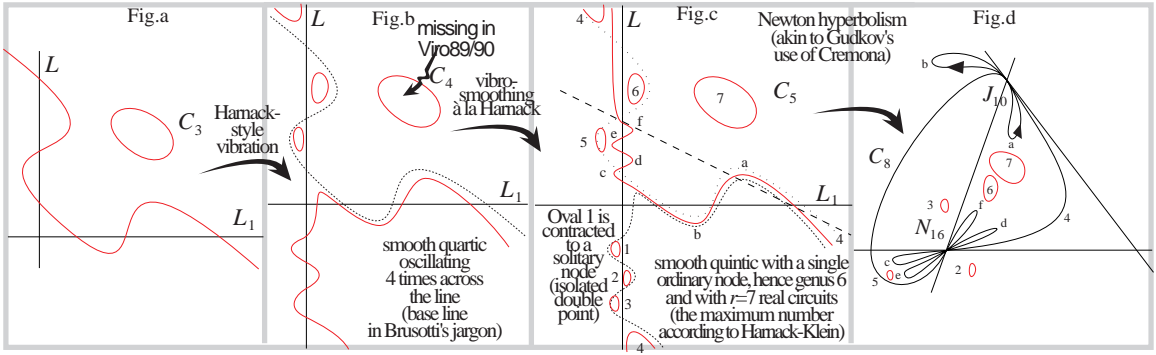,width=152mm}
\captionskipAG
  \caption{\label{ViroDEGREE8_SING_BIS:fig}%
  Viro's 2nd curve (from Viro89/90 plus our correction)}
\figskip
\end{figure}

[07.05.13] Along  these 3 constructions we expected to obtain  all
the 42 (new) $M$-octics of Viro 1980. Yet it seems that this is
not yet the case, because some few schemes assigned the letter
``V'' of Viro 1980's announcement are not yet constructed in our
text based on Viro 89/90 (compare on our
Fig.\,\ref{Degree8-M-curve-TABLE:fig} the V-symbols not yet
covered by squares, rhombs or parallelograms). Probably the other
V-schemes are obtainable by variants of the 2 tricky methods just
exposed. Maybe one can even drift to an art-form of freehand
drawing of such singular octics with petals at two points. Notice
that besides the petaliform circuit both singular $C_8$ used
%%%yet
above have 5 extra ovals whose location can be derived
%%%ad-hoc-ally
a posteriori from Gudkov's congruence.

More precisely Fig.\,\ref{ViroDEGREE8_SING_TRIS:fig}a,b recalls
the 2 ground curves of Viro. One can drag a petal inside to get
two bipetals as on Fig.\,c. However after smoothing we get
Fig.\,d,
%%which violates however
violating B\'ezout (trace the line through the 2 nests of depth
2). Of course this is not a contradiction against Viro method,
B\'ezout being already violated on our liberal singular octic of
Fig.\,c. Let us instead drag the bi-petal of Fig.\,a outside to
get Fig.\,e. Alas this produces only schemes that were already all
obtained by the quadri-ellipse $C_8$ of Viro. Of course the
construction of dragging can be much varied. For instance
Figs.\,g,h produce the new scheme $6\frac{15}{1}$. Of course
Fig.\,g is pure freehand drawing and so our patchwork is  a bit
sloppy.

At this stage the game is to reach the scheme $2\frac{19}{1}$ with
only two outer ovals (compare the diagrammatic of
Fig.\,\ref{Degree8-M-curve-TABLE:fig}). Some few thinking brings
us to Fig.\,\ref{ViroDEGREE8_SING_TRIS:fig}i, which produces
rigorously (without hand drawing) this and the former scheme
$6\frac{15}{1}$. The corresponding schemes are marked by
green-stars on the main-table
(Fig.\,\ref{Degree8-M-curve-TABLE:fig}). Especially interesting is
also the RKM-schemes $3\frac{16}{1}$ and  $7\frac{12}{1}$ which
affords another corruption of Rohlin's maximality conjecture.

\begin{figure}[h]\Figskip
%\vskip-1.2cm\penalty0
\centering
%\hskip-0.7cm\penalty0
\epsfig{figure=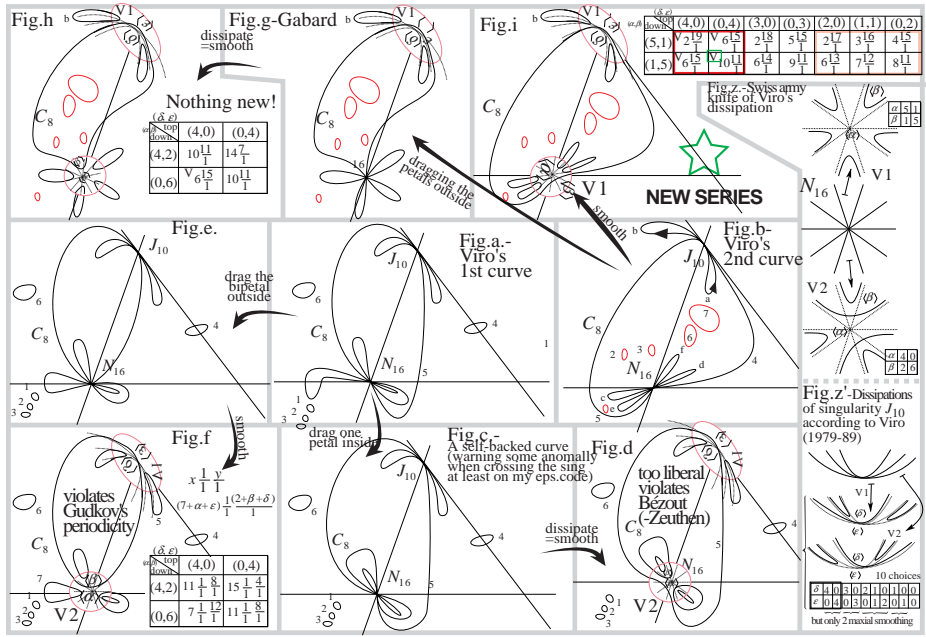,width=122mm}
\captionskipAG
  \caption{\label{ViroDEGREE8_SING_TRIS:fig}%
  Freehand tracing of fundamental octics}
\figskip
\end{figure}

At this moment we nearly have understood Viro's method that one
can realize curves with preassigned topology, e.g. those
Viro-types not yet realized on the table.

\subsection{Some messy ideas of Gabard}

The next idea that came to us is that since we just discovered
another series allied to  Viro's 2nd curve (cf.
Fig.\,\ref{ViroDEGREE8_SING_TRIS:fig}i) there must also be a 2nd
series allied to  Viro's 1st curve. However a priori it seems that
the series so obtained will not be extremely exciting and will
probably coincide with the boring curves of Fig.\,77 of Viro 89/90
(p.\,1129), that are already obtained by the more elementary
device of the quadri-ellipses. Let us however trace them carefully
to check our guess. As above the idea is to close the bi-petal by
the pair of ``paralleling'' braids, yet this time in such a
fashion that the resulting nest of depth 2 is not charged by extra
ovals. This is possible after symmetrizing
%the
Viro's gluing
%%%infinitesimal
 pattern, and leads  to
Fig.\,\ref{ViroDEGREE8_SING_QUATRIS:fig}b. This realizes 3 schemes
marked by little rhombs on the main-table
(Fig.\,\ref{Degree8-M-curve-TABLE:fig}), which (as expected) are
the extremely boring specimens $9(1,10 \frac{1}{1})$, $13(1,6
\frac{1}{1})$, $17(1,2 \frac{1}{1})$ already obtained by the
simpler device (of perturbation of 4 ellipses).

\begin{figure}[h]\Figskip
%\vskip-1.2cm\penalty0
\centering
%\hskip-0.7cm\penalty0
\epsfig{figure=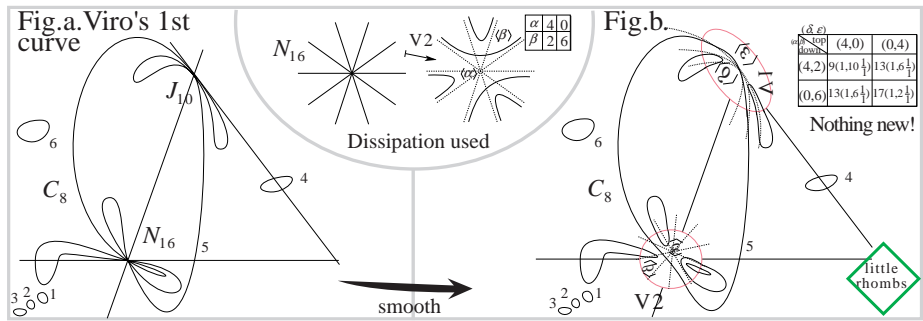,width=122mm}
\captionskipAG
  \caption{\label{ViroDEGREE8_SING_QUATRIS:fig}%
  A variant from Viro's 1st  curve}
\figskip
\end{figure}

At this stage one is slightly puzzled and it is not clear if the
obtention of the remaining Viro's $M$-schemes requires changing of
fundamental curve or smoothing more cleverly  the 2 singular
curves of Viro.

The next idea of our somewhat random search is to consider
Gabard's curve (so-called not because we are lacking in modesty
but because it is just sloppily hand-drawn) depicted on
Fig.\,\ref{ViroDEGREE8_SING_TRIS:fig}g and reproduced as
Fig.\,\ref{ViroDEGREE8_SING_QUINTIS:fig}a) while changing the mode
of smoothing to get Fig.\,b. Alas all those 4 schemes are
%boring
%being
already realized by Viro as smoothing of 4 ellipses. They
are reported by a septagonal star on the main-table
(Fig.\,\ref{Degree8-M-curve-TABLE:fig}). The most interesting
scheme is perhaps $5\frac{6}{1}\frac{9}{1}$ as it flirts nearly
with the $\chi=-16$ row (of the main-table) which is the most
mysterious one containing 4 among the Hilbert-Viro bosons (not yet
known to be realized nor to be prohibited). Those curves having
only one outer ovals, one is inclined to look at a variant of
Viro's 2nd curve where the loop b is dragged inside the singular
circuit (Fig.\,d). Of course doing so
%%%we notice that
the
unicursality of the singular circuit is lost (at least under the
postulate that branches of equal curvature are connected).

Then one produces the smoothing of Fig.\,e where the
$J_{10}$-singularity is smoothed symmetrically so as to create the
maximum number of ovals. However doing this and referring back to
Viro's dissipation list (cf. Fig.\,29, p.\,1103 in Viro 89/90)
%%% or alternatively our Fig.\,\ref{Viro3-15:fig}),
we note that
$\delta+\varepsilon$ cannot be as large as  4 as in the asymmetric
smoothing, but its maximum permissible value is only 3 realized by
the pair $(\delta, \varepsilon)=(3,0)$. Let us however on the
table of Fig.\,e also consider the value $(4,0)$ to look what
monster would result. Actually we obtain the $M$-scheme
$1\frac{20}{1}$. Now recall Petrovskii's estimate $\chi\ge
-\frac{3}{2}k(k-1)=-\frac{3}{2}4\cdot3+1=18=-18$, while our
pseudo-curve has $\chi=1+(1-20)$ and so just respects it. However
it is probably ruled out by the strengthened Petrovskii bound of
Arnold 1972 or ``more'' elementarily by the Arnold congruence mod
4, which can be regarded as a formal consequence of Rohlin's
formula. Of course Gudkov's hypothesis (mod 8) do as well the job,
but is more tricky to prove (Rohlin 1972).

Now again with the idea to attack by surprise the last mysterious
column with $\chi=-16$, let us trace freely a curve with 2
singular circuits by splitting of that of the 2nd Viro curve (cf.
Figs.\,c,f). (After all nobody ever asserted that real curves are
connected and real geometers (say Pl\"ucker-Zeuthen-Klein-Harnack,
etc.) nearly learned us the exact opposite.) Of course during the
process it seems reasonable to destroy one red oval. However on
smoothing the configuration as on Fig.\,g we will have at least 2
outer ovals and the case of 2 ovals is prohibited by Gudkov
hypothesis (or just Arnold). Now the idea would be to split
without loosing an oval while using the symmetrical dissipation of
$J_{10}$. This idea leads to Fig.\,h, whose smoothing Fig.\,i
creates the $M$-scheme $1\frac{7}{1}\frac{12}{1}$ which was never
constructed as yet. Of course since our method is pure free-hand
drawing this does not prove existence of the curve. Yet it is
interesting to vary the parameters to see which kinds of schemes
arise, and actually there is only one maximal companion namely
$5\frac{7}{1}\frac{8}{1}$. Since the latter was first constructed
by Viro's smoothing of coaxial ellipses, some ``principe du
raccord'' (yet another patchwork if you like) gives some very weak
evidence that the scheme $1\frac{7}{1}\frac{12}{1}$ exists
algebraically. Of course upon playing with the dissection of the
singular circuit of Viro's 2nd curve, while keeping in mind that
$\delta=3$ we see that the upper (non empty) oval can contain
either $5,6,7,8,9,10$ ovals. In fact it is convenient to denote by
$U$ the number of upper red ovals inside the $J_{10}$-circuit of
Fig.\,h, i.e. after cellular subdivision of Viro's 2nd curve. This
$U$ can range from $0,1,2,3,4,5$ (a priori), and Fig.\,j tabulates
the resulting schemes after smoothing. We obtain so for one outer
ovals 2 schemes denoted by V, already constructed by Viro80 (via
the most elementary method of 4 ellipses), and for $U=2,4$ two new
schemes denoted by Ga (not yet known to exist), and one scheme
prohibited by Orevkov 2002! Further by choosing instead the lower
dissipation with $(\delta,\varepsilon)=(1,5)$ we get schemes all
realized by the elementary Viro method, so that the naive
principle of propagation could imply that all the former schemes
(1st row) also exists. Of course this would contradict Orevkov's
theorem, and actually the latter can be interpreted as an
obstruction to split Viro's 2nd curve (at least the variant with
introverted the b-loop).

\begin{figure}[h]\Figskip
%\vskip-1.2cm\penalty0
\centering
%\hskip-0.7cm\penalty0
\epsfig{figure=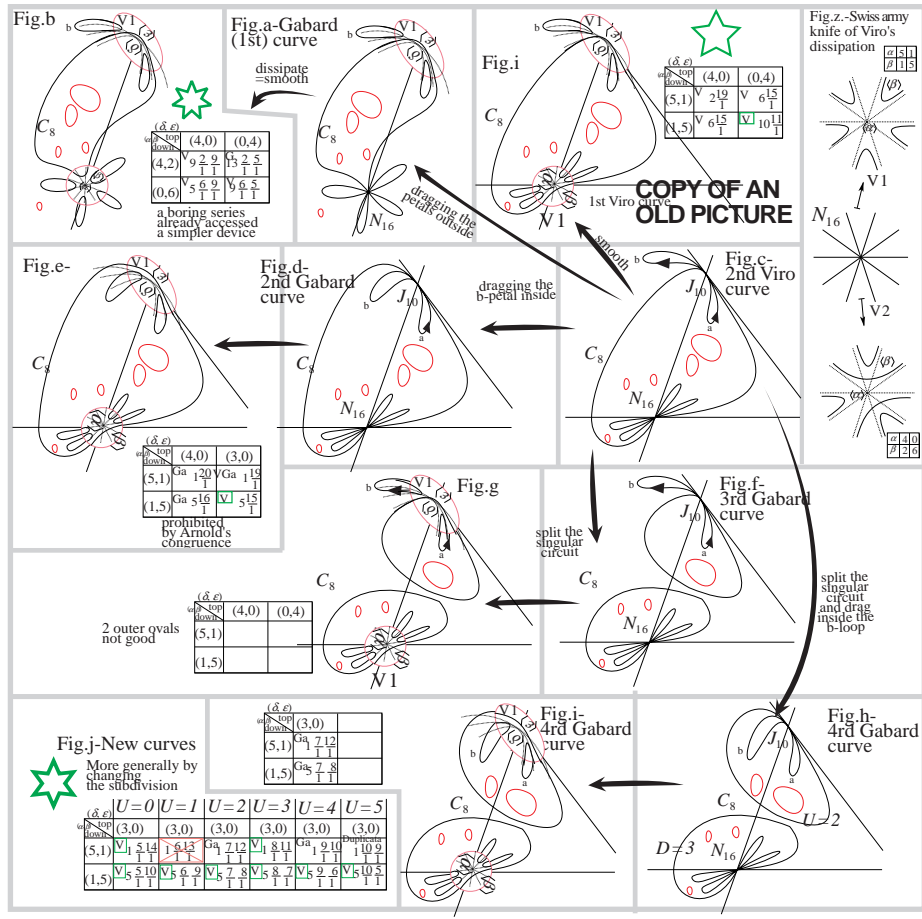,width=122mm}
\captionskipAG
  \caption{\label{ViroDEGREE8_SING_QUINTIS:fig}%
  A variant from  Viro's second curve}
\figskip
\end{figure}

Of course all this is very speculative, and we need to return too
a more pragmatical standpoint.

[05.06.13] On waking this morning (with short hairs), we were
flashed by the idea of looking at the configuration of 3 coaxial
ellipses plus a transverse one. Alas the curves so obtained are
far from maximal and seem to reach at most $(M-6)$ ovals. This is
a big deception.

\begin{figure}[h]\Figskip
%\vskip-1.2cm\penalty0
\centering
%\hskip-0.7cm\penalty0
\epsfig{figure=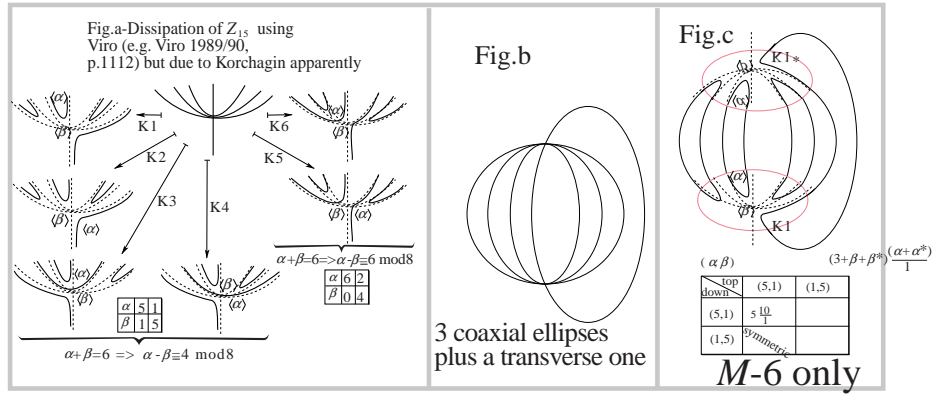,width=122mm} \captionskipAG
  \caption{\label{ViroDEGREE8_GABARD:fig}%
  A variant of  Viro's method due to Gabard
  %(but extremely disappointing)
  (but fails blatantly)} \figskip
\end{figure}

\begin{Scholium} In fact we realize now that Arnold's prose about
the distribution of ovals (in his seminal 1971 paper) might have
been inspired by the allied jargon concerning the distribution of
primes in basic number theory. Probably, the analogy is far
reaching in the sense that up to present knowledge the
distributions realized appears as fairly random, without clear-cut
rule governing the architecture of higher Gudkov pyramids
parametrizing the periodic table of elements (schemes in Rohlin's
terminology).
Notwithstanding, it is only a matter of time to examine deeper the
crystal as to uniformize all Viro, Shustin, and Orevkov
prohibitions while subsuming to one and a sole paradigm, viz.
total reality or perhaps the allied method of deepest penetration
boiling down to B\'ezout for higher order curves.
\end{Scholium}

\subsection{Shustin's constructions}

$\bullet$ Shustin 1985
\cite{Shustin_1985-Indep-removal-and-new-M-curves} (announcement)
and details in 1987--88
(\cite{Shustin_1988-new-M-and-M-1-curves-of-deg-8},
\cite{Shustin_1987/87-a-new-M-curve-of-deg-8}) new constructions
of $6+1=7$ schemes (probably via a variant of Gudkov or Viro) [of
course Viro seems more likely]. In fact, it seems that Shustin's
original proof was somewhat independent of Viro's method, compare
p.\,488 of Shustin 1988
\cite{Shustin_1988-new-M-and-M-1-curves-of-deg-8} where we read:
``The existence of curves of degree 8 with schemes (1)[=the list
of six] was announced by the author in [6](=Shustin 1985
\cite{Shustin_1985-Indep-removal-and-new-M-curves}), where it was
deduced from results on investigation of smoothing of point of
quadratic contact of four non-singular real branches\footnote{So
this seems to be $X_{21}$, yet it looks hard to get all those
schemes via dissipation of the quadri-ellipse.}. Here we give
another proof that was obtained by using Viro's method of gluing
real algebraic curves [2](=Viro 1983
\cite{Viro_1982/84-Gluing-of-plane}).'']. Slightly later, Shustin
found the scheme
$4\frac{5}{3}=4\frac{5}{1}\frac{5}{1}\frac{5}{1}$.
Fig.\,\ref{Degree8-M-curve-TABLE:fig} below shows the exact list
of 7 schemes realized by Shustin, denoted by the letter ``S'',
where the last found is denoted ``S=last''. This is probably the
construction alluded to Orevkov's letter
%%(in Sec.\,\ref{e-mail-Viro:sec}),
(in v.2 of Gabard 2013), and the idea is probably to use octics
with 3 singularities instead of the two used by Viro.

[08.05.13] Now we present the details following Shustin 1988
\cite{Shustin_1988-new-M-and-M-1-curves-of-deg-8}. Again it
suffices to have singular octics while applying the dissipation
method. This time Shustin considers curves with $Z_{15}$
singularities (i.e. 3 branches with contact of order 2 and a
fourth branch transverse to it). The dissipation of this
singularity were apparently classified by Korchagin 1988
\cite{Korchagin_1988-septics-with-one-Z15}, the maximal ones being
depicted on Fig.\,\ref{ViroDEGREE8_SHUSTIN:fig}a. Actually we
shall only employ the smoothings K5, K6. Next, Shustin traces 3
curves but actually the first one already leads to the
construction of the 6 schemes announced by Shustin namely those
marked by the letter ``S'' and an octagonal star on the main-table
(Fig.\,\ref{Degree8-M-curve-TABLE:fig}). Additionally, Shustin's
construction recovers 4 cases claimed by Viro 1980 (namely those
depicted by V? on Fig.\,\ref{ViroDEGREE8_SHUSTIN:fig}b), but which
we were  personally not able to construct (upon reading Viro's
text of 1989/90). Those mixed Viro-Shustin schemes are marked by
the combined symbol ``V~(or\,S)'' plus
%%% a smaller
an octagonal star on the main-table. As a moral, Shustin's method
affords many new types that were inaccessible before.
%fairly delicate.
Interestingly, Shustin propose 2 other ground curves also doted
%%%with
of two $Z_{15}$ singularities, whose smoothings may be worth
exploring, yet it seems that they are not formally required as the
first curve suffices to exhibit all (six) $M$-schemes claimed by
Shustin.

\begin{figure}[h]\Figskip
%\vskip-1.2cm\penalty0
%\centering
\hskip-2.7cm\penalty0
\epsfig{figure=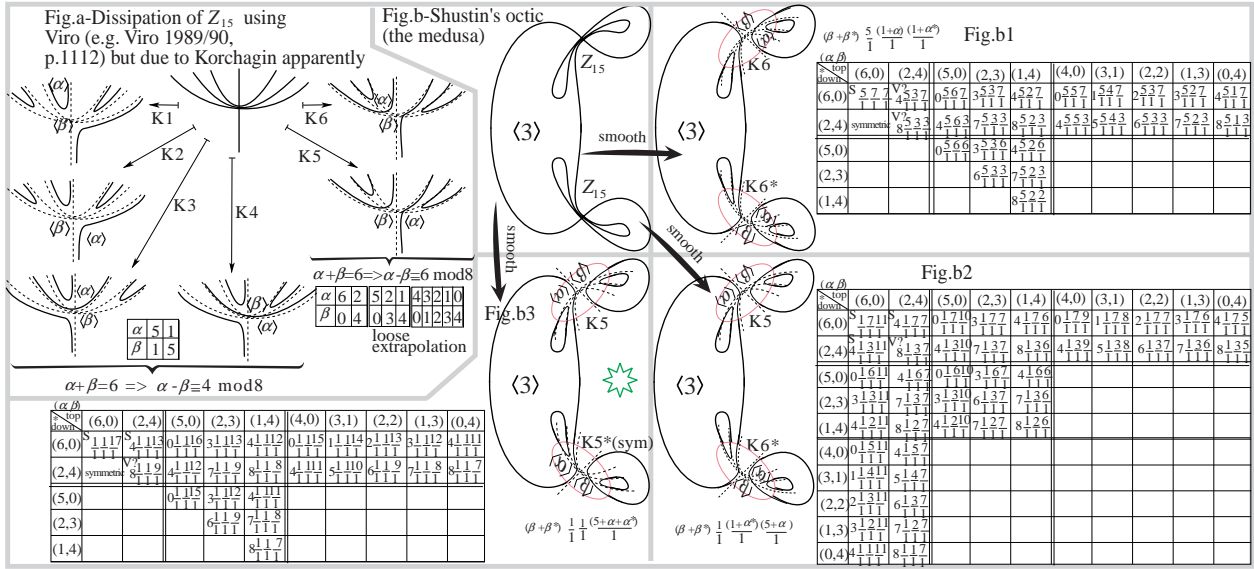,width=172mm} \captionskipAG
  \caption{\label{ViroDEGREE8_SHUSTIN:fig}%
Shustin's series of 6 new $M$-schemes via $Z_{15}$ (twice) on the
  medusa}
\figskip
\end{figure}

At this stage we see that Viro obtained first many $M$-octics by
dissipating $X_{21}$ (quadruple bicontact), then some few others
by smoothing $N_{16}+J_{10}$ (quintuple point\,+\,triple
bicontact), and finally Shustin added to the list several new
schemes by dissipating $Z_{15}$ (triple bicontact with a
transverse branch, alias candelabrum). So philosophically, we see
that Viro's method enjoys two levels of freedom: the choice of the
singularities and the global singular curve which is smoothed.

Further Shustin 1988 (p.\,490--92 of \loccit) gives a detailed
construction of the above singular $C_8$ along the method of
hyperbolism \`a la Huyghens-Newton-Gudkov-Viro and himself. We
shall detail this at the occasion.

It seems however more urgent to inspect what results from the 2nd
(and 3rd) curve of Shustin. His second curve F2 (p.\,489) looks an
apple alike (Fig.\,\ref{ViroDEGREE8_SHUSTIN2:fig}b). As usual the
algorithm is to self-connect the loops with themselves in order to
maximize the number of ovals, selecting appropriately the
dissipation. On the case at hand if we imagine the
$Z_{15}$-singularity as a tree with a (vertical) trunk and 3
branches growing transversally from it, we choose the smoothing
$K1^\ast$ (i.e. K1 symmetrized on the top). On the bottom it seems
harder to find a closing-gluing and actually the one ideally
suited achieves only $\al+\be =5$, cf. again Viro's Fig.\,39 in
Viro 89/90 (p.\,1112) or our Fig.\,a (K7). Hence tolerating a
maximal smoothing ($\al+\be =6$) yet not closing perfectly all the
loops gives Fig.\,c. So we choose for instance $K5^\ast$ where the
star is the symmetrized dissipation of K5. Alternatively, we may
choose a non-maximal dissipation which closes the loops, e.g.
$K7^\ast$ at the bottom, but experimenting a bit (or reading
better Shustin's text especially p.\,490) one sees that this 2nd
curve leads only to $(M-1)$- or even $(M-2)$-curves. So it is not
most exciting for our present purpose of cataloging $M$-curves.

\begin{figure}[h]\Figskip
%\vskip-1.2cm\penalty0
\centering
%\hskip-0.7cm\penalty0
\epsfig{figure=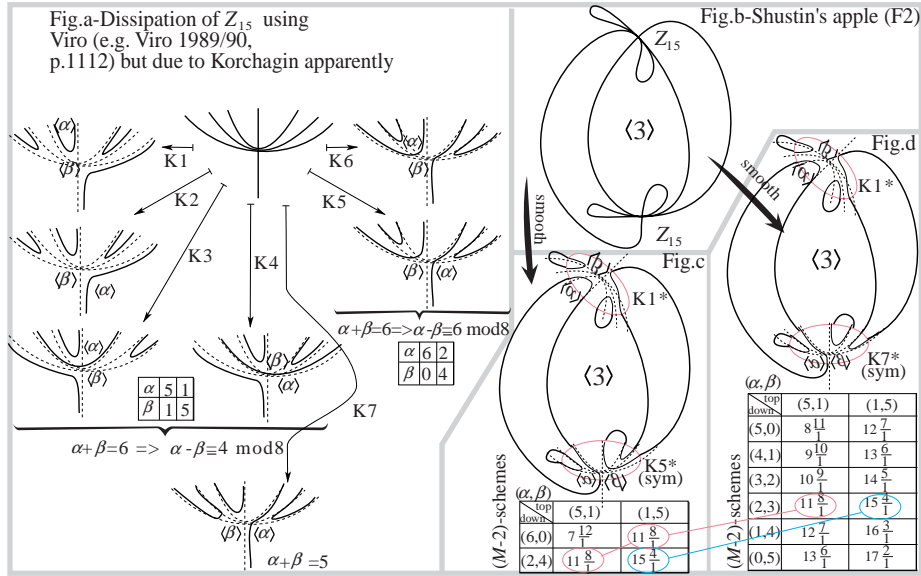,width=122mm}
\captionskipAG
  \caption{\label{ViroDEGREE8_SHUSTIN2:fig}%
  The 2nd Shustin's series}
\figskip
\end{figure}

{\it Added } [02.06.13].---However the $(M-2)$-scheme
$7\frac{12}{1}$ (or its companions $11\frac{8}{1}$,
$15\frac{4}{1}$ obtained on Fig.\,c) is of utmost interest in
relation with Rohlin's maximality conjecture. Indeed the scheme in
question is RKM hence of type~I, and so should kill all its
enlargements. The geography of those extensions are depited by
4-branched asters on Fig.\,\ref{Degree8-(M-i)-curve-TABLE:fig},
and includes for instance the $(M-1)$-scheme
$5\frac{2}{1}\frac{12}{1}=:S_{M-1}$. As the latter stands below
Viro's $M$-scheme $5\frac{2}{1}\frac{13}{1}$ or even
$5\frac{3}{1}\frac{12}{1}$, it may be inferred from
Itenberg-Viro's contraction principle that the scheme $S_{M-1}$ is
very likely to exist. However this would conflict with Rohlin's
maximality conjecture (RMC). Hence we reach a paradox, which can
be solved either by a falsity of Shustin's construction, or of the
contraction principle or finally a disruption of RMC. Actually,
the contraction principle can even be left aside of the token,
just by enlarging directly the $(M-2)$-scheme to the $M$-scheme
constructed by Viro.

Exactly the same comments apply to the scheme $11\frac{8}{1}$ or
$15\frac{4}{1}$ which are also created by Shustin's construction.
Additionally from Fig.\ref{ViroDEGREE8_SHUSTIN2:fig}d, we may
construct more schemes (namely all those of the form
$(8+a)\frac{(11-a)}{1}$, with $0\le a\le 9$), which are not
necessarily RKM, yet we gain no more RKM schemes. To accentuate
the paradox it would be nice to construct the $(M-1)$-schemes
extending the RKM-schemes by hand without reference to the
(nebulous) Itenberg-Viro contraction principle. However even that
is not an absolute prerequisite because it is actually sufficient
to look directly at the $M$-schemes extending out RKM-schemes. So
we get a direct conflict between RMC and Viro's method. For
instance the RKM-scheme $15\frac{4}{1}$ enlarges to
$13\frac{2}{1}\frac{4}{1}$ which in turn enlarges to the
$M$-scheme $13\frac{2}{1}\frac{5}{1}$ constructed by Gudkov or
Viro. So it seems at this stage that there is a clear-cut
corruption of Rohlin's maximality conjecture. We resume the
situation with the following result.

\begin{theorem}\label{RMC:cter-example-via-Shustin}
Shustin's apple construction ({\rm
Fig.\,\ref{ViroDEGREE8_SHUSTIN2:fig}c})
%%%%%produces counter\-examples to
refutes Rohlin's maximality conjecture that a scheme of type~I is
forced to maximality. More precisely any one of the three
$(M-2)$-schemes of degree eight $15\frac{4}{1}$, $11\frac{8}{1}$
or $7\frac{12}{1}$ satisfying the RKM-congruence ($\chi\equiv_8
k^2+4$) are realized algebraically, yet
%%%can be enlarged
enlargeable in the algebraic category via $M$-schemes
constructible \`a la Viro by  dissipating the quadri-ellipse (cf.
{\rm Fig.\,\ref{ViroDEGREE8:fig}d}). Those are respectively for
instance (compare\vadjust{\vskip1pt} {\rm
Fig.\,\ref{Degree8-(M-i)-curve-TABLE:fig}}),
$13\frac{2}{1}\frac{5}{1}$ (constructed by G=Gudkov or V=Viro), or
$9\frac{2}{1}\frac{9}{1}$ (due to Viro), or finally
$5\frac{2}{1}\frac{13}{1}$ (due to Viro).
\end{theorem}

\subsection{Shustin's last construction: $4\frac{5}{1}\frac{5}{1}\frac{5}{1}$}

%A more
Another interesting task is to understand Shustin's last
construction of the $M$-scheme
$4\frac{5}{1}\frac{5}{1}\frac{5}{1}$, which was realized by
Shustin's modification of Viro's method. This is published as
Shustin 1987/87 \cite{Shustin_1987/87-a-new-M-curve-of-deg-8}, but
alas not enough pictures are supplied there. We hope to remedy
this at the occasion.

In fact inspired by Shustin's text we traced the following figure
(Fig.\,\ref{ViroDEGREE8_SHUSTIN_NEW:fig}), which however corrupts
violently Harnack's bound.

[09.05.13] At this stage we hoped to find a more geometric
treatment of Shustin's last curve, in Polotovskii 1988
\cite{Polotovskii_1988---classif-deg-8} but alas not so.

\begin{figure}[h]\Figskip
%\vskip-1.2cm\penalty0
%\centering
\hskip-2.7cm\penalty0
\epsfig{figure=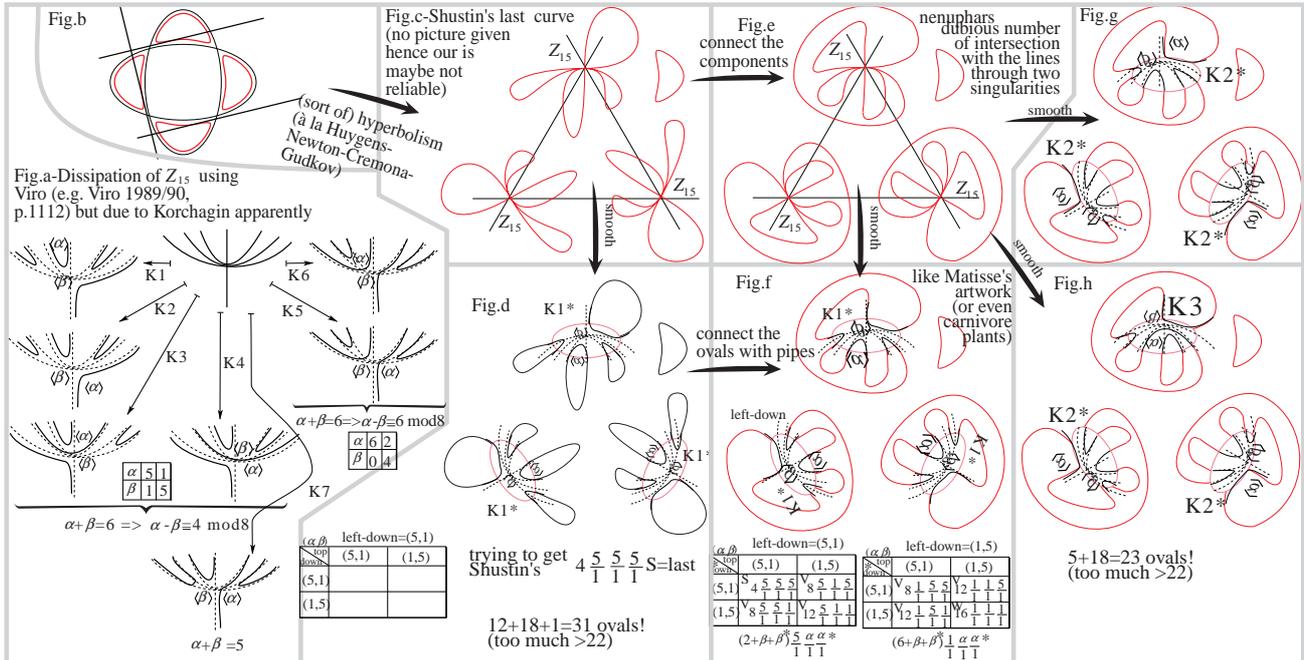,width=172mm}
\captionskipAG
  \caption{\label{ViroDEGREE8_SHUSTIN_NEW:fig}%
Trying to trace Shustin's last (seventh) curve
$4\frac{5}{1}\frac{5}{1}\frac{5}{1}$, but mis-depicted by Gabard
(ohne Gew\"ahr and a lot of contradictions).} \figskip
\end{figure}

[02.07.13] Nearly 2 months later (with interruptions due mother's
health problem) we had the idea that since the curve of Fig.\,d
violates strongly Harnack's bound one should reduce the number of
components by piping together its ovals. This idea suggested us to
trace Fig.\,f, and also Fig.\,e by anticipating the piping atthe
level of the singular curve. Working out the dissipation we find
indeed the curve asserted by Shustin
$4\frac{5}{1}\frac{5}{1}\frac{5}{1}$ as well as the curves
$8\frac{1}{1}\frac{5}{1}\frac{5}{1}$ and
$12\frac{1}{1}\frac{1}{1}\frac{5}{1}$ (both first constructed by
Viro 1980), and finally recover
$16\frac{1}{1}\frac{1}{1}\frac{1}{1}$ (first constructed by Anders
Wiman). All this is excellent (even if somewhat heuristic piping),
if one had not remora with the issue that the line through any 2
of the 3 singularities of the octic $C_8$ of Fig.\,e seems to
intersect the $C_8$ in more than 8 points. Indeed each singularity
being a quadruple point with four branches meeting like the
candelabrum, we see that each singular points contributes already
for 4 intersections so that no extra intersection can occur (say
as on our depiction). Hence the latter must be suitably correct in
order to hope that our interpretation of Shustin's construction is
a tangible one (the correct one).

Despite this little paradox we are presently not able to explain,
we can somewhat cavalier explore the other smoothings of Shustin's
curve (the nenuphar of Fig.\,e) using dissipation K2
(symmetrized). Yet on tracing Fig.\,g one sees (disappointingly)
that the new curve so generated is isotopic to the former one
(Fig.\,f), and since the structural constants $\al, \be$ of the
deformation are the same (Fig.\,a) we convince that no new curve
will emerge from Shustin's nenuphar (nymphae). Of course we could
also combine K1$^\ast$ with K2$^\ast$ dissipations (on different
nenuphars of the curve), yet the same token should kill any hope
to get something new along the way.

Of course one should still analyze other dissipations (like K3,
K4, etc. of Fig.\,a) yet those will not be closing, and so
certainly fail producing $M$-curves. However as already often
illustrated those quasi-maximal curves are still of
%uttermost
interest to appreciate the global architecture of the pyramid.
Alas, we do not feel much motivated doing this work as we are not
sure that our model of Shustin's curve is the correct one. However
on doing it we get Fig.\,h, by choosing just one K3 smoothing on
the top singularity (although we could have prescribed thrice K3
on all three singularities by the independence principle \`a la
Brusotti-Viro-Shustin). However to our little surprise the number
of ovals increases then (due to an unexpected closing)  and so
passes beyond Harnack's bound. A contradiction in mathematics is
obtained! Aber Hallo! Indeed, for the K3-dissipation $\al+\be =6$
too, so we have $3\cdot 6=18$ micro-ovals, plus the $3+2$ traced
on Fig.\,h yielding a total of $23>22$ ovals. So it seems that
Shustin's nymphean curve (Fig.\,e) does not exist or that
(Viro/Korchagin's) theory of dissipation of the candelabrum
$Z_{15}$ is foiled. Of course the most plausible explanation is
that our
%%%interpretation
drawing of Shustin's curve is not the correct one.

[03.07.13] But then what is the correct way to trace Shustin's
singular octic? Of course there exists other ways to pipe together
the ovals of Fig.\,d, and so we have for instance
Fig.\,\ref{ViroDEGREE8_SHUSTIN_NEW2:fig}i and the allied Fig.\,j.
The latter has however $18+6+1=25$ ovals (violating Harnack's
bound by 3 units). Another piping of the ovals (or rather
circuits) is given on Fig.\,k, yet its smoothing Fig.\,l also
violates Harnack. So we need a more radical connection among the
circuits suggesting, e.g., to trace Fig.\,m. Here it seems that
the varied possible smoothings all respect Harnack's bound (as
experimented by Figs.\,n,o,p,q,r). Alas the ground curve (Fig.\,m)
corrupts B\'ezout for line (as the traced line shows 10
interceptions with the presupposed $C_8$). So our our curve
Fig.\,m is still not a {\it bona fide\/} model of Shustin's curve.
Our Fig.\,s still respects Harnack (as $4+18=22$), yet choosing
$(\al,\be)=(6,0)$ on all three singularities gives the scheme
$\frac{21}{1}$ which is prohibited by Gudkov's hypothesis.
(Alternatively it is prohibited by Petrovskii's inequality
(\ref{Petrovskii's-inequalities:thm}) we reads here as
$-18=-\frac{3}{2}\cdot 4(4-1)=-\frac{3}{2}k(k-1)\le\chi$.) So
there is a structural obstruction to the existence of the curve of
Fig.\,m.

\begin{figure}[h]\Figskip
%\vskip-1.2cm\penalty0
%\centering
\hskip-2.7cm\penalty0
\epsfig{figure=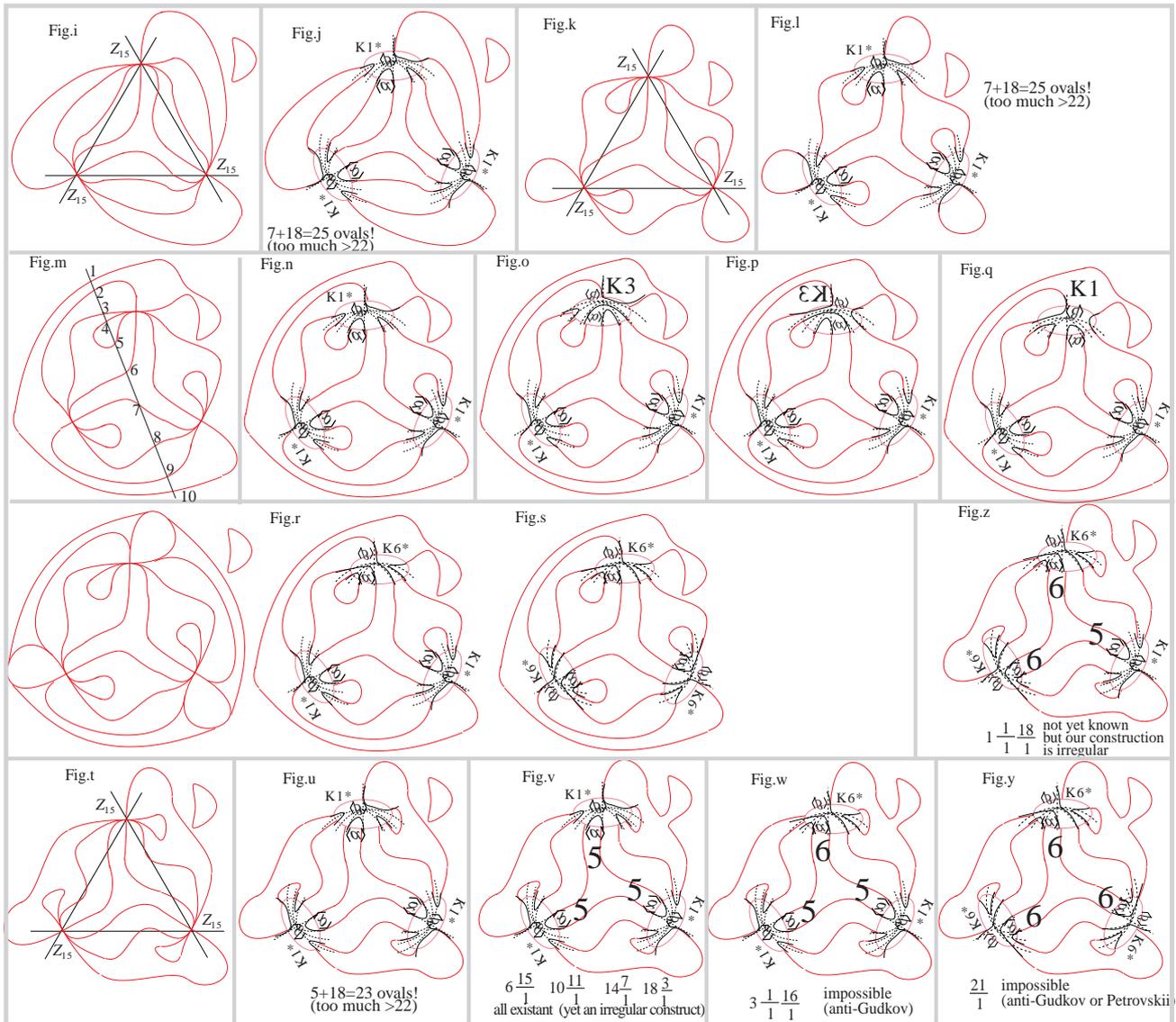,width=172mm}
\captionskipAG
  \caption{\label{ViroDEGREE8_SHUSTIN_NEW2:fig}%
Further attempts to trace Shustin's last (seventh) curve
$4\frac{5}{1}\frac{5}{1}\frac{5}{1}$, but still mis-depicted by
Gabard (ohne Gew\"ahr and a lot of contradictions).} \figskip
\end{figure}

Actually it is clear that our method is very lazy (i.e. purely
topological without any algebro-geometric substance). For instance
we may consider Fig.\,t which smoothed as Fig.\,u produces too
many ovals, so we proceed to the Verschmelzung of Fig.\,v.
Choosing $(\al,\be)=(5,1)$ gives the scheme $6\frac{15}{1}$, while
taking $(\al,\be)=(1,5)$ produces $18\frac{3}{1}$. Changing one
smoothing to K6 (starred=symmetrized) gives Fig.\,w where we
choose the top $\al$ as 6 and the bottom $\al$'s as $5$ (and
always $\al+\be =6$) we get the scheme $3\frac{1}{1}\frac{16}{1}$
which violates Gudkov periodicity (yet not killed by Arnold's
weaker congruence mod 4). So the Fig.\,t with the conjunction of
Fig.\,v is not a viable model of Shustin's curve. Since the above
scheme falls fairly  close to the mystery-scheme
$1\frac{1}{1}\frac{18}{1}$ (not yet known to exist or not) it
seems tempting to reiter a nearby smoothing of Fig.\,y which gives
the impossible scheme $\frac{21}{1}$ when $\al=6$ throughout.
Fig.\,z instead gives as corresponding scheme (for the depicted
distributions of $\al$'s) the following one
$1\frac{1}{1}\frac{18}{1}$, which is precisely the one not yet
known. So we have nearly proved something new but alas not so as
our ground curve (Fig.\,u/v) is constructed by an irregular
(purely topographical) device without control upon the
algebraicity of the picture.

So it is evident that our method is purely heuristic and as yet
not extremely successful. However we cannot exclude that clever
twists of it (experiments) may lead to some new insights (at least
by supplying a topological candidate for an algebraic singular
octic that could produce new schemes by dissipation).

Below Fig.\,\ref{ViroDEGREE8_SHUSTIN_NEW3:fig}a is another failing
attempt to surger the basic curve
Fig.\,\ref{ViroDEGREE8_SHUSTIN_NEW:fig}c. Of course during the
process we noted that this basic curve itself violates B\'ezout
(intersect with the fundamental lines). So we need first to fix
this issue and this ay be arranged by rotating the petals as to
avoid intersection with the 3 ``coordinate-axes''. This gives us
Fig.\,b. Upon surgerying we get Fig.\,c whose smoothing (Fig.\,d)
overwhelms Harnack. This brings us to the idea of connecting the 6
lunes in pair as to gain Fig.\,e, which looks promising. Indeed
its smoothing Fig.\,f yields the desired scheme
$4\frac{5}{1}\frac{5}{1}\frac{5}{1}$ of Shustin, but alas another
smoothing (Fig.\,g) violates Axel Harnack.

\begin{figure}[h]\Figskip
%\vskip-1.2cm\penalty0
%\centering
\hskip-2.7cm\penalty0
\epsfig{figure=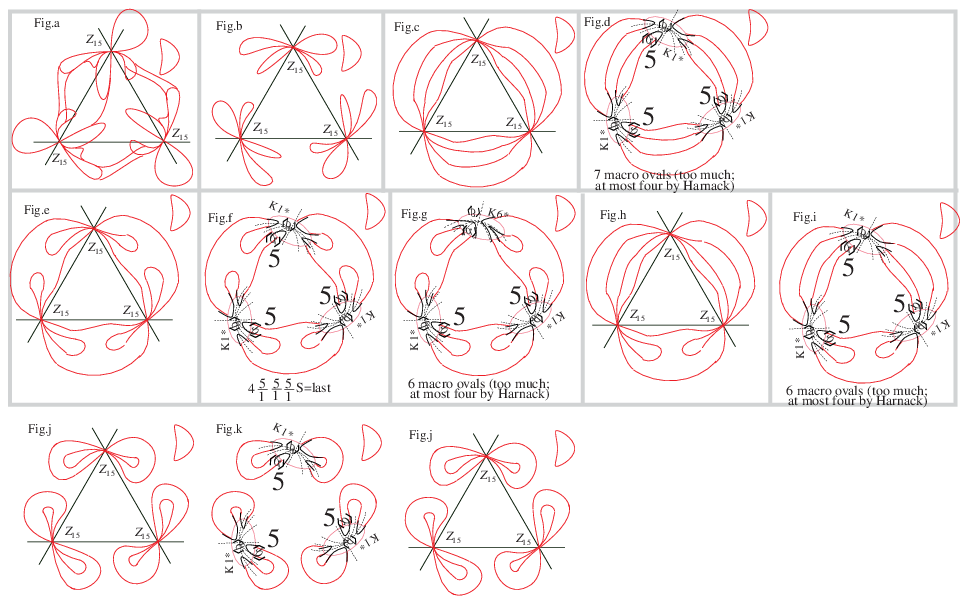,width=172mm}
\captionskipAG
  \caption{\label{ViroDEGREE8_SHUSTIN_NEW3:fig}%
Further attempts to trace Shustin's last (seventh) curve
$4\frac{5}{1}\frac{5}{1}\frac{5}{1}$, but still mis-depicted by
Gabard (ohne Gew\"ahr and a lot of contradictions).} \figskip
\end{figure}

At this moment the situation looks a bit desperate. Of course, we
can connect the 2 extra ovals of Fig.\,g yielding then Fig.\,h,
whose smoothing as Fig.\,i creates again six macro-ovals (too much
for Harnack). The desperation is now complete.

Our question is still how to trace a singular octic whose
dissipation leads to Shustin's scheme
$4\frac{5}{1}\frac{5}{1}\frac{5}{1}$. It seems to us a pity that
Shustin does not supply a picture of this curve and so we  are
relegate to a tedious guessing game. Of course the latter may be
boring yet it could also offer new insights on the cases not yet
settled. However we see that there is rather stringent obstruction
to the manufacture of the divine $C_8$ of Shustin. In a state of
quasi-somnolence (due to the high level of psychological complexes
in front of Shustin's intelligence) we discovered Fig.\,j which
albeit not intrinsically appealing (at least to my intuition)
flashed our attention since we can remember that in Shustin's
vague allusion is made to the union of 3 figures 8. Alas, the
smoothing of Fig.\,k is over-productive yielding 7 macro-ovals (3
more than the four permissible by Axel Harnack).

\subsection{Digressing on Orevkov's hypothetic curve}

[05.07.13] Lacking imagination, we started a random reading of
literature and found some inspiring idea in Orevkov 1999
\cite{Orevkov_1999-Link-theory}, p.\,782, Fig.\,4. There,
Orevkov's Theorem~1.3 states that there is no curve of degree 8 as
shown on Fig.\,\ref{ViroDEGREE8_OREVKOV_HYP:fig}a with
$\al+\be=11$. Under such circumstance it seems of interest to look
at the plethora of curves that would have resulted from
dissipating singularities \`a la Viro. So for instance we get
Fig.\,c but as there is already 22 ovals coming from the brackets
$\la a \ra$, $\la b \ra$ with $a+b=11$ we get a contradiction in
mathematics. (Note that we relabelled Orevkov's $\al,\be $ as
$a,b$ as to avoid a confusion with Viro's parameters for
micro-ovals $\al, \be$ as on Fig.\,b.) So it seems completely
obvious that Orevkov's curve (Fig.\,a) cannot exist and we do not
really understand the interest of his statement (Theorem 1.3 on
p.\,782).

\begin{figure}[h]\Figskip
%\vskip-1.2cm\penalty0
%\centering
\hskip-2.7cm\penalty0
\epsfig{figure=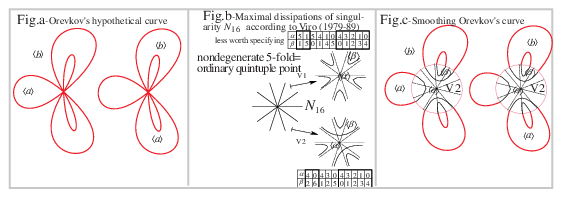,width=172mm}
\captionskipAG
  \caption{\label{ViroDEGREE8_OREVKOV_HYP:fig}%
Further attempts to trace Shustin's last (seventh) curve
$4\frac{5}{1}\frac{5}{1}\frac{5}{1}$, but still mis-depicted by
Gabard (ohne Gew\"ahr and a lot of contradictions).} \figskip
\end{figure}

\section{Decomposing curves:
%hypothetic
loose constructions via freehand drawings}

\subsection{A messy idea \`a la Wiman and decomposing octics}

[06.07.13] In view of the construction of Harnack, Hilbert, Wiman,
etc. it seems that a realist method of construction consist to
split the degree of interest in two and smooth a union of two
curves realizing the given degree. This leads to an art-form
well-known in Russia especially by Polotovskii, Shustin,
Korchagin, Orevkov. For $m=8$ we have the partition $1+7$, $2+6$,
$3+5$ and $4+4$. For $8=4+4$ we have only one $M$-quartic and
effecting a vibration gives a pair of quartics with one oval
maximally intersecting the other along $4\cdot 4=16$ points while
the other 4 ovals are just disjoint replicas. On smoothing $C_4
\cup C_4$ \`a la Brusotti or otherwise (Wiman does not cite
Brusotti) we get Wiman's  $M$-octic with scheme $16
\frac{1}{1}\frac{1}{1}\frac{1}{1}$.

Similar games must be possible with the other partitions.

\subsection{Degree 3+5}

Let us examine $3+5$ first. Here we only one isotopy class of
$M$-quintic (with $r=7$) resp. $M$-cubic (with $r=2$). Let us
assume that a pair of ovals is maximally intersecting along
$3\cdot 5=15$ points. Warning this is a misconception since both
pseudolines of the $C_3$ and $C_5$ have to intersect, hence their
intersection $C_3\cap C_5$ cannot by monopolized by an oval. So
assume rather that both pseudolines are maximally intersecting
along 15 points, but each pair of ovals chosen one from each curve
is disjoint. Smoothing the union $C_3\cup C_5$ gives a curve with
$15+1+6=22$ ovals hence an $M$-octic. Of course knowing its exact
scheme requires knowing more  on how the unique oval of $C_3$
surrounds the $6$ ovals of the $M$-quintic $C_5$. However since
the $C_5$ cannot be nested (unless it is the non-maximal deep
nest) we can infer a priori that the resulting octic scheme will
be simply nested, i.e. of the shape $x \frac{y}{1}$, hence not so
interesting as all those schemes are already realized by Viro's
method (look at the top row of
Fig.\,\ref{Degree8-(M-i)-curve-TABLE:fig}, zoomed as
Fig.\,\ref{Degree8-(M-i)-curve-TABLE_I:fig}).

\subsection{Degree 2+6}

Let us next examine the partition $8=2+6$. Here we imagine one
ovals of the $C_6$ maximally intersecting the conic $C_2$ along 12
points. So one should try to analyze all types of such decomposing
curves. This sort of problems is well-known to experts like
Polotovskii and Orevkov and one sees some direct interconnection
between the isotopic classification of decomposing curves (under
the natural assumption of transversality) and the pure isotopic
classification of a single curve of degree equal to the sum. To
get started let us fix the $C_6$ as being of Harnack's type
$9\frac{1}{1}$. A priori we can imagine that one oval oscillates
across the ellipse $C_2$ as on
Fig.\,\ref{ViroDEGREE8_GABARD_2+6_BOSON:fig}? below. If the
vibrating oval is an outer oval  the resulting scheme is
$(12+8)\frac{1}{1}=20\frac{1}{1}$ which violates Gudkov
periodicity mod 8. If the vibrating oval is the non-empty oval
then we get the (unnested) scheme $9+12+1=22$, which cannot exist
(e.g. by Petrovskii, or Rohlin's formula $(0=)2
(\pi-\eta)=r-k^2$). Finally if the vibrating oval is the unique
empty nested oval we get the scheme $9\frac{12}{1}$ which is also
anti-Gudkov (hence cannot exist).

Okay, but in reality we still have Hilbert's constructions
yielding decomposing curve of ``bidegree'' $(6,2)$ and producing
the interesting (but classical) schemes $1(1,2\frac{17}{1})$, etc.
as depicted on Fig.\,\ref{ViroDEGREE8_GABARD_2+6_BOSON:fig}. This
is fairly exhaustive (i.e. mixing all possible internal versus
external oscillation) yet this still misses the classical scheme
$17\frac{1}{1}\frac{2}{1}$ so that Hilbert's method does not (in
degree 8 as opposed to degree 6) encompass completely Harnack's
one.

\begin{figure}[h]\Figskip
%\vskip-1.2cm\penalty0
%\centering
\hskip-2.7cm\penalty0
\epsfig{figure=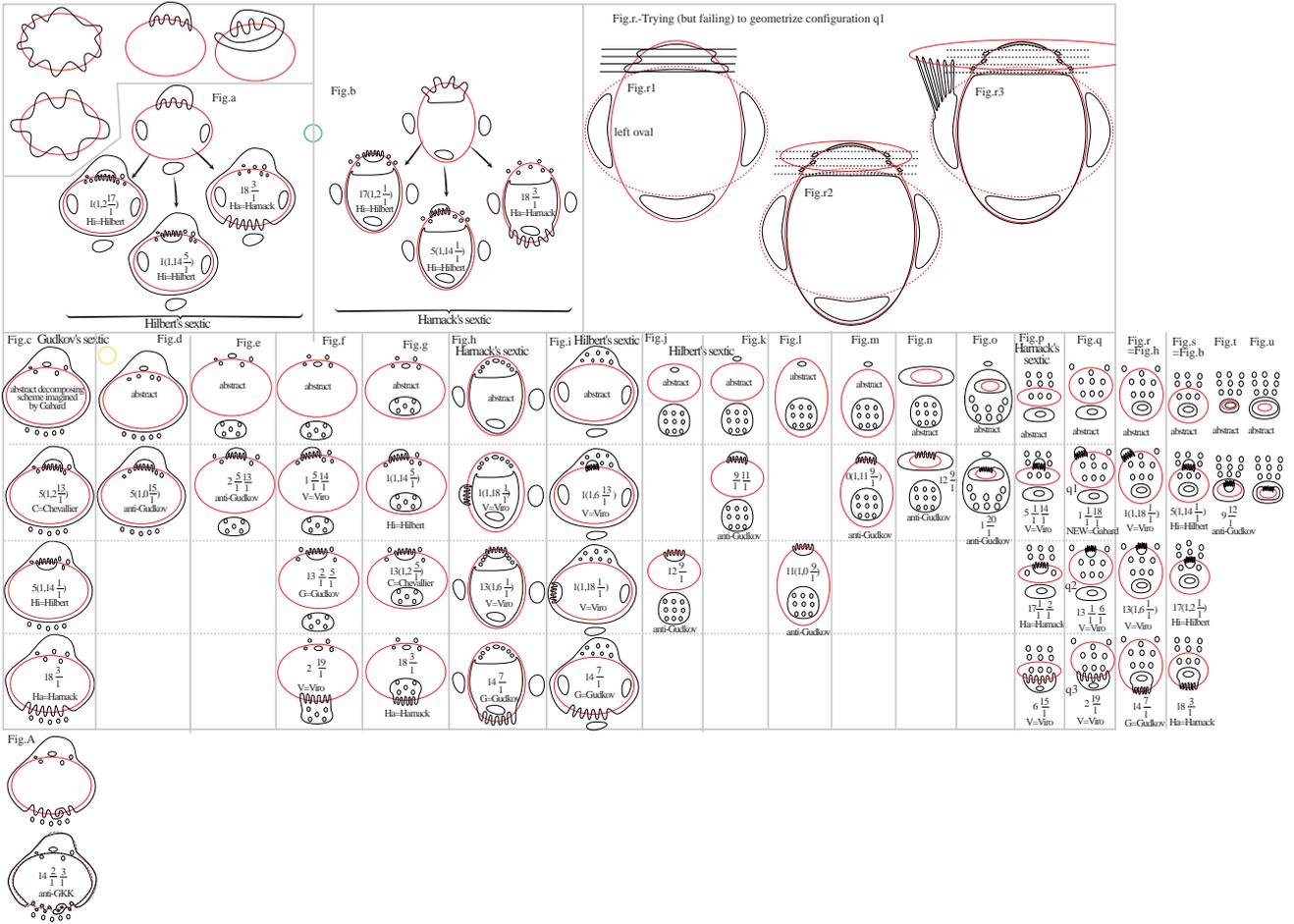,width=172mm}
\captionskipAG
  \caption{\label{ViroDEGREE8_GABARD_2+6_BOSON:fig}%
A flexible Harnack-Hilbert via floppy decomposing octics $2+6$
method yielding 3 schemes of Korchagin (green-colored), 2 schemes
of Chevallier (blue) and one boson only pseudo-realized by Orevkov
(black=Fig.\,b), plus finally one scheme prohibited by Shustin
(red=Fig.\,a). } \figskip
\end{figure}

On doing a similar yet more liberal (i.e. artistic) drawing with a
Gudkov-type sextic we get
Fig.\,\ref{ViroDEGREE8_GABARD_2+6_BOSON:fig}c whose first
smoothing yields the scheme $5(1,2 \frac{13}{1})$ first
constructed by Chevallier. Of course our construction is not a
serious one but maybe can turned to serious. Continuing with this
heuristic method we get Figs.\,d and e which cannot exist as they
create smoothings which are anti-Gudkov (i.e. violates Gudkov's
hypothesis corroborated by V.\,A. Rohlin). However the abstract
qualitative picture Fig.\,f create 3 curves all known to exists
(either first due Gudkov or to Viro) so there is some chance that
Fig.\,f exist (where there is a sort of infinitesimal vicinity of
certain ovals to the ground conic). This concept (of cytoplasmic
curve) is somewhat ill-posed yet fruitful to create the derived
decomposing curves with maximal oscillation. Fig.\,g is likewise
interesting producing another Chevallier's scheme $13(1,2
\frac{5}{1})$, beside more standard birds due to Harnack and
Hilbert respectively. Then it seems that we have exhausted all
possibilities with Gudkov's curve so as to be compatible with
Gudkov hypothesis (which permits only transfer of ovals by quanta
of four-packs). This brought us to Fig.\,h which is Harnack's
curve (of Fig.\,b, i.e. Harnack constructed \`a la Hilbert) with a
transfer of 4 quanta inside. Likewise Fig.\,i shows Hilbert's
sextic with a forced transfer of 4 ovals outside, so as to respect
formally Gudkov hypothesis. Looking at the resulting smoothings we
get curves due to Gudkov and Viro respectively. Thereafter we
consider the series of Figs.\,j,k,l,m and get always anti-Gudkov
curves, hence none of this abstract configuration exist
algebraically. As those 4 options exhaust the distribution of
Hilbert's 2 nests about the ellipse, we would infer that our
discussion is systematic. It remains yet to examine Fig.\,o, whose
outcome is also anti-Gudkov. So the ground ellipse must be inside
the nest and then Figs.\,a,i together give an exhaustive census of
the partition of the inner ovals compatible with the law of 4
quanta imposed by Gudkov periodicity.

Then we shall repeat this game with Harnack's curve
$9\frac{1}{1}$. So we starts (randomly with Fig.\,p) whose
production is pro-Gudkov, and actually includes one scheme of
Harnack and 2 due to Viro. Next we proceed to a delocalization by
a quanta of 4 to get Fig.\,q. Its first smoothing yields the boson
$1\frac{1}{1}\frac{18}{1}$ not yet known to be realized
algebraically (but known to be so pseudo-holomorphically by
Orevkov 2002
\cite{Orevkov_2001/02-classif-flexible-M-curves-degree-8}). So we
nearly made a progress on Hilbert's 16th, but alas not so due to
the purely heuristic character of our method (which is truly just
a flexible artistic avatar of Hilbert's method). Further as the
two other schemes derived from Fig.\,q are classical birds of
Viro, namely $13\frac{1}{1}\frac{6}{1}$ and $2\frac{19}{1}$, we
may get some evidence for the existence of the (bosonic) scheme
$1\frac{1}{1}\frac{18}{1}$. Here by bosonic we just mean hard to
detect!

[[{\it Added} [14.07.13] Actually we can do the depiction slightly
more rationally as on Fig.\,\ref{ViroDEGREE8_GABARD_2+6_BEST:fig}.
At some stage we had also the idea to transfer (package of 4)
ovals in the meander of the oscillation. Unfortunately this is not
possible since this would create another nest but the sextic is
already nested (Fig. a). Here we found a new Korchagin scheme
after transferring 4 ovals in a meander namely
$2(1,14\frac{4}{1})$.

However from there we can via transfers entering once in the red
ellipse and then sorting of it produce delocalization in the
meanders creating first Hilbert's scheme $1(1,14\frac{5}{1})$ and
then the scheme $0(1,14\frac{6}{1})$ prohibited by Shustin 90/91
\cite{Shustin_1990/91-New-restrictions}. So provided the latter
result is correct we have proven the lemma saying that it is
impossible to have a decomposing curve of degree $6+2$ like that
of Fig.\,a. Despite being of modest interest (if one is not an
aficionado of decomposing curves), it is nonetheless a severe
attack on our heuristic method since there is no hope that the
bosonic scheme just obtained as Fig.\,b will really exist.

\begin{figure}[h]\Figskip
%\vskip-1.2cm\penalty0
%\centering
\hskip-2.7cm\penalty0
\epsfig{figure=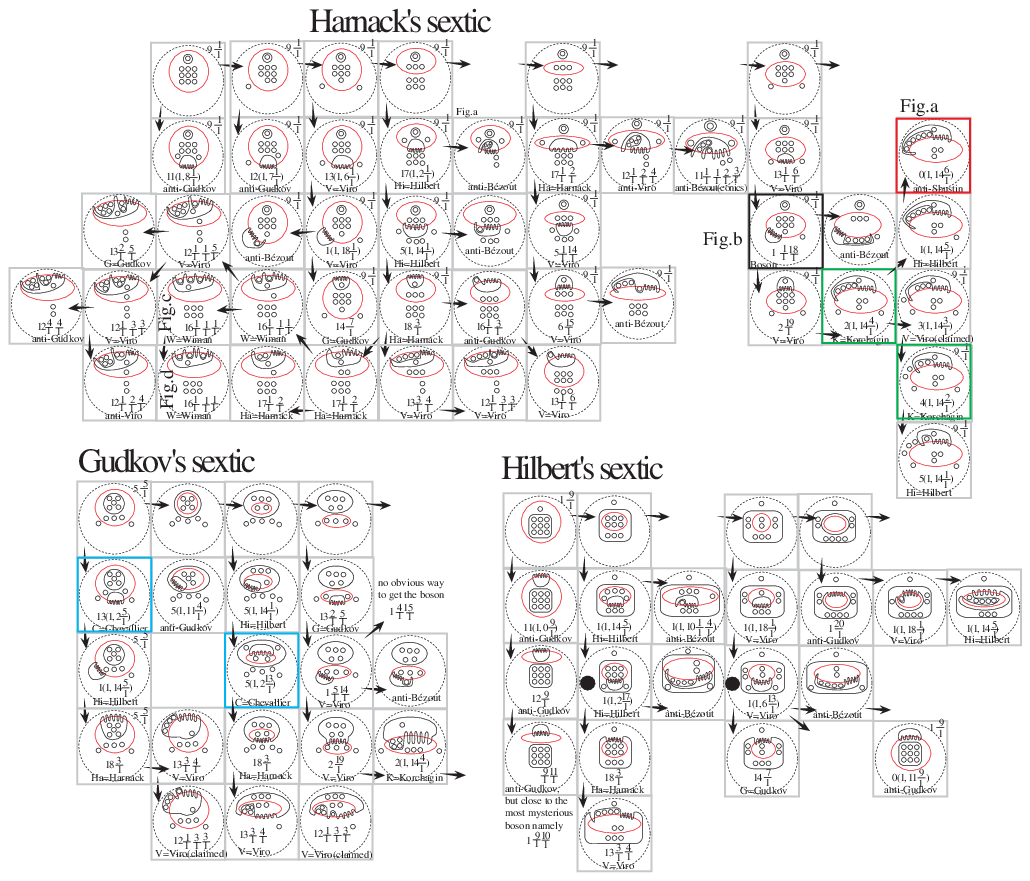,width=172mm}
\captionskipAG
  \caption{\label{ViroDEGREE8_GABARD_2+6_GOOD:fig}%
A flexible Harnack-Hilbert method yielding 2 schemes of Chevallier
and one boson pseudo-realized by Orevkov} \figskip
\end{figure}

[15.07.13] After a lengthy search we arrived at the conclusion
that we explored all decomposing octics of split-degree $2+6$ with
a maximally intersecting oval which is undulating so as to produce
an $M$-curve after perturbation of the nodes. Equivalently the
undulating condition may be expressed by saying that the order of
the 12 intersection points is the same along the ellipse as along
the oval. So we would like to state:

\begin{lemma}
A floppy decomposing curve of bidegree $2+6$ belongs to one of the
isotopy type tabulated on
Fig.\,\ref{ViroDEGREE8_GABARD_2+6_GOOD:fig} plus eventually some
little variants like Figs.\,c and d. Despite some ambiguities we
believe that the resulting list of octic $M$-schemes accessible
via a bi-curve of bidegree $2+6$ is exhaustive, and labelled by
orange colored frames on Fig.\,\ref{SIMPLIFIED-TABLE:fig}.
Actually Fig.\,\ref{ViroDEGREE8_GABARD_2+6_GOOD:fig} has to be
replaced by Fig.\,\ref{ViroDEGREE8_GABARD_2+6_BEST:fig} where 2
additional schemes were discovered.
Likewise we shall see (by virtue of
Fig.\,\ref{ViroDEGREE8_GABARD_1+7_GOOD:fig}) that all octic
$M$-schemes in the light-blue regions of
Fig.\,\ref{SIMPLIFIED-TABLE:fig} are realized through a floppy
split-curve of degree $1+7$ (sometimes referred to as affine
septics, yet this concept is alas not universally defined, compare
a well-known (amicable) controversy between Korchagin-Shustin and
Orevkov, cf. optionally Orevkov 1998
\cite{Orevkov_1998}).\end{lemma}

This was essentially safe that there is an even somewhat more
rational way to present things via
Fig.\,\ref{ViroDEGREE8_GABARD_2+6_BEST:fig} where 2 additional
(2+6)-split $M$-octic schemes are constructed. The underlying idea
is that one starts with simply-nested schemes, e.g. Harnack's
schemes $18\frac{3}{1}$ and then elevates by putting ovals in the
fingers meanders. (Fig.\,a$^\ast$ just shows a useless quantum
transfer to a Gudkov sextic yet producing the same triad as
before.) Then one moves right by a quantum transfer of 4 ovals to
get via Fig.\,b Gudkov's scheme $14\frac{7}{1}$ and mowing upwards
by the finger moves yields a menagerie of 8 schemes moving fairly
linearly on the geographical table. (Fig.\,b$^\ast$ shows an
illegal transfer by the way violating Gudkov periodicity.) Next
one considers the split curve of Fig.\,c yielding the scheme
$6\frac{15}{1}$ and look at its versatile descendance under finger
moves, which one the table correspond to a rectilinear motion with
jump from the 0th pyramid to the 3rd one. Fig.\,d is deduced by a
quantum transfer of 4 ovals. Fig.\,e is deduced by tranquilizing
the vibrating oval and making vibrate the external one, while
Fig.\,f and Fig.\,g should be self-explanatory. Of course albeit
being better organized it is relatively tricky to get convinced
that our search is exhaustive as far as the isotopy type of the
resulting octic $M$-schemes is concerned.

\begin{figure}[h]\Figskip
%\vskip-1.2cm\penalty0
%\centering
\hskip-2.7cm\penalty0
\epsfig{figure=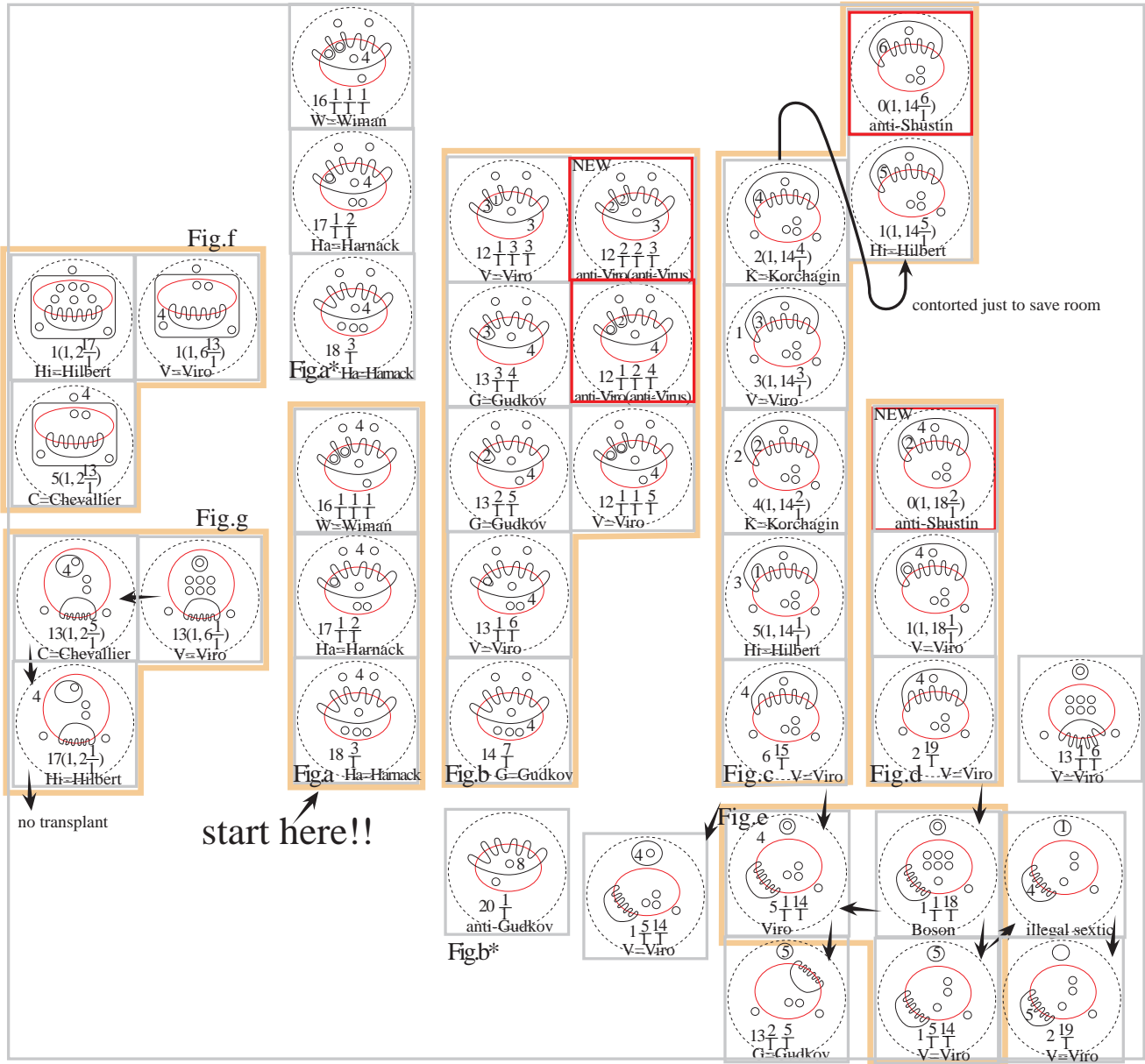,width=172mm}
\captionskipAG
  \caption{\label{ViroDEGREE8_GABARD_2+6_BEST:fig}%
A more rational enumeration yielding 2 supplementary schemes}
\figskip
\end{figure}

(Little idea added [16.07.13]) As we said Fig.\,e contains a new
bosonic octic $M$-scheme namely $1\frac{1}{1}\frac{18}{1}$. Albeit
the latter is not yet known to exist (nor to be prohibited) our
floppy sexy-conic (=split curve of degree $2+6$) is perhaps easier
to prohibit maybe via Thom's genus bound \`a la Mikhalkin (we owe
this idea from an e-mail of Th. Fiedler, cf. eventually
%Sec.\,\ref{e-mail-Viro:sec})
v.2 of Ahlfors. One would like to construct a membrane but we are
not very inspired for the moment.

In particular we see from this geographical report of the
organical game of tracing floppy split octic (in French
``d\'eployement universel du symbole de Gudkov'') that the class
of degree $2+6$ create one boson (not yet known to exist) but
enters twice in conflicts with known prohibitions, namely with
$12\frac{1}{1}\frac{2}{1}\frac{4}{1}$ we enter in conflict with
Viro's imparity law (which is in principle well-established, but
we confess to have not yet studied with enough care its proof to
warranty that the result is true), and beside we enter once in
conflict with Shustin's obstruction of five at the scheme
$(1,14\frac{6}{1})$ (again we have not yet studied Shustin's
proof). Concerning the class of split-curves of degree $1+7$, they
enter frankly in conflict with Viro's obstructions (both the
imparity law as well as the sporadic Viro obstruction when it
comes to the scheme $4\frac{3}{1}\frac{3}{1}\frac{9}{1}$). In
contrast when looking at the 3rd pyramid (encoding schemes with a
subnest), the affine septics are in perfect agreement with the
known construction, and even stronger than that on the blue
regions since it would suggest existence of the bosonic scheme
$14(1,2\frac{4}{1})$.

Finally, albeit there is some overlap between the schemes covered
by both procedures ($2+6$ vs. $1+7$) there is also some
complementary nature of their domain of influences. Of course it
remains to tabulate the geographical position of the split-curves
$3+5$ and $4+4$. What is somewhat surprising is that for $1+7$ we
could give a very regular enumeration whilst for $2+6$ we suffered
under messiness and a chaotical somewhat random search. Perhaps
this is due to the fact that we imposed a too prescription of the
sextic schemes versus letting operate all quantum fluctuations of
ovals compatible with B\'ezout and Gudkov periodicity (optionally
Viro's imparity law). So at the occasion one should try to
reorganize the $2+6$ table.

Finally it may also be observed that the axiom all floppy split
curves does not enter in conflict with Fiedler's obstruction of
four schemes (but as we already noted seriously damage or is
damaged by Viro's imparity law for trinested schemes).

%%%ostensible

\begin{figure}[h]\Figskip
%\vskip-1.2cm\penalty0
%\centering
\hskip-2.7cm\penalty0
\epsfig{figure=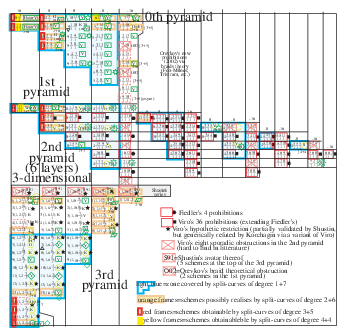,width=172mm} \captionskipAG
  \caption{\label{SIMPLIFIED-TABLE:fig}%
A flexible Harnack-Hilbert method yielding 2 schemes of Chevallier
and one boson pseudo-realized by Orevkov} \figskip
\end{figure}

Let us now examine split curves of degree $3+5$. Reporting on the
map the simple pictures of Fig.\,\ref{ViroDEGREE8_GABARD_3+5:fig}
we get the (little) red framed schemes on
Fig.\,\ref{SIMPLIFIED-TABLE:fig}. The most interesting point is of
course the realizability of the boson $14(1,2\frac{4}{1})$ which
is now realized twice (once in degree $1+7$ and anew in degree
$3+5$). This should perhaps indicate that this boson is most
likely to exist, but of course our viewpoint is so naively
topological that there is no certitude yet. Nonetheless this may
indicate that $14(1,2\frac{4}{1})$ is the less mysterious of all
sic bosons. Actually I should confess that during the enumeration
we first shamefully forgot the schemes below $14(1,2\frac{4}{1})$
and found them only after contemplating the architecture of the
pyramid. So the lesson to keep in mind is always to geometrize the
combinatorics in order to miss nobody.

Finally for bidegree $4+4$, the corresponding schemes are explored
on Fig.\,\ref{ViroDEGREE8_GABARD_4+4:fig}, and their geography
reported by small yellow rectangles on the main-map
(Fig.\,\ref{SIMPLIFIED-TABLE:fig}). The series of schemes so
obtained is a nearly ridiculous collection of 4 schemes, yet
contains Wiman's famous schemes
$16\frac{1}{1}\frac{1}{1}\frac{1}{1}$, which historically was
first obtained through construction of the appropriate $4+4$-slit
curve in the algebro-geometric category (Wiman 1923
\cite{Wiman_1923}) thereby contorting some misconceptions of
Hilbert/Ragsdale. Somewhat more interestingly we also encounter
Viro's scheme $6\frac{15}{1}$ of which it would be nice to know if
Wiman was technologically able to construct it.

It is also interesting to observe that schemes of type $4+4$ form
a subfamily of those of bidegree $2+6$, while those of type $3+5$
constitute a subfamily of those of type $1+7$. Apart from those
observation the resulting architecture looks still a bit
mysterious and our all endeavor shed only minimalist spot of
lights on the overall global Hilbert's 16th problem.

%%%%obtainable

[16h41] Of course one can wonder if our Walt-Disney/Gudkov
flexible depiction can be rigidified in the realm of algebraic
geometry so as to get an existence proof of the (bosonic) scheme
$1\frac{1}{1}\frac{18}{1}$ at the algebraic level. Crudely
speaking it seems to suffice to look at Fig.\,b (right) while
tracing an ellipse enclosing the seven ``upper'' ovals (upperness
being defined w.r.t. to the picture). The difficulty however is to
ensure that 2 remaining outer oval are proximal enough to the
ellipse as to effect the required schematic of Fig.\,q1 (where
``q1'' refers to the first row below Fig.\,q). If feasible then
there is perhaps a perfectly elementary proof of existence of the
bosonic scheme $1\frac{1}{1}\frac{18}{1}$, which is completely at
the level of the Harnack-Hilbert technological level.

Let us state this as follows:

\begin{lemma}
If there exist a decomposing octic of degree $8=6+2$ whose scheme
is Fig.\,\ref{ViroDEGREE8_GABARD_2+6_BOSON:fig}q1 then there
exists a smooth octic with (bosonic) new scheme
$1\frac{1}{1}\frac{18}{1}$, and one wons a Fields medal in
chocolate for a spectacular advance on Hilbert's 16th.
\end{lemma}

So configuration q1 is gold-worth (Goldwert f\"ur
Hilbert/Viro/Orevkov) and trying to geometrize it on the larger
Fig.\,r we meet some evident obstruction on the vibratory model of
Hilbert. First the red ellipse trying to phagocyte all 7 ``micro''
ovals tends to collide with the nonempty oval. Of course this pity
can be salvaged by imagining the ellipse of very large
eccentricity (hence nearly osculating the bottommost horizontal
line). It remains then to take the left (say) macro-oval and to
make it oscillate across the ground ellipse (as heuristically
depicted on Fig.\,r3). This looks structurally hard to do (without
contorting algebraic respectableness) and we see why Orevkov's
bosons $1\frac{1}{1}\frac{18}{1}$ is so hard to construct. Of
course, it can be that there is an obstruction to this scheme.

Further one can of course mentally play with the geometric
parameters of Fig.\,r1 by bringing the bottommost (horizontal)
line closest possible to the crossing of the 2 primitive ellipses,
and then hope that the left oval is close enough to the ground
ellipse of Fig.\,r2 to vibrate across it. This looks a little
puzzle with infinitesimals (we confess hard to believe in).

[07.07.13] What is fairly incredible is that our method (which is
just a heuristic=flexible Harnack-Hilbert-Ragsdale-Brusotti-Wiman
method) is very versatile realizing most of Viro's scheme by the
dissipation of the simplest singularity $A_1$ (ordinary node with
2 real tangents). Hence supposing some intelligence able to
implement at the rigorous level this would constitute a little
attack upon the slogan that Viro's method is structurally stronger
that the classical perturbation method. Albeit vague, our remark
should be precise enough to open a little debate on the point
after more work.

Next we examine Fig.\,r (where a ring is delocalized so that
$\chi$ keeps unchanged), and then considered Fig.\,s. (Alas we did
not realized first that those configurations were already analyzed
as Fig.\,h and Fig.\,b respectively). It seems at this stage that
we have analyzed all possibilities of an elliptical oval located
w.r.t. one of the 3 possible $M$-sextic. Of course this does not
mean that we do have catalogued all decomposing octic of bidegree
$8=6+2$, because a priori the intersection may be not a vibration
but a more complicated ``meander'' as depicted say on the top
Fig.\,\ref{ViroDEGREE8_GABARD_2+6_BOSON:fig}. So for instance
Fig.\,A shows a meander, but it  seems that this causes a loss of
vibratory energy striving us outside the realm of $M$-curves (the
smoothing below being only an $(M-1)$-scheme, actually forbidden
by GKK (the $(M-1)$-avatar of Gudkov periodicity due to
Gudkov-Krakhnov-Kharlamov).

Of course if meander fails to produce $M$-schemes it remains to
investigate the other splitting of 8 as $4+4$, $3+5$, $1+7$.

\subsection{Degree 4+4}

Let us start for evenness psychological simplicity with the
splitting $8=4+4$ (cf. Fig.\,\ref{ViroDEGREE8_GABARD_4+4:fig}).
Here as well-known we recover Wiman's curve
$16\frac{1}{1}\frac{1}{1}\frac{1}{1}$, plus schemes due to
Harnack, and more interestingly Fig.\,f gives a schemes first
cooked by Viro, namely $6\frac{15}{1}$. It would be of didactic
interest to know if this scheme is (rigorously) constructible by
an elementary method \`a la Wiman, thereby circumventing the
intrusion of Viro.

\begin{figure}[h]\Figskip
%\vskip-1.2cm\penalty0
%\centering
\hskip-2.7cm\penalty0
\epsfig{figure=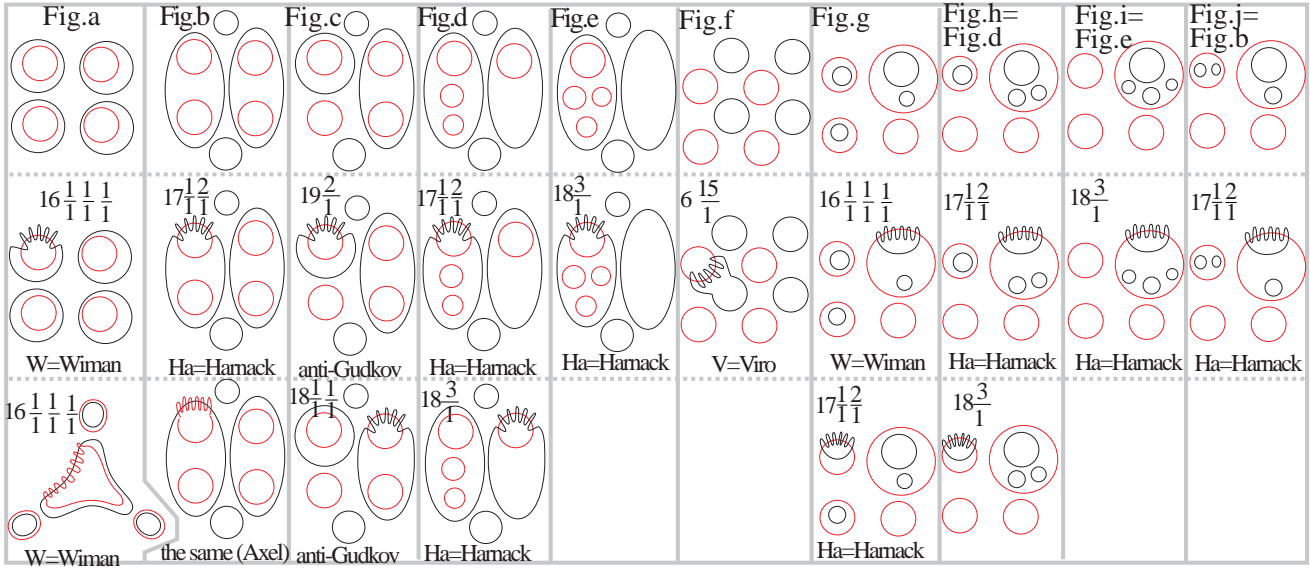,width=172mm}
\captionskipAG
  \caption{\label{ViroDEGREE8_GABARD_4+4:fig}%
Smoothings of decomposable octics of bidegree $4+4$} \figskip
\end{figure}

\subsection{Degree 3+5 (revisited)}

The next case to examine is $8=3+5$. Here we obtain
Fig.\,\ref{ViroDEGREE8_GABARD_3+5:fig}. Fig.\,b,c,d illustrates
the standard phenomenon that when splitting the oscillation into
two circuit we loose vibratory energy failing so to reach the
Harnack maximal case. On the sequel of the map should be
self-explanatory, and it is noteworthy that we get again (added
[14.07.13]) the bosonic scheme $14(1,2\frac{4}{1})$, which we also
realized via a curve of degree $1+7$ (cf.
Fig.\,\ref{ViroDEGREE8_GABARD_1+7_GOOD:fig}). However it could be
that the present schematic view suggest a better technique to
construct the curves as the  individual curves involved in the
splitting have lower degrees hence perhaps easier to control. One
could perhaps try to attack the construction of such a split curve
as an interpolation problem for a cubic given a fixed  $C_5$ in
the background landscape. One could try to start with a Harnack
quintic (e.g. in the model constructed \`a la Hilbert) and then
try to trace the appropriate cubic. Of course the problem looks
violently over-determined and we can only tabulate on a very lucky
stroke to get out. Of course conversely one could imagine as a
reverse engineering telling that any curve $14(1,2\frac{4}{1})$
could be through a dynamical procedure degenerate toward the split
curve depicted (this being maybe reminiscent of a Hilbert-Rohn
type method). Supposing further that one is able to prohibit our
split curve (via say a theory \`a la Orevkov), then one would get
an obstruction upon the bosonic scheme.

\begin{figure}[h]\Figskip
%\vskip-1.2cm\penalty0
%\centering
\hskip-2.7cm\penalty0
\epsfig{figure=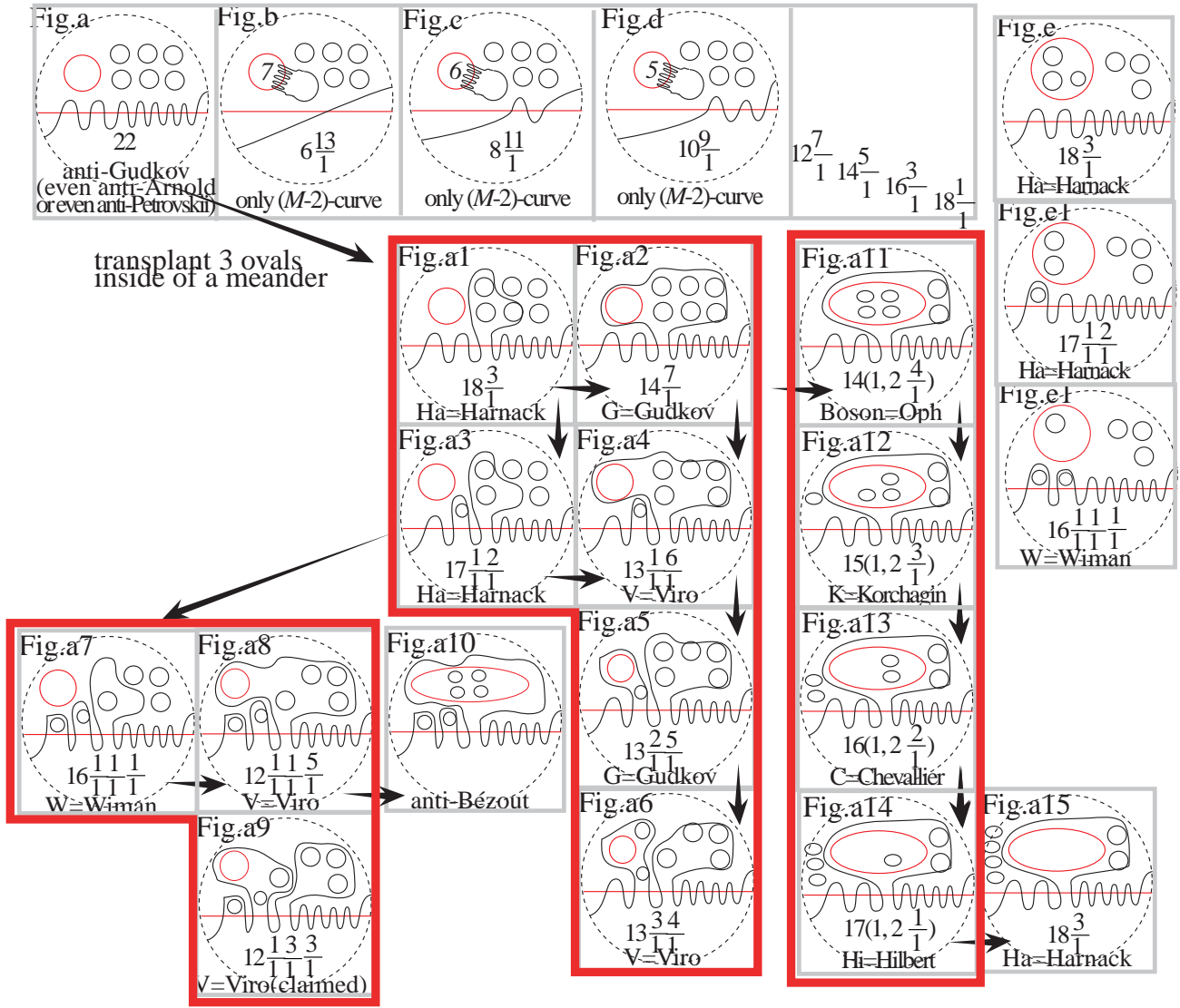,width=172mm}
\captionskipAG
  \caption{\label{ViroDEGREE8_GABARD_3+5:fig}%
Smoothings of split octics of degree $3+5$ (at first view only one
Harnack scheme, but then fairly original Korchagin and Chevallier
schemes, plus even the pseudoholomorphic boson
$14(1,2\frac{4}{1})$).} \figskip
\end{figure}

\subsection{Degree 1+7}

Thereafter we examine $8=1+7$. Here we base the analysis (situs)
upon Viro's census of $M$-septics. This includes precisely the
fourteen $M$-septics schemes $15$, $13\frac{1}{1}$,
$12\frac{2}{1}$, $11\frac{3}{1}$, etc, $2\frac{12}{1}$,
$1\frac{13}{1}$ (where the pseudoline $J$ is omitted from the
symbolism). On looking at the
Fig.\,\ref{ViroDEGREE8_GABARD_1+7:fig} (whose geometric essence is
to say that a maximally dissipable decomposing curve occur when
the  pseudoline is undulating across the line) we see see a fairly
disappointing issue that this method will only created
simply-nested schemes (whose theory is already completely
elucidated through Viro's construction in degree 8). Actually on
taking as septics those with scheme $3\frac{11}{1}$ we get only
$10\frac{11}{1}$ and so this construction do not even realize the
two schemes $6\frac{15}{1}$ and $2\frac{19}{1}$.

\begin{figure}[h]\Figskip
%\vskip-1.2cm\penalty0
%\centering
\hskip-2.7cm\penalty0
\epsfig{figure=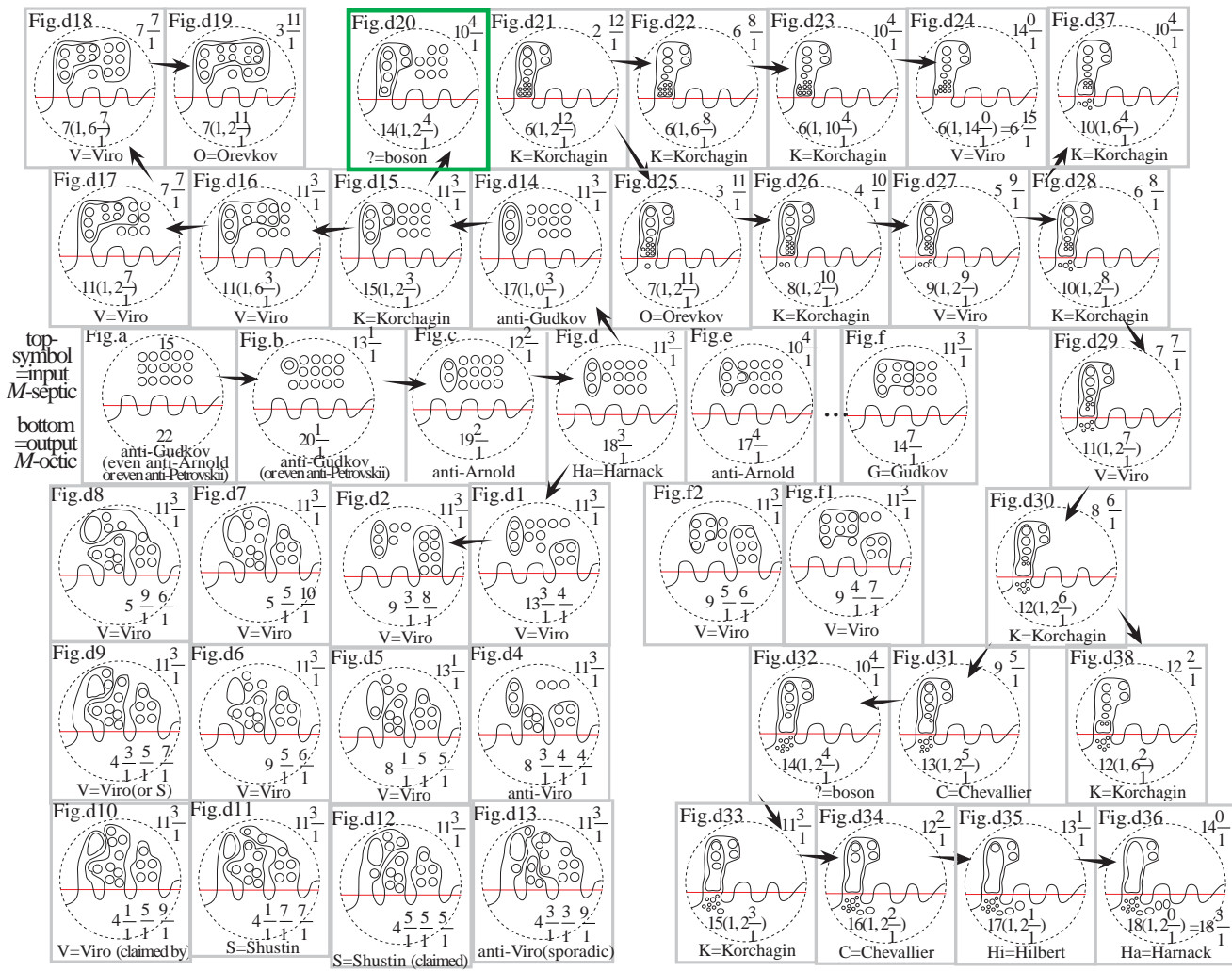,width=172mm}
\captionskipAG
  \caption{\label{ViroDEGREE8_GABARD_1+7:fig}%
Smoothings of decomposable octics of bidegree $1+7$
%(only one
%Harnack scheme)
(cf. next Fig. for a cleaner way)} \figskip
\end{figure}

As a consequence it seems that the method of decomposable curves
(combined with Brusotti's independence of   smoothing for ordinary
nodes) afford only a small list of $M$-schemes, yet some due to
Viro, and 2 of Chevallier, as well as one boson
$1\frac{1}{1}\frac{18}{1}$ not yet known to exist. Hence the
method deserves be investigated more systematically to be sure
that we missed nobody, and then needs to be geometrized in order
to see if the above mentioned boson can be constructed, what
nobody succeeded hitherto to do.

[11.07.13] Of course there is then much more configuration to
analyze and the problem amounts to the classification of affine
$M$-septics. It is notorious that already the case of $M$-sextics
is extremely difficult (initiated by Korchagin/Shustin,
Polotovskii and Orevkov and perhaps fairly close to completion
now). The point is that ovals may be situated in the meanders of
the pseudoline oscillating across the line and this will produce
additional $M$-schemes (or rather pseudo $M$-schemes as we are not
ensured a priori that the configuration exists). However one can
start a qualitative exploration by transferring packages (quanta)
of 4 ovals. So starting from the above Fig.\,d, and transferring
one quanta of 4 outer oval in the meander while smoothing the
resulting $C_7\cup L$ produces the scheme $13
\frac{3}{1}\frac{4}{1}$ which is known to exist by Viro. Of course
continuing in this fashion, we get the $9 \frac{3}{1}\frac{8}{1}$
(also due to Viro), and then $5 \frac{3}{1}\frac{12}{1}$ (also
Virotian), and finally $1 \frac{3}{1}\frac{16}{1}$. The latter is
actually prohibited by Orevkov's braid theory (cf. Orevkov 2002
\cite{Orevkov_2001/02-classif-flexible-M-curves-degree-8}). So our
method is completely hopeless, yet heuristic. Alas it seems that
the originality of our method resides i the fact that we are not
readily considering classical deterministic curve like the famous
Russian decomposing curves. (For a bribe of decomposing curve in
primitive Germany, cf. Mohrmann's 1912 \cite{Mohrmann_1912} lovely
introduction.) In fact we considered on the earlier pictures in
degree $6+2$ and $4+4$
(Figs.\,\ref{ViroDEGREE8_GABARD_2+6_BOSON:fig} and
\ref{ViroDEGREE8_GABARD_4+4:fig} respectively), what could be
called a {\it quantum fluctuating curve}. This is like a
decomposing curve save that both constituents are disjoint yet
with each oval susceptible to oscillates across each of his
neighbors. The power of the method was that we got so series
(triads in general) of classical decomposing curve, which produce
by Brusotti's smoothing some (potentially new) $M$-curves. The
heuristic idea was then that if two of the three smoothing exists
then there is some hope that the 3rd product is also modellizing a
real algebraic curve. Alas when working in odd+odd bidegree both
pseudolines are forced to intersect and so the intersection has to
be monopolized on the pseudolines and we loose the quantum
plurality of matter (where an object admits a nebulosity of things
around it like its descendence). This being said our method seems
to loose its multivalued-ness which was the source of a little
principle of rigidity aping very vaguely the algebraic rigidity of
the real world we are interested in.

However on working more carefully the (quantum) transfers as
Fig.\,d1. and d2 we see that we met an obstruction not allowing
one to reach Viro's scheme $5\frac{3}{1}\frac{12}{1}$, whence a
fortiori not Orevkov's scheme $1\frac{3}{1}\frac{16}{1}$. However
on propagating the quantum fluctuations (suitable transfers of
ovals) we see that we may obtain some schemes due to Shustin
(hence lying somewhat deeper than the first generation of Viro),
like $4\frac{1}{1}\frac{7}{1}\frac{7}{1}$ (Fig.\,d12) or
$4\frac{5}{1}\frac{5}{1}\frac{5}{1}$. All this is pleasant
depiction (morphogenesis or waves with oxygenating bubbles) yet it
seems that there are obstruction to reach by the method more
Orevkov(ian)
 schemes on the first pyramid (of doubly nested schemes) and with
 $\chi=-16$ (i.e. extreme-right portion of the diagram). Of course
 even Viro's schemes on this strip look inaccessible via
 decomposing curve of degree $7+1$. This deserves to be better
 understood at the occasion yet it seems that there is a fairly
simple
%%%experimental
explanation due to the oscillating character of the pseudoline and
the aptitude for at most three meanders to tolerate ovals.
(Otherwise the smoothed octic violates B\'ezout for conics, or
equivalently the doubled quadrifolium
$\frac{1}{1}\frac{1}{1}\frac{1}{1}\frac{1}{1}$ is saturated.)

[13.07.13] Ultimately we noticed that from Fig.\,d5, we changed
the isotopy class of the septics to $13\frac{1}{1}$, and later on
to even $15$ a configuration first constructed by Ragsdale 1906
\cite{Ragsdale_1906} and later rediscovered by Wiman 1923
\cite{Wiman_1923}.

Eventually we had the idea that in order to have a super-nest
(i.e. a subnest in a nest so as to land in the 3rd pyramid), we
only have to put the nest of the septic into the meander as on
Fig.\,d14. This has wrong characteristic but it is a simple matter
to correct this with Fig.\,d15 realizing the scheme
$15(1,2\frac{3}{1})$ first constructed by Korchagin. Then one can
do quantic jumps with packets of 4 ovals up to reach Fig.\,d19
which is a truly remarkable scheme $7(1,2\frac{11}{1})$
 due to Orevkov (essentially the last one constructed up-to-date).
 Even more cleverly, turning back again to Fig.\,d15 (which is
fairly close to the boson $14(1,2\frac{4}{1})$) we realize that it
is enough to drag an outer oval of $C_7$ in the meander and let it
be phagocytozed by the nonempty oval to get the bosonic scheme
$14(1,2\frac{4}{1})$ resisting to present knowledge.
We have proved the modest:

\begin{lemma}
There is no (naive) topological obstruction to the realizability
of the bosonic $M$-scheme $14(1,2\frac{4}{1})$ by  a Brusotti
perturbation of a decomposing curve of degree $7+1$ where the
septics has (reduced) scheme $10\frac{4}{1}$ (pseudoline omitted)
while the line intercepts it maximally as depicted on
Fig.\,\ref{ViroDEGREE8_GABARD_1+7:fig}d20.
\end{lemma}

Dragging the 8 remaining outer ovals into the meander create
$6(1,2\frac{12}{1})$ due to Korchagin. Alas it seems evident that
we will never reach the other boson with a subnest, i.e.
$4(1,2\frac{14}{1})$ because to many outer ovals are created by
the meanders of the oscillation.

Then, starting from Fig.\,d21 (realizing Korchagin's scheme
$6(1,2\frac{12}{1})$), one can imagine a progressive transfer of
the ovals at depth 2 to oval at depth 0 (hence no quanta rule of 4
has to be respected) and we get so Fig.\,d25 which is again
Orevkov's (unique) scheme $7(1,2\frac{11}{1})$. Continuing, this
process sweeps out  all schemes of the fundamental table
(Fig.\,\ref{SIMPLIFIED-TABLE:fig}) yielding successively the
schemes $8(1,2\frac{10}{1})$~(K), $9(1,2\frac{9}{1})$~(V),
$10(1,2\frac{8}{1})$~(K), $11(1,2\frac{7}{1})$~(V),
$12(1,2\frac{6}{1})$~(K), $13(1,2\frac{5}{1})$~(C),
$14(1,2\frac{4}{1})$~(Oph=boson), $15(1,2\frac{3}{1})$~(K),
$16(1,2\frac{2}{1})$~(C), $17(1,2\frac{1}{1})$~(Hi),
$18(1,2\frac{0}{1})=18\frac{3}{1}$~(Ha). As usual, we use the
following abbreviations: Ha=Harnack, Hi=Hilbert, V=Viro,
K=Korchagin, C=Chevallier, O=Orevkov, Oph=Orevkov but only
pseudo-holomorphically. Of course in view of our earlier twist
(Fig.\,d21 up to d24) it is clear that on the main-table we can
move horizontally to the right as well so as to sweep out the
whole portion of the pyramid below Korchagin's
$6(1,2\frac{12}{1})$. So for instance Fig.\,d37 gives Korchagin's
scheme $10(1,6\frac{4}{1})$, and likewise Fig.\,d38 gives
Korchagin's scheme $12(1,6\frac{2}{1})$. Applying the same lateral
dynamics to Chevallier's scheme $13(1,2\frac{5}{1})$ yields merely
Viro's scheme $13(1,6\frac{1}{1})$, yet this might give confidence
in the method. Applying the lateral shift to ``Orevkov's'' boson
$14(1,2\frac{4}{1})$  produces the scheme
$14(1,6\frac{0}{1})=14\frac{7}{1}$ due to Gudkov. Hence this may
give some evidence that the boson in question exists
(algebraically).

At this stage, it seems that we were fairly exhaustive and we
would like a statement about the combinatorial confinement of
$M$-schemes obtained by small perturbation of a decomposing octic
splitting of a line plus a septic.

Before let us make some more basic experiment. Let us start with
Fig.\,\ref{ViroDEGREE8_GABARD_1+7_GOOD:fig}g (which is actually a
perfect replica of the former Fig.\,a). Then as we are interested
mostly in doubly-nested schemes (enclosing 4 mysterious bosons) we
transfer ovals in the meanders, e.g. just one in each meanders to
get Fig.\,g1 which gives the scheme $18\frac{1}{1}\frac{1}{1}$
which is anti-Gudkov. Looking at the map
(Fig.\,\ref{SIMPLIFIED-TABLE:fig}) we see that the closest scheme
is $17\frac{1}{1}\frac{2}{1}$ (Ha) and so we correct to Fig.\,g2
(successive approximation). After that we do quantic jumps of 4
ovals to derive laterally on the right through Figs. g3, g4, g5
but alas cannot so reach the boson $1\frac{1}{1}\frac{18}{1}$
(which we could however find via a $6+2$ splitting, cf.
Fig.\,\ref{ViroDEGREE8_GABARD_2+6_BOSON:fig}). Then looking at the
pyramid again (Fig.\,\ref{SIMPLIFIED-TABLE:fig}) it is clear how
to move down to get $13\frac{2}{1}\frac{5}{1}$ (Fig.\,g6), which
is due to Gudkov. (Probably Gudkov's original construction really
involves this stronger decomposing scheme). This can in turn
derives to the right, yet cannot reach the scheme
$1\frac{2}{1}\frac{17}{1}$ (due to Viro however). Continuing in
the obvious way gives all pictures in the first triangle of
Fig.\,\ref{ViroDEGREE8_GABARD_1+7_GOOD:fig} proving thereby the:

\begin{lemma}
Among all $M$-schemes of the 1st pyramid of
Fig.\,\ref{SIMPLIFIED-TABLE:fig} (consisting of so-called doubly
nested schemes plus those which are simply-nested) only those
located in the sub-triangle where the right edge is deleted are
(potentially) realizable through perturbation of a decomposing
curve of degree $7+1$. In particular it seems that there is no
chance to construct via the splitting $7+1=8$ the four
doubly-nested bosons, namely $1\frac{1}{1}\frac{18}{1}$,
$1\frac{4}{1}\frac{15}{1}$, $1\frac{7}{1}\frac{12}{1}$,
$1\frac{10}{1}\frac{9}{1}$. Actually the first such boson is
constructible via a ``flexible'' sexti-conic of degree $6+2$.
\end{lemma}

Actually it is evident that our messy random table
(Fig.\,\ref{ViroDEGREE8_GABARD_1+7:fig}) can be improved into the
following one (Fig.\,\ref{ViroDEGREE8_GABARD_1+7_GOOD:fig}) whose
architecture is directly adapted to that of the main-pyramid
(Fig.\,\ref{SIMPLIFIED-TABLE:fig}).

\begin{figure}[h]\Figskip
%\vskip-1.2cm\penalty0
%\centering
\hskip-2.7cm\penalty0
\epsfig{figure=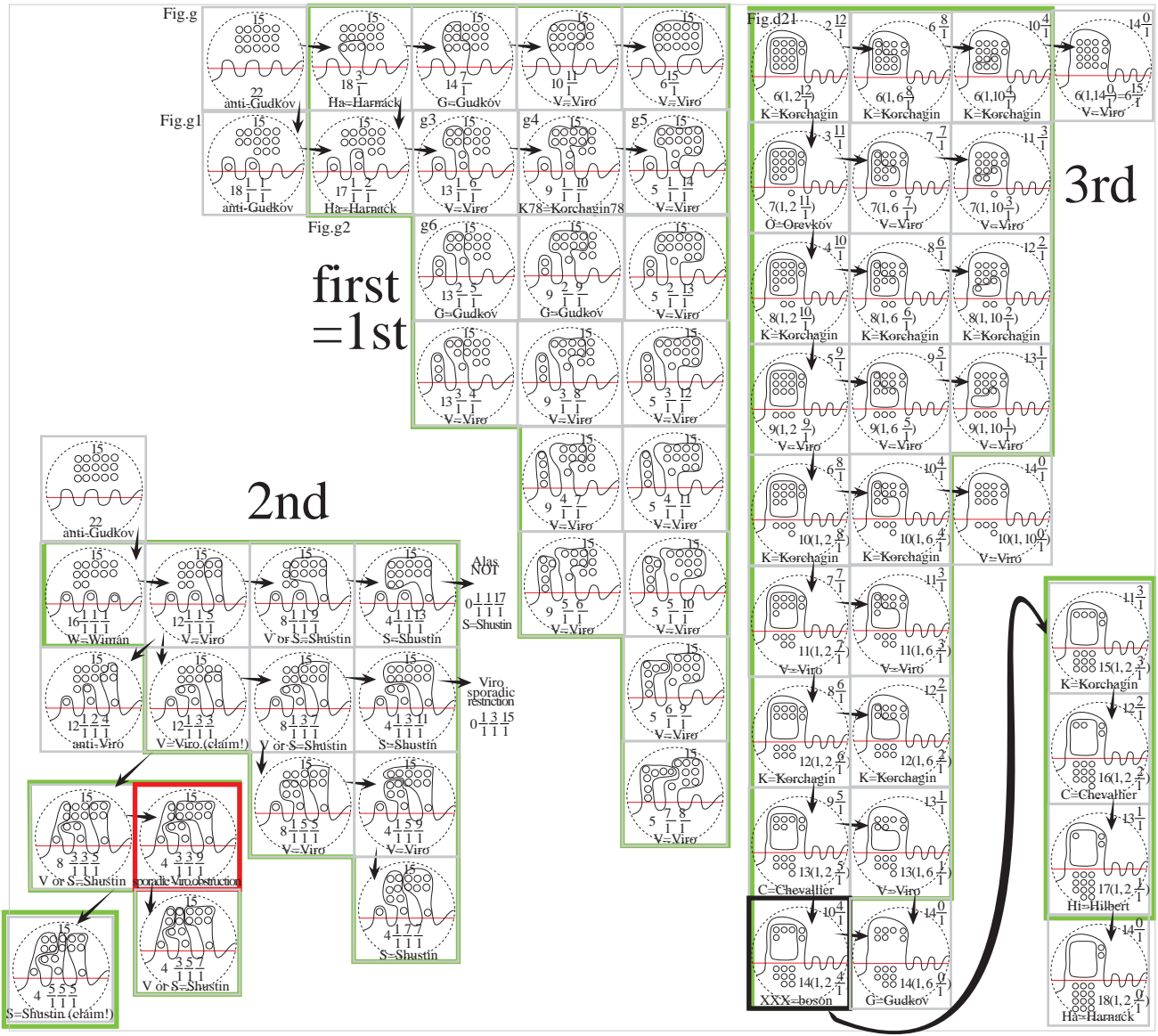,width=172mm}
\captionskipAG
  \caption{\label{ViroDEGREE8_GABARD_1+7_GOOD:fig}%
Smoothings of decomposable octics of bidegree $1+7$: one recovers
many classical schemes, and also one new one yet also enters in
conflict with one Viro sporadic obstruction (of course we suited
to the Fiedler-Viro oddity law for trinested schemes already.)}
\figskip
\end{figure}

As a  moral, Fig.\,\ref{ViroDEGREE8_GABARD_1+7:fig} shows that
potentially Brusotti's classical method is susceptible of
recovering many schemes (in particular many of those of Korchagin
that Viro himself conjectured not to exist). Hence potentially
Brusotti's method could be nearly as puissant as Viro's method,
yet alas apparently nobody ever succeeded to trace those splitting
curves in the algebraic category. Even more importantly, we remark
that the classical Harnack-Brusotti method (with a ground line) is
susceptible of yielding one bosonic scheme (namely
$14(1,2\frac{4}{1})$).

\section{More artwork via
other distributions of singularities}

\subsection{Switching to affine sextics}

{\it Added\/} [17.07.13].---There is a fairly interesting Master
thesis by Daniel Eric Smith 20XX
\cite{Smith_2005-Thesis-under-Korchagin-split-7+1} (one of
Korchagin's student) who exhibits some algebraic models of
$(1+7)$-split octics. Alas the specimens so obtained looks rather
ridiculous and a classification looks fairly out of reach for the
moment. Of course for our purpose of constructing $M$-octics one
does not need to go through a complete census of maximally
interesting oval (in one case of which Smith's work affords no
data) but which is conjectured to be empty by analogy with
Shustin's result that the camel is left unrealized in degree one
less (i.e., $1+6$ corresponding to so-called {\it affine
sextics\/} when the line is interpreted as that lying at
infinity). One finds there (Smith 2005 \loccit) also a lovely
table of affine $M$-sextics known yet to exist. It would be
extremely interesting to adapt this table to degree $1+7$ even
after restricting focus to the comb case. Let us briefly explain
the classification of affine sextics with a maximally intersecting
oval cutting the line transversally across 6 points. First one can
as on Fig.\,\ref{Affine_Sextics:fig}a distinguish the following
configurations termed (by Arnold, Korchagin, ?) the {\it comb\/},
{\it snake\/}, {\it snail\/}, and {\it camel\/} respectively. This
shows all possible isotopic placement of such an oval w.r.t. a
line in $\RR P^2$. Then Fig.\,\ref{Affine_Sextics:fig}b represents
the more global algebraic problem of configuration actually known
to exist. Albeit the problem of classifying such affine sextics
was really started by Korchagin/Shustin (first independently and
then jointly) earlier workers made implicit contribution starting
with the 2 (split)-schemes of Harnack, Gudkov for 2 schemes, Viro
for 9 schemes, Korchagin for 16 new species, Shustin 4 species,
and Orevkov 1998 \cite{Orevkov_1998}, \cite{Orevkov_1998} to two
species (the second of which having been erroneously declared
prohibited in Shustin 1988
\cite{Shustin_1988-Toward-isotopy-class-affine-M-sextic}.).
Fig.\,b$^\ast$ is an attempt to mix the 2nd species of Gudkov (G2)
with V2 the 2nd species of Viro in order to build an interesting
curve of degree 8 (alas it fails seriously to be maximal for
evident reasons). Fig.\,c shows some configurations not
constructed which are either prohibited or perhaps some few which
are not yet known to exist. Especially noteworthy is the deep
collaborative obstruction of Le Touz\'e-Orevkov 2002
\cite{Fiedler-Le-Touzé-Orevkov_2002} where the black framed
species of Fig.\,d is prohibited. Despite of our poor
understanding of this problem, it seems that experts are fairly
close to a complete census of all affine sextics. On the table
above $Ha+G+V+K+S+O=2+2+9+16+4+2=35$ species are constructed.
Actually, prior to Orevkov's constructions we had 33 species
constructed in Korchagin-Shustin 1988
\cite{Korchagin-Shustin_1988/89} and other (according to Orevkov
1998 I \cite{Orevkov_1998} more detailed) constructions of these
33 (affine) curves were given in Korchagin 1996
\cite{Korchagin_1996-smoothing-6-fold}. Further good explanations
are to be found in Orevkov 1998I \cite{Orevkov_1998} where it is
explained that Korchagin-Shustin made some mistake at the
prohibitive level, and that actually Orevkov managed to prohibit
all species except the 33 constructed by Korchagin-Shustin (and
their forerunners like Harnack, Gudkov, Viro) and five species
corresponding to the symbols
$$
A_3(0,5,5), A_4(1,4,5), B_2(1,8,1), B_2(1,4,5), C_2(1,3,6).
$$
Orevkov's 1st note on affine $M$-sextic is devoted to proving
existence of $B_2(1,8,1)$ (which clearly corresponds to
Fig.\,\ref{Affine_Sextics:fig}b1). So it seems that $B$ (or $B_2$)
refers to the snake and the sequence $(1,8,1)$ to the number of
inner ovals when travelled along the cell along the natural sense.
Orevkov's 2nd note is devoted to the construction of $A_3(0,5,5)$
 which corresponds to our Fig.\,b2.
So it seems that the coding of Korchagin-Shustin means that $A$ is
the comb, $B$ is the snake, while the first two entries are the
number of inner ovals ordered by a natural convention of
contiguity between subregion of the cell (as split by the line)
while the last parameter describes the number of outer ovals
(which in principle is predestined so that the total sum is $10$).
The index is somewhat mysterious, but there is surely an
explanation. At any rate the next development is Le
Touz\'e-Orevkov 2002, where the species of Fig.\,d is prohibited.
This should correspond to the symbol, $B_? (1,4,5)$ with $?=2$ the
only choice possible from the above displayed list of five. So if
we decode correctly the symbolism, $A_4(1,4,5)$ corresponds to
Fig.\,d1 (hopefully) while $C_3(1,3,6)$ could be something like
our Fig.\,d2 (as C is probably standing for the snail). In
conclusion there is (in principle) 35 types constructed and only
two remains in doubt. As far as we know the problem did not
progressed anymore since Le Touz\'e-Orevkov 2002
\cite{Fiedler-Le-Touzé-Orevkov_2002} and so the situation is quite
reminiscent of Hilbert's 16th in degree $m=8$.

Further it seems (cf. e.g. an article by Polotovskii, maybe the
end notes of Polotovskii 2000
\cite{Polotovskii_2000-On-the-classif-decompos-7th-degree}) that
decomposing curve do have direct application to Hilbert's original
problem. So it is not impossible (browse also through Orevkov's
texts) that a resolution of the problem of the $35+2$ affine
sextics has some direct impact upon Hilbert in degree $m=8$,
abridged $H(8)$ (for say $M$-curves to simplify a bit).

\begin{figure}[h]\Figskip
%\vskip-1.2cm\penalty0
%\centering
\hskip-2.7cm\penalty0
\epsfig{figure=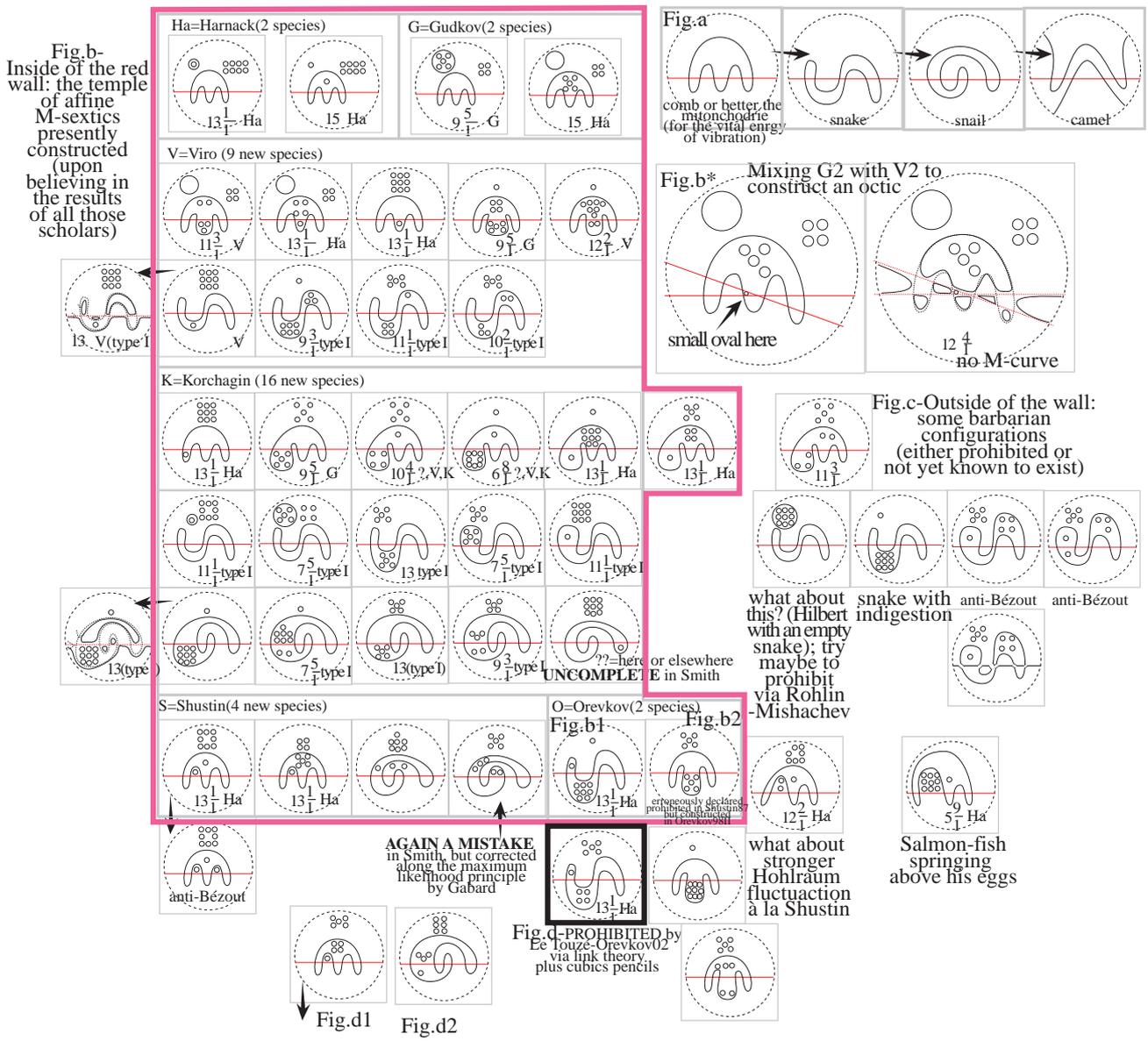,width=172mm} \captionskipAG
  \caption{\label{Affine_Sextics:fig}%
Picture of the fundamental sextic contortions and the (actual, not
definitive?) census of all affine sextics (Harnack, Gudkov,
Korchagin, Shustin, Orevkov). Shustin showed (indirect source via
Smith) that the camel is not realized algebraically (by sextics).}
\figskip
\end{figure}

If so, it would be interesting to know which realizabilty implies
automatically which schemes. So let us assume that one of the 2
bosonic affine $M$-sextic exist. Does this implies existence of a
new $M$-octics?

[18.07.13] Apparently to answer this question, one may look at
Korchagin 1996 \cite{Korchagin_1996-smoothing-6-fold} article
where affine $M$-sextic are used as patches  to be glued in 6-fold
ordinary singularity. For this to be properly understood it is
most convenient to represent an affine sextic in the fundamental
circle, yet instead as above rather in the fashion that the red
lines at infinity corresponds to the fundamental circle with the
usual antipodal identification. Once this is done we directly get
the required patches and so Korchagin 1996 (\loccit) is able to
construct new $M$-nonics starting from certain Viro's quintics
while applying to them quadratic transformation to get a highly
singular nonic to which the gluing method \`a la Viro is applied.

Yet our goal is not
%%%to be
getting sidetracked to the more ambitious case of nonics before
%%having closed
completing the chapter of $M$-octics. Hence in order for affine
$M$-sextics to be useful in the construction of $M$-octics it
seems that one should have a (global) singular octic with a 6-fold
point. Additionally as we impose $m=8$ there is room for a 2-fold
point (alias double point). So basically $8$ is split as $2+6$
(instead of $4+4$ like in the quadri-ellipse basic Viro method, or
as $5+3$ in Viro's more sophisticated avatar involving the 5-fold
singularity $N_{16}$ combined with the triple point $J_{10}$, see
also Shustin's variant involving two quadruple points). At any
rate from our present perspective we could try to use an affine
$M$-sextic as patch for dissipating the 6-fold point $M_{25}$ and
use additionally an ordinary double point $A_1$. To get started we
need only a global singular octic with this prescribed $M_{25},
A_{1}$ singular datum. As we do not feel comfortable with
Huyghens/Newton/Cremona's transformation/hyperbolism we shall
direct appeal to artistic creativity (i.e. free hand drawing) in
the hope that algebra is flexible enough to follow our freewill.

Recall a priori that we know with the last progresses by Orevkov
98/98II \cite{Orevkov_1998_II} 35 types of dissipation of $M_{25}$
so we can expect the method to be quite versatile, and supplying
new progresses in $H(8)$, i.e. Hilbert's 16th in degree $m=8$.

So let us start with the following curve
(Fig.\,\ref{GabardDEGREE8:fig}) obtained by a fairly random
connection between a double and sextuple point which we postulate
as existent as an octic. As we know an ordinary double point
diminish the genus by one, a  triple point counts---when
perturbed---like 3 ordinary points, while a quadruple point
generate $1+2+3=6$ ordinary points. Hence a sextuple points eats
$1+2+3+4+5=15$ units to the genus so our curve has genus
$g=(21-15-1)=5$ hence at most 6 real circuits by Klein's version
of Harnack. Alas it seems that our Fig.\,a was not optimally
chosen because we already consumed 4 circuit for the singular
locus. After some trial we eventually arrive at Fig.\,d with a
single circuit. Alas this curve violates B\'ezout for line yet
maybe this is just a topological model of an isotopy type
admitting also a model satisfying B\'ezout. If not one should try
further maybe by rotating the petals as to be obturating the
vision between both singular points (compare e.g. Fig.\,25 of
Korchagin 1996 \cite{Korchagin_1996-smoothing-6-fold} and then our
idea should be meaningful). This inspired us Fig.\,e which looks
more B\'ezout compatible. Of course as the genus is $g=5$ we have
6 circuits so one can add 5 additional ovals. Presently we do not
know were yet we shall try random locations while trying to fit to
the experimental data already available.

\begin{figure}[h]\Figskip
%\vskip-1.2cm\penalty0
%\centering
\hskip-2.7cm\penalty0
\epsfig{figure=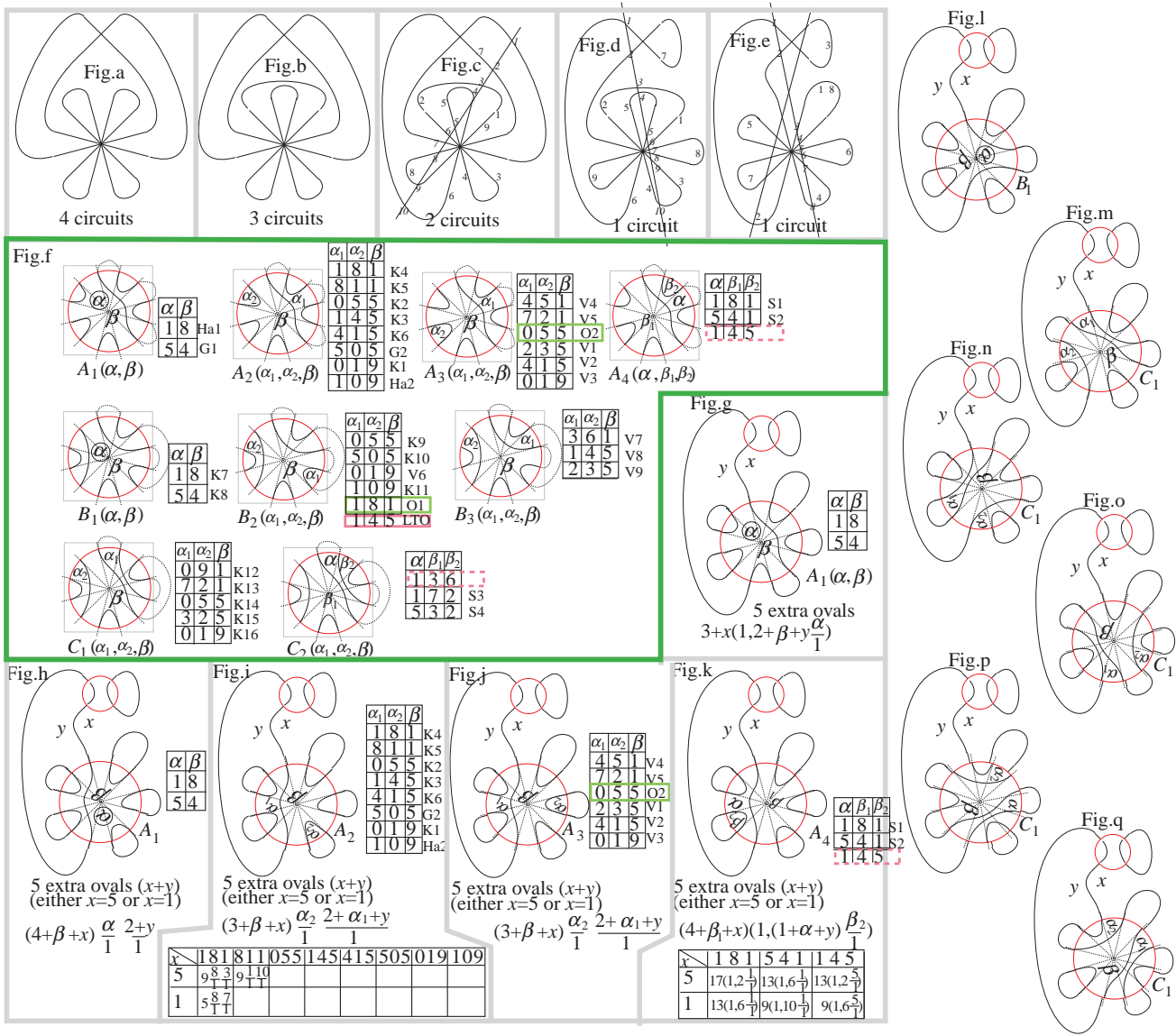,width=172mm} \captionskipAG
  \caption{\label{GabardDEGREE8:fig}%
A freehand singular octic and Korchagin's list of affine
$M$-sextics (in absolute ``circular'' representation)} \figskip
\end{figure}

\subsection{Smoothing $M_{25}$ (=$O_6$=ordinary sextuple point):
after Korchagin, Orevkov, Le Touz\'e}

In order to dissipate the singularity $M_{25}$ (sextuple point) of
Fig.\,d we merely need to repeat the quasi-census of singularity
$M_{25}$ as presented say in Korchagin 1996
\cite{Korchagin_1996-smoothing-6-fold}. This is the same as the
one depicted above except that  the line is represented as the
absolute circle at infinity. So we merely have to copy Korchagin's
series of Figs.\,3 given on p.\,143--144 of \loccit to get our
Fig.\,f. Maybe we just take the initiative for gluing convenience
to kill the square depiction of Korchagin which looks anyway
anecdotic. Korchagin's table turns out to be apparently perfectly
correct (no misprints and permit thereby to correct the 2 obvious
errors that were evident in Smith's table). Of course we just
added to Korchagin's list both $M$-sextics of Orevkov (cf. the
green framed one, namely $B_2(1,8,1)$ and $A_3(0,5,5)$). Finally
we added also the species $B_2(1,4,5)$, albeit prohibited in Le
Touz\'e-Orevkov 2002 \cite{Fiedler-Le-Touzé-Orevkov_2002} since it
is fairly probable that there proof is too involved to be solid.
At least it is always interesting to keep track of such a strange
bird to see what result from it in degree $m=8$. Finally we added
also to Korchagin's list the 2 bosons $A_4(1,4,5)$ and
$C_2(1,3,6)$ (as red dotted frames), again the motivation being to
contemplate what is generated from them in the realm of octics.

\begin{figure}[h]\Figskip
%\vskip-1.2cm\penalty0
%\centering
\hskip-2.7cm\penalty0
\epsfig{figure=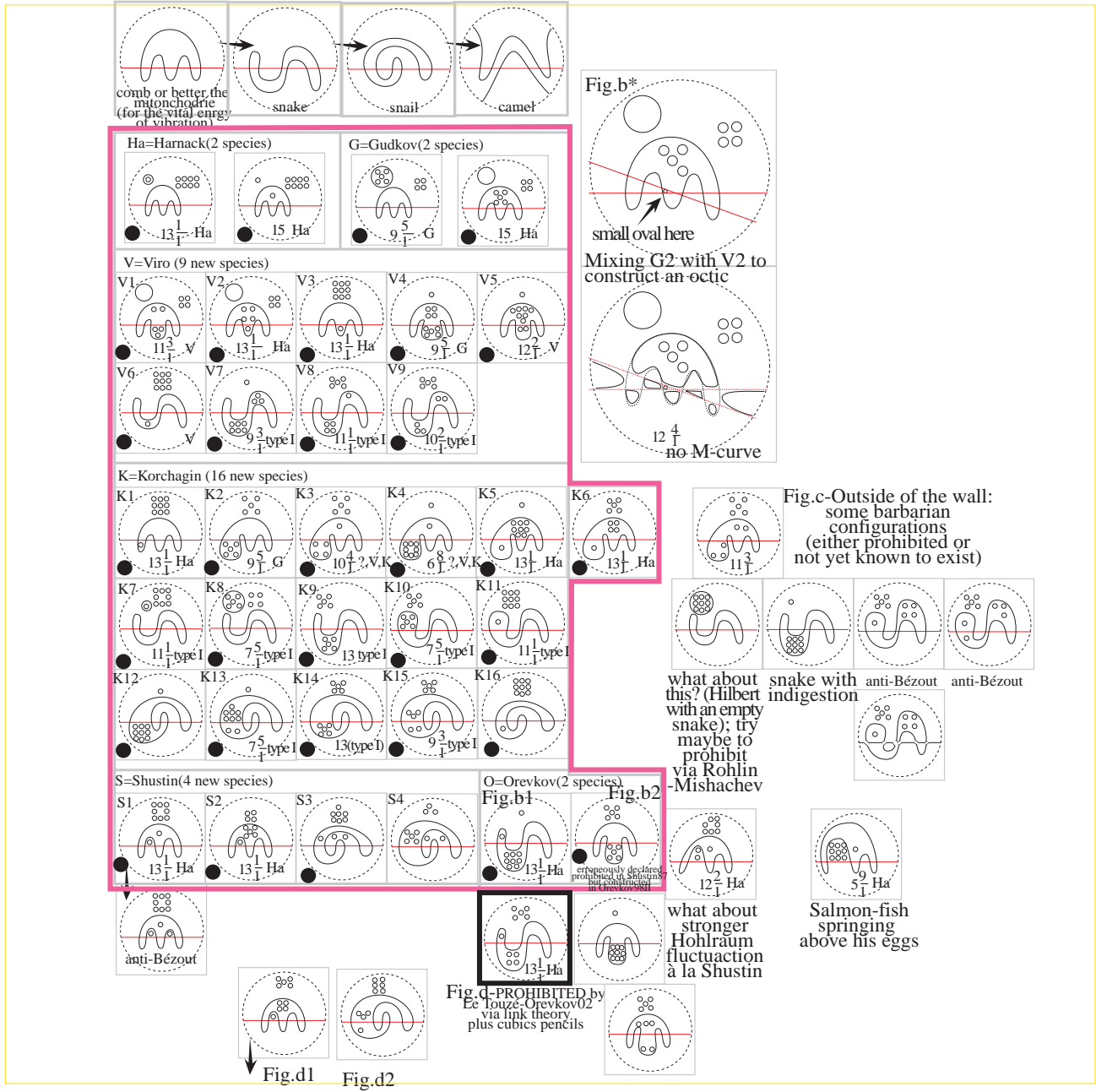,width=172mm} \captionskipAG
  \caption{\label{GabardDEGREE8B:fig}%
 Korchagin's list of affine
$M$-sextics in ``linear'' representation} \figskip
\end{figure}

This being said we are now pared to glue the 35 (known)
dissipations in our floppy singular octic. By choosing the $A_1$
dissipation we get (Fig.\,g) with scheme $3+x (1,2+\be+y
\frac{\al}{1})$. Here $x+y=5$ represent the 5 extra ovals which
according to B\'ezout can only be located as indicated. So the
total number of oval is $7+9+5=21$ only. This is because we did
not choose the optimal dissipation. It suffices however to rotate
the patch to get a maximal formation of ovals (Fig.\,h). The
corresponding scheme has Gudkov symbol $(4+\be+x) \frac{\al}{1}
\frac{2+y}{1}$. For the first listed smoothing $(\al,\be )=(1,8)$,
this specializes to $12+x \frac{1}{1} \frac{2+y}{1}$. Looking at
the table of $M$-octic (Fig.\,\ref{SIMPLIFIED-TABLE:fig}) we see
that Gudkov periodicity forces choosing either $x=5$ or $x=1$, but
the resulting schemes are all well known (i.e. the first 2 colons
of the table). However now the parameter $x$ is nearly fixed and
we may explore other smoothings. For $A_2$ we get Fig.\,i but alas
the first coefficient $3+\be + x$ is at least 5 when $x=1$ and
$\be$ takes its lowest value $1$. The smoothing A3 leads to
Fig.\,j, with the same Gudkov symbol, hence the same defect of
having the 1st coefficient $\ge 5$, as opposed to having it equal
to $1$ as to land in the bosonic zone where little is known (apart
4 construction by Viro and 2 prohibitions by Orevkov, cf.
Fig.\,\ref{SIMPLIFIED-TABLE:fig}). Finally the smoothing A4 leads
to Fig.\,k and even appealing to the exotic dissipation $(1,4,5)$
merely produces Chevallier's scheme $13(1,2\frac{5}{1})$ when
$x=5$, and when $x=1$ we merely get schemes well-known since
Viro's most basic method involving the quadri-ellipses.

So the game is a bit disappointing but it remains of course to
exploit the other dissipation not of comb-type (i.e. class $B$ and
$C$ being respectively the snake and the snail). A priori our
intuition is that those guys will not produce $M$-curves, and this
seems corroborated by Fig.\,l,m,n,o,p,q.

So it seems that we need more artistic imagination when tracing
the ground singular octic. So reminding our former Fig.\,e as
Fig.\,\ref{GabardDEGREE8_2:fig}a, we had first the idea (in order
to kill outer ovals) of looking at Fig.\,b. This seems alas to
have 2 defects. First the curve in question has already 2 circuits
and further the red line exposes the octic to 10 intersections.
Next the morphogenetic brain thinks about Fig.\,c but this does
not satisfy the desideratum of only one outer oval under a
maximally closing dissipation. So it seems that the idea is that
the large circuit has to enclose as many petal possible. Further
still under the desideratum of minimizing the number of outer oval
it is evident that it is more clever if the loop self-connecting
the node is traced introverted as on Fig.\,d, but the latter has
unfortunately 2 circuits. This suggested next Fig.\,e, alas also
with 2 circuits. Then Fig.\,f is lovely for having only 1 circuit
but its maximal closing will have at least 2 outer ovals coming
from the external petals. So the more radical choice is to go to
Fig.\,g but the price to pay is the presence of 2 circuits
(impeding us to add the maximum number of five extra ovals).
Finally after numerous trials we arrived at Fig.\,z. On applying
the dissipation $A_1$ we get an octic violating B\'ezout as it has
2 subnests (Fig.\,z1). Of course we could choose another
dissipation like $A_2$, yet the paradigm of the independence of
smoothing leads rather one to seek an universal object viable
under all possible dissipations. At any rate also using $A_2$
leads to 2 subnests (Fig.\,z2), while rotating more the patch like
on Fig.\,z3, we loose one oval.

\begin{figure}[h]\Figskip
%\vskip-1.2cm\penalty0
%\centering
\hskip-2.7cm\penalty0
\epsfig{figure=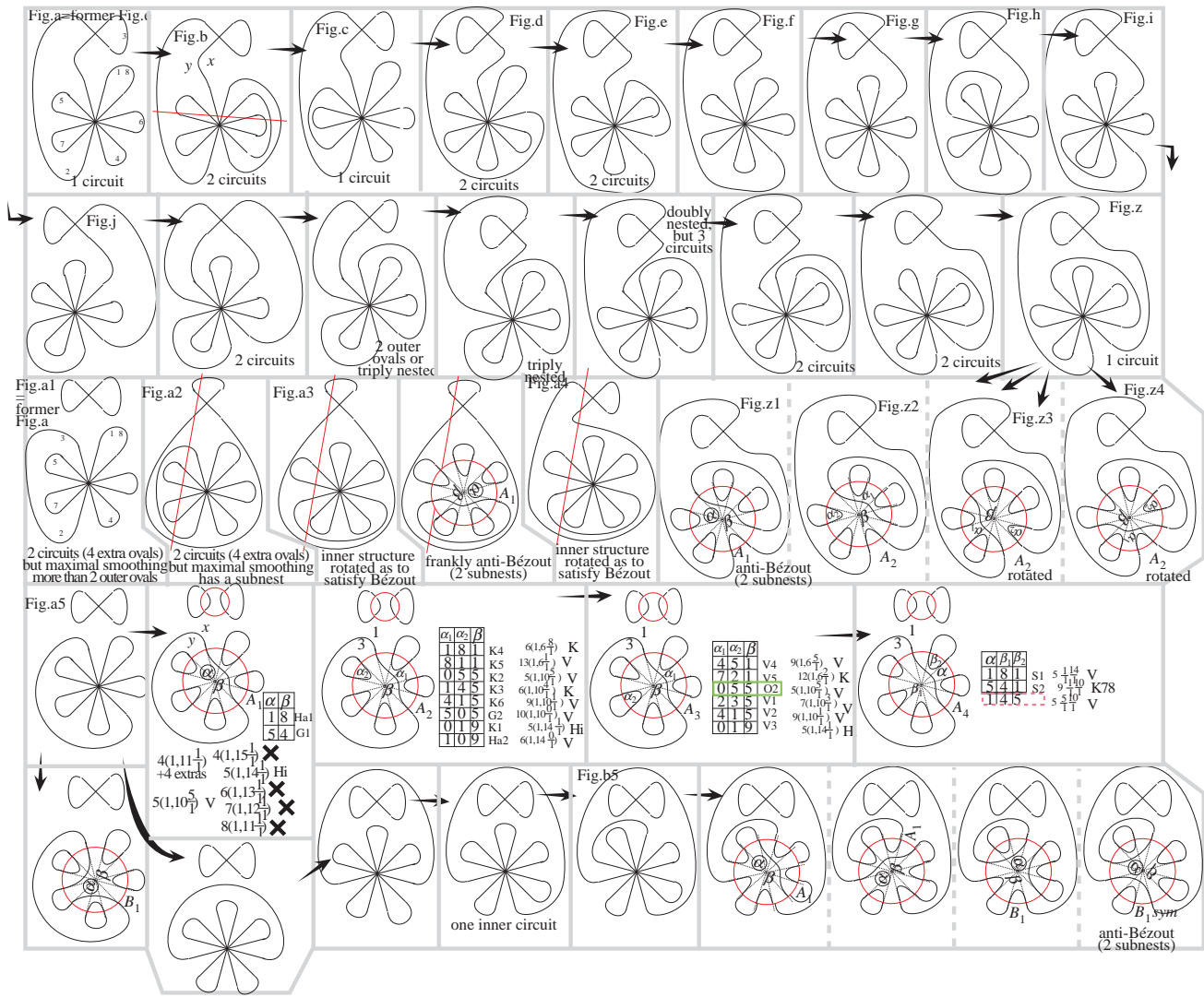,width=172mm} \captionskipAG
  \caption{\label{GabardDEGREE8_2:fig}%
A new freehand browsing through singular octics with smoothings}
\figskip
\end{figure}

[19.07.13] Then we transformed Fig.\,a to Fig.\,a1, and studying
variants arrived at Fig.\,a5 and considered its smoothings. A
first little surprise is that if such a curve exist (parameters
$x,y$ of additional ovals went adjusted as to respect Gudkov
periodicity) then using the $A_3$-smoothing $(7,2,1)$ due to Viro
one gets already the Korchagin's scheme $12(1,6\frac{2}{1})$. In
contrast using the  (semi?) highbrow (sextic) patch $A_3(0,5,5)$
of Orevkov one only recover the basic (quadri-ellipse type) Viro
curve $5(1,10\frac{5}{1})$. Using the smoothing $A_4$ should have
conducted to the most spectacular result especially when using the
hypothetical dissipation $A_4(1,4,5)$, but alas the latter only
produce a basic (quadri-elliptic) Viro scheme, namely
$5\frac{5}{1}\frac{10}{1}$. Of course the other smoothings (snake
and snail like) will not produce $M$-curves (because the
corresponding patches do not have  3 consecutive bumps).

So we need a new curve and we designed Fig.\,b5 after the
evolution law given by the arrow starting for Fig.\,a5. Now the
maximal smoothing does not come from the comb (type $A$) but from
the snake family $B$. Alas the curve so obtained under the
$B_1$-smoothing (symmetrized) violates B\'ezout.

Another idea is to replace the (ordinary) double point by a
solitary double points. And so we consider Fig.\,f1 where the
black dot is the solitary node (the trick being that we placed it
inside of a petal in order to kill an outer oval) so as to land in
the critical bosonic strip where very little is known. As before
the genus of this singular octic is $21-(1+2+3+4+5)-1=5$ so that 6
real circuit are permissible (and 5 of them are not yet traced).
Further in view of the global pattern of the curve it is clear
that the $A$-type (comb) smoothing produces the maximal
$M$-smoothings. Of course according to the independency of
smoothing we choose the solitary node to deform to a little oval,
yet Fig.\,f1A1 shows a scheme with a subnest and an outer nest so
that B\'ezout is contradicted. Hence the configuration Fig.\,f1
cannot exists. (Incidentally one can wonder how difficult it would
be to draw this conclusion if Viro's theory was not available!) So
this is a deception, yet let us surf however the principle of
independency (albeit it might have been proven in this context by
Shustin, cf. e.g. Shustin 1987 \cite{Shustin_1987-versal-deform}
(versal deformation paper)) to look at the more interesting
smoothing $A_2$. This gives Fig.f1A2, where again B\'ezout is
generically foiled except when $\al_2=0$, but then alas $\al_1$ is
not zero so that we get $\ge 2$ outer ovals (and not just one as
desired). After a long sequence of trials we arrived at Fig.\,q1,
whose smoothings along $A_2$ and $A_3$ includes schemes violating
Viro's imparity law. Further, as evident from the scratch, the
smoothing A4 violates B\'ezout.

\begin{figure}[h]\Figskip
%\vskip-1.2cm\penalty0
%\centering
\hskip-2.7cm\penalty0
\epsfig{figure=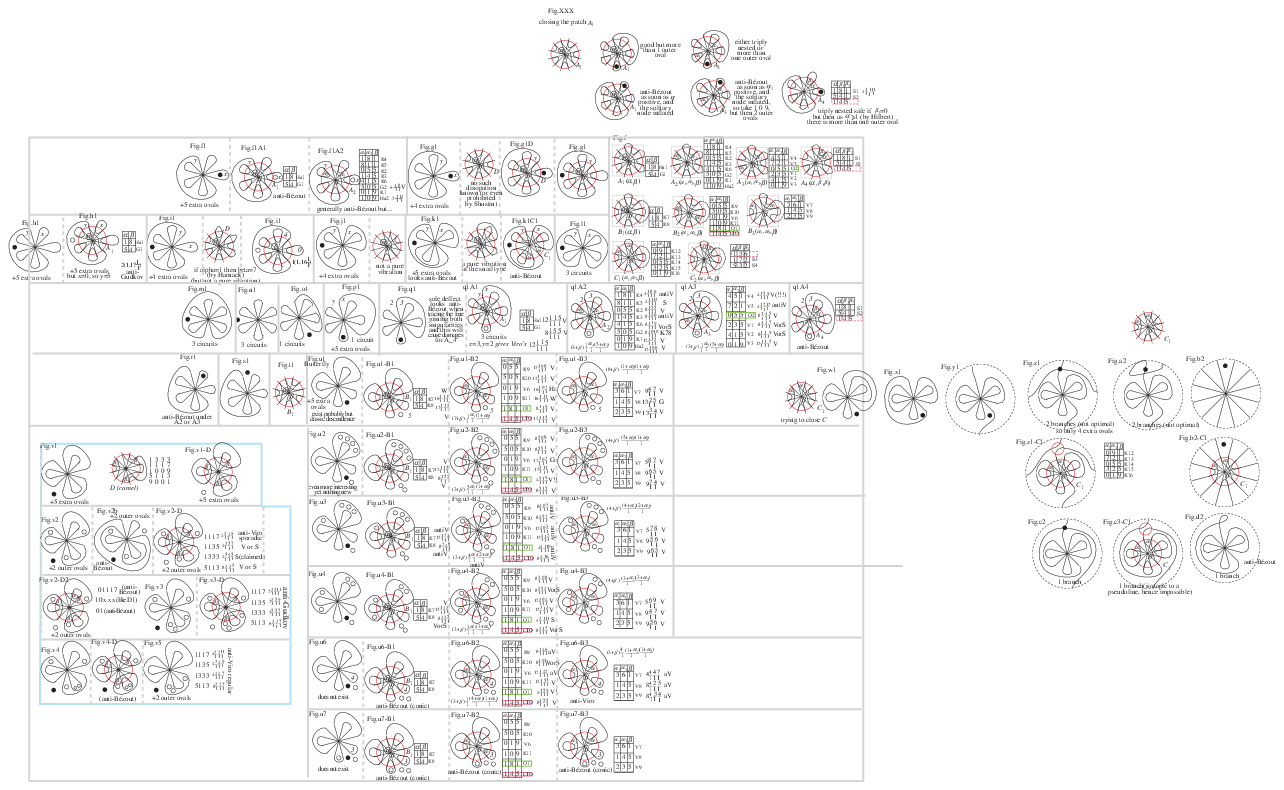,width=172mm} \captionskipAG
  \caption{\label{GabardDEGREE8_3:fig}%
Browsing more through singular octics with a solitary node}
\figskip
\end{figure}

Another idea is to construct the global curve around some fixed
smoothing. So we get from $B_1$ (reproduced as Fig.\,t1) the
Fig.\,u1 where the solitary node is placed so that the line
joining it to the heavy singularity never meet the curve anymore.
Smoothing along $B_1$ gives Fig.\,u1B1 which leads primarily to
the symbol $\frac{1}{1}\frac{1}{1}\frac{1}{1}$, which among
$M$-schemes can only be completed as Wiman's scheme
$16\frac{1}{1}\frac{1}{1}\frac{1}{1}$ so that we can adjust the
number of outer ovals as being 5. So Fig.\,u1
 must be mentally refreshed with 5 outer ovals, and then we can
 explore the full kaleidoscope of all possible smoothings of this
(hypothetical) singular octic. On Fig.\,u1-B2 one sees that
Orevkov's fairly recent smoothing only produce a (fairly)
well-known Viro's scheme (which I could only understood by
Shustin), whilst the smoothing in principle forbidden by Le
Touz\'e-Orevkov produces only the ultra-classical  Viro scheme
$12\frac{1}{1}\frac{1}{1}\frac{5}{1}$. (So maybe the Le
Touz\'e-Orevkov affine-sextic obstruction is wrong!?) Finally
Fig.\,u1-B3 produces only super classical schemes (of G=Gudkov and
V=Viro).

{\it Added} [21.07.13] We first missed to consider as well
singular octics deduced by quantum transfer of 4 ovals. This
produces e.g. Fig.\,u2 where the number of outer ovals is just
diminished by 4. Then there is also Fig.\,u3 which alas runs into
trouble with Viro's imparity law. So if the latter is true an
octic like Fig.\,3 is excluded. Alternatively it is strange that
Fig.\,u3-B3 leads only to viable schemes (first constructed by
Viro) and this phenomenology could imply (if our singular $C_8$
really existed) that there is no independency in the smoothings
(so violating a principle/theorem of Viro/Shustin). Then we can
move to Fig.\,u4, where interestingly Orevkov's smoothing
$B_2(1,8,1)$ yields Shustin's (fairly original) scheme
$4\frac{3}{1}\frac{11}{1}$. Further Le Touz\'e-Orevkov presupposed
anti-smoothing $B_2(1,4,5)$ realizes a even more standard scheme
of Viro (or Shustin). So no obstruction \`a la Le Touz\'e-Orevkov
is detected by patchworking on this curve (Fig.\,u4). Since the
patch is symmetric about a line angled 120 degrees there is no
need to examine the next configuration u5, where 1 oval lye in the
first petal and 3 in the 2nd petal. Finally we could transfer
ovals in the simple petals (e.g., like on Fig.\,u6) but then the
smoothing $B_1$ violates B\'ezout for conics (equivalently the
maximality of the doubled quadrifolium). However it is still
interesting to look at the other smoothing which are less
blatantly foiled (contradicting only with Viro's imparity law).
This little corruption of Viro can be tamed by transferring one
oval in the other petal like on Fig.\,u7, but then each
dissipation foils B\'ezout for conics (due to the presence of 4
nests). It seems evident at this stage that we have exhausted the
possibilities allied to the butterfly type.

So we see that one more efficient method consist to start from the
dissipation and find a global closing of it. Hence of starting of
the curve we suit the curve to the smoothing. So starting from the
A-smoothing we get Fig.\, XXX but alas we fail so to reach the
doubly-nested case with one outer oval.

Of course for reason of symmetry Fig.\,v1 suggests by itself, yet
alas it is not calibrated on a smoothing. Hence starting from the
$C$ smoothing (in any of its incarnation $C_1$ or $C_2$ which have
the same ground skeleton) we may construct the most natural
closing as Fig.\,w1 but it is evident that the latter curve
corrupt B\'ezout (as is especially visible after the smoothing as
we have two nests of respective depths 3 and 2 so forcing 10
intersections with a suitable line chosen to pass through the
deepest ovals). Is there a more clever way to globalize this
infinitesimal smoothing? Of course there is the variant shown as
Fig.\,x1, yet this is is subjected to the same objection. Of
course it could as well globalize as Fig.\,y1, yet this looks an
illegal art-form as it involves a pseudoline (alias branch of odd
degree) not so realizable by an octic curve. The issue seems to be
to look at Fig.\,z1 by introducing again a non-isolated ordinary
node so that another odd branch restore the parity of the degree
(homology class). Alas we found nothing extremely convincing along
closing the snail by a global octic.

[20.07.13] Then coming back to Fig.\,v1 whose maximal dissipation
is depicted right below and clearly identified to the camel type
(which in principle is obstructed by Shustin as we learned in
Smith 200X \cite{Smith_2005-Thesis-under-Korchagin-split-7+1}, see
also Korchagin 1996 \cite{Korchagin_1996-smoothing-6-fold}).
Improvising a short list of dissipation especially the type
$D_1(1,3,3,3)$ which corresponds to a sextic of Harnack type
$9\frac{1}{1}$ would produce the $M$-octic
$7\frac{4}{1}\frac{4}{1}\frac{4}{1}$ which obviously violates
Viro's imparity law, and even more radically Gudkov periodicity
(even in the simplest Arnold version mod 4). So what to conclude?
Either the smoothing $D(1,3,3,3)$ does not exist (as asserted e.g.
in Korchagin 1996) or maybe even the global curve of Fig.\,v1 does
not exist. What about trying another smoothing like $D(1,0,0,9)$.
Patching this into Fig.\,v1 gives the scheme
$7\frac{1}{1}\frac{1}{1}\frac{10}{1}$. Again Gudkov periodicity is
foiled as the first outer coefficient of the symbol as to be a
multiple of four (compare Fig.\,\ref{SIMPLIFIED-TABLE:fig}).

However tacit in our depiction of Fig.\,v1 was the assumption that
the 5 extra ovals are outer and we merely need to delocalize 3
inside the singular circuit to get the first (outer) coefficient
reduced from $7$ to $4$. This suggests Fig.\,v2 which smooth under
$D(1,3,3,3)$ to $4\frac{5}{1}\frac{5}{1}\frac{5}{1}$ a well-known
scheme due to Shustin. So let us pose as an Ansatz that the curve
depicted as v2 exists and that the dissipation $D(1,3,3,3)$ exists
too (albeit this seems to violate a result of Shustin). We would
like to speculate about further smoothing of the camel type yet it
suffice less imaginatively to explore the other more classical
smoothing (but those are nost best suited so fail to produce
$M$-schemes). So we are condemned to adventure in the (in
principle) deserted type of the camel (which according to Smith's
interpretation of a Shustin work does not exist at all). So the
question is: is the Sahara so deserted as to support no camel?

Let us again improvise a list of dissipation right below Fig.\,v1
involving the parameters $(1,3,3,3)$, $(1,1,1,7)$, $(1,1,3,5)$,
etc. The sole condition is that the $\be_i$ have to be odd as to
respect Viro's imparity law and we restrict first attention to the
case of a Harnack curve with one inner oval so $\al=1$. The
complete enumeration of such smoothing is therefore in
lexicographical order $(1,1,1,7)$, $(1,1,3,5)$, $(1,3,3,3)$, and
nothing more. If we move to Gudkov's type of an affine sextic we
get additionally $(5,1,1,3)$, and the Hilbert's sextic leads to
$(9,?,?,?)$ nothing compatible with Viro as the question marks are
odd hence $\ge 1$. Supposing that all these 4 smoothings exists
their injection as gluing into the hypothetical curve v2 (which
looks fairly reasonable and especially esthetic) produces the
schemes listed on Fig.\,v2-D which are all either due to Viro or
Shustin, safe the scheme $4\frac{3}{1}\frac{3}{1}\frac{9}{1}$
which is declared prohibited by Viro's sporadic obstruction.

Of course it remains then also to explore other smoothing of the
camel type with the same ground picture yet with different
locations for the micro-ovals. So have a second type of smoothing
$D_2$ tabulated on Fig.\,v2-D2, but it seems that all of them will
just violate B\'ezout (at least granting existence of our ground
curve Fig.\,v2). So it seems that there is only one class of
dissipation of the camel type (namely $D_1$ abridged as $D$). Then
one can speculate about a curve like Fig.\,v3 yet its smoothing
are anti-Gudkov (see Fig.\,v3-D). However another trading keeping
$\chi$ constant, is to drag the two outer ovals in the one same
nest so as to get Fig.\,v4, but the B\'ezout explodes as seen on
Fig.\,v4-D (if not already apparent on Fig.\,v4 already). Of
course we can also imagine Fig.\,v5, but the resulting scheme are
anti-Viro (imparity law).

Looking back to Fig.\,v1-D one's desideratum could be to get
Wiman's curve $16\frac{1}{1}\frac{1}{1}\frac{1}{1}$ which is more
in the bottom levels of the 2nd pyramid. Of course this would be
possible if we could set all three $\be_i=0$ and $\al=10$. Looking
however at the camel patch it is clear that the tripodal amoebic
region occupied by $\al$ extend through infinity by small
semi-disk, so that this expansion is a cell and therefore the
interior of the traced oval consisting of the 6 branches all
connected together through infinity. Therefore the choice $\al=10$
corresponds to the sextic $\frac{10}{1}$ prohibited by Rohn (or
just via Arnold-Rohlin). Instead a smoothing with parameters
$(\al,\be_1,\be_2,\be_3)=(9,0,0,1)$ is permissible (a Hilbert's
sextic), yet when glued in Fig.\,v1-D would violate Viro's
imparity law, yielding the scheme $15
\frac{1}{1}\frac{1}{1}\frac{2}{1}$ which actually violates Gudkov
periodicity (even in the simple version of Arnold modulo 4). Of
course we could repair it by moving outer ovals inside, but as the
patch is rotationally mobile, this forces one transfer into each
three ovals and so we get again Fig.\,v2 which under smoothing
$D(9,0,0,1)$ gives $12\frac{2}{1}\frac{2}{1}\frac{3}{1}$, which
respects Gudkov but violates Viro (imparity law).

Albeit fairly complicate to get through all this, we see that
there is an intricate interaction between affine sextics (and the
allied infinitesimal dissipation), global singular octics with 2
singularities of type $A_1$ (double point), $M_{25}$ (sextuple
point) and global smooth $M$-octics. Of course for this
interaction to take place it is vital to have a result of
independency of smoothing singularities (that is in principle by
an extension of the usual Riemann-Roch-Severi-Brusotti argument)
to higher singularities (work of Gudkov-Viro-Shustin). A s a
concrete example let us state the following:

\begin{lemma}
Assume that the singular octic of Fig.\,v2 to exist and that the
dissipation $D(1,1,1,7)$ exist as well. Then there is an $M$-octic
with scheme $4\frac{3}{1}\frac{3}{1}\frac{9}{1}$. However the
latter is forbidden by a Viro sporadic obstruction (alas not
extremely well published). Hence either the latter Viro's result
is false, or at least one  of our two hypotheses is erroneous.
\end{lemma}

Of course one can try to attack directly the Fig.\,v2 by arguing
that the line through both singularities has exceeding
intersection, since this line has already $2+6=8$ multiplicity of
intersection and there seems no possible location for the node
from where the fundamental line (through both singularities) could
avoid the rest of the curve.

All this is fairly pleasant yet it seems to sidetrack us in the
realm of trinested $M$-curves (which in principle in totally
settled by the sporadic obstructions of Viro and constructions of
Shustin). In contrast the real golden El~Dorado
%%%%%checked in DICO
involves doubly nested and sub-nested schemes where there is
respectively 4 and 2 bosons still awaiting for some realization
(resp. prohibitions). As yet it seems fairly difficult to get to
doubly nested schemes with only 1 outer ovals with our method
involving a double and sextuple point. Likewise we could not find
much boson in the subnested family.

Another idea is to use $8=3+3+2$ and so to construct an octic by
dissipating 2 triple point and one double point. Of course one can
also exploit $4+2+2$ and look what happens there.

\subsection{8=1+7}

Finally via $8=7+1$ has to be interpreted as smoothing a septuple
point plus a simple (smooth) point. Alas this seems to require a
good understanding of  affine $M$-septics (a complete
understanding looks a bit out of reach for the present days, as we
do not have even settled the degree $6$ case). Yet we can still
try to explore the possibility in a floppy (and sloppy) way so as
to look if there is there some maneuvering room to realize the
bosonic schemes.

Let us be more concrete. Start with a septuple point
(Fig.\,\ref{GabardDEGREE8_septuple:fig}a). Imagine a certain
$M$-smoothing (like Fig.\,1) and imagine a global singular octic
with a unique septuple point (Fig.\,c). It may be noted that the
curve is therefore rational (pencil of lines through the
singularity). Alternatively one can use the genus, and note that a
septuple point diminish the genus by $1+2+3+4+5+6=21$ (imagine a
perturbation of Fig.\,a into a generic arrangement of lines and
count the resulting ordinary nodes). So Fig.\,c is complete, i.e
there is no additional circuits. Gluing the smoothing Fig.\,1 into
Fig.\,c produces Fig.\,c1
 which is $M$-curve, yet one violating Gudkov periodicity even
 mod 4 version \`a la Arnold.
However if we consider the smoothing of Fig.\,2, we get Fig.\,c2,
i.e. the bosonic scheme $1\frac{1}{1}\frac{18}{1}$ not yet known
to exist. Of course this proves little but shows at least that the
$1+7$ method looks better suited (than its $2+6$ companion) to
reach the bosonic strip of doubly-nested schemes with one outer
oval. Of course one deffect of our argument is that our nenuphar
like curve (\F c) does not look a bona fide octics for most lines
through the singularity cut it in nine points. This can be
arranged by looking at \F d. Before doing this let us note that
the affine septics depicted (as \F 2') when smoothed yields the
scheme $18(1, 1\frac{1}{1})$ which violates Gudkov periodicity.

One is next faced to the ingrate duty of listing all dissipations
of a septuple point (at least those of the class A having a
prescribed oscillating pseudoline). This task is difficult and
long yet still manageable at least in qualitative substance.
Bypassing the exact details for the moment it seems that it can
already be inferred (from abstract thinking) that the resulting
scheme (after gluing) will always possess at least 4 outer ovals
for we can only fill 3 ovals (with ovals) without contradicting
B\'ezout for conics. So we will certainly not reach the bosonic
strip by such a construction. However there is perhaps some chance
getting new existence results (constructions) for the 2 subnested
bosons. As one of these bosons is $14(1,2\frac{4}{1})$ it suffices
to consider the dissipation of Fig.\,3 and to glue in Fig.\,d to
get Fig.\,d3 which is the required bosonic scheme. Of course we
notice that we could have directly smoothed out the affine septic
model to get via Brusotti the same scheme. Extracting the exact
philosophy behind  this phenomenology, it seems that the
dissipation of the margarita flower (Fig.\,d) are in bijection
with decomposing curves of degree $1+7$ with a simple slaloming
pseudoline, compare our former
Fig.\,\ref{ViroDEGREE8_GABARD_1+7_GOOD:fig}.) Now the game could
be to reach the other subnested boson $4(1,2\frac{14}{1})$ and
this suggests killing some petals of the flower like on Fig.\,e.
Alas then there is 3 circuits (too much for a rational (unicursal)
curve). Further a line through an inner petal has 3 intersection
outside the septuple point leading hence to a total of 10
intersections (anti-B\'ezout for line alias Euclid? since it the
Euclidean algorithm for polynomials that bound the number of roots
via the degree).

[22.07.13] Next we can at look at Fig.\,f which has the same
defect. Notwithstanding we still consider a suitable smoothing
(Fig.\,f5) which is Wiman's scheme $16
\frac{1}{1}\frac{1}{1}\frac{1}{1}$. The corresponding affine
septics is Fig.\,f6, whose direct smoothing (Fig.\,f7) is only an
$(M-1)$-curve (actually one lying right below Harnack's scheme
$17\frac{1}{1}\frac{2}{1}$). Then Fig.\,h looks ideally suited to
bring us in the bosonic strip of doubly-nested schemes. As shown
by Figs.\,h6--h10 it is clear that we could sweep out all the
bosonic strip (safe for omitting its first item
$1\frac{1}{1}\frac{18}{1}$). This would be a spectacular advance
nearly closing the completion of Hilbert's 16th for $m=8$. Alas
the sole problem seems to be that there is no algebraic model for
Fig.\,h as the latter corrupts B\'ezout upon tracing the singular
line through an inner petal. Via Fig.\,i we can also get the boson
$1\frac{1}{1}\frac{18}{1}$ as shown by Fig.\,i11, but of course
our singular octic is still defective w.r.t. B\'ezout.

\begin{figure}[h]\Figskip
%\vskip-1.2cm\penalty0
%\centering
\hskip-2.7cm\penalty0
\epsfig{figure=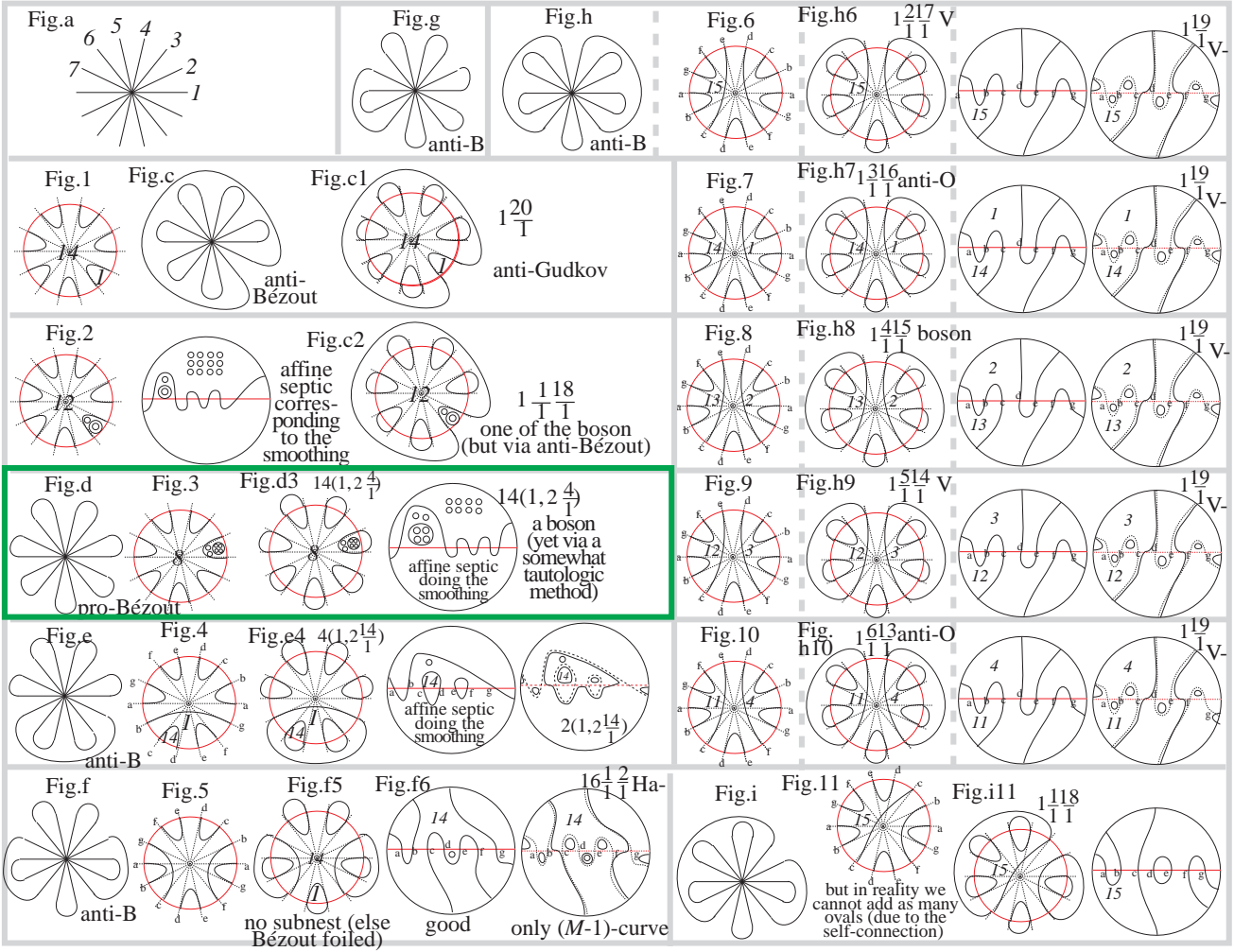,width=172mm}
\captionskipAG
  \caption{\label{GabardDEGREE8_septuple:fig}%
Singular octics with a septuple point with irregular constructions
of all the bosons (and more like schemes prohibited by academician
%%%Stepan
Orevkov)} \figskip
\end{figure}

All this Fig.\,\ref{GabardDEGREE8_septuple:fig} is fairly {\it
avant-gardiste\/}
%%%(vangardish (?)/ vintage)
as it produces all the 6 bosons and even the 2 schemes prohibited
by Orevkov, yet it must be confessed that our method is highly
irregular primarily because the ground singular octics of
Figs.\,c,e,h corrupt B\'ezout for lines. However Fig.\,d is
B\'ezout-regular and still contains a gluing which is bosonic. Of
course it seems however that the margarita curve (Fig.\,d) really
amounts to think directly about the corresponding affine septic.
Hence no real bonus is gained through the power-flower method. To
paraphrase a bit, it seems clear that we have the following
principle:

\begin{Scholium}
The class of $M$-octics occurring through small perturbation of
the margarite curve (Fig.\,d) via
%the
Viro's gluing method coincides exactly with those of decomposing
curves of degree $1+7$ via the Harnack-Brusotti method of
dissipation of ordinary nodes. Those schemes (albeit not being
exactly known) have a range contained in the blue zone of
Fig.\,\ref{SIMPLIFIED-TABLE:fig}. In particular this region
contains the boson $14(1,2\frac{4}{1})$, but not much more.
\end{Scholium}

\subsection{8=2+2+2+2, 8=3+5}

Now we could as well inspect a pair of quadruple ordinary nodes,
and then likewise consider a quintuple plus a triple point. Let us
first look at Fig.\,d as a combination of the quintuple and triple
point. On gluing with maximally closing dissipation like Fig.\,d1
we get Fig.\,d2. The idea is then of course to kill one outer
ovals by injecting micro-ovals inside of it. Alas doing this as on
Fig.\,d3 leads to an affine quintic violating B\'ezout. So it
seems that our sole chance is that the killing of the outer oval
is produced by a global oval of Fig.\,d. Let us calculate the
genus of the curve $C_8$ of Fig.\,d. A quintuple point is
tantamount to $1+2+3+4=10$ ordinary nodes while the triple point
contributes to 3 such nodes. Accordingly, the genus of $C_8$ is
$21-10-3=8$ so that (by Harnack's bound $r\le g+1$) there is 8
extra ovals possible on Fig.\,d. Then we consider the smoothing
maximally closing the ovals of Fig.\,d5, yet since the main
circuit of the patch has already 3 components we can (for a
quintic) only add 4 ovals, and we ascertain with deception that
the resulting gluing (Fig.\,d6) is not maximal. Then we had the
idea to look at Fig.\,e (but the line through the singularities
has excess intersection). This may be settled by Fig.\,e, but then
there are too many free petals inducing many outer ovals impeding
a safe landing in the bosonic strip where there is only one outer
oval. Fig.\,f looks better yet 2 circuits. Fig.\,h is deduced by
creating more petals and looks perfect for having both one circuit
and for respecting B\'ezout. Alas the fairly standard dissipation
of Fig.\,h2 (involving basically an affine quintic of Harnack
type, cf. Fig.\,h3 for the resulting $M$-sextic of Harnack type)
leads after patchworking to an octic violating seriously Gudkov
periodicity (even in the weaker Arnold version mod 4). Our next
idea was to create the second nest form the  8 extra ovals
available suggesting thereby Fig.\,i which smooth indeed to
Fig.\,i3 realizing the boson $1\frac{1}{1}\frac{18}{1}$. At this
stage we could be very happy, if had not realized that the line
through the outer nest just created plus the quintuple point had
at least 9 intersection with our hypothetical octic. Hence the
next idea was to look at Fig.\,j where the outer nest is created
from the petals (and not from the 8 extra oval). The little defect
of this curve is that it has 2 circuits already and therefore the
fundamental part of the patch (i.e. without the micro quantum
ovals) already consume two circuits, so that only 5 are left and
we distributed as on Fig.\,j1 the resulting gluing is only an
$(M-1)$-octic with only 17 inner ovals instead of the 18 desired
ones. To palliate this defect we change a bit the connections to
get Fig.\,k where we have just one circuit and  smoothing
appropriately (Fig.\,k1) yields the bosonic scheme
$1\frac{1}{1}\frac{18}{1}$ via Fig.\,k2. Finally, Fig.\,k3 shows
as series of quantum fluctuation of the 8 indeterminate ovals
creating all schemes of the bosonic strip (of course making so a
big razzia (a clean sweep)
%%%% CHECKED IN DICO
on Hilbert's 16th, of course sometimes entering in conflict with
Orevkov). But of course there is no direct confrontation of our
method with Orevkov's conclusion since the big job is to construct
the singular octic in the algebraic category, yet we know at least
a possible place for where to look for new constructions.
Summarizing:

\begin{figure}[h]\Figskip
%\vskip-1.2cm\penalty0
%\centering
\hskip-2.7cm\penalty0
\epsfig{figure=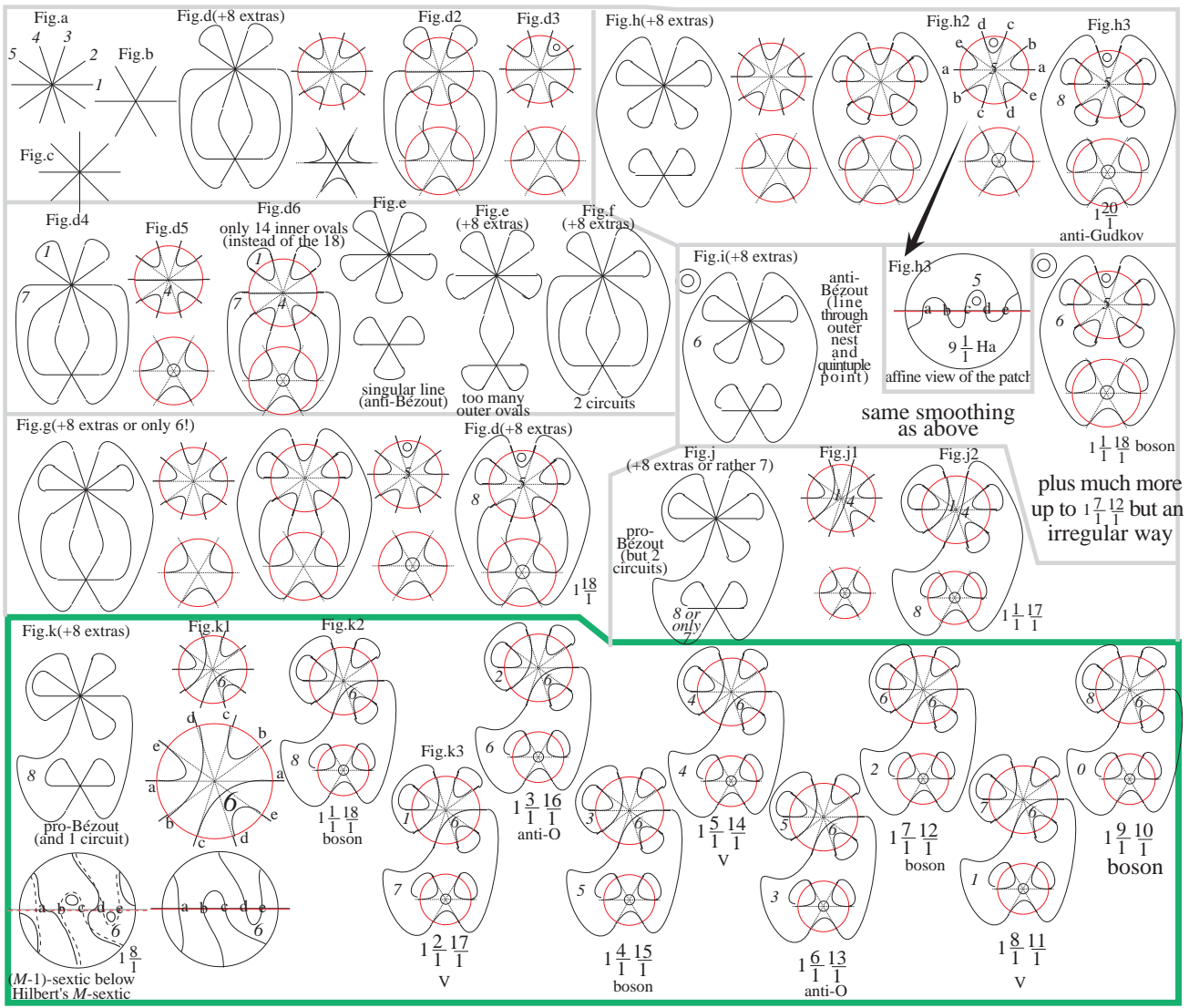,width=172mm} \captionskipAG
  \caption{\label{GabardDEGREE8_5+3:fig}%
Singular octics with a quintuple and triple point while creating
all bosons in a B\'ezout regular fashion (and of course
potentially more like schemes prohibited by academician Stepan
Orevkov)} \figskip
\end{figure}

%%To summarize:

\begin{Scholium}
All of the bosonic $M$-octics are potentially realized by curves
with a quintuple plus a triple point with branches interconnected
%%%%between themselves
as depicted on Fig.\,\ref{GabardDEGREE8_5+3:fig}k. So potentially
there is no more obstruction on Hilbert's 16th than those (15
many) presently known.
\end{Scholium}

{\it Added} [05.10.13].---One can try  exasperating B\'ezout by
passing a conic through the octic singularities along the most
osculating way. For instance we may impose two `horizontal'
tangents at the singularities, plus visiting one simple oval. Then
the multiplicity intersection is 6+4+2=12, still in B\'ezout's
range (2.8=16). Maybe fixing the upper tangent slightly inclined
(angle ca. -10 degrees), then it penetrates inside the double
loop, and to escape it forces 2 intersections. Moreover one more
intersection seems granted by crossing the path joining both
singularities. All this gives 5+2+1+2+4=14 still below B\'ezout.
So it seems that our freehand curve resists the B\'ezout test.

Of course it remains to find algebraic constructions. It remains
also to investigate the composition of two quadruple points.
Philosophically we see that our game (Viro method basically,
sometimes in variants \`a la Shustin, Korchagin, etc.) is
basically a matter of exploring the complement of the discriminant
(smooth curve) by entering through such chambers via the walls of
the discriminant separating the big castle of all curves into
varied chambers. Of course  there is no theological reason a
priori that using even higher strata conditioned by the presence
of several singularities (recall that the smooth locus of the
discriminant is swept out by uninodal curves) we should be able to
access all the chambers. So perhaps Viro's  method can fail to
detect all schemes, but was is certain is that if it works then it
works. This sounds a tautologically but we hope that you know what
we mean.

Let us now work out singularity $4+4$. Then one could also imagine
3 singularity or more. So there is a very exciting global game of
constructing curves with prescribed singularities. The more the
singularity is low the less it eats to the genus and so we can add
several of them. Further it seems that low singularities have an
easy dissipation theory, yet appeal to global curves of a more
complicated nature (high genus). So it seems that one must proceed
to a clever dissection of the problem by dividing in nearly equal
part the difficulties allied to the local and global aspects of
the questions, i.e. to have a relatively easy dissipation theory
while having a fairly easy global singular curve to construct.

[23.07.13] Let us now come back to singular weight $4+4$. After a
quick browse through qualitative pictures, we had the idea to
consider configuration of ellipses like on
Fig.\,\ref{GabardDEGREE8_6:fig}. Alas we had some pain to reach
$M$-curves. Fig.\,g2 is interesting as it seems to contradict the
principle of independence of smoothing since the constructed curve
appears to violate the maximality of the doubled quadrifolium, yet
in reality the bug is of course that the dissipation used is
tantamount to a G\"urtelkurve (quartic with 2 nested ovals) and
the latter do not support any extra infinitesimal ovals being
already a saturated configuration. So Fig.\,g2 is only correct
when one erases all ``3'' parameters and then one recovers the
usual doubled quadrifolium. As to Fig.\,g1 we noticed after a
better inspection that used is an illegal one for it would involve
a quartic with 3 nested ovals (symbol $\frac{3}{1}$) which
obviously violates the principle of saturation of the
G\"urtelkurve $\frac{1}{1}$ (Fig.\,g1a). For sure we can drag the
3 ovals at the center of the patch and then we get the scheme
$16\frac{1}{1}\frac{1}{1}$. Fig.\,g3 shows how to get the
``unnest'' 16, but this dissipation is far from maximal due to
much consanguinity between the behavior of branches at infinity in
the patch, allowing us only to introduce 2 extra ovals. At any
rate it is fairly evident that the quadri-ellipse of Fig.\,g is
not deformable  to any $M$-curve. Next going back to freehand
tracing we manufacture Fig.\,i, which is fairly interesting as it
dissipations includes 3 schemes one of which agressing Orevkov's
link-theoretic prohibition of $1\frac{3}{1}\frac{16}{1}$ and the
other one attacking Viro's sporadic obstruction of
$\frac{3}{1}\frac{3}{1}\frac{13}{1}$. Alas this is not a serious
corruption of the results of those Russian scholars because our
Fig.\,i is not extremely regular with respect to B\'ezout when it
comes to trace the singular line. This invites to consider Fig.\,j
where this B\'ezout trauma is palliated, but alas the price to pay
is the formation of 2 circuits hence reducing the number of boni
ovals. Fig.\,k has the sam defect. Fig.\,l is interesting for
coming quite close to an $M$-scheme prohibited by Shustin (namely
$(1,18\frac{2}{1})$) yet only through an $(M-1)$-approximation
(namely $(1,18\frac{1}{1})$). Next Fig.\,m, though deviating from
the desideratum of having one outer oval only to land in the
bosonic strip, looks excellent otherwise as it respects B\'ezout
and possess only one circuit so that the maximum number of extra
oval permissible by Harnack is gained (namely 9). Let us place
them inside the biggest circuit and get Fig.\,m1 realizing the
Viro scheme $10\frac{11}{1}$. However on exploiting the other
admissible smoothing of Fig.\,m2 we enter in conflict with
B\'ezout since we get 2 subnest in the largest (banana-shaped)
oval.

So we have proven the following paradox:

\begin{lemma}
There is a structural incompatibility between existence of the
curve Fig.\,m (for any quantum placement of the nine extra ovals)
and the Pl\"ucker-Klein-Harnack-Brusotti-Gudkov-Viro-Shustin
principle of independency of smoothing singularities.
\end{lemma}

Yes, but our mistake is simple to detect, namely we used a
dissipation which cannot exit, because the corresponding patch
viewed as an affine quartic yields Fig.\,m4 whose smoothing
(Fig.\,m5) violates B\'ezout (a quintic cannot be nested unless it
is the deep nest $\frac{1}{1}J$). So everything is restored to
normality and we can still expect for an algebraic model of
Fig.\,m. Actually on doing all quantum transfer of ovals
compatible with Gudkov periodicity (and optionally Viro's imparity
law) we get a collection of schemes all permissible (i.e. actually
constructed) in particular two schemes of Korchagin namely
$10(1,6\frac{4}{1})$ and $10(1,2\frac{8}{1})$. Accordingly it may
seems reasonable to guess existence of a singular octic like
Fig.\,m. Alas it seems that postulating its existence leads to no
new result except for reproving some Korchagin's scheme in a
perhaps more elementary fashion (as we truly appeal to a trivial
dissipation theory involving affine quartics). It remains of
course then to construct Fig.\,m via the usual Gudkov-Viro trick
of hyperbolism based on Huyghens-Newton-Cremona (but we are not
not very strong in this game). Incidentally it may be observed
that Fig.\,d is fairly close to Shustin's medusa (depicted on
Fig.\,\ref{ViroDEGREE8_SHUSTIN_NEW:fig}).
%%%%%% CHECK CROSS-LINK).

\begin{figure}[h]\Figskip
%\vskip-1.2cm\penalty0
%\centering
\hskip-2.7cm\penalty0
\epsfig{figure=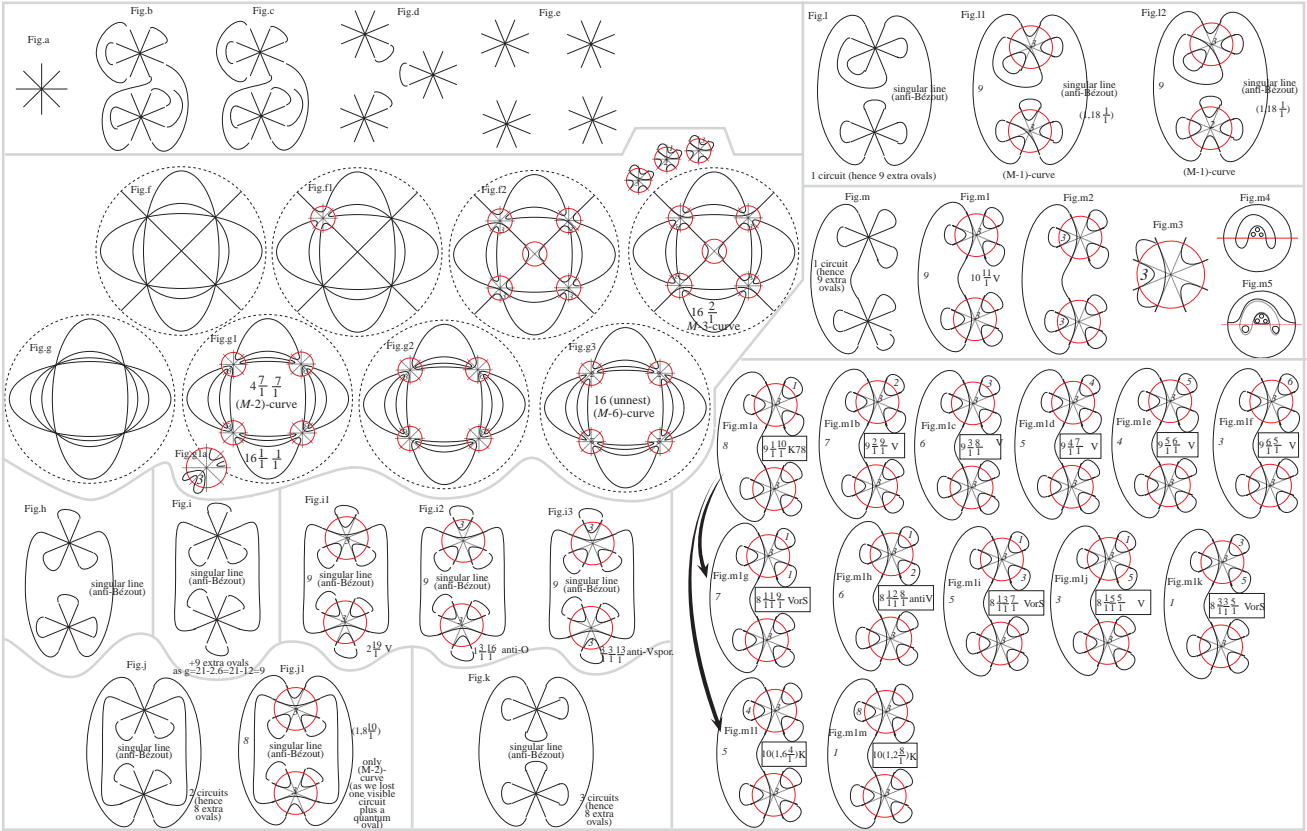,width=172mm} \captionskipAG
  \caption{\label{GabardDEGREE8_6:fig}%
Singular octics from two (or more) quadruple points while creating
2 curves by Korchagin} \figskip
\end{figure}

Of course it is not clear that we exhausted all ways of connecting
the branches of 2 quadruple points, yet it looks hard to reach the
bosonic strip (doubly nested schemes with one outer oval). Of
course we could try to look at what happens with 3 quadruple
points.

Via Fig.\,c we see that we can realize most of Shustin's schemes
(safe those with zero outer ovals). To get to them it suffices
however to start with Fig.\,d which looks even more like a
mushroom and so even more anti-B\'ezout. The latter curve produces
smoothing conflicting with Viro's sporadic obstructions while
sometimes regenerating schemes due to Shustin. From Fig.\,f on, we
consider curve with threefold symmetry and so three quadruple
points. Fig.\,g suitable smoothed (along the unique
$M$-possibility) offers Viro's (non-elementary) scheme
$6\frac{15}{1}$ and quantum fluctuating the 3 outer ovals we may
reach as well schemes of Hilbert, Korchagin and Viro. Fig.\,h is
merely isotopic to the former (yet even dihedrally symmetric).
Fig.\,i offers Harnack's scheme $18\frac{3}{1}$ from which one can
derive through quantum fluctuation to schemes by Hilbert,
Chevallier and Korchagin (yet nothing truly new except that we did
not as yet digested the original construction of Korchagin and
Chevallier as the tend to use the Newton polygon formalism with
which we feel uncomfortable to say the least). Finally Fig.\,j
produce the hard scheme of Shustin
$4\frac{5}{1}\frac{5}{1}\frac{5}{1}$, and in this case B\'ezout
forbids any quantum fluctuation of the 3 outer ovals. So no more
scheme can be derived. Again our construction is perhaps quite
close to Shustin's original which we could not follow due to a
lack of picturing. Phenomenologically, one could imagine that
Fig.\,j appears through perturbation of a doubled tricuspidal
quartic. To be frank and honest, it seems that our picture
(Fig.\,j) violates B\'ezout when tracing the line joining 2
singularities, yet perhaps there is some distortion avoiding this
aberration.

\begin{figure}[h]\Figskip
%\vskip-1.2cm\penalty0
%\centering
\hskip-2.7cm\penalty0
\epsfig{figure=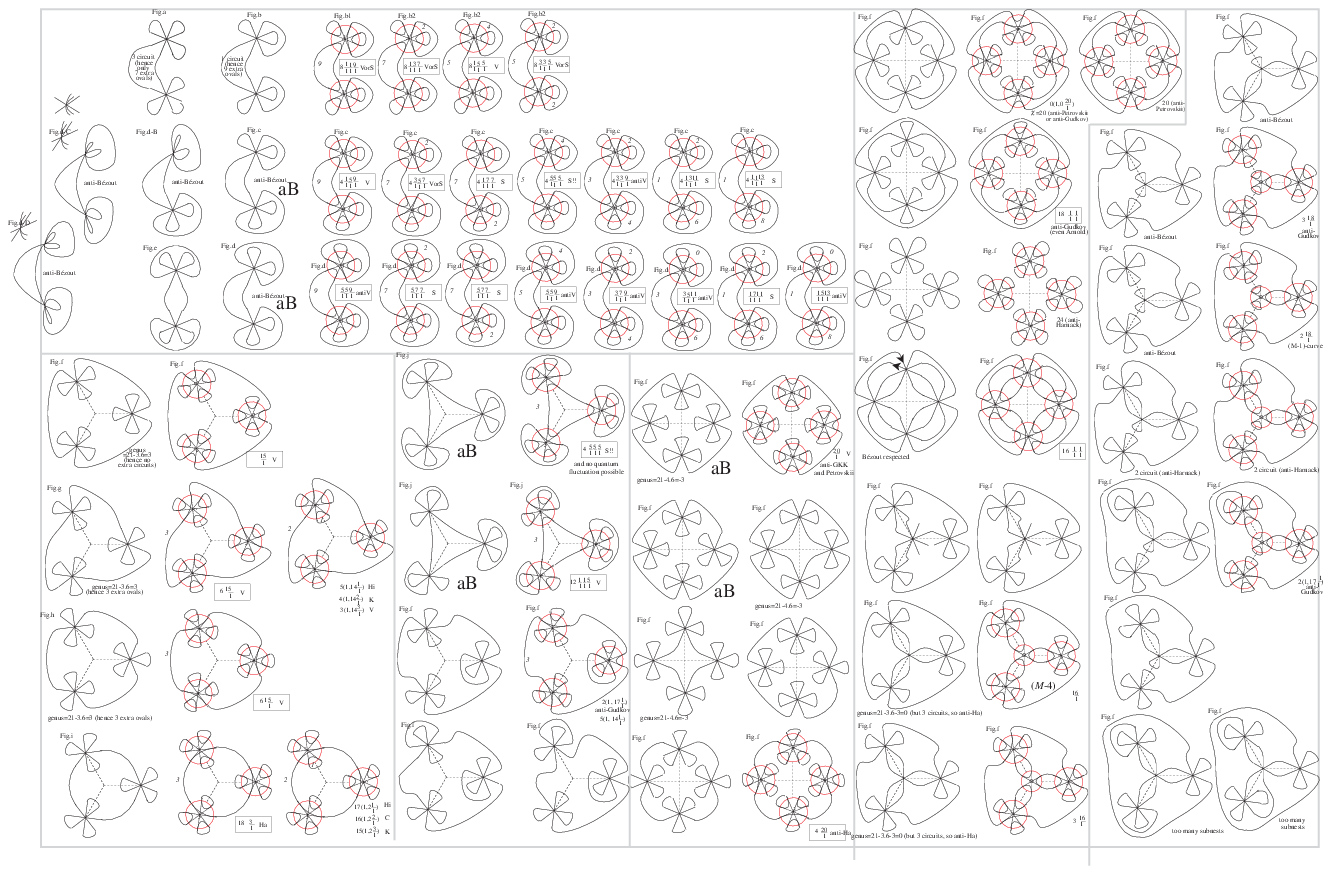,width=172mm} \captionskipAG
  \caption{\label{GabardDEGREE8_7:fig}%
Singular octics from 2 or 3 quadruple points rotationally
invariant} \figskip
\end{figure}

[24.07.13] Perhaps one should as well examine less symmetric
junctures between the three quadruple points. Or one can also
consider a configuration with 4 quadruple points.

[25.07.13] Coming back to Fig.\,d or even Fig.\,c both having rich
descendence (sometimes conflicting with Viro's sporadic
obstructions), yet being themselves optically anti-B\'ezout when
it comes to trace the line through both singularities, we can try
to remedy this defect via Fig.\,d-B. This is an attempt to stretch
the curve as much as we can in order to respect B\'ezout and
converging ultimately to Fig.\,d-C. During this process the
quadruple point loose its ordinariness (i.e. 4 distinct tangents)
so that 2 of the branches acquires  coincident tangents. This new
singularity can be interpreted as an ordinary triple point with a
two-fold (parabolic) branch tangent to one of the 3 linear
branches. Alas it seems that the limiting Fig.\,d-C still corrupts
B\'ezout because the line through the singularities has now a
tangency along two branches and therefore the multiplicity
intersection is about 6 at one singular point, and therefore at
least 12.

Next we meditated more about how to place four quadruple points,
%%%like on Fig.\,??
with 3 of them forming an equilateral triangle
and the fourth one lying in its center. Of course this
configuration seems to lack structural symmetry and so suggested
to us to look at the same with 4 quadruple points distributed on a
square. Then there is 5 quadruple points and thus a conic through
them cut $5\cdot 4=20> 2\cdot 8=16$ so that B\'ezout is surpassed.
In fact it seems that four quadruple points are already prohibited
since considering the pencil of conics through the 4 singularities
one can still pass a conic through any other point of the $C_8$ so
that $4\cdot 4 +1 =17> 2\dot 8$ points are created violating once
more B\'ezout. Of course there is however the little exception of
when the octic split as 4 curves in a pencil of conics in which
case the given distribution of singularity (four quadruple points)
is realized.

Next we had the idea of a central quintuple point plus five triple
points gravitating around. Alas the resulting curves violates
Harnack's bound. And actually the initial singular  curve foils
B\'ezout when tracing the conic through the quintuple  and 4
triple points which has intersection multiplicity at least $5
+4\cdot 3=17> 2\cdot 8$. So let us kill some triple points to get
Fig.\,\ref{GabardDEGREE8_8:fig}b. On smoothing the latter we get
the $M$-scheme $20\frac{1}{1}$, which is anti-Gudkov (or even
anti-Petrovskii). Hence we get as an interesting corollary of
those Russian scholars (whose work is logically founded on either
Rohlin 1952 or Euler-Jacobi-Kronecker interpolation) the following
result not directly imputable to B\'ezout (as far as we can
judge):

\begin{lemma}
There is no singular octic whose real picture is like that of
Fig.\,\ref{GabardDEGREE8_8:fig}b albeit the latter seems perfectly
B\'ezout compatible. Crudely put, there is rigidity of algebraic
curves beyond the (naive) optical level.
\end{lemma}

\begin{figure}[h]\Figskip
%\vskip-1.2cm\penalty0
%\centering
\hskip-2.7cm\penalty0
\epsfig{figure=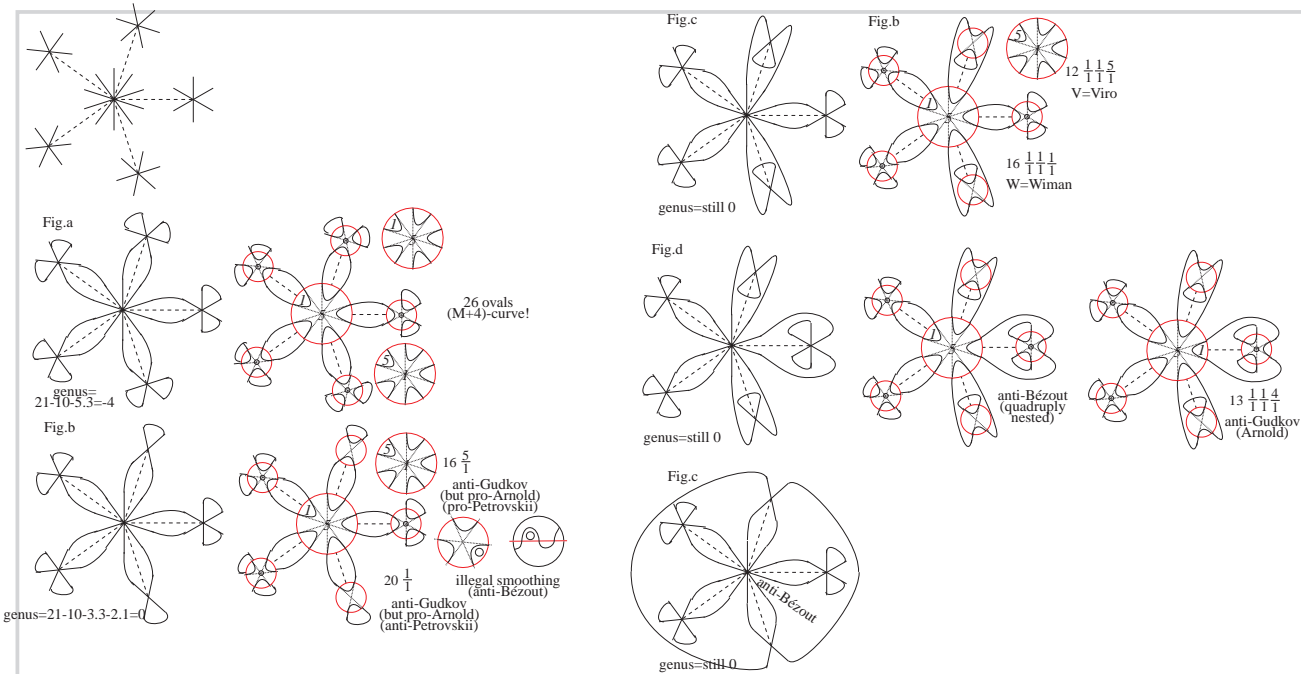,width=172mm} \captionskipAG
  \caption{\label{GabardDEGREE8_8:fig}%
Singular octics from two or three quadruple points invariant under
rotation by $2\pi/3$} \figskip
\end{figure}

Fig.\,b was not  adjusted to Gudkov periodicity and it is a simple
matter to arrange this issue via Fig.\,c where the 2 loops are
introverted. Of course the resulting schemes are not of the most
exciting types, yet one can seek for more introversion yet usually
this runs against B\'ezout's law.

{\it Added} [26.07.13] Then we had the idea to use only a
quintuple point plus a constellation of double points so as to
leave more maneuvering room for closing the petals without quickly
entering in conflict with B\'ezout (see
Fig.\,\ref{GabardDEGREE8_8B:fig}a). Its smoothing (Fig.\,a1) is
alas only an $(M-4)$-curve with 18 ovals. Fig.\,b is a loose essay
to introduce more nodes, the philosophy being that a rational
curve is perhaps easiest to construct {\it ab ovo\/} and having
only one circuit it is not subjected to the extra quantum ovals.
Its smoothing Fig.\,b1 has also only 18 ovals. Fig.\,c is another
way to close the distribution of singularities specified. However
if the 6 extra ovals are lying outside then the smoothing Fig.\,c1
contradicts Gudkov periodicity (proved by Rohlin). In contrast for
other distribution of the outer quantum ovals fluctuating inside
like Fig.\,c2, c3=c4, or c6=c7 we get respectable schemes
originally due to Ha=Harnack, G=Gudkov or V=Viro. Yet the
distribution Fig.\,c5 is forbidden if one believes in Viro's
imparity law. Next one can look at Fig.\,d which has one extra
oval (quantum as its location is not yet decided).  As on Fig.\,d1
let $x$, $y$, $z$ denote the marked position taken by the quantum
oval. If at $x$ then we get the scheme $20\frac{1}{1}$, which is
anti-Gudkov. Otherwise we get 3 schemes with $\chi=18$, hence
violating the primitive version of Gudkov periodicity proved by
Arnold (and a formal consequence of Rohlin's formula).

Hence albeit Fig.\,d is not extremely interesting from the
viewpoints of constructing the 6 missing bosons, it is of some
interest for showing that Viro's method (i.e. the possibility and
independence of smoothing \`a la Pl\"ucker, Klein, Harnack,
Hilbert, Brusotti, Viro, Shustin, etc.) acts actually when
combined with Gudkov as a way of prohibiting  schemes of singular
curves. For instance, we have  proven the following fact.

\begin{lemma}
There is no singular octic of genus 1 whose singular circuit is
isotopic to  Fig.\,\ref{GabardDEGREE8_8B:fig}d, whatever the
location of the one extra oval is.
\end{lemma}

Fig.\,e shows another curve closing the given distribution of
singularities, but its smoothing violates B\'ezout for conics as
the configuration expands the doubled quadrifolium. Fig.\,g is
permissible, yet its smoothings are not of the most exciting type
having in the maximal case at least 3 outer ovals. Finally in view
of the bound of the number of nest (at most three of them else the
curve saturated at the doubled quadrifolium) we were rather
inclined to tolerate at most threefold symmetry suggesting
Fig.\,h, i, j were we successively increased the number of nodes.
Alas Fig.\j smoothed as Fig.\,j1 violated Gudkov' periodicity.
Hence it follows again that Fig.\j cannot exists algebraically
which is not evident optically via say B\'ezout alone. In contrast
a singular scheme like Fig.\,k yields the $M$-scheme
$4\frac{5}{1}\frac{5}{1}\frac{5}{1}$ due to Shustin (albeit we did
not understood his construction), and is therefore not prohibited.

\begin{figure}[h]\Figskip
%\vskip-1.2cm\penalty0
%\centering
\hskip-2.7cm\penalty0
\epsfig{figure=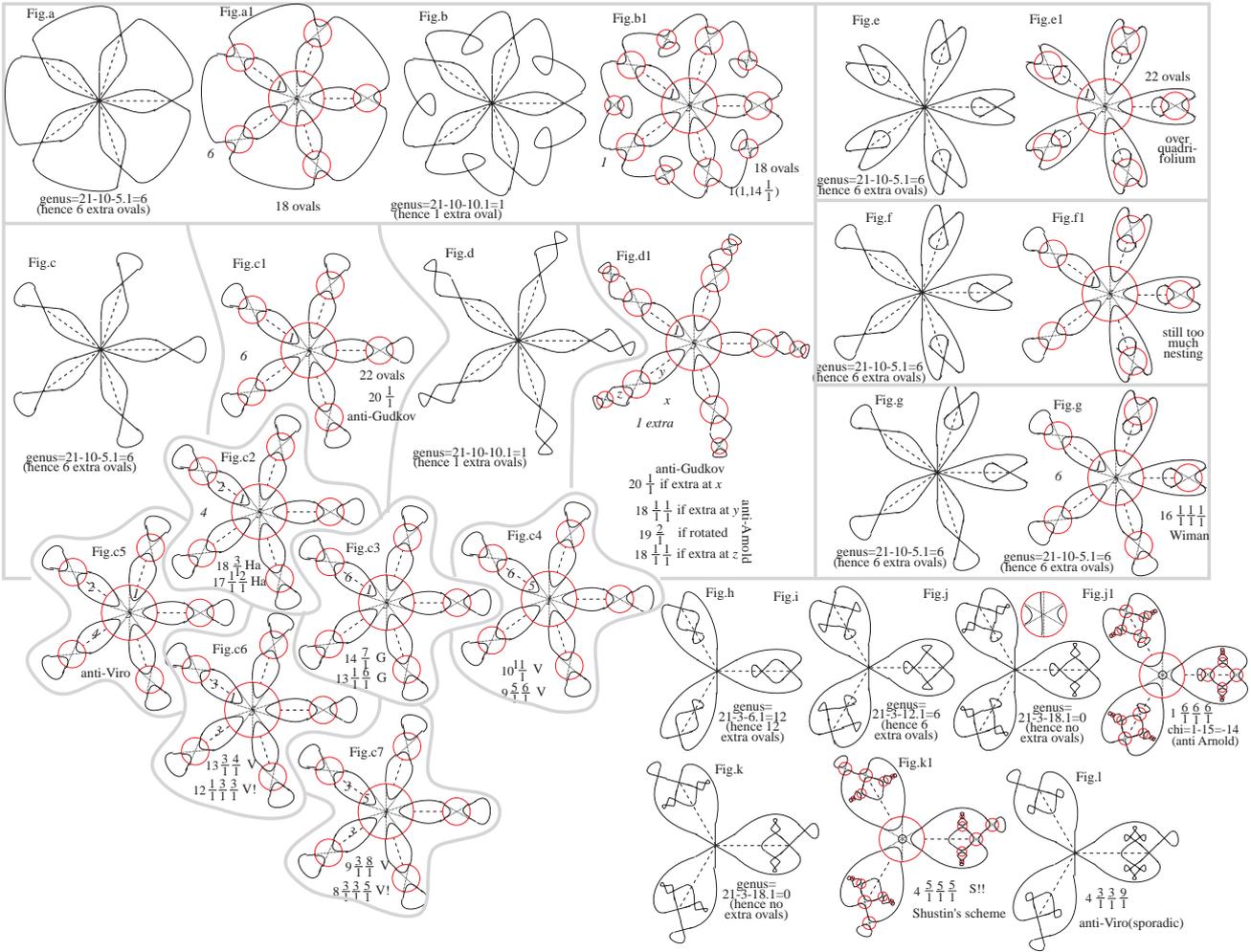,width=172mm} \captionskipAG
  \caption{\label{GabardDEGREE8_8B:fig}%
Singular octics with ...} \figskip
\end{figure}

In some contrast Fig.\,\ref{GabardDEGREE8_8B:fig}l deduced by
transferring two nodes from 2 petals in the first petal, yields
via the unique maximal smoothing the scheme
$4\frac{3}{1}\frac{3}{1}\frac{9}{1}$ in principle prohibited by
Oleg Viro (sporadic obstruction not readily available from the pen
of its discoverer). So granting this as correct we have:

\begin{lemma}
Modulo the truth of Viro sporadic
$4\frac{3}{1}\frac{3}{1}\frac{9}{1}$, there is no rational
singular octic whose real locus is isotopic to
Fig.\,\ref{GabardDEGREE8_8B:fig}l, albeit the latter does not
frankly seem to offend (Monsieur \'Etienne) B\'ezout.
\end{lemma}

[27.07.13] Experimentation shows that our threefold symmetry tends
to lead to the trinested case, while the bosonic strip is
primarily a matter of doubly nested schemes where 4 bosons are
concentrated. This suggests to switch to a simple twofold symmetry
like on Fig.\,\ref{GabardDEGREE8_8C:fig}a. Here we arrange to get
21 nodes  as to drop the genus down to zero (rational curve) for
which one can in principle write down an explicit parametrization.
The Fig.\,b does the job yet create when smoothed the scheme
$\frac{9}{1}\frac{11}{1}$ which is anti-Gudkov. However the
closest Gudkovian approximation is the scheme
$1\frac{9}{1}\frac{10}{1}$ (highly bosonic, i.e nobody ever
succeeded to realize it nor to prohibit it) which could be created
out of the unicursal curve of Fig.\,c where we just traded an
inner node of Fig.\,b for an outgoing node. Of course one can then
successively transplant inner left nodes to the right generating
so Fig.\,c1, c2, c3, \dots, c9, c10 sweeping thereby all the 2
doubly nested bosons via unicursal curves. Hence:

\begin{Scholium}
Maybe there is a simple way to create some of the missing bosons
via a rational curve with 21 nodes. Algorithmically, it is perhaps
a reasonable game to trace by brute computer-force some random
parametrization of degree 8 (yet both components of the
parametrization may be of lower order) so that one of the above
pattern appears
%%%%as the computer tracing.
on the screen (Hilbert's retina). Of course as we have now an
explicit parametrization instead of a implicit (random) equation
of degree 8, we hope that the present problem is somewhat more
tractable electronically that the original setting of Hilbert's
16th involving highly irrational curves.
\end{Scholium}

Of course in the above setting as our curve have only ordinary
nodes they may be thought of as generic immersions of the circle
(of an algebraic character) and thus more or less suited to a
random search.

Fig.\,f shows some spires-like chain of nodes (that could also be
imagined as a cactus like Fig.\,f2 or as the slalom variant of
Fig.\,f3). Alas the resulting smoothed scheme $1\frac{20}{1}$ is
anti-Gudkov, but one can repair this defect by going to Fig.\,g
realizing $2\frac{19}{1}$ instead.

\begin{figure}[h]\Figskip
%\vskip-1.2cm\penalty0
%\centering
\hskip-2.7cm\penalty0
\epsfig{figure=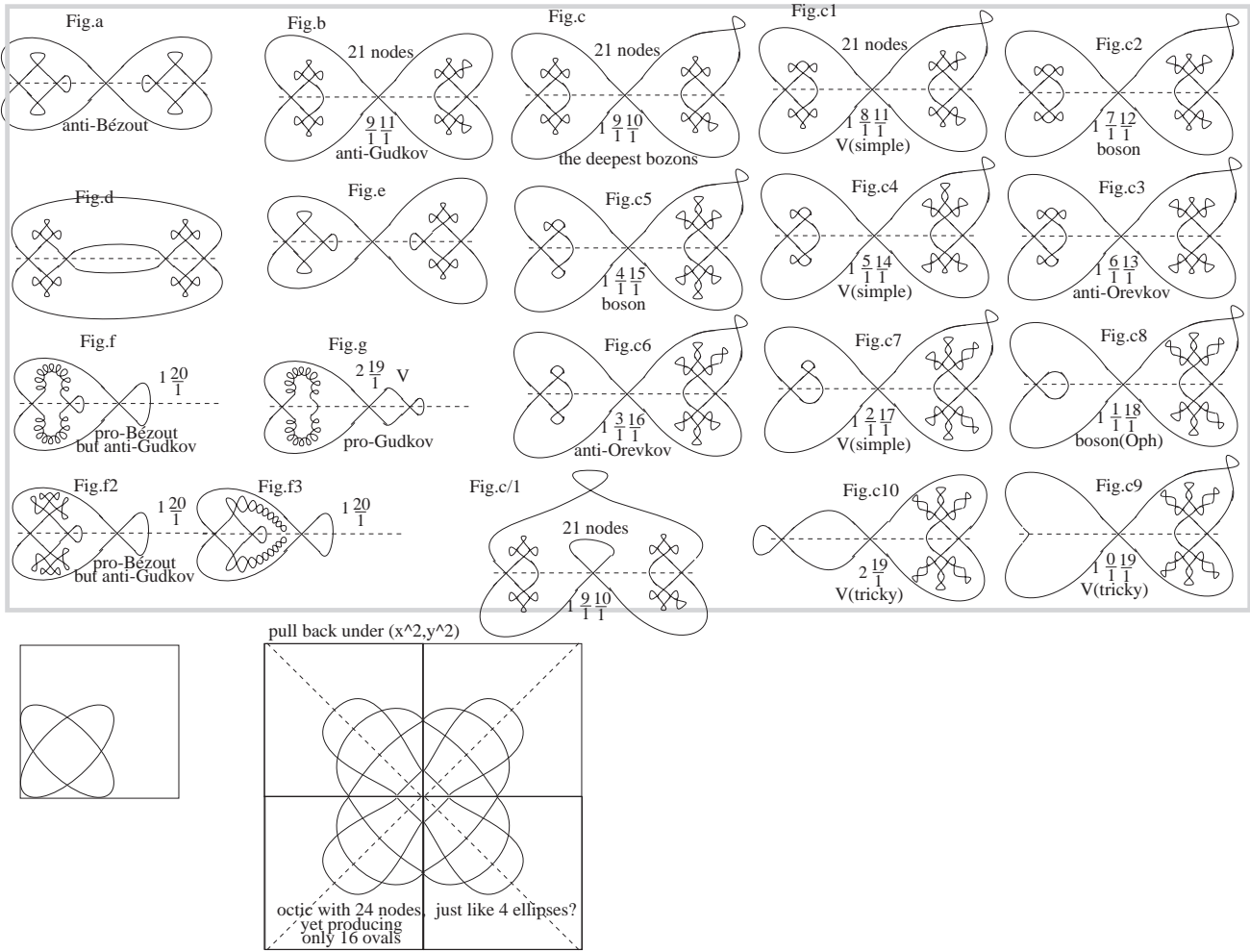,width=172mm} \captionskipAG
  \caption{\label{GabardDEGREE8_8C:fig}%
Singular octics with ...} \figskip
\end{figure}

Of course all this is merely a matter of free-hand tracing of very
hypothetical unicursal curve of degree 8,yet this could be of some
relevance to settling Hilbert's 16th in degree $m=8$, which is the
last one where the combinatorics is still at human size, but where
we have to confess to be faced with serious geometric mysteries
(bosons, sporadic Viro obstructions, Orevkov braid theoretic
obstructions). Of course as a matter of construction it seems that
the above pseudo-construction are not ideally suited as the use
the simplest dissipation theory (for the ordinary node due to
Severi-Brusotti and somehow anticipated by Germans like Klein,
Harnack, etc.) yet with using a fairly complicated global curve
whose range is purely hypothetical. In principle we can increase a
bit the complexity of singularities involved while decreasing a
bit the trickiness of the initial singular curve. We saw already a
lot of such example where using merely ordinary multiple point we
could generate schemes due to Viro yet without appealing to the
dissipation of complex singularities like $X_{21}$ (i.e. quadruple
point with 2nd order tangency between the branches). For a
concrete example yet not perfectly justified see
Fig.\,\ref{GabardDEGREE8_7:fig}g or h.

Further let us observe again that among all of Viro's sporadic
obstructions those impeding $4\frac{3}{1}\frac{3}{1}\frac{9}{1}$
looks to be the more risky one as the scheme is formally covered
by an affine septic.
 So referring to a variant of the dissipation of the septuple
 point one can also easily corrupts this Viro octic obstruction,
which if correct therefore implies an obstruction on one isotopy
type of affine septics   affine

\begin{figure}[h]\Figskip
%\vskip-1.2cm\penalty0
%\centering
\hskip-2.7cm\penalty0
\epsfig{figure=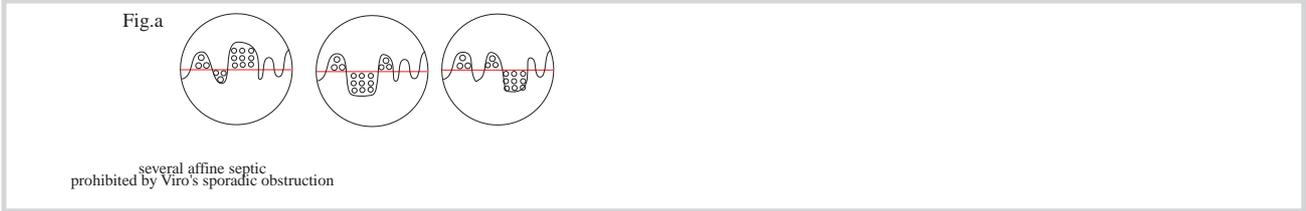,width=172mm}
\captionskipAG
  \caption{\label{GabardDEGREE8_septupleVIRO:fig}%
Several affine septics (indirectly) prohibited by
%%%Stepan
Viro's most stinky sporadic prohibition of
$4\frac{3}{1}\frac{3}{1}\frac{9}{1}$} \figskip
\end{figure}

[28.07.13] Maybe one can give more weight (geometric evidence) to
the curve constructed by assigning a more massive singular points
instead of the many nodes of Fig.\,\ref{GabardDEGREE8_8C:fig}. By
increasing the mass of the singularity we quickly arrive at
Fig.\,\ref{GabardDEGREE8_8D:fig}b where 2 quintuple points are
introduced on both side of an ordinary node which should act as a
splitter leading us in the doubly nested realm (where most bosons
are concentrated). Why two quintuple points? Just because then we
have a lucky stroke configuration as a 5-fold points drops 10 to
the genus so that our curve has genus 0. Yet it should be noted
that only concentrating on the genus we lost from sight the more
basic constraint on the degree $8$ which cannot tolerate an
intersection of 10. So we consider the same configuration yet with
only two {\it quadruple\/} points. Alas the resulting schemes
Fig.\,c,c1,c2 turned to be either anti-B\'ezout, or not
Harnack-maximal. In fact it seems that there is a better
configuration when splitting 8 as $5+3$ yielding Fig.\,d whose
smoothings manage to sweep out virtually all doubly-nested bosons
safe the first species $1\frac{1}{1}\frac{18}{1}$. Of course the
curve Fig.\,d as traced is not B\'ezout permissible yet maybe a
suitable contortion of it can destroy the co-linearity of the 3
singularities. It is evident that each of our depicted flower has
a trunk (stem/stalk/tige) which looks invaginated. It seems clear
that invagination have to stay in front of each other as to
respect B\'ezout for the  line through both singularities. Fig.\,e
supplies a projective depiction. Alas it seems that when the line
penetrate in the hearth shaped region delimited by the two arcs
connecting the triple to the double point that we get again in
trouble with B\'ezout as the line as to escape of the Jordan cell.
So it seems that there is a serious topological obstruction to
existence of a singular octic with the prescribed distribution of
singularities and controlled isotopic type. Of course our failing
attempt to construct decently those bosons does not mean that the
bosons does not exist, yet one could imagine \`a powered
Hilbert-Rohn method stating that any smooth curve (say doubly
nested and with one outer oval, i.e. in the bosonic range)
degenerate toward such a singular curve (with singularities of
masses 2,3,5) and then our argument would give a general
obstruction. Fig.\,f seems to be a solution yet we lost the
splitting double point. Fig.\,g is an attempt to reintroduce the
lost node. Actually Fig.\,f with suitable quantum ovals ($x$ and
$y$) can produce nearly all bosons except the highest one
$1\frac{9}{1}\frac{10}{1}$. Of course this picture can be
affinized (i.e. be put at finite distance) as shown on Fig.\,fA,
which seems a staphylococcus. It is tempting to posit a symmetric
realization of the curve. (Incidentally as in the bosonic strip we
have $\chi=-16$, this Ansatz is compatible with Fiedler's
strengthened version modulo 16 of Gudkov's hypothesis $\chi\equiv
k^2 \pmod 8$.) Next we may examine the varied distribution of
ovals. Fig.\,fA0 shows the case $x=0$, hence with $y=7$ in the
outer oval depicted at infinity. Fig.\,fA1 shows the case $x=1$
where there is one inner oval traced inside the singular circuit.
By symmetry the latter is forced to be self-symmetric (invariant)
but then we get troubles with B\'ezout as the invariant line
through both singularities exhibit already 8 intersections. By the
way the smoothed scheme is prohibited by Orevkov and so we get
some feeling of understanding his highbrow braid-theoretic
obstruction. Fig.\,fA2 with $x=2$ many inner ovals respects
B\'ezout when both are distributed vertically as on Fig.\,fA2-bis,
while the resulting scheme really exists through Viro's simplest
construction. For $x=3$ we get trouble with B\'ezout and the
resulting scheme is a boson. So we feel fairly close to have a
general law emerging. In particular we would say that the first
boson $1\frac{7}{1}\frac{12}{1}$ exist, whilst
$1\frac{4}{1}\frac{15}{1}$ would not. So far so good, but the
sequel make our observational law quite stormy. Indeed as $x=4$,
our naive law  would vote toward existence of the scheme, yet the
scheme is prohibited by Orevkov. (Of course it could be that
Orevkov's result are wrong but this is only weak superstition.)
Next as $x=5$, the odd number of inner ovals forces one to be
invariant and so to possess a supernumerary point of intersection
with the singular line is created, however the scheme is
accessible through Viro's simplest method. So our law is severely
foiled for a second time. When $x=6$, we get the boson
$1\frac{1}{1}\frac{18}{1}$ fairly likely to exist (remember also
its realizability via decomposable curves of degree $2+6$).
Finally as $x=7$, B\'ezout is foiled but the scheme still exist
via a (tricky) Viro method. This is a 3rd corruption of our law.

\begin{figure}[h]\Figskip
%\vskip-1.2cm\penalty0
%\centering
\hskip-2.7cm\penalty0
\epsfig{figure=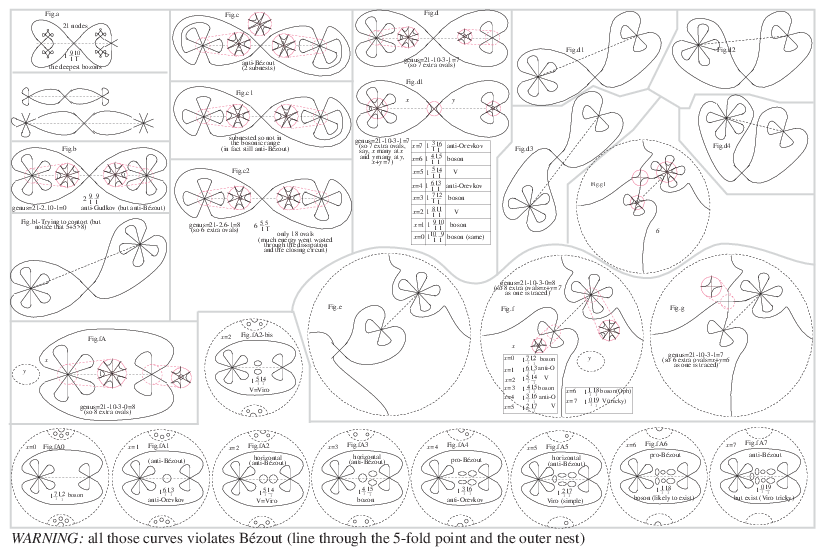,width=172mm} \captionskipAG
  \caption{\label{GabardDEGREE8_8D:fig}%
Singular octics with ...} \figskip
\end{figure}

Of course our law splits in two parts (construction and
prohibition):

(LAW1) If the number of inner ovals is odd then the scheme does
not exist. Of course this part of the law is fairly weak has it
imposes a very peculiar mode of generation of the scheme. By the
way it is contradicted twice by Viro's construtions.

(LAW2) When the number of inner ovals is even then the scheme does
exist. This sense looks logically stronger as it suffices to have
one construction to have a construction. Ye this principle is
foiled once by Orevkov (as $x=4$). So either Orevkov's result is
lase or more likely the octic of Fig.\,fA4
 does not exist even though it looks B\'ezout respectable.

In conclusion very little can be extracted from our naive method,
yet it cannot be excluded that some few new schemes could be
constructed along this strategy.

Now it remains to examine Fig.\,g. One of the smoothing leads back
to Fig.\,f (already analyzed). Actually the singular curve has two
circuits so we do not have the maximum number of quantum ovals
gained and so we certainly fails Harnack maximality. o the
situation does not seem worth investigating any further (and
Fig.\,g1 (Gibraltar's canal)
%%%or detroy???
 shows that as expected we get only 20
ovals).

Let us summarize as follows our weak knowledge:

\begin{lemma}
Among all schemes in the bosonic strip (doubly nested and one
outer oval) all are realized via a singularity of masses $3+5$
safe $1\frac{9}{1}\frac{10}{1}$. This gives weak evidence that the
schemes $1\frac{7}{1}\frac{12}{1}$, $1\frac{3}{1}\frac{16}{1}$,
$1\frac{1}{1}\frac{18}{1}$ do exist algebraically but the middle
term is actually ruled out according to Orevkov's theory. Further
the last boson $1\frac{9}{1}\frac{10}{1}$ albeit not accessible
via the $[3+5]$-method, it is via the $[2+3+5]$-method (see
Fig.\,\ref{GabardDEGREE8_8D:fig}d1), where however the picture is
plagued by robust anti-B\'ezoutism. So our rating agency can give
only low existential evidence for this species to be observed in
nature.
\end{lemma}

[ca. 26.07.13] Maybe (somewhat inspired by a certain Yves
(tailleur de pierre) who learned me yesterday that natural
crystals never stabilizes to fivefold symmetry) we can do the same
game for a sextuple point with a halo of double points (see
Fig.\,\ref{GabardDEGREE8_9:fig}). Suddenly one sees that
arithmetics works then smoother as $21-15-6\cdot 1=0$, so that
more symmetry is gained (while keeping the genus $\ge 0$). We do
not know if there is a direct relation between crystals and
algebraic curves (apart of course in the prose of Alexander
Grothendieck). On tracing Fig.\,c we realized that our
configuration of singularities overwhelms B\'ezout, since the
sextuple and two double points are aligned. After a very erratic
search the sole interesting thing is Fig.\,g5A4 where we recover a
Chevallier scheme using the hypothetical curve of Fig.\,g5.

\begin{figure}[h]\Figskip
%\vskip-1.2cm\penalty0
%\centering
\hskip-2.7cm\penalty0
\epsfig{figure=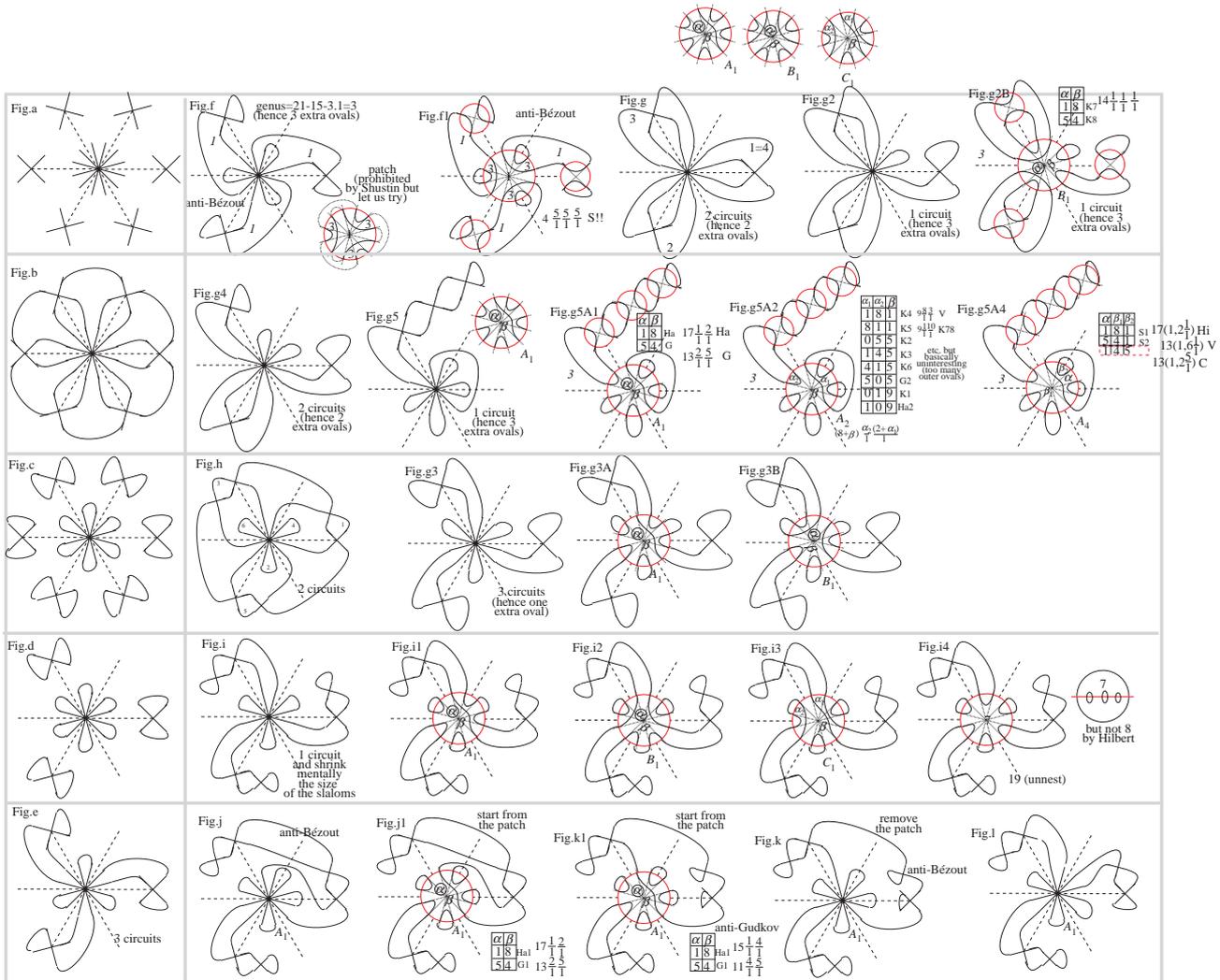,width=172mm} \captionskipAG
  \caption{\label{GabardDEGREE8_9:fig}%
Singular octics with a sextuple point invariant under varied
rotations} \figskip
\end{figure}

On the next figure we try to develop an introverted version of
Fig.\,g5. Philosophically speaking, we thought  initially that
symmetry is likely to be favorable for Harnack maximality yet in
reality it seems rather to be an extra constraint removing some
freedom.

\begin{figure}[h]\Figskip
%\vskip-1.2cm\penalty0
%\centering
\hskip-2.7cm\penalty0
\epsfig{figure=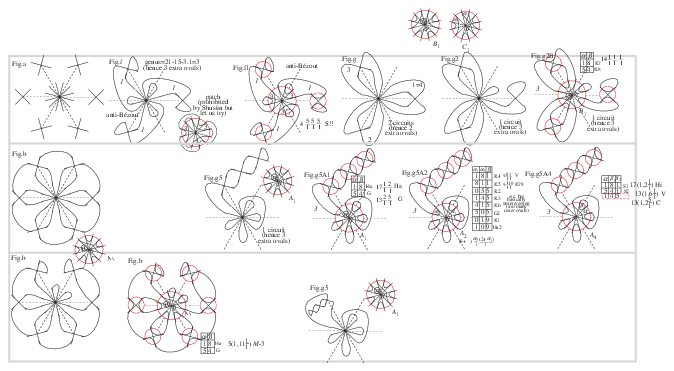,width=172mm} \captionskipAG
  \caption{\label{GabardDEGREE8_10:fig}%
Introverting} \figskip
\end{figure}

\subsection{Reflex de Pavlov  vs. complex de Gromov}

[26.07.13] Among the big scientists of the world, we have two
scientists like Pavlov and Gromov. The latter is a well-known
student of Rohlin () one of our main hero, well known for
semi-deep contributions on Riemannian geometry. Especially
important is Gromov's assertion that mathematics are so trivial
that once after a long effort we have understood the fully story
we do not take the pain to write down the details feeling almost
shameful of the triviality of the truth. This phenomenon we call
the complex de Gromov. Needless to say we feel not affected by the
latter. Needless to say, we feel quite besoffen but it is at such
time (as would say Ahlfors that the mathematical Empire appears in
its full great). It is only then that anodyne details swamp out so
that the true architecture become visible. Actually, our main
philosophical purpose is to give some text to TeX  so as to
picture more tomorrow along the action painting \`a la Jackson
Pollock tomorrow. As we said often in this text it is a pity that
TeX is so painful when it comes to
%%trace
integrates Figures. We, geometers are not linguists willing to
speak a lot, but we are rather pure observers trying to put  so
much pressure on the optical Universe as to force the latter to
deliver its archaic secret about Space, Time, Matter and
Immortality. This sounds pathetic, we confess, but what more
pathetic than masking the true expectation of our endeavor while
dilapidating social funds while bronzing along the coast of the
Aegean sea like Steve Small (annoying Jack Milnor). Science needs
workers and not capitalists. Sorry for all this boring
Naturphilosophie aus the best Bavarian Stock, but it is to remedy
TeX Page-Making (un)skills. H\"atten wir was be{\ss}eres verdient
wenn K.\,Friedrich Gauss zust\"andig gewesen w\"are the Knuth's
progammer
%%gewesen
zu sein? Adobe Illustrator is also not very
efficient for tracing all the pictures we had to trace with much
pain.
%%(more than optically appreciated).

[26.07.13] (Continued). So we must continue our method. Of course
the difficulty is to make a great catalogue as to remember what
has been already explored. Maybe did we already explored a triple
central point plus 3 quadruple points around it (then
genus=$21-3\cdot 6-3=0$).

Then the TeX compilator started to causes problem again.

[13.07.13] Another way to paraphrase our embryonic knowledge
emanating from our naive qualitative pictures of decomposing
curves is as follows.

\begin{Scholium}
Among the $6$ places de r\'esistance impeding the Wehrmacht (or
%rather
the Red-Army) to kill Hilbert's 16th in degree $8$, two of them
(namely $1\frac{1}{1}\frac{18}{1}$ and $14(1,2\frac{4}{1})$) are
accessible through small perturbation of a (hypothetical)
decomposing curve of degree $2+6$ and $1+7$ respectively, whilst
the $4$ other bosons look more mysterious in this respect. We
conjecture therefore (admittedly on hasty evidence)   the ultimate
solution of Hilbert's 16th in degree $m=8$ as materializing both
of those bosons algebraically, but killing the $4$ remaining ones.
\end{Scholium}

In fact our little pictures via decomposing curve might have been
implicitly known to Gudkov (compare e.g. the conjecture he
formulates at the end of the 1974 survey, cf. p.\,??). Albeit
heuristic, our method has the advantage to sweep out nearly all
schemes via a very homogenous method, while the Viro et cie.
method is more random and requiring a mixed patchworking of the
other contributors (Shustin, Korchagin, Chevallier, Orevkov) which
is rigorous yet fairly intricate artwork.

[14.07.13] It seems also worth noticing that the $M$-octic schemes
obtained via our ``floppy'' decomposing curve (our picture are not
algebro-geometric a priori) includes 2 of the 3 schemes
pseudo-holomorphically realized by Orevkov. Maybe this gives some
support that those 2 pseudoholomorphic models can be rigidified
(crystallized) in the algebraic category.

\subsection{Interlude: impact upon the Wiman-Rohlin-Gabard dream of
satellites (le th\'eor\`eme de Riemann rendu synth\'etique)}

[02.07.13] All the little  counterexamples to RMC described
previously in this text seems to affect/jeopardize the grand
programme we drafted in the Introduction of this essay (v.2). Our
thesis was essentially that the phenomenon of total reality should
explain all prohibitions of Hilbert's 16th problem, emphasizing a
great domination of Riemann-Schottky-Klein upon Hilbert. In
particular we expected a phenomenon of stability under satellites,
say of  being of type~I or of being maximal. Actually stronger
than all these conditions is the total reality of a scheme and
this should be the right condition forcing maximality of the
scheme, hence a series of ``criminal'' prohibitions (by killing
all enlargements). So even though Rohlin's maximality conjecture
looks foiled (essentially by Viro), another maximality principle
allied to total reality can still regulate the isotopic
classification of curves (Hilbert's 16th).

\subsection{Temptation of free-hand drawing}

[02.07.13] In view of the failure of Shustin's apple to reach
Harnack-maximality it seems tempting to create more ovals by
splicing the apple in two halves. However most curves so obtained
violates Gudkov's hypothesis validated by Rohlin (et cie), whence
a violent obstruction to dissecting Shustin's apple in two pieces.

\begin{figure}[h]\Figskip
%\vskip-1.2cm\penalty0
\centering
%\hskip-0.7cm\penalty0
\epsfig{figure=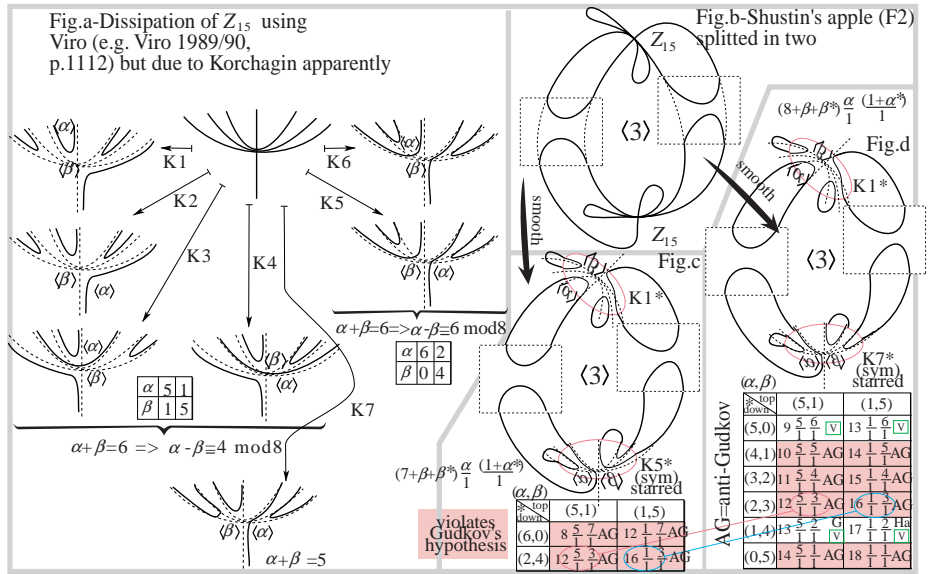,width=122mm}
\captionskipAG
  \caption{\label{ViroDEGREE8_SHUSTIN3:fig}%
  A free split of Shustin's apple, yet violating Gudkov}
\figskip
\end{figure}

\section{Post-Viro neoconstructivism: the latest tricks by
Korchagin, Chevallier, Orevkov}

\subsection{Korchagin's tsunami of 19 new $M$-schemes}

Without  getting too excited by this last scheme of Shustin let us
rather look at the next contribution namely that of Korchagin
which covers some 19 new schemes.

$\bullet$ Korchagin 1988
\cite{Korchagin_1988-dissertation}(=announcement in the Summary of
Candidate's dissertation, Leningrad) and published in 1989
\cite{Korchagin_1989-The-new-M-curves} (constructions of 19 new
schemes, raising Korchagin's score to $20$) probably via a variant
of Viro. Yes, but one variant of his own stock and alas also
difficult to follow as there is no projective picture given, but
just the ``Newton-Viro charts''  which we are not yet acquainted
with. Is it possible to make ``projective'' pictures of
Korchagin's curve, probably but it requires some working aptitude.
TO BE CONTINUED.

\subsection{Shustin's 5 new prohibitions of $M$-schemes}

$\bullet$ Shustin 1990/91 \cite{Shustin_1990/91-New-restrictions}
new prohibitions (of five $M$-schemes of the form $(1,a
\frac{20-a}{1})$ with $a$ adjusted as to verify Gudkov's
hypothesis, cf. Fig.\,\ref{Degree8-M-curve-TABLE:fig} for the
exact values $a=2,6,10,14,18$). This result probably involves a
variant of Hilbert-Rohn [Shustin being a direct student of Gudkov
and arguably among the living best-expert of this method]; but NO
it is rather by a topological method initiated by Viro 1984, whose
first variant seems even to go back to Rohlin, cf. p.\,424 of
\loccit); in fact Shustin claims only to prohibit 5 schemes [and
conjectures that the method employed also prohibit 2 additional
schemes (cf. Remark~6, p.\,443)]. (Omit please, this bracketed
text as this pertains only to $(M-1)$-curves.)

\subsection{21thest century heroes: Chevallier and Orevkov}

$\bullet$ Chevallier 2002 \cite{Chevallier_2002-4-M-curves} (new
constructions of $4$ schemes) probably via a variant of Viro's
method (yes but one requiring a very clever twist); Chevallier's
schemes are $a(1,2\frac{18-a}{1})$ for $a=2,5,13,16$, compare
Fig.\,\ref{Degree8-M-curve-TABLE:fig} for a visualization.

$\bullet$ Orevkov 2002
\cite{Orevkov_2001/02-classif-flexible-M-curves-degree-8}
prohibition of 2 schemes (namely $1\frac{3}{1}\frac{16}{1}$ and
$1\frac{6}{1}\frac{13}{1}$) by his revolutionary techniques
involving braids (Fox-Milnor, Lee Rudolph, quasi-positivity),
pseudoholomorphic curves, etc.

$\bullet$ Orevkov 2002
\cite{Orevkov_2002/XX-New-M-curves-degree-8} construction of one
scheme (namely $7(1,2\frac{11}{1})$), by using apparently
Grothendieck's dessins d'enfants.

\begin{figure}[h]\Figskip
%\vskip-1.2cm\penalty0
%\centering
\hskip-2.7cm\penalty0 \epsfig{figure=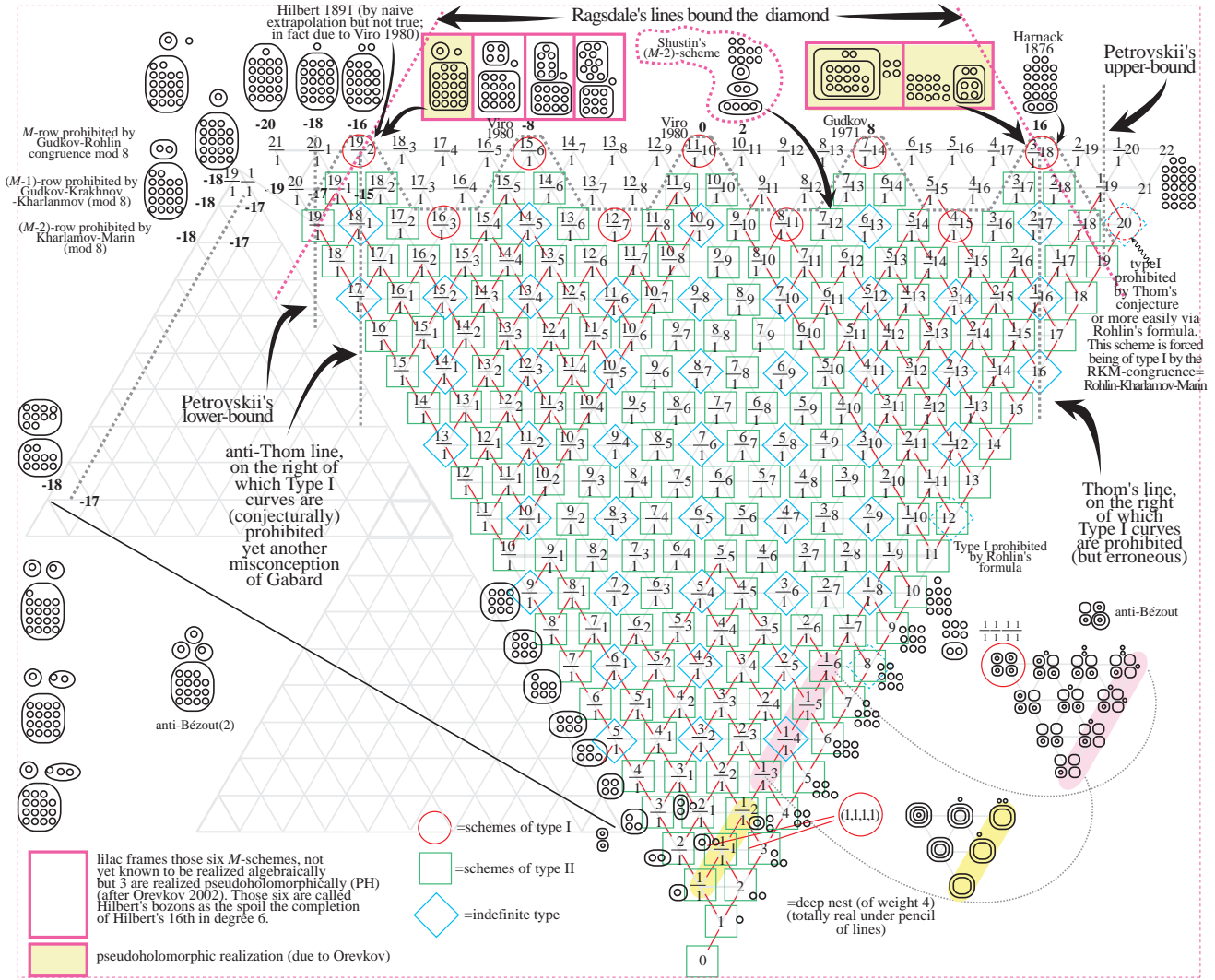,width=172mm}
\captionskipAG
  \caption{\label{Degree8-COPY:fig}%
  Naive pyramid and the six bosons}
\figskip
\end{figure}

At this stage (and the situation does not changed since) there
remains precisely 6 schemes (depicted on
Fig.\,\ref{Degree8-COPY:fig}) which are not yet known to be
realized. For their exact geography in the universe of all
%%%%%%%%%102
104 schemes, cf. again Fig.\,\ref{Degree8-M-curve-TABLE:fig}.
At this stage the proof of the (fragmentary)
theorem~\ref{Hilbert-16th-deg-8-M-curve(Viro-Orevkov):thm} is
completed.
\end{proof}

[30.04.13] We summarize this intricate system of contribution by a
table (Fig.\,\ref{Degree8-M-curve:fig}):

\begin{figure}[h]\Figskip
%\vskip-1.2cm\penalty0
\centering
%\hskip-0.7cm\penalty0
\epsfig{figure=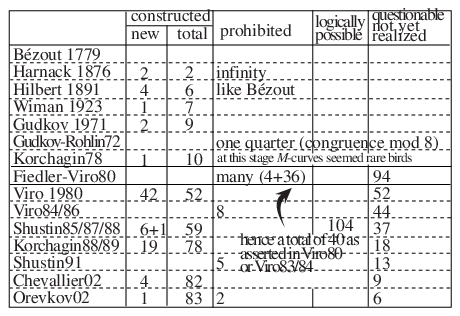,width=82mm} \captionskipAG
  \caption{\label{Degree8-M-curve:fig}%
  Histogram of Hilbert's 16th for $M$-octics}
\figskip
\end{figure}

Alas this does not clarify much the situation and we must of
course do our own exhaustive table of schemes to understand better
what happens. We were guided by the table in Orevkov 2002
\cite{Orevkov_2001/02-classif-flexible-M-curves-degree-8}, but
working out our own map was still necessary to understand better
the architecture.

[01.05.13] To gain a better understanding let us make a table of
all $M$-schemes in degree 8 (such tables are designed in Viro 1980
\cite{Viro_1980-degree-7-8-and-Ragsdale}, Viro 1984/84
\cite{Viro_1983/84-new-prohibitions}, Viro 1986
\cite[p.\,78]{Viro_1986/86-Progress}, and Orevkov 2002
\cite{Orevkov_2001/02-classif-flexible-M-curves-degree-8}), and we
try to stay close to their mode of depiction. Just a little
warning the values of $\chi$ are inverted in comparison to our
convention used on pyramids (e.g. Fig.\,\ref{Degree8-COPY:fig}).
Especially useful is Orevkov's trick to inform the reader of the
original builder of the curve [this trick was already used in Viro
1980 \cite{Viro_1980-degree-7-8-and-Ragsdale}, but in a less
efficient fashion, e.g. why calling Hilbert=(12), and
Harnack=(11)?]. Yet, at some place it seemed to us desirable
(especially in the case of Korchagin) to precise further the date
of fabrication. So for instance K78 means Korchagin 1978
\cite{Korchagin_1978}, where only one scheme was constructed (by a
variant of Brusotti), while in contrast K alone correspond to the
19 schemes constructed later by Korchagin in 1989
\cite{Korchagin_1989-The-new-M-curves}. Further Orevkov's table is
more complete than Viro's, yet contains at some places misprints.
We hope that our version combines the advantages of both
tabulation, while pointing the mistake to be found in Orevkov's.

\begin{figure}[h]\Figskip
%\vskip-1.2cm\penalty0
\centering
%\hskip-0.7cm\penalty0
\epsfig{figure=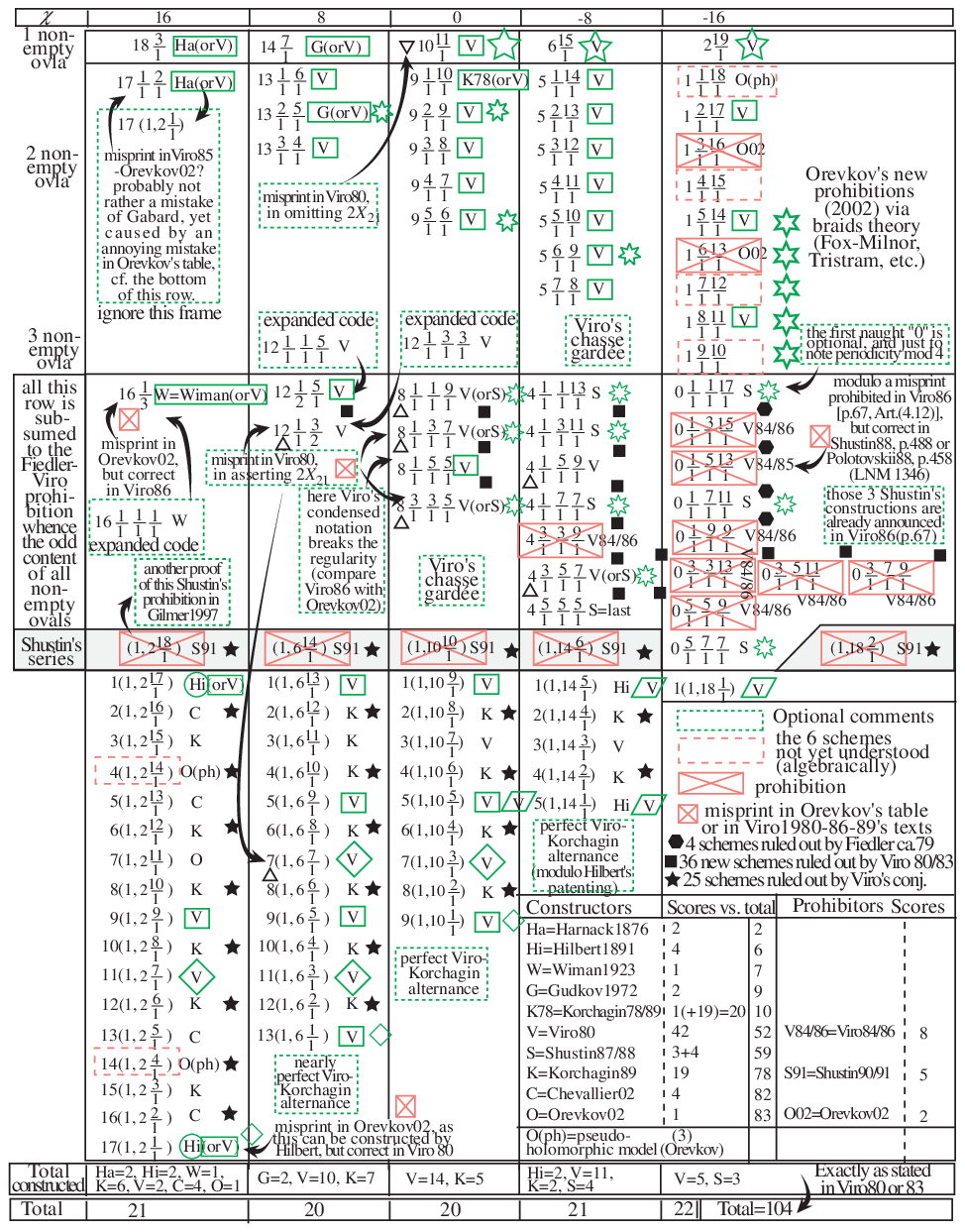,width=122mm}
\vskip-15pt\penalty0
  \caption{\label{Degree8-M-curve-TABLE:fig}%
  The Viro-Orevkov table of $M$-octics: all the
  %%%%%102
  104 Dalmatians
  forming the universe of logically possible schemes after the Fiedler-Viro
  obstruction (announced in Viro
  1980 \cite{Viro_1980-degree-7-8-and-Ragsdale} and detailed in
  Viro 1983/84 \cite{Viro_1983/84-new-prohibitions}).}
\figskip
\end{figure}

Further at the bottom of the first column, it seems that Orevkov
2002 erroneously ascribes to Viro a scheme (namely
$17(1,2\frac{1}{1})$) constructible via Hilbert's method (compare
%%%our Fig.\,\ref{HilbGab2:fig}).
the appropriate Fig. in v.2 of Ahlfors). By the way, this mistake
seems compatible with the fact that on Orevkov's table Hilbert
scores only 3 schemes, while even in Viro 1980 \cite[p.\,568,
Table~1]{Viro_1980-degree-7-8-and-Ragsdale}, Hilbert scored 4
schemes, and incidentally the scheme in question is correctly
ascribed to Hilbert 1891 (cf. 7th line of the first column, in
Viro 1980 \loccit).

Further our map is slightly more practical than Orevkov's as it
also shows the prohibitions. As 83 schemes are constructed and
$8+5+2=15$ are prohibited by V,S and O respectively, yielding a
total of $83+15=98$ so in comparison to the 104 of the
Fiedler-Viro universe,
%%$\Omega$,
6 schemes are left undecided.
%%%%%%%%
\iffalse
%
In fact our universe (i.e. all schemes depicted on our table)
comprises only 102 schemes and not 104 as asserted in Viro 1980
\cite{Viro_1980-degree-7-8-and-Ragsdale}. When this is realized
[and if we remove one credit that Viro stolen at Hilbert (so that
Hi=4, and V=41, not 42!)]\footnote{Omit this bracket as it is a
confusion coming from Orevkov's misprint.} we arrive at $83$
schemes constructed, and $6+5+2=13$ prohibited (so $96$ are fixed)
and $102-99=6$ many remains undecided, which is the present
state-of-the-art in Orevkov 2002
\cite{Orevkov_2001/02-classif-flexible-M-curves-degree-8}. We hope
this being the correct explanation. Note in Shustin 1990/91
\cite{Shustin_1990/91-New-restrictions} the size of the universe
is correctly given as $78$ (constructed), plus 5 prohibited (by
Shustin) and 18 types remain in question. Adding gives
$78+5+18=101$. [not quite right!]
\fi

By the way it should be noted that Viro's table (in Viro 1980
\cite{Viro_1980-degree-7-8-and-Ragsdale}) contains some anomalies,
for instance certain schemes are misplaced depending on the value
of $\chi=p-n$ (e.g. all the last 8 schemes of the series
$(p,n)=(19,3)$ should be moved to $(p,n)=(11,11)$). In this
respect Orevkov's table (2002) is much more reliable. (At any rate
we believe that our table is the most accurate one.) Upon
comparing Viro's 1980 table with ours (or Orevkov's 2002) it seems
that these 8 schemes is the sole inaccuracy in Viro80's table (and
probably this can be ascribed to the typographer who otherwise
would have been much annoyed to split the table in two rows of
equal heights).

Having checked this properly, we see that Viro 1980 (compare his
table) constructs for $\chi=16,8,0,-8,-16$ respectively
$2+10+14+11+5$, hence a total of $42$ schemes (exactly as he
asserts).

%%%%%%
\iffalse As we said it is quite puzzling that Viro estimates the
size of the universe of all $M$-octics to 104. We were fairly
convinced that this was just a misprint in Viro 1980
\cite{Viro_1980-degree-7-8-and-Ragsdale}, but the ``misprint''
reappears in Viro 1983/84
\cite[p.\,416]{Viro_1983/84-new-prohibitions}. I know no
explanation for the moment, except the naive one that Viro
reported the misprint from one paper to the other. The correct
count appears in Risler 1992 \cite{Risler_1992-Viro-BOURBAKI} yet
without serious explanations. By the way the tone of Risler's text
makes it apparent that he did not checked all the details of this
(Russian) census, yet we are of course all in greatest admiration
for those Russian achievements. \fi

At several occasions (during the 1980's) Viro advanced the
following conjecture (1980, 1983 \cite[p.\,416,
2.3.B]{Viro_1983/84-new-prohibitions}, and even 1986
\cite{Viro_1986/86-Progress}):

\begin{conj} {\rm (Viro's conjecture, disproved by Korchagin, but
partially verified on 5 cases by Shustin)}.---If $\al (1, \be
\frac{\ga}{1})$ is the real scheme of an $M$-curve of degree $8$
with $\gamma\neq 0$, then $\al$ and $\ga$ are odd
%%%numbers.
integers.
\end{conj}

This conjecture (wrote Viro 1983/84
\cite[p.\,416/17]{Viro_1983/84-new-prohibitions}) ``arose as a
consequence of unsuccessful attempts to realize the schemes which
are ruled out by it. Despite considerable efforts, no
counterexample to this conjecture has yet been constructed.''
Viro's conjecture would have ruled out all the 25(=9+7+5+3+1)
schemes marked by stars ``$\bigstar$'' on
Fig.\,\ref{Degree8-M-curve:fig}. Viro's conjecture turned out to
be generically false, with many counterexamples due to Korchagin
 ($4+6+4+2=16$ many, the magic number of this theory and
the number chosen by Hilbert, and Rohlin) and 2 counterexamples
due to Chevallier (2002), yet $5$ corroborations of Viro's
intuition by Shustin 1990/91
\cite{Shustin_1990/91-New-restrictions}, and 2 cases which are
still open nowadays. So Viro's conjecture is true with feeble
probability $\ge 5/25=1/5=0.2$ but at most with probability
$7/25=0.28$ [reminding Pl\"ucker's count of the bitangents to a
quartic, or Milnor-Kervaire count of smoothness structures on
$S^7$] in case both bosons are prohibited. Of course Korchagin
(and Chevallier) merely employ a (clever) variant of Viro's
method. It is therefore not clear presently if the 6 bosons will
 succumb under another clever twists of Viro's method, or if
they are new obstructions.

[02.05.13] It is hard to decipher any regularity in the
Viro-Shustin-Korchagin-Orevkov table
(Fig.\,\ref{Degree8-M-curve-TABLE:fig}). The combinatorics of this
table is best appreciated by considering it as constituted of 3
pyramids (triangles), two of them having an edge ``doubled''. The
first pyramid (schemes with 2 nonempty ovals) contains 4
mysterious bosons, namely $1\frac{a}{1}\frac{19-a}{1}$ for
$a=1,4,7,9$. Assuming that all these curves are not realized there
would be some symmetry in the first pyramid. Conversely if the
latter is symmetric, then at least the central schemes with
$a=4,7$ would be prohibited by extending Orevkov's prohibition
symmetrically. However in the second pyramid where the
classification is complete, by virtue of the Viro obstruction
(exposed by Korchagin-Shustin) there is a severe lack of symmetry,
and apparent chaos is reigning. Of course ``Gott w\"urfelt nicht''
(prose of A. Einstein, courtesy of Pharouk Garidi) and so some
hidden symmetry
%%%must exist.
can prevail. One should not exclude the possibility that the
Viro-Shustin con-census on this second pyramid contains
%some
mistakes (this occurred to the best workers, remember Gudkov's
saga 1954/69/72). Finally the 3rd pyramid involving schemes of the
form $a(1,b \frac{c}{1})$ with $a+b+c=20$ is nearly settled,
modulo 2 bosons which looks symmetrical. Like in degree 6, what
could cause symmetry is the phenomenon of eversions. So maybe the
Viro-Shustin asymmetry at the ``top'' of the 2nd pyramid could
represent an obstruction to eversions, or viceversa eversions
could detect an anomaly in the present census. (As we learned from
Kharlamov (letter in
%%%Sec.\,\ref{e-mail-Viro:sec})
v.2 of Ahlfors) and Viro it seems that Gudkov's referee Morosov
suggested so a possible mistake in  Hilbert's immature census.) On
more mature thinking the lack of symmetry is more a defect of our
table than an intrinsic feature of the problem. So it can be an
interesting problem to find a better diagrammatic than our table
respecting some symmetry, and giving full swing to the
%%%%% 6
 Viro's eight (sporadic) prohibitions.

Another crucial philosophical aspect concerns the constructions.
Viro's original method constructed 42 types of $M$-schemes.
Admittedly all what came next to Viro's breakthrough, i.e.
Shustin's 7 schemes, Korchagin 19 (new) schemes, and Chevallier's
4 schemes, are merely variants of Viro's method (yet along very
clever twists). An exception is perhaps Orevkov's construction
(2002) of the scheme $7(1,2\frac{11}{1})$ by a method
%which is
probably fairly distant from Viro's. So the ubiquity of Viro's
method seems very slightly attacked, although it is of course not
clear if Orevkov's scheme cannot be cooked \`a la Viro. {\it
Added\/} [05.10.13].---This comment looks immature, as looking at
Orevkov's paper one sees that is much akin to Chevallier's, which
is pure Viro theory, yet for another singularity than $X_{21}$.

In recent literature , we often read that nearly all objects may
be constructed along Viro's method (cf. e.g. Shustin 2005
\cite{Shustin_2005-Patch-arXiv} where we read (p.\,3): ``In
1979-80, O.~Viro [29,30,31,32] invented a patchworking
construction for real non-singular algebraic hypersurfaces. We
would like to mention that almost all known topological types of
real non-singular algebraic curves are realized in this way.'')
Cf. also Itenberg 2002
\cite{Itenberg_2002-Informal-notes-Constr-real-alg-var} (p.\,3)
where we read: ``Almost all the constructions in topology of real
algebraic varieties since 1979 use the Viro method.''

\subsection{Naive questions on $T$-curves}

For instance via the $T$-construction of Itenberg-Viro(-Orevkov)
is it possible to gain Orevkov's scheme? Maybe experimentally
nobody ever succeeded, but is there a proof that it cannot be
obtained by this recipe? Recall that Haas's result (1997
\cite{Haas_1997-Ragsdale-for-maximal-T-curves}) implies that
maximal $T$-curves respect Ragsdale $\vert \chi\vert \le k^2$, but
this does not answer our question. Maybe in degree 8, the
$T$-construction can be programmed by a machine exploring all
distribution of signs (and triangulations!) so as to assert that
nothing more can be obtained by this device that what is already
tabulated on Fig.\,\ref{Degree8-M-curve-TABLE:fig}. In
contradistinction, it could be just a matter of time until some
machine traces a $T$-curve realizing a new $M$-octic. It would be
actually interesting to know which schemes are realized by
$T$-curves. Usually Itenberg's $T$-construction (which in Itenberg
1994 \cite{Itenberg_1994-Notes-T-curves/1996} is in part ascribed
to Orevkov) is generally presented as a special case of Viro's
method, yet probably is also somewhat stronger (or rather more
flexible) remind for instance Itenberg's breakthrough on
Ragsdale's problem (at the $(M-2)$-level). So it is fairly
probable that the $T$-construction affords more than the 42
schemes derived by Viro's original method (within the limited set
of dissipation envisaged in the 1980 paper).

\begin{ques}
Which $M$-schemes of degree $8$ are realized via the
$T$-construction? How many triangulations and distribution of
signs exist in degree $8$? Were they all listed and analyzed?
\end{ques}

{\it Insertion} [03.04.13] A paper contributing to this question
is De Loera-Wicklin 1998 \cite{De-Loera-Wicklin_1998}. There some
nice tables of schemes realizable via the $T$-construction are
presented with focus on the critical degree $m=8$ (the present
frontier of knowledge). However the quantity of $M$-schemes
created by a random computer search browsing apparently a  million
of $T$-curves in degree 8 (cf. p.\,213 of \loccit) created only 2
types of $M$-curves namely $18\frac{3}{1}$ (cf. box on their
Fig.\,8) and $17 \frac{1}{1}\frac{2}{1}$ (cf. box on their
Fig.\,9). Comparing with our Viro-Orevkov table
(Fig.\,\ref{Degree8-M-curve-TABLE:fig}),  those  schemes are
merely the first ever constructed namely
%those imputable to
by Harnack 1876. Hence,
%%%say the least,
%be frankly honest,
 the De Loera-Wicklin computer
search seems extremely disappointing. As they wrote on p.\,207:
``The obtainable schemes of degree 8 are not yet classified, and
no one knows which schemes are obtainable as $T$-curves.'' The
article contains a basic trick to keep track of all schemes on
only 3 tables which are of course not enough as we have at least
the 4 classes of pyramids depicted on the Viro-Orevkov table.
[Added in proof: I am not sure to understand myself here.] So it
is a good problem to know the number of sheets of paper required
to depict all schemes of degree 8. As a last remark on their
paper, it seems that the authors omit to use the RKM-congruence as
to infer that their scheme on Fig.\,8, with $(\al, \be)=(15,4)$ is
of type~I. Finally on p.\,211, we read a comment that we had some
pain to interpret properly: ``The data for regular triangulations
indicates that Rohlin's comment 20 years ago [18] continues to
hold in view of this new data: nothing so far contradicts the
conjecture that all real schemes not prohibited by known theorems
are of indefinite type.'' This is probably a complex way to
reformulate Rohlin's maximality conjecture that a scheme of type~I
kills all its enlargement. Of course their tables invite one to
play the same game as we did in degree 6, namely smoothing the
union of 4 ellipses to see what can be gained by this more
elementary method. Probably their table is also much weaker than
those compiled by Polotovskii 1988
\cite{Polotovskii_1988---classif-deg-8}.

\subsection{Riemann's gyroscopic total reality}

In sharp contrast it could be that all the 6 last $M$-schemes (or
at least a good portion thereof) are prohibited say by Riemann's
method of total reality (cf. Gabard 2013B
\cite{Gabard_2013B-Riemann's-flirt}) involving pencil of sextics
with $B=M+(m-4)=22+4=26$ basepoints oddly distributed on each oval
(so as to
%%%gain
grant one bonus-intersection on each oval to reach total reality
at $26+22=48=6\cdot 8$). Alas, as yet we were not even able to
recover the two Hilbert/Rohn prohibitions in degree $m=6$ by this
method (of the schemes $11$ and $\frac{10}{1}$ respectively).
However it may be argued that Hilbert's and Rohn's obstructions
derive very simply from Rohlin's formula.
%%%%(\ref{Rohlin-formula:thm}).
In degree $m=6$, Rohlin's formula
affords less obstructions than Gudkov's hypothesis but is it a
general feature? In degree 8, is there any scheme of our table
prohibited by Rohlin's formula? We presume not but this requires a
little exercise. Try for instance cavalier the (still open) scheme
$1\frac{1}{1}\frac{18}{1}$. Then $2(\pi-\eta)=r-k^2=22-16=6$ so
that $\pi-\eta=3$ and as $\pi+\eta=19$ (the number of edges in
Hilbert's tree) the equation is (uniquely) soluble as $\pi=11$ and
$\eta=8$. Further the signs-law affords no constraint (since there
no deep edges available for concatenation). The dream could be
that the dynamics of the electron(s) allied to the totally real
pencil puts some restriction upon complex orientations (via the
dextrogyration principle). Indeed any dividing (so in particular
$M$-) curve appears as a {\it gyroscope\/}\footnote{Coinage of
Emmanuel Boul\'e (ca. April 2013), the cousin of the writer.}
under a holomorphic map of Ahlfors or
Riemann-Schottky-Bieberbach-Grunsky respectively. More
specifically we have the phenomenon of total reality described in
Gabard 2013B \cite{Gabard_2013B-Riemann's-flirt} (as recalled just
above), and under such a sweeping one could try to infer the
structure of complex orientations (as one is able to do in the
trivial case of the deep nest via gyration along the pencil of
lines). From this knowledge one may expect a contradiction with
Rohlin's formula. To be very specific we can prescribe the 4 extra
basepoints either as a tower of 4 points concentrated on one oval
or split them apart in two groups of mass 2. Deciding which
distribution of basepoints is most instructive from the viewpoint
of complex orientation is a puzzle even for the writer.

So meta-principle:

\begin{Scholium} (Riemann's gyroscopic principle).---New
(and probably even old) obstructions in Hilbert's 16th (especially
for $M$-curves) may be derived by a conjunction of Riemann's
gyroscopic principle of total reality combined with Rohlin's
complex orientation formula.
\end{Scholium}

The philosophy is that Rohlin's formula alone explains the
prohibition of Hilbert and Rohn in degree $m=6$, but becomes quite
impuissant in degree 8 (at least if Gudkov's hypothesis is already
imposed as on the Viro-Orevkov table,
Fig.\,\ref{Degree8-M-curve-TABLE:fig}). Yet, perhaps when
%%%% aided
assisted by Riemann's gyroscopic principle (1857) then Rohlin's
formula gain more swing and could rule out one of the 6 bosons not
yet detected (compare with the metaphor by C. Taubes in Bull. AMS
ca. 1994/96 on ``particles'' hard-to-detect, like Higgs' bosons,
when it comes to worry about possible exotic smoothness structure
on the 4-sphere $S^4$).

To ensure total reality of an $M$-octic we have to distribute the
$B=M+(m-4)=22+4=26$ basepoints on the 22 ovals plus the remaining
4 as either a tower (skyscraper) of height 4, or 2 mini-towers
(twin-towers) of height 2. So we have up to continuous deformation
as many total pencil as $22$ (location of the skyscraper) or the
binomial coefficient $\binom{22}{2}=\frac{22\cdot 21}{2}=11\cdot
21=231$ many possible choices of total pencils. Which one (among
those $22+231=253$ many) is the most clever choice in order to
settle the question of the 6 bosons (e.g.
$1\frac{1}{1}\frac{18}{1}$) via Riemann's gyroscopic effect is
hard to predict. One should
%%%%be able to
first
%%%invest money in visualizing
be capable visualizing pencil of sextics as to infer
%%something
valuable information upon complex orientations. Of course it
should also be noted that a priori the real scheme does not
predestine (uniquely) the complex orientations (alias complex
scheme), remember Marin's example in degree $m=7$
%%%(Fig.\,\ref{Marin:fig}).
(cf. appropriate Fig. in v.2 of Ahlfors). It seems evident that
Marin's phenomenon prevails a fortiori in degree 8, and we cannot
expect unique determination of the complex orientations from the
real scheme.

Another remark is in order. We know that some schemes are realized
pseudo-holomorphically (Orevkov 2002), and that Viro's method
without convexity assumption [so called $C$-curves] leads to
pseudo-holomorphic curves (Itenberg-Shustin 2002
\cite{Itenberg-Shustin_2002} and/or 2003
\cite{Itenberg-Shustin_2003-Israel-J-M}). In the same source, it
is also proved that $C$-curves are flexible curves in the sense of
Viro 1986, hence subsumed to all classical obstructions of
topological origins, in particular Rohlin's formula (and even the
Gudkov-Rohlin congruence). Hence there is no chance that Rohlin's
formula alone fixes Hilbert's 16th in degree 8, and the full swing
should come by the adjunction of the totally real pencil
materializing a sort of transverse structure (\`a la Haefliger,
etc.). So Riemann's gyroscopic effect should play a tremendous
role in the final elucidation of Hilbert's 16th problem.

\subsection{The extended Polotovskii tables of octics
(deposited in VINITI? ca. 1985)}

\begin{figure}[h]\Figskip
%\vskip-1.2cm\penalty0
%\centering
%
\hskip-0.7cm\penalty0
\epsfig{figure=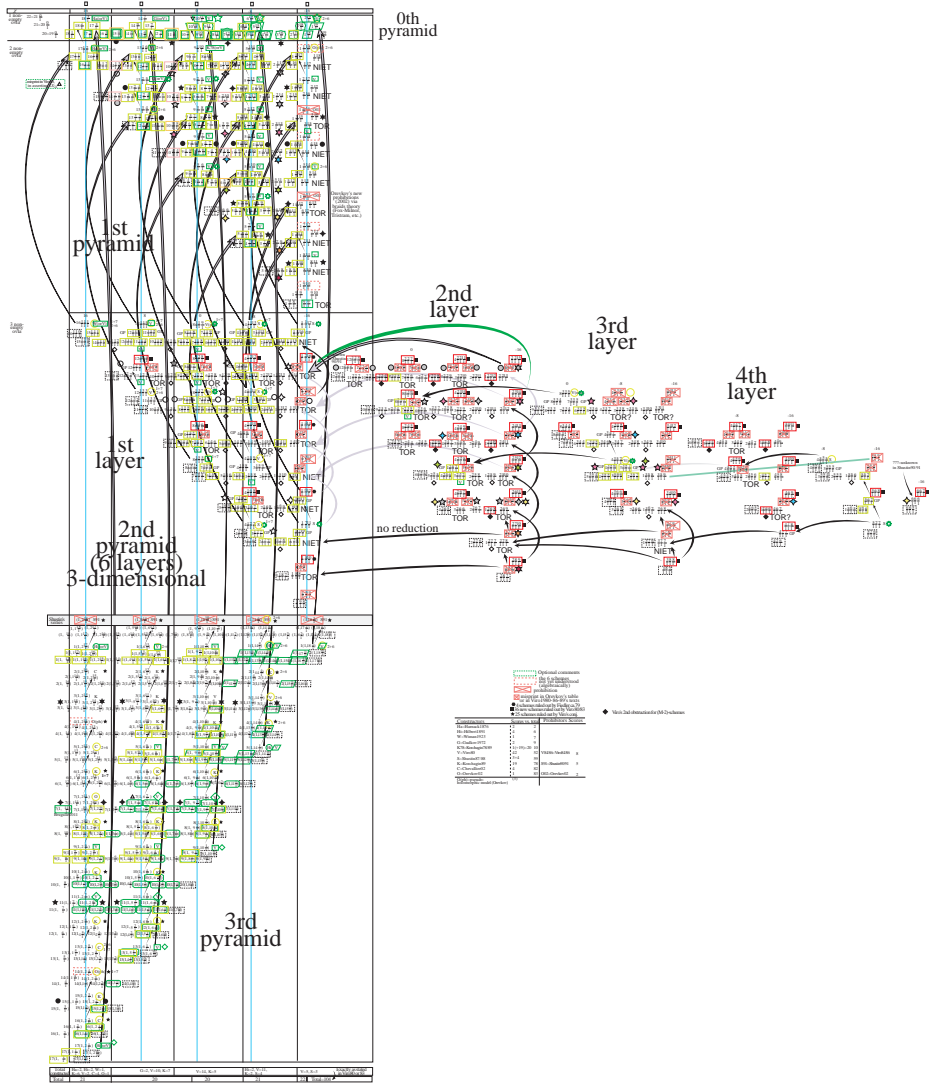,width=122mm}
\captionskipAG
  \caption{\label{Degree8-(M-i)-curve-TABLE:fig}%
  Viro-Polotovskii table for $m=8$ (Viro 1980,
  Polotovskii 1983 published in VINITI): hard-to-read but
enlargements on three sub-plates to be found subsequently
(Figs.\,\ref{Degree8-(M-i)-curve-TABLE_I:fig},\ref{Degree8-(M-i)-curve-TABLE_II:fig},\ref{Degree8-(M-i)-curve-TABLE_III:fig})}
\figskip
\end{figure}

\begin{figure}[h]\Figskip
%\vskip-1.2cm\penalty0
\centering
%
%\hskip-0.7cm\penalty0
\epsfig{figure=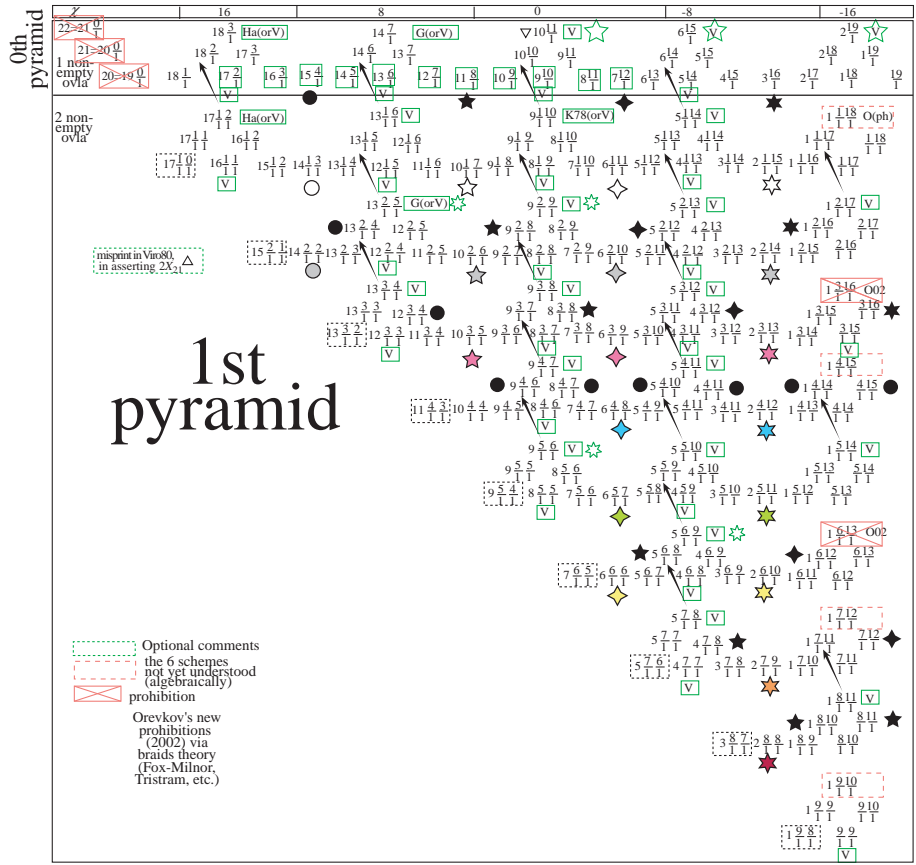,width=122mm}
\captionskipAG
  \caption{\label{Degree8-(M-i)-curve-TABLE_I:fig}%
  Zooming  the 1st pyramid} \figskip
\end{figure}

[10.05.13] To understand better Rohlin's maximality conjecture
(RMC) in degree 8 (where it is still open I think) it would be
valuable to trace an extended table showing also the $(M-1)$- and
$(M-2)$-schemes. [$\bigstar$ {\it Update}.---Meanwhile we think to
have refuted Rohlin's conjecture, cf.
Lemma~\ref{RMC:cter-example-via-Viro-1st-curve=BEAVER}.] All this
must be represented on the same plate in 3D so-to-speak. This is a
fairly good avatar of the degree 6 Gudkov pyramid to degree 8, yet
showing only the 3 top levels with resp. ca. 100, 200 and 400
schemes (recall the trinity of congruences allied to Gudkov {\it
et cie.}).
In fact our multi-pyramid
(Fig.\,\ref{Degree8-(M-i)-curve-TABLE:fig}) was constructed by
taking the Viro-Orevkov table of $M$-schemes
(Fig.\,\ref{Degree8-M-curve-TABLE:fig}) while extending downwards
to their servitude of $(M-1)$- and $(M-2)$-schemes. The plate so
obtained is fairly massive
 and it took us circa one
day just to dress its basic architecture (combinatorial structure)
with the help of Adobe Illustrator. Alas this plate can globally
only be consulted on a computer with moving resolution, so do not
worry if you see nothing on the paper. Understanding the full
architecture of this pyramid
%is of uttermost importance to
%understand
is tantamount settling Hilbert's 16th in degree 8. Especially
exciting is whether Rohlin's maximality conjecture for schemes of
type~I holds true in degree 8. [Same update as before, i.e. cf.
$\bigstar$ above.]
Of course once this global plate will be updated (as to know
exactly which schemes are realized we shall print all individual
pyramids separately as to appreciate on paper-scale those great
achievement of the Russian scholars (on both sides of Ural, i.e.,
Leningrad (Viro) and Gorki (Polotovskii, Shustin, Korchagin), plus
the more recent
%%%workers
contributors like Chevallier and Orevkov.

It is also a good exercise to
%%%understand
appreciate Shustin's disproof of the other sense of Rohlin's
maximality conjecture on the basis of this table
(Fig.\,\ref{Degree8-(M-i)-curve-TABLE:fig}). Indeed the first
Shustin scheme is $10\frac{1}{1}\frac{2}{1}\frac{4}{1}$ which is
of type~II by Arnold's congruence, but from the diagrammatic it is
fairly clear that this scheme is maximal since by Viro's
obstruction there is nothing above it. Actually, right above
Shustin's scheme under examination (i.e.
$10\frac{1}{1}\frac{2}{1}\frac{4}{1}$) we find
$11\frac{1}{1}\frac{2}{1}\frac{4}{1}$ which is actually prohibited
by a deep result of Shustin (cf. Shustin 1990/91
\cite[Thm\,2]{Shustin_1990/91-New-restrictions}), which can be
summarized as follows:

\begin{theorem} \label{Shustin-50-(M-1)-prohibitions:thm}  {\rm (Shustin 1990/91)}.---All trinested $(M-1)$-schemes barred on
Fig.\,\ref{Degree8-(M-i)-curve-TABLE:fig} by a red-cross rectangle
are prohibited. This involves a collection of exactly
%%%%50
fifty  $(M-1)$-schemes (prohibited by Shustin). Basically, this
Shustin obstruction can be summarized as follows. An
$(M-1)$-scheme is prohibited whenever it is immediately dominated
(on Fig.\,\ref{Degree8-(M-i)-curve-TABLE:fig}) by an $M$-scheme
prohibited by Viro's 1st law (imparity law for trinested
$M$-schemes), excepted when it derives from another constructible
$M$-scheme (constructed as a rule either by Viro or Shustin) as
depicted on Fig.\,\ref{Degree8-(M-i)-curve-TABLE:fig} via arrows.
This system of degenerations (each interpretable as a contraction
of an empty oval) explains all fifty obstructions of Shustin, but
leaves open 2 cases namely the $(M-1)$-schemes
$\frac{6}{1}\frac{6}{1}\frac{6}{1}$ and
$12\frac{2}{1}\frac{2}{1}\frac{2}{1}$ which though being
immediately dominated by a Viro-prohibited $M$-scheme, yet
non-dominated by a constructible scheme (for architectonic
reasons, i.e. cf. Fig.\,\ref{Degree8-(M-i)-curve-TABLE:fig}) are
not  known to be prohibited. So those two case are 2 exceptions
toward the synthesizing Shustin's 50 prohibitions (to an uniform
prohibition of 52, like the Boeing B-52). (I am going like a
Boeing.--Joke of Cornelius de Boeck.)
\end{theorem}

Actually, it must be remembered that this deep result of Shustin
is not even logically required to refute RMC, because even if the
scheme $11\frac{1}{1}\frac{2}{1}\frac{4}{1}$  existed
(algebraically) then it would  constitute itself a counterexample
to RMC, by virtue of Klein's congruence and the maximality of the
scheme granted by Viro's imparity law (VIL). Indeed, it must be
noted that this $(M-1)$-scheme cannot be enlarged into a quadruply
nested one (without violating B\'ezout or total reality, i.e.
maximality of the 2nd satellite of the quadrifolium). So it seems
that the sole enlargement possible involves adding nested ovals,
yet those operations (diminishing $\chi$ by $1$) push the scheme
outside of the GKK-range (Gudkov-Krakhnov-Kharlamov) as may be
visualized on the main-Table
(Fig.\,\ref{Degree8-(M-i)-curve-TABLE:fig}).

So we see that Shustin's example can be replicated at several
places whenever we look at an $(M-2)$-scheme where the Gudkov
sawtooth is a broken-line and under the $M$-schemes prohibited by
Viro. So once the geography is fixed it becomes nearly trivial to
get a birdseye view upon all schemes coherently organized into a
pyramid. In particular we can look at a ``broken'' $(M-2)$-scheme
below a Fiedler prohibited scheme, for instance
$\frac{1}{1}\frac{2}{1}\frac{14}{1}$ and Shustin's argument
probably applies as well here (after constructing the scheme via
Viro's method or a suitable variant thereof). If this can be done
we see that from the Germanic angle of view:

\begin{Scholium}
Fiedler's special case of Viro's imparity law certainly suffices
to disprove Rohlin's maximality conjecture and so Klein vache as a
byproduct.
\end{Scholium}

It is also interesting to wonder about Rohlin's maximality
conjecture, i.e. a type~I scheme is maximal. At first sight the
conjecture looks trivially true for geographical reasons. More
precisely in the above pyramid
(Fig.\,\ref{Degree8-(M-i)-curve-TABLE:fig}) we see always Gudkov's
sawtooth undulating piecewise linearly between maxima at
$M$-schemes and minima at $(M-2)$-schemes satisfying the
RKM-congruence $\chi\equiv_8 k^2+4$ (forcing type~I). Hence all
those RKM-schemes (in the depression of the sawtooth) are of
type~I and (at first sight) maximal by virtue of the Gudkov-Rohlin
and GKK congruences (which diagrammatically forbid direct
enlargements of those schemes). However as we shall see there may
be indirect enlargement somewhat more perfidious to visualize on
the planar model of our pyramid. Actually upon looking carefully
at the architecture of this main-table we see that the (unnested)
$(M-2)$-scheme $20$ is RKK and so of type~I (if it existed) and so
would kill the 2 enlargements $21$ and $22$, thereby reproving
nearly Petrovskii's bound $\chi\le 19$ via Rohlin's maximality
principle (RMC). Note yet two objections: first the scheme $20$
does not exist as may be inferred either from Petrovskii's bound
or from Rohlin's formula $2(\pi-\eta)=r-k^2$, hence as the left
side vanishes we have $r=k^2=16$ many ovals and not twenty.
Further the scheme $20$ would kill as well a myriad of other
$(M-1)$-schemes on the row right above it (like $18\frac{2}{1}$,
etc.) as well as larger $M$-schemes which are known to exist by
classical (e.g. Harnack) or neoclassical methods (i.e. Viro).

\subsection{Enlarged plates I, II and III}

In this subsection we just reproduce the Hauptfigur
(Fig.\,\ref{Degree8-(M-i)-curve-TABLE:fig}) on several subplates
viewable at the normal printed size.

Our first plate is the first pyramid.

Next, the second plate is the 2nd pyramid. To depict it at an
acceptable scale, we had to break it like a snake.

\begin{figure}[h]\Figskip
%\vskip-1.2cm\penalty0
\centering
%
%\hskip-0.7cm\penalty0
\epsfig{figure=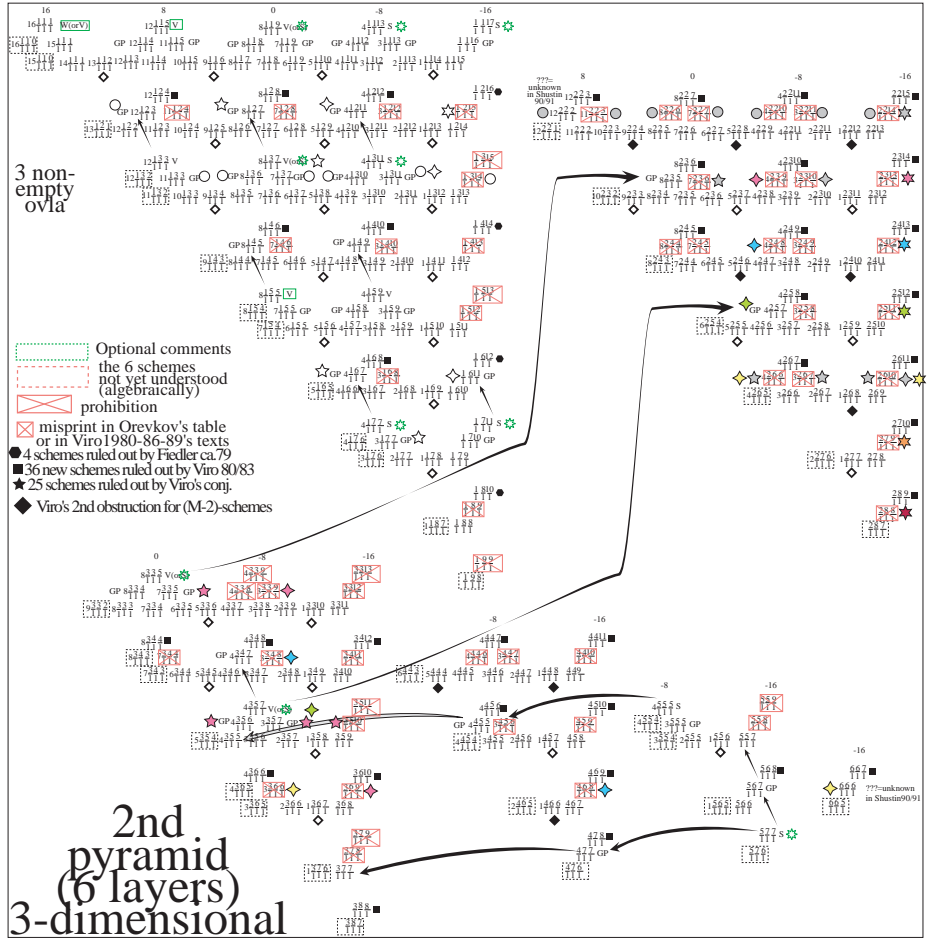,width=122mm}
\captionskipAG
  \caption{\label{Degree8-(M-i)-curve-TABLE_II:fig}%
  Zooming the 2nd pyramid} \figskip
\end{figure}

Finally the 3rd plate is the 3rd pyramid.

\begin{figure}[h]\Figskip
%\vskip-1.2cm\penalty0
\centering
%
%\hskip-0.7cm\penalty0
\epsfig{figure=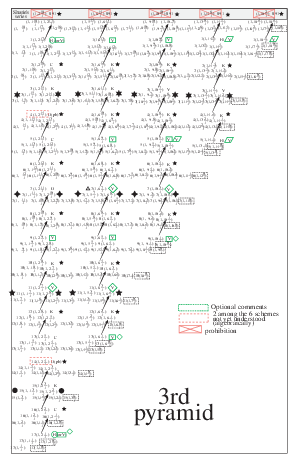,width=102mm}
\captionskipAG
  \caption{\label{Degree8-(M-i)-curve-TABLE_III:fig}%
  Zooming the 3rd pyramid} \figskip
\end{figure}

Propagating the above reasoning, consider a geographical avatar of
the scheme $20$, namely two rows below the $(M-2)$-scheme
$14\frac{2}{1}\frac{2}{1}$ which is RKM. Granting RMC this should
kill all its enlargements, it particular those where the
additional oval separates the 14 outer ovals in two groups. Hence
all schemes of the form $\al \frac{\be}{1}\frac{2}{1}\frac{2}{1}$
where $\al+\be =14$ (or even $15$) are killed. Those schemes are
located in the second multi-pyramid (i.e. the 2nd row of the large
table). This conclusion inferred from RMC is in agreement with
Viro-Fiedler's prohibition via complex orientations. Actually, the
scheme $12\frac{2}{1}\frac{2}{1}\frac{2}{1}$ should likewise be
prohibited, yet this is still unknown as we know since Shustin
90/91 (\ref{Shustin-50-(M-1)-prohibitions:thm}).

On applying the same method to the RKM-scheme
$10\frac{3}{1}\frac{5}{1}$ we get enlargements by looking where
the sub-symbol $\frac{3}{1}\frac{5}{1}$ appears in an
$(M-1)$-scheme. Aided by our map, we locate so the schemes
$7\frac{3}{1}\frac{3}{1}\frac{5}{1}$,
$4\frac{3}{1}\frac{5}{1}\frac{6}{1}$,
$3\frac{3}{1}\frac{5}{1}\frac{7}{1}$. However on behalf of Shustin
90/91, it seems that all those 3 schemes were constructed in
Goryacheva-Polotovskii 1985 \cite{Goryacheva-Polotovskii_1985}
(abridged GP in the sequel). Hence existence of this RKM-scheme
(with principal symbol $(3,5)$) would corrupt Rohlin's maximality
conjecture. So it may be reasonable to expect that our RKM-scheme
does not exist algebraically. ({\it Update} [30.06.13].---With
some more maturity, it seems rather more plausible that RMC is
foiled, as we found even simpler counterexample to it, cf. e.g.
Lemma~\ref{RMC:cter-example-via-Viro-1st-curve=BEAVER}.)

The game can be continued straightforwardly, cf. our map where
sometimes we encounter no obstruction. However sometimes we get
conflicts between RMC and construction claimed by Polotovskii (and
cie.). So for the scheme $6\frac{5}{1}\frac{7}{1}$ we look for the
principal symbol $(5,7)$ and discover on the table the scheme
$4\frac{2}{1}\frac{5}{1}\frac{7}{1}$, and also
$3\frac{3}{1}\frac{5}{1}\frac{7}{1}$. Both those $(M-1)$-schemes
are constructed by GP, hence we get either a conflict with RMC or
a prohibition of the initial RKM-scheme.

Along the same mode-of-thinking the RKM-scheme
$6\frac{6}{1}\frac{6}{1}$ offers an interesting twist. Looking at
$(M-1)$-enlargement of its principal symbol $(6,6)$, we find the
schemes
$4\frac{2}{1}\frac{6}{1}\frac{6}{1}$,
$3\frac{3}{1}\frac{6}{1}\frac{6}{1}$, and finally
$\frac{6}{1}\frac{6}{1}\frac{6}{1}$.
Interestingly the first 2 schemes are prohibited by Shustin 90/91
\cite{Shustin_1990/91-New-restrictions}, while the third is not
known to be realized. Hence there is some evidence that the
original RKM-scheme exists.

Applying the same (heuristic) method to the RKM-scheme
$S:=14\frac{1}{1}\frac{3}{1}$ gives enlargements of the shape $\al
\frac{\be}{1}\frac{1}{1}\frac{3}{1}$ with $\al+\be =14$ (or even
15). Those schemes are located on the first layer of the 2nd
(3-dimensional) pyramid (specifically in its 2nd and 3rd rows).
However some of those schemes (or their direct $M$-enlargements)
are constructed by either Viro (or Shustin's method), for instance
$12\frac{1}{1}\frac{3}{1}\frac{3}{1}$ (is claimed by Viro though
we were incompetent enough to miss this as yet) or
$8\frac{1}{1}\frac{3}{1}\frac{7}{1}$ (claimed by Viro, and which
we were able to manufacture following Shustin). We arrive at an
interesting psychological tension. Several logical issues are
possible.

$\bullet$ First, it could be that the scheme in question (i.e.
$14\frac{1}{1}\frac{3}{1}$) does not exist. (Added [30.06.13]: for
a weak heuristic construction \`a la Viro, cf. our
Fig.\,\ref{ViroDEGREE8_2:fig}f.) This nonexistence looks a priori
quite
%unlikely
improbable as this scheme has nearby two (companion) $M$-schemes
coming  either from Harnack's $17\frac{1}{1}\frac{2}{1}$ or Viro's
$13\frac{1}{1}\frac{6}{1}$ (cf. the Ha and V certificates/patents
of construction on the table) which were both constructed earlier
in this text by Viro's dissipation method
(Fig.\,\ref{ViroDEGREE8:fig}) in its most elementary incarnation
(dissipation of 4 coaxial ellipses).

$\bullet$ Second it could be that Shustin's construction is
erroneous. (Remind that as yet we were not able to realize
$12\frac{1}{1}\frac{3}{1}\frac{3}{1}$ by Viro's method as asserted
in Viro 1980.)

$\bullet$ Third it could be that we located a counterexample to
RMC (Rohlin's maximality conjecture).

So see clearer it would be nice to construct the above scheme $S$
(i.e. $14\frac{1}{1}\frac{3}{1}$). We think that this should be an
easy matter, but let us look at related scheme probably even
easier to construct.

Let us repeat the same method to the RKM-scheme
$S_0:=15\frac{4}{1}$. Its enlargements  have the shape $\al
\frac{\be}{1}\frac{4}{1}$ with $\al+\be =15$ (or even 16). Those
are encountered in the 1st pyramid especially in its fourth row
(while being depicted by black circular bullets). So the same
conflict with Viro's method is obtained. Hence, either Viro's
method is wrong (unlikely but personally we confess to have not
yet understood its mechanism in all details), or Rohlin's
maximality conjecture is false, or finally, the scheme
$S_0=15\frac{4}{1}$ does not exist (algebraically). However since
$S_0$ lies in the depression of the sawtooth between two
$M$-schemes due to Ha=Harnack and G=Gudkov resp., yet most easily
constructed via Viro (namely Ha=$18\frac{3}{1}$ and
G=$14\frac{7}{1}$), it is likely that the scheme in question
exists (albeit not readily obtained by the contraction principle
for empty ovals). [{\it Update} [30.06.13].---For one construction
of this scheme cf.
Fig.\,\ref{ViroDEGREE8_TRICKY_SUITE_XXX_15_4-1:fig}k.]

\begin{figure}[h]\Figskip
%\vskip-1.2cm\penalty0
\centering
%
%\hskip-0.7cm\penalty0
\epsfig{figure=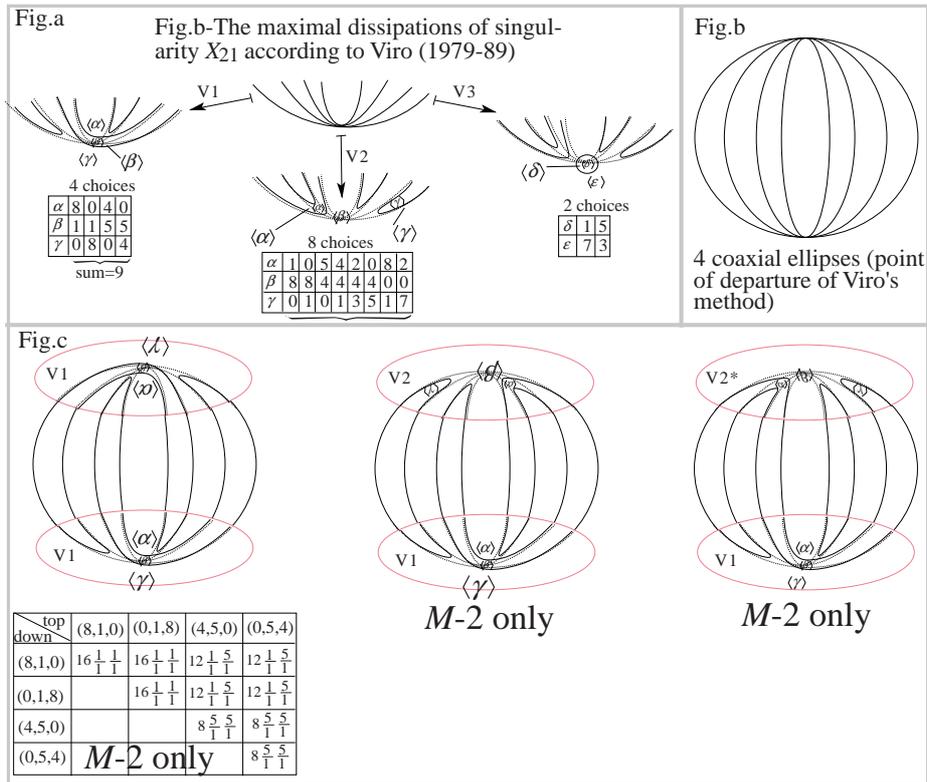,width=122mm} \captionskipAG
  \caption{\label{ViroDEGREE8_BIS:fig}%
  Viro's method for $(M-2)$-schemes}
\figskip
\end{figure}

Let us now try to check more pragmatically this point
(realizability of $S_0$ defined above) to accentuate the paradox.
So we look again at our earlier Fig.\,\ref{ViroDEGREE8:fig} of the
elementary Viro method with 4 ellipses. A priori there is 2
options to create  $(M-2)$-curves. Either employ a non-maximal
dissipation of Fig.\,\ref{ViroDEGREE8:fig}d where V$i$ is glued
with its symmetric V$i^\ast$, or use a maximal dissipation of an
asymmetrical gluing like on Fig.\,\ref{ViroDEGREE8:fig}c. Let us
first explore this 2nd idea. First, look at the left part of
Fig.\,\ref{ViroDEGREE8:fig}c (reproduced below as
Fig.\,\ref{ViroDEGREE8_BIS:fig} for optical convenience) where we
actually already tabulated the possible schemes in the tablet
right-below that figure. Clearly as $\be$ is $\ge 1$ we have (at
least) two nonempty ovals and so our scheme $S_0$ is not realized.
Next look at the middle part of the same Fig.\,c. Then we have a
contorted oval (union of essentially 3 lunes). For the bottom
$\be$ we can only choose $1$ or $5$. The latter being too much for
$S_0$ (where 4 ovals are nested) we choose $\be =1$. But then  the
top $\al$ should be 3, which is however not a permissible value
(cf. Fig.\,a). Let us now examine the 3rd configuration of
Fig.\,c. Again for the bottom $\be$ we are forced to take $\be=1$.
So our desideratum is to choose $\al+\ga=3$ on the top dissipation
V$2^\ast$. However a glance at Fig.\,a shows that $\al+\ga$ can
only be $1,5,9$, and so we fail constructing the desired curve.

Of course all our game looks ancient Russian games \`a la
Viro-Polotovskii
%%%in the late 70's or
of the early 80's. In
%%our background memory there
literature, it is often asserted that Rohlin's maximality
conjecture resisted all assaults of Viro's 1980 Red October
revolution (cf. e.g. Polotovskii 1992
\cite{Polotovskii_1992---classif-decomposing-curve-Rennes}). So it
seems that the scheme $S_0=15\frac{4}{1}$
%%does
should not exist.

Recall at this stage that Viro has another obstruction for
$(M-2)$-curves stating that if the content of a scheme with 3
nonempty ovals is divisible by 4 then two of the nested numbers
are odd and one is even.
%(cf. \ref{Viro-Fiedler-prohibition:thm}
%for the exact statement).
This does not (alas) apply to the case
at hand where there is only one nonempty oval. Yet it is worth
reporting this obstruction via black rhombs on the main-table
(Fig.\,\ref{Degree8-(M-i)-curve-TABLE:fig}). The resulting pattern
is especially delightful of regularity and once more Viro's genius
is baffling our spirit and requires highest admiration. What a
pity just that nobody ever published this table (except perhaps in
Viniti's preprint-series of Gorki not accessible in the west). It
should be noted (from the geography of the
main-table=Fig.\,\ref{Degree8-(M-i)-curve-TABLE:fig}) that Viro's
hypothesis of divisibility by 4 of the content seems
%%%to be
fulfilled precisely in the depressions (of Gudkov's sawtooth)
corresponding to RKM-schemes. As the latter are of type~I, it is
very likely that this Viro obstruction involves again (like the
imparity law) a matter of complex orientations.

This is very interesting but alas does not answer our
%question
query  on the scheme $S_0=15\frac{4}{1}$. So lacking a better idea
let us turn to our first method of construction by using
non-maximal dissipation of $2X_{21}$ (i.e. the coaxial quadruplet
of ellipses). Alas on consulting again Fig.\,55 on p.\,1118 of
Viro 89/90 \cite{Viro_1989/90-Construction}, we realize that only
maximal dissipation are listed there. So what about the other
permissible dissipation? [{\it Update} [01.07.13].---For some
improvised tabulation of those non-maximal dissipations, see our
Fig.\,\ref{ViroDEGREE8_2:fig}, where we nearly got the curve in
question $S_0=15\frac{4}{1}$ through Viro's quadri-ellipse
dissipation.]

\subsection{Shustin's fifty $(M-2)$-obstructions}

[12.05.13] Another valuable piece of information is the table of
Shustin 1990/91 \cite{Shustin_1990/91-New-restrictions} of
$(M-1)$-schemes with 3 nonempty ovals. This reports new
obstruction due to Shustin and constructions made by
Goryacheva-Polotovskii 1985 \cite{Goryacheva-Polotovskii_1985}
(abridged GP on the table or below). Alas it seems that the 6th
and 7th line of the first row contains a misprint. Perhaps one
should read $(1,6,11)$ instead of $(1,7,11)$. So it seems that
there was just a typographical permutation there.

At first sight Shustin's table looks a bit chaotic since an
$(M-1)$-scheme below the Fiedler-Viro $M$-prohibition is
generically prohibited, but there are exceptions to the rule. For
instance the scheme $\frac{1}{1}\frac{6}{1}\frac{11}{1}$ is
constructed by GP, and this may be explained as a degeneration
(better contraction) of Shustin's $M$-scheme
$\frac{1}{1}\frac{7}{1}\frac{11}{1}$. Same remark for
$\frac{4}{1}\frac{7}{1}\frac{7}{1}$ which can be regarded as a
contraction of Shustin's $M$-scheme
$\frac{5}{1}\frac{7}{1}\frac{7}{1}$. Likewise the $(M-1)$-scheme
$\frac{5}{1}\frac{6}{1}\frac{7}{1}$ can  be viewed as contraction
of the same Shustin $M$-scheme. Next
$4\frac{1}{1}\frac{2}{1}\frac{11}{1}$ also looks irregular but
occur as  degeneration of Shustin's scheme
$4\frac{1}{1}\frac{3}{1}\frac{11}{1}$. Hence Shustin's list of
scheme is therefore quite concomitant with the Itenberg-Viro
contraction principle
%%%(\ref{Itenberg-Viro-contraction:conj}),
(see v.2), and the latter could be used to explain regularity of
the sequence (through understanding the architecture of the
pyramid). Precisely, whenever an $(M-1)$-scheme appears as
degeneration of a constructed $M$-scheme then it appears in the
GP-list of Goryacheva-Polotovskii. So for instance it is quite
nice to visualize the 4 possible degenerations of (Viro-Shustin's)
scheme $4\frac{3}{1}\frac{5}{1}\frac{7}{1}$ (diminishing \`a tour
de r\^ole any of the structural constants, i.e  either
$3\frac{3}{1}\frac{5}{1}\frac{7}{1}$,
$4\frac{2}{1}\frac{5}{1}\frac{7}{1}$,
$4\frac{3}{1}\frac{4}{1}\frac{7}{1}$,
$4\frac{3}{1}\frac{5}{1}\frac{6}{1}$). In comparison  Viro's
scheme $8\frac{3}{1}\frac{3}{1}\frac{5}{1}$ admits only 3 possible
contractions. All this and more is best explained by the geometric
view of this 3D-pyramid, when all
%%%levels
6 layers are imagined superimposed.

Further the 1st level of the 2nd pyramid can degenerate upon the
1st pyramid. For instance Viro's $M$-scheme
$8\frac{1}{1}\frac{3}{1}\frac{7}{1}$ can degenerate by contraction
of an empty oval to $8\frac{0}{1}\frac{3}{1}\frac{7}{1}=
9\frac{3}{1}\frac{7}{1}$, which belongs to the 1st pyramid. Albeit
we do not took the pain as yet to construct this very specific
scheme (yet cf. Fig.\,\ref{ViroDEGREE8_2:fig} for a heuristic
construct), its existence looks nearly evident for it may appear
as contraction of many (at least two) other $M$-schemes, namely
$9\frac{3}{1}\frac{8}{1}$ and $9\frac{4}{1}\frac{7}{1}$ (both
obtained via Viro's simplest method with a quadruplet of
ellipses).

In conclusion, the principle of contraction makes all
constructions mentioned in Shustin's table (and implemented by
Polotovskii et al.) look nearly evident. Shustin (or his
typographer?) seems only to miss the $(M-1)$-scheme (immediately)
below Wiman's, i.e. $15 \frac{1}{1}\frac{1}{1}\frac{1}{1}$.
Further (as remarked by Shustin, \loccit) there are 2 interesting
exceptions where the domination principle
%%%%of
via contraction of an empty oval of an $M$-curve does not tell
anything. First there is the $(M-1)$-symbol
$S_{top}:=\frac{6}{1}\frac{6}{1}\frac{6}{1}$ near the summit of
the pyramid materialized by $\frac{6}{1}\frac{6}{1}\frac{7}{1}$
(yet ruled out by Viro's imparity law proved via complex
orientations). So the status of
$\frac{6}{1}\frac{6}{1}\frac{6}{1}$ is extremely puzzling, and
perhaps still open today (a priori unaffected by  the recent
efforts of Chevallier, Orevkov ca. 2002 remodelling other portions
of the pyramid). So if this scheme exists then it is quite likely
to be maximal, yet not of type~I. (For instance the enlargement
obtained by adding one outer oval is not permissible via Gudkov's
hypothesis.) This would be another counterexample to the (reverse)
sense of Rohlin's maximality conjecture. On the other hand, if the
scheme $\frac{6}{1}\frac{6}{1}\frac{6}{1}$ existed, we would by
the contraction principle get the scheme
$\frac{5}{1}\frac{6}{1}\frac{6}{1}$ which looks constructible
being dominated by Shustin's $M$-scheme
$\frac{5}{1}\frac{7}{1}\frac{7}{1}$ (after 2 contractions).

Further another case left open by Shustin is the $(M-1)$-scheme
$12\frac{2}{1}\frac{2}{1}\frac{2}{1}$. It is somewhat more
surprising that this scheme does not fall under Shustin's
prohibition (at least from an architectural viewpoint of a pyramid
builder). Of course if the scheme existed it could degenerate to
$12\frac{1}{1}\frac{2}{1}\frac{2}{1}$, which is  likely to exist
being dominated by $12\frac{1}{1}\frac{2}{1}\frac{3}{1}$ itself
dominated by  Viro's $M$-scheme
$12\frac{1}{1}\frac{3}{1}\frac{3}{1}$ (which alas we were
presently/personally not able to construct despite Viro's
revendication of this territory).

Of course the philosophy of all this is the superiority of
geometry upon arithmetics and symbolism (e.g. Shustin's linear
table). After all, the abstract concept of integer just arose by
visualizing sheep[s]
%shephard
ca. 50'000 BC (evidence of counting according to Boyer-Merzbach
\cite{Boyer-Merzbach_1968/89}).

All this does not, alas,  answers our basic query about
$S_0:=15\frac{4}{1}$ (UPDATE: this is almost surely realized, cf.
either Fig.\,\ref{ViroDEGREE8_SHUSTIN2:fig} for a construction via
Shustin's curve or
Fig.\,\ref{ViroDEGREE8_TRICKY_SUITE_XXX_15_4-1:fig}k for one via
Viro's 1st curve, or even Fig.\,\ref{ViroDEGREE8_2:fig}f for a
derivation from the quadri-ellipse yet with self-guessed accessory
parameters $(\al, \be,\ga)$.) On the table we reported its
$(M-1)$-enlargements by black circles. Those schemes are those of
the 1st pyramid with the symbol $\frac{4}{1}$ occurring as
substring. Further if the additional oval is traced inside the
nonempty one then we get two schemes near the bottom of the 3rd
pyramid, namely $15(1,1\frac{3}{1})$ and $15(1,2\frac{2}{1})$. But
those being dominated by a Korchagin scheme (namely
$15(1,2\frac{3}{1})$) we get again a corruption of Rohlin's
maximality principle (at least when combined with the contraction
principle of Itenberg-Viro, but that is not needed actually).

Likewise we can trace all enlargements of the other RKM-schemes
with one nonempty oval by using other symbols. For instance for
$11\frac{8}{1}$ we get a trajectory of 5-branched stars on the
main-table involving schemes running below Viro's $M$-schemes. So
again a corruption of RMC is derived. The same remark applies to
the schemes $7\frac{12}{1}$ and $3\frac{16}{1}$.
So we arrive at the:

\begin{lemma} Either the four RKM-schemes with one nonempty oval
are prohibited or Rohlin's maximality conjecture (RMC) is false.
(Of course
%a
another dramatic issue would be that RMC is true but
Viro's method is false!)
\end{lemma}

Then we can apply the same method to the RKM-scheme $14
\frac{1}{1}\frac{3}{1}$ and get a series of enlargements in the
2nd pyramid running below $M$-schemes constructed by either Viro
(not understood by us as yet) or Shustin (clearly understood and
depicted earlier in our text, cf.
Fig.\,\ref{ViroDEGREE8_SHUSTIN:fig}).

 So again this RKM-scheme
looks prohibited. One can continue the same game with the
RKM-scheme $10\frac{1}{1}\frac{7}{1}$ and again we find in the 2nd
pyramid symbols running below schemes constructed by Shustin. Idem
for $6\frac{1}{1}\frac{11}{1}$.

\begin{figure}[h]\Figskip
%\vskip-1.2cm\penalty0
\centering
%
%\hskip-0.7cm\penalty0
\epsfig{figure=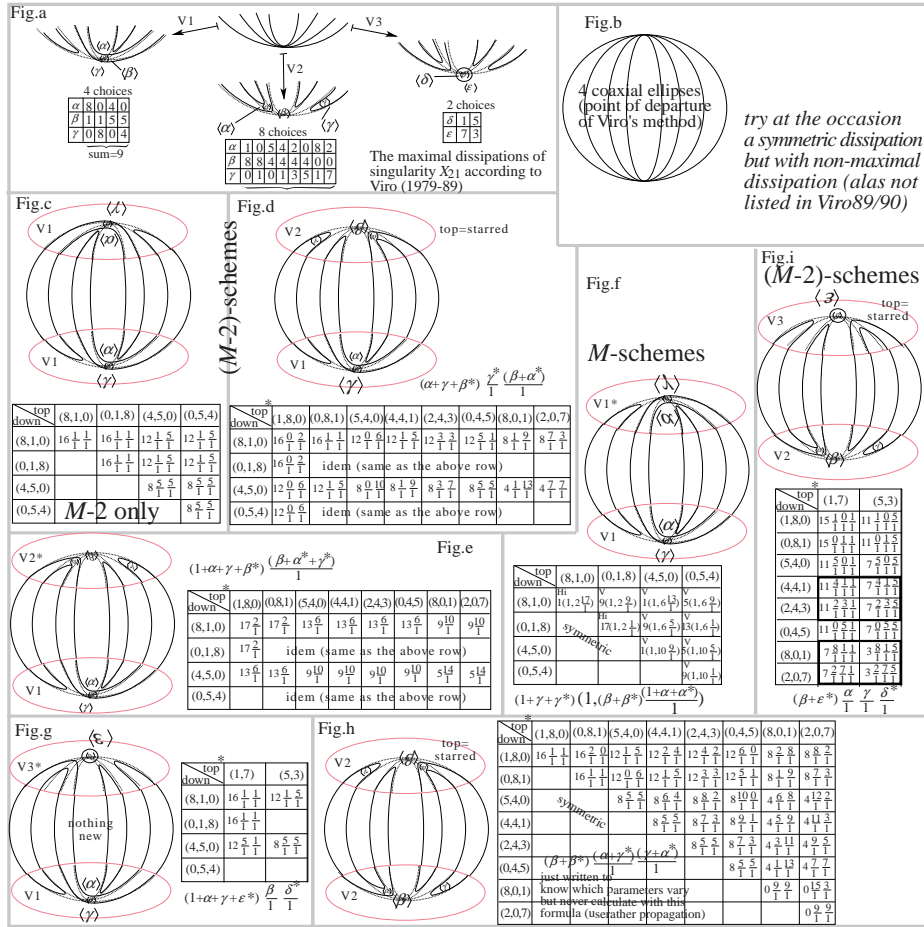,width=122mm} \captionskipAG
  \caption{\label{ViroDEGREE8-(M-2):fig}%
  Viro's method for $(M-2)$-schemes}
\figskip
\end{figure}

Then we arrive at the RKM-scheme $2\frac{1}{1}\frac{15}{1}$. It
admits one $(M-1)$-enlargement in the 2nd pyramid (which is
prohibited by Shustin's obstruction), and by B\'ezout it cannot be
enlarged in the 3rd pyramid. So it seems that the scheme in
question should exist. One can try to realize it via Viro's method
but we failed after a quick try. To proceed more systematically
let us tabulate what is gained by Viro's method for
$(M-2)$-schemes. Compare Fig.\,\ref{ViroDEGREE8-(M-2):fig}. The
schemes so obtained always satisfy the Gudkov congruence mod 8
(i.e. lies 2 stages below certain $M$-curves). This is surprising
but probably related to the fact that we use only the maximal
dissipation. In fact, it seems that one can confine to the
greatest table which is by far the most prolific in creating
schemes. As for $M$-curves the sequence into which Viro's method
creates schemes looks a bit random (for a modest intelligence like
the writer but there is surely some hidden  rule regulating this).

Of course it is fairly evident how to vary the construction to
gain more schemes by using other curves
%%%of
or singularity. The difficulty is just to proceed as
systematically as possible despite the varied choice. Some feeling
for arts it a bit required, or at least patience. As yet the
method is a bit disappointing as it did not generated any
RKM-scheme (lying in the depression of sawtooth). Those are
however the most interesting guys for testing Rohlin's maximality
conjecture. In particular we expect that all RKM-schemes marked by
large symbols (circles, stars with 5,4,6 branches either black or
white colored) are not realized algebraically, safe for the scheme
starred in white by 6 branches (i.e. $2\frac{1}{1}\frac{15}{1}$).
However as yet we could not realize it.

How to realize this scheme via Viro's method? Are the other
schemes known to be prohibited?

%%%% SAXOOO AFTER CHRISTA'S UNFALL AND BEFORE ENSMATH (13.05.13)

[11.05.13] It seems also that upon studying carefully the
combinatorics of this table Ragsdale conjecture is trivially true
in degree 8, while using of course some theoretical gadgets like
Petrovskii/Arnols or maybe Rohlin's formula. This exercise is
worth completing at the occasion (and the truth of this assertion
is mentioned in Rohlin 1978).

\subsection{Hilbert's 16th for all septics (Viro 1979)}

[28.04.13] The degree $m=7$ is the largest presently known where
Hilbert's 16th is completely settled by the effort of a single
hero, O. Ya. Viro 1979. The conceptions of V.\,A. Rohlin who alas
suffered from a first hearth attack ca. 1975 (compare Vershik ca.
1986 for more details about Rohlin's carrier, life and vivisstudes
in the GOULAG) probably contributed to a big extend in Viro's
breakthrough. Let us now state the precise result:

\begin{theorem} \label{Viro-deg-7-isotopy-class-121-schemes:thm} {\rm (Viro 1979, Viro 1980
\cite{Viro_1980-degree-7-8-and-Ragsdale}, 1986
\cite{Viro_1986/86-Progress})}.---Any nonsingular real septics
realize exactly one of the following $121$ schemes (when the
pseudoline is omitted from the symbolism):

$\bullet$ $\al \frac{\be}{1}$ with $\al+\be\le 14$, $0\le \al \le
13$, $1\le \be \le 13$,

$\bullet$ $\al$ with $0\le \al \le 15$,

$\bullet$ $(1,1,1)$ (deep nest).
\end{theorem}

\begin{figure}[h]\Figskip
%\vskip-1.2cm\penalty0
%\centering
%
\hskip-0.7cm\penalty0 \epsfig{figure=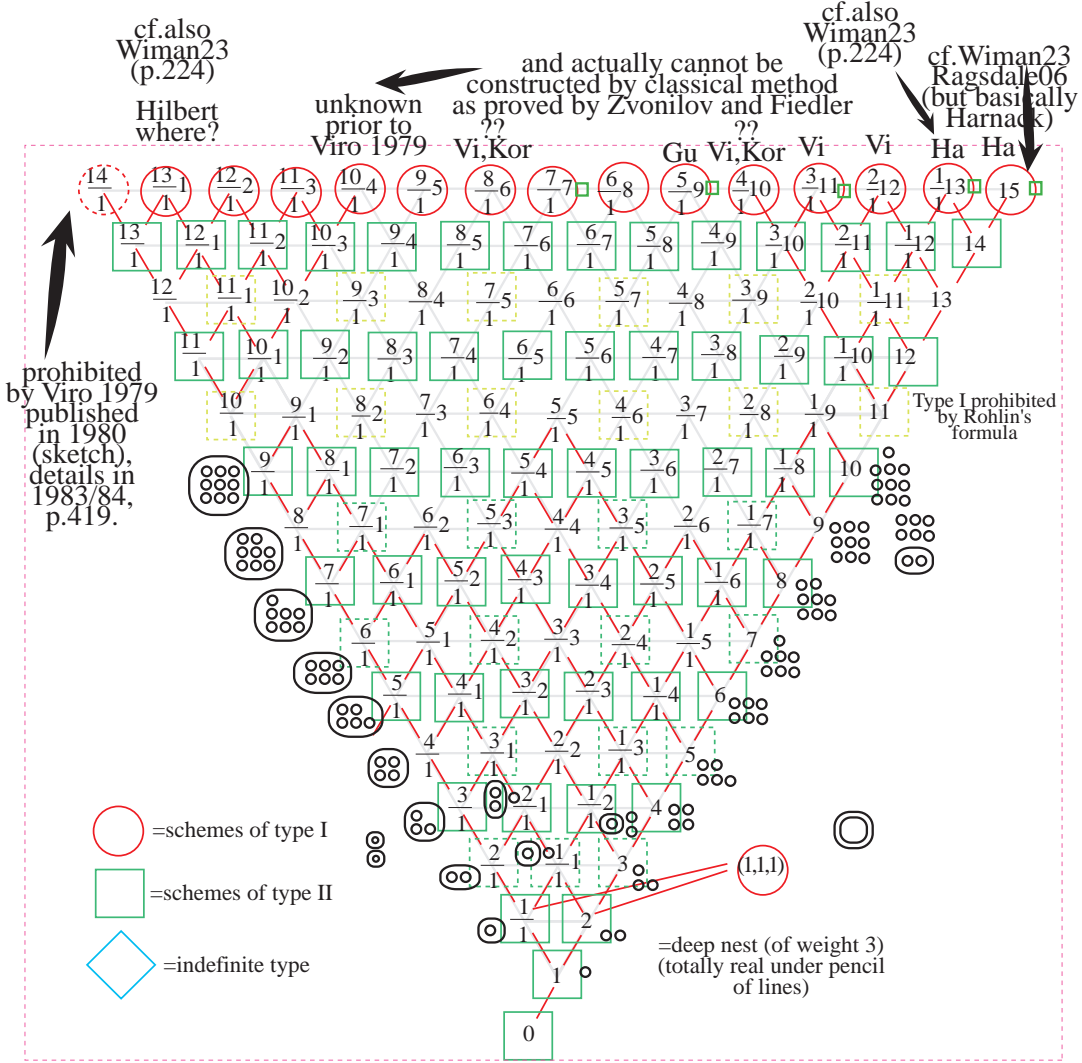,width=142mm}
\captionskipAG
  \caption{\label{Degree7:fig}%
  The Gudkov table in degree 7 due to Viro 1979
  (details published slightly later, notably in Viro 1986
  \cite{Viro_1986/86-Progress})}
\figskip
\end{figure}

It is convenient to visualize this result via a Gudkov pyramid
(Fig.\,\ref{Degree7:fig}). On it we indeed count
$1+2+3\dots+15=\frac{16 \cdot 15}{2}=8\cdot 15=120$ schemes on the
main-triangle and one scheme $0$ must be added getting the total
of 121 isotopy classes listed by Viro. So in degree 7 the Gudkov
pyramid is still nearly planar (upon omitting the deep nest).
Actually adding the deep nest we find a total of 122 schemes.
(Maybe there is a little mistake in Viro at this place, which we
actually cite via Brugall\'e 2005
\cite{Brugallé_2005/07-symm-plane-curves}, so that maybe the
mistake is perhaps due to Brugall\'e). No in fact the explanation
is that the scheme $\frac{14}{1}$ is not realized, since we cannot
take $\be=14$ in (\ref{Viro-deg-7-isotopy-class-121-schemes:thm}).

Of course it would be also useful to enhance this Viro
classification by the types I/II/indefinite in the sense of Klein
1876 (and Rohlin 1978). The first basic point is that with
periodicity 2 all horizontal rows are of type II by Klein's
congruence $r\equiv_2 g+1$. Also by total reality the deep nest
$(1,1,1)$ is of type~I and so are all $M$-schemes by Klein 1876.
Probably this affords the complete list of all type~I schemes.
Also by the Rohlin-Mishachev a dividing curve has at least $r\le
m/2=3.5$ components so we may infer that the scheme $1$ is of
type~II (although this information is not covered by Klein's
congruence). Further according to our ``toutou'' conjecture
%%%\ref{toutou:conj}
(in v.2 of Ahlfors) a scheme of type~II remains of type~II after
addition of a pseudoline (while augmenting its degree by one
unit). First when comparing with the Gudkov table in degree 6
%%%(Fig.\,\ref{Gudkov-Table3:fig})
we see that this principle recover
a good portion of the type~II inferences drawn automatically from
Klein's congruence. But doing transfers at other levels (non
predestined by Klein's congruence) yields a new collection of
schemes of type~I marked by dashed-green squares on
Fig.\,\ref{Degree7:fig}. If this is true we suspect that this and
more information (propagate the 2-by-2 lattice of dashed squares
upwards, cf. yellow-green dashed squares) may be derived from
Mishachev's formula (1975/76 \cite{Mishachev_1975/75}), i.e. the
odd-degree avatar of Rohlin's formula. Once this extension is
effected it is perhaps the case that all remaining schemes are of
indefinite type. All this remains of course to be verified more
seriously but is surely already well-known to Viro, Fiedler, etc
(and perhaps even Rohlin).

 Let us now do some constructions. Let us start with the
sextic $C_6$ of scheme $\frac{2}{1}2$ and of type~II depicted
below while adding a line and smoothing consistently \`a la
Fiedler so as to ensure type~I. The resulting $C_7$ is of type~I
and realize the scheme $\frac{1}{1}3$. This spoils our toutou
conjecture. For cross-reference let us state:

\begin{prop}\label{toutou-conj=FALSE:prop}
The toutou (= the doggy way) conjecture is false, i.e. a scheme of
type II needs not staying of type II after aggregating a
pseudoline.
\end{prop}

\begin{figure}[h]\Figskip
%\vskip-1.2cm\penalty0
%\centering
%
\hskip-0cm\penalty0 \epsfig{figure=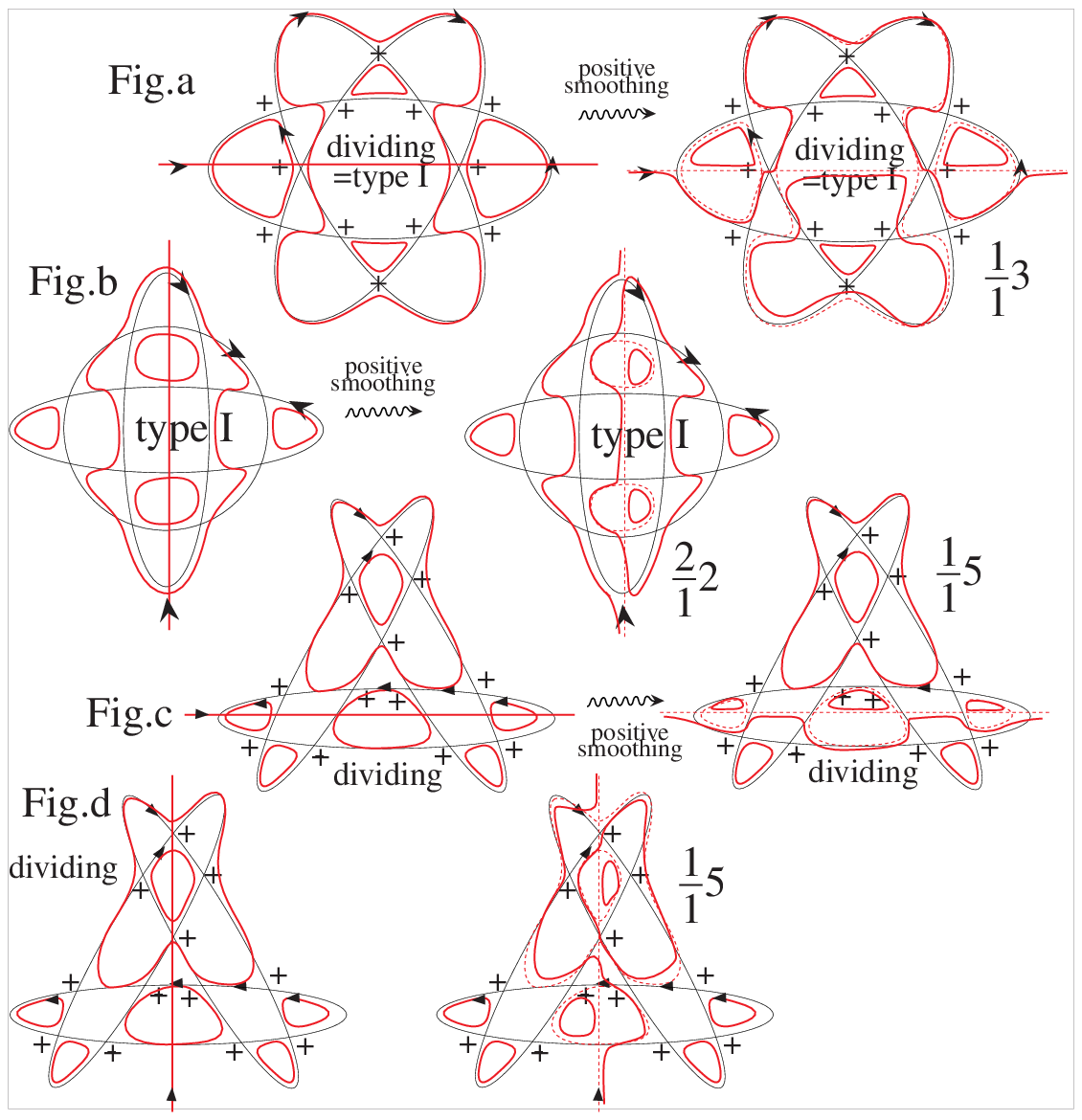,width=122mm}
\captionskipAG
  \caption{\label{V2-12:fig}%
  Corrupting the toutou conjecture}
\figskip
\end{figure}

\begin{proof} In fact our conjecture was motivated by the case of
quintic where the conjecture toutou is true. However in degree 7
it is false as the sextic scheme $\frac{1}{1}3$ is of type~II (cf.
Rohlin's table (in v.2 of Ahlfors)
%%%%Fig.\,\ref{Gudkov-Table3:fig}
as follows from
either Arnold's congruence or Rohlin's formula), yet the same
scheme augmented by a pseudoline is not of type~II since it
contains the dividing representative depicted on
Fig.\,\ref{V2-12:fig}a.
\end{proof}

So in fact we would like to know exactly the distribution into the
3 possible types as Rohlin 1978 was able to do for sextics. The
correct (and complete?) answer seems to be given by Brugall\'e
2005 \cite[p.\,4]{Brugallé_2005/07-symm-plane-curves} upon
assembling the following 3 prohibitive results:

$\bullet$ elementary B\'ezout prohibitions;

$\bullet$ Fiedler's orientation alternating rule (cf. Viro 1986/86
\cite{Viro_1986/86-Progress});

$\bullet$ the Rohlin-Mishachev formula;

while combining at the constructive level with the articles by
Soum 2001 \cite{Soum_2001} and Le Touz\'e 1997
\cite{Fiedler-Le-Touzé_1997-DEA}.

So taking for granted the two lemmas (1.2 and 1.3) stated by
Brugall\'e 2005 \cite{Brugallé_2005/07-symm-plane-curves}, we
arrive at the following corrected map. Hopefully this information
is correct as we had not the patience to check all the details
carefully. Crudely put we see that everything is governed by
Klein's congruence safe the lower rows of blue rhombs which (at
lower altitudes $r=10,8,6,4$) have lacunae probably explained by
Mishachev's formula, while tending to contains an increasing
number of schemes of type~II as $r$ decreases (thermodynamical
interpretation?).

\begin{figure}[h]\Figskip
%\vskip-1.2cm\penalty0
%\centering
%
\hskip-0.7cm\penalty0
\epsfig{figure=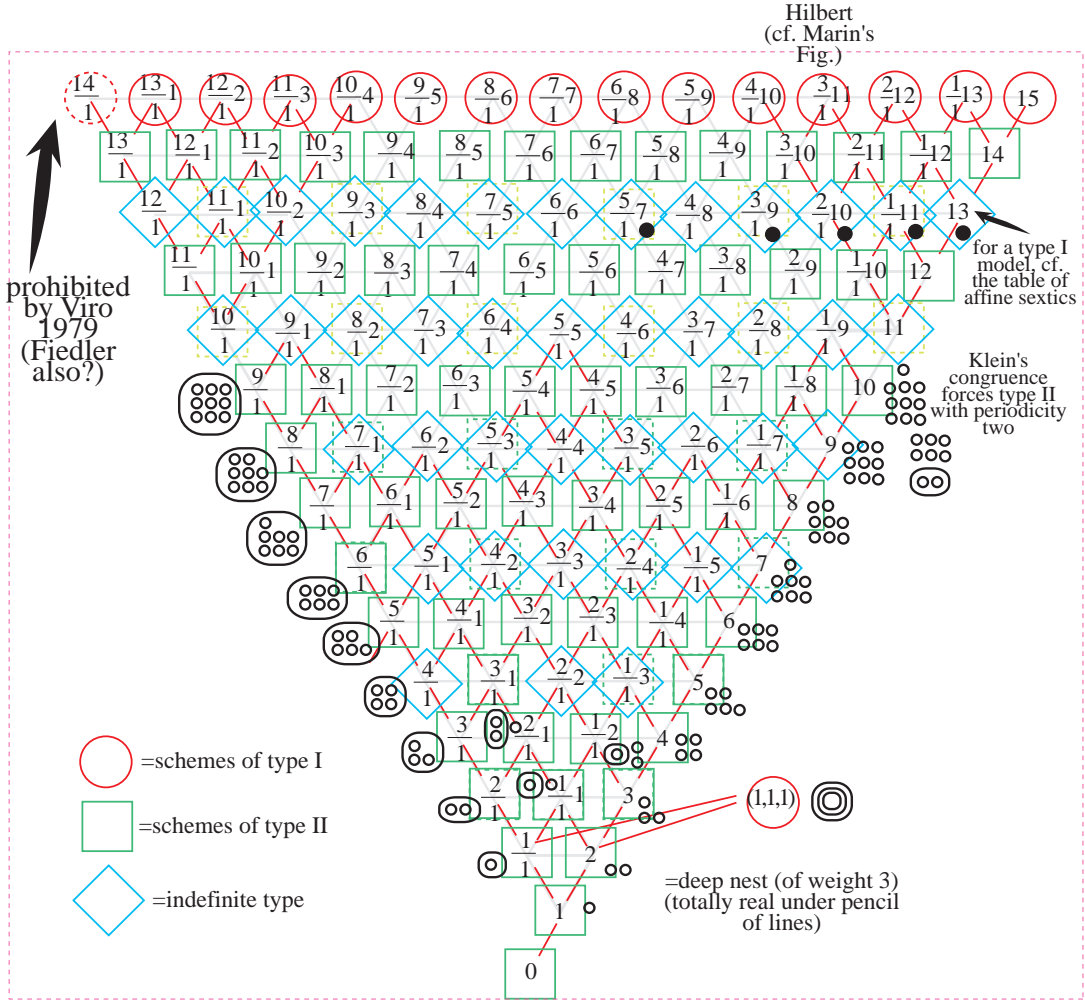,width=142mm} \captionskipAG
  \caption{\label{Degree7Brugalle:fig}%
  The Gudkov table in degree 7 due to Viro 1979
  but decorated \`a la Klein-Rohlin with types (courtesy of
  Le Touz\'e 1997, Soum 2001, Brugall\'e 2005).}
\figskip
\end{figure}

[03.05.13] An incomplete version of this Brugall\'e table was also
published as Fig.\,7 of De Loera-Wicklin 1998
\cite[p.\,204]{De-Loera-Wicklin_1998}. But this contains very
sparse information on degree 7, in sharp contrast with the
complete solution of Brugall\'e.

Of course by virtue of Marin's example
%%%(Fig.\,\ref{Marin:fig})
(cf. v2 of Ahlfors) we cannot like in degree 6 expect the schemes
plus the type data to form a unique rigid isotopy class. Yet
Brugall\'e's result combined with Marin's ``duplication'' of
chambers should produce a lower bound on the number of components
past the discriminant of septics. Probably non body in the world
has the slightest idea of the value of this number $\delta_7$. By
using Brugall\'e only we get a lower bound of ca. the 121 schemes
of Viro, plus some duplication on $13+11+8+6+3=39$ so that we get
160 components past the discriminant. By using Marin's argument
this can be raised to 161. This is certainly  far from sharp and
we speculate that $\delta_7$ could be as high as 500? Probably no
technology like K3 surfaces and their periods tackling the case
$m=6$ is presently available. But can we get an upper bound on
$\delta_7$?

As another little remark it may be observed that Rohlin's
maximality conjecture holds true in degree 7 (as a consequence of
Brugall\'e's table), and even in the strong sense that a scheme is
of type I iff it is maximal.

%%%%%%%%%% TO NEUTRALIZE all 4 MARION's WHEN FINISHED

\subsection{Acknowledgements}

The author wishes to thank the following long list of
geometers (in chronological order of their interaction with
the writer in connection to the present text)

$\bullet$ Felice Ronga (ca. 1997/98 for his explanation of
Brusotti's theorem),

$\bullet$ Claude Weber, Michel Kervaire (for their explanations on
how to classify Klein's symmetric surfaces via the quotient
bordered surface),

$\bullet$ Daniel Coray for the geometric Galois action and all his
logistical support throughout  the years (up to present days),

$\bullet$ Fr\'ed\'eric Bihan for pleasant freshman calculus
discussions about real algebraic geometry,

$\bullet$ Lee Rudolph (ca. 1999 for explaining to us [Claude and
me] what is the natural topological model for a real elliptic
curve with only one ``oval'', namely just a torus $S^1\times S^1$
acted upon by factor permutation $(x,y)\mapsto (y,x)$ fixing
thereby the diagonal circle),

$\bullet$ Alexis Marin, Viatcheslav Kharlamov, Oleg Viro,
Jean-Jacques Risler, Thierry Vust, Michel Kervaire, Pierre de la
Harpe, John Steinig (for their comments and corrections improving
the shape of the article Gabard 2000 \cite{Gabard_2000})

$\bullet$ Ragahavan Narasimhan, Jacek Bochnack (ca. 1999 for {\it
not\/} having been in touch with Ahlfors' result of 1950
\cite{Ahlfors_1950} enabling me some free gestation about thinking
on the problem)

$\bullet$  Manfred Knebusch for his kind interest in the
modest work Gabard 2000 \cite{Gabard_2000},

$\bullet$ Johannes Huisman  for his constant interest
(2001--04--06), and his care about correcting  bugs in both my
Thesis and the article Gabard 2006 \cite{Gabard_2006},

$\bullet$ Sergei Finashin for an exciting discussion in Rennes
2001,

$\bullet$ Jean-Claude Hausmann (ca. 2000/01) for telling me about
the standard surjectivity criterion via the Brouwer degree, which
was decisive to complete Gabard 2006 \cite{Gabard_2006},

$\bullet$ Antonio Costa,  for  his fascinating talks in Geneva,

$\bullet$ Bujalance for his surely over-enthusiastic
Zentralblatt review of my article (Gabard
2006~\cite{Gabard_2006}),

$\bullet$ Fraser-Schoen, whose brilliant work revived my
interest in the theory of the Ahlfors' mapping (ca. the 13
March 2011) at a stage were I was mostly
%involved within the
%theory of so-called
sidetracked by ``non-metric manifolds''
%thanks to efforts of
due to the infectious fascination of Mathieu Baillif and David
Gauld.

%%%added 05.09.13

$\bullet$ Mathieu Baillif (the world leader de la Wissenschaft de
gauche, especially, when it comes to
%%%the legend of
the
Dekrummante) is much acknowledged for all his
%songful
wisdoms,
%and
partnership advices, and
%his
overall pleasant
%vision
approach
%of
to life, science and arts. (Especially noteworthy is the obvious
analogy between some of his minimalist picturing and those by
Polotovskii in the qualitative theory of decomposing curve, se
also the Petit Prince catalogue by Orevkov.)

$\bullet$ Elias Boul\'e, Tarik Garidi, Marie Pittet for frequent
electronic support, and their patience about earing scenes of the
adventure,

$\bullet$
%Grisha Mikhalkin,
Stepan Orevkov, Oleg Viro (2011) for
%their patience about hearing some of my excitation)
their talks and pleasant discussions,

$\bullet$ Marc Coppens (2011--12) for e-mails, and  his
%brilliant
work on the separating gonality (2011 \cite{Coppens_2011})
adumbrating sharper insights on the degree of the Ahlfors function
(or rather the more general allied circle maps). His turning-point
result
%constitutes a turning point in case we
%manage in the future to
appeals to a better conciliation of the analytic theory of Ahlfors
with the algebro-geometric viewpoint.

$\bullet$ (2011/12) Hugo Parlier, Peter Buser, Alexandre Girouard,
Gerhard Wanner und Martin Gander are acknowledged for their recent
e-mail exchanges.

$\bullet$ (Oct. 2012) Daniel Coray for enlarging the capacity of
my TeX-compilator, and for his lovely (Cambridge-style) teaching
about the geometric Galois action in ca. 1999.

$\bullet$ (16 Nov.\;2012) Franc Forstneri\v{c}
%for his
%encouragement and
for pointing out his text with Wold (2012
\cite{Forstneric-Wold_2012}) showing the  state of the art on the
proper holomorphic embedding problem.

$\bullet$ (09 Jan. 2013) Oleg Viro for his excellent answers on
some naive questions on Rohlin's paper
%(cf.
%Sec.\,\ref{e-mail-Viro:sec}).
(cf. v2 of Ahlfors);

$\bullet$ (10 Jan. 2013) Alexis Marin for his invaluable insights
on Klein's intuitions and
%much
more
%%(cf. Sec.\,\ref{e-mail-Viro:sec}).
(cf. v2 of Ahlfors);

$\bullet$ (January-February-March 2013) Viatcheslav Kharlamov,
Stepa Orevkov, Eugenii Shustin, S\'everine Le~Touz\'e, Thomas
Fiedler for all their letters reproduced in
%%Sec.\,\ref{e-mail-Viro:sec}.
(v2 of Ahlfors);

$\bullet$ (ca. 14/16 April 2013) Grisha Mikhalkin for pointing out
that my idea of stability with satellites bears some resemblance
with Wiman 1923 \cite{Wiman_1923} an intuition which turned to to
be perfectly correct when I took the pain to read Wiman's text in
more details about one month later (the 09.05.13). Hence Wiman's
role as a sort of forerunner of Rohlin might be more emphasized
(of course both admit some inspiration from Klein). Wiman's idea
of satellites is very acute relegate all of our priority in
conjecturing a stability under satellites. To pinpoint a bit,
Wiman only stresses  maximality for satellites of his series of
$M$-schemes (cf. Fig.\,\ref{Wiman:fig}) whereas we conjecture it
for any scheme of type~I (so for instance for the double of
Rohlin's schemes in degree 6 satisfying the RKM-congruence).

\section{Bibliographic comments}
%(no exhaustiveness claim is
%made, and please contact the writer if some omission is
%committed.
The writer does not pretend that the following bibliography is
complete (nor that he absorbed all those fantastic contributions
in full details). More extensive bibliographies (overlapping
ours), but covering more material include those of:

$\bullet$ Ahlfors-Sario 1960 \cite{Ahlfors-Sario_1960} (ca. 40
pages times 25 items per pages=1000 entries covering such topics
as the Dirichlet problem, extremal problems, the type problem, the
allied classification theory, etc.);

$\bullet$ Grunsky 1978 \cite{Grunsky_1978} (=562 refs,
including 48 Books).
%%%Joke borrowed in Morrey 1948

Most entries of our bibliography are followed by some comment
explaining briefly the connection to our primary topic of the
Ahlfors map. The following symbolism is used:

$\clubsuit$ serves to point out a special connection to Ahlfors
1950 (especially alternative proofs).

$\spadesuit$ gives other comments (attempting to summarize the
paper contents or to explicit the connection in which we cite
it).

\iffalse

$\P$ signals papers not quoted in the main-body of the text
but
%which looks more-or-less
connected to our topic.
%Further
%we add after this symbol some question perhaps worth pursuing.

\fi

$\bigstar$ marks sources, I could not as yet procure a copy.

\def\Ah47{\textsf{A47}}

\def\A50{\textsf{A50}}

\def\AS60{\textsf{AS60}}
\def\G78{\textsf{G78}}

$\bullet$ the stickers/sigles \AS60, \G78 are assigned when the
source has already been cited in Ahlfors-Sario 1960
\cite{Ahlfors-Sario_1960} resp. Grunsky 1978 \cite{Grunsky_1978}.

$\bullet$ \A50 designates those references citing the paper
Ahlfors 1950 \cite{Ahlfors_1950} (there represents circa 106
articles on ``Google''), and occasionally \Ah47 those quoting
Ahlfors 1947 \cite{Ahlfors_1947}.

{\bf $\heartsuit$n} is something like the indicator of the US
rating agency (to be read ``liked by $n$''). It indicates the
cardinal number {\bf n} of citations of the paper as measured by
``Google Scholar''. The latter machine often misses
cross-citations, especially those in old books, or old articles
with references given in footnotes format. Many sources cited in
Grunsky's book (1978 \cite{Grunsky_1978}) are never cited
electronically. Accordingly, those rating numbers only supply a
statistical idea of the literature ramifications lying beyond a
given entry. Also
 low-citation articles are sometimes the most
 polished product ripe for museum entrance. Forelli 1979
\cite{Forelli_1979}  is typical: self-contained, elegant and
polished proof of Ahlfors result, yet  only rated by 3.
%In
%short, those numbers bear only a low significance and of course
%deviates strongly from the real importance in relation to our main
%interest.

Our bibliography is somewhat conservative with comparatively few
modern references. Our excuse is two-fold: modern expressionism is
sometimes harder to grasp, and recent references are usually well
detected through computer search.
%However, links between the older literature
%are harder to detect by pure computer assisted search.

(Papers are listed in alphabetical, and then chronological order,
regardless of shared co-authorship.)

The primary focus is on  the Ahlfors map and the weaker (but more
general) circle maps. As a such the topic overflow slightly over
the territory of real algebraic geometry. Ahlfors-Sario's book
\AS60 address Riemann surfaces, whereas Grunsky's book \G78
focuses to the case of planar domains. Hence both bibliographies
\AS60, \G78 are quite complementary, and ours is essentially a
fusion of both, but we gradually included more and more recent
contribution. Still additional references are welcome.
%especially
%if there is a little summary pointing out the connection with the
%Ahlfors map.

For conformal maps, it is helpful when browsing the vast
literature to keep in mind the basic question: {\it what result
through which method?}

{\bf Results.} Objects traditionally range along increasing order
of generality through: simply-connected regions,
multiply-connected ones and finally Riemann surfaces. We often add
a humble compactness proviso, as the passage to open objects is
traditionally achieved through the exhaustion trick (going back at
least to Poincar\'e 1883 \cite{Poincare_1883}, and see also Koebe
1907 \cite{Koebe_1907_UbaK1}), and active in recent time (e.g.
Garabedian-Schiffer 1950 \cite{Garabedian-Schiffer_1950}.)

As to the mappings, they may all be interpreted in some way or
another as ramification of RMT (Riemann's mapping to a
circle=disc). We distinguish primarily:

$\bullet$ CM=circle maps (usually not univalent, but multi-sheeted
disc with branch, or winding points=Windungspunkte)

$\bullet$ KNP=Kreisnormierung(sprinzip) (univalent map to a
circular domain)

$\bullet$ SM=slit mappings for various types of them (parallel,
circular, radial, logarithmic spiral, etc.). Those are all allied
to certain natural foliation of the sphere, and some extreme
generality in this respect is achieved in Schramm's Thesis where
any foliation is permitted as support for the slits.

{\bf Methods.} They may be classified
%ranged
%ranked
in two broad classes quantitative vs. qualitative (each having
some branchings):

\noindent$\bigstar$ (Quantitative) variational methods,
including:

$\bullet$ DP=Dirichlet principle (or more broadly speaking,
potential theory=PT, centering around such concepts as the
Green's function, harmonic measures (i.e. harmonic function
with special null/one boundary prescription of the various
contour), etc. Of course, there is a standard yoga between
Dirichlet and Green, so all this is essentially one and the
same method.

$\bullet$ IM=Iterative methods (originators: Koebe and
Carath\'eodory), and by extension this may proliferate up to
including the circle packings.

$\bullet$ EP=extremal problems (e.g. the one of maximizing the
derivative amongst the class of function bounded-by-one) and
leads to the Ahlfors map.

$\bullet$ BK=Bergman kernel (or Szeg\"o kernel), here the
fundamental ideas rest upon Hilbert's space methods, and the
idea of
%orthonormal
orthogonal system. Initially, the method is also inspired by Ritz,
and Bieberbach extremal problem (1914
\cite{Bieberbach_1914-Theo-Praxis-d-konf-Abb}) for the area swept
out by the function. Since the middle 1940's, there were found
several conformal identities among so-called domain functions
(Green's, Neumann's, etc.) and the kernel functions so that
virtually this is now highly connected to DP$\approx$PT. Also the
Ahlfors map is  expressible in term of the Bergman kernel (cf.
e.g., Nehari 1950 \cite{Nehari_1950}) so that this heading is
strongly connected to EP.

$\bullet$ PP=Plateau problem style methods (for RMT, this starts
with the observation of Douglas 1931
\cite{Douglas_1931-Solution}). This strongly allied to DP, albeit
some distinction is useful to keep in mind just for cataloguing
purposes.

\noindent$\bigstar$ (Qualitative) topological methods:

$\bullet$ the {\it continuity method}, as old as Schl\"afli,
(as Koebe notices somewhere) is involved in the accessory
parameters of Schwarz-Christoffel, in Klein-Poincar\'e's
uniformization through automorphic functions, Brouwer
(invariance of the domain), Koebe, etc., e.g. Golusin 1952/57
\cite{Golusin_1952/57})

$\bullet$ Brouwer topological degree and the allied
surjectivity criterion (cf. e.g., Mizumoto 1960
\cite{Mizumoto_1960}, Gabard 2006 \cite{Gabard_2006}). Here
the idea is that there is some topological stability of the
embedding of a curve into its Jacobian via the Abel mapping in
the sense that its homological feature are unsensitive to
variation of the complex (analytic) structure (moduli), and
this enables one to draw universal statement by purely
topological considerations.

Finally we have attempted to manufacture a genealogy map  showing
the affiliation between the authors. The picture turned out to be
so large that TeX prefers reject it at the very end of the file.

%%%\iffalse
\begin{figure}[h]\Figskip
\vskip-1.2cm\penalty0
%\centering
    \hskip-4.2cm\penalty0\epsfig{figure=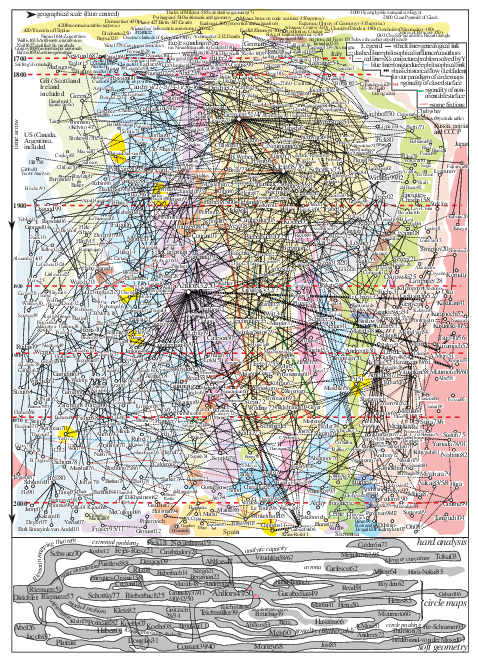,width=205mm}
\captionskipAG
  \caption{\label{Geneal:fig}%
  The authors involved in this bibliography (and genealogy)}
\figskip
\end{figure}
%%%%\fi

\def\fakeref{\hskip-18pt [$\bigstar$]\;\;}

{\small

[15.10.12]  When I reached 884 references, I unfortunatel met the
so-called ``TeX capacity exceeded, sorry.'' obstruction (cf.
Knuth's ``The TeX Book'', p.\,300 for more details). Thus I had to
deactivate some references which are not used for cross-citation,
albeit they clearly belong to our topic.
%Those references are marked by the
%symbol ``\fakeref'' instead of carrying a genuine integer-valued
%label.
%
[16.10.12] This problem was ultimately solved by
%Sir
my advisor Daniel Coray, to whom I express my deepest gratitude
for enlarging the TeX capacity of my compilator.

}

%BIB

{\small

}

\bigskip
%\hspace{+5mm} % To get a little bit of space between the figures
%%%%{\footnotesize
\begin{minipage}[b]{0.6\linewidth} Alexandre
Gabard

Universit\'e de Gen\`eve

Section de Math\'ematiques

2-4 rue du Li\`evre, CP 64

CH-1211 Gen\`eve 4

Switzerland

alexandregabard@hotmail.com
\end{minipage}

%%%%}
%\hspace{-25mm}


\begin{thebibliography}{9}
\bibitem{Gabard-Gauld_2011-Dynamics}
A. Gabard, D. Gauld, {\it Dynamics of non-metric manifolds}, arXiv
2011. [$\spadesuit$ this is just cited for the proof of the
implication: simply-connected $\Rightarrow$ schlichtartig
$\Rightarrow$ orientable, which may be reduced to the five
lemma]~\hfill {\bf $\heartsuit$2}
\end{thebibliography}
\end{document}